\documentclass[11pt]{book}
\usepackage{amsmath}

\def\[#1\]{\begin{equation}#1\end{equation}}
\usepackage{amssymb,amsthm,amsxtra}
\usepackage{amscd}
\usepackage[margin=1in]{geometry}
\usepackage{tikz}
\usetikzlibrary{calc}
\usetikzlibrary{positioning}
\usetikzlibrary{patterns}
\usepackage{hyperref}

\newtheorem{thm}{Theorem}
\newtheorem{cor}[thm]{Corollary}
\newtheorem{lem}[thm]{Lemma}
\newtheorem{prop}[thm]{Proposition}

\theoremstyle{remark}
\newtheorem*{rem}{Remark}
\newtheorem{rems}{Remark}[thm]

\newtheorem{eg}{Example}
\theoremstyle{definition}
\newtheorem{defn}{Definition}

\numberwithin{equation}{chapter}
\numberwithin{thm}{chapter}
\numberwithin{eg}{chapter}
\numberwithin{defn}{chapter}

\newcommand{\sO}{\mathcal O}
\newcommand{\dL}{{\bf L}}
\newcommand{\dR}{{\bf R}}
\newcommand{\ratto}{\dashrightarrow}
\newcommand{\whQ}{{\widehat{Q}}}
\newcommand{\frakm}{{\mathfrak{m}}}

\newcommand{\A}{\mathbb A}
\newcommand{\C}{\mathbb C}
\newcommand{\F}{\mathbb F}
\newcommand{\G}{\mathbb G}
\newcommand{\N}{\mathbb N}
\renewcommand{\P}{\mathbb P}
\newcommand{\Q}{\mathbb Q}
\newcommand{\Z}{\mathbb Z}

\newcommand{\kII}{I\!I}
\newcommand{\kIII}{I\!I\!I}

\DeclareMathOperator{\id}{id}

\DeclareMathOperator{\Hom}{Hom}
\DeclareMathOperator{\Ext}{Ext}
\DeclareMathOperator{\End}{End}
\DeclareMathOperator{\Aut}{Aut}
\DeclareMathOperator{\Isom}{Isom}

\DeclareMathOperator{\Sym}{Sym}
\DeclareMathOperator{\Nm}{Nm}

\DeclareMathOperator{\num}{num}
\DeclareMathOperator{\pt}{pt}
\DeclareMathOperator{\NS}{NS}

\DeclareMathOperator{\Spec}{Spec}
\DeclareMathOperator{\Proj}{Proj}

\DeclareMathOperator{\Pic}{Pic}
\DeclareMathOperator{\Gal}{Gal}
\DeclareMathOperator{\Hilb}{Hilb}
\DeclareMathOperator{\Quot}{Quot}
\DeclareMathOperator{\Alb}{Alb}
\DeclareMathOperator{\Grass}{Grass}

\DeclareMathOperator{\qcoh}{qcoh}
\DeclareMathOperator{\coh}{coh}
\DeclareMathOperator{\perf}{perf}

\DeclareMathOperator{\red}{red}
\DeclareMathOperator{\nr}{nr}

\DeclareMathOperator{\rad}{rad}
\DeclareMathOperator{\gr}{gr}

\DeclareMathOperator{\GL}{GL}
\DeclareMathOperator{\SL}{SL}
\DeclareMathOperator{\PGL}{PGL}
\DeclareMathOperator{\Mat}{Mat}
\DeclareMathOperator{\Map}{Map}
\DeclareMathOperator{\Tr}{Tr}
\DeclareMathOperator{\Br}{Br}

\DeclareMathOperator{\so}{{\mathfrak{so}}}
\DeclareMathOperator{\ad}{ad}

\DeclareMathOperator{\ch}{char}

\DeclareMathOperator{\rank}{rank}
\DeclareMathOperator{\coker}{coker}
\DeclareMathOperator{\im}{im}
\DeclareMathOperator{\ord}{ord}
\DeclareMathOperator{\cone}{cone}

\DeclareMathOperator{\supp}{supp}
\DeclareMathOperator{\Fitt}{Fitt}

\DeclareMathOperator{\Ai}{Ai}

\DeclareMathOperator{\Spl}{{\mathcal S}\!{\it pl}}
\DeclareMathOperator{\Tor}{{\mathcal T}\!{\it or}}
\DeclareMathOperator{\Vect}{{\mathcal V}\!{\it ect}}
\DeclareMathOperator{\sExt}{{\mathcal E}\!{\it xt}}
\DeclareMathOperator{\sHom}{{\mathcal H}\!{\it om}}
\DeclareMathOperator{\sEnd}{{\mathcal E}\!{\it nd}}
\DeclareMathOperator{\Irr}{{\mathcal I}\!{\it rr}}

\let\div\relax
\DeclareMathOperator{\div}{div}

\newcommand{\Mer}{\text{\it ${\mathcal M}\!$er}}

\begin{document}

\title{The (noncommutative) geometry of difference equations\footnote{working title}}
  \author{
Eric M. Rains\\Department of Mathematics, California
  Institute of Technology}

\date{April 22, 2025}
\frontmatter
\maketitle

\chapter{Introduction}
The motivating problem of the present monograph is to understand moduli
spaces of difference (or differential) equations subject to local
constraints telling us about the singularities of the equations, with an
eye to both understanding their basic properties (dimension, rationality,
irreducibility) and natural (and modular!) isomorphisms between them.  It
turns out that this problem is very closely related to a much more
classical problem in algebraic geometry, namely to understand moduli spaces
of {\em sheaves} on smooth projective surfaces, with the one (albeit very
significant) caveat that we need to deform everything to {\em
noncommutative} geometry.  The benefit of this approach is that we gain a
great deal of insight, not to mentions proofs to adapt, from commutative
algebraic geometry; the cost is that we will need to develop the theory of
noncommutative surfaces well beyond what is known in the literature.
(Though this is also a benefit from another perspective: we were led to
discover a number of crucial properties of noncommutative surfaces in the
process of addressing the original problem, without which we would not have
been able to develop the theory nearly as far!)

Although the approach taken here is essentially algebraic (or
algebro-geometric) in nature, the origins of the motivating problem lie in
analysis, and thus we begin by a brief discussion of those aspects.  The
story effectively begins with Riemann's observation that any second-order
Fuchsian equation with three regular singular points is equivalent (up to
changes of variables and rescaling by solutions of first-order equations)
to the equation
\[
z(z-1)f''(z)
+
((a+b+1)z-c)f'(z)
+
abf(z)
=
0,
\]
satisfied by Gauss' hypergeometric function ${}_2F_1(a,b;c;z)$.  In
particular, any such equation is {\em rigid}: it is uniquely determined by
its local behavior near the singular points.  This rigidity phenomenon
persists for other hypergeometric functions, in particular for
${}_nF_{n-1}$ (in which the recurrence for the coefficients is generalized)
and for the Jordan-Pochhammer function (in which Euler's integral
representation is generalized), suggesting that one should investigate
rigid equations more generally.  See for instance
\cite{KatzNM:1996}, which gave an algorithm for
determining when the local system (representation of the fundamental group)
associated to a Fuchsian equation is similarly rigid.

If one allows four singular points, a new structure arises: there is a
two-dimensional family of such equations, as well as a two-dimensional
family of local systems.  Since the fundamental group of
$\C\setminus \{0,1,\lambda\}$ is locally independent of $\lambda$, this
leads to a new operation one can perform: take a Fuchsian equation, apply
the Riemann-Hilbert correspondence to obtain its monodromy representation,
then change $\lambda$ and map back.  This yields an analytic morphism
between the moduli spaces of equations which for infinitesimal changes in
$\lambda$ translates to a certain nonlinear differential equation, the
Painlev\'e VI equation.  (This is the most complicated of a family of
eight\footnote{This was originally viewed as a hierarchy of six equations,
but in the modern version of the classification, the Painlev\'e III
equation is split into three cases, labelled $D_6$, $D_7$, and $D_8$ in
Figure \ref{fig:sakai_good} below; the other Painlev\'e equations
correspond to $D_4$ (PVI), $D_5$ (PV), $E_6$ (PIV), $E_7$ (PII), and $E_8$
(PI) in that diagram.}  second-order nonlinear differential equations
studied by Painlev\'e and his students.)

To give a more explicit example, it will be helpful to degenerate the
situation.  Consider the family of second-order equations that are regular
on $\C$ and have local behavior
\[
v'(z) = \begin{pmatrix} -z^2-t/2+\beta/z &0\\
                          0 & z^2+t/2-\beta/z\end{pmatrix} v(z)
\]
near $\infty$ (i.e., there is a change of basis holomorphic near $\infty$
putting it in that form).   Up to global changes of basis, this is a
2-dimensional family, and can in fact be rationally parametrized, with one
form of a universal equation having the form
\[
v_z
=
\begin{pmatrix}
   -z^2-t/2+u & 2 u(z-x)-2\beta\\
   z+x&z^2+t/2-u
\end{pmatrix}
v,
\]
where $x$ and $u$ are coordinates on the moduli space.  The analogue of the
monodromy-preserving deformation of $\lambda$ is that one can extend
solutions along deformations of $t$; the equation
\[
v_t =
\begin{pmatrix}
  -(z-x)/2 & u \\
    1/2 & (z-x)/2
\end{pmatrix}
v
\]
is compatible with the previous equation as long as $x$ and $u$ satisfy the
nonlinear system
\[
x_t = u - x^2-t/2 \qquad u_t = 2 ux + \beta,
\]
one of the standard forms of the Painlev\'e II equation.  (Note that both
equations are rational in $z$, but {\em analytic} in $t$, since $x$ and $u$
are not in general algebraic.  However, the calculation of the
compatibility relations {\em is} algebraic\dots)

This ``Lax pair'' representation of Painlev\'e II is not unique, however:
the pair
\[
v_z
=
\begin{pmatrix}
  x-\beta/2z & -1/2-u/2z\\
  z+2x^2-u+t & -x+\beta/2z
\end{pmatrix}
v
\qquad
v_t
=
\begin{pmatrix}
  x & -1/2\\
  z & -x
\end{pmatrix}
v
\]
has precisely the same compatibility conditions, but now instead of the
equation in $z$ having a single non-Fuchsian singular point, has two
singular points with one (at 0) Fuchsian.

This description of PII, though not the original description (Painlev\'e
originally discovered it as one of relatively few nonlinear equations for
which the only singularities that depend on the initial conditions are
poles), is more fundamental in many respects; in particular, it explains
the appearance of the equation in many natural asymptotic questions (in
particular in the Tracy-Widom distribution of random matrix
theory \cite{TracyCA/WidomH:1994}).  (The solution of the Riemann-Hilbert
problem that gives rise to the relevant solution of PII itself satisfies a
variant of the above Lax pairs, see \cite{png,iseq2}, where the solution to
the linear problem also appears in the asymptotics of certain phase
transitions.)

Of course, since we obtain multiple such descriptions, this leads to the
natural question of whether there are any others.  This is, in fact the
case, for two different reasons.  The first is that the two Lax pairs above
are actually related by the Fourier-Laplace transform: if one rescales both
equations by suitable first-order equations, then they are related under
the formal transformation $(D_z,z)\mapsto (-z,D_z)$ (and suitable changes
of basis).  If one applies the transformation {\em without} first
rescaling, then the result will be a higher-order differential equation,
and iterating gives arbitrarily high-order Lax pairs for PII.  (It is not
quite trivial that the equation in $t$ survives this process, but this
essentially follows from our general theory below.)  Of course, this tells
us that we should instead be asking for Lax pairs which are ``minimal'',
i.e., cannot be simplified by a formal integral transform, but even there
it turns out there are additional examples (which even for PII have yet to
be worked out explicitly, apart from the fourth-order example given in
Appendix \ref{app:PII} below.).

The Painlev\'e equations turn out to have nice (``integrable''\footnote{The
word ``integrable'' has many slightly different definitions in the discrete
setting, with the most relevant for our purposes being (a) preserves a
Lagrangian fibration, (b) has iterates with subexponential degree growth,
and (c) is the compatibility relation of a Lax pair.  For the discrete
isomonodromy transformations consider below, the first only applies in the
autonomous case (our ``commutative relaxation''), where it is
straightforward, the second is shown to hold in many cases in Section 18.7,
while the third is left as an exercise for the reader.})  discretizations.
We can see one of these already in our PII example.  In the second Lax
pair, note that changing $\beta$ by 2 does not change the local behavior of
the solutions around 0 and $\infty$, and thus it is reasonable to hope that
the two fundamental matrices should be related by a globally meromorphic
(i.e., rational) matrix.  This is, indeed, the case, and in fact applies
for $\beta\mapsto \beta-1$ as long as one is willing to globally scale by
$z^{1/2}$:
\[
\hat{v} = \begin{pmatrix}
      0 & z^{-1/2}/2\\
      -z^{1/2} & (x+\hat{x})z^{-1/2}
      \end{pmatrix}
      v
\]
satisfies the equation with $(\beta-1,\hat{x},\hat{u})$ as long as
\[
x+\hat{x} = \frac{1-\beta}{\hat{u}},
\qquad
u+\hat{u} = 2x^2+t,
\]
a standard form of the {\em discrete} Painlev\'e I equation.  (This is a
nonlinear difference equation in $x$ and $u$, and in a suitable limit as
$\beta\to\infty$ becomes the Painlev\'e I equation.)  Something similar
applies for the first Lax pair.

Note that the Painlev\'e equations are closely related to hypergeometric
equations.  For instance, when $\beta=0$, setting $u=0$ in the second
Lax pair above gives
\[
v_z = \begin{pmatrix} x & -1/2\\
         z+2x^2+t &-x\end{pmatrix}v,
\]
or equivalently $v_1''(z)=2^{-1}(z+t)v_1(z)$, which up to linear change of
variable is nothing other than the Airy equation $f''(z)=zf(z)$.  The
corresponding substitution makes the first Lax pair lower triangular, while
the compatibility condition becomes a Ricatti equation, again related to
the Airy equation, with solution
$(x,u)=(\frac{d}{dt}\log(f(-2^{-1/3}t)),0)$.  Applying the gauge
transformations shifting $\beta$ gives more complicated embeddings of the
Ricatti equation in PII, in which the reducibility of the first Lax pair is
hidden by more and more complicated gauge transformations, while the second
Lax pair has more and more hidden ``apparent'' singularities at $0$.

A similar discrete transformation applies to the hypergeometric equation:
adding $1$ to $a$ does not change the conjugacy classes of the three
generators of monodromy, and thus (since the monodromy is also rigid!) does
not change the monodromy representation, so that one can express the
standard solution ${}_2F_1(a+1,b;c;z)$ in terms of ${}_2F_1(a,b;c;z)$ and
its derivative.  (This is nothing other than one of the standard
``contiguity'' relations.)  This gives rise to a linear difference equation
in $a$ with rational coefficients, and one finds that this equation is {\em
also} rigid in the obvious sense: it is the unique second-order difference
equation with the given finite singularities and behavior near $\infty$.

This can be further generalized by considering $q$-hypergeometric
functions---the $q$-hypergeometric function ${}_2\phi_1(a,b;c;z;q)$
satisfies rigid $q$-difference equations in $a$, $b$, $c$, and $z$---or by
considering other natural families (balanced, well-poised, etc.).  If one
generalizes all the way to {\em elliptic} hypergeometric functions (such as
the generalized beta integrals considered in \cite{dets}), one at first
blush appears to lose rigidity, but this is restored if one notices an
additional symmetry of the situation.  (This symmetry was first {\em
explicitly} noticed in \cite{isomonodromy} for certain equations related to
semiclassical elliptic biorthogonal functions.)  This again can already be
seen for ${}_2F_1$.  If we consider the difference equation satisfied by
${}_2F_1$ under the shift $(a,b)\mapsto (a+1,b-1)$, this at first blush
does not appear to be rigid.  However, the equation has solutions
satisfying the symmetry ${}_2F_1(a,b;c;z)={}_2F_1(b,a;c;z)$, and imposing
the corresponding consistency condition on the difference equation again
makes the equation unique.

If we view this as a difference equation in $t=(a-b)/2$ with $(a+b)/2,c,z$
fixed, then it has the matrix form
\[
v(t+1) = A(t) v(t)
\]
with $A(t)\in \GL_2(k(t))$.  However, since we are interested in symmetric
solutions, we should instead consider the system
\[
v(t+1) = A(t) v(t) \qquad v(-t)=v(t),
\]
at which point we see that $A$ should satisfy the consistency condition
\[
A(-1-t)A(t) = 1.
\]
Note here that there is a subtle issue regarding singularities.  In
general, a singularity of an equation is a local obstruction to having an
invertibly holomorphic fundamental matrix at a point (or on an orbit), and
in the case of a symmetric equation, we are looking for {\em symmetric}
solutions, and thus care about obstructions to having a {\em symmetric}
fundamental matrix.  Thus, for instance, the symmetric difference equation
\[
v(t+1) = - v(z),\qquad v(-t)=v(t),
\]
is singular at $z=-1/2$ since any {\em symmetric} holomorphic solution
near $z=-1/2$ is forced to vanish there.

The condition for an equation to be symmetric can be rephrased as follows.
Without symmetry, we have an action of the infinite cyclic group $\Z$ on
the field $k(t)$ of rational functions (by $t\mapsto t+1$ or $t\mapsto qt$,
or similarly for the field of functions on an elliptic curve), and a
difference equation is nothing other than a 1-cocycle in
$Z^1(\Z;\GL_n(k(t)))$.  The symmetry condition extends the action to one of
the infinite {\em dihedral} group $D_\infty$, and a symmetric equation is a
cocycle in $Z^1(D_\infty;\GL_n(k(t)))$.  A fundamental matrix for a
difference equation is then an explicit trivialization in
$H^1(\Z;\GL_n(\Mer))$ or $H^1(D_\infty;\GL_n(\Mer))$ (i.e., over the field
of meromorphic functions), while the analogue of the discrete isomonodromy
transformation we considered for PII is nothing other than an explicit
gauge equivalence in $H^1(\Z;\GL_n(k(t)))$ or $H^1(D_\infty;\GL_n(k(t)))$.

Since generalizations of hypergeometric functions often satisfy rigid
equations, we would like to understand rigid equations more generally.  In
particular, given a type of equation (differential, difference, symmetric
difference, possibly over a suitable algebraic curve) and a description of
its ``singularities'', we would like to know whether there is a {\em
unique} such equation.  More generally, when the equation is {\em not}
rigid, we would like to understand the moduli space of such equations, and
in particular want to know when there are natural birational maps between
such moduli spaces coming from gauge (i.e., isomonodromy) transformations.

For reasonably simple singularities, it is not too difficult to work with
the moduli spaces infinitesimally (i.e., deformation theory), and thus to
compute expected dimensions.  This would settle the rigidity question, if
it were not for the (a priori quite difficult) problem of determining when
the moduli space is nonempty!  Moreover, the formulas one obtains in such
situations lead to a further mystery, to wit that the expected dimension is
always even.  This suggests that there should be a nondegenerate
alternating pairing on the tangent space, a.k.a., a nondegenerate
$2$-form.  Such a form can be constructed (e.g., in the symmetric elliptic
difference case), leading to the additional question of whether the
resulting 2-form is closed (i.e., whether the moduli space is symplectic).

One difficulty with trying to carry out the above program is a tendency
towards combinatorial explosion: we already have seven types of equation to
consider (differential equations on an arbitrary curve, plus difference
equations on $\G_a$, $\G_m$ or $E$, each of which comes in symmetric
forms), and must then also take into account the various different
possibilities for singularities.  So at the very least, we need to have a
more uniform description for the different types of equations.

Another issue that arises is that if we want to construct a global moduli
space, there is not too much difficulty if we are willing to settle for a
stack (or even an algebraic space), but if we want to obtain an actual Lax
pair, we ideally want the moduli space to be {\em rational} with an
explicit rational parametrization {\em and} a universal family, as we had
in our PII cases.  Moreover, if we want to take a somewhat more classical
approach and have a (quasi-)projective moduli space, we run into the issue
that the usual geometric invariant theory approach tends to change the
moduli problem in unnatural ways, and it can take a great deal of work to
figure out what choices to make in the construction to get a {\em modular}
description of the quotient.

In the Painlev\'e cases, there is, in fact, a natural projective version of
the moduli space.  This was constructed directly for various continuous and
discrete Painlev\'e equations, and eventually systematized by
Sakai \cite{SakaiH:2001}, who observed that every known example was a
rational surface with $K^2=0$ and without $-d$-curves for $d>2$, and
constructed a discrete Painlev\'e equation associated to every such
surface.  Such a surface generically has a unique anticanonical curve
(section of $|-K|$), which is generically smooth of genus 1, and the
associated nonlinear equation is known as the ``elliptic Painlev\'e
equation''.  Unfortunately, Sakai's description, while quite pleasing
aesthetically, does not have any particular relation to the Lax pair
descriptions available for degenerate cases.  (Even in the degenerate cases
for which Lax pairs were known, Sakai's description gives a nice
compactification of the moduli space of equations, but does not come with
any a priori modular description of the compactification.)  An elliptic
analogue of the hypergeometric solutions to Painlev\'e was constructed
in \cite{isomonodromy}, suggesting that indeed the surfaces considered by
Sakai should be interpretable as moduli spaces of symmetric elliptic
difference equations; corresponding Lax pairs were constructed
in \cite{YamadaY:2009} and extended to second-order equations with more
singularities in \cite{ellGarnier}.

Of course, we have already seen that there can be multiple Lax pairs for a
given (discrete) Painlev\'e equation, and thus we would like to understand
how Sakai's classification of Painlev\'e equations splits into a
classification of Lax pairs.  This is a special case of a more general set
of problems: for a given explicit Lax pair, we would like to understand its
generalizations and degenerations.  It is reasonably straightforward to
classify degenerations from elliptic equations to $q$-difference equations
(see the approach of \cite{limits3}), but further degenerations again run into
combinatorial explosion (there are simply too many ways one might take
scaled limits).

The key additional insight to solving the above problems is that in the
differential case, one may work with representations of the algebra of
differential operators (i.e., $D$-modules).  The algebra of differential
operators on a curve may be viewed as a flat deformation of the cotangent
bundle over the curve, and thus in a sense differential equations may be
thought of as sheaves on such a deformation.  This idea was further
bolstered by the construction in \cite{P2Painleve} of elliptic
Painlev\'e via a moduli space of sheaves on a noncommutative $\P^2$.

The difficulty is then to find a suitable construction of noncommutative
surfaces in terms of difference operators.  One of the earliest
constructions in noncommutative geometry was Sklyanin's noncommutative
deformation of $\P^3$ \cite{SklyaninEK:1982,SklyaninEK:1983}, in the form
of a noncommutative graded algebra with the same Hilbert series as $\P^3$.
(The claim about the Hilbert series was conjectural at the time, but proved
in \cite{ArtinM/TateJ/VandenBerghM:1990}.)  In \cite{RosengrenH:2004},
Rosengren observed that the finite-dimensional representations Sklyanin
considered could be interpreted and generalized to representations as
elliptic difference operators, or more precisely {\em symmetric} elliptic
difference operators.  Moreover, each such representation annihilates a
central element of degree 2, and thus this is more naturally a
representation of a noncommutative {\em quadric surface}.  (This
construction is, in fact, reversible, as we will see below.)  In
particular, this gives rise to a module $\Mer$ over the corresponding
graded algebra such that for any other module $M$, $\Hom(M,\Mer)$ is given
by the solutions of a suitable symmetric elliptic difference equation.
Since each degree 2 central element of Sklyanin's algebra is annihilated in
{\em two} such representations, each noncommutative quadric has two such
representations.  Moreover, the commutative quadrics where the two
representations coincide are nothing other than the {\em singular}
quadrics, and thus in the commutative case, we are choosing a ruling of the
quadric.  This suggests more generally that the object we should be looking
at is a noncommutative ruled surface.

Although this is hopeful in one significant respect (we certainly know how
to construct moduli spaces of sheaves on {\em commutative} surfaces!), one
significant remaining issue is that noncommutative algebraic geometry is a
relatively young field (even something as fundamental as ample bundles on
noncommutative ruled surfaces had not been studied!), and thus not only
does one have the problem of proving things, but the more fundamental
problem of figuring out what one {\em should} prove.  This suggests that we
should first investigate modified versions of the questions that relate to
{\em commutative} surfaces.  For the differential case, such a commutative
relaxation is given by the notion of {\em Higgs bundles} (replacing a
connection on a vector bundle $V$ by a morphism $V\to V\otimes \omega_C$),
and there are similar relaxations in the difference cases.

We will see below that such relaxed equations translate naturally to
sheaves on suitable (commutative, projective) ruled surfaces.  Moreover,
the singularities of such relaxed equations can be read off of the way in
which the sheaf meets a certain curve on the ruled surface, which turns out
to always be a section of the anticanonical linear system.  That a choice
of anticanonical curve on a smooth projective surface determines (up to an
overall scalar) a Poisson structure on the surface is quite promising,
since Poisson structures are nothing other than infinitesimal
noncommutative deformations!  Moreover, results
of \cite{TyurinAN:1988,BottacinF:1995,HurtubiseJC/MarkmanE:2002b} tell us
that the moduli space of sheaves on a Poisson surface comes with a natural
Poisson structure, and \cite{BottacinF:2000} that the fibers of restriction
to the anticanonical curve are symplectic! (Note that the fibers being
symplectic was only shown for stable vector bundles, and all results were
only shown to hold in characteristic 0; we will remove both restrictions
below.)

When the anticanonical curve becomes singular, this description of
singularities can become fairly complicated (there are many
finite-dimensional sheaves on the anticanonical curve, even if we include
the complicated condition that they be the restrictions of $1$-dimensional
sheaves on the ambient surface), but there is a trick for simplifying the
problem.  If the sheaf corresponding to our relaxed equation corresponds to
a line bundle on a smooth curve (which is the generic behavior in most
cases), then we can blow up any point of intersection with the
anticanonical curve and lift the map from the curve and thus the sheaf.
Each such blowup reduces the intersection number, and thus eventually we
end up with a sheaf that is {\em disjoint} from the anticanonical curve.
This procedure encodes the original singularity data in two pieces: the
blown up surface (telling us which irreducible singularities appear) and
the first Chern class of the lifted sheaf (telling us their
multiplicities).  In particular, this turns the quite complicated question
of classifying singularities into a much more familiar problem: classifying
(anticanonical) projective surfaces.  (There is a cost, in that some
singularities cannot be resolved in this way, but it turns out that such
singularities can always be simplified by a canonical gauge
transformation.)  Moreover, the blown up surface has additional line
bundles, and it turns out that the process of lifting a sheaf, twisting by
an exceptional line bundle, and then mapping back down is a (canonical)
gauge transformation!  (Note that although line bundles on commutative
surfaces are naturally bimodules over that surface, so twisting gives a new
sheaf on the surface, this fails in the noncommutative case: line bundles
are still bimodules, but the surfaces on the left and right are different!
Of course, this is precisely what we want, given that discrete isomonodromy
transformations are {\em supposed} to move the
singularities\dots)

Part I below studies this commutative relaxation in depth.  (Here, despite
the original analytic motivation being over $\C$, we mostly work over a
general algebraically closed field $k$, eventually generalizing to more
general Noetherian schemes over $\Z$.  This is the natural level of
generality when we switch to noncommutative geometry, and indeed there is
one result (Proposition \ref{prop:irreds_exist} on the existence of
irreducible sheaves) for which our proof in characteristic 0 involves
lifting from finite characteristic.)  Not only does this give significant
insights into what should hold for the noncommutative version, but in a
number of cases noncommutative results {\em reduce} to their commutative
results, and in many others, the commutative proof extends {\it mutatis
mutandis}.  Chapter \ref{chap:sheaves_from_eq_comm} describes how to encode
relaxed equations in certain ruled surfaces, starting with the symmetric
elliptic case\footnote{This corresponds to rational surfaces with smooth
anticanonical curves, so from a geometric perspective is actually the
simplest case!}, then extending to the various degenerations, and discusses
the analogue of the isomonodromy transformations in the commutative
relaxation.  Chapter \ref{chap:classify}, after briefly discussing Poisson
structures from an algebraic perspective, classifies Poisson surfaces, and
thus in particular shows that every nontrivially Poisson surface is
birational (in a way that respects the Poisson structure) to one of the
surfaces arising in Chapter \ref{chap:sheaves_from_eq_comm}.\footnote{This
is not quite correct, as there is an unusual case appearing only in
characteristic 2.}  In Chapter \ref{chap:*!}, we consider the procedure of
lifting to blowups.  The approach via line bundles on smooth curves has two
major deficits: not only are there sheaves in the moduli space not of this
form, but in the noncommutative case, this description essentially never
holds.  It turns out, however, that there is a more functorial description
of a certain ``minimal'' lift to the blowup, which will extend easily to
the noncommutative setting.  In Chapter \ref{chap:poisson_moduli}, we
discuss the result that the moduli space of sheaves on a Poisson surface
inherits a Poisson structure, show that the minimal lift of
Chapter \ref{chap:*!} is Poisson, and consider the question of when a
symplectic leaf is a single point, showing that such rigid $1$-dimensional
sheaves resolve to line bundles on $-2$-curves.

As mentioned above, the translation between sheaves and equations depends
on a choice of ruling.  Surfaces ruled over a curve of positive genus have
a unique ruling, so this is only an issue in the rational case.  Moreover,
since a $-2$-curve on a higher-genus ruled surface must be contained in a
fiber, we only get {\em interesting} rigid sheaves in the rational case.
With this in mind, the rest of Part I is primarily focused on the rational
case.  A representation of a rational surface as an iterated blowup of a
Hirzebruch surface (with choice of ruling as appropriate) induces a basis
of its group of line bundles, and a number of important questions about
divisors can be simplified by using the freedom to choose the blowdown
structure.  Chapter \ref{chap:divs_on_rat} studies this freedom, showing
that there is a (partial, but essentially transitive) action of a Coxeter
group of type $E_n$ (extending the sequence
$A_4,D_5,E_6,E_7,E_8,\tilde{E}_8$) on such blowdown structures, and that
one can adapt some of the standard algorithms of Coxeter theory to
determine whether divisors are nef, effective, or generically integral.
(In particular, one can test whether a given divisor class is represented
by a $-2$-curve, or for a generically integral class find rulings that
minimize the order of the corresponding equation.)  In
Chapter \ref{chap:moduli_of_comm_surfs}, we study problems related to
classification of surfaces: showing that the moduli stacks of rational (or
anticanonical rational) surfaces with blowdown structures are nice, and
studying how the combinatorics of the anticanonical curve splits up the
moduli stack.  (In particular, we give a more or less explicit dictionary
between surfaces and types of singularities.)
Chapter \ref{chap:moduli_of_comm_sheaves} then studies moduli spaces of
$1$-dimensional sheaves, in particular addressing questions such as
rationality or the existence of a universal family, as well as giving an
application: computing the relative Jacobians of rational surfaces.
Finally, in Chapter \ref{chap:moduli_of_eqs_comm}, we bring everything back
to equations, studying the symmetric elliptic case in some detail and
explaining how the various twists by line bundles act on relaxed equations
(as well as natural ``guesses'' for restoring the $q$ parameter in those
descriptions).  We also discuss several examples of using these ideas to
classify degenerations (e.g., classifying hypergeometric functions in the
hierarchy containing ${}_2F_1$, classifying minimal Lax pairs with $2$- or
$4$-dimensional moduli spaces).

For those readers interested in explicit Lax pairs, one should note that
there are essentially no such examples below (the proofs of rationality and
existence of universal families are not usefully effective), but there are
a great number of new {\em existence} results for Lax pairs.  (In
particular, each (discrete) Painlev\'e equation has a family of Lax pairs
parametrized by rational numbers (in fact $\P^1(\Q)$, with $\infty$
corresponding to Sakai's construction), none of which are related by
generalized Fourier transforms.  None of the cases with nontrivial
denominator have been worked out explicitly!)  Since writing down an
explicit Lax pair is in itself a proof of existence, any result along those
lines will be independent of any of the noncommutative geometry results
discussed below.  (See, e.g.,
\cite{ellGarnier}, which constructs second-order Lax pairs of symmetric
elliptic difference type.)

On the other hand, as we already mentioned, the field of noncommutative
algebraic geometry is new enough that most of what we needed to know was
not already in the literature.  The early approach to noncommutative
geometry via flat deformations of graded algebras turns out to be
insufficient for our cases (in the generic case, they only exist for
noncommutative del Pezzo surfaces, or more generally for surfaces on which
some linear combination of components of the anticanonical curve is ample),
and thus we are left with the more modern perspective of identifying a
noncommutative surface with its category of (quasi-)coherent sheaves.  (Of
course, all of our questions are {\em about} sheaves, so this is not an
issue per se, though it does mean that some standard commutative arguments
fail miserably, with a major one being that of taking hyperplane sections
relative to a very ample divisor.)

Although there were constructions of noncommutative projective
planes \cite{ArtinM/TateJ/VandenBerghM:1990,BondalAI/PolishchukAE:1993},
noncommutative ruled surfaces \cite{VandenBerghM:2012}, and even
noncommutative blowups \cite{VandenBerghM:1998}, some of the most basic
questions were still open: e.g., whether blowups in distinct points
commute.  Moreover, as we have already mentioned, the existing
constructions were not particularly well adapted to difference(tial)
operators.  In the symmetric elliptic case, one can extend Rosengren's
observation to give a completely explicit construction: first, of a flat
deformation of the category of all line bundles on the
surface\footnote{This is not quite correct: even in the commutative case,
flatness requires replacing some of the $\Hom$ spaces by proper subspaces
on those fibers where they would otherwise jump.}, and then of a suitable
quotient of the corresponding category of modules (a priori depending on a
choice of ample line bundle, but in fact independent of the choices),
giving a deformation of the category of coherent sheaves on the original
commutative surface.  This again runs into combinatorial explosion if we
try to degenerate, as well as a more fundamental issue that we cannot
expect to have {\em quite} so flat a deformation.  However, it {\em does}
have an interesting generalization to multivariate
operators \cite{elldaha}, and the simplest version of {\em that}
construction turns out to generalize considerably and give us an alternate
version of Van den Bergh's construction of noncommutative $\P^1$-bundles,
in which various morphisms are explicitly given as difference-reflection
(for the affine Weyl group $W(\tilde{C}_1)$, which we encountered above as
the infinite dihedral group $D_\infty$) operators, reducing to difference
operators under a suitable Morita equivalence.

This leads to an interesting (near) dichotomy.  In nearly every case of the
construction, the corresponding algebra has a large center (suitably
defined), and thus the category of sheaves is equivalent to a category of
modules over a certain sheaf of algebras on a {\em commutative} surface.
Indeed, for most components of the relevant moduli stack, every surface in
the family has this form.  Similarly, fully commutative surfaces only
appear on relatively few components of the moduli stack, with the curious
fact being that there is essentially no overlap: the only overlap is a
certain family which is a quaternion algebra in characteristic 0 but
commutative in characteristic 2.  (Moreover, those commutative surfaces are
pathological in other ways, e.g., by having nonreduced Picard group.)  With
this in mind, the surfaces arising from the general construction will be
called ``quasi-ruled'', while those coming from components with truly
noncommutative surfaces will be called ``ruled''.  The classification of
noncommutative quasi-ruled surfaces appears hopeless in general, but the
moduli stack of noncommutative ruled surfaces is fibered in a natural way
over the moduli stack of anticanonical ruled surfaces.

Chapter \ref{chap:daha} discusses this construction in detail, while
Chapter \ref{chap:semicomm} discusses the ``semicommutative'' case, both to
show that we can safely ignore quasi-ruled surfaces, and as a useful tool
in the ruled case (which still has a dense set of semicommutative
instances!).  In particular, we show that if a certain automorphism has
finite order, then the corresponding (quasi-)ruled surface corresponds to a
maximal order on a smooth projective ruled surface (i.e., a sheaf of
algebras such that the generic fiber is a central simple algebra, and
maximal among such sheaves of algebras), and that the maximal order
condition persists under blowups.  (This seems to have been a folklore
result under certain tameness conditions, but those conditions are violated
in many cases of interest to us.)

Although this construction, together with Van den Bergh's construction of
blowing up points on surfaces, is in principle enough to give us all the
noncommutative surfaces we need, there is a fundamental issue: various
natural isomorphisms in the commutative case were not known in the
noncommutative case.  (I.e., that $\P^1\times \P^1$ is a ruled surface in
two different ways, that a one-point blowup of $\P^2$ is a ruled surface,
that a one-point blowup of a ruled surface is a one-point blowup of a
different ruled surface, and that blowups in distinct points commute.)  As
usual, a direct approach would lead to combinatorial explosion; even
something as simple as swapping the rulings on $\P^1\times \P^1$ has 16
different cases to consider.  (Since the resulting formal integral
transformations are of intrinsic interest, we discuss all 16 cases in
Appendix \ref{chap:fourier}.)  Luckily, it turns out that each of the
constructions we are using has a much simpler description that is
particularly well adapted to proving such isomorphisms.

The key observation comes from the elliptic case \cite{generic}, where it
was observed that the {\em derived} category of sheaves has a particularly
simple description.  Indeed, if $X$ is a projective rational surface, then
$\sO_X$ is an exceptional object, and thus there is a semiorthogonal
decomposition $(\sO_X^\perp,\sO_X)$.  One can in general reconstruct the
derived category from the semiorthogonal decomposition, or more precisely
from the two components together with a functor specifying the maps between
them.\footnote{To be precise, this is the derived {\em dg}-category rather
than the more familiar triangulated category, so the gluing functor takes
values in complexes.}  In the case of an anticanonical surface, the ``maps
between them'' turn out to factor through (derived) restriction to the
anticanonical curve, and the noncommutative deformation simply twists this
by a line bundle of degree 0 on that curve.  (One immediate consequence of
this approach is that Serre duality is essentially trivial!)

For many purposes, this would be enough, but in fact we can go further.  To
recover an abelian category from its derived category, it suffices to
specify a $t$-structure on the latter, and it turns out we can give an
inductive description of the $t$-structure for noncommutative
surfaces.\footnote{Unfortunately, the only proof we have that the
dg-category is the derived category of the heart of the $t$-structure
involves using the existing constructions.}  So for each of the
isomorphisms we want to prove, we reduce to (a) constructing a derived
equivalence and (b) checking that the $t$-structure is preserved.

In Chapter \ref{chap:derived}, we discuss this construction in detail, and
in particular discuss when blowups commute, and what happens when they do
not.  Chapter \ref{chap:derived_eqs} then uses the simple description(s) of
the derived category to construct interesting derived equivalences.  The
most fundamental equivalence is an analogue of Cohen-Macaulay duality
(related to the natural duality on equations), though the more {\em
interesting} equivalences are certain deformations of the derived
autoequivalences of commutative elliptic surfaces.  (As part of showing
that duality has cohomological dimension 2, Chapter \ref{chap:derived_eqs}
also discusses how line bundles extend to noncommutative surfaces, and
especially how {\em twisting} by line bundles extends.)
Chapter \ref{chap:families} uses the derived description to give a
relatively simple definition of {\em families} of noncommutative surfaces,
not only in cases with blowdown structures (``split'' cases), but in
somewhat greater generality.  One significant result of this section is
that any coherent sheaf on a family over a Noetherian base is itself
Noetherian.  (This generalizes the usual ``strongly Noetherian'' condition
to allow the surface to vary.)

The description of the derived category makes it trivial to determine the
Grothendieck group of our surfaces, and Chapter \ref{chap:dimens} uses this
to define useful notions such as the dimension of a sheaf (essentially the
degree of the Hilbert polynomial), the N\'eron-Severi lattice
(1-dimensional classes in $K_0$ modulo $0$-dimensional classes), and
effective divisors, as well as showing that our noncommutative surfaces are
irreducible in the sense that any map of line bundles is injective.  (This
chapter also proves most of the remaining conditions required
by \cite{ChanD/NymanA:2013} for the surface to be a ``noncommutative smooth
proper $2$-fold'', with the sole exception being the properties of $\Quot$
schemes shown in Chapter \ref{chap:quot}.) Chapter \ref{chap:birat_noncomm}
then uses this to finish off the proofs that the various isomorphisms from
birational geometry indeed preserve the t-structures, as well as showing
(by an analogue of the arguments of Chapter \ref{chap:divs_on_rat}) that an
analogue of Castelnuovo's criterion holds, and thus we do not need a
separate construction of blowing down.  (To be precise: if a noncommutative
surface is related to a noncommutative plane or noncommutative quasi-ruled
surface by a sequence of blowups and inverses of blowups, then it {\em is}
a noncommutative plane or noncommutative quasi-ruled surface!)

Chapter \ref{chap:effnef} discusses the effective monoid and nef cone (and
gives an algorithm for computing them in the ruled case), and shows that
divisor classes in the interior of the nef cone do, indeed, behave like
ample divisors.  In particular, we can use them to define Hilbert
polynomials just as in the commutative case, enabling us to state (and
prove!) in Chapter \ref{chap:quot} that the analogues of the schemes
$\Quot(M,P)$ remain projective in the noncommutative case.  (This is
surprisingly delicate, mainly because the usual commutative proof uses the
flattening stratification, which is much harder to construct in the
noncommutative case; indeed, our eventual proof that every sheaf has a nice
flattening stratification uses projectivity of the Quot scheme!)  Part of
the argument requires that we understand torsion-free sheaves of rank 1,
and Chapter \ref{chap:quot} also shows that those sheaves are classified by
projective moduli spaces (analogous to the Hilbert scheme of points).

In Chapter \ref{chap:moduli_sheaves_noncomm}, we turn our attention to
general moduli spaces of sheaves.  If we ignore questions of
(semi)stability, then one can obtain quite strong results on the existence
of Poisson structures.  Since the best arguments use derived algebraic
geometry, we defer the discussion to \cite{nc-lagrangian}, and only cite
the result in Section \ref{sec:moduli_simple_noncomm}: in characteristic 0,
the moduli stack of simple objects in the derived category is a gerbe over
a Poisson algebraic space.  (For {\em sheaves}, one can then use a density
argument to extend to finite characteristic.)  The argument also shows that
any {\em smooth} fiber of derived restriction to the anticanonical curve is
symplectic, and thus we spend some effort in
Section \ref{sec:leaves_are_smooth} showing that such fibers are almost
always smooth.  (Indeed, the only case in which they might not be smooth is
when the simple object has trivial restriction to $Q$ and trivial class in
the numerical Grothendieck group, and it is open whether such objects
can even exist!)

We then turn to semistable moduli spaces.  For the most part, the arguments
from the commutative case (see, e.g., \cite{HuybrechtsD/LehnM:1997}) carry
over, with the main exception being an inequality of Le Potier and Simpson,
for which the standard proof involves inductively restricting to hyperplane
sections.  For $1$-dimensional sheaves, there is a workaround sufficient to
take a single hyperplane section, and although the resulting inequality is
significantly weaker than in the commutative case, it turns out to still be
strong enough to show in Section \ref{sec:semistable_noncomm} that
semistable moduli spaces are projective.  (For sheaves of rank $>1$, the
question remains open!)

In Section \ref{sec:painleve_moduli}, with an eye towards the Painlev\'e
equations, we study the special case of $2$-dimensional moduli spaces with
no strictly semistable points.  Any such case induces a derived equivalence
between the moduli space and the noncommutative surface, and we can use our
understanding of the latter derived category to deduce that the moduli
space is again a smooth projective rational surface, and furthermore
describe that surface relatively explicitly.  (The {\em surfaces} are
related by one of the derived equivalences we considered in
Chapter \ref{chap:derived_eqs}, but the two derived equivalences may differ
by autoequivalences.)  Those surfaces are nothing other than the surfaces
studied by Sakai, and in particular we see that the usual second-order Lax
pairs for Painlev\'e equations are related to Sakai's construction of
Painlev\'e equations by such a derived equivalence.

In Section \ref{sec:irreds_exist}, we consider the question of when a
given stable moduli space of 1-dimensional sheaves is nonempty.  The
simplest way to show this is to exhibit an {\em irreducible} sheaf (i.e.,
such that any $1$-dimensional subsheaf has the same Chern class), and such
sheaves are easy enough to find on {\em semicommutative} surfaces, where we
can simply ask for the support on the underlying commutative surface to be
integral.  Smoothness of symplectic leaves lets us lift irreducible sheaves
along discrete valuation rings, and openness of irreducibility makes the
lift irreducible as well, so the result follows essentially immediately by
observing that semicommutative surfaces are dense.  There is, however, one
significant caveat: this density result only holds if one includes surfaces
of finite characteristic!  (In particular, even though ordinary difference
and differential equations are mainly of interest in characteristic 0, the
only proof at the moment that most such families contain irreducible
equations requires knowing that the entire theory works over $\Z$.)

In Section \ref{sec:moduli_are_irred}, we use an idea
of \cite{EllingsrudG/StrommeSA:1993,MarkmanE:2007} from the commutative
case to compute the Chow rings of stable moduli spaces when the surface is
rational and the moduli space is projective, in particular showing that
there is a unique highest-dimensional class, and thus the moduli space is
irreducible.  (The Chow ring calculation is more delicate in the higher
genus cases, but the irreducibility argument carries over.)  We then use
this in Section \ref{sec:integrable} to show that the action of twisting by
line bundles on moduli spaces of 1-dimensional sheaves is algebraically
integrable in a suitable (relative) sense, more specifically showing that
the relative degree of the pullback of a rational function on the moduli
space through such a twist grows at most quadratically in the class of the
line bundle.

In Chapter \ref{chap:diffeq_noncomm}, we again bring everything back to
equations, first carefully studying the relation between sheaves and
equations, and in particular showing that twisting a sheaf by a line bundle
indeed acts on the corresponding equation by a canonical gauge
transformation.  Furthermore, the classification of singularities in the
relaxed/commutative case (which involved relatively explicit computations
of iterated blowups) turns out to carry over immediately to the
noncommutative case, without our needing to generalize the {\em arguments}.

Although this is enough to explain {\em discrete} isomonodromy, {\em
continuous} isomonodromy is more delicate.  In particular, although we show
that the dimension of the space of continuous isomonodromy deformations is
intrinsic to the surface, we construct it in two steps that are not
intrinsic.  (There is, however, a somewhat intrinsic description of what
continuous isomonodromy deformations should be: they come from
infinitesimal deformations of the blown up surface along which sheaves
disjoint from the anticanonical curve naturally extend.)  There is also a
residual mystery: in the Painlev\'e case, continuous isomonodromy can
always be obtained as a limit of discrete isomonodromy, and thus we {\em
should} be able to explain continuous isomonodromy in the noncommutative
setting as a limit of twisting by a line bundle!

Finally, in Chapter \ref{chap:openprobs}, we collect the more significant
of the open problems from Parts I and II.  The most significance of these
comes from the derived equivalences of deformed elliptic surfaces.  In the
differential case, this gives a derived equivalence between a moduli space
of second-order Fuchsian differential equations with four singular points
and a certain noncommutative surface that deforms the corresponding relaxed
moduli space.  This is precisely the sort of derived equivalence that
arises in geometric Langlands theory, and the corresponding structure group
is $\GL_2$, so that this is not simply a trivial instance.  In other words,
for this particular nontrivial subcase of geometric Langlands, we can
generalize the derived equivalence to replace differential equations with
difference, $q$-difference, or even symmetric elliptic difference
equations!  Moreover, our moduli spaces are commutative deformations of
moduli spaces of sheaves on commutative surfaces, and the Poisson structure
on the latter strongly suggests that there should be a corresponding {\em
noncommutative} deformation.  Since the moduli space of sheaves on the
commutative surface has a (Lagrangian) fibration by abelian varieties, it
has induced derived autoequivalences, and it is thus natural to conjecture
that these derived autoequivalences extend to derived equivalences of
corresponding deformations.  (In particular, the derived autoequivalences
correspond to a congruence subgroup of $\SL_2(\Z)$, and that group should
act on the parameter space of the deformation.)  A related problem is to
extend the theory to cover difference equations in which the shift matrix
is contained in an algebraic subgroup of $\GL_n$.

\chapter{Acknowledgements}

Many people have contributed helpful insights (and questions) to
this work (and to work leading up to it) over the years; particular
thanks are due to D. Arinkin, A. Borodin, F. van de Bult, D. Chan,
P. Deift, P. Etingof, T. Graber, M. Ismail, N. Joshi, A. Knutson,
M. Mazzocco, A. Nakamura, A. Okounkov, B. Poonen, B. Pym, S.
Ruijsenaars, J. Stokman, V. Spiridonov, M. Van den Bergh, and X.
Zhu, as well as to the too-many-to-name conferences where many of
these conversations took place.\footnote{Honorable mention goes to
the many people in integrable systems who have asked for
explicit examples, none of which appear below.}  

Special thanks go to P. Diaconis (the best doctoral advisor someone who
refuses to specialize could hope to have), to J. Baik (for the
collaborations where the author first encountered Lax pairs) as well as to
P.  Forrester, both for many collaborations in special functions over the
years and for asking a question many years ago that eventually led the
author to elliptic special functions (and thence to isomonodromy), as well
as to double affine Hecke algebras (a variant of which appear in disguise
in Chapter \ref{chap:daha}).  P. Etingof deserves additional thanks in this
category as well, for many collaborations over the years; even some of the
entirely unrelated projects required the author to learn techniques that
proved useful below.

Much of the author's early work on the noncommutative side
took place on sabbatical (hosted by P. Etingof) at the Massachusetts
Institute of Technology in Fall 2011; although the techniques
developed there were only applicable in the generic case, and have
since been superseded, the approach described below would not exist
without the insights developed in that process.  In addition, this
work was partially supported by grants from the National Science
Foundation, DMS-1001645 and DMS-1500806.

\tableofcontents

\mainmatter

\part{Commutative geometry}
\chapter{Sheaves from difference equations}
\label{chap:sheaves_from_eq_comm}

\section{Symmetric elliptic equations}

As we discussed in the introduction, the analogue of a differential equation
at the top (elliptic) level in the hierarchy of special functions is a {\em
  symmetric elliptic difference equation}, which we should think of as the
pair of equations
\[
v(z+q) = A(z) v(z),
\qquad
v(-z) = v(z),
\]
where $A$ is a matrix of elliptic functions subject to the consistency
condition $A(-q-z)A(z)=1$.  The natural relaxation of this problem is to
forget the difference equation, and simply classify matrices $A$ of
elliptic functions satisfying $A(-q-z)A(z)=1$.  (We will also want to take
into account singularities, but will table that question for the moment.)

Since we plan to relate this to an algebraic geometric object, it will be
helpful to rephrase this original problem in a somewhat more abstractly
geometric way.  Thus we suppose given a smooth genus 1 curve $C_\alpha$ over
an algebraically closed field $k$ (not necessarily of characteristic 0),
along with a translation $\tau_q:C_\alpha\to C_\alpha$ and a hyperelliptic
involution $\eta:C_\alpha\to C_\alpha$, i.e., such that the quotient of
$C_\alpha$ by the involution is isomorphic to $\P^1$.  In the analytic
setting, $C_\alpha$ is $\C/\Lambda$ for some lattice $\Lambda$, $\tau_q$ is
the map $z\mapsto z+q$, and $\eta$ is the map $z\mapsto -q-z$.  We take the
latter choice for $\eta$ so that the problem of classifying $A$ becomes the
following: Classify matrices $A\in \GL_n(k(C_\alpha))$ such that $\eta^*A =
A^{-1}$.

Just as the original problem can be rephrased in terms of $1$-cocycles of
the infinite dihedral group on $\GL_n(k(C_\alpha))$, this question is
itself related to nonabelian cohomology: a matrix $A$ such that
$\eta^*A=A^{-1}$ specifies a $1$-cocycle for the cyclic group
$\langle\eta\rangle$ (of order 2).  Now, the action of $\eta$ allows us to
think of $k(C_\alpha)$ as a Galois extension of the invariant subfield
$k(\P^1)$, and thus we find
\[
H^1(\langle \eta\rangle;\GL_n(k(C_\alpha))) = H^1(\Gal(k(C_\alpha)/k(\P^1)),\GL_n).
\]
It is a classical fact that the latter Galois cohomology set is trivial
(often referred to as Hilbert's Theorem 90, though Hilbert only considered
the case of a cyclic Galois group acting on $\GL_1$), and this translates
to the following fact in our case.  It 

\begin{prop}
Let $L/K$ be a quadratic field extension, and let $A\in \GL_n(L)$ be a
matrix such that $\bar{A}=A^{-1}$, where $\bar\cdot$ is the
conjugation of $L$ over $K$.  Then there exists a matrix $B\in \GL_n(L)$
such that $A = \bar{B} B^{-1}$, and $B$ is unique up to
right-multiplication by $\GL_n(K)$.
\end{prop}

\begin{proof}
  In this case, the argument is particularly simple.  Given any vector
  $w\in L^n$, the vector $v=w+A^{-1}\bar{w}$ satisfies $\bar{v}=Av$.  If
  $x$ generates $L$ over $K$, then we may write
  \[
  w = \frac{\bar{x}}{\bar{x}-x}(w+A^{-1}\bar{w})
  -\frac{1}{\bar{x}-x}(xw + A^{-1}\overline{xw}),
  \]
  and thus the images of $e_1,\dots,e_n$ and $xe_1,\dots,xe_n$ span $L^n$ as
  an $L$-space.  We in particular obtain in this way at least $n$ vectors
  satisfying $\bar{v}=Av$ which are linearly independent over $K$.

  If $B$ is the matrix with those columns, then $B\in \GL_n(L)$ and
  $\bar{B} = A B$, which is what we want.  If $B'$ is another solution of
  the equation, then
\[
\overline{B^{-1}B'} = B^{-1} A^{-1} A B' = B^{-1}B',
\]
and thus $B^{-1}B'\in \GL_n(K)$ as required.
\end{proof}

In our setting, it will turn out to be appropriate to make the
factorization have the form $A = \eta^* B^{-t} B^t$, which simply applies
${-}^{-t}$ to the matrix coming from the above construction.  (In the
noncommutative setting, the most natural correspondence between difference
equations and sheaves is contravariant and holomorphic in $B$.)  The
nonuniqueness (we can still multiply $B$ on the right by any element of
$\GL_n$) is of course still an issue, but it turns out there is a slight
modification which can be made unique.  The first step is to make the
nonuniqueness problem worse by allowing $B$ to be a map between vector
bundles.  Let $\pi_\eta:C_\alpha\to \P^1$ be the morphism quotienting by
the action of $\eta$. Then for any vector bundle $V$ on $\P^1$, and any
meromorphic (and generically invertible) map
\[
B:\pi_\eta^*V\ratto \sO_{C_\alpha}^n,
\]
we obtain a well-defined matrix $\eta^* B^{-t}B^t$,\footnote{Here, by the
transpose $B^t$, we mean the image of $B$ under the functor
$\sHom_{C_\alpha}(-,\sO_{C_\alpha})$.}

and of course any matrix with $\eta^*A = A^{-1}$ can be represented in this
way (just take $V$ to be $\sO_{\P^1}^n$!).  The advantage of allowing $V$
to be a more general vector bundle is that we can then insist that $B$ be
{\em holomorphic} (and thus injective), by absorbing any poles into $V$.
This is still non-unique, since we could freely replace $V$ by any vector
bundle it contains, and still obtain an injective morphism supporting a
factorization of $A$.  However, we can make this unique by imposing a
maximality condition on $V$.

\begin{prop}
  Suppose $B:\pi_\eta^*V_0\to k(C_\alpha)^n$ is an injective map of sheaves, with
  $V_0$ a rank $n$ vector bundle on $\P^1$.  This induces an isomorphism
  $B:\pi_\eta^*(V_0\otimes_{\sO_{\P^1}} k(\P^1))\cong k(C_\alpha)^n$, and the set of
  bundles $V\subset V_0\otimes_{\sO_{\P^1}} k(\P^1)$ such that $BV\subset
  \sO_{C_\alpha}^n$ is nonempty, with a unique maximal element.
\end{prop}

\begin{proof} That the induced map of vector spaces over $k(C_\alpha)$ is an
  isomorphism follows from the fact that it is an injective map of vector
  spaces of the same dimension.  That the set of bundles $V$ is nonempty is
  straightforward, as we have already mentioned (just absorb any poles of
  $B$ into $V$).  Finally, if $V_1$, $V_2$ are vector bundles contained in
  $V_0\otimes_{\sO_{\P^1}}k(\P^1)$ such that $BV_1, BV_2\subset \sO_{C_\alpha}^n$,
  then $V_1+V_2$ is still contained in $V_0\otimes_{\sO_{\P^1}}k(\P^1)$, so
  is torsion-free, thus a vector bundle; and $B(V_1+V_2)=BV_1+BV_2\subset
  \sO_{C_\alpha}^n$.  Since $BV\subset \sO_{C_\alpha}^n$ implies $\deg(\pi_\eta^*V)=\deg(V)\le
  0$, it follows that there is a unique maximal such bundle.
\end{proof}

\begin{rem}
  Although this argument is nonconstructive, one can, in fact, compute the
  maximal $V$.  WLOG, $BV_0\subset \sO_{C_\alpha}$ (it being
  straightforward to clear poles).  If $V_0$ is not the maximal $V$, then
  let $V_1$ be a minimal bundle contained in $V$ and strictly containing
  $V_0$.  Then $V_1/V_0$ is the structure sheaf of a point, and thus there
  is a change of basis of $V_0$ (which, as a bundle on $\P^1$ is a sum of
  line bundles) such that $B|_{\pi^*V_1}$ comes from $B|_{\pi^*V_0}$ by
  dividing a column by the pullback of a homogeneous polynomial of degree
  1.  In particular, $\det(B)|_{\pi^*V_0}$ is a multiple of that
  polynomial.  We may thus proceed as follows: for each zero $x$ of
  $\det(B)|_{\pi^*V_0}$ such that $\bar{x}$ is also a zero (counting with
  multiplicity when $x=\bar{x}$, first do changes of basis on $V_0$ to make
  the maximum number ($n-\rank(B(x))$) of columns vanish at $x$, then do
  further changes to make the maximum number of {\em those} columns vanish
  at $\bar{x}$.  (Note that changes of basis on $V_0$ correspond to column
  operations on $B$ with holomorphic, $\eta$-invariant coefficients.)  Then
  divide those columns by the unique (up to scalars) homogeneous polynomial
  of degree 1 vanishing at $x$ and $\bar{x}$ and iterate.  If we are unable
  to enlarge $V_0$ after checking every $x$, then $V_0$ must have been the
  maximal $V$, and since each enlarging step decreases $\deg(\det(B))$ by
  at least 2, we must eventually terminate.
\end{rem}

To summarize the above considerations, given any matrix $A\in
\GL_n(k(C_\alpha)))$ such that $\eta^*A=A^{-1}$, there is a canonical
factorization $A=\eta^*B^{-t}B^t$ where $B$ is an injective morphism
\[
B:\pi_\eta^*V\to \sO_{C_\alpha}^n
\]
with $V$ a rank $n$ vector bundle on $\P^1$, maximal among those supporting
a map $B$.  This canonical factorization also clarifies issues regarding
singularities.  For instance, as we mentioned in the introduction, the
equation $v(z+q)=-v(z)$ is singular at $-q/2$ as a symmetric equation,
since there are no symmetric solutions which are nonzero and holomorphic at
$-q/2$.  This is nonobvious in terms of $A$, since $A=-1$ has no zeros or
poles here, but becomes clear in terms of $B$, as we find that $B$ must in
fact vanish at every point of the form $-q/2$ (more precisely, at every
fixed point of $\eta$).  We find in general that the points where
$\tau_q^*v=Av$ is singular as a symmetric equation are precisely those
points where $\det(B)=0$.  More generally, the right way to classify
singularities of (symmetric) elliptic difference equations is to consider
the induced cocycle over the ring of ad\`eles; one can show that the
classes of such cocycles are determined by the corresponding elementary
divisors of $B$.

\begin{rem}
  A similar factorization appeared in \cite{isomonodromy}, but the reader
  should be cautioned that they are not quite the same; indeed, the
  factorization of \cite{isomonodromy} involves a partition of the
  singularities in to two subsets, and depends significantly on that
  choice.  It turns out that those matrices (up to transpose) correspond to
  canonical factorizations of cohomologous equations, see Section
  \ref{sec:elldiff}.
\end{rem}

Since $B$ is a map from a pullback, we can use adjunction to relate it to
a map to a direct image: specifying $B$ is equivalent to specifying
\[
\pi_{\eta*}B:V\to \pi_{\eta*}\sO_{C_\alpha}^n.
\]
With this in mind, we can obtain a natural extension of $B$ to a surface
containing $C_\alpha$ (as an anticanonical curve).  Indeed, since $\pi_\eta$
has degree 2, the direct image $\pi_{\eta*}\sO_{C_\alpha}$ is a vector bundle
of degree $2$, and thus we can take the corresponding projective bundle to
obtain a Hirzebruch surface $X = \P(\pi_{\eta*}\sO_{C_\alpha})$.  Note that
$X\cong F_2$, since $\pi_{\eta*}\sO_{C_\alpha}\cong \sO_{\P^1}\oplus
\sO_{\P^1}(-2)$.  Moreover, $X$ contains $C_\alpha$ in a natural way, in such
a way that the induced map from $C_\alpha$ to $\P^1$ is just $\pi_\eta$, and
$C_\alpha$ is anticanonical.  If $\rho:X\to \P^1$ is the corresponding
ruling, and $s_{\min}$ denotes the section of the ruling with minimal
self-intersection ($s_{\min}^2=-2$), then we have a canonical isomorphism
\[
\pi_{\eta*} \sO_{C_\alpha} \cong \rho_*{\cal L} (s_{\min}),
\]
since $\sO_X(s_{\min})$ is the relative $\sO(1)$.  This is just the direct
image under $\rho$ of the restriction map
\[
\sO_X(s_{\min})\to \sO_X(s_{\min})|_{C_\alpha}\cong \sO_{C_\alpha},
\]
where we note that $C_\alpha$ and $s_{\min}$ are disjoint, and we make the
isomorphism canonical by taking the unique global section of
$\sO_X(s_{\min})$ to the unique global section of $\sO_{C_\alpha}$.

In other words, to specify $B:\pi_\eta^*V\to \sO_{C_\alpha}^n$, it is
equivalent to specify its direct image
\[
\rho_*B:V\to \rho_*\sO_X(s_{\min})^n,
\]
where we now think of $B$ as a morphism of sheaves on $C_\alpha\subset X$.
Again using the adjunction between $\rho_*$ and $\rho^*$ gives us a morphism
\[
B:\rho^*V\to \sO_X(s_{\min})^n,
\]
which when restricted to $C_\alpha\subset X$ recovers the original morphism.
Again, we modify this slightly to
\[
B:\rho^*V\otimes \sO_X(-s_{\min})\to \sO_X^n,
\]
which has no effect on the restriction to $C_\alpha$, but is more natural in the
noncommutative setting (and slightly more natural even in the commutative
setting).  In any event, we now have a morphism of vector bundles on the
Hirzebruch surface $X$.

\begin{rem}
Again, this is more constructive than it might appear.  Since $C_\alpha$ is
a double cover of $\P^1$, we can express it as the zero-locus in $F_2$ of
an equation $y^2+a(x,w)y+b(x,w)=0$ where $\deg(a)=2$, $\deg(b)=4$.  The
entries of $B$ are thus expressible as homogeneous polynomials in $y$, $x$,
$w$ of the appropriate degree (assuming an identification of $V$ with a sum
of line bundles).  This expression is not unique, but there is a natural
choice: use the equation of $C_\alpha$ to eliminate all terms of degree
$>1$ in $y$.  But this choice is nothing other than the above extension of
$B$ to a morphism of vector bundles on $F_2$!
\end{rem}

We have thus reduced our classification problem to one of classifying
certain types of morphisms of vector bundles on $F_2$.  If we view $F_2$ as
a family of $\P^1$s over $\P^1$, then our morphism becomes a family of
morphisms of bundles on $\P^1$, each fiber of which is a morphism
$\sO_{\P^1}(-1)^n\to \sO_{\P^1}^n$.  This is precisely the form one
expects for the minimal resolution of a torsion sheaf on $\P^1$ (of degree
$n$), suggesting the following, which we state in greater generality for
future use.

\begin{prop}\label{prop:weak_semiorth_on_ruled}
  Let $\rho:X\to C$ be a ruled surface, with relative $\sO(1)$ denoted by
  $\sO_\rho(1)$, and let $M$ be a coherent sheaf on $X$.  Then the
  following are equivalent.
\begin{itemize}
\item[1.] $M$ is the cokernel of an injective morphism
\[
B:\rho^*V\otimes \sO_\rho(-1)\to \rho^*W
\]
with $V$, $W$ vector bundles of the same rank on $C$.
\item[2.]  $M$ has $1$-dimensional support, $M\otimes \sO_\rho(-1)$ is
  $\rho_*$-acyclic, and $\rho_*M$ is torsion-free.
\end{itemize}
Moreover, if either condition holds, then $B$ is uniquely determined
up to isomorphism by $M$.
\end{prop}

\begin{proof}
  $1\implies 2$: Since $B$ is an injective morphism of vector bundles of
  the same rank, it is an isomorphism on the generic fiber, and thus
  $\supp(M)$ does not contain the generic point of $X$.  It follows that
  $M$ has $\le 1$-dimensional support.  (In fact, $M$ is supported on the
  zero locus of $\det(B)$.)

  Now, since the sheaves $\sO_\rho(d)$ are isomorphic to $\sO_f(d)$ on
  every fiber $f$, we find that $\sO_\rho(d)$ is $\rho_*$-acyclic for $d\ge
  -1$, and has trivial direct image for $d\le -1$.  In particular, we can
  compute the higher direct image long exact sequence associated to the
  short exact sequence
\[
0\to \rho^*V\otimes \sO_\rho(-2)\to \rho^*W\otimes \sO_\rho(-1)\to M\otimes
\sO_\rho(-1)\to 0.
\]
Since $\rho$ has $1$-dimensional fibers, this long exact sequence
terminates after degree 1, and we conclude that $M\otimes \sO_\rho(-1)$ is
$\rho_*$-acyclic.  Similarly, from the untwisted short exact sequence, we
obtain
\[
W\cong \rho_*\rho^*W\cong \rho_*M.
\]
In particular, $\rho_*M$ is torsion-free, and we can recover $B$ as the
kernel of the natural map $\rho^*\rho_*M\to M$.

$2\implies 1$: The condition that $M\otimes \sO_\rho(-1)$ is
$\rho_*$-acyclic implies that if we view $M$ as a family of sheaves on
$\P^1$, then every fiber satisfies $H^1(M_f(-1))=0$.  In particular, every
fiber is $0$-regular in the sense of Castelnuovo and Mumford, and thus $M$
is relatively globally generated \cite{KleimanS:1971}. Since $\rho_*M$ is
torsion-free by assumption, so a vector bundle, it remains only to show
that the kernel of this natural map has the form $\rho^*V\otimes
\sO_\rho(-1)$.  Now, $M$ cannot have any $0$-dimensional subsheaf,
since that would produce a $0$-dimensional subsheaf of $\rho_*M$.  In other
words, $M$ is a pure $1$-dimensional sheaf, and thus has homological
dimension $1$.  In particular, the kernel is a vector bundle (of the same
rank as $W$, since the map is generically surjective), so we can view it as
a flat family of sheaves on $\P^1$.  Since $\rho^*\rho_*M$ and $M$ are
acyclic with isomorphic direct image, it follows that the kernel has
trivial direct image and higher direct image, and thus every fiber of the
kernel has trivial cohomology.  The only sheaves on $\P^1$ with trivial
cohomology are sums of $\sO_{\P^1}(-1)$, and thus the kernel has the form
\[
V'\otimes \sO_\rho(-1)
\]
where $V'$ is a flat family of sheaves on $\P^1$, each fiber of which is a
power of $\sO_{\P^1}$.  In other words, $V'\cong \rho^*V$ for some vector
bundle $V$.
\end{proof}

\begin{rems}
Just as we found $W\cong \rho_*M$, we can also compute $V$ from $M$, since
\[
\rho_*(M\otimes \sO_\rho(-1))\cong V\otimes R^1\rho_*\sO_\rho(-2),
\]
and $R^1\rho_*\sO_\rho(-2)$ is a line bundle on $C$.
\end{rems}

\begin{rems}
  This argument was inspired by the main construction of \cite{BeauvilleA:2000},
  which considered minimal resolutions of sheaves on $\P^n$ for $n>1$; in
  our case, we have a {\em relative} minimal resolution of a family of
  sheaves on $\P^1$.
\end{rems}

\begin{rems}
  This is a special case of a much more general result: any {\em complex}
  $M^{\cdot}\in D^b_{\coh}(X)$ fits into a natural distinguished triangle
  \[
  \rho^*V^{\cdot}\otimes \sO_{\rho}(-1)\to \rho^*W^{\cdot}\to M^{\cdot}\to
  \]
  where $V^{\cdot}$, $W^{\cdot}$ are complexes on $\P^1$.  Indeed, one has
  a natural distinguished triangle
  \[
  N^{\cdot}\to \rho^*R\rho_*M^{\cdot}\to M^{\cdot}\to
  \]
  where $R\rho_*N^{\cdot}=0$, and this implies (checking fiber-by-fiber)
  that $N^{\cdot}\cong \rho^*V^{\cdot}\otimes \sO_{\rho}(-1)$.  (In other
    words, $D^b_{\coh}(X)$ has a semiorthogonal decomposition
    $(\rho^*D^b_{\coh}(\P^1)\otimes
    \sO_{\rho}(-1),\rho^*D^b_{\coh}(\P^1))$.)  For sheaves, this translates
    to a (natural, up to choosing an isomorphism $\omega_{\P^1}\cong
    \sO_{\P^1}(-2)$) five-term sequence
    \[
    0\to \rho^*\rho_*(N_0(-1))(-1)\to \rho^*\rho_* M\to M
     \to \rho^*R^1\rho_*(N_1(-1))(-1)\to \rho^*R^1\rho_*M\to 0.
    \]
\end{rems}

Of course, there remain two conditions to translate into conditions on the
sheaf $M$, namely the constraint on the singularities, and the constraint
that $V$ is maximal.  The former is straightforward: specifying the
elementary divisors of $B$ along $C_\alpha$ is equivalent to specifying the
cokernel of $B$ as a morphism of vector bundles on $C_\alpha$, and thus the
singularities are determined by the restriction $M|_{C_\alpha}$.  (In
particular, we have the overall constraint that $M$ must be transverse to
$C_\alpha$, so that $B$ is generically invertible on $C_\alpha$!)  The
latter is somewhat more subtle, but is not too difficult to deal with.

\begin{prop}
Let $\rho:X\to C$ be a ruled surface, and suppose the sheaf $M$ is given by
a presentation
\[
0\to \rho^*V\otimes \sO_\rho(-1)\xrightarrow{B} \rho^*W\to M\to 0,
\]
where $V$ and $W$ are vector bundles of the same rank.
The morphism $B$ extends to a supersheaf $V\subsetneq V'\subset
V\otimes_{\sO_C} k(C)$ iff $M$ has a subsheaf of the form $\sO_f(-1)$ for
some fiber $f$ of $\rho$.
\end{prop}

\begin{proof}
If $B$ extends to $V'$, then the image of $\rho^*V'\otimes \sO_\rho(-1)$
induces a subsheaf of $M$ isomorphic to
\[
\rho^*(V'/V)\otimes \sO_\rho(-1).
\]
Now, $V'/V$ is $0$-dimensional, so contains a subsheaf of the form $\sO_p$
for some closed point $p\in C$.  This $\sO_p$ itself induces a supersheaf
of $V$, and thus a subsheaf of $M$ of the form
\[
\rho^*(\sO_p)\otimes \sO_\rho(-1)\cong \sO_f(-1),
\]
where $f$ is the fiber over $p$.

Conversely, suppose we have an injective map $\sO_f(-1)\to M$, and let $M'$
be the cokernel.  The higher direct image long exact sequences tell us
\[
\rho_*M'\cong \rho_*M
\qquad
R^1\rho_*M'\cong R^1\rho_*M=0
\qquad
R^1\rho_*(M'\otimes \sO_\rho(-1))\cong R^1\rho_*(M'\otimes \sO_\rho(-1))=0,
\]
and thus $M'$ has a presentation of the form
\[
0\to \rho^*V'\otimes \sO_\rho(-1)\to \rho^*W\to M'\to 0.
\]
Since this construction is functorial, we obtain an injective morphism
\[
\rho^*V\otimes \sO_\rho(-1)\to \rho^*V'\otimes \sO_\rho(-1),
\]
thus an injective morphism $\rho^*V\to \rho^*V'$, and by adjunction,
$V\subset V'$ in such a way that $B$ extends.
\end{proof}

There is a dual condition related to relative global generation.

\begin{prop}
Let $\rho:X\to C$ be a ruled surface, and suppose that $M$ is a pure
$1$-dimensional sheaf on $X$.  If $M$ is $\rho_*$-acyclic, then $M$ is
relatively globally generated iff no quotient of $M$ has the form
$\sO_f(-1)$
for some fiber $f$ of $\rho$.
\end{prop}

\begin{proof}
If $M$ is relatively globally generated, then
\[
\Hom(M,\sO_f(-1))\subset \Hom(\rho^*\rho_*M,\sO_f(-1))\cong
\Hom(\rho_*M,\rho_*\sO_f(-1))=0.
\]

For the converse, consider the natural map $\rho^*\rho_*M\to M$, viewed as
a two-term complex.  The terms in the complex are $\rho_*$-acyclic, and
thus the derived direct image of the complex is
\[
\rho_*\rho^*\rho_*M\cong \rho_* M,
\]
so is exact.  On the other hand, there is a spectral sequence converging to
this result in which we first take the cohomology of the complex before
taking higher direct images.  Since $\rho$ has $1$-dimensional fibers, this
spectral sequence stabilizes at the $E_2$ page, and we thus conclude that
the cohomology sheaves of the complex have trivial direct image and higher
direct image.

We thus conclude that if $M$ is not globally generated, then $M$ has a
surjective morphism to a nonzero sheaf $M'$ with $\rho_*M'=R^1\rho_*M'=0$,
so
\[
M'\cong \rho^*\rho_*(M'\otimes \sO_\rho(1))\otimes \sO_\rho(-1)
\]
Now, $\rho_*(M'\otimes \sO_\rho(1))$ cannot be 0, since that would force
$M'=0$.  It thus admits a surjective map to some $\sO_p$, which induces a
surjective map from $M'$ to a sheaf of the form $\sO_f(-1)$.
\end{proof}

\begin{rem}
In fact, it follows from this that if $M$ is pure $1$-dimensional and
$\rho_*$-acyclic, then it is relatively globally generated iff $M\otimes
\sO_\rho(-1)$ is $\rho_*$-acyclic.  Indeed, a surjection $M\to \sO_f(-1)$
induces a surjection
\[
R^1\rho_*(M\otimes \sO_\rho(-1)) \to R^1\sO_f(-2)\cong \sO_{\pi(f)}
\]
making the former sheaf nontrivial.  (This can also be seen directly from
the five-term sequence discussed above.)
\end{rem}

Similar conditions apply to $\rho_*$-acyclicity and torsion-freeness of
$\rho_*M$.

\begin{prop}
Let $\rho:X\to C$ be a ruled surface, and let $M$ be a pure $1$-dimensional
sheaf on $X$.  Then $\rho_*M$ is torsion-free iff $\Hom(\sO_f,M)=0$ for all
fibers $f$ of $\rho$, and $M$ is $\rho_*$-acyclic iff $\Hom(M,\sO_f(-2))=0$
for all $f$.
\end{prop}

\begin{proof}
For the first condition, we have
\[
\Hom(\sO_f,M)
\cong
\Hom(\rho^*\sO_{\pi(f)},M)
\cong
\Hom(\sO_{\pi(f)},\rho_*M)
\]
Since $\rho_*M$ is torsion-free iff it has no maps from point sheaves, the
first claim follows.

For the second, the same spectral sequence argument based on the complex
$\rho^*\rho_*M\to M$ tells us that if $M'$ is the cokernel of this natural
map, then $R^1\rho_*M\cong R^1\rho_*M'$ and $\rho_*M'=0$.  Since $M'$ is
supported on finitely many fibers as before, it must have a quotient of the
form $\sO_f(-d)$ for some $d>1$, and thus has a nontrivial morphism to
$\sO_f(-2)$.
\end{proof}

If $\Hom(\sO_f,M)=0$ but $\Hom(\sO_f(-1),M)\ne 0$, then any such morphism
is necessarily injective; similarly, if $\Hom(M,\sO_f(-2))=0$ but
$\Hom(M,\sO_f(-1))\ne 0$, then any such morphism is surjective.  We thus
arrive at the final moduli problem: Classify pure $1$-dimensional sheaves
$M$ on $X$ with specified restriction to $C_\alpha$ and such that
$\Hom(M,\sO_f(-1))=\Hom(\sO_f(-1),M)=0$ for all fibers $f$ of the ruling.

We are also interested in understanding when $W$ is trivial, for which we
have the following numerical condition in the rational case.

\begin{lem}\label{lem:cohom_vanish}
Let $\rho:X\to \P^1$ be a Hirzebruch surface, and let $M$ be a sheaf on $X$
with $1$-dimensional support.  Then the following are equivalent:
\begin{itemize}
\item[(a)] $H^0(M)=H^1(M)=0$
\item[(b)] $M$ is $\rho_*$-acyclic and $\rho_*M\cong \sO_{\P^1}(-1)^n$ for
  some $n\ge 0$.
\end{itemize}
\end{lem}

\begin{proof}
  Since $\rho$ has $1$-dimensional fibers, $R^p\rho_*M=0$ for $p>1$; since
  the generic fiber of $\supp(M)$ over $\P^1$ is $0$-dimensional,
  $R^1\rho_*M$ has $0$-dimensional support.  The
  Leray-Serre spectral sequence
\[
H^p(R^q\rho_*M)\implies H^{p+q}(M)
\]
thus implies isomorphisms
\[
H^0(M)\cong H^0(\rho_*M)\qquad H^2(M)\cong H^1(R^1\rho_*M)=0
\]
and a short exact sequence
\[
0\to H^1(\rho_*M)\to H^1(M)\to H^0(R^1\rho_*M)\to 0.
\]
In particular, $H^0(M)=H^1(M)=0$ iff both $\rho_*M$ and $R^1\rho_*M$ have
vanishing cohomology.  In particular, both must be isomorphic to a direct
sum of line bundles $\sO_{\P^1}(-1)$, and since $R^1\rho_*M$ has
$0$-dimensional support, it must be $0$.
\end{proof}

Thus the only way a 1-dimensional sheaf with $H^*(M\otimes
\rho^*\sO_{\P^1}(-1))=0$ could fail to have a presentation of the standard
form is if $\Hom(M,\sO_f(-1))\ne 0$ for some fiber $f$.  (Even if
it fails this last condition, the image of $\rho^*\rho_*M\to M$ gives us a
subsheaf with standard presentation, and thus a subquotient with standard
presentation satisfying maximality.)

One particularly important example of a $1$-dimensional sheaf with standard
presentation is the structure sheaf $\sO_{s_{\min}}$ of the $-2$-curve on
$F_2$.  This has standard presentation
\[
0\to (\rho^*\sO_{\P^1})(-s_{\min})\to \rho^*\sO_{\P^1}\to \sO_{s_{\min}}\to
0
\]
and translating the resulting matrix $B$ back to an equation on $C_\alpha$
gives the equation $v(z+q)=v(z)$, which has no singularities since
$s_{\min}$ is disjoint from $C_\alpha$.

\medskip

When we considered elliptic difference equations above, this was in fact a
simplification: in fact, the difference equations that occur in the theory
of elliptic special functions have {\em theta function} coefficients in
general.  Algebraically speaking, $A$ is really a matrix with coefficients
in ${\cal L}_0\otimes_{\sO_{C_\alpha}}k(C_\alpha)$, where ${\cal L}_0$ is a
degree 0 line bundle equipped with an isomorphism
\[
\eta^*{\cal L}_0\cong {\cal L}_0^*
\]
such that the composition
\[
{\cal L}_0=\eta^*\eta^*{\cal L}_0\cong \eta^*{\cal L}_0^*\cong {\cal L}_0
\]
is the identity.  Again, Hilbert's Theorem 90 allows us to factor ${\cal L}_0$,
though now the nonuniqueness is more significant.  If ${\cal L}_0$ has the
above form, then it can (since we are over an algebraically closed field)
be factored as
\[
\psi:{\cal L}_0\cong \eta^*{\cal L}\otimes {\cal L}^*,
\]
for some line bundle ${\cal L}$.  This line bundle is nonunique in two
respects: the obvious one is that it can be twisted by any power of the
line bundle $\pi_\eta^*{\cal O}_{\P^1}(1)$, but even modulo this, there are
8 possibilities for ${\cal L}$, and for each such choice, 2 possibilities
for the isomorphism $\psi$.  Indeed, we can twist ${\cal L}$ by the degree
$1$ bundle corresponding to any ramification point of $\eta$; if we twist
by them all, this is the same as twisting by $\pi_\eta^*{\cal
  O}_{\P^1}(2)$, except that the isomorphism is multiplied by $-1$.  (In
characteristic 2, the issue of nonuniqueness is somewhat more complicated.
In general, modulo twisting $\pi_\eta^*{\cal O}_{\P^1}(1)$, the
factorizations form a torsor over an abelian group scheme with structure
$\mu_2.\Pic^0(C_\alpha)[2].\Z/2\Z$.)

This nonuniqueness is related to the question of singularities at the four
ramification points: we can no longer canonically distinguish between the
two possible local rank $1$ equations, in order to decide which one should
be viewed as regular.  Since we want to specify the singularity structure,
we should view the factorization of ${\cal L}_0$ as part of the
specification (and indeed, in all of the motivating cases, there is a
natural choice of factorization making the equation regular at the fixed
points of $\eta$ for generic parameters).  With this in mind, we again
obtain a canonical factorization, except now $B$ has the form
\[
B:\pi_\eta^*V\to \pi_\eta^*W\otimes {\cal L}.
\]
Again, two uses of adjunction allow us to extend this to a morphism of
vector bundles on the Hirzebruch surface $\P(\pi_{\eta *}{\cal L})$,
and thus further to the cokernel of that morphism.  The above
considerations extend immediately to the case of general ${\cal L}$.

Note that $\P(\pi_{\eta*}{\cal L})$ is isomorphic to the Hirzebruch surface
$F_1$ iff ${\cal L}$ has odd degree, and is otherwise isomorphic to $F_0$
or $F_2$, with the latter precisely when ${\cal L}$ is a power of
$\pi_\eta^*\sO_{\P^1}(1)$.

We should also note that using the freedom to twist ${\cal L}$ by powers of
$\pi_\eta^*\sO_{\P^1}(1)$, we can assume that ${\cal L}^*$ is represented
by an effective divisor disjoint from the ramification locus.  The
resulting isomorphism ${\cal L}\cong \sO_{C_\alpha}(-D)$ allows us to
translate the problem back to one on $F_2$, but with additional
``apparent'' singularities along $D$.  (These are points where the
obstructions to having symmetric solutions can be gauged away, possibly at
the expense of introducing apparent singularities along other orbits of the
infinite dihedral group.)

\section{Degenerations and generalizations}

Although the generic anticanonical curve on $F_2$ is a smooth genus 1
curve, there is significant scope for degeneration.  We can view $F_2$ as
the minimal desingularization of a weighted projective space, and in this
way anticanonical curves correspond to equations of the form
\[
p_0(x,w) y^2 + p_2(x,w) y + p_4(x,w)=0,
\]
with $p_d(x,w)$ homogeneous of degree $d$.  (Here $w$,$x$,$y$ are
generators of degrees $1$, $1$, and $2$ of a graded algebra, and $F_2$ is
the minimal desingularization of $\Proj(k[w,x,y])$.)  If we want a natural
notion of trivial equation, we need this curve to be disjoint from
$s_{\min}$, forcing $p_0\ne 0$, so that we reduce to considering equations
equations
\[
y^2 + p_2(x,w) y + p_4(x,w)=0.
\]
On any such curve, we have an involution $y\mapsto -p_2(x,w)-y$ which
exchanges the two points on any given fiber.  Given any such curve, the
translation between matrices $B$ on $C_\alpha$ and matrices on $F_2$ is quite
explicit in terms of coordinates: express $B$ in terms of the coordinates,
and use the equation of $C_\alpha$ to eliminate any term of degree $\ge 2$ in
$y$.  The resulting matrix, every coefficient of which is linear in $y$ and
weighted homogeneous, can now be viewed as a matrix on $F_2$, and is a
canonical extension of $B$.

There are several degenerate cases to consider.  In characteristic not 2,
we can complete the square to make $p_2=0$, and the degenerate cases are
classified by the multiplicities of the zeros of $p_4$ (with an additional
case when $p_2=p_4=0$).  These cases all extend to characteristic 2; there
is one extra case in characteristic 2 ($p_2=0$, $p_4$ is not a square) for
which we do not have a natural difference/differential equation
interpretation, so do not consider below.  Note that over a perfect field,
any such curve is equivalent to the curve $y^2=x w^3$.

\begin{itemize}
\item[211:] $C_\alpha$ is integral, with a single node.  Then $C_\alpha$ is
  isomorphic to $\P^1$ with $0$ and $\infty$ identified, and $\eta$ acts on
  this $\P^1$ as $z\mapsto \beta/z$ for some $\beta$.  Up to a change of
  coordinates on $F_2$, $C_\alpha$ has the equation
\[
y^2-xwy+\beta w^4=0,
\]
with $w(z)=1$, $x(z)=z+\beta/z$, $y(z)=z$.  If we choose an automorphism
$\tau_q:z\mapsto qz$, then the above construction applies to relate sheaves
on $F_2$ to symmetric $q$-difference equations
\[
v(qz) = A(z)v(z),
\]
where the symmetry condition is that $A\in \GL_n(k(z))$ satisfies
$A(\beta/z)=A(z)^{-1}$.  Note that this is singular at the node unless
$A(0)=A(\infty)=1$, which in turn happens iff $\det(B)$ is nonzero at the
node.

\item[31:] $C_\alpha$ is integral, with a single cusp.  Then $C_\alpha$ may
  be identified with $\P^1$ such that the cusp maps to $\infty$ and
  $\eta(z)=\beta-z$ for some $\beta$.  Up to a change of coordinates on
  $F_2$, $C_\alpha$ has the equation
\[
y^2-\beta w^2 y+x w^3= 0,
\]
with $w(z)=1$, $x(z)=z(\beta-z)$, $y(z)=z$.  These correspond to symmetric
ordinary difference equations:
\[
v(z+\hbar)=A(z)v(z)
\]
with $A\in \GL_n(k(z))$ such that $A(\beta-z)A(z)=1$.  The equation is
singular at the cusp unless $A(z)=1+O(1/z^2)$ as $z\to\infty$, or
equivalently unless $\det(B)$ vanishes at $\infty$. (The symmetry then
implies $A(z)=1+O(1/z^3)$.)  Here, of course, we can always rescale $z$ to
make $\hbar=1$ when it is nonzero, but it is convenient to retain this
degree of freedom so that we can consider the commutative relaxation
$\hbar\to 0$.

\item[22:] $C_\alpha$ is a union of two smooth components (isomorphic to
  $\P^1$) meeting in two distinct points, and $\eta$ swaps the components.
  In suitable coordinates, $C_\alpha$ has the equation
\[
y^2-xwy=0,
\]
with components $y=xw$, $y=0$; we may view $z=x/w$ as a common coordinate
on the two components.  Then any morphism $B$ on $F_2$ as above specifies a
{\em pair} of morphisms
\[
B_1,B_2:V\to \sO_{\P^1}^n,
\]
agreeing at $0$ and $\infty$.  The corresponding $A$ matrix on $C_\alpha$ is
really a pair of inverse matrices, but we may simply view it as a single
matrix $B_2^{-t}B_1^t$ on one component of $C_\alpha$.  In this way, we
obtain a $q$-difference equation on $\P^1$ without any symmetry condition,
and $B_1$, $B_2$ separate the zeros and poles of the equation.
Singularities at the two nodes of $C_\alpha$ arise when $A(0)\ne 1$ or
$A(\infty)\ne 1$ respectively.

\item[4:] $C_\alpha$ is a union of two smooth components which are tangent at
  a single point, and $\eta$ swaps the components; $C_\alpha$ has equation
  $y^2=w^2 y$, up to changes of coordinates.  Again $B$ specifies a pair of
  morphisms, which now agree to second order at $\infty$.  This corresponds
  to ordinary difference equations without symmetry, which are singular at
  $\infty$ unless $A(z)=1+O(1/z^2)$ as $z\to\infty$.

\item[0:] $C_\alpha$ is nonreduced, with equation $y^2=0$ after a change of
  coordinates.  In this case, the degree 2 morphism $C_\alpha\to \P^1$ is no
  longer generically \'etale, so is not the quotient by an involution.
  However, we can now identify $B$ with a pair of maps
\[
B_0:V\to \sO_{\P^1}^n,\qquad
B_\omega: V\to \omega_{\P^1}^n,
\]
giving a canonical factorization of a meromorphic matrix taking values in
$\omega_{\P^1}$ (a.k.a., a meromorphic Higgs bundle).  Since there is a
canonical connection on $\sO_{\P^1}$, we may use this to interpret the
meromorphic matrix with values in $1$-forms as a meromorphic connection.
This, of course, is just the standard translation between differential
equations and sheaves on $F_2$ arising in the usual theory of Hitchin
systems.
\end{itemize}

\begin{rem}
  The remaining example, the anticanonical curve $y^2=xw^3$ on
  $F_2(\overline{\F_2})$, is an irreducible cuspidal curve as in case $31$
  above, but the degree 2 map to $\P^1$ is not \'etale, so that it wants to
  be simultaneously a difference and differential equation.  One can
  resolve this using the the noncommutative geometry construction below,
  and find that the corresponding algebra of operators alternates between
  differential and difference operators.  The action of these operators is
  far from faithful (there is a kernel in degree 2), which is already an
  issue for ordinary difference and differential equations in finite
  characteristic (where the shift has finite order and differentiation is
  nilpotent).  Of course, for actual applications to integrable systems,
  one only cares about characteristic 0, but from a theoretical
  perspective, there is no reason to exclude finite characteristic.
  (Furthermore, there is one result below, Proposition
  \ref{prop:irreds_exist} about the existence of irreducible equations with
  specified singularities, for which the only proof we have involves
  reduction to finite characteristic!)
\end{rem}

The nonsymmetric difference and differential cases extend to higher genus.

\begin{eg}
Consider a general elliptic difference equation $\tau_q^*v=Av$ with $A\in
\GL_n(k(C))$ for some genus 1 curve $C$.  There is a natural factorization
\[
A = B_\infty^{-t} B_0^t
\]
where $B_0,B_\infty:V\to \sO_{C}^n$ and $V$ is a maximal vector bundle
supporting such a factorization.  (The existence of a meromorphic
factorization is trivial (take $B_\infty=1$, $B_0=A^t$), and implies the
existence of a unique maximal $V$ as above.)  The singularity structure of
$A$ then corresponds in a natural way to the cokernels of $B_0$ and
$B_\infty$ (giving zeros and poles respectively).  The pair
$(B_0,B_\infty)$ extends immediately to a morphism of vector bundles on the
ruled surface $E\times \P^1$: just take $zB_\infty+wB_0$, where $(z,w)$ are
homogeneous coordinates on $\P^1$.  We can then recover the pair as the
restriction of this morphism to the anticanonical curve $zw=0$, a union of
two disjoint copies of $C$.  This is essentially just the construction for
the Sklyanin integrable system (see \cite{HurtubiseJC/MarkmanE:2002a}), the
only difference being that the construction in the literature twists by a
line bundle in order to absorb all of the poles, making $B_\infty=1$, but
making the ruled surface more complicated.  (This corresponds to performing
a sequence of elementary transformations centered at the points where the
sheaf $M$ meets the component $w=0$ of the anticanonical curve.)  More
generally, we should allow difference equations on vector bundles, i.e.,
meromorphic (and meromorphically invertible) maps $A:V\to\tau_q^*V$.  If
there is a holomorphic isomorphism $A_0:V\cong \tau_q^*V$, then one can
divide by $A_0$ to again reduce to a sheaf on $E\times \P^1$.  By Atiyah's
classification of vector bundles on smooth genus 1 curves
\cite{AtiyahMF:1957}, we find that $V\cong \tau_q^*V$ whenever $V$ is a sum
of indecomposable bundles of degree 0 (and this is necessary if $q$ has
infinite order).  This is again an open condition on $V$ (for degree 0
bundles, it is equivalent to semistability), and the isomorphism $A_0$ can
be at least partially globalized (i.e., it exists on some open cover of the
open subset); thus, just as in the rational cases, we can identify large
open subsets of the moduli spaces of sheaves and of difference equations.
(The main distinction is that the identification is no longer canonical,
since $A_0$ is only determined up to scalars.)  Note also that $V$ is
semistable of degree 0 iff there exists some line bundle ${\cal L}$ of
degree 0 such that $H^0(V\otimes {\cal L})=H^1(V\otimes {\cal L})=0$; this
gives an analogue of Lemma \ref{lem:cohom_vanish} for the elliptic case.
\end{eg}

\begin{eg}
Similarly, the Hitchin system corresponds to the analogous factorization
for a meromorphic morphism
\[
A:k(C)\to \omega_C\otimes_{\sO_{C}} k(C),
\]
and gives a sheaf on the surface $\P(\sO_C\oplus \omega_C)$.  As before,
the construction in the literature essentially absorbs the poles of $A$
into the structure of the surface, but this is just a sequence of
elementary transformations.  Once again, the ``true'' moduli space (of
meromorphic connections on vector bundles) and the moduli space of sheaves
can be identified along large open subsets; in this case, the requirement
is that the vector bundle $V$ admit a holomorphic connection.  This is no
longer an open condition (it is equivalent to every indecomposable summand
having degree a multiple of the characteristic
\cite{BiswasI/SubramanianS:2006}), but is implied by open conditions, e.g.,
that $V$ is semistable of degree 0 or stable of degree a multiple of the
characteristic.  It is also implied by the open condition that
$H^0(V\otimes {\cal L})=H^1(V\otimes {\cal L})=0$ for some line bundle
${\cal L}$ of degree $g-1$, though this is no longer equivalent to
semistability.
\end{eg}

In each case, the equation gives rise to a sheaf on a ruled surface with an
embedded curve $C_\alpha$ such that (a) the equation is recovered from the
restriction of a corresponding morphism to $C_\alpha$, and (b) the
singularities of the equation correspond to points where the sheaf meets
$C_\alpha$.  Moreover, in all of our examples, the curve $C_\alpha$ is
anticanonical.  There is a very close relationship between anticanonical
curves on surfaces and Poisson structures: a Poisson structure on $X$ is a
biderivation, so a section of $\omega_X^{-1}$, and any such section
satisfies the Jacobi identity.  In particular, we see that an anticanonical
curve on $X$ determines a Poisson structure on $X$ up to a multiplicative
scalar (which can be fixed by a choice of structure on $C_\alpha$, see
Section \ref{sec:classify_poiss_surf} below).  In fact, we will see in
Chapter \ref{chap:classify} that the above examples hit essentially {\em
  every} birational equivalence class of Poisson surfaces (again with some
exceptions in characteristic 2).  Moreover, as discussed in Section
\ref{sec:comm_poisson_moduli}, this Poisson structure on $X$ induces a
Poisson structure on the moduli space of sheaves on $X$ in such a way that
the fibers of restriction to $C_{\alpha}$ are symplectic.  (In other words,
the moduli spaces of relaxed equations with fixed singularities are
naturally symplectic!)

We also saw in the elliptic case that we could obtain other Hirzebruch
surfaces (still with anticanonical curves) by allowing our equations to be
twisted by a line bundle, or equivalently by allowing apparent
singularities.  As noted, these instances are related by a sequence of
elementary transformations, and in particular by a birational map
respecting the anticanonical curve.  The two surfaces are, in fact, covered
by a common blowup, obtained from $F_2$ by blowing up the apparent
singularities, and the sheaf $M$ can be expressed as the direct image of a
sheaf on the common blowup.  (This is straightforward if $M$ is the direct
image of a line bundle on a smooth curve; for more general sheaves, one
must use the ``minimal lift'' construction of Section
\ref{sec:minimal_lift}.)  One can in fact do something similar for more
general singularities, and thus find that, with the exception of
particularly pathological singularities one can always express $M$ as the
direct image of a sheaf on some iterated blowup in such a way that the
blowup has an induced anticanonical curve {\em disjoint} from the sheaf
associated to $M$.  (Again, we will discuss the details below, where we
will also show that the pathological singularities can always be resolved
by suitable gauge transformations.)

\medskip

There are a few other moduli problems that translate to sheaves on
anticanonical rational surfaces that we want to consider.

\begin{eg}
  Let $C_\alpha$ be a smooth genus 1 curve, and let ${\cal L}$ be a line
  bundle on $C_\alpha$ with $\deg{\cal L}\ge 0$.  Consider the problem of
  classifying matrices $B\in \Mat_n(\Gamma({\cal L}))$ such that
  $\det(B)\ne 0$, up to left- and right- multiplication by constant
  matrices.  For any choice of hyperelliptic involution $\eta$, we can
  encode such matrices via sheaves on the Hirzebruch surface
  $\P(\pi_{\eta*}{\cal L})$, essentially as above.  The only additional
  condition we impose is that the bundle $V$ must also be trivial.
  In the case $\deg({\cal L})=3$, we can also interpret the matrix as one
  over $\Gamma(\sO_{\P^2}(1))$, and on $\P^2$, the matrix is the minimal
  resolution of its cokernel, as in \cite{BeauvilleA:2000}.  Note that
  since in this case we want both bundles to be trivial up to twist, we
  need to impose another numerical condition, which is straightforward to
  determine from the formula for $V$ in terms of $M$ that we gave above.
\end{eg}

\begin{eg}
  On the Hirzebruch surface $F_0=\P^1\times \P^1$, which we can view as the
  smooth quadric $xz=yw$ in $\P^3$, we can apply the previous construction
  to the degenerate anticanonical curve $xz=0$ (and the induced line bundle
  of degree $4$).  This has four components (two fibers and two sections),
  forming a quadrangle.  Any linear combination of the coordinates on
  $\P^3$ is determined by its values at the four points of intersection of
  the quadrangle, and thus to specify a linear matrix $B$, it is equivalent
  to specify a quadruple $(B_0,B_1,B_2,B_3)$ of scalar matrices.  If these
  matrices are invertible, then the restriction of $\coker(B)$ to $xz=0$ is
  determined by its restrictions to the four components, and thus by the
  conjugacy classes of the matrices
\[
B_0^{-1}B_1,B_1^{-1}B_2,B_2^{-1}B_3,B_3^{-1}B_0.
\]
In other words, the problem of classifying sheaves on $\P^1\times \P^1$
with the appropriate kind of presentation and restriction to the quadrangle
is equivalent to the problem of classifying quadruples in $\GL_n(k)$ with
specified conjugacy classes and product $1$.  This is the four-matrix case
of the ``multiplicative Deligne-Simpson problem'',
\cite{Crawley-BoeveyW/ShawP:2006}.
\end{eg}

\begin{eg}
Of course, we can obtain the three-matrix version of the multiplicative
Deligne-Simpson problem by insisting that $B_3=B_0$, or in other words that
the sheaf meet the relevant component in $\sO_p^n$ where $p$ is the point
representing $1$.  If we blow this point up and blow down both the fiber
and the section containing it, we obtain a moduli problem on $\P^2$
concerning sheaves with specified restriction to the triangle $xyz=0$.
\end{eg}

There does not appear to be any way to translate more general multiplicative
Deligne-Simpson problems into the Poisson surface framework.  For the
additive problem, the situation is nicer.

\begin{eg}
  The Hirzebruch surface $F_d$ for $d\ge 1$ has a unique section $s_{\min}$
  of negative self-intersection ($s_{\min}^2=-d$).  Given any $d+2$
  distinct fibers $f_0$,\dots,$f_{d+1}$, the divisor $2s_{\min}+\sum_i f_i$
  is anticanonical, so we may take it as our curve $C_\alpha$.  Given any
  $(d+2)$-tuple of matrices $C_0$,\dots,$C_{d+1}$ with $\sum_i C_i=0$, we
  have a natural corresponding matrix with coefficients in
  $\sO_X(s_{\min}+df)$.  Indeed, we can coordinatize $F_d$ in terms of a
  weighted projective space with a generator $y$ of degree $d$, and
  consider the matrix
\[
B = y + C(x,w)
\]
such that $C(x,w)$ is equal to $C_i$ on $f_i$.  Then the restriction of $B$
to $f_i$ is given by the conjugacy class of $C_i$.  In this way, we obtain
the $(d+2)$-matrix additive Deligne-Simpson problem (classifying $(d+2)$-tuples
of matrices with specified conjugacy classes and sum $0$).
\end{eg}

\begin{rem}
This, of course, is closely related to the problem of classifying Fuchsian
differential equations with specified singularity structure; if we blow up
a point of each fiber $f_i$ then blow down $f_i$, we can in this way
eliminate all components of $C_\alpha$ but the strict transform of $s_{\min}$.
If we choose the centers of the elementary transformations carefully, we
can arrange to end up at $F_2$, and case $0$ above.
\end{rem}

\section{The spectral curve}

Although most of what we will have to say in the relaxed (i.e.,
commutative) case will extend to actual equations, there is one structure
available in the commutative case that does not so extend: since the sheaf
$M$ has $1$-dimensional support, it determines a curve in $X$ which we call
the spectral curve by analogy with the relaxed differential case.  We need
to be somewhat careful here, as the standard notion of the support of a
sheaf can behave badly in families.  Luckily, there is a notion that
remains flat for flat families, namely the $0$-th Fitting scheme.  In
general, if $M$ is a sheaf on a scheme $X$, $\Fitt_0(M)$ is computed from a
presentation $V\to W\to M$ (with $V\to W$ a map of vector bundles) by
taking the zero locus of the exterior power $\wedge^n V\to \wedge^n W$
where $\rank(W)=n$.  (This is independent of the choice of presentation,
and the usual scheme-theoretic support is the reduced subscheme of
$\Fitt_0(M)$, see \cite[\S20.2]{EisenbudD:1995}.)  In particular, we see
that in our case $\Fitt_0(M)$ is nothing other than the curve $\det(B)=0$.

In the typical case (as we will see), $\Fitt_0(M)$ is integral and $M$ is
the direct image of a line bundle on $\Fitt_0(M)$, and thus much of the
structure of $M$ (in particular, the structure of $M|_{C_\alpha}$ and thus
the singularities) can be recovered from the spectral curve.  For instance,
the structure of the singularities are determined by a choice of Chern
class on the blowup $\tilde{X}$ that resolves the intersections between the
spectral and anticanonical curves.  (The analogous resolution will be
extremely important in the nonrelaxed case, and thus one of our goals below
will be to see how to construct it {\em without} reference to the spectral
curve.)  Moreover, since the spectral curve behaves well in flat families,
it determines a family of curves over the moduli space of sheaves, which
tends to interact well with the Poisson structure on the moduli space.  (In
particular, it induces Lagrangian fibrations of the symplectic leaves.)

It is thus worth considering how to compute the spectral curve in general.
In principle, one can do this by following the above algorithm to determine
$B$, and indeed this is the only {\em correct} way to obtain the spectral
curve, but one can obtain a good approximation to the spectral curve
directly from the equation, which can at most differ from the true spectral
curve by some collection of fibers.

The idea is that if we replace $M$ by some $M'$ coming from a non-minimal
factorization of $A$, then $\Fitt_0(M')-\Fitt_0(M)$ is the pullback of a
divisor on $C$, and thus the two divisors become the same after removing
all fibers.  Since the minimization procedure is more involved than the
factorization procedure, this suggests that we should aim to compute that
common result.

This works particularly well in the differential and symmetric difference
cases.

\begin{eg}
  Given a meromorphic Higgs bundle
  \[
  A:V^*\to V^*\otimes \omega_C,
  \]
  the minimal factorization has the form $A = B_\infty^{-t} B_0^t$ with
  \[
  B_0:W\to V,\quad
  B_\infty:W\to V\otimes \omega_C^{-1},
  \]
  corresponding to a matrix $B=B_0y-B_\infty$.  The spectral curve is thus
  \[
  \det(B_0 y - B_\infty) = 0.
  \]
  Since we are ignoring fiber components, we may feel free to multiply this
  by any nonzero section of a line bundle on $C$.  Dividing by
  $\det(B_\infty)$ gives
  \[
  \det(B_0 y - B_\infty)/\det(B_\infty)
  =
  \det((B_0 B_\infty^{-1}) y - 1)
  =
  \det(A^t y - 1)
  =
  \det(A y -1).
  \]
  In other words, the ``horizontal'' part of the spectral curve is nothing
  other than the horizontal part of $\det(Ay-1)$.  (This agrees with the
  usual notion of spectral curve for Higgs bundles.)
\end{eg}

\begin{eg}
  In the ordinary difference case $v(z+\hbar)=A(z)v(z)$, the minimal
  factorization has the form $A=B_\infty^{-t}B_0^t$ with
  $B_0 = B(0,x,w)$, $B_\infty=B(w^2,x,w)$, so that
  \[
  B = B_0 + (B_\infty-B_0)y/w^2.
  \]
  Dividing by $\det(B_\infty)$ tells us that the spectral curve is
  the horizontal part of the divisor given on the affine patch $w=1$ by
  \[
  \det(y-(y-1)A(x)).
  \]
  Similarly, in the $q$-difference case, the spectral curve is the
  compactification of the horizontal part of the divisor of
  $\det(y-(y-x)A(x))$ on the patch $w=1$, while in the elliptic difference
  case the curve is $\det(uA(z)-v)$ (where the anticanonical curve is
  $uv=0$ on $E\times \P^1$).
\end{eg}  

In the symmetric cases, the difficulty is that it is still somewhat
nontrivial to obtain {\em some} factorization (the algorithm is simple, but
there is no obvious easy-to-analyze choice).  The key idea here is that any
symmetric equation can also be interpreted as a nonsymmetric equation, and
the minimal factorization as a symmetric equation is still a
factorization as a nonsymmetric equation (though possibly non-minimal).
Since we can easily compute the spectral curve as a nonsymmetric equation,
we reduce to understanding how the two spectral curves are related.

\begin{eg}
  Consider a symmetric difference equation $v(z+\hbar)=A(z)v(z)$ with
  $A(-z)A(z)=1$.  Given a symmetric factorization
  \[
  A = B(-z)^{-t} B(z)^t,
  \]
  the spectral curve of the corresponding nonsymmetric relaxed equation (on
  the appropriate affine patch) is
  \[
  \det(Y-(Y-1)A(X))
  \propto
  \det((1-Y)B(X) + Y B(-X)).
  \]
  On the other hand, as a symmetric relaxation, the matrix $\hat{B}$ on $F_2$
  is given by $y B_1(x)+B_0(x)$ where $B(z)=B_0(z^2)+zB_1(z^2)$.  (Here we
  have taken the anticanonical curve to be $y^2=x$ on the patch $w=1$.)
  Since
  \[
  (1-Y)B(X) + Y B(-X)
  =
  ((1-2Y)X)B_1(X^2) + B_0(X^2),
  \]
  we see that the nonsymmetric spectral curve is obtained from the
  symmetric spectral curve by the substitution $(y,x)\mapsto
  ((1-2Y)X,X^2)$.  But this can be reversed by taking $X=\sqrt{x}$,
  $Y=1/2-(y/2\sqrt{x})$.  Geometrically, what is going on is that the
  symmetry of $A$ ensures that the horizontal part of the nonsymmetric
  spectral curve is invariant under the involution $(Y,X)\mapsto (1-Y,-X)$,
  and the symmetric spectral curve is nothing other than the quotient under
  that involution.  Similarly, for the symmetric $q$-difference case, the
  two curves are related by $(y,x)\mapsto ((1-1/X^2)Y+(1-X),X+1/X-2)$,
  where the symmetric anticanonical curve is taken to be $y^2+xy=x$.
  Finally, in the symmetric elliptic case, the spectral curve is the
  quotient by $(z,u,v)\mapsto (-z,v,u)$ of the nonsymmetric spectral curve
  on $E\times \P^1$, with an appropriate choice of map from the quotient of
  $E\times \P^1$ to $F_2$.
\end{eg}

\section{Isospectral transformations}

One thing we see immediately from the various formulas for the horizontal
part of the spectral curve is that it depends only on the conjugacy class
of $A$, which leads to the following definition (again standard in the
differential case).

\begin{defn}
  A {\em weak isospectral transformation} between two relaxed differential
  or difference equations represented by matrices $A$, $\tilde{A}$ is an
  invertible meromorphic matrix $C$ (invariant under the involution in the
  symmetric case) such that $\tilde{A}(z) = C(z) A(z) C(z)^{-1}$.
\end{defn}

\begin{rem}
  This definition is somewhat too restrictive in the symmetric case, in
  that there are natural examples of isospectral transformation as
  nonsymmetric equations in which the matrix $C$ is not invariant under the
  involution, with the two most significant examples being
  \[
  \tilde{A} = B^t \eta^* B^{-t} = B^t A B^{-t}
  \]
  and
  \[
  \tilde{A} = A = A A A^{-1}.
  \]
  Note that in both cases, $C$ satisfies a weaker symmetry condition of the
  form $\eta^* C = A^l C$ for some $l$.  The role of $l$ is clearer in the
  nonrelaxed analogue: these gauge transformations change the involution,
  shifting the corresponding degree 2 divisor by $q^l$.  In any event,
  since there is an essentially canonical choice for each $l$ (taking $C$
  to be a power of $A$ or $B^t$ times a power of $A$), we can safely ignore
  these for the present.
\end{rem}

In general, the matrix $C$ is a meromorphic map of vector bundles
$V^t\ratto \tilde{V}^t$, which allows us to describe weak isospectral
transformations directly in terms of $B$: we replace $B$ by $C^{-t}B$
then replace $W$ by the minimal $\tilde{W}$ making $B$ holomorphic.
Factoring both rational maps of vector bundles on $\P^1$ gives us a
commuting diagram
\[
\begin{CD}
  \rho^*V(-1) @<<< \rho^*\hat{V}(-1) @>>> \rho^*\tilde{V}(-1)\\
  @V{B}VV  @V{\hat{B}}VV @V{\tilde{B}}VV\\
  \rho^*W@>>> \rho^*\hat{W} @<<< \rho^*\tilde{W},
\end{CD}
\]
where the horizontal maps are pulled back from $\P^1$.  Taking cokernels of
the vertical maps gives sheaves $M$, $\hat{M}$, $\tilde{M}$ such that each of
$M$ and $\tilde{M}$ is the middle homology of some complex
\[
\rho^*Z(-1)\hookrightarrow \hat{M}\twoheadrightarrow \rho^*Z'
\]
in which $Z$ and $Z'$ are torsion sheaves on $\P^1$.  More generally, if
$M$ and $\hat{M}$ are two sheaves, possibly corresponding to non-minimal
factorizations, then an expression of $M$ as the middle homology of such a
complex induces a weak isospectral transformation between the two relaxed
equations.

\begin{prop}
  Let $M$ and $\hat{M}$ be two $\rho_*$-acyclic, relatively globally
  generated sheaves with 1-dimensional support and torsion-free direct
  image, and let $U$ be the preimage of an open subset of the curve $C$.
  Then an isomorphism $M|_U\cong \hat{M}|_U$ induces an isospectral
  transformation between the corresponding relaxed equations.
\end{prop}

\begin{proof}
  The map $M|_U\to \hat{M}|_U$ can be represented by a morphism of the form
  $M\to \hat{M}\otimes \rho^*\sO_{\P^1}(D)$ for some divisor $D$ supported
  on the complement of the image of $U$.  Moreover, the natural morphism
  $\hat{M}\to \hat{M}\otimes \rho^*\sO_{\P^1}(D)$ certainly induces an
  isospectral transformation (in which the matrix $C$ is a multiple of the
  identity), so that we may reduce to the case of that the local
  isomorphism is a morphism.  (Note that the twist of $\hat{M}$ satisfies
  the same hypotheses.)  We can then factor that morphism through its image
  (which again satisfies the hypotheses) to reduce to the cases of
  surjections and injections.

  For an injection, we have a short exact sequence
  \[
  0\to M\to \hat{M}\to T\to 0,
  \]
  where $T$ must be supported on the complement of $U$.  Since $\hat{M}$ is
  $\rho_*$-acyclic and relatively globally generated, $T$ itself is
  $\rho_*$-acyclic and relatively globally generated, so that we have a
  surjection of the form $\rho^*\rho_*T\to T$.  Pulling back the short
  exact sequence gives
  \[
  0\to M\to \tilde{M}\to \rho^*\rho_*T\to 0.
  \]
  Moreover, the kernel of $\rho^*\rho_*T\to T$ is, as usual, a sheaf of the
  form $\rho^*Z(-1)$, so that we also have a short exact sequence
  \[
  0\to \rho^*Z(-1)\to \tilde{M}\to \hat{M}\to 0.
  \]
  But this is precisely what we need to induce a weak isospectral
  transformation from $M$ to $\hat{M}$.

  For the surjective case, we have
  \[
  0\to T\to M\to \hat{M}\to 0
  \]
  with $\Hom(T,\sO_f)=0$, so that (either using the five-term sequence or
  Cohen-Macaulay duality) we have a short exact sequence of the form
  \[
  0\to T\to \rho^*Z_1(-1)\to \rho^*Z_2\to 0,
  \]
  and can then proceed as in the injective case.
\end{proof}

\begin{cor}
  Any sheaf $M$ with presentation of the standard form admits a weak
  isospectral transformation to a sheaf for which the spectral curve has no
  vertical components.
\end{cor}

\begin{proof}
  If the spectral curve of $M$ contains a fiber $f$, let $\hat{M}$ be the
  kernel of the natural morphism $M\to M\otimes\sO_f$.  Then the inclusion
  $\hat{M}\to M$ is an isomorphism on the complement of $f$, so induces a
  weak isospectral transformation.  Since this operation removes a copy of
  $f$ from the spectral curve, finitely many such steps suffice to produce
  a spectral curve with no vertical components as required.
\end{proof}

Since removing a fiber from the spectral curve also makes its intersection
with the anticanonical curve smaller, we see that this operation makes the
singularities simpler.  In other words, if the spectral curve of a relaxed
equation contains a fiber, then there are less singular equations in the
same isospectral equivalence class, justifying our viewing such cases as
corresponding to apparent singularities (i.e., those singularities that can
be gauged away).

This also leads to a somewhat stricter notion of isospectral
transformation, consisting of those weak isospectral transformations that
do not change the spectral curve.  In the no-fibers case, this leads to a
simpler geometric description.

\begin{prop}
  Suppose that $M$ and $\hat{M}$ are sheaves corresponding to minimal
  factorizations of equations with no apparent singularities.  Then an
  isospectral transformation between $M$ and $\hat{M}$ induces a common
  subsheaf $\tilde{M}$ such that both quotients have $0$-dimensional
  support.
\end{prop}

\begin{proof}
  Such an isospectral transformation induces an isomorphism between $M$ and
  $\hat{M}$ on the complement of a finite union of fibers, which can in
  turn be represented by a morphism $M\to \hat{M}\otimes \rho^*\sO_C(D)$,
  or equivalently $M\otimes\rho^*\sO_C(-D)\to \hat{M}$.  Since this is an
  isomorphism away from a finite subscheme of the spectral curve, it is
  injective with $0$-dimensional quotient, and thus gives the desired
  common subsheaf.
\end{proof}

\begin{rem}
  Of course, the common subsheaf from the theorem is not unique, but once
  we have identified {\em some} common subsheaf, we can then pass to the
  {\em maximal} common subsheaf.
\end{rem}

\begin{cor}
  Any relaxed equation with reduced spectral curve is isospectral to an
  equation such that the corresponding sheaf is the image of a line bundle
  on the normalization of the spectral curve.
\end{cor}

\begin{proof}
  The sheaf $M$ is invertible on the complement of the singular locus of
  the spectral curve, and thus there is an extension of its restriction to
  an invertible sheaf on the normalization.  Since these are isomorphic
  away from a finite set of points, the corresponding relaxed equations are
  isospectral.
\end{proof}

Again, this operation can be viewed as removing apparent singularities, at
least when resolving singularities of the spectral curve along the
anticanonical curve.  (We will see that this indeed makes the singularities
simpler, in that it enables one to resolve the intersections between the
spectral and anticanonical curves.)  We will give a more natural (and more
general) version of this construction in Lemma
\ref{lem:resolving_pseudotwist} below.

It is important to note that for the application to integrable systems, we
are actually interested in {\em families} of isospectral transformations;
at the very least, we have an equation defined over some function field
(typically the function field of some moduli space of equations) and need
to understand the isospectral transformations that are defined over that
field.  So in general we are not interested in isospectral transformations
per se, but rather natural constructions of isospectral transformations.
(Moreover, since the spectral curve does not make sense for actual
equations, we will want such constructions that can be defined without
reference to the spectral curve!)  Our primary example of this are the
``pseudo-twists'' considered below, which by Corollary
\ref{cor:pseudo_twist_is_isospectral} indeed give isospectral
transformations.

\section{Quadratic and other transformations}

Our method for computing the horizontal spectral curve in the symmetric
case involved forgetting the symmetry.  This operation induces a morphism
between moduli spaces, for instance from the moduli space of symmetric
elliptic difference equations to that of all elliptic difference equations.
On the surface $C\times \P^1$, consider the involution $\eta\times
(t\mapsto t^{-1})$.  This preserves our standard choice of anticanonical
curve, and more precisely preserves the corresponding Poisson structure (in
contrast to $\eta\times 1$, say, which {\em negates} the Poisson
structure).  It follows that this involution acts on the corresponding
moduli space, and the symmetric equations are fixed points of this
involution.  (Moreover, the involution preserves the Poisson structure on
the moduli space, so the fixed locus inherits a Poisson structure.)  There
are some difficulties in studying symmetric equations from this
perspective, however.  One is that, as we have seen, the notion of
singularity should really take into account the symmetry, but another is
that when working with moduli spaces, a fixed point merely indicates a
sheaf which is {\em isomorphic} to its image under the symmetry.  Since a
sheaf only determines a matrix up to a choice of basis, not every point of
the fixed locus actually corresponds to a symmetric equation.  (The
situation is not too dire, though: symmetric equations form a component of
the fixed locus.)  Note that the quotient of $E\times \P^1$ by the above
involution is still an elliptic surface (with constant $j$ invariant, and
with two $I_0^*$ fibers in characteristic not $2$), and thus must be blown
down eight times to reach a Hirzebruch surface.  (We can arrange to reach
the usual $F_2$ constructed from $(C,\eta)$, in which case the map from the
elliptic surface blows up each fixed point of $\eta$ twice.)  This reflects
both the fact that the notion of singularity changes and the fact that the
fixed points of the moduli space are the equations which are symmetric {\em
  up to isomorphism}.

Similar comments apply if we try to relate symmetric and nonsymmetric
difference equations in the $q$-difference and ordinary difference cases.
More generally, we could consider any Poisson involution on one of our
Poisson Hirzebruch surfaces.  (Note that in characteristic not 2 such an
involution has isolated fixed points outside $C_\alpha$, so must act as a
nontrivial involution on $F_2$ not fixing any fiber, so up to changes of
variables has the form $(y,x,w)\mapsto (-y,w,x)$.)  We find that the most
general Poisson involution (in characteristic not 2) is again at the
elliptic level, and is simply given by translation by a $2$-torsion point
$p$ of $\Pic^0(C)$; any other Poisson involution on $F_2$ is a degeneration
of this (on some degenerate curve).  Given any symmetric elliptic
difference equation on the isogenous curve $C/\langle p\rangle$, we can
interpret it as an equation on $C$ (typically with twice as many
singularities), and the corresponding sheaf will be invariant under the
Poisson involution.  Since any Poisson involution degenerates this, it in
particular follows that the embedding of symmetric equations in the moduli
space of nonsymmetric equations is a degeneration of this ``quadratic
transformation''.  (So called because at the bottom, differential level,
that is precisely what it is: performing a quadratic change of variables in
the differential equation.)  Once again, the mismatch between the two
notions of singularity and the fact that equations can be symmetric up to
isomorphism without being symmetric is reflected in the fact that the
quotient by the involution is a (singular) del Pezzo surface of degree 4
with an $A_1A_1A_3$ configuration of $-2$-curves.  (One of the $-2$-curves
comes from the original $-2$-curve, while the other four come from fixed
points of the involution, two of which are on the original $-2$-curve.)

For the $q$-difference and differential cases, one may also consider
higher-order transformations, in which one replaces the involution by an
Poisson automorphism of higher (but still finite) order.  Again, there is
an induced Poisson structure on the fixed locus on the moduli space, but
there are issues with understanding singularities at the ramification
points of the automorphism.

There may also be some interesting phenomena related to anti-Poisson
involutions of rational surfaces (which can be identified by the fact that
they are hyperelliptic when restricted to the anticanonical curve).  Though
these remain anti-Poisson on the moduli space, they can be combined with a
natural duality operation on sheaves to again obtain a Poisson involution
on the moduli space.  One example of this is the adjoint operation
$A\mapsto A^{-t}$, see Section \ref{sec:elldiff} below, for which the
fixed points would include relaxed versions of equations with structure
group $\text{SO}$ or $\text{Sp}$.

Note that in each case, although there are subtleties about how to
interpret singularities, it is clear (or will be following our
classification of singularities in Section \ref{sec:singularities}) that we
should consider two fixed points as having the same singularities iff they
lie on the same symplectic leaf for the induced Poisson structure on the
fixed locus, or equivalently if they have the same restriction to a
suitably equivariant sheaf on $C_\alpha$.

\chapter{Birational classification of Poisson surfaces}
\label{chap:classify}
\section{Poisson structures}

As we have seen, a number of types of relaxed equations can be translated
into sheaves on suitable ruled surfaces with anticanonical, or roughly
equivalently, Poisson structures.  Although this will have direct
consequences already for the relaxed case, perhaps the greatest
significance comes from the fact that Poisson structures are essentially
infinitesimal versions of {\em noncommutative} structures, so that this is
our first indication of a possible connection to noncommutative geometry.

It will be helpful to recall some standard facts about Poisson structures.
As usual in geometry, we start with the case of algebras.  Recall that a
Poisson algebra over a commutative ring $R$ is a commutative $R$-algebra
$A$ equipped with an $R$-bilinear map $\{,\!\}:A\times A\to A$ which is a
biderivation (a derivation in each variable) and satisfies the identities
\[
\{f,f\}=0,\qquad
\{f,\{g,h\}\}+\{g,\{h,f\}\}+\{h,\{f,g\}\}=0
\]
for all elements $f$, $g$, $h$.  (Note that we make no assumption about
characteristic here.)  A map $\pi:A\to B$ of Poisson $R$-algebras is called
a Poisson map if
\[
\pi(\{f,g\}) = \{\pi(f),\pi(g)\};
\]
we will also need to consider {\em anti-Poisson} maps, such that
\[
\pi(\{f,g\}) = -\{\pi(f),\pi(g)\}.
\]
(One could more generally consider maps that multiply the Poisson structure
by some other scalar, but these will be less important for us.)

Poisson structures originally arose in classical physics, which we now
understand as coming from the noncommutativity of observables in quantum
mechanics, or, more algebraically, from degenerations of noncommutative
algebras to commutative algebras.  To make this precise, let $R$ be a
discrete valuation ring (dvr) with uniformizer $t$, say $k[[t]]$ (though
the mixed characteristic case works equally well), and let $A$ be a flat
$R$-algebra such that $A_0:=A\otimes_R k$ is commutative.  Then $A$ induces
a binary operation $\{,\}$ on $A_0$ as follows: given $x,y\in A$ with
images $\bar{x},\bar{y}\in A_0$, we define
\[
\{\bar{x},\bar{y}\} := \overline{t^{-1}[x,y]}
\]
If $\bar{x}=\bar{z}$, then we can write $z=x+tv$, and thus
\[
t^{-1}[z,y] = t^{-1}[x,y] + [v,y];
\]
since $A_0$ is commutative, $[v,y]\in tA$ and thus $\overline{[v,y]}=0$.
Moreover, since $[x,y]$ is an alternating biderivation on $A$, $\{,\}$ is
an alternating biderivation on $A_0$; similarly, the fact that commutators
satisfy the Jacobi identity implies that $\{,\}$ satisfies the Jacobi
identity.  In other words, $\{,\}$ makes $A_0$ a Poisson algebra.  (If
$\{,\}$ vanishes identically, one can instead consider the reduction of
$t^{-2}[x,y]$, etc.  This will eventually give a nontrivial Poisson
structure as long as $A$ is not itself commutative.)

This is not quite a universal construction (there can be obstructions to
lifting a Poisson structure to a noncommutative structure, especially once
we have globalized to schemes), but still leads one to hope that any given
(nontrivial!) Poisson structure should lift.  (We will do this for Poisson
structures on surfaces; for Poisson structures on moduli spaces, this is a
very interesting open question!)

One advantage of working with Poisson structures directly is that they are
relatively straightforward to globalize.  By universality of K\"ahler
differentials, to specify an alternating $R$-linear biderivation on $A$, it
is equivalent to specify a homomorphism
\[
\wedge^2\Omega_{A/R}\to A.
\]
of $A$-modules.  Moreover, if $\{,\!\}$ is an alternating biderivation,
then the left-hand side of the Jacobi identity is an alternating
triderivation, which is required to vanish.  These notions localize in a
natural way, so one may thus define a Poisson structure on a scheme $X/S$
to be a morphism
\[
\alpha:\wedge^2\Omega_{X/S}\to \sO_X
\]
of $\sO_X$-modules such that the associated map
\[
\wedge^3\Omega_{X/S}\to \sO_X
\]
(defined locally in the obvious way) is 0.  Note that this induces a
Poisson algebra structure on $\Gamma(U,\sO_X)$ for all open subsets $U$.

One can define Poisson and anti-Poisson morphisms of schemes accordingly.
In scheme-theoretic terms, a Poisson morphism $(X,\alpha_X)\to (Y,\alpha_Y)$ is
a morphism $f:X\to Y$ such that the diagram
\[
\begin{CD}
f^*\wedge^2\Omega_{Y/S}@>f^*\alpha_Y>> f^*\sO_Y\\
   @VVV @| \\
\wedge^2\Omega_{X/S}   @>\alpha_X>> \sO_X
\end{CD}
\]
commutes, and similarly, with appropriate sign changes, for anti-Poisson
morphisms.  Note that in general the identity morphism gives an
anti-Poisson morphism $(X,\alpha)\cong (X,-\alpha)$.  (We will also have
occasion to use the obvious analogues for a general biderivation.)

\begin{prop}
  Let $Y$ be a Poisson scheme, and suppose $f:X\to Y$ is a morphism such
  that $f^*\Omega_{Y/S}\to \Omega_{X/S}$ is an isomorphism (e.g. if $f$ is
  \'etale).  Then there is a unique Poisson
  structure on $X$ making $f$ Poisson.
\end{prop}

\begin{proof}
Clearly, the biderivation on $X$ must be
\[
\wedge^2 \Omega_{X/S}\to \wedge^2\pi^*\Omega_{Y/S}\to \pi^*\sO_Y\cong \sO_X,
\]
which is manifestly alternating.  The isomorphism on differentials in
particular implies that a derivation vanishes iff it vanishes on functions
pulled back from $Y$ (i.e., on pullbacks of sections of $Y$ on open subsets
$U$).  In particular, it follows that the Jacobi triderivation vanishes as
required.
\end{proof}

\begin{prop}
  Let $X$ be a Poisson scheme, and suppose $f:X\to Y$ is a morphism such
  that the natural map $f_*\sO_X\to \sO_Y$ is an isomorphism.
  Then there is a unique Poisson structure on $Y$ making $f$ Poisson.
\end{prop}

\begin{proof}
  The composition
\[
f^*\wedge^2\Omega_{Y/S}\to \wedge^2\Omega_{X/S}\to \sO_X
\]
induces by adjunction a map
\[
\wedge^2 \Omega_{Y/S}\to f_*\sO_X\cong \sO_Y,
\]
i.e. an alternating biderivation.  To check that this satisfies the Jacobi
identity, we may reduce to the case $S=\Spec(R)$, $Y=\Spec(A)$ for an
$R$-algebra $A$, in which case this biderivation is just the induced
biderivation on $\Gamma(X,\sO_X)\cong A$.  But this is clearly a Poisson
bracket.
\end{proof}

\begin{rem}
Examples include a projective birational morphism $\pi:X\to Y$ with $X$
integral and $Y$ normal (and locally Noetherian), as well as the morphism
$\rho:X\to C$ associated to a rationally ruled surface.  (Of course, in the
latter case, the induced Poisson structure is trivial!)
\end{rem}

One difficulty in applying the above definition is that although
alternating biderivations are generally fairly straightforward to construct
(natural constructions of $X$ tend to make $\Omega_{X/S}$ easy to
understand), the map from alternating biderivations to triderivations is
fairly complicated, so that it can be quite nontrivial to show that a given
alternating biderivation satisfies the Jacobi identity.  One helpful
simplification is that satisfying the Jacobi identity is a closed
condition.  This is certainly true insofar as it is the vanishing of a
triderivation, but in fact, there is a somewhat stronger statement one can
make.

\begin{lem}\label{lem:closure_is_Poisson}
  Let $X/S$ be a Poisson scheme, and let $f:X\to Y$ be a morphism over $S$
  such that there is a biderivation on $Y/S$ compatible with the Poisson
  biderivation on $X$.  Then the biderivation on $Y$ restricts to a
  Poisson biderivation on the Zariski closure of $f(X)$.
\end{lem}

\begin{proof}
  This reduces to a statement about algebras: if $g:A\to B$ is a
  biderivation-preserving morphism of $R$-algebras-with-biderivations and
  $B$ is Poisson, then so is $g(A)$.  But this is obvious: the
  compatibility condition tells us that $\{g(x),g(y)\} = g(\{x,y\})$ and
  thus the Poisson bracket restricts to a biderivation on $g(A)$, which
  clearly inherits alternation and the Jacobi identity.
\end{proof}

A {\em Poisson subscheme} of a Poisson scheme $(X/S,\alpha)$ is a locally
closed subscheme $Y/S\subset X/S$ which locally satisfies $\{f,{\cal
  I}_Y\}\subset {\cal I}_Y$ where ${\cal I}_Y$ is the ideal sheaf of $Y$
and $f$ is any local section of $\sO_X$.  In particular, we immediately
obtain an induced Poisson structure on $Y$ such that the corresponding
inclusion morphism is Poisson, and such a structure exists iff $Y$ is a
Poisson subscheme.

\begin{rem}
  We can still define a Poisson subscheme when the biderivation on $X$ fails
  to be alternating or satisfy the Jacobi identity, and the Lemma then
  implies that the Zariski closure of any set of Poisson subschemes is a
  Poisson subscheme.
\end{rem}

An important class of Poisson schemes are the symplectic schemes.  A {\em
  symplectic scheme} is a pair $(X,\omega)$ such that $X/S$ is smooth and
$\omega\in \wedge^2\Omega_{X/S}$ is a closed $2$-form that induces an
isomorphism $\Omega_{X/S}\cong T_{X/S}$.  Call $(X,\omega)$ {\em
  presymplectic} if $\omega$ is nondegenerate, but not necessarily closed.
Then $\omega$ induces a biderivation
\[
\alpha:\Omega_{X/S}\otimes \Omega_{X/S}\xrightarrow{1\otimes\omega}
\Omega_{X/S}\otimes T_{X/S}\to \sO_X,
\]
which, as the inverse of an alternating form, is alternating.  Conversely,
if $\alpha$ is an alternating biderivation which is a nonsingular alternating
form on $\Omega_{X/S}$ at every point, then it determines a presymplectic
structure on $X$.

\begin{prop}
Suppose $(X,\omega)$ is a presymplectic scheme.  Then $(X,\omega)$ is
symplectic iff the corresponding biderivation is Poisson.
\end{prop}

\begin{proof}
Indeed, just as in characteristic 0, $\omega$ induces an isomorphism
\[
\wedge^k\Omega_{X/S}
\cong
\wedge^k\Hom(\Omega_{X/S},\sO_X)
\cong
\Hom(\wedge^k\Omega_{X/S},\sO_X)
\]
for any $k$, which for $k=3$ takes $d\omega$ to the triderivation coming from
the Jacobi identity.
\end{proof}
%

\begin{rem}
  In characteristic 0, it is well-known that a Poisson subscheme of a
  symplectic scheme is necessarily open.  (An ideal of a symplectic local
  ring is Poisson iff it is preserved by all derivations, but this is
  impossible in characteristic 0 for a proper, nontrivial ideal.)  This
  fails in finite characteristic; indeed, in characteristic $p$, the
  principal ideal generated by a $p$-th power will always be Poisson.
\end{rem}

\section{Poisson surfaces}
\label{sec:classify_poiss_surf}

By a {\em Poisson surface}, we will mean a pair $(X,\alpha)$ where $X/S$ is
a projective surface with geometrically smooth and irreducible fibers, and
$\alpha$ is a Poisson structure on $X$ over $S$ which is not identically 0
on any fiber.  Note that since $X$ is a surface,
\[
\wedge^2\Omega_{X/S}\cong \omega_{X/S}
\]
and
\[
\wedge^3\Omega_{X/S}=0,
\]
and thus any section of the anticanonical bundle $\omega_{X/S}^{-1}$
determines a Poisson structure on $X/S$.  Thus a nontrivial Poisson
structure is determined up to a scalar by the divisor of the corresponding
section of $\omega_{X/S}^{-1}$, the {\em associated anticanonical curve},
on the complement of which the surface is symplectic.  The scalar can
itself be identified: if $C_\alpha$ denotes the associated anticanonical
curve, then the short exact sequence (which depends on $\alpha$ only via
$C_\alpha$) 
\[
0\to \omega_X\to \omega_X(C_\alpha)\to \omega_{C_\alpha}\to 0
\]
induces an injection $H^0(\omega_X(C_\alpha))\to H^0(\omega_{C_\alpha})$.
A Poisson structure with associated anticanonical curve $C_\alpha$ is a
nonzero global section of $\omega_X^{-1}(-C_\alpha)$, so an isomorphism
$\sO_X\cong \omega_X^{-1}(-C_\alpha)$.  The inverse of this isomorphism
determines a nonzero global section of $\omega_X(C_\alpha)$, and thus a
nonzero holomorphic differential on $C_\alpha$ in the kernel of the
connecting map $H^0(\omega_{C_\alpha})\to H^1(\omega_X)$.  Conversely, any
nonzero element of that kernel determines a nonzero global section
$\omega_X(C_\alpha)$, the inverse of which is a Poisson structure with
associated anticanonical curve $C_\alpha$.  Note that the differential on
$C_\alpha$ scales inversely with the Poisson structure.

As observed above, one can push Poisson structures forward through proper
birational morphisms; in the case of Poisson surfaces, one can be fairly
precise about lifting as well.  Since everything is local on the base, we
will assume $S=\Spec(\bar{k})$ for some algebraically closed field
$\bar{k}$.  Note that if $\pi:X\to Y$ is a birational morphism of smooth
projective surfaces over $\bar{k}$, then there is a (unique) effective
divisor $e_\pi$ supported on the exceptional locus such that $K_X\sim
\pi^*K_Y + e_\pi$.  (See the discussion before Corollary \ref{cor:invis_on_e_pi}
below.)

\begin{lem}
  Let $(Y,\alpha)$ be a Poisson surface over $\bar{k}$, with associated
  anticanonical curve $C$, and suppose $\pi:X\to Y$ is a birational
  morphism.  Then there is a (unique) Poisson structure on $X$ compatible
  with that of $Y$ iff the divisor $\pi^*C-e_\pi$ is effective.
\end{lem}

\begin{proof}
  Indeed, $\pi^*C-e_\pi$ is the unique anticanonical divisor agreeing with
  $C$ where $\pi$ is an isomorphism, so if there is a Poisson structure on
  $X$, its associated anticanonical curve must be $\pi^*C-e_\pi$.  In
  particular, $\pi^*C-e_\pi$ must be effective.  If it is, then it
  determines a Poisson structure up to scalar multiplication, and that
  Poisson structure agrees with $\alpha$ (again up to scalar multiplication)
  where $\pi$ is an isomorphism.  We may thus eliminate the scalar freedom
  and obtain a Poisson structure agreeing with $\alpha$ at the generic point
  of $X$.
\end{proof}

\begin{prop}
  Let $(X,\alpha_X)$, $(Y,\alpha_Y)$ be Poisson surfaces over $\bar{k}$, and
  suppose $f:X\to Y$ is a rational map which is Poisson at the generic
  points of $X$ and $Y$.  Then there exists a Poisson surface $Z$ and
  birational Poisson morphisms $g:Z\to X$, $h:Z\to Y$ such that $f = h
  \circ g^{-1}$.
\end{prop}

\begin{proof}
  Let $f=h\circ g^{-1}$ be the minimal factorization (i.e., $Z$ is the
  minimal desingularization of the graph of $f$); we need to show that $Z$
  has a compatible Poisson structure.  We can factor each of $g$ and $h$ as
  a product of monoidal transformations, and we can arrange that those
  monoidal transformations centered in points of the anticanonical curves
  come first.  We thus reduce to the case that $X$ and $Y$ are isomorphic
  on some neighborhoods of their respective anticanonical curves, and wish
  to show that $f$ is an isomorphism.  Indeed, otherwise, we have an
  identity
\[
g^*C_X-e_g = h^*C_Y-e_h
\]
of divisors on $Z$.  Since $f$ is an isomorphism on neighborhoods of the
anticanonical curves, $g^*C_X=h^*C_Y$, and thus $e_g=e_h$.  If $g$ is not
an isomorphism, then $e_g$ contains some $-1$ curve, which must be
contracted by $h$, contradicting minimality of $Z$.
\end{proof}

We in particular have a well-behaved notion of a Poisson birational map
between Poisson surfaces.  Since our objective in this note is to
understand how moduli spaces of sheaves on Poisson surfaces are affected by
Poisson birational maps, it will be convenient to give a classification of
Poisson surfaces up to Poisson birational equivalence.  (Compare
\cite{BottacinF:1995,BartocciC/MacriE:2005}, which consider minimal Poisson
surfaces, but without considering birational maps between them.)  Call a
Poisson surface {\em standard} if it is isomorphic as an algebraic surface
to the ruled surface $\P(\sO_C\oplus \omega_C)$ for some smooth curve $C$,
in such a way that the morphism $\sO_C\oplus \omega_C\to \omega_C$
determines a section disjoint from the anticanonical curve.  (These are, of
course, precisely the surfaces we considered in relation to equations!)

There is also an odd case in characteristic 2.  Let $C$ be a smooth genus 1
curve, and let $X$ be the ruled surface corresponding to the nontrivial
self-extension of $\sO_C$.  In characteristic different from 2, such a
surface has a unique anticanonical curve, but in characteristic 2, there is
actually an anticanonical pencil.  Geometrically, we may take $C$ to be in
Weierstrass form $y^2+a_1xy+a_3y = x^3+a_2x^2+a_4x+a_6$, and the vector
bundle to have transition matrix
\[
\begin{pmatrix} 1 & 0 \\ y/x &
  1\end{pmatrix}
\]
between the complement of the identity and a suitable neighborhood of the
identity.  In general, a Poisson structure on $\P(V)$ corresponds to a
global section of $S^2(V^*)\otimes \det(V)\otimes \omega_C^{-1}$, and thus
in our case an anticanonical curve is cut out by a global section of
$S^2(V^*)$.  This has transition matrix
\[
\begin{pmatrix}
  1 & 0 & y^2/x^2\\
  0 & 1 & y/x\\
  0 & 0 & 1
\end{pmatrix},
\]
a nontrivial extension of $\sO_C^2$ by $\sO_C$.  Since
$\dim\Ext^1(\sO_C,\sO_C)=1$, this has precisely two dimensions of global
sections.  If $z$, $w$ are the standard basis of the trivialization of
$V^*$ away from the identity, then the global sections of $S^2(V^*)$ are
spanned by $w^2$ and $z^2+a_1 z w + x w^2$.  Indeed, if $u$, $v$ are the
basis of the other trivialization, then the latter becomes
\[
u^2 + a_1 u v + (a_2 + a_3 y/x^2 + a_4/x + a_6/x^2) v^2,
\]
the coefficients of which are holomorphic at the identity.  Thus the
generic anticanonical curve on this surface has the form $z^2 + a_1 z + c =
x$ for some $c$.  Note that we can replace $c$ by $c+d^2+a_1 d$ for any
$d$, and thus all such anticanonical curves are related by geometric
automorphisms of $X$.  With that in mind, we call a Poisson surface of this
form ``quasi-standard''.  Note that in general, the degree 2 map from
$C_\alpha$ to $C$ induces a 2-isogeny $J(C_\alpha)\to J(C)$ which is dual
to the Frobenius $2$-isogeny.

\begin{thm}\label{thm:poisson_class}
  Any Poisson surface $(X,\alpha)$ over $\bar{k}$ such that $\alpha$ is not
  an isomorphism is Poisson birational to a standard or quasi-standard
  Poisson surface.
\end{thm}

\begin{proof}
  Since $-K_X$ is nontrivial and effective, $X$ has Kodaira dimension
  $-\infty$, so either $X\cong \P^2$ or $X$ admits a birational morphism to
  a ruled surface $\rho:X'\to C$.  We may thus reduce to the case that $X$
  is a ruled surface; if $X\cong \P^2$, simply blow up a point of the
  anticanonical curve $C_\alpha$.  Since $\rho$ admits infinitely many
  sections and $C_\alpha$ has finitely many components, we may choose a
  section that is not a component of the anticanonical curve; by mild abuse
  of notation, we call this curve $C$.

  If $C$ intersects $C_\alpha$, perform an elementary transformation
  (blow up a point, then blow down the fiber) based at a point of
  intersection.  Blowing up the point reduces the intersection
  number by 1, and, since $C$ meets the fiber only in the point
  being blown up, blowing down the fiber does not change the
  intersection number.  Thus by induction, there is a sequence of
  Poisson elementary transformations making $C$ disjoint from
  the anticanonical divisor, so we can reduce to the case that $C$
  was already disjoint from $C_\alpha$.

  In particular, $\omega_X\otimes \sO_C\cong \sO_C$, and thus we may use
  adjunction to compute $\omega_C\cong \sO_C(C)$.  Applying
  $\rho_*$ to the short exact sequence
\[
0\to \sO_X\to \sO_X(C)\to \sO_C(C)\to 0
\]
gives the short exact sequence
\[
0\to \sO_C\to \rho_*(\sO_X(C))\to \omega_C\to 0,
\]
where we observe that
\[
X\cong \P(\rho_*(\sO_X(C))).
\]
We claim that this short exact sequence splits unless the surface is
quasi-standard.  If $g(C)=0$, this is immediate.  If $g(C)=1$, then
$-K_X\sim 2C$; since we also have $-K_X\sim C_\alpha$, we conclude that
$h^0(-K_X)>1$, implying that the sequence splits (or the surface is
quasi-standard).  Finally, if $g(C)\ge 2$, then the fact that
$\chi(\sO_{C_\alpha})=0$ and $C_\alpha$ has a degree 2 map to $C$ (note
that $C_\alpha$ cannot have a fiber as a component, since that would
intersect $C$) implies that $C_\alpha$ cannot be integral.  We can thus
write $C_\alpha = C_1+C_2$; since $C_\alpha.f=2$, where $f$ is the
(numerical) class of a fiber, we must have $C_1.f=C_2.f=1$.  But then $C_1$
and $C_2$ are sections disjoint from $C$, giving the desired splitting of
$\rho_*(\sO_X(C))$.  (Such a splitting is moreover unique, so that
$C_\alpha=2C_1$.)
\end{proof}

In particular, we have the following classification of Poisson surfaces, up
to Poisson birational equivalence:
\begin{itemize}
\item[1.] $X$ is the Hirzebruch surface $F_2$, and $C_\alpha$ is a reduced
  anticanonical curve disjoint from the minimal section of $X$.
\item[2.] $X$ is $C\times \P^1$, with $C$ a smooth genus 1 curve, and
  $C_\alpha$ is a union of two disjoint fibers over $\P^1$.
\item[3.] $X$ is $\P(\sO_C\oplus \omega_C)$ for some smooth curve $C$, and
  $C_\alpha=2C_0$, where $C_0$ is the section corresponding to the morphism
  $\sO_C\oplus \omega_C\to \sO_C$.
\item[4.] $X$ is quasi-standard.
\item[5.] $X$ is a minimal surface with trivial anticanonical bundle
  ($K_3$, abelian, certain surfaces in characteristic 2 or 3 with
  nonreduced $\Pic^0$.)
\end{itemize}
Note that the reduced condition in case 1 is only there to make it disjoint
from case $3$.  Also, as we have seen, case $1$ splits into subcases,
corresponding to the Kodaira symbols $I_0$, $I_1$, $\kII$, $I_2$, $\kII$
(smooth, nodal integral, cuspidal integral, or two components meeting in a
reduced or nonreduced scheme, respectively).

Note that when the surface is not rational, either the anticanonical
divisor is trivial and there are no nontrivial Poisson birational maps, or
the (rational) ruling is uniquely determined, and thus any Poisson
birational map between surfaces of the above canonical form is a
composition of elementary transformations.  In contrast, the rational case
admits a rich structure of birational automorphisms respecting the Poisson
structure.  In particular, the $I_1$ and $I_2$ subcases are in the same
birational equivalence class, as are the $\kII$ and $\kIII$ subcases and
the nonreduced case, though there are significant differences in the
corresponding Hitchin-type systems.

We also note that the surfaces birational to quasi-standard Poisson
surfaces can be characterized as those Poisson surfaces with rational
rulings over genus 1 curves such that the anticanonical curve is integral.

\chapter{Lifting through birational morphisms}
\label{chap:*!}

\section{The minimal lift}
\label{sec:minimal_lift}

Much of what we have to say about the minimal lift construction applies to
arbitrary surfaces, so we will for the moment drop the assumption that our
surfaces are Poisson.  They also involve only local considerations, so we
drop the requirement that the surfaces be projective, but also restrict to
geometric fibers.  Thus for the purposes of this chapters, an ``algebraic
surface'' will mean an irreducible smooth $2$-dimensional scheme over an
algebraically closed field.  In addition, when we say that $\pi:X\to Y$ is
a birational morphism, we will mean that $X$ and $Y$ are algebraic surfaces
and $\pi$ is projective.

Given a birational morphism $\pi:X\to Y$, there is an obvious way
to transport coherent sheaves on $Y$ to coherent sheaves on $X$,
namely the inverse image functor $\pi^*$.  This is very well-behaved
on sheaves of homological dimension $\le 1$ (so in particular for sheaves
with standard presentation on a ruled surface).

\begin{lem}\label{lem:upper_star_in_hd_1}
  Let $\pi:X\to Y$ be a birational morphism.  Then any coherent sheaf $M$
  on $Y$ of homological dimension $\le 1$ is $\pi^*$-acyclic.  Moreover,
  the sheaf $\pi^*M$ is $\pi_*$-acyclic of homological dimension $\le 1$,
  and the natural map $M\to \pi_*\pi^*M$ is an isomorphism.
\end{lem}

\begin{proof}
If $\pi$ is monoidal and $M=\sO_Y$, this is a standard fact about
birational morphisms of surfaces.  The result for general $\pi$ follows by
induction along a factorization of $\pi$ into monoidal transformations, and
the case $M$ locally free follows easily.

In general, choose a locally free resolution
\[
0\to V\xrightarrow{B} W\to M\to 0.
\]
Applying $\pi^*$ gives an exact sequence
\[
\pi^*V\xrightarrow{\pi^*B} \pi^*W\to \pi^*M\to 0,
\]
with $\pi^*V$ and $\pi^*W$ locally free.  Since $B$ is injective, $\pi^*B$
is still injective on the generic fiber; since $\pi^*V$ is locally free, we
conclude that $\pi^*B$ is injective, and thus $M$ is $\pi^*$-acyclic as
required.  The remaining claim follows upon taking the direct image of the
resulting short exact sequence.
\end{proof}

Regarding sheaves of homological dimension $\le 1$, we have the following
characterization.

\begin{prop}\label{prop:hd1_is_pointless}
A coherent sheaf $M$ on an algebraic surface $X$ has homological
dimension $\le 1$ iff it has no subsheaf of the form $\sO_p$ for
$p$ a closed point of $X$.
\end{prop}

\begin{proof}
  We note that $M$ has homological dimension $\le 1$ iff
  $\Tor_2(M,\sO_{\overline{x}})=0$ for all (not necessarily closed) points
  $x\in X$.  If $\dim(\overline{x})=2$, then $\sO_{\overline{x}}=\sO_X$, so
  there is no condition, while if $\dim(\overline{x})=1$, then
  $\overline{x}$ is a divisor on $X$, and its structure sheaf has
  homological dimension $1$.  In other words, we find that $M$ has
  homological dimension $\le 1$ iff $\Tor_2(M,\sO_p)=0$ for all closed
  points $p$.  Now, the minimal resolution of $\sO_p$ in the local ring at
  $p$ is self-dual, and thus we have an isomorphism
\[
\Tor_2(M,\sO_p)\cong \sHom_X(\sO_p,M).
\]
Since
$\sHom_X(\sO_p,M)$ is supported on $p$, $\sHom_X(\sO_p,M)=0$ iff
$\Hom(\sO_p,M)=0$, iff $\sO_p$ is not a subsheaf of $M$.
\end{proof}

\begin{rems}
  One can prove Lemma \ref{lem:upper_star_in_hd_1} directly from this
  characterization, which in fact leads to a stronger statement: a sheaf
  $M$ is $\pi^*$-acyclic as long as it has no subsheaf supported on the
  exceptional locus of $\pi$.
\end{rems}

\begin{rems}
  This characterization extends more easily to the noncommutative setting:
  there is no particularly nice notion of vector bundles on noncommutative
  surfaces, but there {\em is} a good notion of sheaves with
  $0$-dimensional support.
\end{rems}

Since $\pi_*\pi^*M=0$, if $\pi_*$ had a right adjoint $\pi^!$, then
the isomorphism $\pi_*\pi^*M\cong M$ would induce morphisms $\pi^*M\to
\pi^!M$ and $M\to \pi_*\pi^!M$.  Of course, since $\pi_*$ is not
right exact, it cannot have a right adjoint, but duality theory
gives an adjoint to the derived functor.  In our case, we actually
obtain something stronger.

\begin{lem}
Let $\pi:X\to Y$ be a birational morphism of algebraic surfaces.
Then there exists a right exact functor $\pi^!:\coh(Y)\to \coh(X)$,
acyclic on locally free sheaves, such that for any coherent sheaves
$M$, $N$ on $X$ and $Y$ respectively (or bounded complexes of such
sheaves),
\[
\dR\sHom(\dR\pi_*M,N) \cong \dR\sHom_Y(M,\dL\pi^!N).
\]
Moreover, there is a natural isomorphism
\[
N\cong \dR\pi_*\dL\pi^!N.
\]
\end{lem}

\begin{proof}
Since $\pi_*$ is proper, $\pi_*=\pi_!$, and thus $\dR\pi_*$ has a
right adjoint $\pi^!$ on the derived category.  Since $X$ and $Y$ are smooth,
\[
\omega_X
\cong
(\phi\circ \pi)^!\sO_{\bar{k}}[-2]
\cong
\pi^!(\phi^!\sO_{\bar{k}}[-2])
\cong
\pi^!\omega_Y,
\]
where $\phi:Y\to \Spec(\sO_{\bar{k}})$ is the structure morphism
of the $\bar{k}$-scheme $Y$.  But then for any locally free sheaf $V$ on $Y$,
\[
\pi^!V
\cong
\pi^!(\omega_Y\otimes(\omega_Y^{-1}\otimes V))
\cong
\pi^!\omega_Y\otimes \dL\pi^*(\omega_Y^{-1}\otimes V)
\cong
\pi^*V\otimes \omega_X\otimes \pi^*\omega_Y^{-1}
\]
and thus in particular $\pi^!V$ is a sheaf.  Using a locally free
resolution, we conclude that
\[
\pi^!N^\cdot \cong \dL\pi^*N^\cdot \otimes \omega_X\otimes \pi^*\omega_Y^{-1}
\]
for any complex $N^\cdot$.  In other words, the derived functor $\pi^!$
is the left derived functor of the right exact functor
\[
N\mapsto \pi^*N\otimes \omega_X\otimes \pi^*\omega_Y^{-1}.
\]

For the remaining claim, we have
\[
\dR\pi_*\dL\pi^!N
\cong
\dR\sHom_Y(\sO_X,\dL\pi^!N)
\cong
\dR\sHom(\dR\pi_*\sO_X,N)
\cong
\dR\sHom(\sO_Y,N)
\cong
N.
\]
\end{proof}

\begin{rem}
More generally, if $\pi:X\to Y$ is a proper morphism of Gorenstein
varieties, then
\[
\pi^!N\cong \dL\pi^*N\otimes \omega_X\otimes
\pi^*\omega_Y^{-1} [\dim(X)-\dim(Y)];
\]
this is essentially an observation
of Deligne \cite[Prop.~7]{DeligneP:1966}.
Similarly, acyclicity on sheaves of homological dimension 1 and the
isomorphisms $\dR\pi_*\dL\pi^*N\cong N$, $\dR\pi_*\dL\pi^!N\cong N$ hold
for any proper birational morphism between smooth varieties.  The
next lemma encapsulates the key feature of the surface case, which
is much less common for birational morphisms in higher dimension.
\end{rem}

\begin{lem}\label{lem:pi*_dim_1}
Let $\pi:X\to Y$ be a birational morphism of algebraic surfaces.
Then for any sheaf $M$ on $X$, we have $R^i\pi_*M=0$ for $i\ge 2$.
In particular, any quotient of a $\pi_*$-acyclic sheaf is $\pi_*$-acyclic.
\end{lem}

\begin{proof}
Indeed, the fibers of $\pi$ are 1-dimensional, giving the vanishing result,
and the claim about quotients follows by considering the long exact sequence.
\end{proof}

Since $\pi^!$ is a twist of $\pi^*$, we find that sheaves of
homological dimension $1$ are also acyclic for $\pi^!$.  Thus if 
$M$ is a sheaf of homological dimension $1$, we have a natural
isomorphism
\[
M\cong R\pi_*L\pi^!M\cong R\pi_*\pi^!M\cong \pi_*\pi^!M.
\]
Either this or the isomorphism $\pi_*\pi^*M\cong M$ induces by adjunction
a natural morphism
\[
\pi^*M\to \pi^!M,
\]
which will be an isomorphism outside the exceptional locus of $\pi$.

\begin{defn}
Let $\pi:X\to Y$ be a birational morphism of algebraic surfaces,
and let $M$ be a coherent sheaf on $Y$ with homological dimension
$\le 1$.  Then the {\em minimal lift} of $M$ is the image $\pi^{*!}M$
of the natural map $\pi^*M\to \pi^!M$.
\end{defn}

\begin{rem}
In the theory of local systems, there is a functor ``middle extension'',
the image of the natural transformation $j_!\to j_*$ where $j$ is an open
immersion \cite{KatzNM:1990}.  Given the close connection between local
systems and differential equations and the close connection between
(relaxed) differential equations and sheaves, this suggests that there
should be more than just a formal analogy relating this to the minimal
lift.  (One distinction between the two is that our operation is defined on
sheaves, while middle extension is defined on perverse sheaves.)
\end{rem}

To justify the name, we have the following quasi-universal property.

\begin{prop}
Let $N$ be a sheaf of homological dimension $\le 1$ on $Y$, and
suppose $M$ is a $\pi_*$-acyclic sheaf on $X$ with direct image $N$.
Then $M$ has a natural subquotient isomorphic to $\pi^{*!}N$.
\end{prop}

\begin{proof}
The isomorphism $\pi_*M\cong N$ and its inverse induce by adjunction
natural maps
\[
\pi^*N\to M\to \pi^!N.
\]
It remains only to show that the composition agrees with the canonical
map $\pi^*N\to \pi^!N$.  Since we obtained the factors by adjunction,
it follows that the composition
\[
N\to \pi_*\pi^*N\to \pi_*M\to \pi_*\pi^!N\to N
\]
is the identity, which implies that the map $\pi^*N\to \pi^!N$ is adjoint
to the inverse of the canonical morphism $N\to \pi_*\pi^*N$, as required.
\end{proof}

Of course, we also want to know that $\pi^{*!}M$ itself satisfies this!

\begin{prop}
Let $M$ be a sheaf of homological dimension $\le 1$ on $Y$.  Then
$\pi^{*!}M$ is $\pi_*$-acyclic, and there is a natural
isomorphism $\pi_*\pi^{*!}M\cong M$.
\end{prop}

\begin{proof}
Since $\pi^{*!}M$ is a quotient of the $\pi_*$-acyclic sheaf $\pi^*M$,
it is certainly $\pi_*$-acyclic.  Moreover, the composition
\[
M\cong\pi_*\pi^*M\to \pi_*\pi^{*!}M\to \pi_*\pi^!M\cong M
\]
is the identity, and the map
\[
\pi_*\pi^{*!}M\to \pi_*\pi^!M
\]
is injective.  The claim follows.
\end{proof}

The minimal lift behaves well under composition of birational morphisms.

\begin{prop}
Let $\pi_1:X\to Y$, $\pi_2:Y\to Z$ be birational morphisms of
algebraic surfaces.  Then for any coherent sheaf $M$ on $Z$ with
homological dimension $\le 1$,
\[
(\pi_2\circ\pi_1)^{*!}M
\cong
\pi_1^{*!}\pi_2^{*!}M.
\]
\end{prop}

\begin{proof}
Consider the composition
\[
\pi_1^*\pi_2^*M\to
\pi_1^*\pi_2^{*!}M\to
\pi_1^{*!}\pi_2^{*!}M\to
\pi_1^!\pi_2^{*!}M\to
\pi_1^!\pi_2^!M.
\]
The second map is surjective and the third map injective by definition
of $\pi_1^{*!}$, and the first map is surjective since $\pi_1^*$ is
right exact.  It remains only to show that the fourth map is injective,
for which it will suffice to show that $\pi_2^!M/\pi_2^{*!}M$ has 
homological dimension 1.  If not, then the quotient contains the structure
sheaf of a point, and thus has nontrivial direct image under $\pi_2$, 
contradicting the fact that $\pi_2^{*!}$ and $\pi_2^!M$ are
$\pi_{2*}$-acyclic sheaves with isomorphic direct images.

It follows that $\pi_1^{*!}\pi_2^{*!}M$ is the image of the natural map
\[
(\pi_2\circ\pi_1)^*M\cong \pi_1^*\pi_2^*M\to \pi_1^!\pi_2^!M
\cong
(\pi_2\circ\pi_1)^!M
\]
as required.
\end{proof}

We also want to compare the minimal lift to the lifting operation in the
theory of Hitchin systems.

\begin{prop}
Let $\pi:X\to Y$ be a birational morphism of algebraic surfaces,
and let $C\subset X$ be a curve intersecting the exceptional locus
of $\pi$ transversely.  Then for any torsion-free coherent sheaf
$M$ on $C$, viewed as a sheaf on $X$,
\[
\pi^{*!}\pi_*M\cong M.
\]
\end{prop}

\begin{proof}
The transversality assumption implies that $\pi|_C:C\to Y$ has
$0$-dimensional fibers, so $M$ is certainly $\pi_*$-acyclic.
Moreover, if we twist $M$ by the inverse of a sufficiently ample
bundle, we can arrange that $M$ has no global sections, so that the
same applies to its direct image.  But then $\pi_*M$ certainly
cannot have a 0-dimensional subsheaf, so that $\pi_*M$ has homological
dimension $\le 1$.

It follows that $M$ has $\pi^{*!}\pi_*M$ as a subquotient.
Moreover, they are isomorphic away from the exceptional locus, so
the residual sub- and quotient sheaves of $M$ are supported there.
Transversality then makes both sheaves $0$-dimensional, so that the
torsion-free hypothesis makes the subsheaf trivial.  We thus have
a short exact sequence
\[
0\to \pi^{*!}\pi_*M\to M\to T\to 0
\]
where $T$ is $0$-dimensional.  As before, we compute $\pi_*T=0$ from
the long exact sequence, and thus $T=0$ as required.
\end{proof}

\begin{rem}
In particular, if the sheaf $N$ on $Y$ is the direct image of an invertible
sheaf on a smooth curve, then $\pi^{*!}N$ is the corresponding sheaf on the
strict transform.
\end{rem}

Another important example of minimal lifts is the following.

\begin{prop}\label{prop:*!_of_line_bundle_on_curve}
  Let $\pi:X\to Y$ be a birational morphism of algebraic surfaces, and
  suppose that $C$ is a curve on $Y$.  Then $\pi^{*!}\sO_C\cong \sO_{C'}$
  for some curve $C'$.  If $X,Y,\pi$ are Poisson and $C$ is the anticanonical
  curve on $Y$, then $C'$ is the anticanonical curve on $X$.
\end{prop}

\begin{proof}
It suffices to consider the case that $\pi$ is monoidal, with center $p$.
If $p\notin C$, then $\pi^*\sO_C\cong \pi^!\sO_C\cong \sO_{\pi^{-1}C}$,
so the result is immediate.  If $p\in C$, then over the local ring at $p$,
$\sO_C$ has a minimal resolution
\[
0\to (\sO_Y)_p\to (\sO_Y)_p\to \sO_C\otimes (\sO_Y)_p\to 0,
\]
from which we can compute
\[
\pi^{*!}\sO_C\cong \sO_{\pi^*C-e}.
\]

In particular, in the Poisson case, we must have $p\in C_\alpha$, and
$\pi^*C_\alpha-e$ is indeed the anticanonical divisor on $X$.
\end{proof}

\section{Invisible sheaves}
\label{sec:invisible_comm}

Since $\pi^{*!}M$ is naturally a quotient of $\pi^*M$ and a subsheaf of
$\pi^!M$, and adjunction makes the latter sheaves easy to deal with, this
suggests that an investigation of the corresponding kernel and cokernel
might be rewarding.

A key property of those sheaves is that they are in a sense invisible to
$\pi_*$.  To be precise, as we will see, they satisfy the following
definition.  We continue to restrict our attention to smooth surfaces over
algebraically closed fields.

\begin{defn}
  Let $\pi:X\to Y$ be a birational morphism of algebraic surfaces.  A
  coherent sheaf $E$ on $X$ is {\em $\pi$-invisible} if it is acyclic
  with trivial direct image.
\end{defn}

\begin{rem}
In other words, $E$ is $\pi$-invisible iff $\dR\pi_*E=0$.
\end{rem}

We will omit $\pi$ from the notation when it is clear from context.  Note
that an invisible sheaf is necessarily supported (set-theoretically) on the
exceptional locus, since $\pi$ is an isomorphism elsewhere.\footnote{It is
thus tempting to call these {\em exceptional} sheaves, but since we will be
working with exceptional objects in derived categories below, we use
``invisible'' to avoid confusion.}  Also, if $E$ had a $0$-dimensional
subsheaf, that subsheaf would have nontrivial direct image; we thus find
that invisible sheaves have homological dimension $1$.

\begin{prop}
Suppose $M^\bullet$ is a complex of sheaves on $X$ such that
$\dR\pi_*M^\bullet=0$.  Then every homology sheaf of $M$ is invisible.
\end{prop}

\begin{proof}
Since $\pi$ has $\le 1$-dimensional fibers, $R^p\pi_*=0$ for $p\ge 2$.
Thus the hypercohomology spectral sequence
\[
R^p\pi_*(h^q(M^\bullet))\Rightarrow R^{p+q}\pi_*M^\bullet
\]
collapses at the $E_2$ page.  Since the limit of the spectral sequence is
0, every term on the $E_2$ page is 0, and thus
\[
\pi_*(h^q(M^\bullet)) = R^1\pi_*(h^q(M^\bullet))=0.
\]
In other words, $h^q(M^\bullet)$ is invisible for all $q$.
\end{proof}

The most important special case of this for our purposes is the following.

\begin{cor}\label{cor:invisible_ker_coker}
  Suppose $f:M\to N$ is a morphism of $\pi_*$-acyclic sheaves on $X$ such
  that $\pi_*f$ is an isomorphism.  Then $\ker(f)$ and $\coker(f)$ are
  invisible.  Conversely, if $f$ has invisible kernel and cokernel,
  then $\pi_*f$ is an isomorphism; moreover, $\im(f)$ is also
  $\pi_*$-acyclic with $\pi_*M\cong \pi_*\im(f)\cong \pi_*N$.
\end{cor}

\begin{proof}
The first claim is immediate from the proposition (via the other spectral
sequence); for the second, factoring $f$ through its image reduces to the
cases that $f$ is injective or surjective.  But since the cokernel/kernel
is invisible, we obtain an isomorphism between the remaining two derived
direct images, and composing find that $\pi_*f$ is an isomorphism as
required.
\end{proof}

In particular, this implies that the sheaves $\ker(\pi^*M\to \pi^{*!}M)$ and
$\pi^!M/\pi^{*!}M$ are invisible, as we indicated above.

\begin{prop}
The category of $\pi$-invisible sheaves is closed under taking kernels,
cokernels and extensions.  In particular, it is an
abelian category.
\end{prop}

\begin{proof}
  The claim for kernels and cokernels follows from Corollary
  \ref{cor:invisible_ker_coker}.  For extensions
\[
0\to E\to E'\to E''\to 0,
\]
the long exact sequence corresponding to $\dR \pi_*$ immediately tells us
that if two of the sheaves are invisible, then so is the third.
\end{proof}

Another source of invisible sheaves is the following.

\begin{prop}\label{prop:invisible_image}
Suppose $f:M\to N$ is a morphism of sheaves on $X$ such that
$M$ is $\pi_*$-acyclic and $\pi_*N$ has homological dimension $\le 1$.
If $f$ vanishes outside the exceptional locus of $\pi$, then
$\im(f)$ is invisible, and $\pi_*f=0$.
\end{prop}

\begin{proof}
Since $M$ is $\pi_*$-acyclic, so is its quotient $\im(f)$.  Since $f$
vanishes outside the exceptional locus, $\pi_*\im(f)$ is 0-dimensional; but
$\pi_*N$ has no $0$-dimensional subsheaf.
\end{proof}

Invisible sheaves also interact nicely with the lifting operations.

\begin{prop}
Let $E$ be an invisible sheaf on $X$, and $M$ a sheaf on $Y$ with
homological dimension $\le 1$.  Then
\[
\dR\Hom(\pi^*M,E)=\dR\Hom(E,\pi^!M)=0.
\]
\end{prop}

\begin{proof}
This follows immediately from adjunction:
\begin{align}
\dR\Hom(\pi^*M,E)&\cong \dR\Hom(M,\dR\pi_*E)=0,\notag\\
\dR\Hom(E,\pi^!M)&\cong \dR\Hom(\dR\pi_*E,M)=0.\notag
\end{align}
\end{proof}

\begin{prop}
Let $E$ and $M$ be as before.  Then
\[
\Hom(\pi^{*!}M,E)=\Hom(E,\pi^{*!}M)=0.
\]
\end{prop}

\begin{proof}
Indeed,
\[
\Hom(E,\pi^{*!}M)\subset \Hom(E,\pi^!M)=0,
\]
and similarly for $\Hom(\pi^{*!}M,E)$.
\end{proof}

\begin{rem}
We will see below that this actually characterizes those sheaves which are
minimal lifts.
\end{rem}

Invisible sheaves behave well under lifts and direct images.

\begin{lem}
  Suppose $\pi:X\to Y$ and $\phi:Y\to Z$ are birational morphisms of
  algebraic surfaces.  If a sheaf $E$ on $Y$ is $\phi$-invisible, then
  $\pi^*E$, $\pi^!E$, and $\pi^{*!}E$ are $\phi\circ\pi$-invisible.  If a
  sheaf $E$ on $X$ is $\phi\circ\pi$-invisible, then it is
  $\pi_*$-acyclic and $\pi_*E$ is $\phi$-invisible.  Moreover, $\pi_*$
  and $\pi^*$ take projective objects (of the category of invisible
  sheaves) to projective objects, and $\pi_*$ and $\pi^!$ take injective
  objects to injective objects.
\end{lem}

\begin{proof}
First suppose $E$ is $\phi_*$-invisible, and let $E'$ be any of the three
lifts of $E$ to $X$.  Then we have a natural isomorphism
$\dR\pi_*E'\cong E$, and thus
\[
\dR(\phi\circ\pi)_*E'\cong \dR\phi_*\dR\pi_*E'\cong \dR\phi_*E=0
\]
as required.

Now, let $E$ be a $\phi\circ\pi$-invisible sheaf on $X$.  We again have
\[
\dR\phi_*\dR\pi_*E = 0.
\]
Since $\phi$ and $\pi$ have $\le 1$-dimensional fibers, the corresponding
spectral sequence collapses at the $E_2$ page, and thus
\[
R^p\phi_*R^q\pi_*E=0
\]
for $p,q\in \{0,1\}$.  In particular, $\pi_*E$ is invisible, and we need
only show that $E$ is $\pi_*$-acyclic.  But since $R^1\pi_*E$ is supported on
the indeterminacy locus of $\pi^{-1}$, the only way it can have trivial
direct image is to be 0.

The claims about projective and injective objects follow from adjunction.
For instance if $E$ is projective in the category of
$\pi\circ\phi$-invisible sheaves, and $E'$ is any $\phi$-invisible
sheaf, then
\[
\Ext^{p+1}(\pi_*E,E')
\cong
\Ext^{p+1}(E,\pi^!E')
=
0.
\]
\end{proof}

Thus the atomic case (monoidal transformations blowing up a single point)
will be particularly useful.  Here, we can completely characterize the
invisible sheaves.

\begin{lem}
Suppose that $\pi:X\to Y$ is a monoidal transformation, blowing up the
point $p\in Y$ to the exceptional line $e\subset X$.  Then the functor
\[
E\mapsto \Hom(\sO_e(-1),E)
\]
establishes an equivalence between the category of $\pi$-invisible
sheaves and the category of finite-dimensional vector spaces over $\bar{k}$.
\end{lem}

\begin{proof}
  Let $D$ be an effective divisor on $Y$ that has multiplicity $1$ at $p$.
  Then we obtain a short exact sequence
\[
0\to \sO_e(-1)\to \pi^*\sO_D\to \pi^{*!}\sO_D\to 0.
\]
(Use the minimal resolution of $\sO_D$ over the local ring at $p$ to
compute $\pi^{*!}\sO_D$.)  We thus find
\[
R^p\Hom(\sO_e(-1),E)\cong R^{p+1}\Hom(\pi^{*!}\sO_D,E)
\cong R^{p+1}\Hom(\sO_{\tilde{D}},E),
\]
where $\tilde{D}$ is the strict transform of $D$.  Since $\tilde{D}$ is
transverse to the exceptional line, it is transverse to the support of $E$,
and thus the only nonvanishing $\Ext$ group is in degree 1.  We thus
conclude that for any invisible sheaf,
\[
R^p\Hom(\sO_e(-1),E)=0
\]
for $p>0$; in other words, $\sO_e(-1)$ is a projective object in the
category of invisible sheaves.

Since $\sO_e(-1)$ is coherent, with endomorphism ring $k$, it remains only
to show that it generates the category, since then the desired equivalence
follows by Morita theory.  In other words, we need to show that
if $\Hom(\sO_e(-1),E)=0$, then $E=0$.  Since $\tilde{D}$ is transverse to
the exceptional locus, we have a short exact sequence
\[
0\to E\to E(\tilde{D})\to T\to 0
\]
for some $0$-dimensional sheaf $T$.  If $T=0$, then $E$ has support
disjoint from $\tilde{D}$; if nontrivial, it would have $0$-dimensional
support, and thus nontrivial direct image.  Thus $T\ne 0$, and we conclude
that
\[
\chi(E(\tilde{D}))=\chi(T)>0,
\]
so that $E(\tilde{D})$ has global sections, and its direct
image is thus a nontrivial $0$-dimensional sheaf.  It follows from the next
lemma that
\[
\Hom(\sO_e(-1),E)
\cong
\Hom(\sO_e,E(\tilde{D}))\ne 0,
\]
as required.
\end{proof}

\begin{lem}\label{lem:char_of_hd2_and_acyclic}
  Let $\pi:X\to Y$ be the blowup in the closed point $p$ of $Y$, with
  exceptional line $e$, and let $M$ be a coherent sheaf on $X$ of
  homological dimension $\le 1$.  If $M$ is not $\pi_*$-acyclic, then
  $\Hom(M,\sO_e(-2))\ne 0$, while if $\pi_*M$ has homological dimension
  $2$, then $\Hom(\sO_e,M)\ne 0$.
\end{lem}

\begin{proof}
  Since $\pi$ is an isomorphism away from $e$ and $p$, the only way that
  $\pi_*M$ can have homological dimension 2 is if it has a subsheaf
  isomorphic to $\sO_p$.  We can compute
\[
\dR\Hom(\sO_p,\dR\pi_*M)
\cong
\dR\Hom(\dL\pi^*\sO_p,M)
\]
so in particular
\[
\Hom(\sO_p,\pi_*M) \cong \Hom(\pi^*\sO_p,M).
\]
A simple calculation in the local ring at $p$ gives
\[
L_1\pi^*\sO_p\cong \sO_e(-1),
\quad\pi^*\sO_p\cong \sO_e,
\]
giving
\[
\Hom(\sO_p,\pi_*M) \cong \Hom(\sO_e,M),
\]
implying the desired result.

Similarly, if $M$ is not $\pi_*$-acyclic, then $R^1\pi_*M$ is a nontrivial
sheaf supported at $p$, so has a morphism to $\sO_p$.  We compute
\[
\dR\Hom(\dR\pi_*M,\sO_p)
\cong
\dR\Hom(M,\dL\pi^!\sO_p)
\]
which in degree $-1$ gives
\[
\Hom(R^1\pi_*M,\sO_p)
\cong
\Hom(M,L_1\pi^!\sO_p)
\cong
\Hom(M,\sO_e(-2)).
\]
\end{proof}

\begin{cor}
  Let $\pi:X\to Y$ be the blowup in a single point $p\in Y$, with
  exceptional line $e$, and let $M$ be a sheaf of homological dimension $
  \le 1$ on $X$.  If $\Hom(\sO_e(-1),M)=\Hom(M,\sO_e(-1))=0$, then $M\cong
  \pi^{*!}\pi_*M$.
\end{cor}

\begin{proof}
  We first observe that $M$ is $\pi_*$-acyclic with $\pi_*M$ of homological
  dimension $\le 1$.  Indeed, we would otherwise have a nonzero map
  $\sO_e\to M$ or $M\to \sO_e(-2)$, which we could compose with any nonzero
  map $\sO_e(-1)\to \sO_e$ or $\sO_e(-2)\to \sO_e(-1)$.  The second
  composition is necessarily nonzero by injectivity of the second map; the
  first composition could only be zero if the original map had
  $0$-dimensional image.

In particular, we find that we have natural maps
\[
\pi^*\pi_*M\to M\to \pi^!\pi_*M.
\]
The first map is surjective, since its cokernel is invisible, and the
hypotheses imply that $M$ has no nonzero maps to invisible sheaves.
Similarly, the second map is injective since it has invisible kernel.  In
other words, $M$ is the image of the natural map $\pi^*\pi_*M\to
\pi^!\pi_*M$ as required.
\end{proof}

For a more general birational morphism, the exceptional locus is still a
union of finitely many smooth rational curves, which we call the {\em
  exceptional components} of $\pi$.

\begin{lem}\label{lem:no_image_maps_to_except}
  Let $\pi:X\to Y$ be a birational morphism of algebraic surfaces, and let
  $M$ be a nonzero coherent sheaf on $X$ such that $\pi_*M=0$.  Then there
  is some exceptional component $f$ such that $\Hom(M,\sO_f(-1))\ne 0$.
\end{lem}

\begin{proof}
  If $\pi$ is the blowup in a point $p\in Y$, then either $R^1\pi_*M\ne 0$,
  in which case it maps to $\sO_e(-2)$ and thus $\sO_e(-1)$, or
  $R^1\pi_*M=0$, in which case it is invisible, and thus is isomorphic to
  $\sO_e(-1)^r$.

Otherwise, we can factor $\pi=\pi_1\circ \pi_2$ such that $\pi_1$ blows up
a single point.  If $\Hom(M,\sO_e(-1))\ne 0$, where $e$ is the exceptional
locus of $\pi_1$, then we are done.  Otherwise, $M$ is $\pi_{1*}$-acyclic,
and $\pi_{2*}(\pi_{1*}M)\cong \pi_*M=0$.  Moreover, we have a short exact
sequence of the form
\[
0\to \sO_e(-1)^r\to \pi^*_1\pi_{1*}M\to M\to 0.
\]
Since $\Hom(\sO_e(-1),\sO_f(-1))=0$ for $f\ne e$, it will suffice to show
that
\[
\Hom(\pi^*_1\pi_{1*}M,\sO_f(-1))\ne 0
\]
for some exceptional component $f\ne e$. But
\[
\Hom(\pi^*_1\pi_{1*}M,\sO_f(-1))
\cong
\Hom(\pi_{1*}M,\pi_{1*}\sO_f(-1))
\cong
\Hom(\pi_{1*}M,\sO_{\pi_1(f)}(-1))
\]
Every exceptional component of $\pi_2$ is the image of some $f\ne e$, and
thus, by induction, at least one of the latter groups is nonzero.
\end{proof}

\begin{cor}
  Any invisible sheaf has a quotient of the form $\sO_f(-1)$ for some
  exceptional component $f$.
\end{cor}

\begin{proof}
  By the lemma, any invisible sheaf has a nonzero morphism to some
  $\sO_f(-1)$.  The image of such a morphism has the form the form
  $\sO_f(-d)$ for some $d>1$.  Since it is a quotient of an invisible,
  thus $\pi_*$-acyclic, sheaf, it must be $\pi_*$-acyclic.  But then
\[
H^1(\sO_f(-d)) = H^1(\pi_*\sO_f(-d)) = 0,
\]
since it has $0$-dimensional direct image.  It follows that $d=1$.
\end{proof}

\begin{cor}
The category of invisible sheaves is Artinian.
\end{cor}

\begin{proof}
Let $E_1\supsetneq E_2\supsetneq E_3\supsetneq\cdots$ be a descending
chain of invisible sheaves.  Each quotient $E_i/E_{i+1}$ is
a nontrivial invisible sheaf, so surjects on some $\sO_f(-1)$.
In particular, relative to any very ample bundle on $X$, $\deg(c_1(E_i))$
is a strictly decreasing sequence of nonnegative integers.
\end{proof}

\begin{cor}
  Any invisible sheaf admits a filtration in which the successive
  quotients are all of the form $\sO_f(-1)$ with $f$ an exceptional
  component.
\end{cor}

\begin{proof}
Any nontrivial invisible sheaf has a nontrivial map to some $\sO_f(-1)$,
and the kernel is invisible.  Iterating gives a descending chain
of invisible sheaves, which must eventually reach 0.
\end{proof}

This gives us an alternate characterization of minimal lifts.

\begin{thm}\label{thm:minimal_if_no_exceptional}
  Let $\pi:X\to Y$ be a birational morphism of algebraic surfaces, and
  suppose that $M$ is a sheaf on $X$ of homological dimension $\le 1$.  If
  $\Hom(M,\sO_f(-1))=0$ for all exceptional components $f$, then $M$ is
  $\pi_*$-acyclic and $\pi$-globally generated.  If $\Hom(\sO_f(-1),M)=0$
  for all exceptional components $f$, then $\pi_*M$ has homological
  dimension $\le 1$.  Moreover, $M$ is a minimal lift iff
  $\Hom(\sO_f(-1),M)=\Hom(M,\sO_f(-1))=0$ for all exceptional components $f$.
\end{thm}

\begin{proof}
  Suppose that $\Hom(M,\sO_f(-1))=0$ for all exceptional components, and
  consider the natural map
\[
\pi^*\pi_*M\to M.
\]
The direct image of this map is an isomorphism, and thus the usual spectral
sequence shows that the cokernel has trivial direct image.  It follows that
if the cokernel is nonzero, then it maps to $\sO_f(-1)$ for some $f$,
giving a corresponding map from $M$.  We thus conclude that $M$ is
$\pi$-globally generated, and thus $\pi_*$-acyclic (as a quotient of the
$\pi_*$-acyclic sheaf $\pi^*\pi_*M$).

Next, note that if $\Hom(\sO_f(-1),M)=0$ for all exceptional components,
then $\Hom(E,M)=0$ for all invisible sheaves $E$.  So it will suffice to
show that if $\pi_*M$ has homological dimension 2, then $\Hom(E,M)\ne 0$
for some invisible $E$.  Factor $\pi=\pi_1\circ \pi_2$ with $\pi_1$
monoidal.  If $\pi_{1*}M$ has homological dimension 2, then
$\Hom(\sO_e,M)\ne 0$, and thus $\Hom(\sO_e(-1),M)\ne 0$.  Otherwise, by
induction, we have $\Hom(E,\pi_{1*}M)\ne 0$ for some $\pi_2$-invisible $E$,
and thus $\Hom(\pi_1^*E,M)\ne 0$.  Since $\pi_1^*E$ is $\pi$-invisible, the
claim follows.

Now, suppose that both conditions hold.  Then in the composition
\[
\pi^*\pi_*M\to M\to \pi^!\pi_*M
\]
we have already shown the first map to be surjective, and the second map
has invisible kernel, so must be injective.  That $M\cong
\pi^{*!}\pi_*M$ follows.
\end{proof}

Given an exceptional component $f$, let $f^\vee$ denote the linear
combination of exceptional components such that $f^\vee\cdot
g=-\delta_{fg}$ for exceptional components $g$.  (Note that although $X$ is
not assumed projective, it still has a well-behaved intersection theory, at
least where exceptional divisors are concerned, since the exceptional
divisors themselves {\em are} projective.)  Equivalently,
\[
\sO_g(f^\vee)\cong \sO_g(-\delta_{fg}).
\]
By induction on a factorization of $\pi$ into monoidal transformations, we
find that $f^\vee$ is effective.  Define a sheaf $P_f$ by the exact sequence
\[
0\to P_f\to \pi^*\pi_*\sO_X(-f^\vee)\to \sO_X(-f^\vee).
\]

\begin{lem}
The sheaves $P_f$ are projective objects in the category of
$\pi$-invisible sheaves.  More precisely, $P_f$ is the projective cover
of $\sO_f(-1)$.
\end{lem}

\begin{proof}
Since
\[
\Hom(\sO_X(-f^\vee),\sO_g(-1)) \cong H^0(\sO_g(-1-\delta_{fg}))=0,
\]
it follows that $\sO_X(-f^\vee)$ is $\pi_*$-acyclic and $\pi$-globally
generated, so that the exact sequence defining $P_f$ extends to a short
exact sequence.  In addition, we have
\[
\pi_*\sO_X(-f^\vee)\subset \pi_*\sO_X\cong \sO_Y,
\]
and subsheaves of locally free sheaves have homological dimension $\le 1$.

Thus $P_f$ is an invisible sheaf, and we have
\[
R^p\Hom(P_f,\sO_g(-1))
\cong
R^{p+1}\Hom(\sO_X(-f^\vee),\sO_g(-1))
\cong
H^{p+1}(\sO_g(-1-\delta_{fg}))
\]
In particular, $R^p\Hom(P_f,\sO_g(-1))=0$ for $p>0$; since every
invisible sheaf is an extension of sheaves $\sO_g(-1)$, it follows that
$P_f$ is projective.  Moreover,
\[
\dim H^1(\sO_g(-1-\delta_{fg})) = \delta_{fg},
\]
and thus $P_f$ has a unique map to $\sO_f(-1)$, and maps to no other
component; it is thus the projective cover of $\sO_f(-1)$.
\end{proof}

\begin{rem}
Taking the direct image of the short exact sequence
\[
0\to \sO_X(-f^\vee)\to \sO_X\to \sO_{f^\vee}\to 0
\]
gives
\[
0\to \pi_*\sO_X(-f^\vee)\to \sO_Y\to \pi_*\sO_{f^\vee}\to 0,
\]
so that $\pi_*\sO_X(-f^\vee)$ is the ideal sheaf of the $0$-dimensional
subscheme $\pi(f^\vee)\subset Y$, and
\[
P_f\cong L_1\pi^*\pi_*\sO_{f^\vee}\cong L_1\pi^*\sO_{\pi(f^\vee)}.
\]
\end{rem}

Since the category of $\pi$-invisible sheaves has finitely many
irreducible objects, and every irreducible object has a projective cover,
it is isomorphic to a module category.  More precisely, we have the
following, by a straightforward application of Morita theory.

\begin{thm}
  There is a natural equivalence between the category of $\pi$-invisible
  sheaves and the category of finitely generated modules over the
  finite-dimensional $k$-algebra $\End(\bigoplus_f P_f)$, given by
\[
E\mapsto \Hom(\bigoplus_f P_f,E).
\]
Moreover, the Yoneda $\Ext$ groups in this category agree with the $\Ext$ groups
in the category of sheaves.
\end{thm}

\begin{proof}
  The main claim follows from the fact that $\bigoplus_f P_f$ is a
  progenerator of the category (it is coherent, projective, and generates
  the category).  

  By using a projective resolution to compute the Yoneda $\Ext$ groups, we
  see that the remaining claim reduces to showing that the group
  $\Ext^p_{\coh(X)}(P,E)$ vanishes for $p>0$, $P$ projective, and $E$
  invisible.  But this is true for any projective object of the form
  $P_f$, and thus for direct summands of direct sums of such sheaves.
\end{proof}

Since the canonical map $\pi^*\sO_Y\to \pi^!\sO_Y$ is injective (it has
invisible kernel, but $\sO_X$ is torsion-free), it gives rise to a
canonical effective divisor $e_\pi$ representing $K_X-\pi^*K_Y$ (the
divisor class of $\pi^!\sO_Y$), with $e_\pi$ supported on the exceptional
locus.

\begin{cor}\label{cor:invis_on_e_pi}
  Every invisible sheaf is scheme-theoretically supported on a subscheme
  of $e_\pi$.
\end{cor}

\begin{proof}
Every invisible sheaf is a quotient of a sum of sheaves $P_f$, so it
suffices to consider those.  Now, $P_f$ is contained in the kernel $E$ of the
natural map
\[
\pi^*\pi_*\sO_X(-f^\vee)\to \pi^!\pi_*\sO_X(-f^\vee).
\]
Starting with a locally free resolution
\[
0\to V\to W\to \pi_*\sO_X(-f^\vee)\to 0,
\]
we can use the snake lemma to obtain an exact sequence
\[
0\to E\to \pi^!V/\pi^*V\to \pi^!W/\pi^*W.
\]
Since
\[
\pi^!V/\pi^*V
\cong 
(\pi^!\sO_Y/\pi^*\sO_Y)\otimes \pi^*V
\cong
\sO_{e_\pi}(e_\pi)\otimes
\pi^*V,
\]
the result follows.
\end{proof}

\begin{rem}
  The invisible sheaf $\pi^!\sO_Y/\pi^*\sO_Y\cong \sO_{e_\pi}(e_\pi)$ shows
  that this bound on the support is tight.
\end{rem}

The injective objects in the category of invisible sheaves can also be
identified, most easily using the following duality, a special case of
Cohen-Macaulay duality for sheaves of homological dimension 1 and
codimension 1 support.

\begin{lem}
If $E$ is an invisible sheaf, then so is $\sExt^1(E,\omega_X)$, and
$\sExt^1(\sExt^1(E,\omega_X),\omega_X)\cong E$.
\end{lem}

\begin{proof}
  Since $E$ has $1$-dimensional support and homological dimension 1,
\[
\sHom(E,\omega_X)=\sExt^2(E,\omega_X)=0,
\]
and thus $\dR\sHom(E,\omega_X)$ is concentrated in degree 1.  Since
$\omega_X=\pi^!\omega_Y$, we have
\[
\dR\pi_*\dR\sHom(E,\omega_X)
\cong
\dR\sHom(\dR\pi_*E,\omega_Y)
= 0,
\]
and thus $\dR\pi_*\sExt^1(E,\omega_X)=0$ as required.

That this operation is a duality follows by considering how it acts on a
locally free resolution $0\to V\to W\to E\to 0$.
\end{proof}

\begin{rem}
This implies that the $k$-algebra associated to $\pi$ above is isomorphic
to its opposite.
\end{rem}

In particular, we can define injective objects
\[
I_f:= \sExt^1(P_f,\omega_X)\cong \sExt^1(\pi^*\pi_*\sO_X(-f^\vee),\omega_X)
\]
and find that $I_f$ is the injective hull of $\sO_f(-1)$.  Note that the
sheaves $\sO_f(-1)$ are self-dual: considerations of support show that the
dual of $\sO_f(-1)$ has the form $\sO_f(-d)$, and $\sO_f(-1)$ is the only
invisible possibility.  Another self-dual invisible sheaf is the sheaf
$\pi^!\sO_Y/\pi^*\sO_Y$: the dual is most naturally expressed as
$\pi^!\omega_Y/\pi^*\omega_Y$, but twisting by the pullback of a line
bundle has no effect on an invisible sheaf.

\section{Minimal lifts and anticanonical curves}
\label{sec:*!_when_anticanon}

In the case of direct images of line bundles on smooth curves, the purpose
of lifting is to make the lifted sheaf disjoint from the anticanonical
curve.  This is of course not possible for sheaves of positive rank, but we
can still hope to make the restriction to the anticanonical curve simpler.
Thus we would like to understand how the restriction of the minimal lift is
related to the original restriction.  This is not quite the correct
question, as it turns out: it turns out to be easier to twist by a line
bundle and consider
\[
M\otimes \sO_{C_\alpha}(C_\alpha).
\]
The standard resolution
\[
0\to \sO_X\to \sO_X(C_\alpha)\to \sO_{C_\alpha}(C_\alpha)\to 0
\]
shows that
\[
\Tor_p(M,\sO_{C_\alpha}(C_\alpha))
\cong
\sExt^{1-p}_Y(\sO_{C_\alpha},M).
\]
Since the minimal lift of $\sO_{C_\alpha}$ is the structure sheaf of the
new anticanonical curve, our question generalizes to understanding how
$\sExt^*_Y(\pi^{*!}M,\pi^{*!}N)$ and $\sExt^*_Y(M,N)$ are related.

\begin{lem}\label{lem:lift_exts}
Let $\pi:X\to Y$ be a birational morphism of smooth surfaces, and suppose
$M$, $N$ are sheaves on $Y$ with homological dimension $\le 1$.  Then
there is an isomorphism
\[
\sHom_Y(\pi^{*!}M,\pi^{*!}N)\cong \sHom_Y(M,N)
\]
and an exact sequence
\[
0\to \sExt^1_Y(\pi^{*!}M,\pi^{*!}N)
 \to \sExt^1_Y(M,N)
 \to \sHom_Y(E_1,E_2)
 \to \sExt^2_Y(\pi^{*!}M,\pi^{*!}N)
 \to 0,
\]
where $E_1$ and $E_2$ are the $\pi$-invisible sheaves fitting into the
short exact sequences
\begin{align}
&0\to E_1\to \pi^*M\to \pi^{*!}M\to 0\\
&0\to \pi^{*!}N\to \pi^!N\to E_2\to 0.
\end{align}
\end{lem}

\begin{proof}
Using the $E_1$ exact sequence, we obtain the long exact sequence
\[
\cdots\to
R^p\sHom_Y(\pi^{*!}M,\pi^{*!}N)
\to
R^p\sHom_Y(\pi^*M,\pi^{*!}N)
\to
R^p\sHom_Y(E_1,\pi^{*!}N)
\to
\cdots
\]
Moreover, $R^p\sHom_Y(\pi^*M,\pi^{*!}N)\cong R^p\sHom_Y(M,N)$ by
adjunction.

Similarly, since $R^p\sHom_Y(E_1,\pi^!N)=0$, the $E_2$ exact sequence gives
isomorphisms
\[
R^p\sHom_Y(E_1,E_2)
\cong
R^{p+1}\sHom_Y(E_1,\pi^{*!}N).
\]

Combining these gives a long exact sequence
\[
\cdots\to
R^p\sHom_Y(\pi^{*!}M,\pi^{*!}N)
\to
R^p\sHom_Y(M,N)
\to
R^{p-1}\sHom_Y(E_1,E_2)
\to
\cdots,
\]
from which the claim follows immediately.  Note that $\sExt^2_Y(M,N)=0$
since $M$ has homological dimension $\le 1$ by assumption.
\end{proof}

\begin{rem}
Similarly, we have an isomorphism
\[
\Hom(\pi^{*!}M,\pi^{*!}N)\cong \Hom(M,N)
\]
and a long exact sequence
\[
\cdots
\to
\Ext^{p-2}(E_1,E_2)
\to
\Ext^p(\pi^{*!}M,\pi^{*!}N)
\to
\Ext^p(M,N)
\to
\Ext^{p-1}(E_1,E_2)
\to
\cdots
\]
\end{rem}

\begin{cor}
Let $\pi:X\to Y$ be a birational morphism of smooth surfaces, and suppose
$M$ is a sheaf on $Y$ with homological dimension $\le 1$.
Then $\End(\pi^{*!}M)\cong \End(M)$.  In particular, the minimal lift of a
simple sheaf is always simple.
\end{cor}

\begin{cor}
  Let $\pi:X\to Y$ be a Poisson birational morphism of Poisson surfaces,
  and suppose $M$ is a sheaf on $Y$ of homological dimension $\le 1$.  Then
\[
\sHom_Y(\sO_{C_{\alpha_X}},\pi^{*!}M)\cong \sHom_Y(\sO_{C_{\alpha_Y}},M)
\]
and there is an exact sequence
\begin{align}
0\to \sExt^1_Y(\sO_{C_{\alpha_X}},\pi^{*!}M)
 &\to \sExt^1_Y(\sO_{C_{\alpha_Y}},M)\notag\\
 &\to \sHom_Y(\omega_X/\pi^*\omega_Y,\pi^!M/\pi^{*!}M)
 \to \sExt^2_Y(\sO_{C_{\alpha_X}},\pi^{*!}M)
 \to 0.\notag
\end{align}
\end{cor}

For a monoidal transformation, the invisible sheaves associated to a
minimal lift are straightforward to compute.

\begin{prop}\label{prop:four_term_for_monoidal}
  Let $\pi:X\to Y$ be a monoidal transformation, and let $M$ be a sheaf on
  $Y$ of homological dimension $\le 1$.  Then there is an exact sequence
\[
0\to E_1\to \pi^*M\to \pi^!M\to E_2\to 0,
\]
where we have non-canonical isomorphisms
\begin{align}
E_1&\cong \Ext^1(\sO_p,M)\otimes_k \sO_e(-1),\notag\\
E_2&\cong \Hom_k(\Hom(M,\sO_p),\sO_e(-1)).\notag
\end{align}
\end{prop}

\begin{proof}
  There certainly is an exact sequence of the given form, since the kernel
  and cokernel are invisible, thus powers of $\sO_e(-1)$.  From the exact
sequence
\[
0\to E_1\to \pi^*M\to \pi^{*!}M\to 0,
\]
we find
\begin{align}
\Hom(\sO_e(-1),E_1)
\cong
\Hom(\sO_e(-1),\pi^*M)
&\cong
\Hom(\sO_e(-2),\pi^!M)\notag\\
&\cong
\Ext^1(R^1\pi_*\sO_e(-2),M)
\cong
\Ext^1(\sO_p,M).\notag
\end{align}
Similarly,
\[
\Hom(E_2,\sO_e(-1))
\cong
\Hom(\pi^!M,\sO_e(-1))
\cong
\Hom(\pi^*M,\sO_e)
\cong
\Hom(M,\sO_p).
\]
\end{proof}


This has a nice consequence for twists of minimal lifts, which we will
generalize below.

\begin{cor}\label{cor:monoidal_pseudotwist}
  Let $\pi:X\to Y$ be a monoidal transformation, and $M$ a sheaf on $Y$ of
  homological dimension $\le 1$.  Then the sheaves $\pi^{*!}M(\pm e)$ are
  $\pi_*$-acyclic, with direct image of homological dimension $\le 1$.
  Moreover, we have short exact sequences
\[
0\to M\to \pi_*(\pi^{*!}M(e))\to \Ext^1(\sO_p,M)\otimes_k \sO_p \to 0,
\]
and
\[
0\to \pi_*(\pi^{*!}M(-e))\to M\to 
\Hom_k(\Hom(M,\sO_p),\sO_p)\to 0.
\]
In other words, $\pi_*(\pi^{*!}M(e))$ is the universal
extension of $\sO_p$ by $M$, and $\pi_*(\pi^{*!}M(-e))$ is
the kernel of the universal homomorphism from $M$ to $\sO_p$.
\end{cor}

\begin{proof}
Take the short exact sequences
\begin{align}
0\to E_1\to \pi^*M\to \pi^{*!}M\to 0\notag\\
0\to \pi^{*!}M\to \pi^!M\to E_2\to 0\notag,
\end{align}
twist the first by $\sO_X(e)$ and the second by $\sO_X(-e)$, then
take (higher) direct images.  Note that after twisting, the middle terms in
the short exact sequence are still $\pi_*$-acyclic (since they are $\pi^!M$
and $\pi^*M$ respectively).

It remains only to show that $\pi^{*!}M(-e)$ is
$\pi_*$-acyclic and $\pi^{*!}M(e)$ has direct image of
homological dimension 1; the other two claims follow from the facts that
they are a subsheaf or quotient of a well-behaved sheaf.  But
\[
\Hom(\pi^{*!}M(-e),\sO_e(-1))
\cong
\Hom(\pi^{*!}M,\sO_e(-2))
\subset
\Hom(\pi^{*!}M,\sO_e(-1))
=0,
\]
and similarly for
\[
\Hom(\sO_e(-1),\pi^{*!}M(e)).
\]
\end{proof}

Note that if we twist by more than just $\pm e$, then we cannot expect to
obtain a well-behaved direct image.  Indeed, twisting $\sO_X\cong
\pi^{*!}\sO_Y$ by $\sO_X(2e)$ gives (naturally enough) $\sO_X(2e)$,
which fits in an exact sequence
\[
0\to \sO_X(e)\to \sO_X(2e)\to \sO_e(-2)\to 0
\]
so that $R^1\pi_*\sO_e(-2)\cong \sO_p$.  So we cannot iterate the twisting
itself, but must instead iterate the direct-image-of-twist-of-minimal-lift
operation (which we call a {\em pseudo-twist}).  Again, we must be careful
to note that the two pseudo-twists are not inverses to each other.  Indeed,
if $M$ is not torsion, then the pseudo-twist by $\sO_X(e)$ and the
pseudo-twist by $\sO_X(-e)$ change $[M]\in K_0$ by different multiples of
$[\sO_p]$, since
\begin{align}
\dim\Hom(M,\sO_p)
-
\dim\Ext^1(\sO_p,M)
&{}=
\dim\Ext^2(\sO_p,M)
-
\dim\Ext^1(\sO_p,M)\notag\\
&{}=
\dim\Ext^2(\sO_p,M)-\dim\Ext^1(\sO_p,M)+\dim\Hom(\sO_p,M)\notag\\
&{}=
\rank(M).\notag
\end{align}

Even if $M$ is torsion, they need not be inverses, e.g., if $M$ is an
invertible sheaf on a curve singular at $p$.  However, here the situation
is much nicer.  First, if $\pi^{*!}M$ is transverse to $e$, then no matter
how far we twist $\pi^{*!}M$ in either direction, the result will still be
a minimal lift, and thus in the transverse case, the pseudo-twists act as a
group.  More generally, the subsheaf $\Tor_1(\pi^{*!}M,\sO_e(-1))\subset
\pi^{*!}M$ will be a vector bundle on $e\cong \P^1$ having no global
sections.  If we repeatedly pseudo-twist by $\sO_X(-e)$, then each step
either strictly decreases the rank of this vector bundle or strictly
increases its degree.  Thus after a finite number of downwards
pseudo-twists, we obtain a sheaf such that $\pi^{*!}M$ is transverse to the
exceptional curve.  We could also obtain a sheaf with transverse minimal
lift via finitely many upwards pseudo-twists, or by some combination of the
two; the resulting sheaves need not be in the same orbit under the group of
pseudo-twists.

In generalizing pseudo-twists beyond monoidal transformations, we encounter
the problem of determining which twists of a minimal lift are guaranteed to
have well-behaved direct image (i.e., $\pi_*$-acyclic with direct image of
homological dimension $\le 1$).  This appears to be a tricky problem in
general (especially in the presence of exceptional components of
self-intersection $<-2$), but there is one sufficiently large special
case.

\begin{thm}\label{thm:elem_pseudo_twists}
Let $\pi:X\to Y$ be a birational morphism of algebraic surfaces, and
suppose $D$ is a divisor on $X$ which is supported on the exceptional locus
and satisfies $D^2=-1$.  If $M$ is any sheaf of homological dimension
$\le 1$ on $Y$, then $\pi^{*!}M(D)$ is $\pi_*$-acyclic,
with direct image of homological dimension $\le 1$.
\end{thm}

\begin{proof}
Since the intersection pairing on the exceptional locus is integral and
negative definite, there are at most $2n$ exceptional divisors of
self-intersection $-1$, where $n$ is the number of monoidal transformations
in a factorization of $\pi$.  We can moreover construct those divisors
easily enough.  Given a factorization of $\pi$, let $e_i$ be the total
transform of the $i$-th exceptional curve through the remaining $n-i$
blowups; then $e_i^2=-1$, so $\pm e_i$ satisfy the hypotheses.

In other words, we need to show that $\pi^{*!}M(\pm e_i)$ have well-behaved
direct images.  Moreover, we may assume $i=n$, since twisting by $\sO_X(\pm
e_i)$ commutes with blowing down the last $n-i$ points.  Factor
$\pi=\pi_n\circ \pi_{[n]}$, where $\pi_n$ is the monoidal transformation
blowing down the $-1$-curve $e_n$.  Then $\pi^{*!}M(e_n)$ is a quotient of
$\pi_n^!\pi_{[n]}^*M$, so is $\pi_*$-acyclic, while $\pi^{*!}M(-e_n)$ is a
subsheaf of $\pi_n^*\pi_{[n]}^!M$, so has direct image of homological
dimension $\le 1$.

It remains to show that $\pi^{*!}M(-e_n)$ is $\pi_*$-acyclic and
$\pi^{*!}M(e_n)$ has direct image of homological dimension 1.  The key idea
is that not only is $e_n$ effective, but up to twisting by pullbacks of
line bundles, so is $-e_n$.  More precisely, on $\pi_n(X)$, we may choose a
divisor $D$ which is transverse to $\pi^{*!}_{[n]}M$ and meets the
exceptional locus only in $\pi_n(e_n)$ (this certainly exists in a formal
neighborhood of the exceptional locus).  Then the strict transform
$\pi^*D-e_n$ of $D$ is transverse to $\pi^{*!}M$, so that
\[
\pi^{*!}M(e_n)\subset\pi^{*!}(M(D)),
\]
implying that $\pi_*(\pi^{*!}M(e_n))$ has homological
dimension $\le 1$.  Similarly, $\pi^{*!}M(-e_n)$ is
a quotient of $\pi^{*!}(M(-D))$, so is $\pi_*$-acyclic.
\end{proof}

Note that changing the factorization of $\pi$ only permutes the divisors
$e_i$, since they are characterized as the unique divisors of
self-intersection $-1$ which are nonnegative linear combinations of
exceptional components.  Each divisor $e_i$ is the total transform of a
$-1$-curve, the strict transform of which is an exceptional component
$f_i$; this correspondence is canonical, so that we can more naturally
label the orthonormal divisors by exceptional components.  I.e., for any
exceptional component $f$, $e_f$ is the unique effective exceptional
divisor such that $e_f^2=-1$ and $e_f\cdot f=-1$.

\begin{cor}\label{cor:pseudo_twist_is_isospectral}
Let $\pi:X\to Y$ be a birational morphism of algebraic surfaces, and let
$f$ be an exceptional component of $\pi$.  If $M$ is any sheaf of
homological dimension $\le 1$ on $Y$, then there are short exact sequences
of the form
\[
0\to M\to \pi_*(\pi^{*!}M(e_f))\to \sO_{\pi(f)}^{r_1}\to 0
\]
and
\[
0\to \pi_*(\pi^{*!}M(-e_f))\to M\to 
\sO_{\pi(f)}^{r_2}\to 0.
\]
The dimensions $r_1$, $r_2$ are given by
\[
r_1 = c_1(\pi^{*!}M)\cdot e_f,\qquad
r_2 = c_1(\pi^{*!}M)\cdot e_f+\rank(M).
\]
\end{cor}

\begin{proof}
The first short exact sequence comes from the direct image of the sequence
\[
0\to 
E_1(e_f)
\to
\pi^*M(e_f)
\to
\pi^{*!}M(e_f)
\to
0,
\]
and thus $r_1 = -\chi(E_1(e_f))$.  Since $\chi(E_1)=0$ and
$\rank(E_1)=0$, Hirzebruch-Riemann-Roch gives $r_1 = -c_1(E_1)\cdot e_f$.
Since
\[
c_1(E_1) = \pi^*c_1(M)-c_1(\pi^{*!}M),
\]
the formula for $r_1$ follows.

Similarly, the second short exact sequence is the direct image of
\[
0\to \pi^{*!}M(-e_f)\to \pi^!M(-e_f)
 \to E_2(-e_f)\to 0,
\]
so that $r_2 = \chi(E_2(-e_f))=-c_1(E_2)\cdot e_f$.
We find
\[
c_1(E_2) = c_1(\pi^!M)-c_1(\pi^{*!}M)
\]
and
\[
c_1(\pi^!M) = c_1(\pi^*M)+\rank(M)e_\pi.
\]
Since $e_\pi=\sum_f e_f$, we find $e_\pi\cdot e_f=-1$, and the formula
follows.
\end{proof}

Note that just as in the monoidal case, the upwards and downwards
pseudo-twists do not form a group unless $M$ is $1$-dimensional and
$\pi^{*!}M$ is transverse to the exceptional locus.  In addition, if $M$ is
$1$-dimensional and $\pi^{*!}M$ is not transverse to $e_\pi$, then we need
only finitely many pseudo-twists to make it transverse.  This has a
particularly nice consequence in the Poisson case.

\begin{lem}\label{lem:resolving_pseudotwist}
  Let $(Y,\alpha)$ be a Poisson surface with anticanonical curve $C_\alpha$.
  Let $M$ be a pure $1$-dimensional sheaf on $Y$ (i.e., $M$ has
  $1$-dimensional support, and no subsheaf has $0$-dimensional support)
  which is transverse to $C_\alpha$.  Then there exists a Poisson birational
  morphism $\pi:X\to Y$ such that some pseudo-twist $M'$ of $M$ has minimal
  lift disjoint from the anticanonical curve on $X$.
\end{lem}

\begin{proof}
  If $\supp(M)$ is disjoint from $C_\alpha$, then there is nothing to prove.
  Otherwise, let $\pi_1:Y_1\to Y$ be the (Poisson) blow up in a point of
  intersection, and let $M_1$ be a pseudo-twist of $M$ such that
  $\pi^{*!}_1M_1$ is transverse to the exceptional locus.  Then
  $\pi^{*!}_1M_1$ is also transverse to $C_{\alpha_1}$, where $\alpha_1$ is the
  induced Poisson structure on $Y_1$.  Moreover,
\[
c_1(\pi^{*!}_1M_1)\cdot c_1(C_{\alpha_1})
=
(\pi^*c_1(M_1)-r_1e_1)(\pi^*c_1(C_\alpha)-e_1)
<
c_1(M_1)\cdot c_1(C_\alpha)
=
c_1(M)\cdot c_1(C_\alpha).
\]
(That $c_1(M_1)=c_1(M)$ is immediate from the short exact sequences for the
atomic pseudo-twists.)  Thus if we iterate this operation, then after at
most $c_1(M)\cdot c_1(C_\alpha)$ blowups, we will achieve disjointness.
\end{proof}

\begin{rem}
Since one of our motivations is to avoid using support (since this
makes little sense in a noncommutative context), we should note
that ``pure $1$-dimensional'' can be rephrased strictly in terms
of Hilbert polynomials: $M$ has linear Hilbert polynomial, and no
nonzero subsheaf has constant Hilbert polynomial.  Similarly, the
first Chern class of a pure $1$-dimensional sheaf specifies how the
linear term in the Hilbert polynomial depends on the choice of very
ample line bundle.
\end{rem}

If $\pi^{*!}M$ is disjoint from $C_{\alpha_X}$, then we find
\[
M\otimes \sO_{C_{\alpha_Y}}(C_{\alpha_Y})
\cong
\sExt^1_Y(\sO_{C_{\alpha_Y}},M)
\cong
\sHom_Y(\omega_X/\pi^*\omega_Y,\pi^!M/\pi^{*!}M).
\]
Moreover, disjointness has strong consequences for the invisible quotient
$\pi^!M/\pi^{*!}M$.  Since for any exceptional component, $f\cdot
C_{\alpha_X} = f^2+2$, any component with $f^2<-2$ is contained in
$C_{\alpha_X}$, and thus must be disjoint from $\pi^{*!}M$.

\begin{prop}
  Let $\pi:X\to Y$ be a birational morphism of algebraic surfaces.  If $M$
  is a sheaf of homological dimension $\le 1$ on $Y$ such that
  $\supp(\pi^{*!}M)$ does not contain any exceptional components of
  self-intersection $\le -3$, then $\pi^!M/\pi^{*!}M$ is an injective
  object of the category of $\pi$-invisible sheaves.  More precisely,
\[
\pi^!M/\pi^{*!}M
\cong
\bigoplus_f I_f^{-f\cdot c_1(\pi^{*!}M)-\rank(M)f\cdot e_\pi}.
\]
\end{prop}

\begin{proof}
Let $f$ be any exceptional component.  Since
\[
\dR\Hom(\sO_f(-1),\pi^!M)=0,
\]
we have
\begin{align}
\Ext^1(\sO_f(-1),\pi^!M/\pi^{*!}M)
&{}\cong
\Ext^2(\sO_f(-1),\pi^{*!}M)
\cong
H^1(\sExt^1_X(\sO_f(-1),\pi^{*!}M))\notag\\
&{}\cong
H^1(\sO_f(f^2+1)\otimes \pi^{*!}M)
\cong
H^0(R^1\pi_*(\sO_f(f^2+1)\otimes \pi^{*!}M))
\end{align}
If $\pi^{*!}M$ is transverse to $f$, the tensor product is $0$-dimensional,
so acyclic.  Otherwise, if $\sO_f(f^2+1)\otimes \pi^{*!}M$ (a sheaf
supported on $f\cong \P^1$) is not acyclic, then there is a nontrivial
morphism
\[
\sO_f(f^2+1)\otimes \pi^{*!}M\to \sO_f(-2),
\]
so a nontrivial morphism
\[
\pi^{*!}M\to \sO_f\otimes \pi^{*!}M\to \sO_f(-3-f^2),
\]
which is impossible unless $f^2\le -3$.

For the final claim, since $\pi^!M/\pi^{*!}M$ is injective, it is a direct
sum of injective hulls of simple sheaves, and $I_f$ occurs as a summand
with multiplicity $-f\cdot c_1(\pi^!M/\pi^{*!}M)$, since
$c_1(I_f)=c_1(P_f)=f^\vee$.  The claim then follows from
\[
c_1(\pi^!M/\pi^{*!}M)
=
\pi^*c_1(M) + \rank(M)e_\pi-c_1(\pi^{*!}M).
\]
\end{proof}

\begin{rem}
Note that the hypotheses are much weaker than disjointness; in particular,
if there are no exceptional components of self-intersection $\le -3$, then
$\pi^!M/\pi^{*!}M$ is injective as an invisible sheaf for {\em all}
sheaves $M$ of homological dimension $\le 1$.  Conversely, if $f^2\le -3$
for some $f$, then $\pi^!\sO_Y/\pi^*\sO_Y$ is not injective, since
\[
\Ext^1(\sO_f(-1),\pi^!\sO_Y/\pi^*\sO_Y)
\cong
\Ext^2(\sO_f(-1),\sO_X)
\cong
H^0(\sO_f(-1)\otimes \omega_X)^*
\cong
H^0(\sO_f(-3-f^2))^*
\ne 0.
\]
\end{rem}

In the disjoint case, we have
\[
M\otimes \sO_{C_{\alpha_Y}}
\cong
M\otimes \sO_{C_{\alpha_Y}}(C_{\alpha_Y})
\cong
\bigoplus_f \sHom_Y(\omega_X/\pi^*\omega_Y,I_f)^{-f\cdot
  c_1(\pi^{*!}M)}.
\]
Thus to complete our understanding of the disjoint case, we need to
understand the $0$-dimensional sheaves
\[
\sHom_Y(\omega_X/\pi^*\omega_Y,I_f)
\cong
\sHom_Y(P_f,\pi^!\sO_Y/\sO_X).
\]

\begin{prop}
  Let $\pi:X\to Y$ be a birational morphism of algebraic surfaces, let $f$
  be an exceptional component, and set $d=f^\vee\cdot e_\pi$.  Then there
  exists a homomorphism $\sO_{Y,\pi(f)}\to k[t]/t^d$ such that
\[
\sHom_Y(P_f,\pi^!\sO_Y/\sO_X)
\]
is the $\sO_Y$-module induced by the regular representation of $k[t]/t^d$.
\end{prop}

\begin{proof}
  Let $D$ be a smooth divisor on $X$ meeting the exceptional locus only in
  $f$, and that in a single reduced point $p$ (as before, we replace $X$ by
  a formal neighborhood of the exceptional locus if necessary).  Then
  $D=\pi^*\pi_*D-f^\vee$, so that $P_f$ fits into an exact sequence
\[
0\to P_f\to \pi^*\pi_*\sO_X(D)\to \sO_X(D)\to 0.
\]
Since $P_f$ is invisible, the canonical global section of $\sO_X(D)$
lifts to a canonical global section of $\pi^*\pi_*\sO_X(D)$, and we can
thus quotient by the two copies of $\sO_X$ to obtain
\[
0\to P_f\to \pi^*\pi_*(\sO_D(D))\to \sO_D(D)\to 0
\]
Since $\sO_D(D)$ is transverse to the exceptional locus, it
is a minimal lift, and thus we have an exact sequence
\begin{align}
0\to \sExt^1_Y(\sO_D(D),\sO_X)
 &{}\to \sExt^1_Y(\sO_D(D),\pi^!\sO_Y)\notag\\
 &{}\to \sHom_Y(P_f,\pi^!\sO_Y/\sO_X)
 \to \sExt^2_Y(\sO_D(D),\sO_X)
 \to 0\notag
\end{align}
Using the natural locally free resolution of $\sO_D(D)$
turns this into a short exact sequence
\[
0
\to 
\pi_*(\sO_D)
\to
\pi_*(\sO_D\otimes \pi^!\sO_Y)
\to
\sHom_Y(P_f,\pi^!\sO_Y/\sO_X)
\to
0;
\]
the assumptions on $D$ mean that $R^1\pi_*\sO_D=0$, and thus
$\sExt^2_Y(\sO_D(D),\sO_X)=0$.  In particular, the above
short exact sequence must agree with the direct image of the short exact
sequence
\[
0\to \sO_D\to \sO_D\otimes \pi^!\sO_Y\to J\to
0.
\]
Now, since $D$ is smooth and meets the exceptional locus in a single point,
$J$ is the structure sheaf of a jet, so has the form $k[t]/t^d$, where $t$
is a uniformizer of $D$ at $p$.  The claim then follows using the
composition
\[
\sO_{Y,p}\to \sO_{D,p}\to k[t]/t^d
\]
to compute the direct image.
\end{proof}

\begin{rem}
Note that the homomorphism to $k[t]/t^d$ need not be surjective, reflecting
the fact that $\sHom_Y(P_f,\pi^!\sO_Y/\sO_X)$ need not be the structure
sheaf of a jet.
\end{rem}

The sheaf $\pi^!\sO_Y/\sO_X$ also has some relevant structure.  If $\pi$
simply blows up a collection of distinct points on $Y$, then any choice of
ordering on those points gives a natural filtration of $\pi^!\sO_Y/\sO_X$.
Interestingly, we have the same freedom even for much more complicated
birational morphisms.

\begin{prop}\label{prop:sing_filtration}
The invisible subsheaves of $\pi^!\sO_Y/\sO_X$ form a boolean lattice.
The maximal chains in the lattice are in natural correspondence with the
orderings on the exceptional components; if $f_1$, $f_2$,\dots,$f_n$
is such an ordering, then there is a unique maximal chain
\[
0=E_0\subset E_1\subset\cdots E_n=\pi^!\sO_Y/\sO_X
\]
of invisible sheaves such that
\[
c_1(E_i/E_{i-1})=e_{f_i}.
\]
Moreover, for any exceptional component $f$,
\[
\sHom_Y(P_f,E_i/E_{i-1})\cong \sO_{\pi(f)}^{f^\vee\cdot e_{f_i}}.
\]
\end{prop}

\begin{proof}
  The subsheaves of $\pi^!\sO_Y/\sO_X$ are in natural correspondence with
  the subsheaves of $\pi^!\sO_Y$ containing $\sO_X$.  If $M$ is such a
  subsheaf, then $M/\sO_X$ is invisible iff $\pi^!\sO_Y/M$ is
  invisible.  In particular, $\pi^!\sO_Y/M$ cannot have any
  $0$-dimensional subsheaves, so that it has homological dimension $\le 1$.
  Thus $M$ is locally free, and $M\cong \sO_X(D)$ for some exceptional
  divisor $D$.  Moreover, we compute
\[
0 = \chi(M/\sO_X) = D\cdot (e_\pi-D).
\]
In particular, if we express $D$ in terms of the orthonormal basis $e_f$,
then every coefficient must be either 0 or 1.  In other words,
\[
D = \sum_{f\in S} e_f
\]
for some subset $S$ of the set of exceptional components.  Since the sheaf
$\sO_D(D)$ has Euler characteristic 0 and is contained
in the invisible sheaf $\pi^!\sO_Y/\sO_X$, it too is invisible.  We
thus conclude that the invisible subsheaves of $\pi^!\sO_Y/\sO_X$ are in
order-preserving bijection with the sets of exceptional components.
Since $c_1(\sO_D(D))=D$, the claim about Chern classes
of subquotients in maximal chains follows.

If $E$ is a subquotient in some maximal chain, with $c_1(E)=g$, then
\[
E\cong \sO_{e_g}(\sum_{h\in S} e_h)
\]
where $S$ is any set of exceptional components containing $g$.  We thus
need to show that for any such sheaf, $\sHom_Y(P_f,E)$ is
scheme-theoretically supported on $\sO_{\pi(f)}$; that it has length
  $f^\vee\cdot e_g$ follows from a Chern class computation.  Factor
  $\pi=\pi_1\circ\pi_2$ with $\pi_1$ monoidal, having exceptional curve
  $e$.  If $f\ne e$, then
\[
P_f\cong \pi^*_1P_{\pi_1(f)},
\]
so that
\[
\sHom_Y(P_f,E)
\cong
\sHom_Y(P_{\pi_1(f)},\pi_{1*}E).
\]
If $g=e$, then $E\cong \sO_g(-1)$, and thus $\pi_{1*}E=0$.  Otherwise, $E$
has the form $\pi_1^*E'$ or $\pi_1^!E'$ (depending on whether $e\in S$),
where $E'$ is a sheaf of the same form in $\pi_1(X)$.  Either way, the
claim follows by induction.

If $f=e$, then we may choose a curve $D$ as above, so that
\[
\sHom_Y(P_f,E)\cong \sExt^1_Y(\sO_D(D),E)
              \cong \pi_*(\sO_D\otimes E).
\]
Since $E$ is an invertible sheaf on its support and $D$ is transverse to
that support, $\sHom_Y(P_f,E)$ depends only on the support.  Thus if $g\ne
e$, then we may assume $E\cong \pi_1^!E'$, and thus
\[
\sHom_Y(P_f,E)\cong \sHom_Y(\pi_{1*}P_f,E').
\]
Since $\pi_{1*}P_f$ is projective, the claim follows by induction.

Finally, if $f=e=g$, then $\sHom_Y(P_f,E)$ has length 1, so is necessarily
supported on a single point.
\end{proof}

We close the section with another result on resolving sheaves via minimal
lifts.  This does not involve pseudo-twists, but as a result requires
somewhat stronger hypotheses.  Here we recall that $\Fitt_0(M)$ denotes the
$0$-th Fitting scheme of $M$, which for $M$ pure $1$-dimensional is a
canonical divisor representing $c_1(M)$.

\begin{prop}\label{prop:jets_resolvable_by_!*}
  Suppose $(Y,\alpha)$ is a Poisson surface, and $M$ is a pure
  $1$-dimensional sheaf on $Y$ such that the divisor $\Fitt_0(M)$ meets
  $C_\alpha$ in a disjoint union of jets.  Let $\pi:X\to Y$ be the minimal
  desingularization of the blowup of $Y$ in the intersection.  Then $\pi$
  is Poisson, and $\pi^{*!}M$ is disjoint from the induced anticanonical
  curve.
\end{prop}

\begin{proof}
  We can compute $\pi$ by repeatedly blowing up single points of
  intersection, so that $\pi$ is in particular Poisson.  The only issue is
  to show that the jet condition (and thus transversality) is preserved
  under the blowup $\pi:X\to Y$ in a point $p\in C_\alpha\cap\supp(M)$.
If $C_\alpha$ is smooth at $p$, the jet condition in a neighborhood of $p$ is
automatic, and this remains true after blowing up $p$.

Thus suppose $C_\alpha$ is singular at $p$.  Then we have
\[
M\otimes (\sO_X/{\cal I}_p^2)
\cong
M\otimes \sO_{C_\alpha}\otimes (\sO_X/{\cal I}_p^2),
\]
where ${\cal I}_p$ is the ideal sheaf of $p$.  Thus since $M\otimes
\sO_{C_\alpha}$ is a sum of jets, so is $M\otimes \sO_X/{\cal I}_p^2$.  But
then we can compute $\pi^{*!}M$ near $e$ using the minimal resolution of
$M$ over the local ring at $p$, and conclude that $\pi^{*!}M$ is transverse
to the exceptional curve (meeting it in the tangent vectors to the jets
through $p$), and thus is transverse to $\pi^{*!}\sO_{C_\alpha}$.

It follows that we have a short exact sequence
\[
0\to \pi_*\sExt^1_X(\pi^{*!}\sO_{C_\alpha},\pi^{*!}M)
 \to \sExt^1_Y(\sO_{C_\alpha},\pi^{*!}M)
 \to \sHom_Y(\sO_e(-1),\pi^!M/\pi^{*!}M)
 \to 0,
\]
and moreover that
\[
\sHom_Y(\sO_e(-1),\pi^!M/\pi^{*!}M)
\cong
\sHom_Y(\sExt^1_Y(\sO_{C_\alpha},\pi^{*!}M),\sO_p).
\]
In particular,
\[
\pi_*\sExt^1_X(\pi^{*!}\sO_{C_\alpha},\pi^{*!}M)
\]
is again a sum of jets, so that the same holds for
\[
\sExt^1_X(\pi^{*!}\sO_{C_\alpha},\pi^{*!}M).
\]
\end{proof}

\chapter{Poisson aspects of moduli}
\label{chap:poisson_moduli}

\section{Moduli of simple sheaves}
\label{sec:comm_poisson_moduli}

In \cite{BottacinF:1995,BottacinF:2000}, Bottacin showed that the moduli
space of stable vector bundles on a Poisson surface in characteristic 0 has
a natural Poisson structure (extending a result of Tyurin
\cite{TyurinAN:1988}, who constructed the form but did not prove the Jacobi
identity), and identified the symplectic leaves of this structure (and in
particular showed that they are algebraic).  This was extended to general
stable sheaves (in fact, simple sheaves), subject to some smoothness
assumptions, in \cite{HurtubiseJC/MarkmanE:2002b}.

Since we wish to understand how these spaces interact with birational maps,
we immediately encounter a problem: the definition of stability depends on
a choice of ample line bundle.  Since there is no canonical way to
transport ample bundles through birational maps, this adds a great deal of
complexity when trying to understand whether a given operation preserves
stability.  Moreover, since we also want to understand pseudo-twists, we
need to consider the effect of twisting by a line bundle, and though this
preserves the choice of ample bundle, it does {\em not} preserve stability
in general.

One way to avoid this is to replace ``stable'' with ``simple'', having no
nontrivial endomorphisms.  As simplicity is intrinsic, we avoid any
consideration of ample line bundles entirely, yet a result of Altman and
Kleiman means that we still have a reasonably well-behaved moduli space.

Let $X/S$ be locally projective, finitely presented morphism of schemes.  A
{\em family of simple sheaves} on $X$ is a sheaf $M$ on $T\times_S X$, flat
over $T$, such that the fiber $M$ over every geometric point $t\in T$
satisfies $\End(M(t))\cong k(t)$.  We consider two such families $M$, $M'$
equivalent if there is a line bundle ${\cal L}$ on $T$ such that $M'\cong
M\otimes \pi_1^*{\cal L}$; this naturally preserves isomorphism classes of
fibers.  Note that the simplicity condition implies that the line bundle
${\cal L}$ is unique (up to isomorphism) if it exists, since for a simple
family, we have
\[
\Hom_{T\times_S X}(M,M\otimes \pi_1^*{\cal L})\cong \Gamma(T\times_S
X,{\cal L})
\]

Roughly speaking, this moduli problem is represented by an algebraic space.
More precisely, this needs to be extended to the \'etale topology.  Define
a {\em twisted family of simple sheaves} on $X$ parametrized by $T$ to be a
family of simple sheaves parametrized by some \'etale cover $U\to T$ such
that the two induced families on $U\times_T U$ are equivalent.

\begin{thm}\cite{AltmanAB/KleimanSL:1980}
  There is a quasi-separated algebraic space $\Spl_{X/S}$ locally finitely
  presented over $S$ which represents the moduli functor of simple sheaves,
  in the sense that there is a natural bijection between twisted families
  of simple sheaves on $X$ and morphisms to $\Spl_{X/S}$.
\end{thm}

In particular, we obtain such an algebraic space associated to any Poisson
surface, and we wish to show that it is Poisson.  Of course, this
encounters another difficulty: what does it mean for an algebraic space to
be Poisson?  This is less of a difficulty than one might think, since as we
have seen, Poisson structures can be pulled back through \'etale morphisms
(and if a Poisson biderivation descends, the image is clearly Poisson).  In
other words, the notion of a Poisson structure makes perfect sense in the
\'etale topology.  We thus obtain the following definition of a Poisson
structure on an algebraic space.

\begin{defn}
Let ${\cal X}$ be an algebraic space.  A Poisson structure on ${\cal X}$ is
an assignment of a Poisson structure to the domain of every \'etale morphism
$f:U\to {\cal X}$, such that if $g:U'\to U$ is another \'etale morphism,
then $g:(U',\alpha_{f\circ g})\to (U,\alpha_f)$ is a Poisson morphism.
\end{defn}

\begin{rem}
If ${\cal X}$ is a Poisson scheme, then we obtain a Poisson structure in
the above sense by assigning to every \'etale morphism the induced Poisson
structure on its domain.  Thus any Poisson scheme remains Poisson as an
algebraic space.  In the other direction, the same argument lets one make
sense of Poisson structures on arbitrary Deligne-Mumford stacks.
\end{rem}

Thus to obtain our Poisson moduli space, we need simply define a suitably
canonical Poisson structure on the base of every twisted family of simple
sheaves corresponding to an \'etale morphism to $\Spl_{X/S}$.  In fact, it
will be enough to consider only untwisted families, as the required
compatibility of Poisson structures immediately gives us the conditions
needed to descend the Poisson structure from the \'etale cover of the base
of the twisted family.  Note that a family corresponds to an \'etale
morphism iff it is formally universal.

In \cite{BottacinF:2000}, Bottacin showed that the moduli space of stable
vector bundles on a Poisson surface has a well-behaved foliation by
algebraic symplectic leaves: for any bundle $V$, the subscheme
parametrizing stable vector bundles $V'$ with $V'|_{C_\alpha}\cong
V|_{C_\alpha}$ is a symplectic Poisson subscheme (where $C_\alpha$ is the
anticanonical curve).  This description does not quite carry over in full
generality, although it is straightforward enough to fix: with this in
mind, we define the {\em derived restriction} of a sheaf to be the
complex
\[
M|^{\dL}_{C_\alpha}:=\dL i^*M,
\]
where $i:C_\alpha\to X$ is the inclusion morphism.

We will mainly consider the case of a surface over a separably closed
field, which we suppress from the notation.

\begin{thm}\label{thm:poisson}
  Let $(X,\alpha)$ be a Poisson surface which is not symplectic and not
  birational to a quasi-standard Poisson surface.  Then the moduli space
  $\Spl_X$ has a natural Poisson structure, and for any complex
  $M^\bullet_\alpha$ of locally free sheaves on the associated
  anticanonical curve $C_\alpha$, the subspace of $\Spl_X$ parametrizing
  sheaves $M$ with $M|^{\dL}_{C_\alpha}\cong M^\bullet_\alpha$ is a
  (smooth) symplectic Poisson subspace.
\end{thm}

\begin{rems}
With this in mind, we define a {\em symplectic leaf} of $\Spl_X$ to be
the subspace corresponding to some fixed derived restriction, or more
generally a union of components of such a subspace.
\end{rems}

\begin{rems}
  Of course, this continues to hold when $X$ has symplectic fibers (subject
  to the technical condition that $\Pic^0(X)$ is smooth), with $\Spl_X$
  itself smooth and symplectic.  (See, e.g., Chapter 10 of
  \cite{HuybrechtsD/LehnM:1997}.)  The arguments below encounter
  difficulties in the symplectic case, however.  The claim fails when $X$
  is birational to a quasi-standard surface.
\end{rems}

\begin{rems}
  Since the biderivation can be defined over an arbitrary base, and the
  further conditions to be a Poisson structure are closed on the base, we
  conclude that $\Spl_{X/S}$ is Poisson for any Poisson surface over a
  reduced scheme $S$, such that no fiber of $X$ is symplectic.  The claim
  for the symplectic leaves is harder to generalize, mainly since it is no
  longer true in general that conditions of the form
  $M|^{\dL}_{C_\alpha}\cong M^\cdot_\alpha$ are locally closed.  It appears
  the correct condition to place on a family $M^\cdot_\alpha$ over a
  reduced base is that it can be represented by a two-term perfect complex,
  and that
\[
\dim\Hom(M^\cdot_\alpha,M^\cdot_\alpha)
-
\dim\Ext^{-1}(M^\cdot_\alpha,M^\cdot_\alpha)
\]
is constant.  (By the proof of Proposition \ref{prop:symp_dims} below, this
condition is necessary for the tangent space to the symplectic leaf to have
constant dimension.)  With this proviso, the full Theorem most likely holds
over an arbitrary reduced base, and a careful restatement of the constancy
condition should extend to any locally Noetherian base.  (For $X$ rational
(and probably for ruled surfaces over $C$ with $g(C)=1$), the reduced case
is sufficient, since suitable moduli stacks of anticanonical surfaces are
reduced; for $g(C)>1$, though, most components of the moduli functor have
obstructed deformations.)
\end{rems}

\begin{rems}
  This can be further generalized: instead of simple {\em sheaves}, one can
  consider simple {\em complexes}, or more precisely objects in
  $D^b_{\coh}(X/S)$ such that the natural map $\sO_S\to R\End(M,M)\to
  \tau_{\le 0} R\End(M,M)$ is a quasi-isomorphism.  It turns out this
  moduli problem is still represented by an algebraic space (see the
  discussion in Section \ref{sec:moduli_simple_noncomm} below), and
  moreover \cite{nc-lagrangian} that, at least in characteristic 0, the
  algebraic space is Poisson, with symplectic leaves given by derived
  restriction to the anticanonical curve.  (This uses the notion of shifted
  Lagrangian structures on derived stacks, which have only been defined in
  characteristic 0, but of course Lemma \ref{lem:closure_is_Poisson}
  immediately extends the result to the Zariski closure of the
  characteristic 0 locus.)
\end{rems}

In other words, we need to construct a Poisson structure on the base of
every formally universal family of simple sheaves, and show that the
foliation of the base by isomorphism class of derived restrictions to the
anticanonical curve is a foliation by symplectic Poisson subspaces.

The first requirement is to compute the cotangent sheaf.

\begin{lem}
Let $M$ be a formally universal family of simple sheaves on $X$
parametrized by $U$.  Then there is a natural isomorphism
\[
\Omega_U\cong \sExt^1_U(M,M\otimes \omega_X).
\]
\end{lem}

\begin{proof}
Fix an ample line bundle $\sO_X(1)$ on $X$, and for each pair $(m,n)$ of
nonnegative integers, consider the fine moduli space of subsheaves
$V\subset \sO_X(-m)^{\oplus n}$ such that the quotient $M$ is simple and
acyclic for $R\Hom(\sO_X(-m),{-})$, and the induced map
\[
\Hom(\sO_X(-m),\sO_X(-m)^{\oplus n})\to \Hom(\sO_X(-m),M)
\]
is an isomorphism.  On the one hand, this moduli space can be directly
identified as a disjoint union of open subschemes of Quot schemes, so
Theorem 3.1 of \cite{LehnM:1998} describes its cotangent sheaf.  On the
other hand, taking the quotient yields a map to $\Spl_X$ making it a
principal $\PGL_n$-bundle over an open subspace of $\Spl_X$.  We may thus
proceed as in the proof of Theorem 4.2 op.~cit.~to calculate the cotangent
sheaf of that subspace of $\Spl_X$.  Since as $m$, $n$ varies these
subspaces cover $\Spl_X$, the result follows.
\end{proof}

The required biderivation on $U$ is then given (following
\cite{TyurinAN:1988}) by
\[
\Omega_U^{\otimes 2}
\cong
\sExt^1_U(M,M\otimes \omega_X)^{\otimes 2}
\xrightarrow{1\otimes \alpha}
\sExt^1_U(M,M\otimes \omega_X)
\otimes
\sExt^1_U(M,M)
\to
\sO_U,
\]
where the last map is the trace pairing (i.e., Serre duality for
the smooth morphism $X\times U\to U$).  Note the commutative diagram
\[
\begin{CD}
\sExt^1_{U}(M,M\otimes\omega_X)^{\otimes 2}
@>>{1\otimes\alpha}>
\sExt^1_{U}(M,M\otimes\omega_X)
\otimes
\sExt^1_{U}(M,M)
\\
@VV{\alpha\otimes 1}V  @VVV
\\
\sExt^1_{U}(M,M)
\otimes
\sExt^1_{U}(M,M\otimes\omega_X)
@>>> \sO_{U}
\end{CD}
\]

Since the above biderivation is clearly functorial in $U$, and is invariant
under twisting $M$ by any line bundle, so in particular by line bundles
pulled back from $U$, it remains only to show that it is Poisson, and to
verify the claim about symplectic leaves.

Following \cite{HurtubiseJC/MarkmanE:2002b}, there are three cases to
consider, depending on the homological dimension of the simple sheaf.
Simple sheaves of homological dimension 2 are necessarily structure sheaves
of closed points, and thus the corresponding component of
$\Spl_X$ is naturally isomorphic to $X$ itself, and the induced
Poisson structure is the same.  The case of simple sheaves of homological
dimension $<2$ can then be reduced to that of locally free simple sheaves.
This reduction applies equally well in finite characteristic, though we
will need to do some more work (beyond the arguments of
\cite{HurtubiseJC/MarkmanE:2002b}) to identify the symplectic leaves.  We
will consider this case shortly, but will first deal with locally free
sheaves.

\section{Poissonness for locally free sheaves}

The argument of \cite{BottacinF:1995} could most likely be carried over to
finite characteristic, but is sufficiently complicated to make it somewhat
difficult to verify this, especially in characteristic 2.  In light of
Lemma \ref{lem:closure_is_Poisson}, alternation and the Jacobi identity in
characteristic 0 implies alternation and the Jacobi identity on the Zariski
closure of characteristic 0.  But for vector bundles, this is everything!

Now, let $\Vect_X$ denote the subspace of $\Spl_X$ classifying simple
locally free sheaves; a sheaf being locally free is an open condition, so
$\Vect_X$ is open in $\Spl_X$.

\begin{lem}
The algebraic space $\Vect_X$ is smooth and Poisson.
\end{lem}

\begin{proof}
It suffices to prove that the moduli problem is formally smooth;
i.e., any infinitesimal deformation of a vector bundle can be extended.
But the obstructions to deformations of $V$ are classified by
\[
\Ext^2(V,V)\cong \Hom(V,V\otimes\omega_X)^*\cong \Hom(V,V(-C_\alpha))^*.
\]
Since by assumption $\Hom(V,V)$ is generated by the identity, 
the injective map
\[
\Hom(V,V(-C_\alpha))\to \Hom(V,V)
\]
must be 0.

Since this moduli problem is formally unobstructed, for any nonsymplectic
Poisson surface $X$ over a dvr $R$, any simple bundle on $X_k$ extends to
$X$, and thus $\Vect_X\times_{\Spec(R)}\Spec(K)$ is dense in $\Vect_X$.  If
$K$ has characteristic 0, it follows from \cite{BottacinF:1995} that
$\Vect_{X_K}$ is Poisson, and thus so is $\Vect_{X_k}$.
\end{proof}

\begin{rems}
  To be precise, Bottacin only considers stable bundles, but the argument
  there only uses simplicity.
\end{rems}

\begin{rems}
  Note that the smoothness argument fails for symplectic surfaces, as well
  as in the quasi-standard case (where the Poisson structure does not lift
  to characteristic 0).
\end{rems}

The identification of symplectic leaves is somewhat more delicate,
especially since we want to know that they are smooth.  It is helpful to
first consider the case of invertible sheaves.

\begin{prop}\label{prop:Pic_inj}
The natural restriction map
\[
\Pic^0(X)\to \Pic^0(C_\alpha)
\]
is an injective map of group schemes, and the natural Poisson structure is
trivial on $\Pic(X)\subset \Spl_X$.
\end{prop}

\begin{proof}
  It follows from the classification of Poisson surfaces that either $X$ is
  rational, or there is a splitting of the induced map from $C_\alpha$ to
  the base $C$ of the ruling.  If $X$ is rational, then $\Pic^0(X)=1$, and
  there is nothing to prove, while if $C_\alpha\to C$ splits, then the
  composition
\[
\Pic^0(X)\to \Pic^0(C_\alpha)\to \Pic^0(C)
\]
agrees with the natural isomorphism.

Taking derivatives gives an injection
\[
H^1(\sO_X)\to H^1(\sO_{C_\alpha})
\]
implying that
\[
H^1(\alpha):H^1(\omega_X)\to H^1(\sO_X)
\]
is 0.  But this implies for invertible $M$ that the induced map
\[
\Ext^1(M,M\otimes\omega_X)\to \Ext^1(M,M)
\]
is 0, and thus so is the Poisson structure.
\end{proof}

\begin{rem}
  Note that this again fails when $X$ is birational to a quasi-standard Poisson
  surface.  Indeed, in that case, the composition
  \[
  \Pic^0(C)\cong \Pic^0(X)\to \Pic^0(C_\alpha)
  \]
  is the Frobenius isogeny, and thus has nonreduced kernel ($\mu_2$
  or $\alpha_2$, depending on whether $J(C)$ is ordinary or
  supersingular).  In particular, the fibers of restriction to
  $C_\alpha$ can fail to be smooth, and thus Theorem \ref{thm:poisson}
  fails in this case.  Moreover, since the map $H^1(\omega_X)\to
  H^1(\sO_X)$ induced by $\alpha$ is an isomorphism, we find that
  the natural biderivation on the $1$-dimensional scheme $\Pic^0(X)$
  induces a nondegenerate pairing on fibers of the cotangent bundle.
  In particular, the natural biderivation does not give a Poisson
  structure, as it is not even alternating.  Note that we can exclude
  both this case and the symplectic cases on which the conclusion
  of Theorem \ref{thm:poisson} fails by insisting that the kernel
  of $\Pic^0(X)\to \Pic^0(C_\alpha)$ be reduced.
\end{rem}

\begin{cor}
If $X$ is rationally ruled over the smooth curve $C$, then
$h^1(\sO_{C_\alpha})=g(C)+1$.
\end{cor}

\begin{proof}
Indeed, one has a short exact sequence
\[
0\to H^1(\sO_X)\to H^1(\sO_{C_\alpha})\to H^2(\omega_X)\to 0
\label{ses:cor3.5}
\]
where the first map is injective by the proposition.  Since $h^1(\sO_X)=g(C)$,
$h^2(\omega_X)=1$, the claim follows.
\end{proof}

\begin{rem}
  Again, in the quasi-standard case, this fails with $h^1(\sO_{C_\alpha})=1$.
\end{rem}

\begin{lem}\label{lem:leaves_of_vect_are_smooth}
  For any locally free sheaf $V_\alpha$ on $C_\alpha$, the subspace of
  $\Vect_X$ on which $V\otimes \sO_{C_\alpha}\cong V_\alpha$ is locally closed
  and smooth, and the corresponding conormal sheaf is the radical of the
  natural bilinear form on $\Omega_{\Vect_X}$.
\end{lem}

\begin{proof}
  That the subspace is locally closed is standard.
%
%
For smoothness, we need
  to show that any obstruction to a deformation inside this space must
  vanish.  As in \cite{BottacinF:1995,BottacinF:2000}, we obtain a
  self-dual long exact sequence the middle of which is
\[
\to
 \Hom(V_\alpha,V_\alpha)
\to \Ext^1(V,V\otimes\omega_X)\to \Ext^1(V,V)\to \Ext^1(V_\alpha,V_\alpha)\to
 \Ext^2(V,V\otimes \omega_X)\to
\label{eq:le_vect}
\]
Since $\Vect_X$ is smooth, any infinitesimal deformation of $V$ can be
extended, and the question is whether the extension can be chosen in such a
way that $V_\alpha$ remains constant.  Equivalently, any such extension
determines (as an extension of the trivial deformation of $V_\alpha$) a
class in $\Ext^1(V_\alpha,V_\alpha)$, which can be made trivial iff it is
in the image of $\Ext^1(V,V)$, or equivalently iff its image in
$\Ext^2(V,V\otimes\omega_X)$ is trivial.  Since $V$ is simple, the trace
map
\[
\Tr:\Ext^2(V,V\otimes\omega_X)\to H^2(\omega_X)
\]
is an isomorphism (it is dual to the natural map $k\to \End(V)$).  It is
thus equivalent to check whether the trace of the class in
$\Ext^1(V_\alpha,V_\alpha)$ is in the image of $\Tr(\Ext^1(V,V))$.  But by
\cite{ArtamkinIV:1988}, the trace map is the map induced by $V\mapsto
\det(V)$ on the corresponding spaces of deformations.  In particular the
trace of the given class in $\Ext^1(V_\alpha,V_\alpha)$ is the class of the
associated deformation of
\[
\det(V_\alpha)\cong \det(V)|_{C_\alpha}.
\]
We thus reduce to the case that $V_\alpha$, $V$ are invertible sheaves, in which
case it follows by Proposition \ref{prop:Pic_inj} that the deformation of
$V$ must be trivial, so certainly extends.

Now, by smoothness, $U$ is an \'etale neighborhood in the given locally closed
subspace iff we have an exact sequence
\[
0\to {\cal T}_U\to \sExt^1_U(V,V)\to \sExt^1_U(V_\alpha,V_\alpha),
\]
where the first map is the Kodaira-Spencer map and the second map comes
from the relative analogue of the above long exact sequence.  Equivalently,
the dual sequence
\[
\sHom_U(V_\alpha,V_\alpha)\to \sExt^1_U(V,V\otimes\omega_X)\to \Omega_U\to 0
\]
must be exact, and thus the conormal sheaf is equal to the radical, as
required.
\end{proof}

\begin{rem}
  A variant of this argument applies to the moduli space of simple objects
  in $D^b_{\coh}(X)$, as long as $M|_{C_\alpha}\ne 0$.  See Proposition
  \ref{prop:leaves_are_smooth_nc1} below, which gives the noncommutative
  version of this statement, as well as Proposition
  \ref{prop:leaves_are_smooth_nc2}, which gives a variant of the statement
  applying to most cases with $M|_{C_\alpha}=0$.
\end{rem}

But this is precisely what we wanted to show: we have identified a
foliation by locally closed subspaces such that at any point, the cotangent
space to the leaf is the radical of the form given by the Poisson
structure.  In particular, each leaf is a Poisson subspace, and the induced
Poisson structure is nondegenerate, making it a symplectic leaf as
required.

\section{Sheaves of homological dimension $\le 1$}
\label{sec:moduli_hd1}

For general sheaves of homological dimension $\le 1$, we use the reduction
of \cite{HurtubiseJC/MarkmanE:2002b}.  This extends easily to finite
characteristic, but we will need to do some additional calculations to
recover the symplectic leaves.

Let $U\to \Spl_X$ be an \'etale neighborhood with corresponding family $M$
of simple sheaves, and suppose that the fibers of M have homological
dimension $\le 1$.  Let $\pi$ denote the projection $U\times X\to U$.

Twisting $M$ by a line bundle has no effect on the biderivation, so we
may assume that the natural map $\pi^*\pi_*M\to M$ is surjective, and both
$M$ and $M\otimes \omega_X$ are $\pi_*$-acyclic.  Then, as observed in
\cite{HurtubiseJC/MarkmanE:2002b}, one obtains a corresponding locally free
resolution
\[
0\to V\to W\to M\to 0
\]
with $W=\pi^*\pi_*M$ and $V$ simple, such that we can reconstruct $M$ from
$V$.  Indeed, acyclicity implies $\pi_*M$ is locally free, and the
homological dimension condition then ensures that $V$ is locally
free.  Moreover, the induced map
\[
\sHom_U(W,\sO_X)\to \sHom_U(V,\sO_X)
\]
is an isomorphism, and thus given $V$ we can recover $M$ as the cokernel of
the natural injection
\[
V\to \sHom_U(\sHom_U(V,\sO_X),\sO_X).
\]
It follows in particular that $V$ is also a family of simple sheaves,
giving an injective morphism $U\to \Vect_X$, and we need to compare the two
induced biderivations, then understand the symplectic leaves.

For the first part, we use a mild variant of the calculation in
\cite{HurtubiseJC/MarkmanE:2002b}.  The functoriality of long exact
sequences gives a commutative diagram
\[
\begin{CD}
\sHom_U(V,M)  @>>> \sExt^1_U(M,M)\\
@VVV                   @VVV\\
\sExt^1_U(V,V)@>>> \sExt^2_U(M,V)
\end{CD}
\]
Here, the top arrow is surjective, since
\[
\sExt^1_U(W,M)\cong \sHom(\pi_*(M),R^1\pi_*(M))=0,
\]
and the right arrow is injective since
\[
\sExt^1_U(M,W)\cong \sHom(R^1\pi_*(M\otimes\omega_X),\pi_*M)=0.
\]
There is also a natural map $\delta:\sExt^1_U(M,M)\to \sExt^1_U(V,V)$
coming from the above construction (the derivative of the map of moduli
spaces).  Since the top arrow above is surjective, we may represent any
self-extension of $M$ as the pushforward of the resolution $0\to V\to W\to
M\to 0$ by some morphism $\phi:V\to M$.  This gives a self-extension $M'$
as the cokernel in the short exact sequence
\[
0\to V\to W\oplus M\to M'\to 0.
\]
The surjection $W\to M$ then induces a surjection $W\oplus W\to M'$,
so an exact sequence
\[
0\to V'\to W\oplus W\to M'\to 0,
\]
where is a self-extension of $V$ (the image of $M'$ under $\delta$).
Pulling this back through the injection $M\to M'$ gives a short exact
sequence
\[
0\to V'\to V\oplus W\to M\to 0,
\]
so that $V'$ is the pullback through $\phi$ of $0\to V\to W\to M\to 0$.

In other words, given a class in $\sExt^1_U(M,M)$, we can compute its image
under $\delta$ by choosing a preimage in $\sHom_U(V,M)$ and mapping that to
$\sExt^1_U(V,V)$.  In particular, we can put the differential into the
diagram without breaking commutativity.  Note in particular that
commutativity implies that the differential is injective, and in fact
identifies $\sExt^1_U(M,M)$ with the image of the map $\sHom_U(V,M)\to
\sExt^1_U(V,V)$.  We similarly have a commutative square
\[
\begin{CD}
\sHom_U(V,M\otimes\omega_X)  @>>> \sExt^1_U(M,M\otimes\omega_X)\\
@VVV                   @VVV\\
\sExt^1_U(V,V\otimes\omega_X)@>>> \sExt^2_U(M,V\otimes\omega_X)
\end{CD}
\]
Since the connecting maps in this diagram are the duals (up to sign) of the
maps in the original diagram, we see that the dual of the differential also
fits into this diagram; more precisely, we must take the negative of the
dual in order to account for the signs.

Without the diagonal maps, the two squares fit into a commutative cube
induced by $\alpha:\omega_X\to \sO_X$.  We claim that this diagram remains
commutative when we introduce the diagonal maps.  In other words, we need
to show that the composition
\[
\sExt^1_U(V,V\otimes\omega_X)\xrightarrow{-\delta^*}
\sExt^1_U(M,M\otimes\omega_X)\xrightarrow{\alpha}
\sExt^1_U(M,M)\xrightarrow{\delta}
\sExt^1_U(V,V)
\]
agrees with the map induced by $\alpha$.  If we compose both maps with the
connecting map $\sExt^1_U(V,V)\to\sExt^2_U(M,V)$, they agree.  If we can show
that the image of $\sExt^1_U(V,V\otimes\omega_X)\to \sExt^1_U(V,V)$ is
contained in the image of $\delta$, we will be done, since then the error
is in the image of $\delta$, which injects in $\sExt^2_U(M,V)$.  The image
of $\delta$ is the same as the image of $\Hom_U(V,M)$, and thus is the same
as the kernel of the map $\sExt^1_U(V,V)\to \sExt^1_U(V,W)$.  By
commutativity of long exact sequences, it will thus suffice to show that
the map $\sExt^1_U(V,W\otimes\omega_X)\to \sExt^1_U(V,W)$ is 0.  We have a
commutative square
\[
\begin{CD}
\sExt^1_U(W,W\otimes\omega_X)@>\alpha >> \sExt^1_U(W,W)\\
       @VVV                              @VVV\\
\sExt^1_U(V,W\otimes\omega_X)@>\alpha >> \sExt^1_U(V,W)
\end{CD}
\]
The left arrow is an isomorphism, by duality from the fact that
$R^1\pi_*(V)\cong R^1\pi_*(W)$, while the top arrow is $0$ by Proposition
\ref{prop:Pic_inj}, since $\pi_*\omega_X\to \pi_*\sO_X$ is the 0 morphism.
It follows that the bottom arrow is 0 as required.

We have thus shown that the above map of moduli spaces identifies the two
biderivations up to sign, and thus since the biderivation on $\Vect_X$ is
Poisson, so is the biderivation on $U$.  (The minus sign above makes this
an anti-Poisson morphism.)  There are two things remaining for us to do:
first, we need to show that any locally closed subscheme of $U$ on which
$M|^{\dL}_{C_\alpha}$ is fixed is a Poisson subscheme, and second that it
is symplectic.  Of course, we will show both at once if we can identify the
putative leaves of $U$ with the leaves of $\Vect_X$.

Now, $M|^{\dL}_{C_\alpha}$ is represented by the complex $V_\alpha\to
W_\alpha$, where $V_\alpha$ and $W_\alpha$ are the respective restrictions to
$C_\alpha$.  We claim that, just as $V$ determines $M$, $V_\alpha$ determines
$M|^{\dL}_{C_\alpha}$.  Consider the commutative diagram
\[
\begin{CD}
0@>>>\Hom(W,\sO_X)@>>>\Hom(W_\alpha,\sO_{C_\alpha})@>>>\Ext^1(W,\omega_X)@>>>\Ext^1(W,\sO_X)\\
@. @VVV             @VVV                @VVV               @VVV  \\
0@>>>\Hom(V,\sO_X)@>>>\Hom(V_\alpha,\sO_{C_\alpha})@>>>\Ext^1(V,\omega_X)@>>>\Ext^1(V,\sO_X)
\end{CD}
\]
coming from the exact sequence
\[
0\to \omega_X\xrightarrow{\alpha} \sO_X\to \sO_{C_\alpha}\to 0.
\]
Each vertical map fits into a long exact sequence; since $M$ and $M\otimes
\omega_X$ are acyclic, we may apply Serre duality to find that the first
and third vertical maps are isomorphisms, while the fourth map is
injective.  We thus find that the given map $V_\alpha\to W_\alpha$ induces an
isomorphism
\[
\Hom(W_\alpha,\sO_{C_\alpha})\cong \Hom(V_\alpha,\sO_{C_\alpha}).
\]
In particular, $\Hom(V_\alpha,\sO_{C_\alpha})$ is a free
$\End(\sO_{C_\alpha})$-module, and the map $V_\alpha\to W_\alpha$ can be
identified with the natural map
\[
V_\alpha\to \Hom_{\End(\sO_{C_\alpha})}(\Hom(V_\alpha,\sO_{C_\alpha}),\sO_{C_\alpha}).
\]
It follows that the symplectic leaves of $\Vect_X$ pull back to closed
subspaces of the putative symplectic leaves of $U$.

Conversely, we can recover $V_\alpha$ from $M|^{\dL}_{C_\alpha}$ and the
numerical invariants of $M$.  (This also implies that the condition
$M|^{\dL}_{C_\alpha}\cong M^\cdot_\alpha$ is locally closed.)  Since $M$
and $M\otimes \omega_X$ are acyclic, we have a four-term exact sequence
\[
0\to H^{-1}(M|^{\dL}_{C_\alpha})\to H^0(M\otimes \omega_X)\to H^0(M)\to
H^0(M|^{\dL}_{C_\alpha})\to 0,
\]
and thus (locally on the base) we may write
\[
W_\alpha \cong H^0(M|^{\dL}_{C_\alpha})\otimes \sO_{C_\alpha}\oplus
\sO_{C_\alpha}^l
\]
where
\[
l = h^0(M|^{\dL}_{C_\alpha})-h^0(M) = h^0(M|^{\dL}_{C_\alpha})-\chi(M).
\]
Moreover, this isomorphism identifies the map $W_\alpha\to
M|^{\dL}_{C_\alpha}$ with the direct sum of the natural morphism
\[
H^0(M|^{\dL}_{C_\alpha})\otimes \sO_{C_\alpha}\to M|^{\dL}_{C_\alpha}
\]
and the zero morphism
\[
\sO_{C_\alpha}^l\to M|^{\dL}_{C_\alpha}.
\]
But then we can recover $V_\alpha$ as the ``kernel'' of this morphism; more
precisely, the mapping cone of this morphism is a complex representing the
sheaf $V_\alpha$.  In particular, we have
\[
V_\alpha \cong \sO_{C_\alpha}^l\oplus V_{\min},
\]
where
\[
V_{\min}\to H^0(M|^{\dL}_{C_\alpha})\otimes \sO_{C_\alpha}
\]
is the canonical complex quasi-isomorphic to $M|^{\dL}_{C_\alpha}$.  (Of
course, when $\Tor_1(M,\sO_{C_\alpha})=0$, so $M|^{\dL}_{C_\alpha}$ is a
sheaf, $V_{\min}$ is an actual kernel, of the natural (surjective!)
morphism $H^0(M|_{C_\alpha})\otimes \sO_{C_\alpha}\to M|_{C_\alpha}$.)

We thus conclude that the putative symplectic leaf of $\Spl_X$ containing
$M$ can be identified with a Poisson subspace of a symplectic leaf of
$\Vect$.  It remains only to show that this Poisson subspace is open.
(This is automatic in characteristic 0, since then any Poisson subspace of a
symplectic space is open.  The reader only interested in characteristic 0
may still find the direct argument instructive, however.)

Fix a simple sheaf $M$ satisfying the above conditions, with corresponding
locally free sheaf $V$.  In the symplectic leaf of $\Vect_X$ corresponding
to $V_\alpha$, we first impose the following open conditions on the fibers
$V'$:
\begin{itemize}
\item[1.] $V'$ has the same numerical Chern class as $V$.  (This is both
  open and closed.)
\item[2.] $H^0(V')=H^2(V')=0$.  (This is open by semicontinuity)
\item[3.] $\Ext^2(V',\sO_X)=0$ and $\dim\Hom(V',\sO_X)\le
  \rank(W)$. (Also semicontinuity).
\end{itemize}
If $U'$ is an \'etale neighborhood in this open subspace, then
$\sExt^1_{U'}(V',\omega_X)$ is locally free, since the other $\Ext$ groups
are trivial, and thus it has the same rank as
\[
\Ext^1(V,\omega_X)\cong \Ext^1(W,\omega_X)\cong H^1(W)^*
\]
Since $W$ is trivial, we conclude
\[
\rank(\sExt^1_{U'}(V',\omega_X)) = g\rank(W),
\]
where $g=\dim H^1(\sO_X)$ is the genus of the curve over which $X$ is
rationally ruled.  Also, since $V_\alpha$, $W_\alpha$ are fixed and
$h^0(\sO_{C_\alpha})=g+1$, we may compute
\[
h^0(V_\alpha^*) = h^0(W_\alpha^*) = (g+1)\rank(W).
\]
But then the long exact sequence
\[
0\to \sHom_{U'}(V',\sO_X)\to \sHom_{U'}(V',\sO_{C_\alpha})
 \to \sExt^1_{U'}(V',\omega_X)\to\cdots
\]
implies that $\sHom_{U'}(V',\sO_X)$ contains a locally free sheaf of rank
at least $\rank(W)$.  Since we have already imposed an upper bound on the
rank, we conclude that $\sHom_{U'}(V',\sO_X)$ is locally free of rank
$\rank(W)$.  We thus recover an injective map
\[
V'\to \sHom_{U'}(\sHom_{U'}(V',\sO_X),\sO_X).
\]
If we let $M'$ be the cokernel of this map, then it certainly is generated
by global sections, and we can further impose the open conditions that $M'$
and $M'\otimes \omega_X$ are $\pi_*$-acyclic.  We thus obtain a
neighborhood of $V$ in its symplectic leaf such that every sheaf in that
neighborhood comes from a sheaf $M'$.

This completes the proof of Theorem \ref{thm:poisson}.
\qed

\begin{rem}
  One can perform a similar reduction for simple sheaves of homological
  dimension $2$, and obtain an anti-Poisson map from the moduli space of
  structure sheaves of points to the moduli space of ideal sheaves of
  points.  In particular, one can start with {\em any} simple sheaf and
  perform the above reduction twice to obtain a Poisson map to $\Vect_X$.
\end{rem}

\section{Birational maps and moduli spaces}

In this section, $\pi:X\to Y$ will be a Poisson birational morphism of
Poisson surfaces, which are now over a more general base scheme $S$.
For convenience of notation, we will silently identify sheaves on $S$ with
their pullbacks to $X$ or $Y$ as appropriate.

Although our goal is to show that $\pi^{*!}$ respects the Poisson
structure (in a suitable sense), we will in fact mainly focus on $\pi_*$
instead; since $\pi^{*!}$ is an inverse to $\pi_*$, it will be easy to
derive Poissonness of $\pi^{*!}$ from Poissonness of $\pi_*$, but the latter
will also let us deal with pseudo-twists.

The first issue is that direct images do not in general preserve flatness
of families.  This is, of course, just a relative version of semicontinuity
questions, but is particularly easy to deal with in our case.

\begin{lem}
  Let $M$ be a coherent sheaf on $X$, flat over $S$.  If $R^1\pi_*M$ is
  flat over $S$, then so is $\pi_*M$, and moreover the natural map
\[
\pi_*M\otimes N\to \pi_*(M\otimes N)
\]
is an isomorphism for all coherent sheaves $N$ on $S$.  In particular, this holds if
every fiber of $M$ is $\pi_*$-acyclic.
\end{lem}

\begin{proof}
For any coherent sheaf $N$ on $S$, we have an isomorphism
\[
\dR\pi_*M\otimes^{\dL} N\cong \dR\pi_*(M\otimes^{\dL} N)\cong \dR\pi_*(M\otimes N).
\]
Since $\pi$ has fibers of dimension $\le 1$, the spectral sequence for
$\dR\pi_*M\otimes^\dL N$ has entries in only two rows.  We thus obtain
isomorphisms
\[
\Tor_{p+2}(R^1\pi_*M,N)\cong \Tor_p(\pi_*M,N)
\]
for $p>0$, an isomorphism
\[
R^1\pi_*M\otimes N\cong R^1\pi_*(M\otimes N),
\]
and an exact sequence
\[
0\to \Tor_2(R^1\pi_*M,N)\to \pi_*M\otimes N\to \pi_*(M\otimes N)\to
\Tor_1(R^1\pi_*M,N)\to 0.
\]
The claim follows immediately.
\end{proof}

It follows that $\pi_*$ induces a morphism from an open subspace of $\Spl_X$
to $\Spl_Y$.  More precisely, we have the following.

\begin{lem}
  Let $\Spl_{X,\pi}\subset \Spl_X$ be the subspace of simple sheaves $M$
  which are $\pi_*$-acyclic and have simple direct image.  Then
  $\Spl_{X,\pi}$ is an open subspace, and $\pi_*$ induces a surjective
  morphism from $\Spl_{X,\pi}$ to $\Spl_Y$.
\end{lem}

\begin{proof}
  That $\pi_*$-acyclicity is an open condition follows from the fact that
  for any $S$-flat family $M$, the functor $R^1\pi_*$ commutes with taking
  fibers; thus the $\pi_*$-acyclic locus is just the complement of the
  image of the closed subscheme $\Fitt_0(R^1\pi_*M)$ under the proper map
  $X\to S$.  Then $\pi_*M$ is a flat family, and simplicity is an open
  condition on flat families \cite{AltmanAB/KleimanSL:1980}.

  For surjectivity, note that point sheaves in $\Spl_Y$ are direct images
  of point sheaves in $\Spl_X$, while any other simple sheaf is the direct
  image of its minimal lift.
\end{proof}

We will need one more fact in the proof of Poissonness, which we will also
use in the next section, so separate out from the proof.

\begin{lem}\label{lem:useful_factorization}
Let $M$ be a $\pi_*$-acyclic sheaf on $X$ with direct image of homological
dimension $\le 1$.  Then there is a commutative diagram
\[
\begin{CD}
\pi^!\pi_*M\otimes \pi^*\omega_Y @>1\otimes \pi^*\alpha_Y>> \pi^!\pi_*M\\
 @VVV @AAA\\
  M\otimes \omega_X @>1\otimes \alpha_X>> M,
\end{CD}
\]
where the first vertical map factors as
\[
\pi^!\pi_*M\otimes \pi^*\omega_Y
\to
\pi^*\pi_*M\otimes \omega_X
\to
M\otimes \omega_X.
\]
\end{lem}

\begin{proof}
Since $\pi$ is an isomorphism outside the exceptional locus, the diagram
can only fail to commute on the exceptional locus.  The failure to commute
is thus measured by a morphism
\[
\pi^!\pi_*M\otimes \pi^*\omega_Y \to \pi^!\pi_*M
\]
that vanishes outside the exceptional locus.  We can thus apply Proposition
\ref{prop:invisible_image} to conclude that this morphism has invisible
image.  Since $\pi^!\pi_*M$ has no invisible subsheaf, the image is 0,
and thus the diagram commutes.
\end{proof}

\begin{thm}\label{thm:direct_image_is_Poisson}
  Let $\pi:X\to Y$ be a Poisson birational morphism of Poisson surfaces.
  The direct image functor $\pi_*$ defines a Poisson morphism
  $\pi_*:\Spl_{X,\pi}\to \Spl_Y$.
\end{thm}

\begin{proof}
  We have already shown that it defines a morphism, so it remains only to
  show that it respects the Poisson structure.  Consider a sheaf $M$ in
  $\Spl_{X,\pi}$.  If $\pi_*M$ has homological dimension $2$, then $\pi_*M$
  must be a point, and either $M$ is a point or $M$ is supported on the
  exceptional locus.  On the subspace $X\subset \Spl_X$ parametrizing point
  sheaves, $\pi_*$ is just $\pi$, so is certainly Poisson.  If $M$ is
  supported on the exceptional locus, then $\pi_*M$ is supported on the
  anticanonical curve, so that for such sheaves, $\pi_*$ maps to a single
  point with trivial Poisson structure, so again is Poisson.

  We may thus restrict our attention to the case that $\pi_*M$ has
  homological dimension $\le 1$.  Now, the adjunction between $\pi^!$ and
  $\pi_*$ induces natural maps
\begin{align}
  \Ext^1(\pi_*M,\pi_*M\otimes \omega_Y)&\cong \Ext^1(M,\pi^!\pi_*M\otimes \pi^*\omega_Y)\notag\\
  \Ext^1(\pi_*M,\pi_*M)&\cong \Ext^1(M,\pi^!\pi_*M),\notag
\end{align}
with inverse given by $\pi_*$.  (To be precise, given a $\pi_*$-acyclic
sheaf $N$, there is a natural map $\pi_*:\Ext^1(M,N)\to
\Ext^1(\pi_*M,\pi_*N)$ given by taking the direct image of the extension.)
In particular, the Poisson structure on a neighborhood of $\pi_*M$ in
$\Spl_Y$ is obtained by composing the above isomorphisms with the map
\[
\Ext^1(M,\pi^!\pi_*M\otimes \pi^*\omega_Y)
\xrightarrow{\Ext^1(M,1\otimes \pi^*\alpha_Y)}
\Ext^1(M,\pi^!\pi_*M).
\]
By Lemma \ref{lem:useful_factorization}, this factors through the natural map
\[
\Ext^1(M,M\otimes \omega_X)
\xrightarrow{\Ext^1(M,1\otimes \alpha_X)}
\Ext^1(M,M),
\]
and the composition
\[
\Ext^1(M,M) \to \Ext^1(M,\pi^!\pi_*M)\cong \Ext^1(\pi_*M,\pi_*M)
\]
is just $\pi_*$, i.e., the differential of the morphism
$\pi_*:\Spl_{X,\pi}\to \Spl_Y$.  In other words, the map
$\Omega_{\Spl_Y}\to \Omega^*_{\Spl_Y}$ induced by the Poisson structure is
the direct image of the corresponding map on $\Spl_{X,\pi}$, and thus
$\pi_*$ is Poisson.
\end{proof}

\begin{rem}
  The analogous statement holds for the moduli space of simple objects in
  the derived category, with the condition on the domain being that the
  derived direct image is simple.
\end{rem}

\begin{cor}
The functors $\pi^*$ and $\pi^!$ induce Poisson morphisms from $\Spl_Y$ to
$\Spl_{X,\pi}$.
\end{cor}

\begin{proof}
Indeed, $\pi^*$ and $\pi^!$ are injective, and preserve simplicity and
flatness.  Since the images of $\pi^*$ and $\pi^!$ satisfy the hypotheses
of the Theorem, $\pi_*$ is Poisson on both images; since it is also an
isomorphism on both images, the inverses are Poisson.  But this is
precisely what we want to prove.
\end{proof}

\medskip

For the minimal lift, we again have an issue with flatness.  In this case,
however, we can completely control the corresponding flattening
stratification.  Let $\Spl_{Y,\le 1}$ denote the subspace parametrizing
sheaves of homological dimension $\le 1$.

\begin{lem}
  Suppose $\pi$ is monoidal, blowing up the point $p\in Y$.  Then for each
  integer $m\ge 0$, the subspace of $\Spl_{Y,\le 1}$ parametrizing sheaves
  $M$ with $\dim(\Ext^1(\sO_p,M))=m$ is a Poisson subspace, and contains
  every symplectic leaf it intersects.
\end{lem}

\begin{proof}
Since $\pi$ is Poisson, $p$ is contained in the anticanonical curve
$C_\alpha$, and thus
\[
M\otimes \sO_p\cong (M\otimes \sO_{C_\alpha})\otimes \sO_p
\]
is constant on symplectic leaves, so determines a stratification by locally
closed Poisson subspaces.  Since
\[
\dim(H^0(M\otimes \sO_p))
=
\dim(\Hom(M,\sO_p))
=
\dim(\Ext^1(\sO_p,M))+\rank(M),
\]
the same applies to $\dim(\Ext^1(\sO_p,M))$.
\end{proof}

Now, given any flat family of sheaves on $Y$ with homological dimension
$\le 1$, we obtain a stratification of the base of the family by taking the
flattening stratification of the minimal lift of the family.  Normally this
depends on a choice of relatively very ample bundle, but one can always
refine such a stratification by imposing the open and closed conditions
that the numerical Chern classes be constant, and the result will then be
independent of the relatively very ample bundle.  It is this canonical
stratification that we will mean; it refines the usual decomposition of $Y$
by numerical Chern class.
 
\begin{cor}
  If $\pi$ is monoidal, then the stratification induced on $\Spl_{Y,\le 1}$
  by $\pi$ agrees with the stratification by $\dim\Ext^1(\sO_p,M)$ and
  numerical Chern class.
\end{cor}

\begin{proof}
Indeed, for any coherent sheaf $M$ of homological dimension $\le 1$, we
have the equation
\[
[\pi^{*!}M]=\pi^*[M]-\dim\Ext^1(\sO_p,M)[\sO_e(-1)]
\]
in the Grothendieck group.  It follows that $\pi^{*!}M$ has constant
numerical Chern class iff $M$ has constant numerical Chern class and
$\dim\Ext^1(\sO_p,M)$ is constant.
\end{proof}

\begin{thm}
  Let $\pi:X\to Y$ be a Poisson birational morphism of Poisson surfaces.
  Every stratum of the stratification induced on $\Spl_{Y,\le 1}$ by $\pi$
  is a Poisson subspace, and the minimal lift induces a Poisson isomorphism
  from each stratum to an open subspace of $\Spl_{X,\le 1}$.
\end{thm}

\begin{proof}
  First suppose that $\pi$ is monoidal.  Then we have already shown that
  the stratification is Poisson.  On each stratum, the minimal lift
  preserves flatness (and simplicity), so defines a morphism from the
  stratum to $\Spl_{X,\le 1}$.  That the image is open follows from
  semicontinuity and the fact that minimal lifts are characterized by the
  vanishing of $\Hom(\sO_e(-1),M)$ and $\Hom(M,\sO_e(-1))$.

  We thus obtain an isomorphism between each stratum and the corresponding
  open subspace, and it remains only to show that it is Poisson.
  Equivalently, we need to show that $\pi_*$ is Poisson on the image of
  $\pi^{*!}$, but this is immediate from Theorem
  \ref{thm:direct_image_is_Poisson}.

  More generally, if we factor $\pi=\pi_1\circ\pi_2$ with $\pi_1$ monoidal,
  then the stratification induced by $\pi$ refines the stratification
  induced by $\pi_2$, and is identified via $\pi_2^{*!}$ with the
  stratification induced by $\pi_1$ on $\pi_1(X)$.  The claim then follows
  easily by induction.
\end{proof}

\begin{cor}
  Let $\pi:X\to Y$ be a birational morphism of Poisson surfaces.  Then for
  every stratum of the stratification induced by $\pi$ on $\Spl_{Y,\le 1}$,
  any pseudo-twist operation defines a Poisson morphism on the open
  subspace of the stratum where the pseudo-twist is simple.
\end{cor}

\begin{proof}
Indeed, a pseudo-twist is a composition of the Poisson morphism
$\pi^{*!}$, the Poisson morphism ${-}\otimes\sO_X(\pm e_f)$, and
the (partially defined) Poisson morphism $\pi_*$.
\end{proof}

\begin{cor}
  Let $(Y,\alpha)$ be a Poisson surface, let the subscheme $J\subset C_\alpha$
  be a disjoint union of jets, and let $\pi:X\to Y$ be the minimal
  desingularization of the blowup of $Y$ along $J$.  For any sheaf $M_\alpha$
  which is a direct sum of structure sheaves of subschemes of $J$,
  $\pi^{*!}$ is a symplectomorphism on the symplectic leaf
  $M\otimes \sO_{C_\alpha}=M_\alpha$ of $\Spl_{Y,\le 1}$.
\end{cor}

\begin{proof}
Indeed, $M\otimes \sO_{C_\alpha}=M_\alpha$ iff $\pi^{*!}M$ is disjoint from
$C_\alpha$, and has the correct numerical Chern class.
\end{proof}

\begin{rem}
In particular, if $C_\alpha$ is smooth, then any symplectic leaf consisting
of sheaves transverse to $C_\alpha$ is symplectomorphic to an open subspace of
some $\Spl_X$ with $\pi:X\to Y$ a Poisson birational morphism: simply choose
$\pi:X\to Y$ so that $\pi^{*!}M$ is disjoint from $C_\alpha$.
\end{rem}

\medskip

The fact that inverse images and minimal lifts behave nicely with respect
to the Poisson structure has an interesting consequence.  Let $f:X\ratto Y$
be a Poisson birational map.  Then we can use the above facts to
essentially embed $\Spl_X$ as a Poisson subspace of $\Vect_Y$.  More
precisely, suppose $U\to \Spl_{X\le 1}$ is a Noetherian \'etale
neighborhood (in the remaining case, when $U$ is an \'etale neighborhood of
$X$, simply consider the corresponding family of ideal sheaves).  Then
there is a corresponding morphism $U\to \Vect_Y$ (\'etale to its image)
such that the two induced Poisson structures are the same.

To see this, let $Z$ be a Poisson resolution of $f$, so that $f$ factors
through Poisson birational morphisms $g:Z\mapsto X$ and $h:Z\mapsto Y$.
Then if $M$ is the original flat family of sheaves on $X$, $g^* M$ is a
flat family of sheaves on $Z$, with the same Poisson structure.  We may
then apply the construction of \cite{HurtubiseJC/MarkmanE:2002b} (see
Section \ref{sec:moduli_hd1} above) to turn this into a flat family of
simple locally free sheaves on $Z$ (negating the Poisson structure in the
process).  Let $V$ be this family.  Now, if we twist $V$ by a suitable line
bundle, we may arrange to have
\[
\Hom(V,\sO_e(-1))=0
\]
for every $h$-exceptional component $e$. Indeed, we have
\[
\Hom(V,\sO_e(-1))
\cong
\Hom(V|_e,\sO_e(-1))
\cong
H^1(V|_e(-1))^*
\]
so we need simply twist by a relatively ample bundle that makes
$\bigoplus_e V|_e(-1)$ acyclic.  But this makes $V$ a minimal lift under
$h$!  Indeed, since $V$ is torsion-free, we also have $\Hom(\sO_e(-1),V)=0$
for all $e$, and thus Theorem \ref{thm:minimal_if_no_exceptional} applies.
In particular, it follows that $h_*V$ is simple, and induces an
anti-Poisson map from $U$ to $\Spl_{Y,\le 1}$.  Applying the construction
of \cite{HurtubiseJC/MarkmanE:2002b} again gives the required Poisson
morphism $U\to \Vect_Y$.  Each step of the above process is reversible, so
the maps $U\to \Spl_{X\le 1}$ and $U\to \Vect_Y$ have isomorphic images
(algebraic subspaces of $\Spl_{X\le 1}$, $\Vect_Y$ respectively).

\bigskip

In addition to images and lifts, there is one more natural morphism we want
to consider.  We have already discussed the dualization functor
\[
{-}^D:=\sExt^1_Y({-},\omega_Y)
\]
in the context of invisible sheaves, but of course much of what we have
said applies to any pure $1$-dimensional sheaf.  Duality also interacts
nicely with the inverse image and minimal lift functors.

\begin{prop}
Let $\pi:X\to Y$ be a birational morphism of smooth projective surfaces.
Then for any pure $1$-dimensional sheaf $M$ on $Y$, $\pi^*M$, $\pi^!M$,
and $\pi^{*!}M$ are pure $1$-dimensional sheaves, and we have
\begin{align}
(\pi^*M)^D&\cong \pi^!(M^D)\notag\\
(\pi^{*!}M)^D&\cong \pi^{*!}(M^D)\notag\\
(\pi^!M)^D&\cong \pi^*(M^D).\notag
\end{align}
\end{prop}

\begin{proof}
This is essentially by inspection using a locally free resolution
\[
0\to V\to W\to M\to 0
\]
of $M$, and the corresponding resolution
\[
0\to W^*\to V^*\to M^D\to 0
\]
of $M^D$.
\end{proof}

Since duality interacts nicely with these Poisson morphisms, it is natural
to suspect that it is itself Poisson. This is not quite the case: it is in
fact anti-Poisson.

\begin{prop}
Let $(X,\alpha)$ be a Poisson surface, and let $\Spl^1_X$ be the subspace
parametrizing pure $1$-dimensional sheaves.  Then ${-}^D$
is an anti-Poisson involution on $\Spl^1_X$.
\end{prop}

\begin{proof}
Recall that at a point $M\in \Spl^1_X$, the pairing on the sheaf of
differentials is given by the composition
\[
\Ext^1(M,M\otimes \omega_X)\otimes \Ext^1(M,M\otimes \omega_X)
\to
\Ext^2(M,M\otimes \omega_X^2)
\to
\Ext^2(M,M\otimes \omega_X)
\to
k.
\]
Now, ${-}^D$ induces an isomorphism
\[
\Ext^1(M,M\otimes \omega_X)
\to
\Ext^1(M^D\otimes \omega_X^{-1},M^D)
\to
\Ext^1(M^D,M^D\otimes \omega_X).
\]
Indeed, a class in $\Ext^1(M,M\otimes \omega_X)$ corresponds to an
extension
\[
0\to M\otimes \omega_X\to N\to M\to 0.
\]
Since $N$ is an extension of pure $1$-dimensional sheaves, it is also pure
$1$-dimensional, and thus we can dualize the entire exact sequence to
obtain
\[
0\to M^D\to N^D\to (M\otimes \omega_X)^D\to 0,
\]
and note that $(M\otimes \omega_X)^D\cong M^D\otimes \omega_X^{-1}$.

Similarly, the product of two such extensions is given by a four-term exact
sequence of pure $1$-dimensional sheaves, which we can again dualize.  If
$\phi$, $\psi$ are two classes in $\Ext^1(M,M\otimes \omega_X)$, then we have
\[
(\phi\psi)^D = \psi^D\phi^D = -\phi^D\psi^D.
\]
Duality respects the trace map, so all told gives an anti-Poisson involution.
\end{proof}

\section{Rigidity}
\label{sec:rigidity_comm}

In studying moduli spaces of differential and difference equations,
one particularly interesting question is when those moduli spaces
consist of a single point.  For instance, the hypergeometric
differential equation has been known since Riemann to be the unique
second-order Fuchsian equation with three singular points (and any
specific equation is determined by the exponents at the singularities).
In a generalized Hitchin system context, the order and singularity
structure precisely specifies a symplectic leaf in the corresponding
Poisson moduli space.  We are thus interested in understanding which
sheaves are rigid, i.e., the only point of their symplectic leaf.
More precisely, we say that a sheaf $M$ of homological dimension
$\le 1$ is {\em rigid} if it is isomorphic to any sheaf of homological
dimension $\le 1$ with the same numerical Chern class and (derived)
restriction to $C_\alpha$.  (This is stronger than the symplectic
leaf condition, since the symplectic leaf by definition includes
only simple sheaves.)

A slightly easier question is infinitesimal rigidity: at which sheaves is
the tangent space to the symplectic leaf $0$-dimensional?  Of course,
a sheaf is infinitesimally rigid iff the map
\[
\Ext^1(M,M\otimes \omega_X)\to \Ext^1(M,M)
\]
is $0$, or equivalently if the image is $0$-dimensional.  As usual, it will
be convenient to embed this in an Euler characteristic; we also extend to
pairs of sheaves.  Thus for two coherent sheaves $M$, $N$ on the Poisson surface
$(X,\alpha)$, we define
\[
\chi_\alpha(M,N)
:=
\dim(\Hom(M,N))
-
\dim(\im(\Ext^1(M,N\otimes \omega_X)\to \Ext^1(M,N)))
+
\dim(\Hom(N,M)).
\]
Note that duality immediately gives $\chi_\alpha(M,N)=\chi_\alpha(N,M)$;
it also gives
\[
\dim(\Hom(N,M))=\dim(\Ext^2(M,N\otimes \omega_X)),
\]
making $\chi_\alpha$ look more like an Euler characteristic.  (In fact, it
{\em is} an Euler characteristic: it is the Euler characteristic of the
tangent complex of the {\em derived} symplectic leaf.)

\begin{prop}\label{prop:symp_dims}
  The $\alpha$-twisted Euler characteristic $\chi_\alpha(M,N)$ depends only
  on the numerical Chern classes $c_*(M)$, $c_*(N)$ and the derived
  restrictions $M|^{\dL}_{C_\alpha}$ and $N|^{\dL}_{C_\alpha}$.
\end{prop}

\begin{proof}
We first note that we have a long exact sequence
\[
0\to C^{-1}
 \to A^0
 \to B^0
 \to C^0
 \to A^1
 \to B^1
 \to C^1
 \to A^2
 \to B^2
 \to C^2
 \to 0
\]
where
\begin{align}
A^p &:= \Ext^p(M,N\otimes\omega_X),\notag\\
B^p &:= \Ext^p(M,N),\notag\\
C^p &:= \Ext^p(M|^{\dL}_{C_\alpha},N|^{\dL}_{C_\alpha}).\notag
\end{align}
Indeed, since the complex $N\otimes\omega_X\to N$ represents
$N\otimes^{\dL} \sO_{C_\alpha}$, to construct such a sequence, it will
suffice to give a quasi-isomorphism
\[
\dR\Hom(M,N\otimes^{\dL} \sO_{C_\alpha})
\cong
\dR\Hom(M|^{\dL}_{C_\alpha},N|^{\dL}_{C_\alpha}).
\]
If $i:C_\alpha\to X$ is the natural closed embedding, then
\[
N\otimes^{\dL} \sO_{C_\alpha} \cong \dR i_*\dL i^* N,
\]
and thus
\[
\dR\Hom(M,N\otimes^{\dL} \sO_{C_\alpha})
\cong
\dR\Hom(M,\dR i_* \dL i^* N)
\cong
\dR\Hom(\dL i^* M,\dL i^* N),
\]
Since $\dL i^*M\cong M|^{\dL}_{C_\alpha}$, the claim follows.

Then, by definition,
\[
\chi_\alpha(M,N) = \dim(B^0) -\dim(\im(A^1\to B^1)) + \dim(A^2),
\]
and exactness gives
\[
\dim(\im(A^1\to B^1))
=
\dim(A^1)-\dim(C^0)+\dim(B^0)-\dim(A^0)+\dim(C^{-1}),
\]
so that
\[
\chi_\alpha(M,N)
=
\dim(A^0)-\dim(A^1)+\dim(A^2)
+
\dim(C^0)-\dim(C^{-1}).
\]
The $C$ dimensions manifestly only depend on $M|^{\dL}_{C_\alpha}$ and
$N|^{\dL}_{C_\alpha}$, while
\[
\dim(A^0)-\dim(A^1)+\dim(A^2)
=
\dim\Ext^2(N,M)-\dim\Ext^1(N,M)+\dim\Hom(N,M)
\]
is the usual $\Ext$ Euler characteristic, and thus only depends on the
numerical Chern classes of $M$ and $N$ (and can be computed via
Hirzebruch-Riemann-Roch).
\end{proof}

\begin{rem}
For convenient reference, we note that
\begin{align}
\dim\Ext^2(N,M)&{}-\dim\Ext^1(N,M)+\dim\Hom(N,M)\notag\\
{}={}&
-\chi(\sO_X)\rank(M)\rank(N)
+\chi(M)\rank(N)
+\rank(M)\chi(N)\notag\\
&-c_1(M)\cdot c_1(N)
+\rank(M) K_X\cdot c_1(N).
\end{align}
If $\Tor_1(M,\sO_{C_\alpha})=\Tor_1(N,\sO_{C_\alpha})=0$, then $\chi_\alpha(M,N)$
is given by the above together with the contribution
$\dim\Hom(M|_{C_\alpha},N|_{C_\alpha})$ from the restriction to the
anticanonical curve.  Particularly nice is the case $\rank(M)=\rank(N)=0$,
in which case
\[
\chi_\alpha(M,N)
=
-c_1(M)\cdot c_1(N)
+
\dim\Ext^0(M|^{\dL}_{C_\alpha},N|^{\dL}_{C_\alpha})
-
\dim\Ext^{-1}(M|^{\dL}_{C_\alpha},N|^{\dL}_{C_\alpha}).
\]
\end{rem}

In the case $M=N$, we call $\chi_\alpha(M,M)$ the {\em index of rigidity} of
$M$, following \cite{KatzNM:1996}; note
\[
\chi_\alpha(M,M)
=
2\dim(\End(M))-\dim(\im(\Ext^1(M,M\otimes \omega_X)\to \Ext^1(M,M))).
\]
In particular, if $M$ is simple, then it is infinitesimally rigid iff
$\chi_\alpha(M,M)=2$.

The $\alpha$-twisted Euler characteristic, and thus the index of rigidity,
behaves nicely under direct images and minimal lifts.  Actually, the middle
term (which is the one of greatest interest in any case) has the weakest
hypotheses to ensure good behavior.

\begin{prop}
  Let $\pi:X\to Y$ be a Poisson birational morphism of Poisson surfaces,
  and suppose $M$ and $N$ are coherent sheaves on $X$, at least one of
  which is $\pi_*$-acyclic with direct image of homological dimension $\le
  1$.  Then
\begin{align}
\dim(\im(\Ext^1(\pi_*M,\pi_*N\otimes \omega_Y)\to
  \Ext^1(\pi_*M,\pi_*N)))\notag\qquad &\\
\le 
\dim(\im(\Ext^1(M,N\otimes \omega_X)\to \Ext^1(M,N))).&
\end{align}
\end{prop}

\begin{proof}
  By symmetry, we may assume that $N$ is $\pi_*$-acyclic with direct image
  of homological dimension $\le 1$.  Then, using adjunction, we observe
  that the maps
\[
\Ext^1(\pi_*M,\pi_*N\otimes \omega_Y)
\to
\Ext^1(\pi_*M,\pi_*N)
\]
and
\[
\Ext^1(M,\pi^!\pi_*N\otimes \pi^*\omega_Y)
\to
\Ext^1(M,\pi^!\pi_*N)
\]
have isomorphic images.  By Lemma \ref{lem:useful_factorization}, the
latter factors through the map
\[
\Ext^1(M,N\otimes \omega_X)
\to
\Ext^1(M,N),
\]
giving us a bound on the image.
\end{proof}

\begin{prop}
  Let $\pi:X\to Y$ be a birational morphism of algebraic surfaces, and
  suppose $M$ and $N$ are coherent sheaves on $X$ such that $M$ is
  $\pi_*$-acyclic, $\pi_*N$ has homological dimension 1, and either $M$ has
  no invisible quotient or $N$ has no invisible subsheaf.  Then
  $\dim\Hom(\pi_*M,\pi_*N)\ge \dim(\Hom(M,N))$.
\end{prop}

\begin{proof}
It suffices to show that the natural map $\Hom(M,N)\to \Hom(\pi_*M,\pi_*N)$ 
is injective.  Suppose $f:M\to N$ is in the kernel.  Then it is 0 outside
the invisible locus, so has invisible image, giving both an invisible
quotient of $M$ and an invisible subsheaf of $N$.
\end{proof}

Combining hypotheses gives the following.

\begin{cor}
  Let $M$, $N$ be coherent $\pi_*$-acyclic sheaves on $X$ with direct
  images of homological dimension $\le 1$.  If (1) $M$ has no
  invisible quotient or $N$ has no invisible subsheaf, and (2) $N$ has
  no invisible quotient or $M$ has no invisible subsheaf, then
  $\chi_\alpha(M,N)\le \chi_\alpha(\pi_*M,\pi_*N)$.
\end{cor}

\begin{cor}
  Let $M$ be a coherent $\pi_*$-acyclic sheaf on $X$ with direct image of
  homological dimension $\le 1$.  Suppose further that either $M$ is
  simple, $M$ has no invisible quotient, or $M$ has no invisible
  subsheaf.  Then $\chi_\alpha(M,M)\le \chi_\alpha(\pi_*M,\pi_*M)$.
\end{cor}

\begin{proof}
If $M$ has no invisible quotient or no invisible subsheaf, this is just
a special case of the proposition.  Those hypotheses were only used to
prove the inequalities on $\Hom$ spaces, which are automatic when $M$ is
simple.
\end{proof}

\begin{prop}
  If $M$, $N$ are coherent sheaves on $Y$ of homological dimension $\le 1$,
  then the minimal lift respects the $\alpha$-twisted Euler characteristic:
  $\chi_\alpha(\pi^{*!}M,\pi^{*!}N)=\chi_\alpha(M,N)$.
\end{prop}

\begin{proof}
From the remark following Lemma \ref{lem:lift_exts}, we know that
\[
\Hom(\pi^{*!}M,\pi^{*!}N)\cong \Hom(M,N)
\]
and
\[
\Ext^1(\pi^{*!}M,\pi^{*!}N)\subset \Ext^1(M,N).
\]
Dually, the map
\[
\Ext^1(M,N\otimes \omega_Y)
\to
\Ext^1(\pi^{*!}M,\pi^{*!}N\otimes \omega_X)
\]
is surjective.  But then the maps
\[
\Ext^1(\pi^{*!}M,\pi^{*!}N\otimes \omega_X)
\to
\Ext^1(\pi^{*!}M,\pi^{*!}N)
\]
and
\[
\Ext^1(M,N\otimes \omega_Y)
\to
\Ext^1(M,N)
\]
have isomorphic images.
\end{proof}

\begin{prop}
Let $\pi:X\to Y$ be a birational morphism of Poisson surfaces, and $M$ a
coherent sheaf on $Y$ of homological dimension $\le 1$.  If $M'$ is any
pseudo-twist of $M$, then $\chi_\alpha(M,M)\le \chi_\alpha(M',M')$.
\end{prop}

\begin{proof}
  Since
  $\chi_\alpha(M,M)=\chi(\pi^{*!}M,\pi^{*!}M)=\chi_\alpha(\pi^{*!}M(\pm
  e_f),\pi^{*!}M(\pm e_f))$, it remains to show that the last index of
  rigidity is nonincreasing under direct images.  That this holds for the
  $\Ext^1$ term follows from the fact that $\pi^{*!}M(\pm e_f)$ is
  $\pi_*$-acyclic with direct image of homological dimension $\le 1$.  For
  the $2\dim(\End(M))$ term, consider the composition
\[
\End(M)\cong \End(\pi^{*!}M)\cong \End(\pi^{*!}M(\pm
e_f))\to \End(M').
\]
This is an isomorphism away from the point $\pi(f)$, and thus
any endomorphism in the kernel must vanish away from this point.  In
particular, any such endomorphism would have $0$-dimensional image, which
must be 0 since $M$ has homological dimension $\le 1$.
\end{proof}

\medskip

In the Hitchin system context, the main question is which pure
$1$-dimensional sheaves are rigid.  A first step is the following.

\begin{lem}
If $M$ is a simple 1-dimensional sheaf on the Poisson surface $Y$,
transverse to the anticanonical curve, then $M$ is rigid iff it is
infinitesimally rigid and $\Fitt_0(M)$ is
an integral curve.
\end{lem}

\begin{proof}
If $M$ has integral support and is infinitesimally rigid, then since
it is simple, we have $\chi_\alpha(M,M)=2$.  It follows that for any
other sheaf $N$ of homological dimension $\le 1$ with the same
numerical Chern class and restriction to $C_\alpha$, we have
\[
\chi_\alpha(M,N)=2.
\]
But this implies either $\Hom(M,N)>0$ or $\Hom(N,M)>0$.  Since $M$
is integral, any subsheaf is either 0 or has $0$-dimensional quotient.
Thus if $f:M\to N$ is a nonzero morphism, then either $f$ is
injective, or $f$ has $0$-dimensional image.  The latter cannot
happen, since $N$ has homological dimension $\le 1$, and the former
implies that $f$ is an isomorphism, by comparison of Chern classes.
Similarly, if $f:N\to M$ is a nonzero morphism, then the cokernel
is $0$-dimensional, but then comparison of Chern classes shows that
the kernel must also be $0$-dimensional, again a contradiction
unless $f$ is an isomorphism.  Either way, we obtain an isomorphism between $M$ and $N$, and thus $M$ is rigid.

Conversely, if $M$ is rigid, then it is certainly infinitesimally rigid
(since symplectic leaves in $\Spl_Y$ are smooth), so only the support needs
to be controlled.  Twisting by a line bundle has no effect on rigidity, so
we may assume that $M$ has no global sections.  If $\Fitt_0(M)$ is not
integral, then we may choose a nonzero subsheaf $N$ with strictly smaller
support, and a point $p$ of the divisor $\Fitt_0(M)-\Fitt_0(N)$ which is
not on $C_\alpha$.  For some sufficiently ample divisor $D$ on $Y$, the
twist $N'=N(D)$ will have a global section.  Define a
sequence of sheaves
\[
M(D)\cong M_{D\cdot c_1(M)}\supset M_{D\cdot
  c_1(M)-1}\supset\cdots\supset M_0,
\]
all containing $N'$, in the following way.  To
obtain $M_{k-1}$ from $M_k$, choose a nonzero map
$M_k/N'\to \sO_p$;
such a map exists since $p$ is in the support of $M_k/N'$.  Then $M_{k-1}$
is the kernel of the induced map $M_k\to \sO_p$.

This process has no effect on the first Chern class, and each step
reduces the Euler characteristic by 1, so that $M_0$ has the same
Chern class as $M$.  In addition, twisting by $\sO(D)$ as no
effect on the restriction to $C_\alpha$ (by transversality), and each
inclusion in the sequence is an isomorphism away from $p$.  So
rigidity implies that $M\cong M_0$.  But $M$ has no global sections,
while $M_0\supset N'$ does.
\end{proof}

\begin{rem}
Note that the argument that a map from $M$ to $N$ must be an
isomorphism also shows that any pure $1$-dimensional sheaf with
integral Fitting scheme is simple.  The argument from constancy of
the index of rigidity is essentially that of \cite[Thm.~1.1.2]{KatzNM:1996}.
\end{rem}

In many cases, the above argument ruling out non-integral support can be
interpreted as twisting by a suitable invertible sheaf on the support, and
noting that the resulting sheaf is nonisomorphic.  We could try to do
something similar in cases with integral support and positive genus; the
main technicality is that the action of invertible sheaves on torsion-free
sheaves has nontrivial stabilizers.  However, if we attempt to identify the
part of $\Ext^1(M,M)$ coming from such twisting, we are led to the
following result.

\begin{lem}
Let $M$ be a pure $1$-dimensional sheaf on $X$ transverse to the
anticanonical curve.  Then there is a natural injection
\[
H^1(\sHom(M,M))\to \im(\Ext^1(M,M\otimes\omega_X)\to \Ext^1(M,M)).
\]
\end{lem}

\begin{proof}
  The local-global spectral sequence for $\Ext^*(M,M)$ collapses at the
  $E_2$ page; indeed, this happens for $\Ext^*(M,N)$ as long as $M$ has
  homological dimension 1 and $N$ has $\le 1$-dimensional support.  We thus
  obtain the following short exact sequence
\[
0\to H^1(\sHom(M,M))\to \Ext^1(M,M)\to H^0(\sExt^1(M,M))\to 0.
\]
Similarly, the spectral sequence for $\Ext^*(M,M\otimes\omega_X)$
collapses, giving a corresponding exact sequence and a commutative diagram
\[
\begin{CD}
H^1(\sHom(M,M)\otimes\omega_X)@>>> \Ext^1(M,M\otimes\omega_X)\\
@VVV @VVV\\
H^1(\sHom(M,M))@>>>\Ext^1(M,M)
\end{CD}
\]
But
\[
H^1(\sHom(M,M)\otimes \sO_{C_\alpha})
=
0,
\]
by dimensionality, and thus the map
\[
H^1(\sHom(M,M)\otimes\omega_X)\to H^1(\sHom(M,M))
\]
is surjective.  The claim follows immediately.
\end{proof}

\begin{rem}
Note that the pairing between $\Ext^1(M,M\otimes\omega_X)$ and
$\Ext^1(M,M)$ restricts to the trivial pairing
\[
H^1(\sHom(M,M)\otimes\omega_X)\otimes H^1(\sHom(M,M))
\to
H^2(\sHom(M,M)\otimes\omega_X)
=
0.
\]
Since the overall pairing is perfect, it follows that
\[
\dim(\im(\Ext^1(M,M\otimes\omega_X)\to \Ext^1(M,M)))
\ge
2h^1(\sHom(M,M)).
\]
This inequality can be strict, e.g., if the first Fitting scheme of $M$ is
$0$-dimensional and not contained in $C_\alpha$.  (Indeed, in that case,
there is another sheaf with the same invariants and smaller endomorphism
sheaf.)
\end{rem}

\begin{thm}\label{thm:rigid_is_minus2}
  Let $M$ be a simple, pure $1$-dimensional sheaf on the Poisson surface
  $X$, transverse to the anticanonical curve.  Then $M$ is rigid iff there
  exists a Poisson birational morphism $\pi:Y\to X$ and a $-2$-curve $C$ on
  $Y$ disjoint from the anticanonical curve, such that $M\cong
  \pi_*\sO_C(d)$ for some integer $d$.
\end{thm}

\begin{proof}
We first show that sheaves of the form $\pi_*\sO_C(d)$ are rigid.
Since $C$ is transverse to the exceptional locus, $\sO_C(d)$ is
the minimal lift of its direct image.  Since $\sO_C(d)$ is certainly
simple, and we can compute
\[
\chi_\alpha(\sO_C(d),\sO_C(d)) = -C^2 = 2,
\]
we conclude that it is infinitesimally rigid, and thus so is its direct
image.  Since the image has integral Fitting scheme, the result follows.

  Now, suppose that $M$ is rigid, with (integral) support $C_0$, and let
  $\psi:\tilde{C}_0\to C_0$ be the normalization of $C_0$.  Then we observe
  that
\[
\sHom(M,M)\subset \psi_*\sO_{\tilde{C}_0},
\]
with quotient supported on the singular locus of $C_0$.  Indeed, if $M'$
denotes the quotient of $\psi^*M$ by its torsion subsheaf, then $M'$ is
torsion-free, so invertible, on the curve $\tilde{C}_0$.  This operation is
functorial, and an isomorphism on the smooth locus of $C_0$, so
\[
\sHom(M,M)\subset \psi_*\sHom(M',M')=\psi_*\sO_{\tilde{C}_0}.
\]
Now, $\dim\Hom(M,M)=h^0(\psi_*\sO_{\tilde{C}_0})=1$, and thus we find that
\[
h^1(\sHom(M,M))\ge h^1(\psi_*\sO_{\tilde{C}_0})=h^1(\sO_{\tilde{C}_0}),
\]
with equality only if $\sHom(M,M)=\psi_*\sO_{\tilde{C}_0}$.  Since $M$ is
infinitesimally rigid, the Lemma gives $h^1(\sHom(M,M))=0$, and thus
$\sHom(M,M)=\psi_*\sO_{\tilde{C}_0}$ and $h^1(\sO_{\tilde{C}_0})=0$.

In particular, $\tilde{C}_0$ is a smooth rational curve, and $M$ is the direct
image of a torsion-free sheaf on $\Spec(\sHom(M,M))\cong \tilde{C}_0$.
Since $M$ is the direct image of an invertible sheaf on a smooth curve,
there exists a Poisson birational morphism $\pi:Y\to X$ such that
$\pi^{*!}M$ is disjoint from $\pi^{*!}\sO_{C_\alpha}$.  But then
\[
\chi(\pi^{*!}M)=\chi_\alpha(\pi^{*!}M)=\chi_\alpha(M)=2,
\]
so that $c_1(\pi^{*!}M)^2=-2$, and disjointness gives $c_1(\pi^{*!}M)\cdot
K_Y=0$.  Since the Fitting scheme of $\pi^{*!}M$ is the image of
$\tilde{C}_0$ under a birational morphism, it follows that $\pi^{*!}M$ is
an invertible sheaf on a $-2$-curve as required.
\end{proof}

Note that if $X$ is not rational, then the only possible $-2$ curves on $X$
are components of fibers of the rational ruling of $X$.  Thus rigidity in
Hitchin-type systems is essentially only a phenomenon of the rational case,
which is an additional reason we will focus on the rational case below.
There is a weaker notion of rigidity that is perhaps more natural for
higher genus cases.  In higher genus, the first-order equations without
singularities are already non-rigid, forming a $2g$-dimensional family.
Since the tensor product of regular first-order equations is a regular
first-order equation, these equations actually form a group, and this group
acts on every symplectic leaf.  This suggests that one should ask about
{\em projectively} rigid equations, for which this group acts transitively
on the symplectic leaf.  (Since the group acts by scalar gauge
transformations, the resulting equations will then be unique up to scalar
gauge.)  It seems unlikely, however, that there are any nontrivial examples
of such equations (i.e., for which the equation is not first-order), so we
have not pursued this question any further.

\chapter{Divisors on rational surfaces}
\label{chap:divs_on_rat}

\section{Sheaves as equations revisited}

Using the above theory, we can give a further simplification of our
translation of the problem of classifying (relaxed) equations to sheaves:
apart from certain ill-behaved cases that can always be resolved by
isospectral transformations, there is a (Poisson) blowup of the original
ruled surface on which the minimal lift of the relaxed equation is disjoint
from the anticanonical curve.  We can generalize this somewhat further by
noting (following the proof of Theorem \ref{thm:poisson_class}) that we can
blow down any -1-curve disjoint from the sheaf without losing any
information, as long as we keep track of the marked section on the original
ruled surface.

More precisely, let $(X,\alpha)$ be a smooth, projective Poisson surface
(with $\alpha$ neither 0 nor an isomorphism), suppose $\rho:X\to C$ is a
rational ruling of $X$ (i.e., a map to a smooth projective curve such that
the generic fiber is a smooth curve of genus 0), and let $s:C\to X$ be a
section of $\rho$ which is transverse to the anticanonical curve
$C_\alpha$.  If we repeatedly blow up points where the strict transform of
$s$ meets the anticanonical curve, we will eventually reach a surface
$\tilde{X}$ in which the strict transform of $s$ is disjoint from the
anticanonical curve, at which point we may repeatedly blow down $-1$-curves
disjoint from $s$ and eventually reach a standard Poisson surface, with $s$
the marked section disjoint from $C_\alpha$.  Thus given a pure
$1$-dimensional sheaf on $X$, the direct image of its minimal lift to
$\tilde{X}$ will have a natural interpretation as a relaxed equation.

This interpretation of course depends on all of the data, but it turns out
to depend only mildly on the choice of section.  To understand this
dependence, it suffices to consider the case that $X$ is already a standard
Poisson surface, but the chosen section is not the marked section.  It is
easy to see that the two anticanonical curves are naturally isomorphic
(they are certainly birational, and agree along the fibers being modified),
so that this will produce a new relaxed equation on the same curve.  The
fact that we should be taking the minimal lift to $\tilde{X}$ is a
complicating factor, but we can resolve this by taking a suitable
non-minimal factorization.  This extends by a quotient sheaf with trivial
direct image to the base curve; since this property survives the direct
image of minimal lift procedure, starting with a non-minimal factorization
has no effect on the final relaxed equation.  In particular, we may feel
free to compose $B$ with the natural map $W\to W(\rho^{-1}D)$ where the
divisor $D$ on $C$ is the image of the divisor on $C_\alpha$ induced by the
chosen section.  This makes $B$ vanish at the points we are blowing up,
making it straightforward to compute a resolution of the minimal lift.  The
resulting resolution is actually the inverse image of a resolution on the
new standard surface, and thus the direct image is again straightforward to
compute.  We find that the new matrix $B'$ on $C_\alpha$ is nothing other
than a scalar multiple of the original matrix $B$.  More precisely, the
sheaf $\sO_{s_X}$ maps to the structure sheaf of a section on $X'$ which
thus corresponds to a first-order equation, and $B'$ is the multiple of $B$
by the corresponding scalar.

In other words, changing the section changes the relaxed equation by a
scalar gauge transformation (the relaxed analogue of multiplying the
solution vector $v$ by the solution of a first-order equation), and thus
without the section, we still have a well-defined ``projective'' relaxed
equation.  (I.e., if we think of the equations as having structure group
$\GL_n$, we can at least recover the image in $\PGL_n$.)  Equivalently, we
may think of the data without section as determining a {\em twisted}
equation, and the choice of section as specifying a first-order twisted
equation to gauge by in order to undo the twist.  This suggests a slight
generalization: rather than specify the section itself, we should specify a
line bundle on the section.  Of course, any such line bundle corresponds to
a line bundle on the base curve, so we can undo this by twisting our sheaf
by the inverse line bundle.

The choice of rational ruling, on the other hand, is either completely
trivial or extremely significant, depending in large part on the birational
equivalence class of surfaces we are considering.  For a higher genus
surface (i.e., $H^1(\sO_X)\ne 0$), there {\em is} no choice of rational
ruling: since the Albanese torsor of $\P^1$ is a point, the natural map
from $X$ to its Albanese torsor contracts the fibers of any rational
ruling, and thus the map from $X$ to its image in $\Alb^1(X)$ is a
canonical rational ruling on $X$ with any other rational ruling equivalent
to the canonical ruling.

In contrast, in the rational case, a given rational surface can have
multiple (or even infinitely many!) rational rulings.  The simplest
instance of this is $\P^1\times \P^1$, for which the two maps to $\P^1$
give distinct rational rulings.  Given an anticanonical curve on
$\P^1\times \P^1$, each ruling gives its own interpretation of sheaves up
to scalar gauge.  Since the order of equation corresponding to a sheaf may
be computed by the intersection of its Chern class with the class of a
fiber of the ruling, we see that the two interpretations of a sheaf will in
general have different orders.  In fact, they need not even be the same
{\em kind} of equation!  In $\P^1\times \P^1$, consider the anticanonical
curve $x_1y_1y_2^2=0$.  Relative to the first ruling, the horizontal part
of the anticanonical curve is nonreduced, and thus after resolving the
choice of normalizing section, the anticanonical curve will remain
nonreduced, and thus the relaxed equation will be a Higgs bundle.  On the
other hand, relative to the {\em second} ruling, we have two distinct
horizontal components which will be tangent after the resolution (there is
only one possible point of intersection, but the intersection multiplicity
must be 2), and thus {\em that} relaxed equation is of ordinary difference
type.  In all, there are 16 different combinatorial types of anticanonical
curve in $\P^1\times \P^1$, each of which gives a transformation on
relaxed equations.  In the non-relaxed case, these can be thought of as
formal integral transforms (with the above example corresponding to the
Mellin or $z$-transform), which are given in Appendix \ref{chap:fourier}
below.

One reason this is significant is that the moduli space itself is
independent of the choice of ruling (and in principle on the choice of
anticanonical curve, though generally we will have blown up enough points
to make that unique), and thus different choices of ruling will give the
same moduli spaces (i.e., the same spaces of initial conditions for
isospectral problems).  Thus in the rational case we expect any given
isospectral (or isomonodromic) integrable system to have many different
such interpretations, and we will want to understand how to determine the
range of possibilities.

\section{Blowdown structures on rational surfaces}
\label{sec:blowdown_structures_comm}

With this in mind, we want to understand the set of possible ways to blow a
rational surface down to a Hirzebruch surface.  It turns out to be useful
to record slightly more information than just a birational morphism to a
Hirzebruch surface, and we thus make the following definition.

\begin{defn}
Let $X$ be a rational surface with $X\not\cong \P^2$.  A (Hirzebruch)
blowdown structure on $X$ is a chain $\Gamma$ of morphisms
\[
X=X_m\to X_{m-1}\to\cdots\to X_0\to \P^1,
\]
such that for $1\le i\le m$, the morphism $X_i\to X_{i-1}$ is the blowup in a
single point of $X_{i-1}$, while the morphism from $X_0\to \P^1$ is a geometric
ruling.  Two blowdown structures will be considered equivalent if they fit
into a commutative diagram
\[
\begin{CD}
X@>>> X_{m-1}@>>>\cdots @>>> X_0@>>> \P^1\\
 @|  @VVV \cdots @. @VVV @VVV\\
X@>>> X'_{m-1}@>>>\cdots @>>> X'_0@>>> \P^1
\end{CD}
\]
\end{defn}

\begin{rems}
  Related structures have been studied in the case that the rational
  surface can be blown down to $\P^2$, \cite{SakaiH:2001}; one also
  considers the related notion of an exceptional configuration (essentially
  the analogue for $\P^2$ of the notion of numerical blowdown structure
  below) \cite{LooijengaE:1981}.  Our considerations below are somewhat
  more general (since not every surface blows down to $\P^2$, especially
  including $F_2$), but of course closely related; for instance
  \cite{LooijengaE:1981} already saw the appearance of the root system
  $E_{m+1}$.
\end{rems}

\begin{rems}
  One can make a similar definition for higher genus rationally ruled
  surfaces, but there is much less of interest to say in that case, so we
  will defer that discussion to the noncommutative case.
\end{rems}

Note that in addition to keeping track of a factorization of the birational
morphism, we also keep track of the ruling at the end (but only up to
$\PGL_2$).  Most of the time, of course, the latter provides no
information, since most Hirzebruch surfaces have a unique geometric ruling
(and thus have a unique blowdown structure), with the lone exception being
$\P^1\times \P^1$.

One reason for including the above information in the blowdown structure is
that it allows us to associate a basis of $\Pic(X)$ to any blowdown
structure.  For each monoidal transformation $X_i\to X_{i-1}$, we have an
exceptional curve $e_i$ on $X_i$, and the total transform of this curve
gives us a divisor on $X$, which we also denote $e_i$.  Since
\[
\Pic(X)\cong \Pic(X_0)\oplus \bigoplus_i \Z e_i,
\]
it remains only to give a basis of $\Pic(X_0)\cong \Z^2$.  One basis
element is obvious, namely the class $f$ of the fibers of the ruling.  For
the other, there is also an obvious choice, namely the class $s_{\min}$ of
a section with minimal self-intersection.  This turns out not to be the
best choice for our purposes, however, as it gives us a countable infinity
of different intersection forms to consider.  A slightly different basis
greatly reduces the number of cases.  Define a divisor class
\[
s:=s_{\min} + \lfloor -s_{\min}^2/2\rfloor f.
\]
Since $s_{\min}\cdot f=1$, $f^2=0$, we find that $s\cdot f=1$, and $s^2$ is
either 0 or $-1$, depending on whether $s_{\min}^2$ was even or odd.  With
this in mind, we call a blowdown structure even or odd depending on the
parity of $s_{\min}^2$.  (Note that if $X_0$ comes from a line bundle on a
hyperelliptic genus 1 curve as in the case of twisted elliptic difference
equations, then this choice of basis element essentially corresponds to
choosing the bundle to have degree $1$ or $2$, as then $\sO_X(s)$ agrees
with the relative $\sO(1)$.)

The basis we obtain has one of two possible intersection forms, depending
on parity.  In the even case, we have
\[
s^2 = 0,\quad s\cdot f=1,\quad f^2=0,
\quad s\cdot e_i=f\cdot e_i=0,\quad
e_i\cdot e_j=-\delta_{ij},
\]
while in the odd case, we have the same, except $s^2=-1$.
The expansion of the canonical class in the basis again only depends on
parity:
\[
K_X = \begin{cases}
-2s-2f+\sum_{1\le i\le m} e_i & \text{$\Gamma$ even}\\
-2s-3f+\sum_{1\le i\le m} e_i & \text{$\Gamma$ odd}.
\end{cases}
\]
And of course in either case we find $K_X^2 = 8-m$.  When $X_0=F_1$, we can
blow it down to $\P^2$, suggesting an alternate basis in the odd case:
replace $f$ by $h=s+f$, the class of a line in $\P^2$.  This gives an
orthonormal basis (up to signs) for $\Pic(X)$, with $K_X=-3h+s+\sum_i e_i$,
but makes the effective cone look rather strange when $X_0=F_{2d+1}$ for
$d>0$.  When working in this basis, it is natural to denote $s$ by $e_0$
for symmetry purposes.

We order the basis by $s,f,e_1,\dots,e_m$ or $h,e_0,e_1,\dots,e_m$ in order
for the following to hold.

\begin{lem}
  Let $X$ be a rational surface with blowdown structure $\Gamma$.  Then the
  expansion of any effective divisor in the corresponding basis is
  lexicographically nonnegative.
\end{lem}

\begin{proof}
  Let $D$ be an effective divisor on $X$, which we may clearly assume
  nonzero.  The coefficient of $s$ in $D$ is $f\cdot D$, which is
  nonnegative since $f$ has a representative which is integral of
  nonnegative self-intersection, so is nef.  If the coefficient of $s$ is
  0, then let $i$ be the largest index such that the image of $D$ in $X_i$
  is still a divisor.  That image must be a multiple of $f$ or $e_i$ as
  appropriate, and since $f$ and $e_i$ are effective, that multiple must be
  positive.
\end{proof}

Note that we can recover the blowdown structure from the corresponding
basis for $\Pic(X)$: blow down $e_m$, then the image of $e_{m-1}$, etc.,
and construct a map $X_0\to \P^1$ using $f$.  Of course, not every basis
with the correct intersection form will correspond to a blowdown structure,
but we will eventually give an algorithm for determining when a given basis
(expressed in terms of some original blowdown structure) also corresponds
to a blowdown structure.  In any event, we will define a {\em numerical}
blowdown structure to be a basis of $\Pic(X)$ having the same intersection
form as an even or odd blowdown structure on $X$.

The surface $X_1$ was obtained by blowing up a point of $X_0$, and we find
that the fiber containing that point becomes a pair of $-1$ curves on
$X_1$, of divisor classes $e_1$, $f-e_1$.  We thus obtain an alternate
blowdown structure on $X_1$ by blowing down $f-e_1$, producing $X'_0$
differing from $X_0$ by an elementary transformation.  The basis elements
$f$ and $e_i$ for $i\ge 2$ are unchanged by this transformation, but $e_1$
and $s$ are transformed as follows.
\[
(s',e'_1)
=
\begin{cases}
(s-e_1,f-e_1) & \text{$\Gamma$ even}\\
(s+f-e_1,f-e_1) & \text{$\Gamma$ odd}
\end{cases}
\]
Note that this swaps the even and odd cases, and if we perform the
transformation twice, we end up back at the original blowdown structure.

Another natural way to transform a blowdown structure is to rearrange
blowups.  If the morphism $X_{i+2}\to X_i$ blows up two distinct points of
$X_i$, then we can perform the blowdown in the other order, thus swapping
the basis elements $e_i$ and $e_{i+1}$.  Unlike the elementary
transformation case, this operation is {\em not} always legal, as when
$X_{i+2}\to X_{i+1}$ blows up a point of $e_i$, there is no longer any
choice in how to reach $X_i$.  However, when it applies, it has a
particularly nice action on the basis: it is simply the reflection with
respect to the intersection form in the divisor class $e_i-e_{i+1}$.
Similarly, when $X_0\cong \P^1\times \P^1$, we obtain another blowdown
structure by changing to the other ruling.  This swaps the basis elements
$s$ and $f$, and is the reflection in the divisor class $s-f$.

In this way, we obtain a collection of $m$ reflections in the even case,
$m-1$ in the odd case; if we perform an elementary transformation, the two
sets mostly overlap, but we obtain a total of $m+1$ different reflections
in this way (assuming $m\ge 2$).  The corresponding vectors are linearly
independent, and are given by
\[
s-f,f-e_1-e_2,e_1-e_2,\dots e_{m-1}-e_m
\]
in the even case and
\[
s-e_1,f-e_1-e_2,e_1-e_2,\dots,e_{m-1}-e_m
\]
in the odd case.  Note that each one of these vectors is orthogonal to $K$;
we can see this either by direct calculation or by noting that the
expansion of $K$ depends only on the parity of the blowdown structure, so
had better be invariant under the above reflections.

\begin{lem}
Suppose $m\ge 2$.  Then the above sets of vectors give a basis for the
orthogonal complement of $K$ in $\Pic(X)$.  With respect to the negative of
the intersection form, they form the set of simple roots of a Coxeter group
of type $E_{m+1}$.
\end{lem}

\begin{proof}
  In either case, we have $m+1$ vectors, while $\Pic(X)$ has rank $m+2$, and
  thus we obtain bases of the orthogonal complement over $\Q$.  Since the
  bases are obviously saturated (they are essentially triangular with unit
  diagonal), the first claim follows.

  That the vectors are simple roots for a Coxeter system follows from the
  fact that their inner products are nonpositive (i.e., the intersections
  are nonnegative).  To identify the system, note that the $e_i-e_{i+1}$
  roots are the simple roots of a Coxeter group of type $A_{m-1}$,
  adjoining $f-e_1-e_2$ extends this to $D_m$, and adjoining $s-f$ or
  $s-e_1$ as appropriate extends one of the short legs of the $D_m$ Dynkin
  diagram, giving a diagram of type $E_{m+1}$ (i.e., star-shaped diagram
  with legs of lengths $2$, $3$, and $m-2$).
\end{proof}

\begin{rem}
Note the small $m$ cases
\begin{align}
\notag E_3&=A_1\times A_2\\
\notag E_4&=A_4\\
\notag E_5&=D_5\\
\notag E_9&=\tilde{E}_8,
\end{align}
with $E_6$, $E_7$, $E_8$ as expected.
When $m=1$, we have only the root $s-f$ or $s-e_1$ as appropriate, and when
$m=0$, we have only $s-f$ in the even case, and no roots in the odd case.
\end{rem}

With this in mind, we refer to the given vectors as the simple roots for
the (numerical) blowdown structure.  The corresponding simple reflections
clearly give an action of $W(E_{m+1})$ on the set of numerical blowdown
structures.

\begin{lem}
Suppose $\Gamma$ is a blowdown structure for the rational surface $X$, and
let $\sigma$ be a simple root for $\Gamma$, with corresponding reflection
$r_\sigma$.  If $\sigma$ is ineffective, then $r_\sigma\Gamma$ is a
blowdown structure.
\end{lem}

\begin{proof}
  To be precise, we mean here that if the numerical blowdown structure
  $\Gamma$ comes from a blowdown structure, then so does $r_\sigma\Gamma$,
  as long as the divisor class $\sigma$ is ineffective.

  Using elementary transformations as appropriate, we may reduce to the
  cases $\sigma=s-f$ and $\sigma=e_i-e_{i+1}$.  If $s-f$ is ineffective on
  $X$, it is certainly ineffective on $X_0$, but then $X_0$ must be
  $\P^1\times \P^1$ (if $X_0\cong F_{2d}$ for $d>0$, then $s_{\min}=s-df$
  is effective), and we have already seen that the reflection gives a
  blowdown structure.  Similarly, if $e_i-e_{i+1}$ is ineffective on $X$,
  it is ineffective on $X_{i+1}$, which implies that $X_{i+1}\to X_{i-1}$
  blows up two distinct points of $X_i$, so that we need merely blow up the
  points in the opposite order.
\end{proof}

\begin{rems}
  The converse is also true: $\sigma$ is the negative of a simple root of
  $r_\sigma\Gamma$, and thus cannot be effective.
\end{rems}

\begin{rems}
The roots of the $A_{m-1}$ subsystem act without changing the Hirzebruch
surface $X_0$, by permuting the distinct points being blown up.  Similarly,
the simple reflections of the $D_m$ subsystem leave the rational ruling
invariant.  If we combine those reflections with the action of the
elementary transformation, we obtain a group of type $C_m$ acting on the
different ways to blow down to a Hirzebruch surface compatibly with the
given ruling.
\end{rems}

Given a blowdown structure $\Gamma$, call an element $w\in W(E_{m+1})$ {\em
  ineffective} if there exists a word $w = r_1 r_2\cdots r_l$ with each
$r_i$ a simple reflection such that the corresponding simple root is
ineffective for the relevant blowdown structure, $r_{i+1}r_{i+2}\cdots
r_l\Gamma$.  (That this numerical blowdown structure comes from an actual
blowdown structure follows by an easy induction.)  In particular, if $w$ is
ineffective, then $w\Gamma$ is a blowdown structure.

We thus need to understand the effective simple roots.  By a ``$-d$-curve''
on a rational surface, we mean a smooth rational curve of self-intersection
$-d$.

\begin{lem}\label{lem:fixed_of_simple}
  Let $X$ be a rational surface with blowdown structure $\Gamma$ and
  $K_X^2<8$.  Then any effective simple root $\sigma$ can be decomposed as
  a nonnegative linear combination of $-d$-curves with $d>0$.  There is at
  most one fixed component of $\sigma$ not orthogonal to $\sigma$, and that
  component has self-intersection $\le -2$, with equality only if $\sigma$
  is a $-2$-curve.
\end{lem}

\begin{proof}
  If $\sigma=e_i-e_{i+1}$ is effective, then it is the total transform of a $-2$
  curve on $X_{i+1}$.  At each later step in the blowing up process, either
  we blow up a point not on the total transform, in which case the
  decomposition is unchanged, or we blow up a point on the total
  transform, in which case we acquire an additional component $e_j$, and
  the component(s) containing the center of the monoidal transform have
  their self-intersection decreased by $1$.  Thus by induction every
  component is a $-d$-curve for some $d>0$.  We also see that $\sigma$ is
  uniquely effective, so every component is fixed, but, except for the strict
  transform of the original $-2$ curve, orthogonal to $\sigma$.  The case
  $\sigma=f-e_1-e_2$ follows by elementary transformation.
  
  For the remaining case, assume for convenience that $\Gamma$ is even, so
  the remaining root is $s-f$.  This is effective precisely when $X_0\cong
  F_{2d}$ with $d>0$, when we can write it in the form
\[
s-f = s_{\min} + (d-1)f.
\]
This same decomposition applies to $X_m$, so that the only fixed components
are those of the total transform of $s_{\min}$, to which the previous
calculation applies.  On the other hand, we could choose the fibers in this
decomposition to all pass through the point blown up on $X_1$, obtaining a
decomposition
\[
s-f = s_{\min} + (d-1)(f-e_1)+(d-1)e_1,
\]
or
\[
s-f = (s_{\min}-e_1) + (d-1)(f-e_1)+de_1
\]
on $X_1$, the latter when the point being blown up is on $s_{\min}$.  The
components of this decomposition are all $-e$-curves for varying $e>0$, and
as before this property is preserved on taking the total transform to $X_m=X$.
\end{proof}

Since reflections in ineffective simple roots take blowdown structures to
blowdown structures, we can define a groupoid (the {\em strict} groupoid of
blowdown structures on $X$) as follows: the objects are the blowdown
structures on $X$, while the morphisms are given by the actions of
ineffective elements of $W(E_{m+1})$.

\begin{thm}\label{thm:two_orbits_of_blowdowns_comm}
If $X$ is a rational surface with $K_X^2<8$, then the strict groupoid of
blowdown structures on $X$ has precisely two isomorphism classes, one for
each parity of blowdown structure.  (If $K_X^2=8$, the groupoid has only
one isomorphism class.)
\end{thm}

In other words, any two blowdown structures on $X$ with the same parity are
related by a sequence of reflections in ineffective simple roots.  For
$K_X^2<8$, elementary transformations imply that we need only consider the
even parity case.

\begin{proof}
  If $K_X^2=8$, this is obvious, as either there is only one blowdown
  structure or $X\cong \P^1\times \P^1$, and there are two blowdown
  structures related by reflection in $s-f$.  Now, suppose $K_X^2=7$, and
  fix an even blowdown structure on $X$.  To blow $X$ down to a Hirzebruch
  surface, we must blow down a $-1$ curve, which in particular gives a
  divisor class $D$ such that $D^2=D\cdot K_X=-1$.  There are only three
  such divisor classes, namely $D\in \{e_1,s-e_1,f-e_1\}$.  Only the case
  $D=e_1$ could blow down to an even Hirzebruch surface, since the other
  two cases have classes of odd self-intersection in their orthogonal
  complements.  In other words, a surface with $K_X^2=7$ blows down to a
  unique even Hirzebruch surface, and thus the even blowdown structures on
  $X$ are bijective with the even blowdown structures on this Hirzebruch
  surface.  (Note also that $f-e_1$ is always a $-1$ curve, while $s-e_1$,
  related to it by a simple reflection, is either a $-1$ curve or
  decomposes as $s-e_1=(s-f)+(f-e_1)$, depending on whether $s-f$ is
  effective.)

  For $K_X^2<7$, we may induct on $m$, and thus reduce to the question of
  showing that every $-1$ curve on $X$ can be moved to $e_m$ by a sequence
  of reflections in ineffective simple roots.  Again, we may assume we are
  starting from an even blowdown structure, conjugating by elementary
  transformations as appropriate.  Let
\[
E = n s + d f - \sum_i r_i e_i
\]
be the class of the given $-1$ curve, and recall that $E$, being effective,
must be lexicographically positive, so in particular $n\ge 0$.

If $E\cdot \sigma<0$ for some simple root, then the root must be ineffective,
since otherwise $E$ would be a fixed component of self-intersection $-1$
not orthogonal to the simple root.  In particular, we can always perform
the corresponding reflection, and this makes the vector
$(n,d,-r_1,\dots,-r_{m-1})$ lexicographically smaller.  Since the
reflections preserving $n$ form a finite group (of type $D_m$), we conclude
that after finitely many reflections in ineffective simple roots, we will
obtain a divisor such that $E\cdot \sigma\ge 0$ for every simple root $\sigma$.
We may also assume $E\cdot e_m\ge 0$, since otherwise $E=e_m$.

For such a divisor, we have the inequalities
\[
d\ge n\ge r_1+r_2;\quad r_1\ge r_2\ge\cdots\ge r_m\ge 0.
\]
Since $(E-e_m)\cdot K_X = 0$, we may also express $E-e_m$ as a linear
combination of simple roots:
\[
E-e_m = n(s-f) + (n+d)(f-e_1-e_2)
+(n+d-r_1)(e_1-e_2)+\!
\sum_{2\le k\le m-1}\! (1+\sum_{k<i} r_i)(e_i-e_{i+1}),
\]
clearly a nonnegative linear combination.  Since
\[
E\cdot (E-e_m) = -1 - r_m <0,
\]
we obtain a contradiction.
\end{proof}

\begin{rems}
We can adapt this to an algorithm for testing whether a given class is a
$-1$-curve (and thus whether a given numerical blowdown structure comes
from an actual blowdown structure): reflect in simple roots with $E\cdot
\sigma<0$ until one of the roots is effective, $E\cdot f<0$, or $E=e_m$.
Then $E$ is a $-1$-curve iff the last termination condition holds.  This is
the first of many variants we will see of the standard
reduce-to-the-fundamental-chamber algorithm from Coxeter theory.
\end{rems}

\begin{rems}
Of course, we could obtain a groupoid with a single isomorphism class by
including morphisms of the form $w\epsilon$ where $\epsilon$ is the
elementary transformation, but this is somewhat inconvenient, since the
morphisms no longer correspond directly to elements of a group.
\end{rems}

\begin{rems}
  An analogous statement holds for higher genus ruled surfaces, as a
  special case of the even more general result proved for noncommutative
  surfaces in Proposition \ref{prop:formal_-1_can_be_blown_down} below.
  (That Proposition implicitly assumes a choice of anticanonical curve, but
  this is not needed in the commutative case.)
\end{rems}

Any $-2$-curve on $X$ is a (real) root of the root system $E_{m+1}$.  More
precisely, we have the following.

\begin{prop}
  Suppose $D$ is the class of a $-2$ curve.  Then there exists a blowdown
  structure on $X$ for which $D$ is a simple root.
\end{prop}

\begin{proof}
Let $\Gamma$ be an even blowdown structure on $X$, and write
\[
D = ns + df - \sum_{1\le i\le m} r_i e_i
\]
as before.  Since $D\cdot K_X=0$, we can expand $D$ as a linear combination
of simple roots
\[
D = n(s-f) + (n+d)(f-e_1-e_2)
+(n+d-r_1)(e_1-e_2)+\!
\sum_{2\le k\le m-1}\! (\sum_{k<i} r_i)(e_i-e_{i+1}),
\]
and thus find as before that $D\cdot \sigma<0$ for some simple root $\sigma$.
Once more, $D$ would have to be a fixed component of $\sigma$ if $\sigma$ were
effective, and thus either $D=\sigma$ or $\sigma$ is ineffective.  As before, in the
latter case, reflecting makes $D$ lexicographically smaller, so this
process must terminate.
\end{proof}

\begin{rem}
  Again, this translates to an algorithm for testing whether a given
  divisor class is represented by a $-2$-curve, which is formally very
  similar to the algorithm of \cite{KatzNM:1996} for testing whether a
  local system is rigid.
\end{rem}

In the case of a surface with a chosen anticanonical curve, there is a
related groupoid with more morphisms and nontrivial stabilizers.  Call a
simple reflection $s\in S(E_{m+1})$ {\em admissible} for the blowdown
structure $\Gamma$ (and the anticanonical curve $C_\alpha$) if the
corresponding simple root is either ineffective or has intersection 0 with
every component of $C_\alpha$.  (This implies that it is a $-2$-curve
disjoint from $C_\alpha$; if it were reducible, some component would have
negative intersection with $C_\alpha$.)  Although $s\Gamma$ is no longer a
blowdown structure when $s$ is effective but admissible, we define a
modified action as follows.  If $s$ is ineffective, then $s\cdot
\Gamma:=s\Gamma$, while if $s$ is effective but admissible, then $s\cdot
\Gamma:=\Gamma$.

More generally, call an element $w\in W(E_{m+1})$ admissible for $\Gamma$
if there exists a word
\[
w = s_1 s_2 \cdots s_l
\]
for $w$ such that for each $1\le i\le l$, $s_i$ is admissible for
\[
s_{i+1}\cdot s_{i+2}\cdots s_{l-1}\cdot s_l\cdot \Gamma.
\]
In this case, we also call the given word admissible.

\begin{prop}
If $w$ is admissible for $\Gamma$, then every reduced word for $w$ is
admissible for $\Gamma$.  Moreover, if
\[
w = s_1\cdots s_l = s'_1\cdots s'_{l'}
\]
are two admissible words representing $w$, then
\[
s_1\cdot s_2\cdots s_{l-1}\cdot s_l\cdot \Gamma
=
s'_1\cdot s'_2\cdots s'_{l'-1}\cdot s'_{l'}\cdot \Gamma.
\]
\end{prop}

\begin{proof}
  We first note that if $s$ is admissible for $\Gamma$, then it is also
  admissible for $s\cdot \Gamma$; either $s$ is ineffective and remains so,
  or $s$ is a $-2$ curve, and $s\cdot\Gamma=\Gamma$.  Either way, we find
  $s\cdot s\cdot\Gamma = \Gamma$.  In other words, if a reflection occurs twice
  in a row in an admissible word, we can remove the pair without affecting
  admissibility or the final blowdown structure.

  Since any word can be transformed into a reduced word by a sequence of
  braid relations and removal of repeated reflections, and any two reduced
  words are related by a sequence of braid relations, it remains only to
  show that the claim holds for braid relations.  In other words, given a
  braid relation in $W(E_{m+1})$, we need to show that either both sides
  are inadmissible or both sides are admissible and produce the same
  blowdown structure.

  Let $s$, $t$ be distinct simple reflections.  If both are inadmissible,
  there is nothing to prove, so suppose that $s$ is admissible for
  $\Gamma$.  Then we observe that if $t$ is inadmissible for $\Gamma$, then
  it is also inadmissible for $s\cdot\Gamma$.  Thus only the case that $s$
  and $t$ are both admissible need be considered.  If the braid relation is
  $st=ts$, then we need merely check that the relation holds in each of the
  four cases ($s$ effective or not, $t$ effective or not).

  Thus suppose the braid relation is $sts=tst$.  Let $r_s$, $r_t$ be the
  corresponding simple roots.  If $r_s+r_t$ is ineffective, then
  either $r_s$, $r_t$ are both ineffective (and the braid relation follows
  from the fact that the action agrees with the linear action) or precisely
  one (say $r_t$) is effective.  But then we find that $r_s$ is effective
  in the blowdown structure $t\cdot s\cdot \Gamma$, and thus
\[
s\cdot t\cdot s\cdot \Gamma = t\cdot s\cdot \Gamma = t\cdot s\cdot t\cdot\Gamma.
\]
If $r_s+r_t$ is effective and $r_t$ is effective, then $r_s$ is also
effective.  Indeed, $(r_s+r_t)\cdot r_t=-1$, and thus any representative of
$r_s+r_t$ contains $r_t$ as a component, implying $r_s+r_t-r_t$ effective.
Thus in this case, we have
\[
s\cdot t\cdot s\cdot \Gamma = t\cdot s\cdot t\cdot\Gamma = \Gamma,
\]
since the blowdown structure never changes.

  Finally, we have the case $r_s+r_t$ effective but $r_s$, $r_t$ are
  ineffective.  Relative to the blowdown structure $t\cdot \Gamma$, $r_s$ is
  effective, and thus
\[
s\cdot t\cdot s\cdot t\cdot\Gamma
=
t\cdot s\cdot t\cdot t\cdot\Gamma
=
t\cdot s\cdot \Gamma
\]
(or both sides are undefined, if $r_s$ is inadmissible for $t\cdot \Gamma$).
If both sides are defined, then $s$ is admissible for both blowdown
structures, and thus
\[
t\cdot s\cdot t\cdot \Gamma
=
s\cdot s\cdot t\cdot s\cdot t\cdot\Gamma
=
s\cdot t\cdot s\cdot t\cdot t\cdot\Gamma
=
s\cdot t\cdot s\cdot \Gamma,
\]
and we are done.
\end{proof}

If we use admissible elements in place of effective elements in defining
the groupoid of blowdown structures, the resulting ``weak'' groupoid has
nontrivial stabilizers, a conjugacy class of reflection subgroups of
$W(E_{m+1})$.

\begin{prop}
  The stabilizer of $\Gamma$ in the weak groupoid of blowdown structures is
  the reflection subgroup of $W(E_{m+1})$ generated by reflections in
  $-2$-curves disjoint from the anticanonical curve.
\end{prop}

\begin{proof}
Given a $-2$-curve $v$ disjoint from $C_\alpha$, let $w$ be an effective element
of $W(E_{m+1})$ such that $v$ is a simple root $\sigma$ in $w\Gamma$.  Then 
$r_\sigma$ stabilizes $w\Gamma$, so $r_v=w^{-1} r_\sigma w$ stabilizes
$\Gamma$.

Conversely, consider an admissible reduced word $w$ stabilizing $\Gamma$.
If every reflection in $w$ is ineffective, then $w$ acts linearly, and
since it stabilizes a basis, we have $w=1$.  Otherwise, we can write
\[
w = w_1 r w_2
\]
where $r$ is an effective but admissible simple reflection, $w_2$ is
ineffective, and $\ell(w)=\ell(w_1)+\ell(w_2)+1$.  Since $r$ is effective,
it stabilizes $w_2\cdot \Gamma=w_2\Gamma$.  We thus conclude that we can factor
\[
w = w_1 w_2 (w_2^{-1}r w_2)
\]
where both $w_1w_2$ and $w_2^{-1}r w_2$ are admissible elements stabilizing
$\Gamma$.  The second factor is a reflection in the simple root
corresponding to $r$ in $w_2\Gamma$ (which is admissible, so a $-2$-curve
disjoint from $C_\alpha$), while the first factor has length strictly
smaller than $\ell(w)$.  Thus, by induction, $w$ can be written as a product
of reflections in $-2$-curves disjoint from $C_\alpha$.
\end{proof}

Note from \cite{HarbourneB:1997} that the $-2$-curves disjoint from
$C_\alpha$ can be determined in the following way: restriction to $C_\alpha$
gives a natural homomorphism $\Pic(X)\to \Pic(C_\alpha)$, and the $-2$-curves
are precisely the simple roots in the system of positive roots in the
kernel of this homomorphism.  This is easy to see from our perspective, as
it reduces to checking when a simple root of $E_{m+1}$ is a $-2$-curve
disjoint from $C_\alpha$.

Given an anticanonical rational surface $X$, there is a natural
combinatorial invariant of blowdown structures, namely how the components
of $C_\alpha$ (which we fix an ordering of) are expressed in terms of the
corresponding basis.  That is, if we fix an ordered decomposition
\[
C_\alpha = \sum_i c_i C_i
\]
where the $C_i$ are the distinct components of $C_\alpha$ (so each
$c_i>0$), then given any blowdown structure, we may associate the sequence
of pairs $(c_i,v_i)$ where $v_i\in \Z^{m+2}$ is the image of $C_i\in
\Pic(X)$ under the isomorphism $\Pic(X)\cong \Z^{m+2}$ corresponding to
$\Gamma$.  (If $C_\alpha$ is integral, this invariant simply distinguishes
between even and odd blowdown structures.)  If we fix $X$ and an ordering
of the components of $C_\alpha$, this invariant takes on only finitely many
values as we vary $\Gamma$.  In fact, something much stronger holds: if we
take the union over all anticanonical surfaces and all orderings of their
anticanonical components, then for any lower bound on the self-intersections of
the anticanonical components, there are only finitely many choices for this
combinatorial invariant.

Indeed, if we put a lower bound on the self-intersections of the components
of $C_\alpha$, then this implies a lower bound on the self-intersections of
any $-d$-curve on $X$ (since any $-d$-curve with $d>2$ has negative
intersection with $C_\alpha$, so is a component).  In particular, this gives
only finitely many possible Hirzebruch surfaces that $X$ can be blown down
to.  On a given Hirzebruch surface, there are only finitely many
combinatorially distinct decompositions of anticanonical curves, and as we
blow up points, the change in invariant only depends on the set of
components containing the point being blown up.

In particular, this combinatorial type splits the weak groupoid of blowdown
structures into finitely many groupoids.  We observe that each of these
groupoids is a quotient groupoid $G/H$ for some $G\subset W(E_{m+1})$,
where $H$ is the group generated by reflections in $-2$-curves disjoint
from $C_\alpha$.  Indeed, whether a simple root is admissible only depends on
the combinatorial type, and thus the admissible elements of $W(E_{m+1})$
preserving the combinatorial type form a group.  This group is certainly
contained in the stabilizer of the sequence of vectors corresponding to the
components of $C_\alpha$, and itself contains a reflection group.

\begin{prop}
  Suppose $\rho$ is a positive root which is orthogonal to every component of
  $C_\alpha$.  Then the corresponding reflection is admissible.
\end{prop}

\begin{proof}
  It suffices to consider the case that $\rho$ is {\em simple} among the
  root system of positive roots orthogonal to every component of
  $C_\alpha$.  Then we claim that there is a blowdown structure in which
  $\rho$ is a simple root of $E_{m+1}$.  If $\rho$ is already simple in
  $E_{m+1}$, this is immediate.  Otherwise, let $\sigma$ be a simple root
  such that $\sigma\cdot \rho<0$.  If $\sigma$ is ineffective, we may
  reflect in $\sigma$ and proceed by induction.  Otherwise, $\sigma$ is
  effective, and some component $c$ of $\sigma$ which is not a component of
  $C_\alpha$ satisfies $c\cdot \rho<0$.  Since $f\cdot \rho\ge 0$ for every
  positive root, we have $c\ne f$, and thus by the proof of Lemma
  \ref{lem:fixed_of_simple}, $c$ must be a fixed component of $\sigma$.  If
  $c^2=-2$, then it is orthogonal to every component of $C_\alpha$, but
  then the fact that $c\cdot \rho<0$ contradicts simplicity of $\rho$
  unless $\rho=c$, in which case reducing to a simple root of $E_{m+1}$ is
  straightforward.

  Otherwise, by the classification of fixed components of $-2$-curves, we
  find that $c=e_i$ for some $i$, and that $e_j\cdot \sigma=0$ for $j\ge
  i$.  Since the only positive roots satisfying $\rho\cdot e_i<0$ are those of
  the form $e_i-e_j$ for some $j>i$, we conclude that $\rho\cdot
  \sigma=(e_i-e_j)\cdot \sigma=0$, a contradiction.
\end{proof}

\begin{rem}
  In general, the full stabilizer need not be a reflection subgroup.  For
  instance, let $X=X_8$ be a rational elliptic surface with an
  anticanonical curve of Kodaira type $I_3^*$ (corresponding to the root
  system $\tilde{D}_7$).  Since a subsystem of type $D_7$ in $E_8$ has
  trivial stabilizer, we conclude that the stabilizer must be contained in
  the translation subgroup of $W(E_9)=W(\tilde{E}_8)$ in this case.  Any
  translation will add some multiple of $-K_X$ to the different components,
  preserving the property that the relevant linear combination is $-K_X$;
  it follows that the stabilizer contains a corank $7$ subgroup of the
  translation subgroup of $E_8$.  In other words, the stabilizer is
  isomorphic to $\Z$, so is certainly not a reflection subgroup!  (We can
  also directly verify in this case that the generator of the stabilizer is
  admissible.)  The stabilizer also fails to be a reflection subgroup of
  $W(E_9)$ when the anticanonical curve has Kodaira type $I_7$, and in one
  of the two ways it can have Kodaira type $I_8$.  In the $I_8$ case it is
  again isomorphic to $\Z$, while in the $I_7$ case it is isomorphic to
  $W(\tilde{A}_1)\times \Z$.
\end{rem}

\begin{cor}
  If $X$ is a rational surface with an integral anticanonical curve, then
  the weak groupoid of blowdown structures on $X$ is a union of two
  isomorphic quotient groupoids of the form $W(E_{m+1})/H$ where $H$ is the
  group generated by reflections in the $-2$ curves of $X$.
\end{cor}

\begin{rem}
  The reader should be cautioned that a reflection subgroup of an infinite
  Coxeter group need not have finite rank.  Indeed, an example was given in
  \cite[Ex.~2.8]{HarbourneB:1985} of a rational surface with (nodal)
  integral anticanonical curve and infinitely many $-2$-curves, and thus
  the stabilizers in the corresponding groupoid have infinite rank.
\end{rem}

\section{Nef divisors}
\label{sec:nef_comm}

Given a rational surface and blowdown structure, one natural question which
arises is whether a given vector corresponds to an effective divisor class,
or one with an integral representative.  For the latter, it will be helpful
to also have an answer to the question of which vectors correspond to nef
divisor classes (i.e., those having nonnegative intersection with every
effective divisor).  This is complicated in general, but in the
anticanonical case, is quite tractable.  Note that as with the above
algorithms for recognizing $-2$- and $-1$-curves, the algorithms below only
depend on (a) the decomposition of $C_\alpha$ in some initial choice of
blowdown structure, and (b) the kernel of the natural homomorphism
$\Pic(X)\to \Pic(C_\alpha)$.  (The latter is not particularly tractable to
compute in general, but in fact all we really need is the ability to test
membership in the kernel.)

One answer to this question was given in \cite{LahyaneM/HarbourneB:2005}:
the monoid of effective divisors (assuming $K_X^2<8$) is generated by the
$-d$-curves with $d>0$ and $-K_X$.  (We will give a constructive version of
this statement in Section \ref{sec:effdivs} below.)  In principle, this
gives a way of testing whether a vector is nef: simply check that it has
nonnegative intersection with $-K_X$ and every $-d$-curve.  Of course, this
is not an actual algorithm, for the simple reason that a rational surface
can have infinitely many smooth curves of negative self-intersection.  (In
fact, this is the {\em typical} behavior!)

We do, however, obtain the following.  Given an anticanonical rational
surface $X$ and a blowdown structure $\Gamma$ for $X$, define the {\em
  fundamental chamber} to be the monoid in $\Pic(X)$ consisting of classes
having nonnegative intersection with every simple root for $\Gamma$.
Note also that there are only finitely many $-d$-curves on $X$ with $d>2$,
since any such curve must be a component of $-K_X$.

\begin{prop}
  Suppose $X$ is an anticanonical rational surface with $K_X^2\le 6$, and
  let $\Gamma$ be a blowdown structure on $X$.  Let $D$ be a divisor class
  in the fundamental chamber of $\Gamma$.  Then $D$ is nef iff $D\cdot
  C_\alpha\ge 0$, $D\cdot e_m\ge 0$, and $D$ has nonnegative
  self-intersection with every $-d$-curve with $d>2$.  If $C_\alpha^2\ge
  0$, then we can omit the condition $D\cdot C_\alpha\ge 0$.
\end{prop}

\begin{proof}
  That nef divisors satisfy the given conditions is
  obvious, so it remains to show that the given conditions imply that $D$
  is nef.  Since we have assumed $D\cdot C_\alpha\ge 0$, it remains to verify
  that it has nonnegative intersection with every $-d$-curve for $d>0$.  We
  have also assumed this for $d>2$, so only the cases of $-2$- and
  $-1$-curves remain.  Any $-2$-curve is a positive root of $E_{m+1}$, so
  is a nonnegative linear combination of simple roots.  By assumption, $D$
  has nonnegative intersection with every simple root, so nonnegative
  intersection with every positive root.  Similarly, we saw above that any
  $-1$-curve can be written as $e_m$ plus a nonnegative linear combination
  of simple roots, so again $D$ has nonnegative intersection with every
  $-1$-curve.

  If $K_X^2\ge 0$, then we can write $K_X$ as a nonnegative linear combination
  of simple roots and $e_m$, so can omit the corresponding condition.
\end{proof}

\begin{rem}
  Something similar holds when $K_X^2=7$, except that we must also assume
  $D\cdot (f-e_1)\ge 0$.  In any event, we can readily write down the
  effective and nef monoids when $K_X^2=7$.  Indeed, if $X_0\cong F_{2d}$ and
  $X_1\to X_0$ blows up a point of $s_{\min}$ (which we arrange to occur if
  $d=0$), we have
\begin{align}
\text{Eff}(X_1) &= \langle s-df-e_1,f-e_1,e_1\rangle\notag\\
\text{Nef}(X_1) &= \langle f,s+df,s+(d+1)f-e_1\rangle,\notag
\end{align}
while if $d>0$ and $X_1\to X_0$ blows a point not on $s_{\min}$, we have
\begin{align}
\text{Eff}(X_1) &= \langle s-df,f-e_1,e_1\rangle\notag\\
\text{Nef}(X_1) &= \langle f,s+df,s+df-e_1\rangle.\notag
\end{align}
Since $s-df=s_{\min}$, it is clear that the putative generators for
$\text{Eff}(X_1)$ are effective.  Similarly, we find that the putative
generators for $\text{Nef}(X_1)$ have nonnegative self-intersection, and
can be represented by integral divisors, so are nef.  Since in each case
the two bases are dual to each other, they must actually be the effective
and nef monoids.  (The corresponding bases for the monoids relative to an
odd blowdown structure can of course be obtained by an elementary
transformation.)  Note that in either case, $\text{Eff}(X_1)$ is a
simplicial cone generated by $-e$-curves with $e<0$.  Similarly, for
$K_X^2=8$, the effective cone is generated by $s_{\min}$ and $f$.
\end{rem}

Since nef divisors are effective (\cite[Cor.~II.3]{HarbourneB:1997}), we
also conclude that any class $D$ satisfying the above hypotheses is
effective.  In fact, we can do better: we can give an explicit effective
divisor representing $D$.

\begin{prop}\label{prop:decomp_nef}
With hypotheses as above, $D$ can be written as a nonnegative linear
combination of $-d$-curves and $-K_X$.
\end{prop}

\begin{proof}
Suppose first that $D\cdot e_m=0$, so that $D$ is the total transform of a
divisor on $X_{m-1}$.  If $m>2$, then this divisor on $X_{m-1}$ is itself
in the fundamental chamber, and has nonnegative intersection with
$-K_{X_{m-1}}$.  We can thus, by induction, decompose it into $-d$-curves
and copies of the anticanonical curve on $X_{m-1}$.  The total transform of
a $-d$-curve is either a $-d$-curve or the sum of a $-(d+1)$-curve and
$e_m$, while the total transform of the anticanonical curve on $X_{m-1}$ is
$C_\alpha+e_m$, and thus the decomposition on $X_{m-1}$ induces a
decomposition on $X_m$ which again has the desired form.

Similarly, if $m=2$ and $D\cdot e_2=0$, then we still find that $D$ is the
total transform of a nef divisor on $X_1$, and thus obtain the desired
decomposition of $D$ by expanding it in the basis of the simplicial cone
$\text{Eff}(X_1)$.

Finally, suppose $D\cdot e_m>0$.  Then we claim that $D-C_\alpha$ satisfies the
original hypotheses.  Indeed, if $C$ is a $-d$-curve for $d>2$, then
\[
(D-C_\alpha)\cdot C = D\cdot C + (d-2)> D\cdot C,
\]
while if $\sigma$ is a simple root, then $(D-C_\alpha)\cdot \sigma = D\cdot
\sigma\ge 0$.  In addition, $(D-C_\alpha)\cdot e_m = D\cdot e_m-1\ge 0$.
Finally, we have
\[
(D-C_\alpha)\cdot C_\alpha = D\cdot C_\alpha - C_\alpha^2.
\]
Either $C_\alpha^2<0$, so the inequality becomes stronger, or $C_\alpha^2\ge
0$, and the inequality is redundant.  Either way, $D-C_\alpha$ satisfies all
of the hypotheses, and we obtain an explicit decomposition of the given form.
\end{proof}

\begin{lem}
If $D$ is a nef divisor on the rational surface $X$ with $K_X^2\le 6$, then
there exists a blowdown structure (of either parity) such that $D$ is in
the fundamental chamber.  Moreover, the representation of $D$ in the basis
corresponding to such a blowdown structure depends only on the parity.  In
addition, if $e$ is a $-1$-curve with $e\cdot D=0$, then the blowdown
structure can be chosen in such a way that $e_m=e$.
\end{lem}

\begin{proof}
  Choose a blowdown structure of the desired parity on $X$.  If $D$ is not
  already in the fundamental chamber, then there exists a simple root
  $\sigma$ such that $D\cdot \sigma<0$.  Since $D$ is nef, $\sigma$ cannot
  be effective, and thus we can apply the corresponding reflection.  Either
  $\sigma$ is in the subsystem of type $D_m$ (which can only occur finitely
  many times in a row, since that subgroup is finite), or it decreases
  $D\cdot f$.  The latter is nonnegative since $f$ is effective, and thus
  the process will terminate after a finite number of steps.

For uniqueness, suppose $D$ is in the fundamental chamber of both $\Gamma$
and $\Gamma'$, two blowdown structures of the same parity.  Then there
exists an ineffective element $w\in W(E_{m+1})$ such that
$\Gamma'=w\Gamma$.  If $w=1$, then we are done; otherwise, there is an
ineffective simple root $\sigma$ of $\Gamma$ such that $w\sigma$ is
negative (the last root in some reduced word for $w$).  Since $D$ is in the
fundamental chamber for both $\Gamma$ and $\Gamma'$, it has nonnegative
intersection with every positive root of either blowdown structure.  Thus
$D\cdot \sigma\ge 0$ since $\sigma$ is positive for $\Gamma$, and $D\cdot
(-\sigma)\ge 0$ since $-\sigma$ is positive for $\Gamma'$.  In other words,
$D\cdot \sigma=0$, and thus the reflection in $\sigma$ does not change the
expansion of $D$ in the standard basis.  The claim follows by induction on
the length of the reduced word for $w$.

Finally, if $\Gamma$ is any blowdown structure with $e_m=e$, then $D\cdot
(e_{m-1}-e_m)=D\cdot e_{m-1}\ge 0$, and thus the algorithm for putting $D$
in the fundamental chamber will never try to apply the corresponding
reflection, so will never change $e_m$.
\end{proof}

\begin{rem}
  Uniqueness is of course a standard fact from Coxeter theory when we
  restrict to $D$ in the root lattice, and the above argument is adapted
  from the standard one.
\end{rem}

This then gives us the desired algorithm for testing whether a divisor is
nef: First check that it has nonnegative intersection with every $-d$-curve
with $d>2$ and with $C_\alpha$, then repeatedly attempt to reflect in simple
roots with $D\cdot \sigma<0$.  If at any step we have $\sigma$ effective, $D\cdot
f<0$ or $D\cdot e_m<0$, then $D$ is not nef; otherwise, we terminate in the
fundamental chamber, and conclude that $D$ is nef.

\section{Effective divisors}
\label{sec:effdivs}

A similar algorithm works for testing whether a divisor $D$ is effective.
We assume $K_X^2<7$, since otherwise the effective cone is simplicial, so
testing whether $D$ is effective is just linear algebra.

Again, we start by choosing any blowdown structure for $X$, and if at any
step in the process we obtain a divisor with $D\cdot f<0$, we halt with the
conclusion that our divisor was ineffective.  We perform the following
steps, as specified.
\begin{itemize}
\item[1.] If there exists a component $C$ of $C_\alpha$ such that
  $C^2,D\cdot C<0$, then replace $D$ by $D-C$, and repeat step 1.
\item[2.] If $D\cdot e_m<0$, then replace $D$ by $D+(D\cdot e_m)e_m$
  and go back to step 1.
\item[3.] If $D$ is in the fundamental chamber, conclude that the
  original divisor was effective.  Otherwise, choose the lexicographically
  smallest simple root such that $D\cdot \sigma<0$.  If $\sigma$ is effective,
  replace $D$ by $D-\sigma$ and go back to step 1; otherwise, replace $\Gamma$
  by $r_\sigma\Gamma$ and go back to step 2.
\end{itemize}

To see that this algorithm works, we note as before that $f$ is nef, so any
divisor with $D\cdot f<0$ is not effective.  Whenever we replace $D$ by
$D-C$ in step 1, $C$ is an integral curve of negative self-intersection
intersecting $D$ negatively.  But then $D$ is effective iff $D-C$ is
effective; one direction is obvious, while if $D$ is effective, then $C$ is
a fixed component of $D$.  The same argument applies in step 2, while in
step 3 when $\sigma$ is effective, either $\sigma$ is irreducible (so again
the argument applies) or we have $D\cdot c<0$ for some fixed component $c$
of $\sigma$.  We must have $c^2\ge -2$, else $c$ would have been removed in
step $1$; and similarly $c\ne e_m$.  But then the classification of fixed
components of effective simple roots lets us find a lexicographically
smaller simple root having negative intersection with $D$.

Since we terminate at a nef divisor in the fundamental chamber, this
algorithm also gives us an explicit decomposition of $D$ as a nonnegative
linear combination of $C_\alpha$ and $-d$-curves.  In this context, we note
the following.

\begin{prop}
  Let $X$ be an anticanonical rational surface with $K_X^2<8$.  Either
  every representative of $-K_X$ is integral, or some representative is a
  nonnegative linear combination of $-d$-curves with $d\ge 1$.
\end{prop}

\begin{proof}
If some representative of $-K_X$ is reducible, then we can write
\[
-K_X=D_1+D_2
\]
for nonzero effective divisors $D_1$, $D_2$, and it suffices to show that
each $D_i$ is linearly equivalent to a nonnegative linear combination of
$-d$-curves.  Since $D_1$ is effective by assumption, we can write
\[
D_1 \sim m(-K_X) + \sum_j c_j C_j
\]
where each $C_j$ is a $-d$-curve for some $d\ge 1$ and all coefficients are
nonnegative.  This implies $D_1+mK_X$ is effective, and thus $-D_2=D_1+K_X =
D_1+mK_X+(m-1)C_\alpha$ is effective, unless $m=0$.  In other words, $D_1$
has a decomposition as required.
\end{proof}

\begin{cor}
  Let $X$ be an anticanonical rational surface, and suppose $K_X^2\notin
  \{0,1,8,9\}$.  Then the effective monoid of $X$ is generated by the
  integral curves of negative self-intersection.
\end{cor}

\begin{proof}
  If $K_X^2<0$, then $C_\alpha$ has negative self-intersection, and is either
  integral or redundant.  For $1<K_X^2<8$, we note that $\Gamma(-K_X)$
  corresponds to a codimension $m$ subspace of $\Gamma(-K_{X_0})$.  On a
  Hirzebruch surface, either every anticanonical curve is reducible (i.e.,
  on $F_d$ for $d>2$), or the reducible anticanonical curves form a
  codimension $2$ subvariety of the $8$-dimensional projective space of all
  anticanonical curves.  We are imposing $m\le 6$ linear conditions on this
  projective variety, and thus obtain a nonempty set of anticanonical
  curves on $X$ which are reducible on $X_0$ and thus reducible on $X$.
\end{proof}

\begin{rem}
  Similarly, if $K_X^2=1$ but $X$ has a $-2$-curve, then some anticanonical
  curve is reducible (and in fact contains that $-2$-curve).  Also, in any
  case the anticanonical divisor is not needed to generate the {\em
    rational} effective cone when $K_X^2=1$, since then $-2K_X =
  (-2K_X-e_7)+e_7$ is a sum of effective divisors.
\end{rem}

We can also adapt the algorithm to compute $h^0(\sO_X(D))=\dim
H^0(\sO_X(D))$ for an effective divisor.  Indeed, every step of the
algorithm either removes a fixed component of $D$ or leaves $D$ alone and
changes the blowdown structure, and thus the algorithm terminates with a
nef divisor $D'$ with a natural isomorphism
\[
H^0(\sO_X(D'))\cong H^0(\sO_X(D)).
\]
Thus to compute the dimensions of effective linear systems, it remains only
to compute the dimensions of linear systems corresponding to nef divisors
in the fundamental chamber.  So let $D$ be such a divisor class and suppose
$m=0$ or $D\cdot e_m>0$, since otherwise we may as well consider $D$ as a
divisor on $X_{m-1}$.  Note that $(K_X-D)\cdot f=-2-(D\cdot f)<0$, so that
$K_X-D$ is ineffective and thus $h^2(\sO_X(D))=0$ by Serre duality.

If $D\cdot C_\alpha>0$, then it follows from
\cite[Thm.~III.1(ab)]{HarbourneB:1997} that $h^1(\sO_X(D))=0$, and thus
we can use Hirzebruch-Riemann-Roch to compute
\[
h^0(\sO_X(D)) = \chi(\sO_X(D)) = \frac{D\cdot (D+C_\alpha)}{2}+1.
\]
This in particular holds whenever $m<8$, since then either $D\cdot
C_\alpha>0$ or $D=0$.

If $m\ge 8$ and $D\cdot C_\alpha=0$, then from the proof of Proposition
\ref{prop:decomp_nef}, we find that $D-C_\alpha$ is also nef.  Now consider
the short exact sequence
\[
0\to \sO_X(D-C_\alpha)\to \sO_X(D)\to \sO_X(D)|_{C_\alpha}\to 0.
\]
From \cite[Thm.~III.1(d)]{HarbourneB:1997}, we find that the natural inclusion
\[
H^0(\sO_X(D-C_\alpha))\subset H^0(\sO_X(D))
\]
is an isomorphism iff the line bundle $\sO_X(D)|_{C_\alpha}$ is
nontrivial.  Since $h^0(\sO_{C_\alpha})=1$, we conclude that
\[
h^0(\sO_X(D))
=
\begin{cases}
h^0(\sO_X(D-C_\alpha))+1 & \sO_X(D)|_{C_\alpha}\cong \sO_{C_\alpha}\\
h^0(\sO_X(D-C_\alpha)) & \text{otherwise.}
\end{cases}
\]
If $m>8$, then $D-C_\alpha$ is a nef divisor with $(D-C_\alpha)\cdot C_\alpha>0$,
so we reduce to the previous case.  If $m=8$, then $D=r C_\alpha$ for some
$r\ge 1$, and thus $C_\alpha$ must be nef.  We deduce that either $C_\alpha$ is
integral or every component of $C_\alpha$ is a $-2$-curve.  Moreover, the
above recurrence tells us that in this case,
\[
h^0(\sO_X(rC_\alpha)) = \lfloor r/r'\rfloor +1,
\]
where $r'$ is the order of the bundle $\sO_X(C_\alpha)|_{C_\alpha}$ in
the group $\Pic(C_\alpha)$.  (In particular, $h^0(\sO_X(rC_\alpha))=1$ if this
bundle is not torsion.)

\begin{rem}
If
\[
D = ns+df-\sum_i r_i e_i
\]
relative to some even blowdown structure, then
\[
\chi(\sO_X(D))
=
(n+1)(d+1)-\sum_i \frac{r_i(r_i+1)}{2}.
\]
This of course corresponds to the fact that $H^0(\sO_X(D))$ is a
subspace of $H^0(\sO_X(ns+df))$ cut out by the appropriate number of
linear conditions.  (If $X$ blows up $m$ distinct points of $X_0$, the
conditions are simply that the curve have multiplicity $r_i$ at the $i$-th
point.)  In principle, one could determine $h^0(\sO_X(D))$ (and in
particular test whether $D$ is effective) using linear algebra, but the
above approach scales better, and largely separates out the combinatorial
influences from the algebraic influences.
\end{rem}

This algorithm easily extends to allow us to compute the Betti numbers of any
line bundle $\sO_X(D)$: if $D$ is ineffective, then $h^0(\sO_X(D))=0$ while
otherwise we may compute it using the above algorithm, and similarly for
$h^2(\sO_X(D))=h^0(\sO_X(K_X-D))$, allowing us to deduce the remaining
Betti number from the Euler characteristic.

\section{Integral divisors}

By Lemma II.6 and Theorem III.1 of \cite{HarbourneB:1997}, there is a
relatively short list of possible ways that a nef class can fail to be
generically integral.  (The integral classes which are not nef are
precisely the $-d$-curves for $d\ge 1$ and the anticanonical divisor, when
this is integral and has negative self-intersection, and we already know
how to recognize those.)  Although the description given there is purely
geometric, it turns out to be easy enough to recognize the different cases
in terms of the representation of the divisor in a fundamental chamber.
Since this representation is unique, we can (and will) figure out how each
case is represented by placing various geometrically motivated constraints
on the blowdown structure, and checking that the result is in the
fundamental chamber.

\begin{rem}
  Note that in characteristic 0, ``generically integral'' and ``generically
  smooth'' are the same on an anticanonical rational surface: a generically
  integral divisor class on a rational surface has at most one base point,
  and if it does, meets $C_\alpha$ at that point with
  multiplicity 1.  Bertini's theorem implies that the generic
  representative is smooth away from the base point, and the intersection
  with $C_\alpha$ implies smoothness there.
\end{rem}

A linear system $|D|$ is a {\em pencil} if $|D|\cong \P^1$ and the base
locus is $0$-dimensional.

\begin{lem}
  Let $D$ be a divisor on $X$, and suppose $\Gamma$ is an even blowdown
  structure such that $D$ is in the fundamental chamber.  Then $|D|$ is a
  pencil iff one of the following three cases occurs.
\begin{itemize}
\item[(a)] $D=f$.
\item[(b)] $D=2s+2f-e_1-\cdots-e_7$, $D$ is nef, and
  $\sO_X(2s+2f-e_1-\cdots-e_7-e_k)|_{C_\alpha}\not\cong \sO_{C_\alpha}$ for
  $8\le k\le m$.
\item[(c)] $D=r(2s+2f-e_1-\cdots-e_8)$, where
$\sO_X(2s+2f-e_1-\cdots-e_8)|_{C_\alpha}$ is a line bundle of exact order
$r$ in $\Pic(C_\alpha)$.
\end{itemize}
\end{lem}

\begin{proof}
  A pencil is certainly generically integral (lest $X$ be reducible), so
  nef.  Per \cite[Lem.~II.6]{HarbourneB:1997}, there are three possibilities:
$D^2=0$, $D\cdot C_\alpha=2$; $D^2=D\cdot C_\alpha=1$; or $D^2=D\cdot
C_\alpha=0$.

In the first case, the generic fiber of $D$ has arithmetic genus 0, so $D$
is the class of a rational ruling.  It follows that there exists a blowdown
structure such that $D=f$, and we readily verify that $f$ is in the
fundamental chamber.  (If the blowdown structure we end at is odd, simply
perform an elementary transformation, and note that this preserves the
meaning of $f$.)

For the case $D^2=D\cdot C_\alpha=0$, $D$ gives a quasi-elliptic fibration of
$X$, and we can choose a blowdown structure in which we first blow down any
$-1$-curves contained in fibers.  After doing so, we end up at a relatively
minimal quasi-elliptic surface, which must be $X_8$ for the blowdown
structure.  The only isotropic vectors in $\Pic(X_8)$ are the multiples of
the canonical class, and thus $D=r(2s+2f-e_1-\cdots-e_8)$ for some $r$;
again, this is in the fundamental chamber.  For this to be a pencil, it
must not have any fixed component, so $\sO_X(D)|_{C_\alpha}\cong
\sO_{C_\alpha}$ and $\sO_X(-K_{X_8})|_{C_\alpha}$ has order dividing $r$.
If the order strictly divides $r$, then $\sO_X(D)$ will have more than
$2$ global sections.

For the case $D^2=D\cdot C_\alpha=1$, the linear system is again
quasi-elliptic, now with a base point.  The base point is on the
anticanonical curve, namely the unique point such that
\[
\sO_X(D)|_{C_\alpha}\cong {\cal L}_{C_\alpha}(p).
\]
The fibers of $D$ transverse to $C_\alpha$ are either integral or contain a
single $-1$-curve, while the fiber not transverse to $C_\alpha$ contains
$C_\alpha$.  The residual divisor $D-C_\alpha$ has arithmetic genus $1-r$,
where $r=C_\alpha^2-1$, and thus has at least $r$ connected components, each
of which has negative self-intersection (since it is orthogonal to a class
of positive self-intersection).  It follows that every component has
self-intersection $-1$ and arithmetic genus 0, and thus contains a
$-1$-curve.  We may thus choose a blowdown structure in which we first blow
down those $-1$-curves until eventually reaching $X_7$ and
$D=2s+2f-e_1-\cdots-e_7$.  This is a pencil on $X_7$ precisely when it is
nef, and remains a pencil on $X$ as long as we never blow up the base
point.
\end{proof}

\begin{prop}\label{prop:non_integral2}
  Suppose $D$ is a nef divisor class with $D\cdot C_\alpha\ge 2$,
  and let $\Gamma$ be an even blowdown structure for which $D$ is in the
  fundamental chamber.  Then $|D|$ is generically integral unless $D=rf$ for
  some $r>1$.
\end{prop}

\begin{proof}
  By \cite[Thm.~III.1(a)]{HarbourneB:1997}, $D$ is base point free, so is
  generically integral unless it is a strict multiple of a pencil.
\end{proof}

The case $D\cdot C_\alpha=0$, which is the most interesting for us in any
event, is the next easiest case to handle.  If $\sO_X(D)|_{C_\alpha}\not\cong
\sO_{C_\alpha}$, then $D$ can only be integral if $D=C_\alpha$ and
$C_\alpha$ is integral (\cite[Thm.~III.1(d)]{HarbourneB:1997}).

\begin{thm}\label{thm:non_integral}
  Let $X$ be an anticanonical rational surface, let $D$ be a nef divisor
  class on $X$ such that $\sO_X(D)|_{C_\alpha}\cong \sO_{C_\alpha}$, and let
  $\Gamma$ be an even blowdown structure such that $D$ is in the
  fundamental chamber.  Then $|D|$ is generically integral unless one of the
  following two possibilities occurs.
\begin{itemize}
\item[(a)] $D=r(2s+2f-e_1-\cdots-e_8)$, and
  $\sO_X(2s+2f-e_1-\cdots-e_8)|_{C_\alpha}$ is a line bundle of order $r'$
  strictly dividing $r$.  Then the generic representative of $D$ is a
  disjoint union of $r/r'$ curves of genus 1, of divisor class
  $r'(2s+2f-e_1-\cdots-e_8)$.
\item[(b)] $D=r(2s+2f-e_1-\cdots-e_8)+e_8-e_9$ with $r>1$, and
  $\sO_X(2s+2f-e_1-\cdots-e_8)|_{C_\alpha}\cong \sO_{C_\alpha}$.  Then the
  generic representative of $D$ is the union of $r$ divisors of class
  $2s+2-e_1-\cdots-e_8$ (all of genus 1) and a $-2$-curve of class
  $e_8-e_9$.
\end{itemize}
\end{thm}

\begin{proof}
  $D$ is generically integral unless it factors through a pencil or has a
  fixed component.  The first case is precisely (a) above, by the
  classification of pencils.  The second case is described in
  \cite[Thm.~III.1(c)]{HarbourneB:1997}: $D$ has a unique fixed component
  $N$ which is a $-2$-curve, and $D-N$ is a strict multiple of a pencil $P$
  with $P\cdot N=1$.  In particular, there exists a blowdown structure such
  that $P$ is the total transform of some antipluricanonical pencil on
  $X_8$, and ``pluri'' can be ruled out by the fact that $P\cdot N=1$.
  Now, $N$ cannot be contracted by the map $X\to X_8$, since $P$ is still
  base point free on $X_8$; thus $N$ is a rational curve, and since $N\cdot
  P=1$, must be a $-1$-curve.  We can thus further insist that the map
  $X_8\to X_7$ blows down $N$.  Since $N$ is a $-2$-curve on $X$, the map
  $X\to X_8$ blows up a point of $N$ exactly once, and we may insist that
  this is the first point blown up after reaching $X_8$; i.e., that $N$ is
  already a $-2$-curve on $X_9$.  But then $N=e_8-e_9$, and $D$ has the
  claimed form, which we verify is in the fundamental chamber.
\end{proof}

\begin{rems}
For multiplicative Deligne-Simpson problems, a rather more complicated
irreducibility condition was given in \cite{Crawley-BoeveyW/ShawP:2006}.  In
particular, the above result gives a much stronger statement in the case of
$3$- and $4$-matrix multiplicative Deligne-Simpson problems, and it is
natural to wonder if a similarly strong result holds in general.
\end{rems}

The remaining case $C_\alpha\cdot D=1$ can be dealt with in one of two ways.
The easiest is to blow up the intersection with $C_\alpha$, and consider the
strict transform $D'$ of $D$ on $X_{m+1}=:X'$, a divisor class which is
generically disjoint from the new anticanonical curve.  The above
algorithms tell us how to determine the generic decomposition of such a
divisor class: first use the algorithm for testing effectiveness to write
it as a sum of (fixed) $-2$-curves and a nef class in some fundamental
chamber, then use the above result to decompose the latter class.  The
generic decomposition of $D'$ on $X'$ corresponds directly to the generic
decomposition of $D$ on $X$, since $D'$ is a strict transform, and thus
this procedure computes the generic decomposition of $D$.

We can also work out what the nonintegral cases look like in the
fundamental chamber, again using a result of Harbourne,
\cite[Thm.~III.1(b)]{HarbourneB:1997}.  We omit the details.

\begin{prop}\label{prop:non_integral1}
  Suppose $D$ is a nef divisor class on $X$ such that $D\cdot C_\alpha=1$,
  and let $\Gamma$ be a blowdown structure for which $D$ is in the
  fundamental chamber.  Then $D$ is generically integral except in the
  following two cases.
\begin{itemize}
\item[(a)] For some $1\le i\le m$, $D\cdot e_i=0$ and
$\sO_X(D-e_i)|_{C_\alpha}\cong \sO_{C_\alpha}$.  The fixed part of $D$ is
the total transform of the minimal such $e_i$.
\item[(b)] $D=r(2s+2f-e_1-\cdots-e_8)+e_8$ and
  $\sO_X(2s+2f-e_1-\cdots-e_8)|_{C_\alpha}\cong \sO_{C_\alpha}$.  The fixed
  part of $D$ is the total transform of $e_8$.
\end{itemize}
\end{prop}

We also mention a necessary condition for a divisor to be integral, related
to the theory of Coxeter groups of Kac-Moody type.  The convention there is
to consider both ``real'' roots (i.e., roots in the usual sense) and
``imaginary'' roots.  The latter are defined as integral vectors whose
orbit intersects the fundamental chamber in a nonnegative linear
combination of simple roots.

\begin{prop}
  Any integral divisor $D$ such that $D\cdot K_X=0$ is a positive root
  (real or imaginary).
\end{prop}

We have already shown this for $-2$-curves (i.e., that $-2$-curves are
positive real roots), while for nef curves it is a consequence of the
following more general fact.  Note that since $e_m\cdot C_\alpha=1$, the fact
that the simple roots are a basis of $C_\alpha^\perp$ implies that together
with $e_m$, they form a basis of $\Pic(X)$.

\begin{prop}
Let $X$ be an anticanonical rational surface with $K_X^2<7$, and
$D$ a nef divisor class on $X$.  Then for any blowdown structure
on $X$, $D$ is a nonnegative linear combination of the simple roots
and $e_m$.  In fact, if we write (for an even blowdown structure)
\[
D = a(s-f)+b(f-e_1-e_2)+\sum_{1\le i<m}c_i(e_i-e_{i+1})+c_m e_m,
\]
then we have the inequalities
\begin{align}
c_2&\ge c_3\ge\cdots\ge c_m\ge 0\notag\\
c_2&\ge b\ge a\ge 0.\notag\\
c_2&\ge c_1\ge 0\notag
\end{align}
\end{prop}

\begin{proof}
The divisor class $e_i$ is effective for all $i$ (since it is the
total transform of a $-1$-curve on $X_i$), and thus $D\cdot e_i\ge
0$.  Taking $i\ge 3$, we conclude that
\[
c_2\ge c_3\ge \cdots \ge c_m.
\]
To see that $c_m\ge 0$, we note that
\[
c_m = (c_m e_m)\cdot C_\alpha = D\cdot C_\alpha\ge 0.
\]
Similarly, the classes $f+s-e_1-e_2$, $s$, and $f$ are effective on
$X_2$, thus on $X$, so taking inner products with $D$ shows
\[
c_2\ge b\ge a\ge 0.
\]
Finally, the effective
classes $f+s-e_2$, $f+s-e_1$ tell us that
\[
c_2\ge c_1\ge 0,
\]
finishing the proof.
\end{proof}

\begin{rems}
Of course, we can weaken the hypothesis ``$X$ anticanonical'' to ``$D\cdot
K_X\le 0$'', since the latter fact was the only way in which we used the
anticanonical curve.
\end{rems}

\begin{rems}
If $D$ is in the fundamental chamber, the same sequences of coefficients
will be convex.
\end{rems}

\chapter{Moduli of surfaces}
\label{chap:moduli_of_comm_surfs}

\section{General surfaces}

One benefit of considering blowdown structures is that it makes the moduli
problem of rational surfaces much better behaved.  Of course, the standard
approach of choosing an ample bundle also works, but in many cases, there
{\em is} no natural ample bundle, making it hard to control the symmetry;
in contrast, as we have seen, working with blowdown structures gives us a
(rational) action of the Coxeter group $W(E_{m+1})$.

One technical detail is that, since anticanonical surfaces can have
nontrivial (and even positive-dimensional) automorphism groups, the moduli
problem is most naturally represented by a {\em stack} rather than a
scheme.  This adds less complication than one might think at first glance,
as one can accomplish a great deal with stacks without needing to know what
they are; see the discussion in \cite[\S3D]{HarrisJ/MorrisonI:1998}.  In
particular, for our purposes, what we really want to know is that there is
some additional structure we can impose that rigidifies the problem (i.e.,
kills the automorphism group) and is {\em smooth} in the sense that for any
family of surfaces with blowdown structure, the compatible additional
structures are classified by a scheme which is smooth over the base of the
family.  If the surfaces-with-additional-structure are themselves
classified by a scheme, then we can view that scheme as a smooth cover of a
moduli stack of surfaces-without-additional-structure.  This generalizes
somewhat: we may feel free to have multiple notions of additional structure
that only work when certain open conditions are satisfied, as long as those
open conditions cover any family of surfaces.  Any property of schemes that
descends through smooth morphisms makes sense for such stacks, and can be
verified on the smooth covering scheme.  Similarly, one can define the
dimension of the stack as the difference between the dimension of the
smooth cover and the (locally constant) dimension of the fibers.

To construct the moduli stack (an Artin stack) of rational surfaces with
blowdown structures, we first need to construct the moduli stack of
Hirzebruch surfaces.  This is of course essentially just the moduli stack
of rank 2 vector bundles on $\P^1$ (more precisely the moduli stack of
$\PGL_2$-bundles), so is a standard construction, but it will be useful to
keep in mind the details.  (In particular, the construction we use is not
the usual construction for the moduli stack of vector bundles; the extra
structure we use to rigidify the moduli problem has a simpler
interpretation as a structure on $\P(V)$: a basis of
$\Gamma(\sO_{\P^1}(1))$ and a section of given self-intersection.  The
reader who prefers not to think about stacks may feel free to include this
as part of the data in the moduli problems below, and work entirely with
the resulting moduli schemes, adding appropriate constants to the various
statements about dimensions.)

For any integer $d\ge 0$, we have
\[
\Ext^1(\sO_{\P^1}(d+2),\sO_{\P^1})\cong 
H^1(\sO_{\P^1}(-d-2))\cong 
k^{d+1},
\]
and thus the non-split extensions of $\sO_{\P^1}(d+2)$ by $\sO_{\P^1}$ are
classified up to automorphisms of the two bundles by points of the
corresponding $\P^d$.  By a standard construction, this gives rise to a
canonical extension
\[
0\to \sO_{\P^1}\boxtimes \sO_{\P^d}(1) \to V \to \sO_{\P^1}(d+2)\boxtimes
  \sO_{\P^d}\to 0
\]
of sheaves on $\P^1\times \P^d$, each fiber of which is the corresponding
non-split extension of $\sO_{\P^1}(d+2)$ by $\sO_{\P^1}$.

Let $S_k$ denote the locally closed subspace of $\P^d$ on which the fiber
is isomorphic to $V_{d,k}:=\sO_{\P^1}(k+1)\oplus \sO_{\P^1}(d+1-k)$; this
gives a stratification of $\P^d$ by $S_k$ for $0\le k\le d/2$.  Each
stratum can itself be identified as a moduli space of global sections of
the $V_{d,k}$, namely the moduli space of {\em saturated} global sections
(i.e., generating a subbundle), modulo the action of $\Aut(V_{d,k})$.
Since the generic global section is saturated, we have
\[
\dim(S_k) = \dim(\Gamma(V_{d,k}))-\dim(\Aut(V_{d,k}))
          = d-\max(d-2k-1,0);
\]
in other words, $\dim(S_k)=2k+1$ except that $\dim(S_{d/2})=d$.

Since $\dim(\Gamma(V_{d,k}))=d+2$ is independent of $k$, this gives a flat
map to the moduli problem of vector bundles of the form $V_{d,k}$ for $0\le
k\le d/2$.  Taking the relative $\P$ of the bundle gives a smooth map to the
moduli problem of Hirzebruch surfaces; since every Hirzebruch surface
arises in this way for sufficiently large $d$, we find that the moduli
problem of Hirzebruch surfaces is represented by an algebraic stack.  (Note
that if $V$, $V'$ are nonsplit extensions, then $\Hom(V,V')$ is given by a
locally closed subset of global sections of $V'$, so $\Isom(V,V')$ is
indeed a scheme as required, and admits a quotient by $\G_m$ to give
$\Isom_{\P^1}(\P(V),\P(V'))$.)   Note that the stabilizers have the form
$\P\Aut(V_{d,k})\rtimes \Aut(\P^1)$.

This stack has two components (even and odd Hirzebruch surfaces), with
generic fibers isomorphic to $F_0$ and $F_1$ respectively.  In general,
$F_d$ has codimension $d-1$ in the corresponding stack (simply compare
automorphism group dimensions).  Note that the smooth cover corresponding
to $V_d$ classifies pairs $(X,\sigma)$, where $X$ is a Hirzebruch surface
and $\sigma:\P^1\to X$ is an embedding with $\im(\sigma)\cdot f=1$,
$\im(\sigma)^2=d+2$; the map to the moduli stack simply forgets $\sigma$.

To blow up, we proceed as in \cite{HarbourneB:1988}, based on an idea of
Artin.  (There, Harbourne constructed the moduli stack of blowups of
$\P^2$; the extension to Hirzebruch surfaces is straightforward.)  Let
${\cal X}_0$ denote the moduli stack of Hirzebruch surfaces, and let ${\cal
  X}_1$ denote the corresponding universal family of rational surfaces,
with structure maps $\pi_0:{\cal X}_1\to {\cal X}_0$, $\rho:{\cal X}_1\to
{\cal X}_0\times \P^1$.  (This last is something of an abuse of notation;
what we really mean is the $\P^1$-bundle over ${\cal X}_0$ over which the
vector bundles were constructed.  Though this was a product over $\P^d$, we
are quotienting by $\Aut(\P^1)$.) Consider now the problem of classifying
surfaces with $K_X^2=7$.  Such a surface is uniquely determined by a pair
$(X_0,p)$ where $X_0$ is a Hirzebruch surface and $p\in X_0$ is a closed
point.  But points on a surface are classified by the universal surface, so
that the moduli space of rational surfaces with blowdown structure such
that $K_X^2=7$ is precisely ${\cal X}_1$.

Now, extend this to a sequence of stacks ${\cal X}_i$ and morphisms
$\pi_i:{\cal X}_{i+1}\to {\cal X}_i$ for all $i\ge 0$ in the following way.
Using the morphism $\pi_{i-1}$, we may construct the fiber product ${\cal
  X}_i\times_{{\cal X}_{i-1}} {\cal X}_i$ and then blow it up along the
diagonal.  (This is a scheme over ${\cal X}_{i-1}$, so there is no
difficulty in defining the blowup.)  Call the resulting blowup ${\cal
  X}_{i+1}$, and let $\pi_i$ be the morphism induced by the first
projection from the fiber product.

\begin{prop}\label{prop:stack_of_all_surfaces}
  The stack ${\cal X}_i$ represents the moduli problem of rational surfaces
  with blowdown structure and $K_X^2=8-i$.  The universal surface over this
  stack is $\pi_i:{\cal X}_{i+1}\to {\cal X}_i$, and the blowdown structure
  is induced by the maps
\[
{\cal X}_{i+1}\to
{\cal X}_i\times_{{\cal X}_{i-1}}{\cal X}_i
\to
{\cal X}_i\times_{{\cal X}_{i-2}}{\cal X}_{i-1}
\to
\cdots
\to
{\cal X}_i\times_{{\cal X}_0} {\cal X}_1
\to
{\cal X}_i\times \P^1
\]
In addition, for $m\ge 1$, each divisor class $f$, $e_1$,\dots,$e_m$
is represented by a divisor on the universal surface, and there exists a
line bundle of first Chern class $2s$.
\end{prop}

\begin{proof}
This is a simple induction: ${\cal X}_i$ is the universal surface over
${\cal X}_{i-1}$, so also the moduli space of pairs $(X_{i-1},p)$.  To
obtain the universal surface over ${\cal X}_i$, we need to blow up $p$ on
the corresponding fiber, and this is precisely what blowing up the diagonal
does for us.

The claim about divisors is clear for $e_1$,\dots,$e_m$, since these are
just the total transforms of the corresponding exceptional curves.
Similarly, $f-e_1$ is always a $-1$-curve on $X_1$, so gives rise to a
divisor on the universal surface.  Finally, the line bundle
$(\rho^!\sO_{\P^1})^{-1}$ has class $2s$.
\end{proof}

\begin{rem}
  For $m=0$, there is a small difficulty having to do with the fact that
  the $\P^1$ could be twisted; once $m=1$, we have guaranteed that the
  universal $\P^1$ in the construction has a point.  (Of course, the
  anticanonical bundle on the base $\P^1$ is always defined, so lifts to a
  bundle of class $2f$.)  Similarly, although the original construction on
  $\P^d$ comes with a section of the Hirzebruch surface, ${\cal X}_0$
  forgets that section, and the induced automorphisms can act nontrivially
  on the relative $\sO(1)$.
\end{rem}

\begin{cor}
  The moduli stack of rational surfaces with blowdown structure has two
  irreducible components for each value of $K_X^2\le 8$, and each component
  is smooth of dimension $10-2K_X^2$.
\end{cor}

\begin{proof}
  This is clearly true for ${\cal X}_0$ (since the generic Hirzebruch
  surface has a $6$-dimensional automorphism group and the stack is covered
  by open substacks for which some $\P^k\times \PGL_2$-bundle is isomorphic
  to $\P^d$), and each map $\pi_i$ is smooth of relative dimension 2 (being
  a family of smooth projective surfaces).
\end{proof}

\begin{rems}
  Note that the two components are naturally isomorphic for $K_X^2\le 7$:
  just apply the standard elementary transformation.  Also, the formula for
  the dimension holds for $\P^2$ as well, since
  $\dim\Aut(\P^2)=8=2K_{\P^2}^2-10$.
\end{rems}

\begin{rems}
  Essentially the same argument gives an analogous result for higher genus
  rationally ruled surfaces: there are two irreducible components for each
  value of $K_X^2\le 8-8g$, and each component is smooth of dimension
  $10(1-g)-2K_X^2$.  Indeed, the moduli stack of smooth curves of genus $g$
  is smooth of dimension $3g-3$, while for any given curve the moduli
  stack of $\PGL_2$-bundles is formally smooth of dimension $3g-3$, so that
  we obtain a smooth stack of dimension $6g-6$ when $K_X^2=8-8g$, with the
  remaining cases following as in the rational case.
\end{rems}

The action of simple reflections on blowdown structures clearly extends to
give birational automorphisms of these stacks (since each simple root is
clearly generically ineffective).  The action is undefined when the root is
effective, leading us to wonder what those substacks look like.  It turns
out that any given positive root (simple or not) is effective on a closed
substack of codimension 1.  This is a special case of a more general fact
about flat families of sheaves, which we will have occasion to use again.

\begin{lem}\label{lem:tau_func}
  Let $\pi:X\to S$ be a projective morphism of schemes, and suppose $M$ is a
  coherent sheaf on $X$, flat over $S$.  Suppose moreover that
  $R^p\pi_*M=0$ for $p>1$, and the fibers of $M$ have Euler characteristic
  $0$.  Then the locus $T\subset S$ parametrizing fibers with global
  sections has codimension $\le 1$ everywhere.  Moreover, where
  $T\subsetneq S$, it is a Cartier divisor.
\end{lem}

\begin{proof}
  By \cite{KnudsenFF/MumfordD:1976}, the derived direct image of $M$ can be
  represented by a perfect complex on $S$ starting in degree 0.  Since the
  higher direct images vanish, we can replace the degree $1$ term by the
  kernel of the map to the degree $2$ term to obtain a two-term perfect
  complex.  The Euler characteristic condition implies that the two terms
  have the same rank everywhere, and $R^1\pi_*M$ is supported on the zero
  locus of the determinant of the appropriate map, so has codimension
  $\le 1$.  Semicontinuity implies that the fibers of $M$ have global
  sections precisely along the support of $R^1\pi_*M$.  That this is a
  Cartier divisor follows from the construction of
  \cite{KnudsenFF/MumfordD:1976}: the determinant is the canonical global
  section of the bundle $\det \dR\Gamma(M)^{-1}$.  (Note that the reference
  only shows that the fibers have global sections on the zero locus of the
  canonical global section; the argument above gives the converse as well.)
\end{proof}

\begin{cor}
  For any positive (real) root of $W(E_{m+1})$, the corresponding divisor
  class is effective on a codimension $1$ substack of ${\cal X}_m$.
\end{cor}

\begin{proof}
  Indeed, a positive root has $D^2=-2$, $D\cdot C_\alpha=0$, and thus
  $\chi(\sO_X(D))=0$.  Since there exist surfaces for which no positive
  root is effective and surfaces for which every positive root is
  effective, the substack is nonempty, and not all of ${\cal X}_m$, so has
  codimension $1$.
\end{proof}

\begin{rem}
Similarly, that a surface in ${\cal X}_9$ admits {\em some} anticanonical
curve is a codimension $1$ condition.
\end{rem}

\begin{cor}
  The substack of rational surfaces containing a
  $-e$-curve for some $e\ge 3$ is a countable union of closed substacks of
  codimension $\ge 2$.
\end{cor}

\begin{proof}
  A $-d$-curve on $X_m$ for $d>2$ is either the strict transform of a
  $-d$-curve on $X_0$ or arises from a $-(d-1)$ curve on some $X_k$ for
  $k<m$.  We have already noted that $X_0$ contains $-3$-curves (or worse)
  in codimension $\ge 2$.  For the other case, we observe by induction
  (using the $d=2$ case above) that $X_k$ containing a $-2$-curve is a
  countable union of codimension $1$ conditions (since any effective root
  is a sum of $-2$-curves, and in the absence of $-3$-curves or worse, is
  uniquely effective).  If we blow up a point of an effective root, the
  result will necessarily contain a $-3$-curve, and any $-3$-curve not
  already present must arise in this way.
\end{proof}

\begin{rem}
  Similarly, for $d>3$, we have a $-d$-curve or worse on a countable union
  of locally closed substacks of codimension $\ge d-1$.  This fails to be
  closed since the $-2$-curve we are turning into a $-d$-curve could
  degenerate into a reducible effective root.  So, for instance, the closed
  substacks corresponding to $-4$-curves also contain configurations of
  self-intersection $-4$ consisting of two $-3$-curves connected by a chain
  of $-2$-curves.
\end{rem}

\medskip

If we try to extend the action of simple roots to the entire moduli space,
we encounter a problem along these codimension $\ge 2$ substacks.  For
simplicity in exhibiting the problem, we consider blowups of $\P^2$.
Define maps $p_1$, $p_2$, $p_3:\A^1\to \P^2$ by
\[
p_1(u) = (0,0),\qquad
p_2(u) = (u,0),\qquad
p_3(u) = (0,u^2),
\]
and define two family of surfaces parametrized by $\A^1$: $Y_u$ is the blowup
of $\P^2$ in $p_1$, $p_2$, then $p_3$, while $Z_u$ is the blowup of $\P^2$
in $p_3$, $p_2$, then $p_1$.  (At each step, the maps not already used
extend to the blowup, so give a well-defined point at which to blow up.)
For $u\ne 0$, we have $Y_u\cong Z_u$, since the points are distinct, so the
different blowups commute.  On the other hand $Y_0$ is the blowup of $F_1$
in two {\em distinct} points of the $-1$-curve, while $Z_0$ blows up the
same point of the $-1$-curve twice.  But then $Y_0\not\cong Z_0$, since
$Z_0$ contains a $-2$-curve, and $Y_0$ does not.

Since every rational surface locally looks like $\P^2$, we see that there
is no way to extend the action of $S_3$ on a three-fold blowup to include
the surfaces where the three-fold blowup introduces a $-3$-curve.
This is essentially the only difficulty, however.

\begin{thm}
  Let ${\cal X}^{{\ge}{-}2}_m$, $m\ge 2$, denote the stack parametrizing
  pairs $(X,\Gamma)$ where $X$ is a rational surface with $K_X^2=8-m$ not
  containing any $-d$-curves for $d>2$ and $\Gamma$ is a blowdown structure
  on $X$.  There is a natural action of the Coxeter group $W(E_{m+1})$ on
  this stack, which for every simple root is given by the usual action on
  blowdown structures, where the root is ineffective.
\end{thm}

\begin{proof}
  Since $m\ge 1$, we may use an elementary transformation to identify the
  two components of ${\cal X}^{{\ge}{-}2}_m$, and thus have both an odd and an
  even blowdown structure on $X$.

  Since $X$ has no $-d$-curves for $d>2$, an odd blowdown structure maps it
  to $F_1$, so that we can proceed on to $\P^2$.  Moreover, every
  infinitely near point that gets blown up as we proceed to $X$ is a jet.
  More precisely, the blowdown structure on $X$ is determined by (a) a
  union $F_{m+1}$ of jets on $\P^2$ and (b) a filtration $F_i$ of the
  structure sheaf of this union such that each quotient is the structure
  sheaf of a point.  When a given simple root $e_i-e_{i+1}$ is ineffective,
  the corresponding quotient $F_{i+1}/F_{i-1}$ is supported on two distinct
  points, and the reflection makes the other choice of $F_i$.  This extends
  immediately to the locus where the two points agree; the jet condition
  ensures that the degree 2 scheme parametrizing choices of $F_i$ is
  separable.  In particular, we find that the action of $S_{m+1}$ extends
  to the full stack.

  Similarly, an even blowdown structure corresponds to (a) a choice of
  $F_0$ or $F_2$, (b) a union of jets on this surface (disjoint from the
  $-2$-curve), and (c) a corresponding filtration.  The same argument tells
  us that the corresponding $S_m$ acts.  Since the two subgroups cover the
  full set of simple roots, we obtain
  the desired action of $W(E_{m+1})$.  (The braid relations hold because
  they hold generically.)
\end{proof}

\begin{rems}
  The stack ${\cal X}^{{\ge}{-}2}_m$ is not quite a substack of ${\cal
    X}_m$, in general, as it may be necessary to remove a countable
  infinity of closed substacks.  For instance, we can obtain a $-3$-curve
  in ${\cal X}_9$ by blowing up a point of any $-2$-curve, and thus each of
  the infinitely many positive roots of $E_{8+1}$ produces a different
  component we must remove.
\end{rems}

\begin{rems}
  Note that the simple reflections act trivially on the locus where the
  corresponding simple root is effective (thus a $-2$-curve).  Since that
  locus has codimension $1$, we see that the simple reflections act {\em as
    reflections} on ${\cal X}^{{\ge}{-}2}_m$.
\end{rems}

\begin{rems}
  Note that although the group acts on ${\cal X}^{{\ge}{-}2}_m$, this action
  {\em cannot} extend to the universal surface.  That is, the action
  preserves the isomorphism class of the surface, but the isomorphisms on
  generic fibers degenerate as we approach the bad fibers.  The problem
  here is that the universal surface is itself a moduli space of surfaces,
  but those surfaces could contain $-3$-curves; in other words, the generic
  isomorphisms degenerate precisely on the corresponding $-2$-curves on the
  fiber.
\end{rems}

\begin{rems}
  This also helps quantify the sense in which the moduli stack of surfaces
  is badly behaved when we do not introduce the blowdown structure: even if
  we exclude $-d$-curves for $d>2$, it is the quotient of an Artin stack by
  a discrete group which is infinite when $K_X^2\le 0$.
\end{rems}

\begin{rems}
  Finally, the reader should be cautioned that although we have divisors
  corresponding to the standard bases of $\Pic(X)$ (including $s$, since
  $s$ is represented by a canonical divisor on $X_0$ in the odd case),
  these choices of divisor are not compatible with the action of the
  Coxeter group.  (For instance, the reflection in $e_1-e_2$ changes the
  representation of $f$ from $(f-e_1)+e_1$ to $(f-e_2)+e_2$.)  In
  particular, although the various line bundles are taken to isomorphic
  bundles under the group action, those isomorphisms are not canonical.
\end{rems}

\section{Anticanonical surfaces}

When trying to extend the above construction to anticanonical surfaces, we
encounter the difficulty that the dimension of the anticanonical linear
system varies with the rational surface, and this variation depends in
subtle ways on the configuration of $-d$-curves on the surface with $d>2$.
It turns out that this is not too serious an obstruction to constructing
the linear system, however.

\begin{prop}\label{prop:linsys}
  Let ${\cal L}$ be a line bundle on a family $X$ of smooth projective
  varieties over a $d$-dimensional locally Noetherian integral base $S$,
  such that on every fiber, $H^p({\cal L})=0$ for $p\ge 2$.  Then the
  moduli functor $|{\cal L}|$ defined by taking $|{\cal L}|(T)$ to be the
  set of effective divisors on $X|_T$ with $\sO_X(D)\cong {\cal L}$ is
  represented by a locally projective scheme $|{\cal L}|$ over $S$, of
  dimension at least $d+\chi({\cal L})-1$ everywhere locally.  If the base
  is smooth and the dimension is equal to $d+\chi({\cal L})-1$, then
  $|{\cal L}|$ is a local complete intersection.
\end{prop}

\begin{proof}
  This certainly holds (and with equality for the dimension) when ${\cal
    L}$ is acyclic, since then Grauert's theorem tells us that the direct
  image of ${\cal L}$ is a vector bundle, and $|{\cal L}|$ is just the
  corresponding projective bundle.

  More generally, we may as well assume $S$ is Noetherian and affine.  In
  particular, there is an effective divisor $D$ such that ${\cal L}(D)$ is
  acyclic on every fiber, so that we may directly construct $|{\cal
    L}(D)|$.  The long exact sequence of cohomology associated to
  \[
  0\to {\cal L}\to {\cal L}(D)\to {\cal L}\otimes \sO_D(D)\to 0
  \]
  then tells us that for $p>0$, $H^p({\cal L}\otimes \sO_D(D))\cong
  H^{p+1}({\cal L})$ on every fiber, and thus that ${\cal L}\otimes
  \sO_D(D)$ is acyclic.  The dual vector bundle to $\Gamma({\cal L}\otimes
  \sO_D(D))$ then defines a system of equations on $|{\cal L}(D)|$, cutting
  out those sections containing $D$.  Call the resulting closed subscheme
  $Y$.  Since $Y$ is an intersection of $h^0({\cal L}\otimes \sO_D(D))$
  hypersurfaces, it has dimension at least
  \[
  \dim |{\cal L}(D)| - h^0({\cal L}\otimes \sO_D(D))
  =
  d-1+\chi({\cal L}(D))-\chi({\cal L}\otimes \sO_D(D))
  =
  d-1+\chi({\cal L})
  \]
  everywhere locally.  The restriction to $Y$ of the universal divisor on
  $|{\cal L}(D)|$ has $D$ as a fixed locus, which we may subtract to obtain
  a universal divisor for $|{\cal L}|$ as required.

  The local complete intersection property follows from the fact that $Y$
  is obtained from a smooth scheme by intersecting the same number of
  hypersurfaces as its codimension.
\end{proof}

\begin{rem}
  Since automorphisms of ${\cal L}$ act trivially on $|{\cal L}|$, we can
  apply this even when ${\cal L}$ is merely an isomorphism class of line
  bundles, with the one caveat being that this will make $|{\cal L}|$ a
  family of Brauer-Severi varieties.
\end{rem}

In particular, we may define the stack ${\cal X}^\alpha_m$ of anticanonical
surfaces to be the linear system $|-K|$ on the universal surface over
${\cal X}_m$.  Of course, we would like to know that this is irreducible,
which will require some more work.

For Hirzebruch surfaces, things are not too difficult to control, as in
each case we can write $-K=D_0+D_1$ where $D_0$ is the divisor of fixed
components and $D_1$ is acyclic:
\begin{align}
F_0&: -K = 0+(2s+2f)\notag\\
F_1&: -K = 0+(2s+3f)\notag\\
F_2&: -K = 0+(2s+2f)\notag\\
F_{2d+1}&: -K = (s-df)+(s+(d+3)f),\quad d\ge 1\notag\\
F_{2d} &: -K = (s-df)+(s+(d+2)f),\quad d\ge 2\notag
\end{align}
In particular, we find in each case that $|-K|$ is smooth over the given
locally closed substack, of relative dimension $\chi(D_1)$.  This has the
following curious effect: for all $k\ge 3$, the dimension of the stratum of
${\cal X}^\alpha_0$ corresponding to $F_k$ is $0$, while for $F_0$, $F_1$,
$F_2$ the dimension is $2$, $2$, $1$ respectively.  (Recall that ${\cal
  X}_0$ itself has dimension $-6$, since $F_0$ and $F_1$ have 6-dimensional
automorphism groups.)  Since ${\cal X}^\alpha_0$ has dimension at least $2$
everywhere, we find that it has precisely two irreducible components, as we
would have expected, both of which are l.c.i.

\begin{rems}
  This is already a distinct departure from the general case, since now all
  $-d$-curves for $d>2$ are a codimension $2$ phenomenon, not just the
  $-3$-curves.  It is clearer why this should be so for blowups: to obtain
  a $-d$ curve for $d>2$, we need simply blow up a point of a $-(d-1)$
  curve.  Since we already need to blow up a point of the anticanonical
  curve, this is a codimension 0 condition on the blowup!  Thus really the
  question is when the anticanonical curve is reducible, and this is
  codimension 2 (either the curve is reducible on $X_0$, or we must blow up
  a singular point).
\end{rems}

\begin{rems}
  Note that this moduli problem is not formally smooth.  For instance, if
  $X=F_4$ and $C_\alpha$ contains $s_{\min}$ with multiplicity $2$, then the
  anticanonical section $\alpha$ extends to an anticanonical section on an
  open subset of ${\cal X}_0$.  This exhibits a subspace of the tangent
  space to $(X,C_\alpha)$ as a direct sum of the tangent space to $X$ in
  ${\cal X}_0$ and the tangent space to $\alpha$ in
  $\P(H^0(\omega_X^{-1}))$.  But this subspace is larger than the generic
  tangent space!
\end{rems}

To understand ${\cal X}^\alpha_m$ in general, we may proceed by induction in $m$.

\begin{thm}\label{thm:antican_moduli_lci}
  The moduli problem of classifying triples $(X,C_\alpha,\Gamma$), where
  $X$ is a rational surface with $K_X^2=8-m$, $C_\alpha\subset X$ is an
  anticanonical curve, and $\Gamma$ is a blowdown structure, is represented
  by an Artin stack ${\cal X}^\alpha_m$.  This stack is a local complete
  intersection of dimension $m+2$, with one irreducible component for each
  parity of blowdown structure, with both components integral.  If $m\ge
  1$, the two components are canonically isomorphic.
\end{thm}

\begin{proof}
  We have already shown this for $m=0$, and we also note that the
  anticanonical curve is generically smooth in that case.  Now, the fibers
  of the natural forgetful map ${\cal X}^\alpha_m\to {\cal X}^\alpha_{m-1}$ are
  straightforward to determine: a point in the fiber just indicates which
  point of the anticanonical curve was blown up.  In other words, ${\cal
    X}^\alpha_m$ is the universal anticanonical curve over ${\cal
    X}^{\alpha}_{m-1}$.  By induction, the latter has two irreducible
  components, both integral, and the generic fiber of either component has
  smooth anticanonical curve.  In particular, we find that the fibers of
  ${\cal X}^\alpha_m$ over ${\cal X}^\alpha_{m-1}$ are all $1$-dimensional, and
  generically integral, so that each component of ${\cal X}^\alpha_{m-1}$
  has integral preimage.  Moreover, we immediately find that $\dim({\cal
    X}^\alpha_m)=m+\dim({\cal X}^\alpha_0)=m+2$ as required.

  That the components are isomorphic for $m>0$ follows by elementary
  transformation as before, and the local complete intersection property
  follows from the fact that
  $
  2m-6 + \chi(-K_X)-1 = m+2.
  $
\end{proof}

\begin{rem}
  Much of this statement is special to the rational case.  In particular,
  for genus 1 ruled surfaces, the typical anticanonical curve has to
  distinct components, and thus every time we blow up a point, we increase
  the number of components of the moduli stack.  For higher genus, not only
  do we typically have multiple anticanonical components, but we always
  have nonreduced components, and thus the moduli stack becomes nonreduced
  as soon as we blow up a point.  Both pathologies are compatible with
  being l.c.i., but the argument fails, as in either case, the dimension of
  the linear system is larger than the Euler characteristic would suggest.
\end{rem}

The above construction showing that $W(E_{m+1})$ cannot in general act in
the presence of $-3$-curves works equally well in the anticanonical case
(as long as the surface we start with has $K_X^2\ge 3$, so that there is an
anticanonical curve containing any three points).  Thus we introduce the
substack ${\cal X}^{\alpha,{\ge}{-}2}_m\subset {\cal X}^{\alpha}_m$ as
before, by excluding all $-d$-curves for $d>2$.  This is actually a
substack in this case, since as we noted above, once we bound the minimal
section of the Hirzebruch surface, there are only finitely many possible
configurations of components of the anticanonical curve.  Of course, not
having any $-d$-curves for $d>2$ is a very strong condition to impose on an
anticanonical surface: in particular, for $K_X^2<0$, it forces $C_\alpha$
to be integral.  (Indeed, otherwise the components of $C_\alpha$ are smooth
rational curves, at least one of which has negative intersection with
$C_\alpha$, so is a $-d$-curve for $d>2$.)

\begin{prop}
  For $0\le m<8$, ${\cal X}^{\alpha,{\ge}{-}2}_m$ is a $\P^{8-m}$-bundle
  over ${\cal X}^{{\ge}{-}2}_m$, and thus has two smooth components,
  isomorphic if $m>0$.  For $m>8$, ${\cal X}^{\alpha,{\ge}{-}2}_m$ can be
  identified with a closed substack of ${\cal X}^{{\ge}{-}2}_m$.
\end{prop}

\begin{proof}
  For $m>8$, the anticanonical curve is integral if it exists, and since it
  has negative self-intersection is rigid.  Thus there is at most one
  anticanonical curve on a rational surface without $-d$-curves for $d>2$.
  We have seen that this is a closed codimension 1 condition for $m=9$,
  while for $m>9$, it combines the closed conditions that the image in
  ${\cal X}^{{\ge}{-}2}_{m-1}$ is in ${\cal X}^{\alpha,{\ge}{-}2}_{m-1}$ and
  that the point being blown up lies on the anticanonical curve.

For $m<8$, the anticanonical divisor is nef on a surface without
$-d$-curves for $d>2$, and since $C_\alpha^2=8-m>0$, the corresponding line
bundle is acyclic.  But then the linear system is a $\P^{8-m}$-bundle as
required.
\end{proof}

\begin{rems}
  The cases with $m<8$ are nothing other than (degenerate) del Pezzo
  surfaces.  There is, of course, a lot one can say about those surfaces,
  but they are less useful for present purposes as they tend not to have
  sheaves disjoint from the anticanonical curve.  It is, however, worth
  noting that the action of $W(E_6)$ in the $m=5$ case is the familiar one
  on the 27 lines (i.e., $-1$-curves) on a cubic surface.
\end{rems}

\begin{rems}
The case $m=8$ is more subtle, as the surface could have a unique
anticanonical curve, or could have a $1$-parameter family of anticanonical
curves (making it an elliptic surface with no multiple fibers).
\end{rems}

In any case, since the choice of $C_\alpha$ is independent of the blowdown
structure, the action of $W(E_{m+1})$ extends immediately to ${\cal
  X}^{\alpha,{\ge}{-}2}_m$.  As before, the action does not extend to the
universal surface (it must act linearly on $\Pic(X)$, so does not respect
the effective cone), but it turns out that there is a strong sense in which
it {\em does} act on the {\em line bundles} on the universal surface.

Given any vector $v\in \Z s+\Z f+\sum_i \Z e_i$, we have a corresponding
line bundle ${\cal L}_v$ on the universal surface over ${\cal
  X}^{\alpha,{\ge}{-}2}_m$.  (In general, we only knew this when the
coefficient of $s$ was even, but the assumptions imply that an odd blowdown
structure reaches $F_1$, where $s$ is canonically a divisor, and the claim
for even blowdown structures follows by elementary transformation.)  Of
course, the space of global sections of this bundle can vary wildly with
the surface, and can similarly vary if we replace $v$ by $wv$ for any
element $w\in W(E_{m+1})$.  These are essentially the same phenomenon,
however.  The main problem with the global sections of ${\cal L}_v$ is that
we can have sections of ${\cal L}_v$ on a given surface that do not extend
to neighboring surfaces in the moduli space.  This can be fixed by taking
the direct image sheaf rather than the fiberwise global sections.  This can
cause problems in general, however, which are characterized by the
following result (a strong (albeit specialized) form of semicontinuity).

\begin{lem}
Let $\pi:X\to S$ be a projective morphism, and suppose that $M$ is a
sheaf on $X$, flat over $S$, such that every fiber of $M$ has $H^p=0$ for
$p\ge 2$.  Then for any sheaf $N$ on $S$, we have isomorphisms
\[
\Tor_{p+2}(R^1\pi_*M,N)\cong \Tor_p(\pi_*M,N)
\]
for $p>0$, along with a short exact sequence
\[
0\to \Tor_2(R^1\pi_*M,N)\to \pi_*M\otimes N\to \pi_*(M\otimes \pi^*N)\to
\Tor_1(R^1\pi_*M,N)\to 0
\]
and an isomorphism
\[
R^1\pi_*M\otimes N\cong R^1\pi_*(M\otimes \pi^*N).
\]
In particular, $\pi_*M$ is flat iff $R^1\pi_*M$ has homological dimension
$\le 2$, the fibers of $\pi_*M$ inject in the corresponding spaces of
global sections of $M$ iff $R^1\pi_*M$ has homological dimension $\le 1$,
and the injection is an isomorphism iff $R^1\pi_*M$ is flat.
\end{lem}

\begin{proof}
  As in the proof of Lemma \ref{lem:tau_func}, we find that $\dR\pi_*M$ is
  represented by a two-term perfect complex on $S$.  If $V^0\to V^1$ is
  this complex, then we have a four-term exact sequence
\[
0\to \pi_*M\to V^0\to V^1\to R^1\pi_*M\to 0,
\]
and thus any flat resolution of $\pi_*M$ extends to a flat resolution of
$R^1\pi_*M$.  The claim follows upon tensoring this resolution with $N$ and
observing that
\[
\dR\pi_*M\otimes^\dL N\cong \dR\pi_*(M\otimes^\dL\pi^*N)\cong
\dR\pi_*(M\otimes \pi^*N),
\]
with the last isomorphism following from the fact that $M$ is flat.
\end{proof}

\begin{rem}
  As an example, consider the divisor class $-2K_X$ on the moduli stack of
  anticanonical Hirzebruch surfaces.  This acquires cohomology when the
  surface has a $-d$-curve for any $d\ge 3$; since this locus has
  codimension $2$, $R^1\pi_*(\omega_X^{-2})$ has homological dimension $\ge
  2$.  As a result, the fibers of the direct image sheaf do not inject in
  the spaces of global sections of the fibers.  Similarly, the
  anticanonical bundle itself fails this criterion in the presence of a
  $-d$-curve for $d\ge 4$.
\end{rem}

  By the last claim of the Lemma, we can compute $R^1\pi_*M$ fiberwise, and
  the main contribution comes from hypersurfaces: those where a given
  positive root becomes effective, and those where ${\cal L}_v|_{C_\alpha}$
  has a global section.  Near a generic point of such a hypersurface, we
  find that $\pi_*{\cal L}_v$ is flat and injects fiberwise in the space of
  global sections of ${\cal L}_v$, since $R^1\pi_*{\cal L}_v$ is a flat
  sheaf on a hypersurface, so has homological dimension $\le 1$.  Although
  there could in principle be problems coming from intersections of the
  hypersurfaces, this at least suggests the following result; note that by
  the previous remark, we cannot allow worse than $-2$-curves, even if we
  did not care about the $W(E_{m+1})$ action.

\begin{thm}
  The direct image of any line bundle ${\cal L}_v$ is a flat sheaf ${\cal
    V}_v$ on ${\cal X}^{\alpha,{\ge}{-}2}_m$, and the action of $W(E_{m+1})$
  extends to these sheaves.  More precisely, for any element $w\in
  W(E_{m+1})$, we have an isomorphism
\[
w^*{\cal V}_v\cong {\cal V}_{wv},
\]
defined up to scalar multiplication, and the isomorphisms are compatible,
again up to scalar multiplication.  Moreover, the multiplication map
\[
{\cal V}_v\times {\cal V}_{v'}\to {\cal V}_v\otimes {\cal V}_{v'}\to {\cal
  V}_{v+v'}
\]
induced by
\[
{\cal L}_v\otimes {\cal L}_{v'}\cong {\cal L}_{v+v'}
\]
has no zero divisors.
\end{thm}

\begin{proof}
  First note that if $v$ is not generically effective, then ${\cal V}_v=0$,
  since then no global section of ${\cal L}_v$ on a fiber can extend to an
  open substack of the moduli space.  The generically effective divisors
  form a cone invariant under the action of $W(E_{m+1})$, so the various
  claims are immediate outside this cone.  A generically effective divisor
  will have $v\cdot f\ge 0$, so $(-C_\alpha-v)\cdot f\le -2$, and thus
  $-C_\alpha-v$ cannot be effective.  We thus conclude that $H^2({\cal
    L}_v)=0$ for such a divisor, which is all we need to apply the Lemma.

  It will suffice to show that whenever $v$ is generically effective, the
  group acts and $R^1\pi_*{\cal L}_v$ has homological dimension $\le 1$.
  Indeed, this implies flatness of ${\cal V}_v$, as well as the fact that
  multiplication has no zero-divisors, the latter since the map
\[
\Gamma({\cal L}_v)\times \Gamma({\cal L}_{v'})
\to
\Gamma({\cal L}_{v+v'})
\]
is injective on every fiber.

Now, using an elementary transformation as necessary, we may suppose our
blowdown structure is odd and consider $X$ as an $m+1$-fold blowup of
$\P^2$.  In the corresponding basis of $\Pic(X)$, we have
\[
v = nh - \sum_{0\le i\le m} r_i e_i.
\]
If $r_i=v\cdot e_i<0$ for any $i$, then we have a short exact sequence
\[
0\to {\cal L}_{v-e_i}\to {\cal L}_v\to {\cal L}_v|_{e_i}\to 0.
\]
On the generic fiber, the quotient is a sheaf of negative degree on the
smooth rational curve $e_i$, and thus has no global sections; on the
general fiber, the quotient has $1$-dimensional support, so that the Lemma
applies.  We thus find that $\pi_*({\cal L}_v|_{e_i})=0$ and (since that
certainly injects!) that $R^1\pi_*({\cal L}_v|_{e_i})$ has homological
dimension $\le 1$.  It follows that
\[
{\cal V}_{v-e_i}\cong {\cal V}_v,
\]
and $R^1\pi_*(\sO_X(v))$ has homological dimension $\le 1$ iff
$R^1\pi_*(\sO_X(v-e_i))$ has homological dimension $\le 1$.  By
induction, if we set
\[
v' = nh - \sum_{0\le i\le m} \max(r_i,0)e_i,
\]
then (since this operation respects the action of $S_{m+1}$) it suffices to
prove the claim for $v'$.

Thus suppose $r_i\ge 0$ for $0\le i\le m$, and consider the short exact
sequence
\[
0\to {\cal L}_v\to {\cal L}_{nh}\to Q\to 0.
\]
Since ${\cal L}_{nh}$ is acyclic, ${\cal V}_{nh}$ is flat, and we have
exhibited ${\cal V}_v$ as a subsheaf of this flat sheaf.  Moreover,
the fibers ${\cal V}_v$ inject in $\Gamma({\cal L}_v)$ iff they inject in
the corresponding fibers of ${\cal V}_{nh}$.

Now, $\pi_*Q$ is the kernel of a two-term perfect complex (since $Q$ has
$1$-dimensional support), and is thus a subsheaf of a locally free sheaf.
In particular, $\pi_*Q$ is torsion-free, and thus the map ${\cal V}_v\to
{\cal V}_{nh}$ is determined by its action on the generic fiber.  This
action is clearly $S_{m+1}$-covariant, and thus so is ${\cal V}_v$.
Since the morphism ${\cal V}_v\to {\cal V}_{nh}$ determines the injectivity
condition, we also conclude that if $R^1\pi_*{\cal L}_v$ has
homological dimension $1$, then so does $R^1\pi_*{\cal L}_{w v}$ for any
$w\in S_{m+1}$.

A similar calculation with an even blowdown structure shows that the
corresponding $S_m$ acts on the bundles, and preserves the homological
dimension condition.  The one technicality is that the ambient bundle
${\cal L}_{ns+df}$ need not be acyclic, but we can use $S_{m+1}$-invariance
(conjugated by an elementary transformation) to assume $n\ge d$.

We thus now have full $W(E_{m+1})$-covariance, so that it suffices to prove
the homological dimension claim for $v$ in the fundamental chamber.  Of
course, if $v=0$, then ${\cal L}_0=\sO_X$, and the claim is obvious, so
suppose $v\ne 0$.  If $v\cdot C_\alpha>0$, then ${\cal L}_v$ is acyclic, and
we are done.  Otherwise, consider the short exact sequence
\[
0\to {\cal L}_{v+K}\to {\cal L}_v\to {\cal L}_v|_{C_\alpha} \to 0
\]
Generically, ${\cal L}_v|_{C_\alpha}$ is a nontrivial degree 0 sheaf on a
smooth genus 1 curve, and thus we again find
\[
\pi_*({\cal L}_v|_{C_\alpha})=0
\]
and thus ${\cal L}_v$ satisfies the homological dimension condition iff
${\cal L}_{v-C_\alpha}$ satisfies the homological dimension condition.
\end{proof}

We should note a couple of things here.  First, the argument shows that on
a blowup of $\P^2$, ${\cal V}_v\subset {\cal V}_{(v\cdot h)h}$ whenever
$v\cdot h>0$, with locally free quotient, and similarly ${\cal V}_v\subset
{\cal V}_{(v\cdot f)s+(v\cdot s)f}$ relative to an even blowdown structure.
This can be generalized to the noncommutative case, by first constructing
suitable noncommutative analogues of the ambient bundles (as spaces of
symmetric difference operators), imposing suitable conditions on the
generic fiber, and then using an analogue of the above argument to prove
flatness.  Of course, this only works when the anticanonical curve is
integral, but can be made largely explicit for smooth anticanonical
curve, see \cite{generic}.

Next, the result allows us to construct a flat family of categories with a
nice action of $W_{E_{m+1}}$.  The objects of the categories are the
vectors $v\in \Z s+\Z f+\sum_i \Z e_i$, while the morphisms from $v$ to
$v'$ are given by ${\cal V}_{v'-v}$, with the natural multiplication maps.
The dimensions of the $\Hom$ spaces in this category are constant as we
vary the choice of anticanonical surface, and the group acts in the obvious
way.  Again, this deforms to the noncommutative case as an equally flat
family of categories, see \cite{generic}, with one additional parameter (a
point of $\Pic^0(C_\alpha)$) controlling the noncommutative deformation.

In fact, one can generalize even further: if one replaces the $\Hom$ spaces
by their symmetric powers, there is an additional deformation parameter
(again a point of $\Pic^0(C_\alpha)$), thus giving flat deformations of
symmetric powers of surfaces with smooth anticanonical curve.  (See
\cite{elldaha}, which also gives a conjectural extension to a deformation
of Hilbert schemes of points.)

Unfortunately, since this fails when the anticanonical curve is
nonintegral, we will not be able to use this approach in the general
noncommutative setting.

\medskip

In the case $C_\alpha$ integral, we can give a direct construction of the
substack of surfaces with anticanonical curve isomorphic to $C_\alpha$.
(Presumably this can be extended to general curves, but it is unclear what
the precise conditions will be.)

For any connected projective curve $C$ (integral or not) of arithmetic
genus 1, there is a natural moduli problem mapping flatly to the locally
closed substack of ${\cal X}^{\alpha}_m$ where $C_\alpha\cong C$, namely the
problem of classifying triples $(X,\phi,\Gamma)$ where $\phi:C\to X$ embeds
$C$ as an anticanonical curve.  Given such a triple, the restriction
morphism $\phi^*:\Pic(X)\to \Pic(C)$ gives us a sequence of (isomorphism
classes of) bundles $\phi^*(s)$, $\phi^*(f)$ and $\phi^*(e_i)$ for $1\le
i\le m$.  The classes $\phi^*(e_i)$ have degree 1, $\phi^*(f)$ has degree
2, and $\phi^*(s)$ has degree $1$ or $2$ depending on whether $\Gamma$ is
odd or even.

\begin{lem}
  The triple $(X,\phi,\Gamma)$ is determined up to isomorphism of $X$ by
  the classes $\phi^*(s)$, $\phi^*(f)$, and $\phi^*(e_i)$, $1\le i\le m$.
\end{lem}

\begin{proof}
  We may assume $X$ and $C$ are defined over an algebraically closed field,
  so that the given classes actually correspond to line bundles.  Since
  $K_X+f$ and $-f$ are both ineffective, we conclude that
  $H^0(\omega_X(f))=H^2(\omega_X(f))=0$, and Euler characteristic
  considerations imply $H^1(\omega_X(f))=0$.  It follows that we have a
  natural isomorphism $\Gamma(X,\sO_X(f))\cong \Gamma(C,\phi^*(f))$, so that
  we can recover the induced map $\rho:C\to \P^1$ from $\phi^*(f)$ (up to
  $\PGL_2$).  Similarly, $\dR\rho_*(\omega_X\otimes {\cal L}_s)=0$, so that
  we have an isomorphism $\rho_*{\cal L}_s\cong \rho_*\phi^*(s)$ allowing
  us to recover $X_0$ as the projective bundle of $\rho_*(\phi^*(s))$, as
  expected.

  Now, suppose we have reconstructed the surface $X_k$, and consider the
  direct image of $\phi^*(e_{k+1})$ on the corresponding anticanonical
  curve $C_k$.  This can be identified with the direct image of the sheaf
  $\sO_X(e_{k+1})|_{C_\alpha}$ on $X$, and thus by Corollary
  \ref{cor:pseudo_twist_is_isospectral} fits into an exact sequence
\[
0\to \sO_{C_k}\to \phi_{k*}\phi^*(e_{k+1})\to \sO_{p_{k+1}}\to 0
\]
where $p_{k+1}$ is the point of $C_k$ that gets blown up on $X_{k+1}$.
We have
\[
\Gamma(X_k,\phi_{k*}\phi^*(e_{k+1}))
=
\Gamma(X_{k+1},\sO_{X_{k+1}}(e_{k+1})|_{C_\alpha}).
\]
Since $e_{k+1}$ is a $-1$-curve on $X_{k+1}$, we find that
$\sO_{X_{k+1}}(e_{k+1})$ is acyclic and uniquely effective, while
\[
H^p(\sO_{x_{k+1}}(e_{k+1}-C_\alpha)) \cong
H^{2-p}(\sO_{x_{k+1}}(-e_{k+1}))^* = 0.
\]
It follows that $\phi_{k*}\phi^*(e_{k+1})$ is uniquely effective, so
determines $p_{k+1}$ and thus $X_{k+1}$.
\end{proof}

For $C$ integral, we readily see that any sequence of invertible sheaves of
the correct degrees will give rise to a valid triple, and thus we can
identify the moduli space of triples with the product
\[
\Pic^2(C)\times \Pic^2(C)\times \Pic^1(C)^m
\qquad
\text{or}
\qquad
\Pic^1(C)\times \Pic^2(C)\times \Pic^1(C)^m,
\]
depending on parity.  This is a principal $\Aut(C)$-bundle over the
corresponding substack of ${\cal X}^{\alpha}_m$, and since the choice of
embedding of $C$ is independent of the choice of blowdown structure, the
action of $W(E_{m+1})$ extends.  Of course, this extension is just the
restriction of the obvious linear action on $\Pic(C)^{m+2}$!  Note also
that if we allow $C$ to vary over the moduli stack of smooth genus 1
curves, then this construction gives a dense open substack of ${\cal
  X}^\alpha_m$.

For nonintegral curves, there will certainly be an additional constraint on
the degree vectors of the invertible sheaves, since those degrees can be
read off from the combinatorial type of $(X,\Gamma)$ (i.e., the
representations of the components of $C_\alpha$ in the standard basis).  And,
of course, there is the additional difficulty that the moduli stack of
curves has many more pathologies once one allows nonintegral curves,
especially since we do not want to impose any stability conditions.  (For
instance, there are reduced but reducible curves of arithmetic genus 1 that
are not even Gorenstein.)

Along these lines, we note the following constraint on the anticanonical
curve of a rational surface, for later use in the noncommutative case
(which admits the same anticanonical curves).

\begin{lem}\label{lem:Ca_is_num_conn}
  Let $X$ be a rational surface, and suppose we can write $-K_X=A+B$ with
  $A$, $B$ nonzero effective divisors.  Then $H^1(\sO_A)=H^1(\sO_B)=0$ and
  $A\cdot B=2h^0(\sO_A)=2h^0(\sO_B)>0$.
\end{lem}

\begin{proof}
  Since $A$ and $B$ are nonzero effective divisors, we find that
  $H^0(\sO_X(-A))=H^0(\sO_X(-B))=0$, and thus by duality
  $H^2(\sO_X(-A))=H^2(\sO_X(K_X+B))$=0.  The long exact sequence associated to
  the natural presentation of $\sO_A$ has a piece
\[
H^1(\sO_X)\to H^1(\sO_A)\to H^2(\sO_X(-A)),
\]
and thus $H^1(\sO_A)=0$, with $H^1(\sO_B)=0$ following similarly.  We thus
have $h^0(\sO_A)=\chi(\sO_A)=\frac{1}{2}A\cdot (-K_X-A)=\frac{1}{2}A\cdot
B$, so the remaining claim follows.
\end{proof}

\begin{rem}
We also find $h^1(\sO_X(-A))=-\chi(\sO_X(-A))=\frac{1}{2}A\cdot B-1$.
\end{rem}

This has the following interesting consequence; this was established for
anticanonical curves in $\P^2$ and $\P^1\times \P^1$ in
\cite[Cor.~5.7]{ArtinM/TateJ/VandenBerghM:1991}, but given the above Lemma, the
proof carries over directly.

\begin{prop}\label{prop:Pic0_acts}
  Let $(X,C_\alpha)$ be an anticanonical rational surface.  Then there is a
  natural action of the group scheme $\Pic^0(C_\alpha)$ on $C_\alpha$,
  which for an invertible sheaf ${\cal Q}\in \Pic^0(C_\alpha)$ takes a
  point $p\in C_\alpha$ to the unique point $p'$ such that ${\cal
    I}_{p'}\cong {\cal Q}\otimes {\cal I}_p$.
\end{prop}

\begin{rem}
  Similarly, the proof of Proposition 5.9 and Lemma 5.10 op.~cit.~tells us
  that for any ${\cal Q}\in \Pic^0(C_\alpha)$, the corresponding action
  $\tau_{\cal Q}$ fixes every singular point of $C_\alpha$, stabilizes any
  irreducible component of $C_\alpha$, and for any invertible sheaf ${\cal
    L}$ one has
\[
\tau_{\cal Q}^* {\cal L}\cong {\cal L}\otimes {\cal Q}^{-\chi({\cal L})}.
\]
The last part of Proposition 5.9 op.~cit.~suggests that if the component
$C$ occurs with multiplicity $m$, then $\tau_{\cal Q}$ restricts to the
identity on $(m-1)C$.
\end{rem}

\section{Partitioning the moduli stack}
\label{sec:types_comm}

Since the type and singularities of a difference or differential equation
depends on how the anticanonical curve decomposes and interacts with the
map to $\P^1$, we would like to understand the corresponding decomposition
of ${\cal X}^\alpha_m$.

The simplest part of the decomposition is by whether the curve is smooth or
of multiplicative or additive degeneration.  More precisely, we define the
``Picard type'' of an anticanonical rational surface to be one of the
symbols $e$, $*$, or $+$, depending on whether $\Pic^0(C_\alpha)$ is an
elliptic curve, multiplicative group, or additive group.

\begin{lem}
  The anticanonical rational surfaces of Picard type $*$ or $+$ (with
  blowdown structure) form an irreducible codimension $1$ closed substack
  ${\cal X}^{\alpha,\bar{*}}_m$ of ${\cal X}^\alpha_m$.
\end{lem}

\begin{proof}
  If $m>0$, then we may use the line bundle $\sO_{C_\alpha}(e_m)$ to induce
  a map from $C_\alpha$ to a weighted projective space with generators of
  degrees $1$, $2$, $3$, and the image of $C_\alpha$ will be the
  Weierstrass model of $\Pic^0(C_\alpha)$; then ${\cal
    X}^{\alpha,\bar{*}}_m$ is cut out by the vanishing of the discriminant:
  an invariant hypersurface of degree 12.

  For $m\le 0$, we need merely observe that if we blow up a point of
  $C_\alpha$, this will not change the Picard type, and thus the condition
  remains closed of codimension 1.
\end{proof}

\begin{rem}
  We could also use the nef divisor $f$ on $X_0$ or $h$ on $\P^2$ and
  reduce to known results on hyperelliptic curves of genus 1 or cubic
  plane curves.  The Weierstrass case is particularly nice over $\Z[1/6]$,
  however, since then we obtain well-defined functions $a_4$, $a_6$ on the
  $G_m$-bundle over ${\cal X}^\alpha_m$ in which we have chosen a Poisson
  structure on $X$ (since that corresponds to a choice of nonzero
  holomorphic differential).  Then ${\cal X}^{\alpha,\bar{*}}_m$ is cut out in
  those coordinates by the equation $-64 a_4^3-432 a_6^2=0$.
\end{rem}

The above works over $\Z$, but the additive case behaves differently in
characteristic $2$ and $3$, and thus we have a slightly weaker result.

\begin{lem}
  Over an algebraically closed field or a $\Z[1/6]$-algebra, the
  anticanonical rational surfaces of Picard type $+$ (with blowdown
  structure) form an irreducible codimension $1$ closed substack
  ${\cal X}^{\alpha,+}_m$ of ${\cal X}^{\alpha,\bar{*}}_m$.
\end{lem}

\begin{rem}
  Over $\Z[1/6]$, the additive Weierstrass curves are cut out from the full
  stack of Weierstrass curves by the equations $a_4=a_6=0$, but over $\F_3$
  the equations have degrees $2$, $12$ and over $\F_2$ they have degrees
  $1$, $12$; the same holds for hyperelliptic and cubic models.
\end{rem}

One issue that arises with the stack ${\cal X}^{\alpha,\bar{*}}_m$ is that
it does not distinguish $0$ and $\infty$.  For instance, in the irreducible
case, there are two branches of $C_\alpha$ at the node, and the limit as we
approach the two points is different.  We can resolve this by taking a
suitable double cover of the moduli space.  Indeed, the Weierstrass model
has a unique singular point, which we may take to be $(0,0)$, giving a
curve of the form
\[
y^2 + a_1 xy = x^3+a_2 x^2.
\]
The two branches at $\infty$ are given by the roots of $b_1^2+a_1 b_1-a_2$,
and thus the double cover has a model of the form
\[
y^2 + c_1 xy = x^3,
\]
where $c_1=a_1+2b_1$ is uniquely determined if we have chosen a nonzero
holomorphic differential on $C_\alpha$.  The double cover is ramified along
${\cal X}^{\alpha,+}_m$, and the equation for the corresponding substack is
$c_1=0$, now valid over $\Z$.  If $C_\alpha$ has multiplicative type with
multiple components, then the double cover distinguishes the two branches
at any given one of the nodes, and that identification can then be carried
along the polygon.

We similarly define the ``combinatorial type'' of an anticanonical rational
surface with blowdown structure to be the (multi)set of pairs $(C_i,m_i)\in
\Pic(X_m)\times \Z$, where $C_i$ ranges over the components of the
anticanonical curve and $m_i$ is the corresponding multiplicity.  Again, if
we have two components with the same divisor class and multiplicity, we
will usually want to distinguish them.  We thus choose an ordering of the
components, which gives us a somewhat different way to describe the type:
if we have $c$ components, then we have a morphism $\phi:\Z^c\to \Pic(X_m)$
taking $e_i$ to $C_i$, and a vector $\mu\in \Z^c$ giving the
multiplicities, with the property that $\phi\mu = C_\alpha$.

The ``type'' of a rational surface is then defined to be the symbol
$(c,\phi,\mu)_\sigma$ where $(c,\phi,\mu)$ gives the combinatorial type and
$\sigma$ gives the Picard type.  Note that most of the time the Picard type
is determined by the combinatorial type: if $\max_i\mu_i>1$, then the
Picard type is $+$, while if $c>3$ and $\max_i\mu_i=1$, then the Picard
type is $*$.  So only the cases with $c\le 3$ and no multiplicities have
any ambiguity.  (And, of course, the only case with Picard type $e$ is
$(1,(-K),(1))_e$.)

It is, in principle, straightforward to determine the set of possible types
for any given $m$, as it is easy to write down a complete (though countably
infinite) list for $m=0$, and each point being blown up either lies on a
single component or on a set of mutually intersecting components.  One also
has the following.  Define the ``dimension'' of a type to be
\[
\dim((c,\phi,\mu)_\sigma)
=
\begin{cases}
  m+2 & \text{$c=1$, $\sigma=e$}\\
  m+1 & \text{$c=1$, $\sigma=*$}\\
  m & \text{$c=1$, $\sigma=+$}\\
  m-1 & \text{$c=2$, $\mu=(11)$, $\sigma=+$}\\
  m-2 & \text{$c=3$, $\mu=(111)$, $\sigma=+$}\\
  2m-6 + \sum_{1\le i\le c} \phi_i\cdot (\phi_i-K)/2, & \text{otherwise.}
\end{cases}
\]

\begin{lem}
  For any type $(c,\phi,\mu)_\sigma$, the anticanonical rational surfaces
  of that type form an irreducible Artin stack of dimension
  $\dim((c,\phi,\mu)_\sigma)$, which for $\sigma\ne e$ maps to a locally
  closed substack of the double cover of ${\cal X}^{\alpha,\bar{*},m}$.
\end{lem}

\begin{proof}
  We proceed by induction on $m$.  For $m=-1$ (i.e., $\P^2$), the result is
  straightforward: since we have labelled the components, we may construct
  the stack inside the product of the linear systems on $\P^2$, and the
  various conditions we want to impose or exclude are closed, with known
  codimension; subtracting $8=\dim\Aut(\P^2)$ then gives the desired result
  for the dimension.  For most cases with $m=0$, we can proceed similarly;
  when the combinatorial type includes a section of nonpositive
  self-intersection $-d$, the surface must be $F_d$, while if it otherwise
  contains a component having negative intersection with $s-f$, it must be
  $F_0$ or $F_1$.  The only remaining cases are
  \[
  (1,(2s+2f),(1))_{e/*/+},\quad\text{and}\quad
  (2,(s+f,s+f),(1,1))_{*/+}.
 \notag
  \]
  The anticanonical curve on such a surface has a natural hyperelliptic
  model, and the surface itself is determined by this model together with a
  class in $\Pic^0$ (which is trivial iff the surface is $F_2$), and thus
  the stack has dimension 1 more than the corresponding stack of
  hyperelliptic curves of genus 1 (which is 1-dimensional in the first
  case, $-1$-dimensional in the second).

  For $m>0$, the type of $X_m$ determines the type of $X_{m-1}$ as well as
  the set of components the point being blown up lies on.  Irreducibility
  is then almost immediate, as the fibers of the stack parametrizing $X_m$
  over the stack parametrizing $X_{m-1}$ are either
  single points or open subsets of curves.  (Here we use the fact that we
  have chosen branches at the singular points, and thus when blowing up a
  point on a surface of type $(2,\phi,(1,1))_*$, we can distinguish the two
  points.)

  It remains only to show that change in the dimension of the type is the
  same as the change in dimension of the classifying stack.  If neither
  type falls into one of the exceptions in the dimension formula, then the
  change in the dimension of the type is indeed $1$ or $0$ depending on
  whether the point being blown up is on $1$ or $2$ components.  It is then
  straightforward to deal with the exceptions case-by-case.
\end{proof}

In particular, this tells us how many parameters a given type has, modulo
the effect of automorphisms.

Of course, simply knowing the possible types of degenerations of surfaces
(or of difference equations) is only part of the story: in general, we
would like to understand the limiting relations between the different
types.  That is, we would like to know which types appear in the closure in
${\cal X}^\alpha_m$ of the locally closed substack corresponding to a
particular type.  Note that this is somewhat more than just a poset, as a
given type may appear in multiple ways in the closure whenever there are
ambiguities in the labelling.  In particular, we expect the answer to be an
ordered category rather than simply a poset.  (There is a further issue, in
that there are cases in small characteristic in which the intersection of
two closures of types is not a union of types, see below.)

We have not been able to answer this question completely, but do have a
couple of necessary conditions, the simpler of which is as follows.  Define
a (na\"{i}ve) category structure on the set of types as follows.  Let
$(c_1,\phi_1,\mu_1)_{\sigma_1}$, $(c_2,\phi_2,\mu_2)_{\sigma_2}$ be a pair
of types.  If $\sigma_1>\sigma_2$ relative to the order $e>*>+$, then there
are no morphisms between the types, while otherwise a morphism is given by
a linear transformation $\psi:\Z^{c_2}\to \Z^{c_1}$ with nonnegative
coefficients such that $\phi_2=\phi_1\circ \psi$, $\psi(\mu_2)=\mu_1$.
Note that any endomorphism in this category is an automorphism, so this is
indeed an ordered category.

\begin{lem}
  A flat morphism from a dvr to ${\cal X}^\alpha_m$ induces a morphism from
  the type of the (surface corresponding to the) special fiber to the type
  of the generic fiber.
\end{lem}

\begin{proof}
  This is certainly true for the Picard type, so it suffices to consider
  the combinatorial type.  Making a quasi-finite flat base change as
  necessary, we may assume that the irreducible components of the
  anticanonical curve on the generic fiber are geometrically irreducible,
  and then choose an ordering on the components to induce a well-defined
  combinatorial type.  The corresponding family of anticanonical curves
  defines a divisor on the family of surfaces, not containing any fiber.
  As a result, the same applies to each irreducible component of that
  divisor.  Those are in one-to-one correspondence with the components of
  the generic anticanonical curve, and thus (by restriction) induces a
  divisor on the special fiber corresponding to each component of the
  generic anticanonical curve.  Each such divisor is contained in the
  anticanonical curve on the special fiber, and is thus a sum of geometric
  components of that curve.  The linear transformation such that $\psi_i$ is
  the linear combination of components corresponding to the restriction of
  $C_i$ gives a morphism of combinatorial types as required.
\end{proof}

We call a morphism in the category of types ``effective'' if it arises from
a family over a dvr, and ``strongly effective'' if there is a family over a
dvr such that the special fiber is the {\em generic} surface of the given
type.  It is easy to see that if $\psi_1$, $\psi_2$ are morphisms and
$\psi_1$ is strongly effective, then $\psi_2\circ \psi_1$ is (strongly)
effective iff $\psi_2$ is (strongly) effective.

The one technique we have (other than producing an explicit deformation)
for proving strong effectiveness is the following.

\begin{prop}
  Let $\psi:T_1\to T_2$ be a morphism such that for any factorization
  \begin{align}
  \begin{CD}
  T_1 @>\psi_1>> T' @>\psi_2>> T_2
  \end{CD}
  \notag
  \end{align}
  of $\psi$ with $T'\not\cong T_2$, $\dim(T')<\dim(T_2)$.  Then $\psi$ is
  strongly effective.
\end{prop}

\begin{proof}
  Let $T_2=(c,\psi,\mu)$, and let $Y$ be the fiber product $\prod_i
  |\psi_i|$ of linear systems on the universal surface over ${\cal X}_m$.
  There is a natural map from $Y$ to ${\cal X}^\alpha_m$ given by taking
  the anticanonical curve to be $\sum_i \mu_i D_i$ where $D_i$ is the fiber
  of $|\psi_i|$.  It follows from Proposition \ref{prop:linsys} that $Y$
  has dimension at least $2m-6 + \sum_{1\le i\le c} \phi_i\cdot
  (\phi_i-K)/2$.  Moreover, if the combinatorial type does not force the
  Picard type, then imposing the corresponding condition introduces one or
  two more hypersurfaces as appropriate.  We thus obtain a stack $Y_{T_2}$
  of dimension at least $\dim(T_2)$ everywhere locally.

  Now, it is easy to see that this construction defines a functor on the
  category of types, and thus we have a morphism $Y_{T_1}\to Y_{T'}\to
  Y_{T_2}$ for any factorization as described.  Moreover, the moduli stack
  of surfaces of type precisely $T$ embeds as an open substack of
  $Y_T$ for each $T$.

  Now, if $\psi$ were not strongly effective, then the substack
  corresponding to $T_1$ would necessarily meet some component of $Y_{T_2}$
  not containing the substack corresponding to $T_2$.  Any such component
  has dimension at least $\dim(T_2)$, and thus the type $T'$ of its generic
  fiber has at least that dimension.  But the existence of surfaces of type
  $T'$ in $Y_{T_2}$ implies that there is a morphism $T'\to T_2$, and we
  thus obtain a factorization of precisely the form we excluded.
\end{proof}

The following special cases are fairly straightforward, if occasionally
tedious; we omit the details.

\begin{cor}
  Any morphism between multiplicative types is strongly effective.
\end{cor}

\begin{cor}
  Any morphism to a type $(c,\phi,\mu)$ with $c\le 3$, $\max_i \mu_i=1$ is
  strongly effective.
\end{cor}

\begin{cor}
  For $m\le 1$, any morphism is strongly effective.
\end{cor}

Unfortunately, this is not the full story in general.  To see this,
consider the moduli stack of singular del Pezzo surfaces of degree 1.  In
our terms, this is the quotient stack ${\cal X}^{\alpha,\ge -2}_7/W(E_8)$,
and thus in particular inherits a decomposition from the one we have
constructed on ${\cal X}^{\alpha,\ge 2}_7$.  We can thus gain insight into
the geometry of this decomposition by considering the quotient.  In
characteristic 0, a singular del Pezzo surface has an equation of the form
\[
y^2 = x^3 + a_4(x,w)x + a_6(x,w)
\]
where $a_4$, $a_6$ are homogeneous of the given degree.  If we blow up the
base point of the anticanonical linear system, we obtain an elliptic
surface, and the type (mod $W(E_8)$) of the anticanonical curve $w=0$ can
be read off from the Kodaira symbol of the corresponding fiber of the
elliptic surface.  In particular, we can read off a parametrization of the
surfaces of any given type from Tate's algorithm (see specifically the
discussion in \cite[\S IV.9]{SilvermanJH:1994}), and then perform an
elimination to determine the equations satisfied by the coefficients of
$a_4$, $a_6$ in general.  (There is a minor issue, in that the scheme
corresponding to type $I_8$ is reducible, but it is easy enough to factor.)
We thus obtain a collection of $22$ ideals corresponding to the different
possible special fibers, and it is easy enough to both determine the
containment relations between the different ideals and to verify that the
intersection of any two of the closures is a union of closures of types.
We thus obtain a stratification of the moduli stack of del Pezzo surfaces
{\em in characteristic 0}, which we may denote pictorially as in Figure
\ref{fig:sakai_good}.\footnote{These Hasse diagrams were adapted from the
\TeX\ code for the analogous diagram in \cite{JoshiN/NakazonoN/ShiY:2016}.}
(Here we specify the root system of the $-2$-curves rather than the Kodaira
symbol.)  The open stratum $A_0^e$ is $m+2=9$-dimensional, and each
successive column decreases the dimension by $1$.  For surfaces with
$K^2=0$, one obtains the same diagram, with dimensions increased by $1$
(and affine rather than finite root systems), since ${\cal X}^{\alpha,\ge
  -2}_8$ is the universal $\Pic^0(C_\alpha)$ over ${\cal X}^{\alpha,\ge
  -2}_7$, and the map on types coming from blowing down is clearly
bijective.

\begin{figure}[t]
\centering
\begin{tikzpicture}[scale = 1]
\begin{scope}
\coordinate (P11s) at (0,0);
\coordinate (P11e) at ($(P11s)-(0.8,0)$);
\coordinate (P10s) at ($(P11e)-(0.6,0)$);
\coordinate (P12e) at ($(P11s)+(0.6,0)$);
\coordinate (P12s) at ($(P12e)+(0.8,0)$);
\coordinate (P13e) at ($(P12s)+(0.6,0)$);
\coordinate (P13s) at ($(P13e)+(0.8,0)$);
\coordinate (P14e) at ($(P13s)+(0.6,0)$);
\coordinate (P14s) at ($(P14e)+(0.8,0)$);
\coordinate (P15e) at ($(P14s)+(0.6,0)$);
\coordinate (P15s) at ($(P15e)+(0.8,0)$);
\coordinate (P16e) at ($(P15s)+(0.6,0)$);
\coordinate (P16s) at ($(P16e)+(0.8,0)$);
\coordinate (P17e) at ($(P16s)+(0.6,0)$);
\coordinate (P17s) at ($(P17e)+(0.8,0)$);
\coordinate (P18e) at ($(P17s)+(0.6,0)$);
\coordinate (P18s) at ($(P18e)+(0.8,0)$);
\coordinate (P19e) at ($(P18s)+(0.6,0)$);
\coordinate (P19s) at ($(P19e)+(0.8,0)$);
\coordinate (P21s) at (0,-1);
\coordinate (P21e) at ($(P21s)-(0.8,0)$);
\coordinate (P20s) at ($(P21e)-(0.6,0)$);
\coordinate (P22e) at ($(P21s)+(0.6,0)$);
\coordinate (P22s) at ($(P22e)+(0.8,0)$);
\coordinate (P23e) at ($(P22s)+(0.6,0)$);
\coordinate (P23s) at ($(P23e)+(0.8,0)$);
\coordinate (P24e) at ($(P23s)+(0.6,0)$);
\coordinate (P24s) at ($(P24e)+(0.8,0)$);
\coordinate (P25e) at ($(P24s)+(0.6,0)$);
\coordinate (P25s) at ($(P25e)+(0.8,0)$);
\coordinate (P26e) at ($(P25s)+(0.6,0)$);
\coordinate (P26s) at ($(P26e)+(0.8,0)$);
\coordinate (P27e) at ($(P26s)+(0.6,0)$);
\coordinate (P27s) at ($(P27e)+(0.8,0)$);
\coordinate (P28e) at ($(P27s)+(0.6,0)$);
\coordinate (P28s) at ($(P28e)+(0.8,0)$);
\coordinate (P29e) at ($(P28s)+(0.6,0)$);
\coordinate (P29s) at ($(P29e)+(0.8,0)$);
\coordinate (P31s) at (0,-2);
\coordinate (P32e) at ($(P31s)+(0.6,0)$);
\coordinate (P32s) at ($(P32e)+(0.8,0)$);
\coordinate (P33e) at ($(P32s)+(0.6,0)$);
\coordinate (P33s) at ($(P33e)+(0.8,0)$);
\coordinate (P34e) at ($(P33s)+(0.6,0)$);
\coordinate (P34s) at ($(P34e)+(0.8,0)$);
\coordinate (P35e) at ($(P34s)+(0.6,0)$);
\coordinate (P35s) at ($(P35e)+(0.8,0)$);
\coordinate (P36e) at ($(P35s)+(0.6,0)$);
\coordinate (P36s) at ($(P36e)+(0.8,0)$);
\coordinate (P37e) at ($(P36s)+(0.6,0)$);
\coordinate (P37s) at ($(P37e)+(0.8,0)$);
\coordinate (P38e) at ($(P37s)+(0.6,0)$);
\coordinate (P38s) at ($(P38e)+(0.8,0)$);
\coordinate (P39e) at ($(P38s)+(0.6,0)$);
\coordinate (P39s) at ($(P39e)+(0.8,0)$);
\coordinate (P41s) at (0,-3);
\coordinate (P42e) at ($(P41s)+(0.6,0)$);
\coordinate (P42s) at ($(P42e)+(0.8,0)$);
\coordinate (P43e) at ($(P42s)+(0.6,0)$);
\coordinate (P43s) at ($(P43e)+(0.8,0)$);
\coordinate (P44e) at ($(P43s)+(0.6,0)$);
\coordinate (P44s) at ($(P44e)+(0.8,0)$);
\coordinate (P45e) at ($(P44s)+(0.6,0)$);
\coordinate (P45s) at ($(P45e)+(0.8,0)$);
\coordinate (P46e) at ($(P45s)+(0.6,0)$);
\coordinate (P46s) at ($(P46e)+(0.8,0)$);
\coordinate (P47e) at ($(P46s)+(0.6,0)$);
\coordinate (P47s) at ($(P47e)+(0.8,0)$);
\coordinate (P48e) at ($(P47s)+(0.6,0)$);
\coordinate (P48s) at ($(P48e)+(0.8,0)$);
\coordinate (P49e) at ($(P48s)+(0.6,0)$);
\coordinate (P49s) at ($(P49e)+(0.8,0)$);
\node at ($(P10s)-(0.4,0)$){$A_0^e$};
\node at ($(P21s)-(0.4,0)$){$A_0^*$};
\node at ($(P22s)-(0.4,0)$){$A_1^*$};
\node at ($(P23s)-(0.4,0)$){$A_2^*$};
\node at ($(P24s)-(0.4,0)$){$A_3$};
\node at ($(P25s)-(0.4,0)$){$A_4$};
\node at ($(P26s)-(0.4,0)$){$A_5$};
\node at ($(P27s)-(0.4,0)$){$A_6$};
\node at ($(P28s)-(0.4,0)$){$A_7$};
\node at ($(P29s)-(0.4,0)$){$A_8$};
\node at ($(P18s)-(0.4,0)$){$A_7'$};
\node at ($(P32s)-(0.4,0)$){$A_0^+$};
\node at ($(P33s)-(0.4,0)$){$A_1^+$};
\node at ($(P34s)-(0.4,0)$){$A_2^+$};
\node at ($(P35s)-(0.4,0)$){$D_4$};
\node at ($(P36s)-(0.4,0)$){$D_5$};
\node at ($(P37s)-(0.4,0)$){$D_6$};
\node at ($(P38s)-(0.4,0)$){$D_7$};
\node at ($(P39s)-(0.4,0)$){$D_8$};
\node at ($(P47s)-(0.4,0)$){$E_6$};
\node at ($(P48s)-(0.4,0)$){$E_7$};
\node at ($(P49s)-(0.4,0)$){$E_8$};
%
%
%
%
\draw [->, thick] (P21s)--(P22e);
\draw [->, thick] (P22s)--(P23e);
\draw [->, thick] (P23s)--(P24e);
\draw [->, thick] (P24s)--(P25e);
\draw [->, thick] (P25s)--(P26e);
\draw [->, thick] (P26s)--(P27e);
\draw [->, thick] (P27s)--(P28e);
\draw [->, thick] (P28s)--(P29e);
\draw [->, thick] (P32s)--(P33e);
\draw [->, thick] (P33s)--(P34e);
\draw [->, thick] (P34s)--(P35e);
\draw [->, thick] (P35s)--(P36e);
\draw [->, thick] (P36s)--(P37e);
\draw [->, thick] (P37s)--(P38e);
\draw [->, thick] (P38s)--(P39e);
\draw [->, thick] (P47s)--(P48e);
\draw [->, thick] (P48s)--(P49e);
\draw [->, thick] ($(P27s)+(0,0.2)$)--($(P18e)-(0,0.2)$); 
\draw [->, thick] ($(P10s)-(0,0.2)$)--($(P21e)+(0,0.2)$); 
\draw [->, thick] ($(P21s)-(0,0.2)$)--($(P32e)+(0,0.2)$); 
\draw [->, thick] ($(P22s)-(0,0.2)$)--($(P33e)+(0,0.2)$); 
\draw [->, thick] ($(P23s)-(0,0.2)$)--($(P34e)+(0,0.2)$); 
\draw [->, thick] ($(P24s)-(0,0.2)$)--($(P35e)+(0,0.2)$); 
\draw [->, thick] ($(P25s)-(0,0.2)$)--($(P36e)+(0,0.2)$); 
\draw [->, thick] ($(P26s)-(0,0.2)$)--($(P37e)+(0,0.2)$); 
\draw [->, thick] ($(P27s)-(0,0.2)$)--($(P38e)+(0,0.2)$); 
\draw [->, thick] ($(P28s)-(0,0.2)$)--($(P39e)+(0,0.2)$); 
\draw [->, thick] ($(P36s)-(0,0.2)$)--($(P47e)+(0,0.2)$); 
\draw [->, thick] ($(P37s)-(0,0.2)$)--($(P48e)+(0,0.2)$); 
\draw [->, thick] ($(P38s)-(0,0.2)$)--($(P49e)+(0,0.2)$); 
\draw [->, thick] ($(P26s)-(0.1,0.35)$)--($(P47e)+(0.1,0.35)$);
\draw [->, thick] ($(P27s)-(0.1,0.35)$)--($(P48e)+(0.1,0.35)$);
\draw [->, thick] ($(P28s)-(0.1,0.35)$)--($(P49e)+(0.1,0.35)$);
\draw [->, thick] ($(P18s)-(0.1,0.35)$)--($(P39e)+(0.1,0.35)$);
\end{scope}
\end{tikzpicture}
\caption{The natural stratification of ${\cal X}_7^{\alpha,\ge -2}/W(E_8)$
  over $\Z[1/6]$.}
\label{fig:sakai_good}
\end{figure}
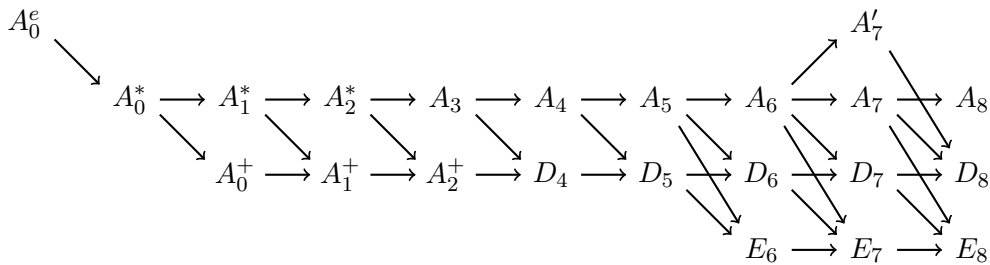
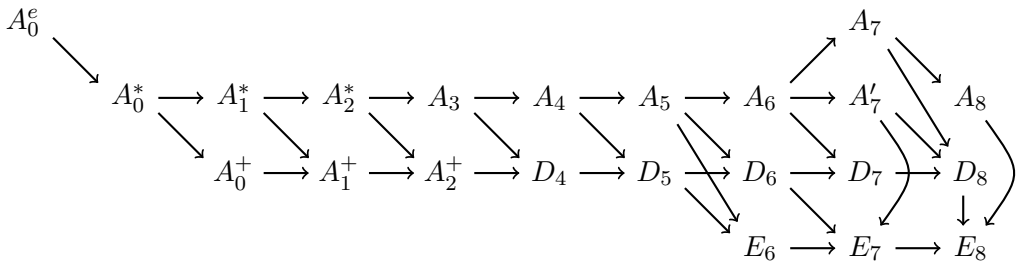
\begin{figure}[t]
\centering
\begin{tikzpicture}[scale = 1]
\begin{scope}
\coordinate (P11s) at (0,0);
\coordinate (P10e) at ($(P11s)-(0.8,0)$);
\coordinate (P10s) at ($(P10e)-(0.6,0)$);
\coordinate (P11e) at ($(P11s)+(0.6,0)$);
\coordinate (P12s) at ($(P11e)+(0.8,0)$);
\coordinate (P12e) at ($(P12s)+(0.6,0)$);
\coordinate (P13s) at ($(P12e)+(0.8,0)$);
\coordinate (P13e) at ($(P13s)+(0.6,0)$);
\coordinate (P14s) at ($(P13e)+(0.8,0)$);
\coordinate (P14e) at ($(P14s)+(0.6,0)$);
\coordinate (P15s) at ($(P14e)+(0.8,0)$);
\coordinate (P15e) at ($(P15s)+(0.6,0)$);
\coordinate (P16s) at ($(P15e)+(0.8,0)$);
\coordinate (P16e) at ($(P16s)+(0.6,0)$);
\coordinate (P17s) at ($(P16e)+(0.8,0)$);
\coordinate (P17e) at ($(P17s)+(0.6,0)$);
\coordinate (P18s) at ($(P17e)+(0.8,0)$);
\coordinate (P18e) at ($(P18s)+(0.6,0)$);
\coordinate (P19s) at ($(P18e)+(0.8,0)$);
\coordinate (P19e) at ($(P19s)+(0.6,0)$);

\coordinate (P21s) at (0,-1);
\coordinate (P20e) at ($(P21s)-(0.8,0)$);
\coordinate (P20s) at ($(P20e)-(0.6,0)$);
\coordinate (P21e) at ($(P21s)+(0.6,0)$);
\coordinate (P22s) at ($(P21e)+(0.8,0)$);
\coordinate (P22e) at ($(P22s)+(0.6,0)$);
\coordinate (P23s) at ($(P22e)+(0.8,0)$);
\coordinate (P23e) at ($(P23s)+(0.6,0)$);
\coordinate (P24s) at ($(P23e)+(0.8,0)$);
\coordinate (P24e) at ($(P24s)+(0.6,0)$);
\coordinate (P25s) at ($(P24e)+(0.8,0)$);
\coordinate (P25e) at ($(P25s)+(0.6,0)$);
\coordinate (P26s) at ($(P25e)+(0.8,0)$);
\coordinate (P26e) at ($(P26s)+(0.6,0)$);
\coordinate (P27s) at ($(P26e)+(0.8,0)$);
\coordinate (P27e) at ($(P27s)+(0.6,0)$);
\coordinate (P28s) at ($(P27e)+(0.8,0)$);
\coordinate (P28e) at ($(P28s)+(0.6,0)$);
\coordinate (P29s) at ($(P28e)+(0.8,0)$);
\coordinate (P29e) at ($(P29s)+(0.6,0)$);

\coordinate (P31s) at (0,-2);
\coordinate (P31e) at ($(P31s)+(0.6,0)$);
\coordinate (P32s) at ($(P31e)+(0.8,0)$);
\coordinate (P32e) at ($(P32s)+(0.6,0)$);
\coordinate (P33s) at ($(P32e)+(0.8,0)$);
\coordinate (P33e) at ($(P33s)+(0.6,0)$);
\coordinate (P34s) at ($(P33e)+(0.8,0)$);
\coordinate (P34e) at ($(P34s)+(0.6,0)$);
\coordinate (P35s) at ($(P34e)+(0.8,0)$);
\coordinate (P35e) at ($(P35s)+(0.6,0)$);
\coordinate (P36s) at ($(P35e)+(0.8,0)$);
\coordinate (P36e) at ($(P36s)+(0.6,0)$);
\coordinate (P37s) at ($(P36e)+(0.8,0)$);
\coordinate (P37e) at ($(P37s)+(0.6,0)$);
\coordinate (P38s) at ($(P37e)+(0.8,0)$);
\coordinate (P38e) at ($(P38s)+(0.6,0)$);
\coordinate (P39s) at ($(P38e)+(0.8,0)$);
\coordinate (P39e) at ($(P39s)+(0.6,0)$);

\coordinate (P41s) at (0,-3);
\coordinate (P41e) at ($(P41s)+(0.6,0)$);
\coordinate (P42s) at ($(P41e)+(0.8,0)$);
\coordinate (P42e) at ($(P42s)+(0.6,0)$);
\coordinate (P43s) at ($(P42e)+(0.8,0)$);
\coordinate (P43e) at ($(P43s)+(0.6,0)$);
\coordinate (P44s) at ($(P43e)+(0.8,0)$);
\coordinate (P44e) at ($(P44s)+(0.6,0)$);
\coordinate (P45s) at ($(P44e)+(0.8,0)$);
\coordinate (P45e) at ($(P45s)+(0.6,0)$);
\coordinate (P46s) at ($(P45e)+(0.8,0)$);
\coordinate (P46e) at ($(P46s)+(0.6,0)$);
\coordinate (P47s) at ($(P46e)+(0.8,0)$);
\coordinate (P47e) at ($(P47s)+(0.6,0)$);
\coordinate (P48s) at ($(P47e)+(0.8,0)$);
\coordinate (P48e) at ($(P48s)+(0.6,0)$);
\coordinate (P49s) at ($(P48e)+(0.8,0)$);
\coordinate (P49e) at ($(P49s)+(0.6,0)$);
\node at ($(P10s)-(0.4,0)$){$A_0^e$};

\node at ($(P21s)-(0.4,0)$){$A_0^*$};
\node at ($(P22s)-(0.4,0)$){$A_1^*$};
\node at ($(P23s)-(0.4,0)$){$A_2^*$};
\node at ($(P24s)-(0.4,0)$){$A_3$};
\node at ($(P25s)-(0.4,0)$){$A_4$};
\node at ($(P26s)-(0.4,0)$){$A_5$};
\node at ($(P27s)-(0.4,0)$){$A_6$};
\node at ($(P18s)-(0.4,0)$){$A_7$};
\node at ($(P28s)-(0.4,0)$){$A_7'$};
\node at ($(P29s)-(0.4,0)$){$A_8$};

\node at ($(P32s)-(0.4,0)$){$A_0^+$};
\node at ($(P33s)-(0.4,0)$){$A_1^+$};
\node at ($(P34s)-(0.4,0)$){$A_2^+$};
\node at ($(P35s)-(0.4,0)$){$D_4$};
\node at ($(P36s)-(0.4,0)$){$D_5$};
\node at ($(P37s)-(0.4,0)$){$D_6$};
\node at ($(P38s)-(0.4,0)$){$D_7$};
\node at ($(P39s)-(0.4,0)$){$D_8$};
\node at ($(P47s)-(0.4,0)$){$E_6$};
\node at ($(P48s)-(0.4,0)$){$E_7$};
\node at ($(P49s)-(0.4,0)$){$E_8$};

\draw [->, thick] (P21s)--(P21e);
\draw [->, thick] (P22s)--(P22e);
\draw [->, thick] (P23s)--(P23e);
\draw [->, thick] (P24s)--(P24e);
\draw [->, thick] (P25s)--(P25e);
\draw [->, thick] (P26s)--(P26e);
\draw [->, thick] (P27s)--(P27e);

\draw [->, thick] (P35s)--(P35e);
\draw [->, thick] (P36s)--(P36e);
\draw [->, thick] (P37s)--(P37e);
\draw [->, thick] (P38s)--(P38e);

\draw [->, thick] (P47s)--(P47e);
\draw [->, thick] (P48s)--(P48e);

\draw [->, thick] ($(P27s)+(0,0.2)$)--($(P17e)-(0,0.2)$);
\draw [->, thick] ($(P18s)-(0,0.2)$)--($(P28e)+(0,0.2)$);

\draw [->, thick] ($(P24s)-(0,0.2)$)--($(P34e)+(0,0.2)$);
\draw [->, thick] ($(P25s)-(0,0.2)$)--($(P35e)+(0,0.2)$);
\draw [->, thick] ($(P26s)-(0,0.2)$)--($(P36e)+(0,0.2)$);
\draw [->, thick] ($(P27s)-(0,0.2)$)--($(P37e)+(0,0.2)$);
\draw [->, thick] ($(P28s)-(0,0.2)$)--($(P38e)+(0,0.2)$);

\draw [->, thick] ($(P36s)-(0,0.2)$)--($(P46e)+(0,0.2)$);
\draw [->, thick] ($(P37s)-(0,0.2)$)--($(P47e)+(0,0.2)$);
\draw [->, thick] ($(P38e)+(0.3,-0.3)$)--($(P48e)+(0.3,0.3)$);

\draw [->, thick] ($(P26s)-(0.1,0.35)$)--($(P46e)+(0.1,0.35)$);
\draw[->, thick] ($(P28s)-(0.2,0.3)$) .. controls ($(P38s)+(0.3,0)$) .. ($(P48s)+(-0.2,0.3)$);
\draw[->, thick] ($(P29s)-(0.2,0.3)$) .. controls ($(P39s)+(0.3,0)$) .. ($(P49s)+(-0.2,0.3)$);

\draw [->, thick] ($(P10s)-(0,0.2)$)--($(P20e)+(0,0.2)$); 
\draw [->, thick] ($(P21s)-(0,0.2)$)--($(P31e)+(0,0.2)$); 
\draw [->, thick] ($(P22s)-(0,0.2)$)--($(P32e)+(0,0.2)$); 
\draw [->, thick] ($(P23s)-(0,0.2)$)--($(P33e)+(0,0.2)$); 
\draw [->, thick] (P32s)--(P32e);
\draw [->, thick] (P33s)--(P33e);
\draw [->, thick] (P34s)--(P34e);

\draw [->, thick] ($(P18s)-(0.1,0.35)$)--($(P38e)+(0.1,0.35)$);

\end{scope}
\end{tikzpicture}
\caption{The na\"{i}ve poset of types in ${\cal X}_7^{\alpha,\ge -2}/W(E_8)$.}
\label{fig:sakai_bad}
\end{figure}

Now, the components of the anticanonical curve of a del Pezzo surface of
degree 1 all have negative self-intersection, so are rigid; as a result, if
there is a morphism at all between two such types, it is unique.  We find
in particular that the corresponding subcategory of the category of types
is simply the natural poset of root subsystems of $E_8$ (as extended by the
Picard type), ordered by inclusion.  However, the resulting diagram (Figure
\ref{fig:sakai_bad}) is slightly different (essentially from
\cite{SakaiH:2001}, except that we have added the missing arrow from $A_7$
to $D_8$; note that $A_7$ and $A_7'$ have swapped positions in the diagram
to avoid arrows from the top row to the bottom row in each case).  In
particular, we find that although the combinatorics suggests that there
should be degenerations
\[
A'_7\to E_7,\quad
A'_7\to E_8,\quad
 D_8\to E_8,\quad
 A_8\to E_8,
\]
the corresponding morphisms are not effective.  Note that the three cases
corresponding to vertical arrows cannot possibly be {\em strongly}
effective, since the dimensions of the corresponding substacks are the same.

In fact, the situation is even worse than this suggests: any surface of
type $E_8$ in characteristic 3 can be obtained as the reduction mod 3 of a
surface of type $A_8$ over a suitable $3$-adic field.  (Similar statements
hold for the other three cases in characteristic 2.)  Even if one restricts
ones attention to equicharacteristic deformations, there are still issues:
in characteristic 3, one can obtain {\em some} surfaces of type $E_8$ as
limits from type $A_8$.  As a result, our decomposition of the moduli stack
of del Pezzo surfaces is neither a stratification in characteristic 3 nor
over $\Z$, as the closure of type $A_8$ meets type $E_8$ in a proper
substack.

  Indeed, over a field, a surface with an $A_8$ singularity has the form
  \[
  y^2 +(a_{11}t+a_{10}u) xy + a_{30}u^3 y = x^3,
  \]
  up to changes of basis in $x$ and $y$ (but with no such changes of basis
  in $t$ and $u$ required).  The discriminant of this surface has the form
  \[
  c_1^2 t^3 u^9 + 3 c_1c_2 t^2 u^{10} + 3 c_2^2 t u^{11} + c_3 u^{12}
  \]
  for suitable functions $c_1$, $c_2$, $c_3$ of the parameters.  The
  discriminant of an $E_8$ surface in characteristic $\notin \{2,3\}$ has
  degree precisely 2 in $t$, but such a discriminant is not in the Zariski
  closure of the above set of discriminants, and thus there can be no
  degeneration from $A_8$ to $E_8$ in such cases.  On the other hand,
  consider the del Pezzo surface
\[
y^2 = x^3 - \frac{243t^2-54tu-u^2}{4}x^2 - \frac{3u^2(27t-5u)}{2}x
-\frac{tu^4(27t-4u)}{4}
\]
over the rationals.  This has an $A_8$ singularity at $y=x=u=0$, but
modulo $3$ becomes the surface
\[
y^2 = x^3 + u^2 x^2 + t u^5,
\]
which now has an $E_8$ singularity at $y=x=u=0$.  A characteristic 3
surface with an $A_8$ singularity at $u=0$ has discriminant of
the form $u^9(at^3+bu^3)$, and thus cannot degenerate to the above surface
of discriminant $-tu^{11}$.  However, the family of surfaces
\[
y^2 = x^3 + v^3(tv-u)^2 x^2 + u^2(u-vt)(u-v^3t) x + u^4 t (u+v^3t)
\]
over $\Spec(\F_3[v])$ generically has an $A_8$ singularity at $u=0$, but at
$v=0$ becomes the surface
\[
y^2 = x^3 + u^4 x + u^5 t
\]
with an $E_8$ singularity.  Every smooth fiber of the latter surface has
$j$-invariant $0$, so is a supersingular curve.  (There is one more
isomorphism class of surfaces with an $E_8$ singularity, namely the
quasielliptic surface $y^2=x^3+u^5 t$, but this is a degeneration of the
$j=0$ case.)

Let us consider what a surface of type $A_8$ would look like if we
considered it relative to a smooth anticanonical curve $C'$.  To obtain
type $A_8$, two things must happen: all roots of $A_8$ must vanish in
$\Pic(C')$, but also every root {\em not} in the subsystem must {\em not}
vanish.  Indeed, if some root not in $A_8$ were to vanish, then every root
of $E_8$ would vanish, and the surface would have type $E_8$ instead.
Since the lattice $\Lambda_{A_8}$ has index 3 in $\Lambda_{E_8}$, we see
that the image of $\Lambda_{E_8}$ in $\Pic^0(C')$ is a $3$-torsion
subgroup.  But in characteristic not $3$, it is impossible to degenerate a
nontrivial $3$-torsion point to the identity.  (In the $3$-adic case, we
can take the $3$-torsion point to be in the kernel of reduction, while in
the equicharacteristic case, we may take the special fiber of $C'$ to be
supersingular.)

Of course, the above argument is quite ad hoc, but it turns out that the
obstruction generalizes.  First, we note that the claim about degeneration
of torsion points holds for curves of arbitrary genus.

\begin{lem}\label{lem:order_r_closed}
  Let $R$ be a dvr with field of fractions $K$ and residue field $k$.
  Let $C_R$ be a smooth curve over $R$, and let ${\cal L}$ be a line bundle on
  $C_R$ such that ${\cal L}_K$ has exact order $r$ in $\Pic(C_K)$.  If $r\in R^\times$,
  then ${\cal L}_k$ also has exact order $r$.
\end{lem}

\begin{proof}
  We first note that ${\cal L}_k^r\cong \sO_{C_k}$, since this is a closed
  condition.  But the assumption on $r$ implies that $\Pic^0(C_R)[r]$ is
  \'etale, and thus the closed subscheme $\Pic^0(C_R)[d]$ is also open for
  any divisor $d|r$.  It follows that ${\cal L}_k^d\not\cong \sO_{C_k}$ for every
  proper divisor $d|r$.
\end{proof}

\begin{rem}
  This also holds if $R$ is equicharacteristic and $C_k$ is ordinary, since
  then the complement of $\Pic^0(C_R)[p^l]$ in $\Pic^0(C_R)[p^{l+1}]$ is
  again both closed and open.  The claim fails in the remaining cases with
  $\ch(k)|r$, however.
\end{rem}
  
Now, a key fact about the curve $C'$ we used above was that it was
orthogonal to all of the roots in the relevant root systems.  In
particular, we can view $C'$ as a curve on the surface obtained by
contracting the $-2$-curves.  Moreover, $C'$ is an {\em ample} divisor on
that singular surface, and the claim boiled down to showing that the image
of the Picard group of the minimal desingularization in $C'$ contained a
copy of $\Lambda_{A_8}^\perp/\Lambda_{A_8}$.

Note that in Lemma \ref{lem:order_r_closed}, we only consider the Picard
groups of the two fibers; as a result, to apply the result, we need only
understand surfaces over fields.  Moreover, base changing to a finite
extension of $R$ has no effect on the order of $L_K$ or $L_k$, and thus we
can take a limit to a valuation ring in which both the residue field and
the field of fractions are algebraically closed.  The objective, therefore,
is to find more general conditions in which the restriction of $\Pic(X)$ to
some curve in $X$ contains nontrivial torsion points.

With that in mind, let $Y$ be a normal surface over an algebraically closed
field $k$, with minimal desingularization $\tilde{Y}$.  To fix ideas,
suppose for the moment that we have a rational map $\phi:Y\to \P^1$ such
that the locus of indeterminacy is a single smooth point of $Y$ and
$\phi^*\sO_{\P^1}(1)$ is ample.  Then blowing up the corresponding point of
$\tilde{Y}$ gives a surface $\tilde{X}$ with a morphism to $\P^1$.  Now,
let $D$ be a divisor class on $\tilde{X}$ such that the restriction of $D$
to the generic fiber of $\tilde{X}$ is principal.  If we choose a function
$f$ with that divisor (which is uniquely determined modulo $k(\P^1)^*$),
then $D-\div(f)$ is certainly linearly equivalent to $D$, but now has
trivial restriction to the generic fiber.  This implies that $D-\div(f)$ is
a supported on a finite set of fibers, and is thus a linear combination of
components of fibers.

Now, suppose $C$ is a component of a fiber.  If $C$ does not meet the
exceptional curve of $\tilde{X}\to \tilde{Y}$, then its image in $Y$ is
disjoint from the generic fiber of $\phi$.  Since $\phi^*\sO_{\P^1}(1)$ is
ample, this implies that the image of $C$ must be a single point, and is
thus one of the singular points of $Y$.  We thus see that the components of
fibers split into two classes: those which are contracted in $Y$, and those
which meet the exceptional curve.  Since the exceptional curve meets each
fiber precisely once, we see that every fiber contains exactly one
component meeting the exceptional curve.

Suppose $C_1$,\dots,$C_n$ are the contracted components, and let $C$ be any
other component.  If $F$ is the fiber containing $C$, then $F-C$ is a sum
of components not meeting the exceptional locus, and thus we have
\[
F-C\in \Z\langle C_1,\dots,C_n\rangle.
\]
Now, suppose we are given a divisor class $D\in \Pic(\tilde{Y})$ such that
$rD\in \Z\langle C_1,\dots,C_n\rangle$ for some integer $r\ge 1$, and
suppose that $D$ has trivial restriction in the Picard group of the generic
fiber of the pencil.  Then the pullback of $D$ to $\tilde{X}$ is a linear
combination of components of fibers, and therefore has an expression of the
form
\[
D \sim \sum_i a_i C_i + \phi^*Z
\]
where $a_i\in \Z$ and $Z$ is a divisor on $\P^1$.  But this implies that
\[
(D-\sum_i a_i C_i)\cdot C_j = 0
\]
for all $j$.  By Mumford's criterion for contractibility
\cite{MumfordD:1961}, the intersection form of $C_1$,\dots,$C_n$ is
negative definite, and we thus find that
\[
D=\sum_i a_i C_i\in \Z\langle C_1,\dots,C_n\rangle.
\]
In other words, the map from $\Q\langle C_1,\dots,C_n\rangle\cap
\Pic(\tilde{Y})$ to the Picard group of the generic fiber has kernel
precisely $\Z\langle C_1,\dots,C_n\rangle$.  Thus if $\Z\langle
C_1,\dots,C_n\rangle$ is not {\em saturated} in $\Pic(\tilde{Y})$, the
restriction map from $\Pic(\tilde{Y})$ induces a torsion point as required.

The main difficulty in applying the above argument in general is the
requirement that the pencil have a single base point, as this requires in
particular that $\langle C_1,\dots,C_n\rangle^\perp\subset \Pic(\tilde{Y})$
contains smooth (and non-rigid) curves of self-intersection $1$.  We would
thus like to extend the argument to deal with larger base loci.  The
difficulty, of course, is that fibers of the pencil can then have multiple
components meeting the base locus.

Let $Y$,$\tilde{Y}$ be as above, but now let ${\cal L}_0$ be an arbitrary
very ample line bundle on $Y$.  By Bertini's theorem, the corresponding
linear system contains smooth curves, so let us fix such a curve $C_0$.  A
further application of Bertini to $C_0$ lets us choose a curve $C_1$ in the
linear system meeting $C_0$ in $C_0^2$ distinct points, and we can then
choose $C_2$ meeting $C_0$ in the complement of $C_0\cap C_1$, so that the
linear system spanned by $C_0$, $C_1$, $C_2$ is base-point-free.  Let
$P_{12}$ denote the pencil through $C_1$ and $C_2$, and let $P$ denote the
pencil through $C_0$ and the generic point of $P_{12}$.  Now, although $P$
has base points, they are in general defined over an extension field of the
field $k(\P^1)$ over which $P$ is defined.  In fact, we have the following.

\begin{lem}
  The base points of $P$ are defined over the separable closure of
  $k(\P^1)$, and form a single orbit under the action of the absolute
  Galois group $\Gal(k(\P^1))$.
\end{lem}

\begin{proof}
  That the splitting field of the base scheme of $P$ is separable follows
  from the fact that $C_0\cap C_1$ is reduced, and thus the same holds if
  we replace $C_1$ by the generic fiber of $P_{12}$.  The linear
  system $P_{12}|_{C_0}$ is base-point-free and thus determines a morphism from
  $C_0$ to $\P^1$.  The base points of $P$ are then just the preimage of
  the generic point under this morphism, and thus the base scheme is
  $\Spec(k(C_0))$.  Since $k(C_0)$ is a field, transitivity is immediate.
\end{proof}

Now let $\tilde{X}$ be obtained from $\tilde{Y}$ by blowing up the base
locus of $P$.  If $D$ is a divisor on $\tilde{Y}$ {\em defined over $k$}
which is principal on the generic fiber of $P$, then, just as before, its
preimage in $\tilde{X}$ is linearly equivalent to an integer linear
combination of components of fibers of $P$.  Let $D'$ be such a
representative, and let $\sigma$ be in the absolute Galois group of
$k(\P^1)$.  Then $\sigma D'-D'\sim \sigma D-D$ is principal, but disjoint
from the generic fiber, and is thus the divisor of a function pulled back
from the base of the pencil.  But then since $k(\P^1)$ has trivial Brauer
group, we can add such a divisor to $D'$ to make it Galois-stable.
In particular, $D'$ is a sum of Galois-orbits of components of fibers.

The key insight now is that the components which are contracted on $Y$ are
contained in $k$-rational fibers.  If $C$ is a component of such a fiber
which is {\em not} contracted on $Y$, then its image in $\tilde{Y}$ meets
the base locus of $P$.  But then by transitivity of the Galois group on the
base points, we conclude that the image on $\tilde{Y}$ of the sum of the
Galois orbit of $C$ actually {\em contains} the base locus.  Then the
resulting Galois stable curve is simply the strict transform of the
corresponding fiber of $P$ on $Y$, and may thus be expressed as a linear
combination of the preimage of the fiber and the components of the preimage
of the singular points.

Thus if $C_i$ are the contracted components, we obtain an expression
\[
D' = \sum_i a_i C_i + \sum_j F_j + \sum_l D_l
\]
where each $F_j$ is a fiber and each $D_l$ is an effective divisor
supported on a non-rational fiber.  But then
\[
D'\cdot C_j = \sum_i a_i C_i\cdot C_j,
\]
and we can argue as before.  We thus obtain the following result.

\begin{prop}\label{prop:Pic_injectivity}
  Let $Y$ be a normal surface over an algebraically closed field $k$, with
  minimal desingularization $\tilde{Y}$ and exceptional curves
  $C_1$,\dots,$C_n$.  Let $\Lambda=\Z\langle C_1,\dots,C_n\rangle\subset
  \Pic(\tilde{Y})$, and let $\Lambda^+\subset \Pic(\tilde{Y})$ be the
  subgroup consisting of line bundles such that ${\cal L}^r\in \Lambda$ for
  some $r\in \Z$ which is invertible in $k$.  If $C_0$ is a smooth very
  ample curve on $Y$, then the restriction map $\Pic(\tilde{Y})\to
  \Pic(C_0)$ induces an exact sequence
  \[
  0\to \Lambda\to \Lambda^+\to \Pic(C_0)
  \]
\end{prop}

\begin{proof}
  The above argument shows that this holds for the generic curve in some
  pencil containing $C_0$, and then the result follows by Lemma
  \ref{lem:order_r_closed}.
\end{proof}

In the following result, note that by \cite{ArtinM:1974} any flat family of
surfaces with sufficiently nice rational singularities admits a uniform
minimal desingularization over an \'etale cover of the base.  In our
specific application, we start with a family of smooth surfaces, so there
is no issue.

\begin{prop}
  Let $R$ be a dvr and $Y/\Spec(R)$ be a projective scheme such that the
  fibers are normal rational surfaces, and suppose that
  $\tilde{Y}/\Spec(R)$ is a fiberwise minimal desingularization of $Y$.  Let
  $C_1$,\dots,$C_n$ be the components of the exceptional locus of
  $\tilde{Y}_k$, and $C'_1$,\dots,$C'_m$ the components of the exceptional
  locus of $\tilde{Y}_K$.  Let $a_{ij}$ be the multiplicity of $C_j$ in the
  special fiber of the closure of $C'_i$, and let $\vec{a}_i\subset \Z^n$
  be the corresponding collection of vectors.  Then the quotient
  \[
  (\Q\langle \vec{a}_1,\dots,\vec{a}_m\rangle \cap \Z^n)
  /
  (\Z\langle \vec{a}_1,\dots,\vec{a}_m\rangle)
  \]
  is trivial if $\ch(k)=0$ and a $p$-group if $\ch(k)=p$.
\end{prop}

\begin{proof}
  Both $Y_k$ and $Y_K$ have the property that their generic hyperplane
  section relative to some sufficiently very ample divisor class is smooth,
  and thus there is a hyperplane defined over $R$ such that both fibers of
  the corresponding section $C_0$ are smooth.  For any vector in $\Q\langle
  \vec{a}_1,\dots,\vec{a}_m\rangle\cap \Z^n$, let $D$ be the corresponding
  linear combination of $C_1,\dots,C_n$.  Since $Y$ has rational fibers,
  the invertible sheaf $\sO_{Y_k}(D)$ extends uniquely from $Y_k$ to $Y$,
  and some power of $\sO_{Y_K}(D)$ will have divisor class in $\Z\langle
  C'_1,\dots,C'_m\rangle$.  The restriction $\sO_{C_0}(D)$ is trivial on
  the special fiber of $C_0$, and thus by Lemma \ref{lem:order_r_closed} is
  trivial (or has $p$-power order) on the generic fiber.  The result
  follows from Proposition \ref{prop:Pic_injectivity}.
\end{proof}

To apply this to anticanonical surfaces, we need to understand when
configurations of anticanonical components can be contracted.

\begin{lem}
  Let $X$ be an anticanonical rational surface, and let $C_1$,\dots,$C_l$
  be a sequence of anticanonical components such that the intersection
  matrix of $C_1,\dots,C_l$ is negative definite and not every component
  appears.  Then the curve $\cup_i C_i$ is contractible.
\end{lem}

\begin{proof}
  First note that if any of $C_1$,\dots,$C_l$ are $-1$-curves, then we can
  simply blow down that curve and apply the Lemma to $X_{m-1}$.  We may
  thus assume that $C_i^2\le -2$ for $1\le i\le l$.  Now, let
  $m_1,\dots,m_s$ be the multiplicities of the anticanonical components,
  and consider the divisor class $Z^+=\sum_{1\le i\le l} m_i C_i$.  Then
  for $1\le i\le l$,
  \[
  Z^+\cdot C_i = (C_\alpha-\sum_{l<j\le s} m_j C_j)\cdot C_i
  = (2+C_i^2) - \sum_{l<j\le s} m_j (C_j\cdot C_i)
  \le 0,
  \]
  while
  \[
  (Z^+)^2 = \sum_{1\le i\le l}
    m_i(2+C_i^2) - \sum_{\substack{1\le i\le l\\ l<j\le
        s}} m_im_j (C_j\cdot C_i)
    <0,
  \]
  since $C_\alpha$ is connected.  It follows that the {\em fundamental
    cycle} \cite{ArtinM:1966} of $\cup_i C_i$ satisfies $Z=\sum_i r_i C_i$
  with $1\le r_i\le m_i$, and thus in particular both $Z$ and $C_\alpha-Z$
  are effective.  It follows from Lemma \ref{lem:Ca_is_num_conn} that every
  connected component of $Z$ has arithmetic genus 0, and then
  \cite{ArtinM:1962,ArtinM:1966} tells us that every connected component of
  $\cup_i C_i$ contracts to a rational singularity.  The result essentially
  follows immediately.  (Artin assumes the configuration is connected, but
  the proof carries over.)
\end{proof}

\begin{cor}
  Let $\psi:T_1\to T_2$ be a morphism of types, suppose $C_1,\dots,C_l$ is
  a contractible configuration of anticanonical components of $T_1$, let
  $\Lambda_1\subset \Pic(X_m)$ be the corresponding sublattice, and let
  $\Lambda_2\subset \Lambda_1$ be the sublattice coming from those
  components of $T_2$ which are supported on $C_1,\dots,C_l$ (relative to
  $\psi$).  If $\Lambda_1/\Lambda_2$ has torsion of degree prime to the
  characteristic, then $\psi$ is ineffective.
\end{cor}

This in particular explains the obstructions for del Pezzo surfaces of
degree 1 (and corresponding obstructions for rational surfaces with $K^2=0$
and no $-d$-curves with $d>2$, the subject of \cite{SakaiH:2001}), as well
as giving obstructions in other cases.  For instance, let $T_2$ be the type
of (quartic del Pezzo) surfaces with anticanonical curve decomposition:
\[
C_\alpha
=
(s-e_1-e_4)+
(f-e_1-e_2)+
(s-e_2-e_3)+
(f-e_3-e_4)+
(e_1)+
(e_2)+
(e_3)+
(e_4),
\]
and let $T_1$ correspond to surfaces with
\[
C_\alpha
=
2(s-f)+4(f-e_1-e_2)+3(e_1-e_2)+6(e_2-e_3)+5(e_3-e_4)+4(e_4).
\]
Since the components of $T_1$ are linearly independent, it is easy to
verify that there is a morphism $\psi:T_1\to T_2$.  On the other hand, the
only components of $T_2$ supported on the complement of $e_4$ are the first
four.  The corresponding sublattice of $\Z^5$ is not saturated (it contains
$2e_2-2e_4$ but not $e_2-e_4$), and thus such a morphism can only
correspond to a degeneration when the special fiber has characteristic 2.
Here we should note that if we blow up a generic point on each component of
the generic fiber (and then blowup a corresponding point downstairs), then
there are no subsets of the anticanonical components to which the
obstruction applies.  Thus to apply this obstruction fully, one must in
principle consider not only the types themselves but also all possible
blowdowns.  Since there are infinitely many possible blowdowns, it is
unclear if this is actually computable!

\medskip

For surfaces with $K^2=0$, although the computer calculation only showed
that Figure \ref{fig:sakai_good} was valid in characteristic 0, it is not
too hard to extend it to $\Z[1/6]$.  If $E$ is a fixed elliptic curve, then
del Pezzo surfaces (with blowdown structure) with anticanonical curve $E$
are classified by maps $\phi\in \Hom(\Lambda_{E_8},E)$.  If $R\subset E_8$
is an indecomposable subsystem, then there is an anticanonical curve of
type $R$ if every root in $\ker\phi$ is contained in $R$.  In particular,
if $\Lambda_R$ is saturated in $\Lambda_{E_8}$, then the generic point in
$\Hom(\Lambda_{E_8}/\Lambda_R,E)$ will give rise to such a surface.
Moreover, we may verify by direct computation that the generic surface with
such an anticanonical fiber has non-constant $j$-invariant, at which point
dimension considerations tell us that the generic surface arises from a map
$\Hom(\Lambda_{E_8}/\Lambda_R,E)$.  It follows easily that if $R\subset R'$
with $\Lambda_R$ saturated, then the corresponding degeneration is strongly
effective in any characteristic.  In particular, every arrow in Figure
\ref{fig:sakai_good} remains strongly effective over any field.

One is thus left to consider degenerations from the unsaturated cases
$A'_7$, $A_8$, $D_8$.  Since $\Lambda_{E_8}/\Lambda_{D_8}\cong \Z/2\Z$, and
$\Lambda_{E_8}/\Lambda_{A_8}\cong \Z/3\Z$, such surfaces are classified by
points of order precisely 2 or 3 on $E$ as appropriate, and thus over
$\Z[1/6]$ cannot be degenerated.  Similarly, $A'_7$ surfaces with a fiber
isomorphic to $E$ are labelled by a $2$-torsion point and a point in $E$,
and thus degenerate to the surface of type $D_8$ labelled by the same
$2$-torsion point, but not to surfaces of type $E_7$ or $E_8$ over
$\Z[1/6]$.

This also makes it relatively straightforward to construct the exotic
degenerations in small characteristic, as one need simply construct an
appropriate family of curves $E$ with a $2$- or $3$-torsion point
degenerating to the identity.  (This is how the above examples of $A_8\to
E_8$ degenerations were produced.)

\section{Surfaces and singularities of equations}
\label{sec:singularities}

Now that we have introduced the decomposition of the moduli stack into
types, it is natural to ask what a given type implies about the
corresponding difference/differential equations.  The type of a blown up
surface translates into information about how the direct image of a sheaf
with given Chern class meets the anticanonical curve on $F_2$, and we
argued above that this controls the singularities of the equation.  Though
it is clear that the sheaf must meet the anticanonical curve on $F_2$ at
any singular point of the equation, and that the equation is singular at
any point in the support of the intersection, we have not so far been clear
about how the local behavior of the equation at the singularity is recorded
in the sheaf structure of the intersection.  For any given equation, of
course, one can simply turn it into a sheaf and keep track of the different
blowups; however, when classifying degenerations of a given elliptic
scenario (or generalizations of a given non-elliptic scenario), one
generally obtains simply a list of possible types of surface, and would
then like to translate that into information about the singularities of the
equation.

We restrict our attention to the case that the sheaf on $F_2$ can be
separated from the anticanonical curve by a sequence of blowups and minimal
lifts.  This is a relatively mild restriction, since Lemma
\ref{lem:resolving_pseudotwist} tells us that any (relaxed) equation is
isospectral to one that can be so resolved, so that the cases we are
excluding involve apparent singularities.  We similarly assume that the
support of the sheaf does not contain a fiber, as again this leads to
apparent singularities.

More precisely, since singularities are local in nature, we are interested
in singularities at some fixed point $x\in \P^1$, and only need to be able
to separate the sheaf from the anticanonical curve in a neighborhood of
that fiber.  We are thus given a blowup of $F_2$ in which every point blown
up is on the same fiber, and a sheaf $M$ on that blowup such that
$x\notin \rho(\pi(\supp(M)\cap C_\alpha))$.
If $M$ has first Chern class $ds+d'f-r_1e_1-\cdots-r_me_m$, then the
intersection of $\pi_*M$ with the anticanonical curve near $\rho^{-1}x$ is
independent of $d'$, and thus the question is to translate the remaining
data into local information about the equation.  Equivalently, we are given
the intersection numbers of $c_1(M)$ with the various components of the
preimage of $x$.

Again, since singularities are local, we should actually replace the
surface by a formal neighborhood of the fiber.  Helpfully, the intersection
numbers are still well-defined even when $M$ is only a sheaf on the formal
completion of the fiber, since the various components are still proper
curves.  Now, on the formal neighborhood, $M$ has a natural direct sum
decomposition, with one summand for each point of intersection.  This
remains true as we blow up points, and thus eventually we end up with the
case of a spectral curve meeting an exceptional curve $e$ in a single point
not on the anticanonical curve.  If the spectral curve has a single branch,
then we can parametrize it to identify the given component of $M$ as a
torsion-free sheaf on the image of a map from $\Spec(k[[t]])$ to the formal
completion along the fiber.  More generally, we may use the fact that
$\supp(M)$ is {\em generically} smooth (at least in characteristic 0!) to
reduce to the case that $M$ is locally the direct image of the structure
sheaf of $\Spec(k[[t]])$ (which is necessarily true when $M$ has smooth
support).

This leaves questions of limits, but these are issues in any event, since
we need to allow $\supp(M)$ to acquire additional ($-2$-curve) components.
So in general what we need to understand is not the precise condition
corresponding to $M|_{C_\alpha}$ but rather its Zariski closure!  We may
thus feel free to impose various nonempty open conditions.  If $|c_1(M)|$
is not generically integral, then we can write $M$ globally as an
extension; if it is and $\dim|c_1(M)|>0$, then in sufficiently large
characteristic the generic curve in $|c_1(M)|$ (after separating from
$C_\alpha$) meets each exceptional curve $e$ in an \'etale subscheme.
Indeed, $|c_1(M)|$ has no base points (since it is generically disjoint
from $C_\alpha$) and thus induces a linear family of polynomials on $e$ with no
common divisor.  Such a family is generically irreducible (any
factorization over $k(u)$ induces a factorization over $k[u]$ and degree
considerations force one factor to be constant), and thus the only way it
can fail to cut out an \'etale subscheme is if it is inseparable, i.e., if
$c_1(M)\cdot e$ is a multiple of the characteristic.  Thus (away from such
characteristics) over the algebraic closure of the generic point, $M$ has
$c_1(M)\cdot e$ factors in a formal neighborhood of $e$, with each factor
meeting $e$ with multiplicity 1.  (In particular each factor is smooth.)
We can also feel free to avoid any countable set of possible intersection
points (e.g., the points where $e$ meets other exceptional components or
the strict transform of a fiber).  Note that in the rigid case, where this
fails, we can still determine the local conditions $M$ has to satisfy by
taking the Zariski closure of the conditions on low-order coefficients of
$B$.

It is instructive to look at the case of a smooth point of $C_\alpha$
(which is always the case for the last separating blowup!).  We may
parametrize the formal neighborhood by $\Spec(k[[x,y]])$ such that
$C_\alpha$ is given by $x=0$.  The isomorphism class of $M|_{C_\alpha}$ is
determined by a partition $\lambda$ giving the lengths of the different
summands of $M|_{C_\alpha}$, and $M$ itself can then be written as the
cokernel of an $m\times m$ matrix $\Psi$ reducing modulo $x$ to the
diagonal matrix with entries $y^{\lambda_i}$, where $m$ is the number of
(nonzero) parts of $\lambda$.  We thus see that $\det(\Psi)$ is contained
in the product ideal $\prod_i (y^{\lambda_i},x).$ In fact, if $k=\bar{k}$,
one generically has a factorization $\det(\Psi) = \prod_i p_i(x,y)$ where
$p_i\in (y^{\lambda_i},x)$.  Indeed, any element of the product ideal has
Newton polygon bounded by the polygon with slopes $1/\lambda_i$, and thus
having such a factorization is an open condition on the product ideal,
which is clearly nonempty for $\Psi$ (consider diagonal matrices!).  We
thus generically have $m$ factors, each of which can be parametrized by
expressing $x$ as a power series in $y$.  There are, of course nongeneric
cases, such as the matrix
\[
\Psi=
\begin{pmatrix}
  y^2 & x\\
  x & y
\end{pmatrix},
\]
where the Newton polygon is no longer generic.  Such cases correspond to
precisely the sorts of behaviors excluded above; in this particular
instance, the pathology is that the strict transform of the spectral curve
meets the two exceptional components in their point of intersection.  (One
also has more degenerate cases in which the minimal lift of $M$ contains an
exceptional component, another type of apparent singularity.)

When applying this to the problem of classifying relaxed equations, there
is one complication, however: we in general want to consider sheaves on
more general rational surfaces with blowdown structure, where there is no
longer a canonical association to equations.  Of course, as we discussed
above, we still have well-defined equations up to scalar gauge, and we can
fix that freedom by choosing a section transverse to the anticanonical
curve, blowing up additional points as necessary to separate that section
from $C_{\alpha}$, then acting by $W(E_{m+1})$ to find a blowdown structure
where that section is $s-f$.  It is generally straightforward to determine
how the choice of section affects the gauge; in any case where it is not
straightforward, one can simply choose another section of the same degree
and work on the surface where both sections are resolved, at which point
the resulting scalar gauge transformation can be read off from the
singularities of the resulting first-order equation.  In particular, once
we have understood singularities in the $F_2$ case, we will be able to
understand the general case, up to an overall scalar gauge transformation
freedom (which we can also understand).

Note also that once we have restricted $M$ to be the image of a structure
sheaf of a locally smooth curve, we may feel free to blow down any
component of the fiber which is a $-1$-curve not meeting $M$.  After doing
so, we obtain a surface such that every point we blow up after the first is
a point on the previous exceptional curve, as otherwise the final
exceptional locus would contain $-1$-curves disjoint from $M$.  Moreover,
any time we blow up a point which is on just one anticanonical component,
the type of surface that results is not affected by the choice of point.
In particular, we may as well stop blowing up just before we blow up a
smooth point of the anticanonical curve, as in that case the remainder of
the resolution will simply blow up that point multiple times.  (For
bookkeeping purposes, it is convenient to blow that point up {\em once}, so
that its location is recorded in the surface.)  We thus restrict our
attention to types of surfaces in which every blowup but the first blows up
a point of an exceptional component and every $-1$-curve but the last is an
anticanonical component.  Call such a type ``minimal''.  Minimal types
correspond to the case of unibranch spectral curves meeting the last
exceptional anticanonical component in a single point (with multiplicity
1); any minimal type is thus the first in a sequence of types obtained by
repeatedly blowing up that point of intersection.

In this way, we reduce to the following problems:
\begin{itemize}
  \item[(1)] Given a minimal type, what do the resulting maps from
    $\Spec(k[[t]])$ to $F_2$ look like?  (And can this be inverted?)  What
    about the higher-order versions of the minimal type?
  \item[(2)] Given a map from $\Spec(k[[t]])$ to $F_2$, what is the local
    structure of the corresponding equation?
  \item[(3)] What is the effect of taking limits from the typical case of a
    minimal type to more special cases?
\end{itemize}
We will not consider (3) here (except in Section \ref{sec:elldiff} for the
elliptic case), as for most applications it is the typical behavior that
matters.  Presumably, the only effect of taking such limits is to make
changes on the order of the $o()$ terms below.

Of the remaining problems, the simpler to deal with is (2).  Supposing for the
moment that the fiber of interest is $x=\infty$ (which will be the main
interesting case in any event), then the map to $F_2$ takes the form
$(y,x,w)=(Y(t),1,W(t))$ for suitable power series $Y$ and $W$.  Moreover,
$W(t)$ cannot be identically 0, since then the image would be contained in
the fiber.  The key point now is that we may view the formal neighborhood
of the fiber as a $\P^1$-bundle over $\Spec(k[[w]])$, and perform our usual
calculations to convert between sheaves on the $\P^1$-bundle and morphisms
on the anticanonical curve.  In our case, the sheaf is the structure sheaf
of an affine curve, and is thus simply $k[[t]]$ viewed as a $k[[w]]$-module
via the map $w\mapsto W(t)$.  Since this is a totally ramified extension,
$k[[t]]$ is a free $k[[w]]$-module with basis $(1,t,\dots,t^{\ord(W)-1})$,
and the action $M_t$ of multiplication by $t$ in this basis is the
companion matrix of the minimal polynomial of $t$ over $k[[w]]$.  The
associated morphism of vector bundles on $F_2$ is then simply $B:=y -
Y(M_t)\in \Mat_n(k[[w]][y])$.  Note that in characteristic 0 (or greater
than $\ord(W)$), we may change variables so that $W(t)=t^{\ord(W)}$.

This order $\ord(W)$ has a direct geometric interpretation: since $W=0$ is
a fiber, $\ord(W)=c_1(M)\cdot f$.  As a result, the combinatorics let us
compute the generic value of $\ord(W)$ by writing the pullback of $f$ as a
sum of components.  In particular, we see that $\ord(W)$ is generically
equal to the sum of (a) the multiplicity of the last exceptional curve $e$
in this expression (since we can generically decompose into curves with
$c_1(M)\cdot e=1$) and (b) the multiplicity with which $\supp(M)$ meets the
strict transform of the fiber $W=0$.  Any point in which the spectral curve
meets the strict transform of $W=0$ must also be a point of $e$, and thus
generically the contribution from (b) is 0.  Note that since we are
excluding irrelevant blowups, we can never resume blowing up points on the
strict transform of $W=0$ once we stop.  We thus see that the strict
transform has class $f-e_1$ or $f-e_1-e_2$, since any further blowups would
make it have negative intersection with the anticanonical curve.  Moreover,
in either case any further blowups give exceptional curves disjoint from
the strict transform, so irreducibility actually forces the additional
contribution to vanish.  In particular, we find that $\ord(W)\ge
c_1(M)\cdot e_1$ with equality unless (a) we are blowing up a smooth
unramified point once or (b) we are blowing up a (possibly singular)
ramified point.  Moreover, in case (b), we have $\ord(W)\ge c_1(M)\cdot
(e_1+e_2)$ except when we blow up the ramified point exactly twice.  Note
that in terms of the parametrization, the first inequality is trivial, as
$c_1(M)\cdot e_1 = \min(\ord(W),\ord(Y))$, and the genericity condition
becomes $\ord(W)\le \ord(Y)$ except when the strict transform has type
$f-e_1-e_2$.  Moreover, since $f-e_1-e_2$ is a $-2$-curve, having a strict
transform of this form is itself a codimension 1 condition (for any fixed
type of surface), and thus again does not happen generically.

If the point $(0,1,0)$ we are blowing up is a typical (i.e., not ramified
for the involution) smooth point of the fiber, then the anticanonical curve
(formally) locally has the form $y^2-y=0$, and thus the spectral and
anticanonical curves meet with multiplicity $\ord(Y)$.  In the generic case
$\ord(W)=1$, so WLOG $W=t$, this gives $B$ in the form $y-Y(t)$, so a zero
of multiplicity $\ord(Y)$, which in turn gives rise to a zero of $A$ of
that multiplicity and a pole of the same multiplicity at the image under
the involution.  (The ramified case is more complicated, as in that case
the anticanonical curve in odd characteristic has the typical form $y^2=w$,
forcing $\ord(W)>1$ when the multiplicity of intersection is $>1$.  When
the multiplicity is 1, the resulting equation has $A=-1$ at the
ramification point, while the typical higher order case is $2\times 2$ with
both a zero and pole there.  One has the same overall structure in
characteristic 2, with the only difference being that the wild ramification
allows $A$ to remain holomorphic in some higher order cases.)

The next simplest case to consider is actually the differential case, where
we assume characteristic 0.  Here the anticanonical curve is $y^2=0$, and
the reduced curve has a natural parametrization $(0,1,1/z)$, so that in
terms of the parameter we have $z=W(t)^{-1}=t^{-\ord(W)}$.  The resulting
Higgs bundle then has the form $A(z)^{-t} = -Y(M_t)|_{w=1/z}$.  Note that
we may view the resulting equation as one over the field of Puiseux series;
if $W(t)=t^a$, $\ord(Y(t))=b$, then the equation has the form
\[
\frac{v'(z)}{v(z)} = z^{-b/a} \sum_{0\le l} c_l z^{-l/a-2},
\]
where $t^b/Y(t) = \sum_{0\le l} c_l t^l$.  (Here the factor of $z^{-2}$
comes from the fact that we need to work relative to a differential which
is holomorphic at $z=\infty$.)  This should be compared with the
classification of the local behavior of differential equations
\cite{vandenEssenA/LeveltAHM:1992}, though we will see below that the
combinatorics of the surface actually encodes more subtle information about
which coefficients are nonzero.

The nonsymmetric $q$-difference case is the next simplest, as the
anticanonical curve is $y^2-xwy=0$, and thus we can easily parametrize both
branches:
\[
(y,x,w) = (0,1,1/z)\quad\text{and}\quad (y,x,w) = (z,1,z);
\]
here we have parametrized so that the singular point is $\infty$ on one
branch and $0$ on the other branch.  In this case the equation has matrix
equivalent to
\[
A(t) = \frac{Y(t)}{Y(t)-W(t)}\bigg|_{t\mapsto M_t}.
\]
(Here we use the fact that we only care about the conjugacy class of $M_t$
over $k[[w]]=k[[1/z]]$, and $M_t$ and its transpose are conjugate.)

The nonsymmetric difference case is analogous, with the parametrizations
now being
\[
(y,x,w) = (0,1,1/z)\quad \text{and}\quad (y,x,w)=(z^2,1,z),
\]
of the anticanonical curve $y^2-w^2y=0$, giving equation
\[
A(t) = \frac{Y(t)}{Y(t)-W(t)^2}\bigg|_{t\mapsto M_t}.
\]

The symmetric $q$-difference case is only slightly more complicated.  The
two branches on the anticanonical curve of the corresponding double cover
correspond to the parametrizations
\[
(y,x,w) = (\eta z/(z^2+\eta)^2,1,z/(z^2+\eta))\quad\text{and}\quad
(y,x,w) = (z^3/(z^2+\eta)^2,1,z/(z^2+\eta))
\]
of the irreducible anticanonical curve $y^2-xwy+\eta w^4 = 0$.
The key point is that although $w$ is a quadratic function of $z$, the
corresponding extension of $k[[w]]$ is unramified.  We may thus let $Z(t)$
be the non-holomorphic solution of
\[
Z(t)+\frac{\eta}{Z(t)} = \frac{1}{W(t)}
\]
in $k((t))$ and view $M_t$ as a matrix over $k[[1/z]]$ via $w\mapsto
z/(z^2+\eta)$, to obtain an equation of the form
\[
A(t)
=
\frac{\eta Z(t)/(Z(t)^2+\eta)^2-Y(t)}
     {Z(t)^3/(Z(t)^2+\eta)^2-Y(t)}\bigg|_{t=M_t}.
\]

The symmetric difference case is more complicated, as now the
anticanonical curve is ramified over the base of the
ruling, with natural parametrizations (assuming characteristic 0)
\[
(y,x,w) = (1/z^3,1,1/z^2)\quad\text{and}\quad (y,x,w) = (-1/z^3,1,1/z^2)
\]
of the anticanonical curve $y^2 = w^3 x$.  Note that as before, since the
characteristic is 0 we may assume that $W(t)=t^a$.  Let $Y(t) = \sum_{b\le
  l} c_l t^l$, with $c_b\ne 0$.  Then the eigenvalues of $Y(M_t)$ are the
Puiseux series of the form
\[
\sum_{b\le l} c_l z^{-2l/a},
\]
one for each of the $a$ roots $z^{1/a}$, and thus the eigenvalues of $A(z)$
are the different values of
\[
\frac{1/z^3-\sum_{b\le l} c_l z^{-2l/a}}
     {-1/z^3-\sum_{b\le l} c_l z^{-2l/a}}.
\]
Note that if $a$ is odd, then the eigenvalues are all $\pm 1$ at
$z=\infty$, depending on the sign of $3-2b/a$.

\medskip

For (1), the answer is in one sense trivial: given a (spectral) map from
$\Spec(k[[t]])$ to $F_2$, we simply repeatedly blow up the image of the
closed point until the closed point is no longer on the anticanonical
curve.  This is reasonably straightforward, but as stated tells us little
about the inverse problem.  After each blowup, we obtain a patch of $X_m$
containing the image of the special fiber, and a morphism from that patch
to $X_{m-1}$.  Apart from the parametrized spectral curve itself, the other
main piece of information on each patch is the equation of the
anticanonical curve.  A key observation is that after the first few
blowups, there are at most two (local) components of the anticanonical
curve, and those components (which may have multiplicity in $C_\alpha$)
meet with multiplicity 1.  In other words, after a few steps, we can
coordinatize the resulting patch so that the anticanonical curve has the
equation $u^{m_1}v^{m_2}=0$, where the most recent exceptional curve is
$u=0$.  (If $m_2=0$, there is some ambiguity in the coordinate $v$, which
we reduce by insisting that the point being blown up is the origin.)  Thus
the question factors into understanding the first few steps and
understanding this general scenario.  Note that since in the minimal case
we are assuming that we never blow up a point of an exceptional curve which
is not an anticanonical component, we may assume that $m_1>0$.

Let $u(t)$, $v(t)\in t k[[t]]$ be the parametrization of the spectral curve
in this patch.  Note that we may assume that $u(t)^{m_1}v(t)^{m_2}$ is not
identically 0, since otherwise the spectral curve would not be transverse
to the anticanonical curve.  After blowing up the origin, there are three
possibilities for the new affine patch.  The new exceptional curve is
$\hat{u}=0$, and the other relevant component (if any) of the anticanonical
curve is $\hat{v}=0$.
\begin{itemize}
  \item[(a)] If $\ord(u)>\ord(v)$, then the new patch has coordinates
    $\hat{u}(t) = v(t)$, $\hat{v}(t)=u(t)/v(t)$, with
    $\hat{m}_1=m_1+m_2-1$, $\hat{m}_2=m_1$.
  \item[(b)] If $\ord(u)<\ord(v)$, then the new patch has coordinates
    $\hat{u}(t) = u(t)$, $\hat{v}(t)=v(t)/u(t)$, with
    $\hat{m}_1=m_1+m_2-1$, $\hat{m}_2=m_2$.
  \item[(c)] If $\ord(u)=\ord(v)$, then the new patch has coordinates
    $\hat{u}(t)=u(t)$, $\hat{v}(t)=v(t)/u(t)-\alpha$, where $\alpha=\lim_{t\to
    0} v(t)/u(t)$, $\hat{m}_1=m_1+m_2-1$, $\hat{m}_2=0$.
\end{itemize}
Note that if $m_2=0$, then (b) and (c) give the same combinatorial type, so
we merge them into a case (bc), essentially allowing $\alpha=0$ in (c).
Geometrically, these cases correspond to the case that the spectral curve
is tangent to $u=0$, tangent to $v=0$, or tangent to neither, respectively.
Every time we hit case (c) or (bc) we pick up a parameter $\alpha$, which
is in $k^*$ when $m_2\ne 0$ and $k$ when $m_2=0$.

If $m_1=1$, $m_2=0$, then the intersection with the anticanonical curve
depends only on the order of vanishing $\ord(u)$.  Note here that the
choice of point to blow up on the exceptional curve was made at the {\em
  previous} step, and thus there are no more parameters (since this blowup
separates the spectral and anticanonical curves).  We find the following by
an easy induction.
\begin{itemize}
  \item[(1)] If $m_1,m_2>0$, then the intersection of the image of the
    structure sheaf with $C_\alpha$ depends only on
    $
    u(t)+o(t^{m_1\ord(u)+(m_2-1)\ord(v)})$ and
    $v(t)+o(t^{(m_1-1)\ord(u)+m_2\ord(v)})$.
  \item[(2)] If $m_1>0$, $m_2=0$, then it depends only on
    $
    u(t)+o(t^{m_1\ord(u)})$ and $v(t) + o(t^{(m_1-1)\ord(u)})$.
\end{itemize}
In each case, we either have $(m_1,m_2)=(1,0)$ or have enough information
to determine in which case (a--c) we are in, and in case (c) what the value
of the parameter is.  Moreover, in each case, we find that we know the new
coordinates to at least the correct precision (using the fact that we know
a ratio to the same {\em relative} precision as we know the numerator and
denominator).  (Of course, this data is in general more than we need to
determine the intersection, since already for $(m_1,m_2)=(1,0)$ we only
care that the leading coefficient of $u$ not vanish.)

Note in particular that if $(m_1,m_2)=(1,1)$, then the combinatorial type
depends only on $\ord(u)$ and $\ord(v)$.  Moreover, we see that the blowing
up process is essentially just performing the subtraction form of Euclid's
algorithm for computing $\gcd(\ord(u),\ord(v))$, and the additional data
determining the intersection is just the constant term of
\[
u(t)^{\ord(v)/\gcd(\ord(u),\ord(v))}v(t)^{-\ord(u)/\gcd(\ord(u),\ord(v))}.
\]
(The minimal case has $\gcd(\ord(u),\ord(v))=1$.)

The nonreduced cases (i.e., $\max(m_1,m_2)>1$) are more complicated to deal
with, and we can only give a satisfying description in characteristic 0.
The advantage of working in characteristic 0 is that we can take roots of
power series.  In particular, we know $u$ to at least enough precision to
know its leading term, and thus we may choose a new coordinate $s$ such
that $u(s) = s^{\ord(u)}$, and find that we now need to specify only the
coordinate $v(s)$, to precision $o(s^{(m_1-1)\ord(u)+m_2(\ord(v))})$.  If
$m_2>0$, we could instead reparametrize so that $v(s)=s^{\ord(v)}$, so that
the only data is $u(s)+o(s^{m_1\ord(u)+(m_2-1)\ord(v)})$.  Note that
knowing
\[
v(s) + o(s^{(m_1-1)\ord(u)+m_2(\ord(v))})
\]
tells us
\[
v(s)^{1/\ord(v)} + o(s^{1+(m_1-1)\ord(u)+(m_2-1)(\ord(v))}),
\]
and thus tells us the inverse function to the same precision.  In other
words, when both parametrizations are defined, the information we obtain is
the same in either case.

To proceed further, we need to understand better how the two
parametrizations are related.  To begin with, let $H\subset \N$ be a
submonoid (i.e., $0\in H$ and $H$ is closed under addition).  A power
series {\em supported on $H$} is then an element of $k[[t]]$ in which the
coefficient of $t^l$ is 0 unless $l\in H$.

\begin{lem}
  For any submonoid $H\subset \N$, the power series supported on $H$ form a
  ring, in which the series with nonzero constant terms are units.  If
  $f\in k[[t]]$ is supported on $H$ with $f(0)=0$, then for any other element
  $g\in k[[t]]$, $g\circ f$ is supported on $H$.
\end{lem}

\begin{proof}
  That the power series supported on $H$ form a ring is trivial, and thus
  any polynomial in such a power series is also supported on $H$.  Since
  the set of power series supported on $H$ is closed under formal limits,
  the claim regarding composition follows.  That the reciprocal of a series
  with nonzero constant term is still supported on $H$ then follows by
  plugging in $f-f(0)$ into the appropriate geometric series.
\end{proof}

More surprisingly, we have the following.

\begin{lem}
  For any submonoid $H\subset \N$, the power series with $f(0)=0$,
  $f'(0)\ne 0$ and $f(t)/t$ supported on $H$ form a group under composition.
\end{lem}

\begin{proof}
  Suppose $g$ is such a series, with
  \[
  g(t) = \sum_{n\in H} c_n t^{n+1},
  \]
  $c_0\ne 0$.  Then
  \[
  g(f(t)) = \sum_{n\in H} c_n (f(t))^{n+1}
  = \sum_{n\in H} c_n t t^n (f(t)/t)^{n+1}.
  \]
  Since $(f(t)/t)^{n+1}$ is supported on $H$ and $t^n$ is supported on $H$,
  we find that these series are indeed closed under composition.

  It remains only to show that the compositional inverse of such a series
  is again of this form.  Let $f(t)$ be such a series, and let $g(t)$ be
  its compositional inverse.  Suppose $g(t)/t$ is not supported on $H$, and
  let $n$ be the minimal nonnegative integer not in $H$ such that the
  coefficient of $t^{n+1}$ in $g(t)$ is nonzero.  We may thus write
  \[
  g(t) = g_0(t) + c_n t^{n+1} + o(t^{n+1}),
  \]
  where $g_0(t)/t$ is supported on $H$.  We then have
  \[
  t = g(f(t)) = g_0(f(t)) + c_n t^{n+1} (f(t)/t)^{n+1} + o(t^{n+1}),
  \]
  so that
  \[
  1-\frac{g_0(f(t))}{t} = c_n t^n (f(t)/t)^{n+1} + o(t^n).
  \]
  The left-hand side is supported on $H$, and thus the coefficient of $t^n$
  on the left vanishes, while on the right it equals
  \[
  c_n f'(0)^{n+1}\ne 0,
  \]
  giving the desired contradiction.
\end{proof}

In particular, the monoid generated by the exponents of the nonzero terms
of $s^{-\ord(v)}v(s)$ in the $u(s)=s^{\ord(u)}$ parametrization will be the
same as that coming from the other parametrization, when both
parametrizations are valid, and this monoid will be unchanged under blowup
as long as we remain in scenarios (a) and (b).  Unfortunately, scenario (c)
can change the monoid in somewhat unpredictable ways.  Luckily, we can
enlarge the monoid without losing information about the blowup.

Call a subset $H\subset \N$ ``arithmetic'' if it contains $0$ and for all
$n\in H$, also contains $n+\gcd(\{l:l\in H|l\le n\})$.  Such a set is
certainly a monoid, and we can still consider the arithmetic set generated
by the exponents of nonzero terms of a power series.  Note that an
arithmetic set has a minimal set of generators: $d_1$ is the minimal
nonzero element, $d_2$ the minimal element not a multiple of $d_1$, $d_3$
the minimal element not a multiple of $\gcd(d_1,d_2)$, etc.

We associate an arithmetic set to our configuration $u(t),v(t)$ as follows.
If $m_1,m_2>0$, then we change coordinates so that $u(s)=s^{\ord(u)}$, and
take the arithmetic set generated by $\gcd(\ord(u),\ord(v))$ and the
exponents of the nonzero terms of $s^{-\ord(v)}v(s)$.  (Here, by
convention, we suppose that all terms past the required precision are
nonzero.)  If $m_1>0$, $m_2=0$, then we take the same coordinate change,
but now the arithmetic set is generated by $\ord(u)$ and the exponents of
nonzero terms of $v(s)$.

For $m_1,m_2>0$, we find that blowing up has no effect on the arithmetic
set, while for $m_2=0$, scenario (bc) replaces the arithmetic set by
$(H-\ord(u))\cap \N$,
while scenario (a) replaces it by
$
H\cup \gcd(\ord(u),\ord(v))\N$.
Moreover, for $m_2=0$, we can detect which of (a) or (bc) we are in knowing
only $\ord(u)$ and $H$; scenario (a) is the case that $H$ contains a
nonzero element less than $\ord(u)$.

\medskip

Let us now see how this information gets used for the different types of
equations.  The $q$-difference cases are now the simplest, as the
anticanonical curve can now already be written as $uv=0$ in suitable
(formal) coordinates.  For the nonsymmetric $q$-difference case, the two
coordinates are $y$ and $w-y$, so the surface determines $a:=\ord(Y(t))$,
$b:=\ord(W(t)-Y(t))$, and
\[
\gamma:=\lim_{t\to 0}\frac{Y(t)^b}{(W(t)-Y(t))^{a}}.
\]
If we write $Y(t) = \alpha t^a (1+O(t))$, $W(t)-Y(t)=\beta t^b (1+O(t))$,
then we have three cases.  If $a<b$, then $W(t) = \alpha t^a (1+O(t))$, so
we may as well reparametrize to make $\alpha=1$, $W(t)=t^a$, and thus the
equation has the symbolic form
\[
\frac{v(qz)}{v(z)} = \gamma^{1/a} z^{-1+b/a} (1+O(z^{-1/a})),
\]
where we translate Puiseux series to matrices as above.  (Note that in the
non-minimal case $\gcd(a,b)>1$, the exponent only has denominator
$a/\gcd(a,b)$, but there must be {\em some} lower-order term with
denominator $a$, as otherwise $Y,W$ would not give a parametrization!)
Similarly, in the case $a>b$, we have $W(t) = \beta t^b (1+O(t))$,
and again may take $\beta=1$ to give
\[
\frac{v(qz)}{v(z)} = \gamma^{1/b} z^{1-a/b} (1+O(z^{-1/b})).
\]
Finally, if $a=b$, then $W(t) = (\alpha+\beta) (1+O(t))$, and the equation is
\[
\frac{v(qz)}{v(z)} = \gamma (1+O(z^{-1/a})).
\]
This is essentially the same as the formal classification of $q$-difference
equations given in \cite{RamisJ-P/SauloyJ/ZhangC:2013}.  We obtain the same
three cases in the symmetric $q$-difference case, as the needed
reparametrizations do not affect the leading terms of the numerator and
denominator of $A$.  Note that in contrast to the case of a ``finite''
singularity (i.e., at a smooth point of the anticanonical curve), the
non-minimal singularities here correspond to higher-dimensional matrices,
and similarly for the other cases considered below.

For the differential case, we start with $(m_1,m_2)=(2,0)$ with $u=y$,
$v=w$.  The above description of the spectral curve thus has the form
\begin{align}
  Y(t) &= u(t) = t^a\\
  W(t) &= v(t) = \sum_{0\le l\le a-b} c_l t^{l+b} + o(t^a),
\end{align}
with $c_0\ne 0$, where the corresponding arithmetic set is determined by
the degrees where the $\gcd$s of the exponents up to that point (with $a$)
drop.  (In the minimal case, the $\gcd$ actually reaches $1$ before we hit
$o(t^a)$; in general, we do not care where the $\gcd$ drops past there.)
Moreover, our genericity assumptions imply $b\le a$, and thus we could
instead reparametrize to have
\[
Y(t) = t^a(\sum_{0\le l\le a-b} e_l t^l+o(t^{a-b})),\quad
W(t) = t^b,
\]
where the $e_l$ are determined from $c_l$ and generate the same arithmetic
set.  This gives an equation of the form
\[
\frac{zv'(z)}{v(z)} = \sum_{0\le l\le a-b} f_l z^{(a-l)/b-1} +o(1),
\]
where again the arithmetic set is the same.  Note that we can record the
combinatorial information in the equivalent form of specifying the leading
exponent and each later exponent where the common denominator increases.
Also, although the translation from the parameters of the surface to the
coefficients $f_l$ is fairly complicated, there is a one-to-one
correspondence between those parameters and the potentially nonzero
coefficients of the Puiseux series, such that each $f_l$ depends only the
parameters up to that point, and is degree 1 in the parameter corresponding
to $l$, unless the parameter lies in $k^*$ in which case it is proportional
to an appropriate power of the parameter.  In particular, if we are
considering two irreducible singularities that share the first few steps of
the blowup, then the corresponding Puiseux series must agree in precisely
the first $N$ terms (among those allowed to be nonzero), where $N$ is the
number of shared parameters.

For instance, suppose we are given an irreducible singularity of the form
\[
\frac{zv'(z)}{v(z)}
=
f_0 z^3 + f_6 z^2
+f_{12}z+f_{14}z^{2/3}+f_{15}z^{1/2}+f_{16}z^{1/3}
+f_{17}z^{1/6}+f_{18} + o(1).
\]
Assuming $f_0$, $f_{14}$, $f_{15}$ are nonzero, then the combinatorial data
is determined by the sequence $(3,2/3,1/2)$, and the corresponding moduli
space can be identified with $(k^*)^3\times k^5$.  The corresponding
sequence of $(m_1,m_2)_{a/b/c}$ is:
\begin{align}
&(2,0)_a(1,2)_b(2,2)_b{\bf (3,2)_c}{\bf (4,0)_{bc}}{\bf (3,0)_{bc}}\notag\\
&(2,0)_a(1,2)_b{\bf (2,2)_c}\notag\\
&(3,0)_a{\bf (2,3)_c}{\bf (4,0)_{bc}(3,0)_{bc}(2,0)_{bc}(1,0)_{bc}},\notag
\end{align}
where we have indicated the steps introducing parameters in bold.  We
should recall here that the subscript on a tuple $(m_1,m_2)_{a/b/c}$
determines the {\em next} point to be blown up.  Note that the three
parameters from $k^*$ are
\[
-\frac{1}{f_0}, \frac{f_0^6}{f_{14}}, \frac{f_{14}^9}{9 f_0^{14}f_{15}^2},
\]
and the remaining parameters in $k$ are even more complicated functions of
the coefficients.  The anticanonical curve on the resulting surface has the
decomposition
\begin{align}
&2(s+f-e_1-e_2-e_3-e_4)+(e_1-e_2)+2(e_2-e_3)+3(e_3-e_4)+4(e_4-e_5)+3(e_5-e_6)\notag\\
{}+{}&2(e_6-e_7-e_8-e_9)+(e_7-e_8)+2(e_8-e_9)\notag\\
{}+{}&3(e_9-e_{10}-e_{11})+2(e_{10}-e_{11})+4(e_{11}-e_{12})+3(e_{12}-e_{13})+2(e_{13}-e_{14})+(e_{14}-e_{15}).\notag
\end{align}
Note that although setting
a coefficient other than $f_0$, $f_{14}$, $f_{15}$ to 0 has no effect on
the combinatorics, setting one of the critical coefficients to 0 can make
significant changes.  For instance, the generic subcase with $f_{15}=0$ has
one fewer blowup, so this does not correspond to a degeneration of surfaces
as considered above.

The one complication in the nonsymmetric difference case is that we do not
start in a configuration $u^{m_1}v^{m_2}=0$, since the two branches are
tangent.  We must blow up twice in order to achieve this, and this leads to
several cases.  Note first that blowing up a general point of $e_1$ yields
an equation
\[
\frac{v(z+1)}{v(z)} = 1+\alpha/z + o(1/z),
\]
$\alpha\ne 0$, which leaves the case that we blow up the triple
intersection.  This yields an anticanonical decomposition
$(s+f-e_1-e_2)+(s+f-e_1-e_2)+(e_1)+2(e_2)$, and there are four
possibilities for the next blowup: a general point $\lambda$ of $e_2$, or
the intersection with one of the other three components.  The resulting
patches, with the equation of the anticanonical curve and the map to the
original coordinates are:
\begin{itemize}
\item[(1)] $u^2=0$, where $y=u^2 (v+\lambda)$, $w=u$.
\item[(2)] $u^2 v = 0$, where $y=u^2 v$, $w=u$.
\item[(3)] $u^2 v = 0$, where $y=u^2 (v+1)$, $w=u$.
\item[(4)] $u^2 v =0$, where $y=u^2 v$, $w=uv$.
\end{itemize}
In case (1), if we take $u(t) = t^a$, $v(t) = \sum_{b\le l\le a} c_l t^l+o(t^a)$,
then the resulting equation has the symbolic form
\[
\frac{f(z+1)}{f(z)}
=
\frac{\lambda+v(z^{-1/a})}{\lambda-1+v(z^{-1/a})} + o(z^{-1})
=
\frac{\lambda}{\lambda-1} + \sum_{b\le l\le a} d_l z^{-l/a} + o(z^{-1}),
\]
where $d_b\ne 0$ and the common denominator increases at the same places as
$v(z^{-1/a})$.  In case (2), we now know $v(t)$ to relative precision
$o(t^a)$, and now the equation is
\[
\frac{f(z+1)}{f(z)}
=
\frac{v(z^{-1/a})}{v(z^{-1/a})-1}(1+o(z^{-1}))
=
\sum_{b\le l\le a+b} d_l z^{-l/a} + o(z^{-b/a-1}),
\]
where again the common denominator increases at the same points as
$v(z^{-1/a})$.  Case (3) is analogous, and simply gives the dual equation.
Finally, for case (4), we want to parametrize so that $u(t)v(t)=t^a$,
and have $v(t) = \sum_{b\le l\le a} c_l t^l + o(t^a)$ with $c_b\ne 0$,
giving the equation
\[
\frac{f(z+1)}{f(z)}
=
\frac{1}{1-v(z^{-1/a})}
=
1 + \sum_{b\le l\le a} d_l z^{-l/a} + o(z^{-1}),
\]
now with $a>1$, since $e_3$ is a multiple component of the fiber in this
case.  Again, these agree with the formal classification of difference
equations \cite{TurrittinHL:1960,PraagmanC:1983}, apart from the
description of the combinatorial data.

The symmetric difference case is, naturally, even more complicated, as now
we may need to blow up three times to have just a pair of transverse
branches to consider.  The simplest case is that we blow up a general point
of $e_1$, which gives an equation of symbolic form
\[
\frac{f(z+1)}{f(z)} = \exp(\alpha/z+O(1/z^3)).
\]
(Here we ignore the effect that the shift has on the symmetry condition.)
Similarly, blowing up the point of tangency then a general point of $e_2$
gives an equation of symbolic form
\[
\frac{f(z+1)}{f(z)} = -\exp(\alpha/z+O(1/z^3)).
\]
Otherwise, we are blowing up the triple intersection, and again have four cases.
\begin{itemize}
\item[(1)] $u^2=0$, where $y=u^3 (v+\lambda)^2$, $w=u^2(v+\lambda)$.
\item[(2)] $u^2 v = 0$, where $y=u^3 v^2$, $w=u^2 v$.
\item[(3)] $u^2 v = 0$, where $y=u^3 (v+1)^2$, $w=u^2(v+1)$.
\item[(4)] $u^2 v =0$, where $y=u^3 v$, $w=u^2v$.
\end{itemize}

In case (1), we reparametrize so that $u(t)^2(\lambda+v(t))=t^{2a}$,
corresponding to $t=z^{-1/a}$, giving symbolic equations of the form
\[
\frac{f(z+1)}{f(z)}
=
\frac{1-\sqrt\lambda}{1+\sqrt\lambda}
\exp(\sum_{1\le l\le a} c_l z^{-l/a}+o(1/z)).
\]
Of course, the symmetry means that these equations come in pairs,
corresponding to the fact that the full equation has even order.
Case (3) is the next simplest, as we may still reparametrize so that
$u(t)^2(1+v(t))=t^{2a}$, giving again pairs of symbolic forms typified by
\[
\frac{f(z+1)}{f(z)}
=
\frac{\sqrt{1+v(z^{-1/a})}-1}{\sqrt{1+v(z^{-1/a})}+1}
=
\sum_{b\le l\le a+b} c_l z^{-l/a} + o(z^{-1-b/a})
\]
For (2) and (4), the natural parametrization is $z^{-2}=u(t)^2v(t) =
t^{2a+b}$, where $u$ and $v$ have order $a$ and $b$ respectively, and are
known to relative precision $o(t^a)=o(z^{-2a/(2a+b)})$.  Case (4) then
gives an equation of the form
\[
\frac{f(z+1)}{f(z)}
=
\exp(\sum_{0\le l\le a} c_l z^{-(b+2l)/(2a+b)} + o(1/z)).
\]
If $b$ is odd, this equation satisfies the symmetry condition, while if $b$
is even, then the overall order is even and the symbolic equations come in
pairs.

\medskip

It is worth noting in each case how the corresponding canonical
isomonodromy transformations (the noncommutative analogue of pseudo-twists)
behave.  Note that when performing the isomonodromy transformation
corresponding to twisting by $e_l$, the effect is to perform the $d$-th
iterate of the twist by $e_m$, where $d$ is the multiplicity of $e_m$ in
the pullback of $e_l$, since twisting by a component of $C_\alpha$ has no
effect.  We thus need only consider the twist by $e_m$.  The result is then
to replace the sheaf by the sheaf of functions vanishing at that point, and
is thus easily seen to correspond to gauging by $M_t$ (which we may view
symbolically as the appropriate power of $z$).  In the $q$-difference
cases, gauging by $z^{-1/a}+o(z^{-1/a})$ multiplies the leading coefficient
by $q^{-1/a}$; since the true parameter is the $a$-th power of the leading
coefficient, this shifts the parameter by $q$ as expected.  In the
nonsymmetric ordinary difference case, gauging by $z^{-1/a}+o(z^{-1/a})$
multiplies $A$ by $1-1/az+o(1/z)$, and similarly for the symmetric case.
And of course in the differential case gauging by $z^{-1/a}$ adds
$-1/a+o(1)$ to $zf'/f$.

In the differential case, we of course also have continuous isomonodromy
transformations, which in the Puiseux form correspond to equations
\[
\frac{f_u}{f} = \frac{a}{l+a} z^{l/a+1}
\]
giving
\[
\frac{d}{du} A = z^{l/a}.
\]
In particular, for a unibranch type with $l$ parameters, we obtain $l-1$ such
deformations.  (There is also a deformation changing the location
of the singularity.)

In the nonsymmetric difference case, if we write the symbolic equation as
\[
\frac{f(z+1)}{f(z)} = z^{b/a} \exp(g(z^{-1/a}))
\]
then the equation
\[
f_u(z) = \frac{a}{a+l} z^{1+l/a} f(z)
\]
gives
\[
\frac{d}{du} g(z^{-1/a}) = z^{l/a}(1+o(1/z))
\]
again giving continuous deformations for every parameter but the last.
The symmetric difference case is analogous; the only change is that the
continuous equation needs to preserve the symmetry, but this is easy to
arrange.

If we take into account the symmetries (e.g., the $\PGL_2$ symmetry in the
differential case), we find in general that the net number of {\em local}
continuous isomonodromy deformations is 0 if the anticanonical curve is
reduced, and otherwise can be expressed as $1-C_\alpha\cdot
(C_\alpha-C_\alpha^{\text{red}})$.  This follows by an easy induction:
blowing up a point of the smooth locus of $C_\alpha$ subtracts $e_m$ from
$C_\alpha^{\text{red}}$, and similarly for blowing up a point on two
components, while blowing up a point on a single component with nontrivial
multiplicity leaves $C_\alpha^{\text{red}}$ alone.  One thus reduces to the
cases in which $C_\alpha$ first becomes nonreduced: $C_\alpha=2(s+f)$,
which has $-3$ deformations due to the symmetries,
$C_\alpha=(s+f-e_1-e_2)+(s+f-e_1-e_2)+(e_1-e_2)+2(e_2)$, with no continuous
deformations, and
$C_\alpha=(2s+2f-2e_1-e_2-e_3)+(e_1-e_2-e_3)+(e_2-e_3)+2(e_3)$, likewise.

The fact that the result depends only on the geometry suggests that there
should be a geometric explanation of these continuous isomonodromy
deformations, just as the discrete isomonodromy deformations come from
twists by line bundles.  Of course, this assumes that the local
deformations can in fact be glued together to form global deformations,
though the results of \cite{OrmerodCM/RainsEM:2017b} suggest that this is
indeed the case, and indeed we establish this in general in Section
\ref{sec:cont_isomonodromy} below.

\chapter{Moduli of sheaves}
\label{chap:moduli_of_comm_sheaves}

\section{Rational surfaces}

Let $(X,C_\alpha)$ be an anticanonical rational surface.  Say a coherent
sheaf on $X$ has {\em integral support} if its $0$-th Fitting scheme is an
integral curve on $X$, and it contains no $0$-dimensional subsheaf.

\begin{thm}\label{thm:mod_irr}
  Let $(X,C_\alpha)$ be an anticanonical rational surface over an
  algebraically closed field of characteristic $p$, let $D$ be a divisor
  class with generic representative an integral curve disjoint from
  $C_\alpha$, and let $r$ be the largest integer such that $D\in r\Pic(X)$.
  Then the moduli problem of classifying sheaves $M$ on $X$ with integral
  support, $c_1(M)=D$, $\chi(M)=x$, and $M|_{C_\alpha}=0$ is represented by
  a quasiprojective variety $\Irr_X(D,x)$ of dimension $D^2+2$, with a
  symplectic structure induced by any choice of nonzero holomorphic
  differential on $C_\alpha$, and the natural map to $|D|$ (given by
  $M\mapsto \Fitt_0(M)$) is a Lagrangian fibration where the fibers are
  smooth.  Moreover, $\Irr_X(D,x)$ is unirational if the generic
  representative of $D$ has no cusp, separably unirational if $p=0$ or
  $\gcd(x,r,p)=1$, and rational if $x\bmod r\in \{1,r-1\}$.  Finally, if
  $\gcd(x,r)=1$, then there exists a universal sheaf over $\Irr_X(D,x)$.
\end{thm}

\begin{proof}
  Quasiprojectivity follows from the standard GIT construction: for any
  choice of stability condition, a sheaf with integral support is stable.
  The symplectic structure (in particular smoothness) follows from the fact
  that sheaves with integral support are simple (have no nonscalar
  endomorphisms) together with Theorem \ref{thm:poisson}.  This requires a
  choice of Poisson structure on $X$, or equivalently a choice of nonzero
  holomorphic differential on $C_\alpha$; the symplectic structure on the
  moduli space scales linearly with the choice of differential.

  Now, the typical sheaf in the moduli space corresponds to a pair $(C,M)$
  where $C$ is an integral curve of class $D$ (and disjoint from
  $C_\alpha$) and $M$ is a torsion-free sheaf on $C$.  If $g=D^2/2+1$, then
  $\Gamma(\sO_X(D))$ has dimension $g+1$, and thus the integral curves in
  the linear system form an open subset of a $\P^g$.  We also compute that
  $C$ has arithmetic genus $g$, and thus the fiber over the point
  corresponding to $C$ is a compactification of $\Pic^{x+g-1}(C)$, so has
  dimension $g$.  Moreover, we have a commutative diagram
  \[
  \begin{CD}
    \Ext^1_C(M,M)\otimes \Ext^1_C(M,M)@>>> \Ext^2_C(M,M)
    @VVV @VVV\\
    \Ext^1_X(M,M)\otimes \Ext^1_X(M,M)@>>> \Ext^2_X(M,M)
  \end{CD}
  \]
  and thus the pairing on the tangent space at $M$ to its fiber factors
  through $\Ext^2_C(M,M)$, which vanishes as long as $M$ is invertible, so
  that the tangent space is not only half-dimensional, but isotropic, and
  thus Lagrangian as required.  (When $M$ is not invertible, the tangent
  space to the fiber is too large.)

  If $x=1$, then $\deg(M)=g$, and thus the generic such sheaf has a unique
  global section.  The quotient by the corresponding trivial subsheaf is a
  $0$-dimensional sheaf of degree $g$, thus giving a birational
  correspondence with the punctual Hilbert scheme $X^{[g]}$.  Since
  symmetric powers of rational surfaces are rational varieties
  \cite{MattuckA:1968}, it follows that $\Irr_X(D,1)$ is rational.

  More generally, since $D/r$ is a primitive element of the Picard lattice
  of $X$, there exists a divisor $D'$ such that $D\cdot D'=r$.  In
  particular, we can twist by powers of this divisor to obtain isomorphism
  $\Irr_X(D,x)\cong \Irr_X(D,r+x)$ for any $x$.  Similarly, the duality
  morphism $M\mapsto \sHom(M,\omega_C)$ on invertible sheaves can be
  defined globally (and extended to torsion-free sheaves) by $M\mapsto
  \sExt^1(M,\omega_X)$, and gives an isomorphism $\Irr_X(D,x)\cong
  \Irr_X(D,2g-2-x)$.  Since $2g-2=D^2$ is a multiple of $r^2$, the
  rationality claim follows.

  For unirationality, note that since $2g-2$ is a multiple of $r^2$, $g-1$
  is a multiple of $r$, and thus $\Irr_X(D,2-g)$ (classifying sheaves of
  degree 1) is rational.  Since the generic sheaf is an {\em invertible}
  sheaf on the generic curve, we can take its $d$-th power and thus obtain
  a rational map $\Irr_X(D,2-g)\to \Irr_X(D,d+1-g)$ for any $d$.  Since the
  generic curve $C_{\text{gen}}$ is integral, this map is dominant unless
  $C_{\text{gen}}$ is cuspidal and $d$ is 0 in $k$.  This implies
  unirationality in the noncuspidal case for any $d\ne 0$; the case $d=0$
  reduces to the case $d=g-1$ by twisting.  When $\gcd(d,p)=1$, the
  multiplication by $d$ map is separable, and again we may feel free to add
  multiples of $r$ to make this happen.

  Finally, we note that the obstruction to the existence of a universal
  sheaf on $\Irr_X(D)$ is given by a class in $H^2(\Irr_X(D),\G_m)$, which
  can be constructed as follows.  There certainly exists a universal sheaf
  \'etale-locally, and we may twist by a line bundle to ensure that this
  universal sheaf is acyclic; i.e., that $x\ge g$.  Then, since the fibers
  of the universal sheaf are simple, the endomorphism ring of the direct
  image descends to an Azumaya algebra on $\Irr_X(D)$, the class of which
  is the desired obstruction.  This Azumaya algebra has degree $x$, and
  thus the obstruction has order dividing $x$.  Since this is true for any
  twist of sufficiently large degree, we also find that the order of the
  obstruction divides $x+r$, and thus that it divides $\gcd(x,r)$.  In
  particular, if $x$ and $r$ are relatively prime, there is no obstruction.
\end{proof}

\begin{rems}
  The generically cuspidal case can of course only occur in finite
  characteristic, but can certainly occur there, say if $D$ is the class of
  a fiber in a rational quasi-elliptic surface.
\end{rems}

\begin{rems}
  Often in the literature, one restricts ones attention to the subscheme
  where $C$ is not just integral but smooth, making the fibers of the
  Lagrangian fibration abelian varieties.  Of course, this is problematical
  in finite characteristic, where there may not be any smooth curves in the
  linear system.  In addition, since singularity is a codimension 1
  condition, this removes an entire hypersurface from the moduli space,
  based on a condition which is rather unnatural from the difference
  equation perspective.  (Indeed, as we mentioned, difference equations
  correspond most naturally to sheaves on noncommutative surfaces, and
  there the notion of support fails altogether.  In contrast, the failure
  of integrality corresponds to reducibility of the equation in a suitable
  sense.)  For instance, in the generic $2$-dimensional case, both the
  surface and the moduli space are elliptic surfaces, and there are 12
  fibers where the support is singular.  Similarly, there are 12 points of
  the moduli space where the sheaf is not invertible on its support, again
  an odd condition in terms of difference equations.
\end{rems}

\begin{rems}
  It seems likely that the condition $\gcd(x,r)=1$ for the existence of a
  universal sheaf is necessary.  Some condition is needed, see the remark
  following Proposition \ref{prop:ell_surf_rat} below.
\end{rems}

\begin{rems}
  Note that although the proof implicitly gives a parametrization in the
  rational case, the parametrization is likely to be relatively
  complicated, if for no other reason than that the parametrization of
  symmetric powers is already fairly complicated.
\end{rems}

The rational case $x=1$ is particularly nice for another reason: although
the definition of stability generally requires the choice of an ample
bundle, it turns out that when $\chi=1$, this choice is irrelevant.  One
finds in this case $M$ is stable iff any nonzero quotient of $M$ has
positive Euler characteristic (and there are no strictly semistable
sheaves).  We thus find that $\Irr_X(D,1)$ extends naturally to a
projective moduli space.  This space is no longer symplectic, but since
every sheaf in the space is stable, thus simple, it still inherits a
Poisson structure, with smooth symplectic leaves determined by the
quasi-isomorphism class of the complex $M\otimes^{\dL} \sO_{C_\alpha}$.  In
particular, the open subvariety where $M|_{C_\alpha}=0$ is still smooth and
symplectic.

For our purposes, the most natural case is $x=D\cdot f$.  Indeed, the sheaf
corresponding to a standard matrix-form difference equation comes from a
sheaf on a Hirzebruch surface with presentation
\[
0\to \rho^*V\otimes \sO_\rho(-1)\to \sO_X^n\to M_0\to 0.
\]
If we twist by $-f$, then both sheaves in the resolution have vanishing
cohomology, and thus $H^*(M_0(-f))=0$; conversely, by Lemma
\ref{lem:cohom_vanish}, any sheaf with $H^*(M_0(-f))=0$ at the least has a
canonical subsheaf with a presentation of the above form.  (In Chapter
\ref{chap:sheaves_from_eq_comm}, we imposed the additional open conditions
$\Hom(M,\sO_f(-1))=\Hom(\sO_f(-1),M)=0$ for all $f$; ignoring those
conditions gives us a natural partial compactification.)

Of course, we do not have a sheaf on a Hirzebruch surface, but rather a
sheaf on some blowup of the Hirzebruch surface.  However, we have the
following fact, by the same spectral sequence argument as Lemma
\ref{lem:cohom_vanish}.

\begin{lem}
Let $\pi:X\to X_0$ be a birational morphism of smooth projective surfaces,
and let $M$ be a $1$-dimensional sheaf on $X$.  Then $H^0(M)=H^1(M)=0$
iff $M$ is $\pi_*$-acyclic and $H^0(\pi_*M)=H^1(\pi_*M)=0$.
\end{lem}

In other words, a sheaf $M$ on $X$ induces a matrix-form difference
equation (up to twisting by $-f$) iff $H^0(M)=H^1(M)=0$.  (Again, it could
fail to be the natural sheaf associated to a difference equation, but this
can be avoided by imposing the additional conditions
$\Hom(M,\sO_g(-1))=\Hom(\sO_g(-1),M)=0$ for any smooth rational curve $g$
contained in a fiber.)  We are thus led to consider the space
$\Irr_X(D,0)$.  Once again, the stability condition turns out to be
independent of the choice of ample bundle: a $1$-dimensional sheaf $M$ on
$X$ with $\chi(M)=0$ is stable iff any proper nontrivial subsheaf has
negative Euler characteristic, and similarly for semistability.  Since we
need semistable sheaves, we do not immediately inherit a Poisson structure,
although this will certainly exist on the complement of the semistable
locus.  (In most cases, this can be resolved by deforming the ample bundle
to break any ties arising in semistability.)

It remains only to consider the condition $H^0(M)=H^1(M)=0$.  By Lemma
\ref{lem:tau_func}, this is the complement of a codimension $1$ condition
on any family of $1$-dimensional sheaves, cutting out a Cartier divisor.
Of course, this is only well-defined outside the semistable locus, but the
generic sheaf is integral, and thus stable, so we still obtain a
well-defined divisor on the projective moduli space.  On the integral
locus, this divisor is just the canonical theta divisor in the relative
$\Pic^{g-1}$, while in general, it is the zero locus of a canonical
global section of a canonical line bundle $\det \dR\Gamma(M)^{-1}$.

We can obtain a whole family of such divisors by noting that for any vector
$v\in D^\perp$, we can twist $M$ by $\sO_X(v)$ without affecting the Euler
characteristic, thus obtaining rational automorphisms of the projective
moduli space (these are only rational, since twisting can affect stability;
but this is an automorphism on the integral locus).  In particular, we
obtain in this way a canonical global section of $\det
\dR\Gamma(M(v))^{-1}$, which we call a ``tau function'' by analogy with
\cite{ArinkinD/BorodinA:2009}.  It is of course a misnomer to call it a
function (just like a theta function is not an algebraic function), but it
is at the very least a convenient way of describing divisors on the moduli
space.  (Similarly, the divisor on the moduli stack of surfaces where a
given divisor class is a $-2$-curve can also be viewed as a tau function,
as can theta functions themselves.)

The main significance of the $\chi=0$ case is that any standard matrix-form
relaxed equation can be turned into a non-relaxed equation by simply
changing $q$ (or $\hbar$), and vice versa.  So although the true moduli
space is different in the non-relaxed case, this gives an open subset where
they can be identified.  This gives some results for free (e.g., it tells
us that when $r=1$, the true moduli space is rational), and is in any event
rather suggestive.  In particular, we will see below that this
identification of sheaves with $R\Gamma=0$ between the two cases extends to
a corresponding identification of subcategories of the derived categories.

\medskip

We now turn to low-dimensional cases.  The moduli space is $0$-dimensional
when $D$ is a $-2$-curve, and in that case, $M$ is uniquely determined by
its Chern classes and the constraint that it be disjoint from $C_\alpha$;
moreover, its image on $X_0$ is similarly determined by its bidegree and
its intersection with $C_\alpha$ (see Theorem \ref{thm:rigid_is_minus2}
above).  Since the moduli spaces are symplectic, the next interesting case
is the $2$-dimensional case $D^2=D\cdot K_X=0$.  The only such divisors in
the fundamental chamber are the classes $-rK_{X_8}$ for integer $r$, such
that $\omega_{X_8}|_{C_\alpha}$ has exact order $r$ in $\Pic(C_\alpha)$.
Now, in that case, $X_8$ is itself an quasi-elliptic surface, and $D$ is
the class of a fiber.  The corresponding moduli spaces are just the
relative Picard varieties of this quasi-elliptic surface.

\begin{prop}\label{prop:ell_surf_rat}
  Let $\psi:X\to \P^1$ be a relatively minimal rational quasi-elliptic
  surface.  Then for any integer $x$, the relative $\Pic^x$ of $X$ over
  $\P^1$ is a rational surface.
\end{prop}

\begin{proof}
  This is essentially a result of
  \cite[Prop.~5.6.1]{CossecFR/DolgachevIV:1989}.  To be precise, that
  Proposition shows that a relatively minimal quasi-elliptic surface is
  rational iff it has at most one multiple fiber (and that ``tame'') and
  its relative Jacobian is rational.  Since the relative Jacobian of the
  relative $\Pic^x$ is isomorphic to the original relative Jacobian, and
  the relative $\Pic^x$ cannot turn non-multiple fibers into multiple
  fibers (or tame multiple fibers into wild multiple fibers), the claim
  follows.
\end{proof}

\begin{rem}
  It is worth noting that although the moduli space here is rational, there
  is in general no universal sheaf, even if we restrict to an open subset
  (or the function field) of the moduli space.  In particular, when $x=0$,
  it follows from Bhatt's appendix to \cite{KrashenD/LieblichM:2008} that
  the obstruction can be identified with the class of $X$ as a torsor over
  the relative Jacobian.  In particular, we find in that case that the
  obstruction is nontrivial whenever $r>1$.
\end{rem}

Since our surfaces are anticanonical, it makes sense to ask which rational
surface one obtains in this way.  That is, if we start with a blowdown
structure on $X$ and a section of the anticanonical linear system, is there
a natural way to choose a blowdown structure on the (minimal proper regular
model of the) relative $\Pic^x$ such that we can compute the new
anticanonical curve and the new morphism $\Z^{10}\to \Pic(C_\alpha')$?  Note
that on a (quasi-)elliptic rational surface, the line bundle
$\omega_X|_{C_\alpha}$ has finite order (say $r$), and thus determines a
subgroup of degree $\gcd(x,r)$.  Since the corresponding bundle for the
relative $\Pic^x$ has order $r/\gcd(x,r)$, the resulting surface should
depend only on the composition $\Z^{10}\to\Pic(C_\alpha)$ with the quotient
by this subgroup.  We will show this in the case $\gcd(x,r)=r$, and give an
explicit description of the surface, in the following section.

There are additional $2$-dimensional cases in which the generic sheaf does
not have irreducible support as defined above, but in which one still has a
fibration by generically abelian varieties.  The key point is that for any
pair $(r,d)$ of integers with $r>0$, $\gcd(r,d)=1$, the moduli space of
stable vector bundles of degree $d$ and rank $r$ on an elliptic curve $E$
can be identified with $\Pic^d(E)$ by \cite{AtiyahMF:1957}.  If $X_8$ is a
quasi-elliptic rational surface with no multiple fibers, then a stable
sheaf on $X_8$ with $c_1=rC_8$ and $\chi=d$ must be supported on a single
fiber, as otherwise it would split as a direct sum of sheaves with
different slopes.  We can similarly rule out the case that its
scheme-theoretic support is nonreduced, and thus see that it must be a
torsion-free sheaf on a single fiber, which must therefore be stable.  So
again we get a $2$-dimensional moduli space.  This suggests that the usual
Lax pairs for discrete (or continuous) Painlev\'e equations should have
analogues for arbitrary $d/r\in \Q$, and indeed this is the case (see
Theorem \ref{thm:Painleve_moduli_spaces}, as well as Section
\ref{sec:Painleve_Lax} for the case $d=1$).

Past the $2$-dimensional cases, it no longer makes much sense to ask which
variety we obtain (since birational geometry is extremely complicated, even
for $4$-folds).  It is fairly straightforward to write down divisor classes
giving such moduli spaces, however; for instance the $4$-dimensional moduli
spaces correspond (in the fundamental chamber of an even blowdown
structure) to one of
\[
2s+3f-\sum_{1\le i\le 10}e_i
\qquad\text{or}\qquad
4s+4f-2\sum_{1\le i\le 7} e_i-e_8-e_9
\]
The former is always generically integral (assuming of course that the
corresponding line bundle is trivial on $C_\alpha$), while the latter is
generically integral unless $e_8-e_9$ is a $-2$-curve.  Both of these cases
extend to a sequence of moduli problems of dimension $2g$ for arbitrary
$g>1$.  The first case extends to the general problem of classifying
second-order problems with simple singularities (related to generalized
Garnier-type systems), while the second case generalizes the ``matrix
Painlev\'e equations'' of \cite{KawakamiH:2015}.

In general, for any dimension bigger than $2$, there are only
finitely many possibilities for the representative of $D$ in the
fundamental chamber.  Indeed, we may write any divisor $D$ in the
fundamental chamber in the form
\[
D = D_7 +
c(2s+2f-e_1-\cdots-e_8)
-\sum_i \lambda_i e_{8+i},
\]
where $D_7\in \langle s,f,e_1,\dots,e_7\rangle$
is in the fundamental chamber and $\lambda$ is a partition with all parts
at most $c$.  The constraint $D\cdot C_\alpha=0$ becomes
\[
|\lambda| = D_7\cdot C_\alpha
\]
and we have
\[
D^2 = D_7^2 + 2 c D_7\cdot C_\alpha - \sum_{1\le i}\lambda_i^2
    \ge D_7^2 + c D_7\cdot C_\alpha.
\]
Since every generator of the fundamental chamber for $m=7$ has positive
intersection with $C_\alpha$, there are only finitely many pairs $(D_7,c)$
such that $D_7\ne 0$ and $D_7^2 + c D_7\cdot C_\alpha\le d-2$, where $d$ is
our desired value for $D^2+2$.  (The pairs with $D_7=0$ correspond to
the $D^2=0$ case.)

\section{(Quasi-)elliptic surfaces}

Let $(X,\Gamma,C_\alpha)$ be an anticanonical rational surface with
blowdown structure $\Gamma$, and suppose that $X$ is (quasi-)elliptic, so
that the divisor $rC_\alpha$ is the class of a fiber of a genus 1 pencil
for some positive integer $r$.  Since $-rK$ is the class of a pencil, we
see that the line bundle $\omega_X^r|_{C_\alpha}$ is trivial.  Moreover, if
$\omega_X|_{C_\alpha}$ had order $s$ strictly dividing $r$, then
$sC_\alpha$ would already have linearly independent sections, contradicting
the assumption on $rC_\alpha$.  We thus find that in this scenario,
$\omega_X|_{C_\alpha}$ has exact order $r$.  Per Proposition
\ref{prop:Pic0_acts}, such a torsion bundle induces an automorphism of
$C_\alpha$ of order $r$, and we can quotient by this automorphism to obtain
a new curve $C'_\alpha$.  The pair $(X,\Gamma)$ is determined from the map
$\Lambda_{10}\to \Pic(C_\alpha)$; if we compose with the degree-preserving
map $\Pic(C_\alpha)\to \Pic(C'_\alpha)$, we obtain a new map
$\Lambda_{10}\to \Pic(C'_\alpha)$.  This new map has the same combinatorial
structure as the original map (twisting by $\omega_X|_{C_\alpha}$ preserves
degrees, so the automorphism preserves components), and thus itself arises
from a unique triple $(X',\Gamma',C'_\alpha)$.

\begin{thm}\label{thm:Picr}
  Suppose $C_\alpha$ is reduced.  Then the surface $X'$ constructed in this
  way is the minimal proper regular model of the relative $\Pic^r$ of $X$,
  in such a way that the fiber corresponding to $rC_\alpha$ is $C'_\alpha$.
\end{thm}

\begin{proof}
  Assume for the moment that the sublattice of $\Lambda_{E_8}$
  corresponding to $C_\alpha$ is saturated; this excludes only two cases,
  namely $\Lambda_{A_7}\subset \Lambda_{E_7}\subset \Lambda_{E_8}$ and
  $\Lambda_{A_8}\subset \Lambda_{E_8}$, which we will discuss below.  (The
  remaining unsaturated sublattice $\Lambda_{D_8}$ corresponds to a
  nonreduced $C_\alpha$.)  Together with the hypothesis that $C_\alpha$ is
  reduced, this is equivalent to assuming that $\Pic(X)$ is generated by
  the $-1$-classes meeting each component of $C_\alpha$ positively.  In
  particular, this ensures that the homomorphism $\Pic(X)\to \Pic(C_\alpha)$
  is determined by its restriction to such classes.  A $-1$-class is always
  uniquely effective, the transversality condition implies that the
  corresponding curve meets $C_\alpha$ in a single smooth point, and that
  point in turn determines the image in $\Pic^1(C_\alpha)$.  In particular,
  we can reconstruct $X$ from $C_\alpha$ and the configuration of points in
  which $-1$-classes meet $C_\alpha$, and similarly for $X'$.  As a result,
  to prove the theorem, we will simply need to show that the minimal proper
  regular model of the relative $\Pic^r$ has the correct special fiber, and
  has the relevant $-1$-classes, meeting $C'_\alpha$ in the correct points.

Rather than study the minimal proper regular model directly, we instead
consider the corresponding moduli space of semistable sheaves with first
Chern class $c_1(M)=-rK$ and Euler characteristic $\chi(M)=r$.  Unlike the
cases $\chi\in \{-1,0,1\}$ discussed earlier, in this case the stability
condition depends nontrivially on the choice of ample divisor $\sO_X(1)$.
There are only finitely many divisors with $D$, $-rK-D$ effective (i.e.,
subdivisors of fibers); we may thus choose the ample divisor in such a way
that for any such divisor,  either $D\in\Z K$ or
\[
\frac{D\cdot \sO_X(1)}{K\cdot \sO_X(1)}
\notin
\Z.
\]
This ensures that any semistable sheaf not supported on the special fiber
will be stable; the various inequalities are forced to be strict by
integrality.  This will also force any semistable sheaf supported on the
special fiber to be $S$-equivalent to a sum of stable sheaves supported on
$C_\alpha$, see below. (We should also note that for any fixed $r\ge 1$,
there are only finitely many divisor classes with both $D$ and $-rK-D$
effective on {\em some} triple $(X,\Gamma,C_\alpha)$, and could thus in
principle choose the ample divisor in a uniform way over any family of
(quasi-)elliptic surfaces.)

Since stable sheaves are simple (and stability is an open condition), we
find that the moduli space has an open subset which agrees with an open
subset of the moduli space of simple sheaves.  Any sheaf in the open subset
has support disjoint from $C_\alpha$, and thus that open subset has a natural
symplectic structure.  In particular, we find that this open subset is
smooth.  It is also minimal, in that it cannot contain any $-1$-curve of a
smooth compactification (e.g., the desired minimal proper regular model).
(More precisely, one can define intersections of line bundles with
projective curves in quasiprojective surfaces; that any quasiprojective
symplectic surface is minimal follows by noting that $0=K\cdot e=-1$ for
any curve $e$ which can be blown down.)  And, of course, it has a natural
fibration over an affine line (induced by that of the complement of
$C_\alpha$ in $X$) such that the smooth locus of the generic fiber is
$\Pic^r$ of the corresponding fiber of $X$.  In other words, this open
subset of the moduli space is precisely the minimal proper regular model of
the relative $\Pic^r$ of $X\setminus C_\alpha$.  It is thus natural to
conjecture that the Zariski closure of this open subset is the full minimal
proper regular model.  We call this Zariski closure the {\em main
  component} of the semistable moduli space (in fact, it is typically the
smallest, but most interesting, component).

It will thus be necessary to understand the remainder of the moduli space.
A key observation is that if $M$ is semistable and supported on $rC_\alpha$,
then $M$ is $S$-equivalent to a sheaf scheme-theoretically supported on
$C_\alpha$.  Indeed, $\omega_X|_{C_\alpha}$ is in the identity component of
$\Pic(C_\alpha)$, and thus twisting by $\omega_X$ preserves semistability.
If for some $l>0$, $M$ is supported on $(l+1)C_\alpha$ but not on $lC_\alpha$,
then we have a nonzero morphism $M\otimes \omega_X\to M$ between semistable
sheaves of the same slope, and thus $M$ is $S$-equivalent to the sum of the
image and the cokernel of this morphism.  This makes $M$ $S$-equivalent to
a sheaf supported on $lC_\alpha$, and we may proceed by induction in $l$.

\begin{lem}\label{lem:antican_stable}
Let $C_\alpha$ be an anticanonical curve on a rational surface $X$.  If $M$
is a semistable sheaf supported on $C_\alpha$ with $c_1(M)=rC_\alpha$,
$\chi(M)=r$, then $M$ is $S$-equivalent to a sum of stable sheaves with
$c_1(M)=C_\alpha$, $\chi(M)=1$.
\end{lem}

\begin{proof}
  We first claim that the slope 0 sheaf ${\cal O}_{C_\alpha}$ is stable.
  We need to show that any quotient sheaf has positive Euler
  characteristic, and can easily reduce to the case of a torsion-free
  quotient, i.e., ${\cal O}_D$ for some curve $D\subset C_\alpha$.  But by
  Lemma \ref{lem:Ca_is_num_conn}, we have $\chi({\cal O}_D)=h^0({\cal
    O}_D)>0$ as required.  It follows immediately that the ideal sheaf of a
  point is stable, since all of its subsheaves are subsheaves of ${\cal
    O}_{C_\alpha}$, so have negative Euler characteristic.  Then, by
  duality on $X$, we obtain stability of the sheaf
\[
{\cal O}_{C_\alpha}(p)
:=
\sExt^1_X({\cal I}_p,\omega_X),
\]
the unique nontrivial extension of ${\cal O}_p$ by ${\cal O}_{C_\alpha}$.

It will thus suffice to show that $M$ is $S$-equivalent to a sum of sheaves
of the form ${\cal O}_{C_\alpha}(p)$.  In fact, it will suffice to construct
a nonzero homomorphism from $M$ to some ${\cal O}_{C_\alpha}(p)$: since
${\cal O}_{C_\alpha}(p)$ is stable of the same slope as $M$ (true regardless
of the choice of ample divisor!), such a morphism is necessarily
surjective.  The kernel of the surjection will then remain semistable of
the same slope, and we may proceed by induction.  The same argument shows
that $\Hom(M,{\cal O}_{C_\alpha})=0$ (since ${\cal O}_{C_\alpha}$ is stable of
smaller slope than $M$), and thus by duality (note that $C_\alpha$ is
Gorenstein with trivial dualizing sheaf), $H^1(M)=0$.

For any point $p\in C_\alpha$, consider the short exact sequence
\[
0\to M'_p\to M\to M\otimes {\cal O}_p\to 0.
\]
If the map $H^0(M)\to H^0(M\otimes {\cal O}_p)$ fails to be surjective,
or, equivalently,
\[
\Ext^1(M\otimes {\cal O}_p,{\cal O}_{C_\alpha})
\to
\Ext^1(M,{\cal O}_{C_\alpha})
\]
fails to be injective, then $\Hom(M'_p,{\cal O}_{C_\alpha})\ne 0$, and any
such morphism induces a nontrivial extension of $M\otimes {\cal O}_p$ by
${\cal O}_{C_\alpha}$ together with a nonzero morphism from $M$ to this
extension.  (Indeed, considerations of Hilbert polynomials show that the
image of this morphism has first Chern class $C_\alpha$.)  Since
$\dim\Ext^1({\cal O}_p,{\cal O}_{C_\alpha})=1$, this extension has the form
${\cal O}_p^n\oplus {\cal O}_{C_\alpha}(p)$, and thus $M$ has a nonzero
morphism to ${\cal O}_{C_\alpha}(p)$.

If the map $H^0(M)\to H^0(M\otimes {\cal O}_p)$ is always surjective, then
the map $H^0(M)\otimes {\cal O}_{C_\alpha}\to M$ is surjective on fibers, and
thus surjective.  But since both sheaves have the same first Chern class and
$M$ has the larger Euler characteristic, this is impossible!
\end{proof}

\begin{rems}
  Note that one can reconstruct $p$ from the sheaf ${\cal O}_{C_\alpha}(p)$,
  since the latter has a unique global section.  It follows that there are
  at most $r$ distinct points $p_i$ admitting morphisms $M\to {\cal
    O}_{C_\alpha}(p_i)$.  Since the cokernel of the natural morphism
  $H^0(M)\otimes {\cal O}_{C_\alpha}\to M$ is supported (set-theoretically)
  on those points, we conclude that the natural morphism is injective, and
  the cokernel is a $0$-dimensional sheaf of degree $r$, from which we can
  read off the $S$-equivalence class of $M$.
\end{rems}

\begin{rems}
  If we replace semistability by the weaker condition that any nonzero
  quotient of $M$ has positive Euler characteristic, we may still conclude
  by the same argument that $M$ has a surjective morphism to some sheaf of
  the form $\sO_{C_\alpha}(p)$.  If we further replace $\chi(M)=r$ by
  $\chi(M)<r$, then $H^0(M)\to H^0(M\otimes \sO_p)$ can never be
  surjective, so that $M$ has a surjective morphism to {\em every} sheaf of
  the form $\sO_{C_\alpha}(p)$.
\end{rems}

\smallskip

We thus conclude that the portion of the moduli space classifying sheaves
supported on the special fiber consists (up to $S$-equivalence) of sums of
sheaves ${\cal O}_{C_\alpha}(p)$, and need to know which of these sheaves
lie on the main component.  The key additional constraint comes from the
observation that if $M$ is supported on $X\setminus C_\alpha$, then
$M\otimes \omega_X\cong M$.  Twisting by $\omega_X$ induces an automorphism
of the full semistable moduli space, and it follows that this automorphism
must act trivially on the main component.  In other words, if $M$ is
$S$-equivalent to
\[
\bigoplus_{1\le i\le r} {\cal O}_{C_\alpha}(p_i),
\]
the multiset of points $p_i$ must be permuted by the action of $\omega_X$.
Since this action is free of order $r$ on the smooth locus, we find that
the $S$-equivalence classes fixed by the automorphism consist of sheaves
\[
\bigoplus_{1\le k\le r} {\cal O}_{C_\alpha}(p)\otimes \omega_X^k
\]
with $p$ in the smooth locus, together with sums
\[
\bigoplus_{1\le i\le r} {\cal O}_{C_\alpha}(p_i)
\]
in which each $p_i$ is a singular point of $C_\alpha$.  In our case, since
$C_\alpha$ is reduced, the latter gives only finitely many points.  The first
family of sheaves is manifestly classified by the smooth locus $C'_\alpha$,
with the closure of $C'_\alpha$ containing in addition only the sheaves
${\cal O}_{C_\alpha}(p)^r$ with $p$ singular.  (There could in principle be
isolated additional points in the main component, but we will see below
that this cannot happen.)

The main difficulty at this point is that it is very difficult to determine
tangent spaces to GIT quotients at semistable points (especially so in our
case, since we only want the tangent vectors coming from a particular
component).  To get around this, we will consider one more moduli space.

The condition that a simple sheaf is invertible on its support is open (we
can express it as the condition that the first Fitting scheme is empty), as
is the condition that it be semistable.  If $C$ is any fiber of the genus 1
fibration on $X$, then $\sO_C$ has a unique global section, and thus any
invertible sheaf on $C$ is simple.  We thus obtain an algebraic space
parametrizing semistable invertible sheaves on fibers of $X$.  As before,
any semistable invertible sheaf not supported on the special fiber
$rC_\alpha$ is stable, and thus away from the special fiber, we recover the
N\'eron model of the relative $\Pic^r$.  This fails on the special fiber
for the simple reason that a given $S$-equivalence class can occur more
than once.  By the above classification of $S$-equivalence classes, we find
that the $S$-equivalence class of $M$ is determined by the orbit under
twisting by $\omega_X$ of the invertible sheaf $M|_{C_\alpha}$; in
particular, $M$ is stable iff $M|_{C_\alpha}\cong {\cal O}_{C_\alpha}(p)$ for
some point $p$ of the smooth locus.  Thus each $S$-equivalence class is
represented by $r$ distinct points of the algebraic space parameterizing
semistable invertible sheaves.

The point is that we can compute tangent spaces in this algebraic space:
\begin{align}
\dim\Ext^1_X(M,M)
&=
\dim\Hom_X(M,M)
+
\dim\Hom_X(M,M\otimes \omega_X)\notag\\
&=
\dim\Hom_{rC_\alpha}(M,M)
+
\dim\Hom_{rC_\alpha}(M,M\otimes \omega_X)\notag\\
&=
\dim\Gamma(\sO_{rC_\alpha})
+
\dim\Gamma(\omega_X|_{rC_\alpha})\notag\\
&=
2,
\end{align}
and thus the ($2$-dimensional) algebraic space is smooth at these points.
Since twisting by $\omega_X$ acts without fixed points on the special
fiber, it preserves tangent spaces, and thus the corresponding subset of
the semistable moduli space is smooth.  In particular, we conclude that the
main component of the semistable moduli space is smooth on the locus
represented by invertible sheaves, i.e., on the smooth locus of $C'_\alpha$.

The minimal desingularization of the main component is thus a proper
regular model of the relative $\Pic^r$, so blows down to the minimal proper
regular model.  It follows that if we simply remove the singular points
from the main component, the result maps to the minimal proper regular
model.  Now, the special fiber of the main component has the same number of
components as the special fiber of the minimal proper regular model (which
must have the same Kodaira type as $C_\alpha$,
\cite[Thm.~5.3.1]{CossecFR/DolgachevIV:1989}).  Thus the only way the surfaces
can fail to be isomorphic is if the map from the minimal desingularization
of the main component blows down a component of the original special fiber.
(Indeed, we must blow down as many components as the minimal
desingularization introduces.)  Now, any section of $\Pic^r(X\setminus
C_\alpha)$ extends to a $-1$-curve on the minimal proper regular model, which
must in particular meet the special fiber in a point of the smooth locus.
It follows that if the corresponding curve in the main component meets the
special fiber in a point of the smooth locus, the corresponding component
cannot be contracted.  We will see that (under the additional saturation
hypothesis) any component is met by some $-1$-curve, giving the desired
isomorphism.

Now, let $e$ be any $-1$-class on $X$ which is transverse to $C_\alpha$, and
consider the corresponding $\tau$-divisor $\tau(e)$.  This certainly
determines a well-defined curve in the complement of the special fiber (and
any non-integral fibers), and we claim that its closure in the main
component meets the special fiber in a single point, which lies in the
smooth locus.  Indeed, if $M$ is supported on the special fiber and
$\Gamma(M(-e))\ne 0$, then we find that $M$ is $S$-equivalent to
$\sO_{rC_\alpha}(e)$.  Indeed, we may state the condition as
$\Hom(\sO_{rC_\alpha}(e),M) = \Hom(\sO_X(e),M) \ne 0$.  Since the image is
both a quotient of the semistable sheaf $\sO_{rC_\alpha}(e)$ and a subsheaf
of the semistable sheaf $M$, the image is also semistable, and can be
extended to Jordan-H\"older filtrations of both $M$ and $\sO_{rC_\alpha}$.
Since $M$ and $M\otimes \omega_X$ are $S$-equivalent, this is enough to
completely determine the $S$-equivalence class of $M$ as required.

In particular, we find that $\tau(e)$ meets the special fiber in the image
of $e\cap C_\alpha$ in $C'_\alpha$. (In particular, we may choose $e$ so that
this point lies in any desired component of $C'_\alpha$.)  It remains only to
show that $\tau(e)$, or rather the corresponding invertible sheaf $\det
R\Gamma(M(-e))$, is a $-1$-class on the minimal proper regular model, and
that this correspondence between $-1$-classes extends to a homomorphism
preserving the intersection pairing.

If $\sO_C(e)$ is stable for every fiber $C$, then $\tau(e)$ consists
precisely of sheaves of that form, and is thus a rational curve as
required.  Since it meets the generic fiber (and thus the anticanonical
curve) in a single point, we conclude that it is a $-1$-curve.  More
generally, adding a component of a nonspecial fiber to $e$ does not change
how $\tau(e)$ meets the special fiber or any integral fiber, and in this
way we can arrange for $\sO_C(e)$ to be stable for all $C$.  In particular,
we find that any class $\det R\Gamma(M(-e))$ obtained in this way is the
sum of the class of a $-1$-curve and a linear combination of components of
nonspecial fibers.

It remains to see that this correspondence extends to a homomorphism and
preserves the intersection pairing.  Both of these are closed conditions on
the (irreducible) moduli stack, so we may impose any dense conditions we
desire.  In particular, we may assume that $C_\alpha$ is smooth and every
nonspecial fiber of $X$ is integral, so that $\tau(e)$ is a $-1$-curve for
every $-1$-class $e$.  Let $e'$ be another $-1$-curve on $X$.  Then
$\tau(e')\cdot\tau(e)$ may be computed as the degree of $\det
R\Gamma(M(-e'))|_{\tau(e)}$, or equivalently as the degree of $\det
R\Gamma(\sO_C(e-e'))$ as $C$ varies over fibers of the genus 1 fibration on
$X$.  Now, consider the natural presentation
\[
0\to \sO_X(e-e'+rK)\to \sO_X(e-e')\to \sO_C(e-e')\to 0
\]
Since $(e-e')\cdot K=0$ and $X$ has no $-2$-curves, $e-e'$ is
ineffective, and similarly for $e-e'+rK$, $e'-e+K$, $e'-e+(1-r)K$.
Thus $R\Gamma(\sO_C(e-e'))$ is represented by the complex
\[
H^1(\sO_X(e-e'+rK))\to H^1(\sO_X(e-e')),
\]
since the other cohomology groups vanish.  This map depends linearly on the
original choice of map $\sO_X(rK)\to \sO_X$, and thus the desired degree
may be computed as the common dimension of the two cohomology groups, which
by Hirzebruch-Riemann-Roch is equal to $-1-(e-e')^2/2=e\cdot e'$ as
required.

Finally, to see that this extends to a homomorphism, we note that the
intersection form on the ten $-1$-curves $s-e_1$, $f-e_1$,
$e_1$,\dots,$e_8$ has determinant $-1$, so the corresponding $\tau$
divisors span $\Pic(X')$, just as the original divisors span $\Pic(X)$.
Since we may use the intersection form to expand any element of $\Pic(X)$
in that basis, we conclude that the $\tau$-divisor map is linear.

We excluded two cases above, in which $C_\alpha$ has Kodaira symbol $I_8$ or
$I_9$.  In the latter case, we can easily see that any Jacobian fibration
with an $I_9$ fiber has Weierstrass form $y^2+txy+a_3y=x^3$ over the
algebraic closure, with $a_3\ne 0$.  Any two such surfaces are
geometrically isomorphic (in a nonunique way), and thus the claim follows
immediately.  Similarly, in the bad $I_8$ case, the corresponding Jacobian
fibration must be the desingularization of the blow-up in the identity of a
surface
\[
y^2+txy = x^3+a_2 x^2+a_4 x
\]
with $a_4\ne 0$; two such surfaces are geometrically isomorphic iff they
have the same value of $a_2^2/a_4$.  There are two $-1$-curves on this
surface that do not meet the singular point, which meet the corresponding
$G_m$ in the points $\lambda$, $1/\lambda$ where
\[
\frac{(\lambda+1)^2}{\lambda} = \frac{a_2^2}{a_4}.
\]
In particular, the surface is determined up to (again nonunique)
isomorphism by the two points of intersection, and again the claim follows.
\end{proof}

\begin{rems}
  This is extended to general $\Pic^d$ in Corollary \ref{cor:Picr} below,
  which also removes the condition that the special fiber be reduced,
  though the description of the resulting surface is somewhat less
  explicit.  Note that even though the result is purely commutative, the
  argument involves working with noncommutative surfaces, and is based on
  showing that the obstruction to having a universal family over the
  relative $\Pic^d$ is given by a noncommutative rational surface with $K^2=0$.
\end{rems}

\begin{rems}
Note that in the good cases, we prove a slightly stronger fact: not only is
the special fiber of $\Pic^r(X)$ isomorphic to $C'_\alpha$, but the the
isomorphism we construct is compatible with the isomorphism $X'\cong
\Pic^r(X)$.  This still holds in the $I_8$ and $I_9$ cases, but
the above calculation does not suffice.
\end{rems}

\medskip

In \cite{EOR}, three families of algebras were considered based on
deformations of multiplicative Deligne-Simpson problems.  When the
noncommutative parameter $q$ is a root of unity, these algebras have large
center, which was conjectured to be a certain explicit del Pezzo surface.
It turns out that this conjecture can be reduced to a special case of
Theorem \ref{thm:Picr}.  (Given our comment above about how the proof of
Corollary \ref{cor:Picr} works, this suggests that the algebras of
\cite{EOR} are closely related to the noncommutative rational surfaces we
study below; this is true, but surprisingly difficult to show, mainly
because the construction of \cite{EOR} involves explicit presentations,
while the construction below is far more abstract.)

Recall that we mentioned above that multiplicative Deligne-Simpson problems
with three or four matrices could be translated to moduli problems about
sheaves on rational surfaces.  Indeed, given matrices $g_1$, $g_2$, $g_3$
with $g_1g_2g_3=1$, we may define a morphism
\[
B:\sO_{\P^2}(-2)^n\to \sO_{\P^2}(-1)^n
\]
by $B=x+g_1y+g_1g_2z$, which as usual can be recovered from its cokernel.
Moreover, we can recover the conjugacy class of $g_1$ from the restriction
$\coker(B)|_z$, and similarly for $g_2$ and $g_3$.  (The support of
$\coker(B)|_z$ is the set of eigenvalues of $g_1$, and the Jordan block
structure can be read off from the lengths of the direct summands supported
at each eigenvalue.)

In the first case of interest, $g_1$, $g_2$, $g_3$ are $3\times 3$ matrices
with respective characteristic polynomials
$(t-x_{i1})(t-x_{i2})(t-x_{i3})$. If $\prod_{ij} x_{ij}\ne 1$, the
corresponding Deligne-Simpson problem has no solution (since
$\det(g_1g_2g_3)\ne 1$), but if $\prod_{ij}x_{ij}=\zeta_n$, then we can
consider the case of $3n\times 3n$ matrices with the same characteristic
polynomials.  The corresponding (affine) moduli space is then precisely the
desired center (of the $E_6$ algebra of \cite{EOR}).  Resolving the
singularities in the usual way blows up 9 points and thus gives us a
elliptic surface, and we immediately see that $\sO_{C_\alpha}(C_\alpha)$
has order $n$.  It thus follows that the desired moduli space is nothing
other than the semistable version of the relative $\Pic^0$ of this elliptic
surface.  This is not {\em quite} the problem we considered above, but of
course we could fix this by twisting by a suitable line bundle.  In the
generic case, $mC_\alpha+e_8$ is ample for $m\gg 0$, so if we take it as
the ample bundle determining the stability condition, then twisting by it
has no effect on semistability, and has the same effect as twisting by
$e_8$, so changes the degree by $r$.  In general, this bundle is only nef,
but we can make it ample by a suitable deformation, at the cost of
resolving some singularities of the moduli space (corresponding to
semistable sheaves becoming stable).  We thus see that the moduli space is
given by quotienting by the automorphism corresponding to
$\sO_{C_\alpha}(C_\alpha)$, which is nothing other than the automorphism
multiplying every parameter by $\zeta_n$, so that the quotient map raises
all of the parameters to the $n$-th power.

Now, for the sheaves in this open subset to have presentations involving
trivial bundles, we want $H^0(M)=H^1(M)=0$, and thus the moduli space is
the complement of a fiber and section on the relative Jacobian.  The fiber
corresponds to the original $n$-tuple fiber, and has the same Kodaira type
(since relative Jacobians preserve Kodaira types of tame multiple fibers,
\cite[Thm.~5.3.1]{CossecFR/DolgachevIV:1989}); since the anticanonical
divisor on that surface was integral (nodal), we see that the fiber being
removed is an integral nodal curve.  Moreover, we could blow down the
section before removing it, and in this way obtain a del Pezzo surface of
degree 1 with a nodal fiber of type $I_3$ removed.  The problem of
identifying this del Pezzo surface reduces to the problem of identifying
the corresponding elliptic surface, and thus to a special case of Theorem
\ref{thm:Picr}.  For $n=1$, this agrees with the del Pezzo surface for
which explicit equations were given in \cite{EOR}; for $n>1$, the
conclusion of Theorem \ref{thm:Picr} settles the conjecture made there (to
wit that the formula for the equation of the moduli space need only be
modified by taking $n$-th powers of the input).

For the $E_7$ case, we now have 2, 4, and 4 eigenvalues on each leg of the
triangle, and thus most naturally obtain sheaves on $X_9$.  However,
the sheaves of interest have Chern class
\[
n(4h-2e_0-2e_1 -e_2-e_3-e_4-e_5 -e_6-e_7-e_8-e_9,
\]
where the anticanonical curve has decomposition
\[
-K_X=(h-e_0-e_1)+(h-e_2-e_3-e_4-e_5)+(h-e_6-e_7-e_8-e_9),
\]
and this divisor class is not in the fundamental chamber.  Reflecting in
$h-e_0-e_1-e_2$ gives
\[
n(3h-e_0-e_1-e_3-e_4-e_5-e_6-e_7-e_8-e_9)
\]
with decomposition
\[
-K_X = (e_2)+(h-e_2-e_3-e_4-e_5)+(2h-e_0-e_1-e_2-e_6-e_7-e_8-e_9),
\]
and thus we can blow down $e_2$ (after permuting the blowups) to again obtain a
moduli problem on an elliptic surface.  (Again, the affine surface we
obtain is most naturally del Pezzo of degree 1, but we can blow down a
$-1$-curve at infinity to obtain a del Pezzo surface of degree 2.)

For the $E_8$ case, the problem originates on $X_{10}$ with Chern class
\[
n(6h-3e_0-3e_1-2e_2-2e_3-2e_4-e_5-e_6-e_7-e_8-e_9-e_{10})
\]
and anticanonical decomposition
\[
-K_X = (h-e_0-e_1)+(h-e_2-e_3-e_4)+(h-e_5-e_6-e_7-e_8-e_9-e_{10}).
\]
Again, this is not in the fundamental chamber, and indeed the component
$h-e_0-e_1$ of $C_\alpha$ is a $-1$ curve disjoint from our sheaf, so may
be blown down.  Doing so makes the image of $h-e_2-e_3-e_4$ itself a
$-1$-curve, so that we may blow {\em that} component down to again reach an
elliptic surface.

There is also a $D_4$ case considered in \cite{OblomkovA:2004}.  In that case,
there are four matrices, so one starts with $\P^1\times \P^1$, and obtains
a moduli problem with Chern class
\[
n(2s+2f-e_1-e_2-e_3-e_4-e_5-e_6-e_7-e_8)
\]
and anticanonical decomposition
\[
-K_X=(s-e_1-e_2)+(f-e_3-e_4)+(s-e_5-e_6)+(f-e_7-e_8),
\]
so again an elliptic surface.  The only difference from the $E_6$ case is
that we need to impose {\em two} vanishing conditions, so are removing two
$-1$-curves in addition to the anticanonical curve of the moduli space, and
thus naturally obtain a del Pezzo surface of degree 2, which can be further
blown down to an affine cubic surface with three lines at infinity as
specified in \cite{OblomkovA:2004}, or even to an affine del Pezzo surface of
degree 4.

There is one technical issue we glossed over above: the above moduli
problems are only the {\em expected} behavior of solutions, but in general
one should actually allow the multiplicities of the eigenvalues to differ.
(This is because the actual problem only specifies the minimal polynomials
of $g_1$, $g_2$, $g_3$.)  Thus, for instance, in the $E_6$ case, the actual
constraint on the Chern class is only that $D\cdot h=3n$ and that $D$ is
disjoint from the three anticanonical components.  But it is easy to see
that any other choice will have negative self-intersection, and thus its
support must contain a $-2$-curve.  Removing such $-2$-curves either leaves
nothing (making the solution rigid) or a sheaf of self-intersection 0,
which must be a sheaf of the above form (so that the $-2$-curves don't
affect the solution to the Deligne-Simpson problem).

\medskip

In \cite{Crawley-BoeveyW/ShawP:2006}, multiplicative Deligne-Simpson
problems were related to Coxeter groups, in this case to groups with
arbitrary star-shaped Dynkin diagrams.  (In particular, $E_{m+1}$ has a
star-shaped diagram, and the corresponding Deligne-Simpson problems have a
quadratic and a cubic minimal polynomial.)  At least in the three- and
four-matrix cases, we can see these Coxeter groups from the rational
surface perspective.  In the three-matrix case, these are precisely the
reflection subgroups of the stabilizer of the decomposition of $C_\alpha$,
see \cite{LooijengaE:1981}; the simple roots of the subsystem are
\[
h-e_{11}-e_{21}-e_{31},\qquad\text{and}\qquad e_{ij}-e_{i(j+1)},
\]
where the first subscript on the $e$'s indicates which component of
$C_\alpha$ they intersect.  In the four-matrix case, the description is
slightly more subtle: a four-matrix Deligne-Simpson problem corresponds to
a sheaf on $\P^1\times \P^1$, and swapping the rulings changes the
Deligne-Simpson problem; thus in addition to stabilizing the decomposition
of $C_\alpha$, we also want to stabilize the root $s-f$.  In addition to
the obvious $A_l$-subsystems, we have a simple root
$s+f-e_{11}-e_{21}-e_{31}-e_{41}$.  In any event, even if we assume that
these reflections generate the full stabilizer (rather than just the
subgroup generated by reflections), the algorithms still become more
complicated.  Indeed, we have already seen that there are several kinds of
elliptic pencil on such a surface (even if we exclude multiple fibers),
distinguished by the components that get blown down on the relatively
minimal model; as a result, it is more difficult to identify whether a
class is integral based merely on its image in the fundamental chamber of
this smaller Coxeter group.

Of course, as we remarked above, the stabilizer of the anticanonical
decomposition can fail to be a reflection group (and that reflection group
can apparently fail to be of finite rank).  There are a few cases
\cite{LooijengaE:1981} in which the stabilizer is a reflection group and
has been explicitly identified, specifically those with nodal (and thus
polygonal) anticanonical curve, having at most $5$ components.  (Looijenga
also remarks that the $6$ component case is probably feasible; and we have
already seen that the stabilizer can fail to be a reflection group when
there are $7$ or more components.) The $2$ component case has a star-shaped
diagram with one very short leg, corresponding to the relation between
nonsymmetric $q$-difference equations and three-matrix multiplicative
Deligne-Simpson problems in which one of the minimal polynomials is
quadratic (discussed in Section \ref{sec:degen} below).

\chapter{Moduli of equations}
\label{chap:moduli_of_eqs_comm}

\section{Symmetric elliptic difference equations}
\label{sec:elldiff}

We now wish to translate the above theory back to the realm of difference
equations, and in particular see what it suggests should be true for the
non-relaxed case.  From a geometric perspective, the simplest case is that
of symmetric elliptic difference equations, since then not only is the
surface smooth, but so is the anticanonical curve.  If $C$ is a smooth
genus 1 curve, then the above considerations tell us that matrix-form
symmetric difference equations on $C$ twisted by a line bundle are in
natural correspondence with triples $(X_0,\phi,M_0)$ where $X_0$ is a
Hirzebruch surface (in particular with specified map to $\P^1$), $M_0$ is a
sheaf on $X_0$ with $H^0(M_0)=H^1(M_0)=0$, and $\phi:C\to X_0$ embeds $C$
as an anticanonical curve.  (As we mentioned above, this is not quite
correct, as these sheaves also include degenerate cases with apparent
singularities.)  Moreover, the pairs $(X_0,\phi)$ are classified by
elements of $\Pic^2(C)\times \Pic^2(C)$ or $\Pic^1(C)\times \Pic^2(C)$,
depending on the parity of the Hirzebruch surface, or equivalently
depending on the parity of the degree of the twisting line bundle.

Since $C\cong C_\alpha$ is smooth, $M_0|_{C_\alpha}$ is a direct sum of
structure sheaves of jets.  If $X$ is the minimal desingularization of the
blowup of $X_0$ in those jets, then Proposition
\ref{prop:jets_resolvable_by_!*} tells us that the minimal lift $M$ of
$M_0$ to $X$ is disjoint from $C_\alpha$.  This encodes the {\em locations}
of the singularities of the equation in the surface $X$, while the
multiplicities of the singularities are encoded in the first Chern class of
$M$.

More precisely, the surface is determined by the classes $\phi^*(s)$,
$\phi^*(f)$ and $\phi^*(e_i)$ for $1\le i\le m$; since $C$ is smooth, the
bundles $\phi^*(e_i)$ can be identified with line bundles $\sO_C(p_i)$ for
points $p_i$.  Now, consider the extension
\[
B:\rho^*V\to \rho^*W(s)
\]
of our original $B$ to $X_0$.  The Chern class of $M_0$ is given by that of
$\det(B)$, so has the form $ns+df$ where $n=\rank(W)$ is the order of the
difference equation.  At any point $p\in C_\alpha$, we can view $B$ as a
matrix over the local ring $\sO_{C_\alpha,p}$.  Up to left- and
right-multiplication by invertible matrices, we can diagonalize $B$, and
then define a partition $\lambda(B;p)$ by letting $\lambda_j(B;p)$ for
$j\ge 1$ be the number of diagonal elements contained in $\mathfrak{m}^j$.
(This partition can be computed in a more basis-independent way: one can
readily recover it from the valuations of the $\gcd$s of the coefficients
of the exterior powers $\wedge^k B$.)  Local computations then tell us that
$\lambda_1(B;p)$ is the rank of $M_0$ at $p$, and $\lambda_j(B;p)$ is the
rank after blowing up $p$ $j-1$ times.  If $e_{p,j}$ denotes the $j$-th
class in the sequence $e_1$,\dots,$e_m$ such that $\phi^*(e_i)\cong
\sO_C(p)$, then we find
\[
c_1(M) = c_1(M_0) - \sum_{p,i} \lambda_i(B;p) e_{p,i}.
\]
(Note that we must blow up $p$ at least as many times as there are parts of
$\lambda$ in order to make the resulting sheaf disjoint from the
anticanonical curve.)

We can also describe the invariants $\lambda(B;p)$ in terms of the original
shift matrix $A$.  If $p$ is not fixed by $\eta$, then we can again
diagonalize $A$ over the local ring at $p$ (by left- and right-
multiplication).  The resulting equivalence classes are given by weights of
$\GL_n$, i.e., nonincreasing sequences of integers.  We then find that
$\lambda(B;p)$ is determined by the positive coefficients of this weight;
the negative coefficients of the weight appear in $\lambda(B;\eta(p))$.
When $p$ is fixed by $\eta$, the situation is more complicated (as we saw
in our discussion of singularities); up to the relevant equivalence
relation (left multiplication by invertible matrices over the local ring,
right multiplication by symmetric invertible matrices over the local
field), $B$ is a direct sum of matrices
\[
\begin{pmatrix}1\end{pmatrix},\qquad
\begin{pmatrix}u\end{pmatrix},\qquad
\begin{pmatrix}1 & u\\0 & u^e\end{pmatrix}, e>1.
\]
(This needs to be adjusted slightly when the equation is twisted by a line
bundle.)  The second case is a singularity of order 1, but corresponds to
an eigenvalue $-1$ of $A$ (assuming the characteristic is not 2); the third
cases have order $e$, but appear to have order $e-1$.  (In characteristic
2, this phenomenon is worse: singularities of order $e$ appear to have
order $\max(e-2,0)$ or $\max(e-4,0)$ depending on whether $C$ is ordinary
or supersingular.)

Of course, a sheaf $M$ on $X$ disjoint from $C_\alpha$ need not come from a
maximal morphism $B$.  The condition
$\Hom(M,\sO_g(-1))=\Hom(\sO_g(-1),M)=0$ for every component $g$ of a fiber
is simplified by disjointness, since we need only consider those $g$ which
are disjoint from $C_\alpha$.  In particular, $g$ must be a $-2$-curve, and
since it is contained in a fiber, must be a root of the $D_m$ subsystem.
The subquotient corresponding to the standard representation of a
difference equation will differ from $M$ by a number of copies of sheaves
$\sO_g(-1)$, and thus in particular correspond to strictly semistable
points of the moduli space.  As a result, the specific extension classes
will be irrelevant, and we thus obtain precisely one point of the moduli
space for each difference equation that arises.  (This, of course, is
special to the case $\chi(M)=0$; in general, the extension matters, and
corresponds to a choice of filtration of the fiber at $p$ of the ambient
vector bundle, coming from Proposition \ref{prop:sing_filtration} with the
obvious ordering on exceptional components at $p$.  Without vertical
components, this filtration agrees with the obvious one coming from the
diagonalization of $B$.)

The constraint on the difference equations underlying such semistable
points is that their Chern class must differ from the specified Chern class
by a nonnegative sum of $-2$-curves disjoint from $C_\alpha$.  This can be
translated in terms of the local data $\lambda(B;p)$ as follows.
Subtracting roots of the form $e_i-e_j$ simply replaces the partitions
$\lambda(B;p)$ by partitions of the same size which cover it in the
dominance ordering.  Roots of the form $f-e_i-e_j$ either subtract 1 from
the first parts of both $\lambda(B;p)$ and $\lambda(B;\eta(p))$ or (when
$p=\eta(p)$) subtract $1$ from the first two parts of $\lambda(B;p)$.  For
$p\ne \eta(p)$, the conditions combine to say that the relevant weight of
$\GL_n$ (related to the conjugate partitions) becomes smaller in dominance
order, with something similar in the case of ramification points.

One special case we should note is that when $D=s-f$ (assuming this is
effective), then as usual, the moduli space is a point (the sheaf
$\sO_{s-f}(-1)$), and the corresponding difference equation is just the
trivial equation $v(z+q)=v(z)$.  More generally, if $s-f$ is effective and
$D\cdot (s-f)<0$, then the corresponding difference equation will have a
block-triangular structure such that the first or last block is trivial.
(More precisely, one has such a block-triangular structure only up to
isomonodromy transformation; indeed, the theory of semiclassical orthogonal
polynomials \cite{MagnusAP:1995} or semiclassical biorthogonal functions
\cite{isomonodromy} involve precisely such equations, with the triangular
structure becoming less and less apparent as the degree of the polynomial
grows.)

\section{Isomonodromy}

Of course, a sheaf on an anticanonical surface $X$ with $C_\alpha\cong C$
does not determine a difference equation unless we also choose a blowdown
structure (more precisely, a choice of blowdown structure modulo the action
of ineffective roots of the $A_m$ subsystem).  In other words, a given
moduli space of sheaves corresponds to many different moduli spaces of
difference equations, one for each blowdown structure on the surface.  In
addition, we can also twist by line bundles and apply the duality
$\sExt^1(-,\omega_X)$.  The latter is canonical, but to make sense of the
former requires a choice of blowdown structure.  Thus in the generic
situation, we have an action of $\Aut(C)\times (W(E_{m+1})\times
Z_2)\ltimes \Z^{m+2}$, where the cyclic group $Z_2$ acts by duality; in the
nongeneric situation, there is a partial action taking into account the
usual issues with effective reflections.  Since the group $W(E_{m+1})$
simply acts on the set of ways of interpreting sheaves, it certainly
respects the Poisson structure on the moduli space, and the construction of
the Poisson structure implies that it is preserved by twisting.  Duality is
anti-Poisson, and $\Aut(C)$ acts on the Poisson structure in the same way
it acts on holomorphic differentials.  In particular, $\Pic^0(C)\subset
\Aut(C)$ preserves the Poisson structure, and hyperelliptic involutions are
anti-Poisson; in the $j=0$ and $j=1728$ cases, we also have automorphisms
multiplying the Poisson structure by other roots of unity.  In terms of
moduli spaces of unrelaxed difference equations, each of these operations
will change the parameters (the twisting line bundle, the points with
allowed singularities, $q$), but should give birational maps between the
corresponding moduli spaces.  The action on $q$ is essentially forced:
Poisson maps should preserve $q$, while anti-Poisson maps should negate $q$
(and for $j=0$, $j=1728$, $\Aut(C)$ acts as one would expect on $q\in
\Pic^0(C)$).

Note that the subgroup $D^\perp\subset \Z^{m+2}$ acts as a (large) abelian
group of rational Poisson automorphisms of the moduli spaces $\Irr(D,x)$.
(The element $K_X$ acts trivially, of course, as does any effective class
in $D^\perp$.) In particular, they give rise to a discrete integrable
system acting on a rational variety, which relative to the fibration by
$\supp(M)$ acts by translation within each fiber (a torsor over the
Jacobian of the support).  Similarly, we will describe a $q$-twisted
version of this action, which an analogue of (higher-order) discrete
Painlev\'e equations, a non-autonomous version of translation on an abelian
variety.

Of course, in our setting, we only have a relaxation of the true moduli
spaces of difference equations, but will still gain insight by looking at
how simple reflections, twists by basis elements, and duality act; in each
case, there will be an obvious way to take into account the shifting by
$q$.  Since the simplest operations do not preserve the triviality
condition $H^0(M)=H^1(M)=0$ even generically (since they do not preserve
the condition $\chi(M)=0$), we need to allow $W$ to be nontrivial.  We thus
note (per \cite{generic}, or from the discussion in Section
\ref{sec:sheaves_from_eq_nc}) that in the (analytic) difference equation
case, $B$ corresponds to a morphism
\[
B:\pi_{\eta'}^*V\to \pi_{\eta}^*W\otimes {\cal L}
\]
of bundles on $\C/\Lambda$, where $\eta$ is the involution $z\mapsto -q-z$,
$\eta'$ is the involution $z\mapsto -z$, and ${\cal L}$ is the twisting
line bundle.  We do, however, assume $V$ maximal where convenient, since
this is in any event the main case of interest (and it is easy enough to
figure out what goes wrong when maximality fails).

The simplest operation is twisting by $\sO_X(f)$, which simply twists
the bundles $V$ and $W$ by $\sO_{\P^1}(1)$.  This is also easy to extend to
an action on difference equations: the only change is that since $V$ and
$W$ are pulled back through different degree $2$ maps, we must absorb the
difference into ${\cal L}$, thus changing the twisting bundle by the
element of $\Pic^0(C)$ corresponding to $q$.  (In other words, twisting by
$f$ changes $\phi^*(s)$ by $q$.)

\medskip

Although the operation $M\mapsto M(s)$, or equivalently
$M_0\mapsto M_0(s)$, is just as natural in terms of
sheaves, the translation to morphisms of vector bundles on $C$ is quite a
bit more subtle, for the simple reason that twisting does not respect the
resolution we are using.  Now, we can write the original morphism
$B:\pi_\eta^*V\to \pi_\eta^*W\otimes \phi^*\sO_X(s)$ in the form
\[
B = B_0(x,w) y_0 + B_1(x,w) y_1
\]
where $x,w$ are homogeneous coordinates on
$\P^1\cong\pi_\eta(C)\cong\rho(X)$, and $y_0$, $y_1$ generate
$\pi_{\eta*}\phi^*\sO_X(s)\cong \rho_*\sO_X(s)$, viewed as a graded module
over $k[x,w]$.  Here, we write $V$ and $W$ as sums of line bundles, and
the degree of a given coefficient of $B_0$ and $B_1$ depends on the
corresponding difference of degrees of line bundles.  (In the untwisted
case, we used $\phi^*\sO_X(s_{\min})$, but this just differs by twisting by
$f$; i.e., we assume that the twisting bundle has degree 2 when its degree
is even.) Now, since $B$ comes from the standard resolution of $M_0$, we
can twist by $\sO_X(s)$ and take direct images to obtain a short exact
sequence
\[
0\to V\xrightarrow{(B_0,B_1)} W\otimes
\rho_*\sO_X(s)\xrightarrow{(B'_1,-B'_0)} \rho_*(M_0(s))\to 0
\]
where $B'_1$, $B'_0$ are suitable morphisms of sheaves on $\P^1$.
The sheaf $\rho_*(M_0(s))$ is torsion-free, since 
otherwise $\Hom(\sO_f(-1),M_0)$ would be nonzero for some fiber $f$,
contradicting maximality of $V$.  In particular, we have
\[
W' \cong \rho_*(M_0(s)).
\]
We also have a commutative diagram
\[
\begin{CD}
0@>>> \rho^*V @>>> \rho^*W\otimes \rho^*\rho_*\sO_X(s)@>>>
\rho^*\rho_*(M_0(s)) @>>>0\\
@. @| @VVV @VVV @.\\
0@>>> \rho^*V @>>> \rho^*W(s)@>>>
M_0(s) @>>>0
\end{CD}
\]
with exact rows; since each sheaf in the bottom sequence is $\rho$-globally
generated, the vertical morphisms are surjective, and have isomorphic
kernels.  We thus find that $V'$ fits into an exact sequence
\[
0\to \rho^*(V')(-s)\to \rho^*W\otimes \rho^*\rho_*\sO_X(s)\to
\rho^* W(s)\to 0.
\]
It follows that
\[
V'\cong \begin{cases} W& s^2=0\\
                      W(-f) & s^2=-1
\end{cases}
\]
We moreover find that (apart from this twisting), the map from $V'$ to $W'$
is simply given by $B'_1y_1+B'_0y_0$.  To avoid the issue with twisting, we
will compute $M(s+f)$ in the odd case.

We now observe that, viewing this as a morphism $B'$ on $C$, we have
\[
B' \eta^*B - \eta^*(B'\eta^*B)
= 
(B'_0B_1-B'_1B_0)(y_0\eta^*y_1-y_1\eta^*y_0)
=
0,
\]
since $B'_0B_1=B'_1B_0$ by construction.  Since $B'\eta^*B$ is
$\eta^*$-invariant, we can use it to modify the factorization of $A$ to
obtain
\[
A = (\eta^*B)^{-t}B^t = (B')^t (\eta^*B')^{-t},
\]
and find that the new $A$ has the form
\[
A' = (\eta^*B')^{-t} (B')^t = (B')^{-t} A (B')^t.
\]
It is somewhat more natural to express the inverse of this operation.

\begin{prop}
  Let $A(z)$ be a twisted elliptic matrix with $\eta^*A = A^{-1}$, twisted
  by a line bundle of degree $\delta$, and let $M$ be the corresponding
  sheaf.  Then the twisted sheaf $M(-s-(2-\delta) f)$ corresponds to the
  matrix $B^t A B^{-t}=B^t \eta^*B^{-t}$, where $B$ comes from the minimal
  factorization of $A$.
\end{prop}

This operation need only be modified very slightly (since conjugation
should become a gauge transformation) to make sense for difference
equations: if we start with the system
\[
v(z+q) = B(-q-z)^{-t} B(z)^t v(z)\qquad v(-z)=v(z),
\]
which we can write in the form
\[
v(-z) = v(z)\qquad B(-q-z)^t v(-q-z) = B(z)^t v(z),
\]
we simply want the equations satisfied by $w(z)=B(z)^t v(z)$, namely
\[
w(-q-z)=w(z)\qquad B(-z)^{-t}w(-z)=B(z)^{-t} w(z).
\]
The only nonobvious point is that this new equation is symmetric with
respect to a slightly different involution; we will see this phenomenon
naturally arising in the noncommutative setting.  (Here twisting by $s$
changes $\phi^*(f)$ by $q$; it also changes $\phi^*(s)$ when $s^2=-1$.  In
general, the rule is that twisting by $D$ changes $\phi^*$ by $(D\cdot
{-})q$; this is the only $W(E_{m+1})$-invariant rule compatible with what
we have so far seen.)  This can also be sidestepped by choosing a point
``$q/2$'' such that $2(q/2)=q$, and replacing the above $w$ by $w(z) =
B(z-q/2)^t v(z-q/2)$.

\medskip

The next operation we consider is duality, as this can also be computed on
the Hirzebruch surface.  Indeed, we have seen that the minimal lift
operation commutes with the canonical duality, so we just need to
understand the dual of $M_0$.  Applying $\dR\sHom(-,\omega_X)$ to the
standard presentation
\[
0\to \rho^*V(-s)\xrightarrow{B} \rho^*W\to M_0\to 0
\]
gives
\[
0\to \sHom(\rho^*W,\omega_X)\xrightarrow{B^t}
\sHom(\rho^*V,\omega_X)(s)\to \sExt^1(M_0,\omega_X)\to 0.
\]
Now, we have
\[
\sHom(\rho^*W,\omega_X) \cong \rho^*\sHom(W,\sO_{\P^1})\otimes \omega_X
\]
and $\omega_X\cong \sO_X(-2s-(4-\delta)f)$, where $\delta\in \{1,2\}$ is
the degree of the twisting bundle.  We thus see that this is a presentation
of the alternate kind we just considered, and can thus determine the
corresponding relaxed difference equation.  We thus find that
$\sExt^1(M,\omega_X)(2f)$ corresponds to the matrix
\[
A' = B (\eta^*B)^{-1} = A^t = \eta^* A^{-t}.
\]
In other words, dual sheaves correspond (up to the involution) to adjoint
difference equations.  More precisely,
if we both dualize and act by $\eta$, we obtain the adjoint equation
\[
w(z+q) = A(z)^{-t} w(z);
\]
the dual sheaf itself corresponds to
\[
w(z-q) = A(z-q)^t w(z),
\]
which is of course precisely the same equation viewed as a $-q$-difference
equation.  (The negation of $q$ comes from the fact that duality negates
the Poisson structure.)  Note that the additional $2f$ twist is precisely
what we need in order for the $\chi(M(-f))=0$ condition to be preserved by
duality.

We next turn to twists by $e_i$.  From Corollary
\ref{cor:pseudo_twist_is_isospectral}, we find that the action of such
twists on $M_0$ has the following form.  If $M'\cong M(-e_i)$, then $M'$ is
acyclic for $\pi:X\to X_0$, and its direct image $M'_0$ fits into a short
exact sequence
\[
0\to M'_0\to M_0\to \sO_p^r\to 0
\]
where $p$ is the point of $X_0$ lying under $e_i$, and the morphism $M_0\to
\sO_p^r$ is suitably canonical (with $r = c_1(M)\cdot e_i$).  Local
computations let us describe this morphism in the elliptic case.  First, if
$e_i$ arises from the first time we blow up $p$, then it is just the
canonical morphism
\[
M_0\to \Hom_k(\Hom(M_0,\sO_p),\sO_p).
\]
More generally, the structure of $B$ over the local ring induces a natural
increasing filtration $F_l$ of $\rho^*W\otimes\sO_p^n$, induced by tensor
product from the filtration
\[
F^+_l = \im(u^{1-l}B)\cap \rho^*W
\]
where $u$ is a uniformizer.  If $e_i$ arises from the $l$-th time we blow up
$p$, then the corresponding morphism $M_0\to \sO_p^r$ is induced by the
morphism
\[
\rho^*W\to \rho^*W\otimes\sO_p\cong F_{\infty}\to F_{\infty}/F_l.
\]

Once we have identified the map $M_0\to \sO_p^r$, we then need to
understand how the new $A'$ is related to the original $A$.  We first note
that $M'_0$ is indeed $\rho_*$-acyclic.  Otherwise, there is a nonzero
morphism from $M'_0$ to $\sO_f(-2)$ for some fiber $f$, so that the short
exact sequence defining $M'_0$ pushes forward to an extension
\[
0\to \sO_f(-2)\to F\to \sO_p^r\to 0.
\]
Since $\Hom(M_0,\sO_f(-2))=0$, this is a non-split extension, so $p\in f$;
but then $F\cong \sO_f(-1)\oplus \sO_p^{r-1}$, contradicting the fact that
$\Hom(M_0,\sO_f(-1))=0$ (since $M_0(-s)$ is
$\rho_*$-acyclic).  We also find that $W'=\rho_*M'_0$ is torsion-free, since
$\Hom(\sO_f,M'_0)\subset \Hom(\sO_f,M_0)=0$ for any fiber $f$.

Since $M'_0$, is $\rho_*$-acyclic, we can take the direct image of the
defining extension to obtain a short exact sequence
\[
0\to W'\xrightarrow{D_1} W\to \sO_{\rho(p)}^r\to 0,
\]
for a suitable morphism $D_1$, and since $\rho$ is flat the inverse image is
also exact:
\[
0\to \rho^*W'\to \rho^*W\to \sO_f^r\to 0,
\]
where $f$ is the fiber containing $p$.  We thus obtain a map of short exact
sequences
\[
\begin{CD}
0@>>> \rho^*V(-s) @>B>> \rho^*W @>>> M_0@>>> 0\\
@. @VVV @VVV @VVV @.\\
0@>>> \sO_f(-1)^r @>>> \sO_f^r @>>> \sO_p^r @>>> 0,
\end{CD}
\]
in which the second and third vertical maps are surjective.  The snake
lemma tells us that the natural four-term exact sequence
\[
0\to \rho^*V'(-s)\to \rho^*W'\to M'_0\to \sO_f(-1)^{r'}\to 0
\]
(the cokernel has the form $\sO_f(-1)^{r'}$ since this is the only extent
to which $M'_0$ can fail to be $\rho$-globally generated) is related to a
four-term exact sequence
\[
0\to \rho^*V'(-s)\to \rho^*V(-s)\to \sO_f(-1)^r\to \sO_f(-1)^{r'}\to 0.
\]
This is
the twist of the inverse image of an exact sequence
\[
0\to V'\xrightarrow{D_2} V\to \sO_p^r\to \sO_p^{r'}\to 0.
\]
Of course, if $r'\ne 0$, then $M'_0$ is not globally generated, so we
should really replace $M'_0$ by the generated subsheaf; in this case, the
twisting operation will not be invertible, but a finite amount of such
twisting will suffice to remove any components of $M'_0$ supported on
fibers.  In any event, the new $B$ can be written as
\[
B' = D_1^{-1} B D_2,
\]
and thus, since $D_1$ and $D_2$ are $\eta$-invariant,
\[
A' = D_1^t A D_1^{-t}.
\]
Again, this conjugation should become a gauge transformation: if
$w(z) = D_1(z)^t v(z)$, then $w$ satisfies the equation
\[
w(z+q) = D_1(z+q)^t A(z) D_1(z)^{-t} w(z).
\]
This gauge transformation has the effect of shifting the singularity at $p$
by $q$ (as we expect from how twisting should affect $\phi^*$); in terms of
the invariants $\lambda(B;p)$, it moves the appropriate part to the
partition corresponding to $p-q$.  The case $r'\ne 0$ corresponds to a
situation in which the shifted singularity ends up cancelling an existing
singularity at $p-q$.  (In particular, we recover our earlier observation
that the sheaves with components supported on fibers correspond to
equations with apparent singularities, that is to say singularities which
can be removed by a suitable gauge transformation.)

\begin{rem}
  A similar construction (unfortunately called ``elementary
  transformations'') for sheaves on $\P^2$ was given in
  \cite{VinnikovV:1990}.
\end{rem}

We note that since $M\otimes \sO_{C_\alpha}=0$, twisting by
\[
\sO_X(C_\alpha)\cong \sO_X(2s+2f-\sum_i e_i)
\]
has no effect on the sheaf.  Since for difference equations, the various
twisting operations all change various parameters by multiples of $q$, this
cannot quite be true for difference equations; instead, twisting by the
canonical class simply shifts $z$ by $q$.  (This is the gauge
transformation by the matrix $A$ itself.)

\begin{rem}
  With the above constructions in mind, we can also easily identify the
  various gauge transformations of \cite{isomonodromy} (called isomonodromy
  transformations there) with twists; in particular, the gauge
  transformations with matrices described by \cite[Thm.~4.6]{isomonodromy}
  are twists by $s+f-\sum_{1\le i\le m+3} e_i$.
\end{rem}

\medskip

It remains to understand how the group $W(E_{m+1})$ acts.  The $S_m$
subgroup is of course easy to understand, as it simply changes the order in
which we blow up the distinct singular points.  Indeed, since it does not
change the final Hirzebruch surface, we should not expect it to have any
effect on the interpretation of the sheaf.

To understand the subgroup $W(D_m)$, it will be enough to understand how
the elementary transformation acts, as it conjugates the two $S_m$
subgroups in $W(D_m)$.  We suppose now that $X_0$ is even, since of course
the odd to even elementary transformation is just the inverse.  And of
course, since the elementary transformation does not change the blowdown
structure past $X_1$, it suffices to consider the direct image $M_1$ of $M$ on
$X_1$.  Let $e$ be the exceptional curve on $X_1$ over $X_0$, and let
$X'_0$ be the transformed Hirzebruch surface.  The nature of the
minimal lift operation implies that we have a short exact sequence of the form
\[
0\to \sO_e(-1)^r\to \pi^*M_0\to M_1\to 0.
\]
If we take the direct image under $\pi'$ (the map that blows down the
complement $e'$ of $e$ in its fiber), then we see that $\sO_e(-1)$ is
acyclic with direct image of the form $\sO_f(-1)$, and thus
\[
0\to \sO_f(-1)^r\to \pi'_*\pi^*M_0\to \pi'_*M_1\to 0
\]
is exact.  This, of course, is precisely the situation we encounter with
non-maximal splittings; in particular, we can compute the new matrix $A'$
equally well from either $\pi'_*M_1$ or $\pi'_*\pi^*M_0$.  Now,
\[
\rho'_*\pi'_*\pi^*M_0 \cong \rho_*\pi_*\pi^*M_0\cong \rho_*M_0,
\]
and thus $W'\cong W$.  Similarly, the nonmaximal bundle $V''$ can be
computed (up to a scalar) by
\begin{align}
\rho'_*(\pi'_*\pi^*M_0(-s'-f))
&=
\rho'_*\pi'_*(\pi^*M_0(e_1-s-f))\notag\\
&=
\rho_*\pi_*(\pi^*M_0(e_1-s-f))
\cong
\rho_*(M_0(-s-f)),\notag
\end{align}
and thus $V''\cong V\otimes \sO_{\P^1}(-1)$.  We furthermore find that
the corresponding map $\rho'_*(B''(s'+f))$ on $\P^1$
factors as
\[
V\otimes \sO_{\P^1}(-1)\xrightarrow{\rho_*(B(s))}
W\otimes \rho_*\sO_{X_0}(s-f)
\xrightarrow{1\otimes \psi}
W\otimes \rho'_*\sO_{X_0}(s')
\]
where $\psi:\rho_*\sO_{X_0}(s-f)\to \rho'_*\sO_{X_0}(s')$ is the image of the
natural map $\sO_{X_1}(s-f)\to \sO_{X_1}(s-e_1)$ on $X_1$.  But then, as a
morphism on $C$, $B''$ is just the composition with the corresponding map
of vector bundles on $C$.  (If we want, we can then compute the true $B'$
by restoring maximality.)

That is, if ${\cal L}_0$ is the original (degree 2) twisting bundle, then
the new twisting bundle is ${\cal L}_1:={\cal L}_0(-p)$
where $p\in C$ is the point corresponding to $e_1$; and $B''=\psi B$ where
$\psi$ is the unique (up to scalars) global section of
\[
\sO_C(-p)\otimes\pi_\eta^*\sO_{\P^1}(1),
\]
essentially a degree 1 theta function vanishing at $\eta(p)$.  We then find
that $A' = \eta^*\psi^{-1}\psi A$.  In other words, the elementary
transformation simply multiplies the shift matrix by a ratio of two degree
1 theta functions, preserving the symmetry.

\begin{rem}
  If we were using a line bundle on a marked section to identify our
  equations with untwisted elliptic equations, then the action of $W(D_n)$
  is actually trivial, since it does not change the normalizing sheaf.
  This is consistent with the above action since the normalizing equation
  is acted on by the same scalar gauge transformation as any other
  equation, and thus the effect is cancelled when we divide by the
  normalizing equation.
\end{rem}

\medskip

The remaining simple reflection is much more subtle, as can be seen in
particular by the fact that it does not preserve the rank of the equation.
Indeed, if $D=ns+df-\sum_i r_i e_i$ is the original class, then after
reflecting in $s-f$, we obtain an equation of class $D=ds+nf-\sum_i r_i
e_i$.  Since this swaps the order of the equation and a measure of its
degree (relative to ordinary multiplication), this suggests that this
operation should correspond to some sort of generalized Fourier
transformation.  This is in fact the case, and the transform is essentially
that of Spiridonov and Warnaar \cite{SpiridonovVP/WarnaarSO:2006}; see the
discussion in \cite{generic} or Appendix \ref{chap:fourier}.

\section{Elliptic hypergeometric equations}

As we mentioned in the Introduction, we were led to consider symmetric elliptic
difference equations by their appearance in two contexts: as the equations
satisfied by elliptic hypergeometric integrals, and as equations related to
elliptic biorthogonal functions (and Painlev\'e theory).  We should
therefore explain how these equations fit into the current framework.

In \cite{dets}, Spiridonov and the author computed the explicit matrix $A$
for the difference equation satisfied by the ``order $m$ elliptic beta
integral''.  For generic parameters, these equations are nonsingular at the
ramification points, and it is thus straightforward to compute their
singularity structure.  The order $m$ elliptic beta integral satisfies an
elliptic equation of order $m+1$ with $2m+4$ ``simple zeros'', i.e., points
where $A$ is holomorphic and $\det(A)$ vanishes once, as well as
two points where $A$ vanishes identically.  This gives a sheaf of Chern class
\[
(m+1)s+(m+2)f-(m+1)e_1-(m+1)e_2-\sum_{3\le i\le 2m+6}e_i
\]
on a blowup of $F_2$ (we can recover the coefficient of $f$ by degree
considerations once we have found all the singularities).  If we perform an
elementary transformation, swap $e_1$ and $e_2$, then again perform an
elementary transformation (i.e., if we reflect in $f-e_1-e_2$), this gives
us a sheaf on a blowup of $F_0$ or $F_2$ with Chern class
\[
(m+1)s+f-\sum_{3\le i\le 2m+6}e_i,
\]
reflecting the fact that the equation given in \cite{dets} had two
singularities introduced precisely in order to make it elliptic rather than
twisted.  Since the first two blowups are independent of the remaining
blowups, we can move those to the end, then ignore them.  Thus the most
natural sheaf-relaxation of this equation has Chern class
\[
(m+1)s+f-\sum_{1\le i\le 2m+4}e_i
\]
on a surface with $K_X^2=4-2m$, relative to an even blowdown structure.
Now, this vector is actually a positive (real) root for $E_{2m+5}$, and
thus (since by construction the sheaf is supported on the complement of the
anticanonical curve) is generically the class of a $-2$-curve.  In
particular, the sheaf, and thus the difference equation, is rigid; this of
course explains why it was even possible to write down the equation
explicitly.

To verify that the vector is a positive root, we can of course apply the
usual algorithm.  The only simple root that has negative intersection with
the class is $s-f$, and thus the first step is the generalized Fourier
transformation mentioned above.  This gives a first-order equation of class
\[
s+(m+1)f-\sum_{1\le i\le 2m+4}e_i,
\]
at which point the action of $D_{2m+4}$ suffices to transform it to the
trivial equation $s-f$ (assuming sufficiently general parameters).  Since
first-order equations have explicit meromorphic solutions given by elliptic
Gamma functions \cite{RuijsenaarsSNM:1999} and the generalized Fourier
transform is at least formally an integral transformation, we see that
the above $m+1$-st order equation should have a solution expressed as an
integral involving elliptic Gamma functions.  This is, of course, hardly
surprising considering that the equation arose as the equation of an
integral, and indeed, we recover the elliptic beta integral in this way.

Now, suppose one starts with the trivial equation and performs some
sequence of elementary transformations in various points and
Spiridonov-Warnaar transformations.  This will have the effect of replacing
the original $-2$-curve $s-f$ by some image under $W(E_{m+1})$, which will
still be the class of a $-2$-curve, and thus corresponds to a rigid
equation.  (Indeed, for the relaxation, rigidity is a property of the sheaf
on $X$, so is independent of the blowdown structure.)  In terms of
solutions, this starts with $1$, and performs some sequence of the
operations ``multiply by a symmetric product of elliptic Gamma functions''
and ``apply the Spiridonov-Warnaar transform''.  One thus expects that the
result of such a sequence of operations will always satisfy a {\em rigid}
difference equation.  This is not quite an algebraic statement, as there
are analytic issues in showing that the Spiridonov-Warnaar transform takes
solutions of the original equation to solutions of the transformed equation
(mainly having to do with the possibility of boundary terms in the analogue
of integration by parts).  This can be shown, however, as well as the
converse, namely that any rigid elliptic equation has a solution given by
an iterated integral corresponding to its reduction to $s-f$, see \cite[\S
  13.5]{generic}. (This appears related to the notion of Bailey
chains/trees, see \cite{SpiridonovVP:2002,SpiridonovVP:2004} for the
elliptic case.)

\section{Semiclassical biorthogonal functions}

The other main motivating family of equations are those of
\cite{isomonodromy}, the equations satisfied by certain functions which are
biorthogonal with respect to the order $m$ elliptic beta integral.  These
equations are no longer explicit (though the residues can be expressed as
multivariate integrals of products of elliptic Gamma functions), but it is
still quite feasible to determine their singularities.  These start out
twisted, so we need to make one of the 16 compatible choices of twisting
data, but there is a natural choice making the equation nonsingular at the
ramification points, at least for generic parameters.  (Alternately, as in
the elliptic beta integral case, \cite{isomonodromy} gives a
well-controlled elliptic version of the difference equation, corresponding
to the non-elliptic version by a pair of elementary transformations.)  We
find that we are in the even case, and have a second-order equation with
$2m+6$ simple singularities.  The corresponding Chern class is thus
\[
2s+(m+1)f-\sum_{1\le i\le 2m+6} e_i.
\]
When $m=0$, this is rigid (and indeed all of the multivariate integrals
arising as coefficients can be explicitly evaluated); this is of course one
of the rigid cases we just saw, corresponding to the fact that the order 0
elliptic beta integral admits hypergeometric biorthogonal functions.  (This
generalizes the fact that the Jacobi polynomials, orthogonal with respect
to the usual beta integral, are hypergeometric.)  Otherwise, the class is
in the fundamental chamber, so nef, and is easily checked to be generically
integral.  Since the coefficients of the Chern class are relatively prime,
we find that the corresponding moduli space of sheaves is rational, and
thus the same is true for the moduli space of difference equations (at
least for the components where the moduli spaces are birational; of course,
then $q$-deformed twisting should make this true for the difference
equations in any component).

The case $m=1$ is of particular interest, as in that case the
$2$-dimensional moduli space is an open subvariety of an elliptic surface,
the relative Jacobian of the original $X_8$ (which is isomorphic to $X_8$
given the choice of a section).  More precisely, it is obtained from the
relative Jacobian by removing both the Jacobian of $C_\alpha$ and the
divisor where the bundle fails to be trivial.  As we mentioned, this is
just the theta divisor, so corresponds to the identity section of the
relative Jacobian.  Since a section of a rational elliptic surface is a
$-1$-curve (it is smooth, rational, and meets $C_\alpha$ in a single
point), we could just as well blow down that section before removing it.
This gives us an alternative interpretation of the moduli space as the
complement of $C_\alpha$ in a del Pezzo surface of degree 1: the section
blows down to a point of $C_\alpha$, so there is nothing else to remove.
The $q$-deformed twisting operations discussed above must still act as
rational maps on this del Pezzo surface, and one can check directly that
those actions all factor through blowing up the point of $C_\alpha$
corresponding to $q$.  In other words, the moduli space of difference
equations of this type is precisely the sort of rational surface studied by
\cite{SakaiH:2001}, and the twist operations correspond directly to the
elliptic Painlev\'e equation.  Since the divisor class is anticanonical,
its expansion in the standard basis is invariant under the action of
$W(E_9)$ on blowdown structures, and thus we obtain an action of $W(E_9)$
on the family of moduli spaces of difference equations.  The translations
(which is how Sakai defined the elliptic Painlev\'e equation in
\cite{SakaiH:2001}) act the same way on the parameters as the twist
operations, and thus must actually act in the same way on the moduli
spaces.  In other words, there is a large abelian subgroup of
$\Z^{10}\rtimes W(E_9)$ that acts trivially on the moduli space; this
almost certainly is special to the $m=1$ case.  In particular, from a
difference equation perspective, the correct way to generalize the elliptic
Painlev\'e equation is clearly the twisting action rather than the Coxeter
group action, as the former is what acts by gauge (i.e., isomonodromy)
transformations in general.  (In fact, the relation between these
interpretations is mediated by a derived equivalence, see \cite{generic},
or the generalization in Theorem \ref{thm:weird_langlands} below.)

For general $m$, the fact that $B$ is a morphism between rank 2 bundles and
becomes rank 1 at the singular points (at least generically) allows us to
define a number of rational functions on the moduli space.  A typical
example is the following: when $W$ is trivial, the $1$-dimensional
subspaces $\ell_i:=\im(B(p_i))$ determine eight points in a common
projective line, and thus we can take the cross-ratio of any four of them
to obtain a rational function on the moduli space.  This function has a
particularly nice divisor, and can formally be written as a ratio of tau
functions.  More precisely, if for any $v\in \Z^{2m+8}$, we define
$\tau(v)$ to be the tau function which vanishes when $H^0(M(v-f))\ne 0$, then
\[
\chi(\ell_1,\ell_2,\ell_3,\ell_4)
\propto
\frac{\tau(f-e_1-e_3)\tau(f-e_2-e_4)}
     {\tau(f-e_1-e_4)\tau(f-e_2-e_3)}.
\]
(This is really just a statement about the divisor of the left-hand side.)
This follows from the observation that (for generic parameters)
$\ell_i=\ell_j$ precisely when twisting by $-e_i-e_j$ changes $W$ from
$\sO_{\P^1}^2$ to $\sO_{\P^1}(-2)\oplus \sO_{\P^1}$, an easy consequence of
our description of how twisting affects $W$.  (It is not quite clear how
this should extend to the locus $\tau(0)$ where $W$ itself is not trivial,
but it is clear that whatever order it has along $\tau(0)$ will not depend
on $i$ or $j$.)  If we take into account the behavior when
$\im(B(p_i))=\im(B(p_j))$, this suggests a more precise statement
\[
\chi(\ell_1,\ell_2,\ell_3,\ell_4)
\propto^?
\frac{\theta(e_1-e_3)\tau(f-e_1-e_3)\theta(e_2-e_4)\tau(f-e_2-e_4)}
     {\theta(e_1-e_4)\tau(f-e_1-e_4)\theta(e_2-e_3)\tau(f-e_2-e_3)},
\]
where for any positive root $r$, $\theta(r)$ is the divisor on the moduli
stack of surfaces that vanishes where the root is a $-2$-curve.  This is
very suggestive of the formula for the cross-ratio given in
\cite{isomonodromy}, in which there are multivariate integrals taking the
places of factors
\[
\frac{\tau(f-e_i-e_j)}{\theta(f-e_i-e_j)};
\]
and further suggests that we should have an equation of the form
\begin{align}
&\theta(e_1-e_2)\tau(f-e_1-e_2)\theta(e_3-e_4)\tau(f-e_3-e_4)\notag\\
{}-{}&
\theta(e_1-e_3)\tau(f-e_1-e_3)\theta(e_2-e_4)\tau(f-e_2-e_4)\notag\\
{}+{}&
\theta(e_1-e_4)\tau(f-e_1-e_4)\theta(e_2-e_3)\tau(f-e_2-e_3)=^?0
,\notag
\end{align}
though it is as yet unclear how to make precise sense of such a statement.
(The main issue is producing suitable canonical isomorphisms between tensor
products of bundles of the form $\det \dR\Gamma$, since $\theta$ and $\tau$
are canonical global sections of such bundles.)  Note in particular that
when $m=1$, if we could make sense of the above relation, then it would
induce an entire $W(E_9)$-orbit of relations, which are precisely the
relations satisfied by a tau function for the elliptic Painlev\'e equation,
as given in \cite[Thm.~5.2]{KajiwaraK/MasudaT/NoumiM/OhtaY/YamadaY:2006}.

\begin{rem}
  In \cite{recur}, a four-term variant (related to Pl\"ucker relations
  for pfaffians) of this $W(E_9)$-invariant system of recurrences
  was introduced, corresponding to a slightly different family of
  multivariate hypergeometric integrals.  Is there a corresponding
  model (presumably involving relative Pryms rather than relative
  Jacobians) from a geometric perspective?  This would presumably
  involve a component of the fixed locus of an involution of the
  form $M\mapsto \iota^*\sExt^1(M,\omega_X)$, where $\iota$ is an
  involution on an anticanonical surface of the form $X=X_8$ that
  acts as a hyperelliptic involution on the relevant anticanonical
  curve.  (The latter condition ensures that the combined involution is
  Poisson.)
\end{rem}

\medskip

The case $D=2rs+2rf-\sum_{1\le i\le 8} re_i$ is also likely to have
interesting behavior (assuming it is generically integral, i.e., that
$\sO_{C_\alpha}(C_\alpha)$ has exact order $r$ on $X_8$).  In this
case, $X_8$ is an elliptic surface on which $C_\alpha$ appears as the
underlying curve of an $r$-fold section, but the moduli space is still an
open subset of the relative Jacobian.  As we saw, the moduli space of
matrices $B$ is always rational in the $2$-dimensional case, and the same
reasoning as before tells us that it is an affine del Pezzo surface of
degree 1.  Again, for any $q$, we have an induced action of $\Z^{10}\rtimes
W(E_9)$ as birational maps on this family of del Pezzo surfaces.  Since
Theorem \ref{thm:Picr} tells us the parameters of these del Pezzo surfaces,
we can control the action of the birational maps, and find that the action
factors through a suitable one-point blowup, just as in the case $r=1$.
Note, however, that the lack of a universal family on this moduli space
means that the most obvious way of producing a corresponding Lax pair will
not work.

The situation is somewhat better if we take the same $D$ but ask for the
Euler characteristic to be some $d$ with $\gcd(r,d)=1$; there, the moduli
space is again expressible as a rational surface with $K^2=0$, so that one
again recovers elliptic Painlev\'e.  (For the relaxed case, this follows
from a relative version of the standard derived autoequivalences of
elliptic curves; in general, this requires derived equivalences of
noncommutative surfaces, see Theorem \ref{thm:Painleve_moduli_spaces}
below.)

\section{Isomonodromy and apparent singularities}

Although we have mainly focused on the irreducible case, we should note
that even in generically irreducible cases, the moduli space can contain
reducible fibers (or even fibers with nonreduced support), and because we
were careful to develop the theory without relying too much on the notion
of ``support'', everything mostly carries over, with the main issues being
that stability can fail, and that components of Chern class $f$ can
disappear on taking the direct image to the ruled surface, or can make the
vector bundles degenerate to have torsion.  Even in that case, one can
still make much of the theory work.

Fix a symmetric elliptic difference equation (with elliptic curve $E$, and
which we suppose is untwisted)
\[
v(q+z) = A(z)v(z)
\]
supported on a vector bundle $\pi^* W^*$, and for $x\in E$, $\lambda \in
W^*_{\pi(x)}$ (nonzero), $l\in \N$, define a subbundle $W_{l,x,\lambda}$ by
a short exact sequence
\[
0\to W_{l,x,w}\to W\to \sO_{\pi(x)}\oplus \cdots\oplus \sO_{\pi(x+(l-1)q)}\to 0,
\]
where the map $W\to \sO_{\pi(x+iq)}$ is given by
\[
A(x+(i-1)q)\cdots A(x+q)A(x)\lambda,
\]
which makes sense as long as $A$ is regular on the orbit of $x$.  (If $A$
is regular on the orbit of $y$, then $A(y)$ induces an isomorphism between
$(\pi^*W)_y^*$ and $(\pi^*W)_{y+q}^*$.)  Let
$C_{l,x,\lambda}:W_{l,x,\lambda}\to W$ be the natural inclusion.
Then we can gauge by $C_{l,x,\lambda}$ to obtain a new equation
with shift matrix
\[
A_{l,x,\lambda}(z) = C_{l,x,\lambda}(z+q)^t A(z) C_{l,x,\lambda}(z)^{-t}.
\]
A priori, this has singularities at $x,\dots,x+nq$ (and their images under
$z\mapsto -q-z$), but the choices of linear functionals above ensure that
the singularities at $x+q,\dots,x+(l-1)q$ cancel, and thus one has only a
pole at $x$ (with residue of rank 1) and a zero at $x+lq$ (with $A(x+lq)$
having kernel of dimension 1).

Since $W_{l+1,x,\lambda}\subset W_{l,x,\lambda},W_{l,x+q,A(x)\lambda}$, we
see that $A_{l+1,x,\lambda}$ is a pseudo-twist of $A_{l,x,\lambda}$, and
$A_{l,x+q,A(x)\lambda}$ is a pseudo-twist of $A_{l+1,x,\lambda}$.  In other
words, the pair of pseudo-twists (in the two introduced singular points)
takes $A_{l,x,\lambda}$ to $A_{l,x+q,A(x)\lambda}$ for all $n>0$,
$x,\lambda$.  It follows that if $\lambda=v(x)$ for some solution of the
symmetric elliptic difference equation, then iterating the pseudo-twists
gives $A_{l,x+kq,v(x+kq)}$.  In other words, any solution of the original
equation on the orbit of $x$ induces a sequence of isomonodromic equations
with apparent singularities.

Another way of saying this is that for all $l>0$, there is a map from the
space of solutions of the original equation to the space of solutions of
the corresponding isomonodromy equation.  More precisely, there is an {\em
  embedding} of the space of {\em projective} solutions of the original
equation to the space of solutions of the corresponding isomonodromy
equation.  (This is because $W_{l,x,\lambda}$ is invariant under rescaling
$\lambda$.)  That we lose information about scalar multiples should not, of
course, be surprising, since the isomonodromy equation is actually
intrinsic to the geometry, while the interpretation of the linear equation
requires a choice of normalizing section.  (We could, of course, twist by a
line bundle instead, but then the shift matrix only identifies adjacent
fibers up to a scalar!)

To understand the geometry here, we need to look at what happens to $B$
rather than $A$.  But this is straightforward: the same conditions that
ensure that the intermediate singularities we are introducing cancel
ensures that we can eliminate those singularities by changing $V$.  In
particular, we see that $V_{l,x,\lambda}$ is the kernel of an appropriate
map
\[
V\to \sO_{\pi'(-q-x)}\oplus\cdots\oplus \sO_{\pi'(-(l-1)q-x)}
\]
In the relaxed case $q=0$ (which only works for $l$ less than the order of
the equation, and under the assumption that the relevant linear functionals
are linearly independent), we are decreasing $\chi(V)$ by $l-1$ and
$\chi(W)$ by $l$, and thus we obtain an equation
\[
  [M_{l,x,\lambda}] = [M]+[\sO_f(-l-1)]
\]
in the Grothendieck group of $X$.  Since we also have a morphism
$M_{l,x,\lambda}\to M$, we deduce that $M_{l,x,\lambda}$ is the extension
of $M$ by $\sO_{f}(-l-1)$, where the fiber meets the anticanonical curve in
$x$.  (Something similar of course holds in the noncommutative setting,
where the fact that the sheaf of class $\sO_f(-l-1)$ meets the
anticanonical curve in $x$ implies that it also meets the anticanonical
curve in $-lq-x$.)

The appearance of a sheaf of class $f$ here is, of course, to be expected,
since it does not make sense to even talk about the corresponding equation
without first choosing a ruling.  What this does tell us, however, is that
the relation between sheaves and equations is much more intrinsic than it
may have appeared from the somewhat ad hoc nature of our above
constructions: the equation associated to a sheaf on a noncommutative
surface (once we have made sense of this!) can be recovered up to scalars
from the behavior under pseudo-twists of its extensions by sheaves of the
form $\sO_f(-l-1)$.

Some further comments about this construction: First, we could have
computed our composition of pseudo-twists in the opposite order, i.e., by
passing through $A_{l-1,x+q,A(x)\lambda}$ instead of $A_{l+1,x,\lambda}$.
There is a significant caveat here, however: since $A_{0,x,\lambda}=A$ (we
are dealing with an extension by $\sO_f(-1)$), commutativity breaks for
$l=1$, and thus one must either work in the order we used above, or follow
a suitable infinitesimal deformation of $A_{1,x,\lambda}$ through the
pseudo-twists.  Second, this construction gives a more direct explanation
of why the semiclassical (bi)orthogonal function solutions can be expressed
via hypergeometric functions: these correspond to extensions by sheaves of
the form $\sO_{s-f}(-l-1)$.  Applying a pair of elementary transformations
turns this into $\sO_{s-e_1-e_2}(-l-1)$, which after taking the Fourier
transform gives $\sO_{f-e_1-e_2}(-l-1)$.  Blowing down $e_1$, $e_2$ gives
an extension by $\sO_f(-l-1)$ of a sheaf of Chern class
\[
(m+1)s+f-e_3-e_4-\cdots-e_{2m+6},
\]
corresponding to the equation of the elliptic beta integral of order $m$ (and
degenerations thereof).  Finally, any isomonodromy transformation
of the original equation extends to an isomonodromy transformation of the
equation with removable singularities, and thus identifies the solutions of
the two projective difference equations, recovering the fact that
isomonodromy transformations in the geometric sense correspond to gauge
transformations of the associated equations.

An analogous construction works in the differential case.  The key idea is
that we can work locally near the original singularity (say $0$).  It is
instructive to first consider the case of a trivial equation of order $n$,
which after the gauge transformation to introduce an apparent singularity
at $u$ near $0$ has the form
\[
v'(z,u) = (l(z-u)) \Pi_1 v(z,u)
\]
where $\Pi_1$ is an idempotent of trace $1$ (corresponding to gauging by
$(1-\Pi_1)+(z-u)^{-l}\Pi_1$).  This can be completed to a Lax pair by taking
\[
v_u(z,u) = (-l/(z-u)) \Pi_1 v(z,u)
\]
(i.e., gauging the trivial Lax pair by the same matrix) and thus in the
trivial case the image of the leading term of the matrix we are gauging by
indeed remains constant.

More generally, if $v'(z)=A(z)v(z)$ is regular at $0$, then we can
solve the equation in a neighborhood of $0$ and thus get an expression
of the form
\[
A(z) = M'(z) M(z)^{-1}
\]
where $M(z)$ is a holomorphic matrix which we may assume satisfies
$M(0)=1$.  In particular, gauging the Lax pair
\[
v'(z,u) = A(z) v(z,u)\qquad v_u(z,u) = v(z,u)
\]
by $M(z)$ give the trivial Lax pair, and thus gauging by $M(z)
((1-\Pi_1)+(z-u)^{-l}\Pi_1)$ gives the above Lax pair with apparent
singularity.  Any other Lax pair with the same kind of apparent singularity
can be obtained by a further invertibly holomorphic gauge transformation,
and thus the gauge transformation introducing the desired apparent
singularity to the nontrivial equation has the form
\[
C(z,u) = M(z) ((1-\Pi_1)+(z-u)^{-l}\Pi_1) F(z,u)
\]
for some invertibly holomorphic matrix $F$.  We thus see that the image of
the leading term of the gauge matrix introducing the apparent singularity
is nothing other than $M(u) \im(\Pi_1)$, i.e., the projective solution of
$v'(z)=A(z)v(z)$ with initial condition $\im(\Pi_1)$.  (A modified version
works for $l<0$, since for that range the image of the leading term of
$C^{-t}$ is $M(u)^{-t} \im(\Pi_1^t)$.  This of course is related by duality
to the $l>0$ case, and thus something similar applies in the difference
case.)

In general, we can apply the generalized Fourier transform to move the sub-
or quotient sheaf of class $f$ to have class $s$, or indeed any class
corresponding to a rational curve of self-intersection 0, and in each case
the isomonodromy transformations moving the auxiliary intersections will
correspond to a projective version of a linear difference (or differential)
equation.  This generalizes further: if $D_r$ is the class of a rational
curve with $D_r^2\ge 0$, then we may consider extensions by (or of) a sheaf
of class $D_r$, and look at the isomonodromy transformations that move the
additional $D_r^2+2$ intersections with the anticanonical curve.  If we
also include permutations of those intersections, we obtain an action of
the affine Weyl group $\tilde{A}_{D_r^2+1}$ on the projective space of
extensions, which extends to a $\PGL$-valued cocycle of
$\tilde{A}_{D_r^2+1}$ with coefficients in an $S_{D_r^2+2}$-cover of the
function field of the projective space of sheaves of class $D_r$
(specifying the $D_r^2+2$ points of intersection).  This is likely to be
particularly interesting in the case when the original sheaf is rigid, as
then the resulting system of partial difference equations is a compatible
family of rigid ordinary difference equations.

\section{Degenerations}
\label{sec:degen}

As one might expect, the story becomes more complicated once the
anticanonical curve becomes singular.  The simplest case is when the
anticanonical curve on $X$ is still integral; in that case, the
considerations of the previous section carry over with little change.  The
main constraint is that the symmetric (ordinary and $q$-) difference
equations must have only finite singularities, since blowing up the node or
cusp will introduce a new component to the anticanonical curve.  This can
actually be violated in a mild way: if the difference equation is twisted
by a line bundle, we can use the singularity to single out a global section
of the bundle (modulo scalars), and in this way obtain an untwisted
equation with only a mild singularity at $\infty$.  (E.g., in the
$q$-difference case, the matrix $A$ will no longer be $1$ at $\infty$, but
will still be a multiple of the identity.)  (Equations with more
complicated singularities at $\infty$ might correspond to sheaves that
cannot be separated from the anticanonical curve by a suitable blowup,
though this can always be fixed by a finite number of pseudo-twists, and we
have seen that these correspond to gauge transformations.)

For more degenerate cases, we can still be guided by what happens in the
elliptic case.  For the $W(E_{m+1})$ action, the $A_{m-1}$ subsystem merely
permutes the singularities, while elementary transformations change the
twisting bundle and multiply the shift matrix by a corresponding
meromorphic section.  The reflection in $s-f$, in contrast, can have a more
significant effect on the nature of the equation.  

We have already considered how the different equations look on $F_2$, and
something similar applies to $F_0$ or $F_1$.  Indeed, since elementary
transformations should not affect the type of equation, but can introduce
or contract fibers which are components of $C_\alpha$, the rule is quite
simple: contract any components of class $f$, and then recognize the curve from
the same list of possibilities as for $F_2$.  (More precisely, choose any
section of the ruling which is transverse to $C_\alpha$, and perform a
sequence of elementary transformations moving that section to the
$-2$-curve of an $F_2$; the resulting anticanonical curve will be disjoint
from $C_\alpha$, and differs from the original only in self-intersections of
components and the contraction of fibers.)  Since this rule treats sections
and fibers differently, it is clear that the result can depend on the
choice of ruling on $F_0$.

Moreover, it is in general not possible to avoid this issue.  One might be
tempted to adapt the algorithms for checking integrality to use only
elements of the group that stabilizes the decomposition of $C_\alpha$, but
this encounters two significant problems: the group need not be a
reflection group, and the reflection subgroup need not have finite rank.
Either possibility denies us any kind of useful ``fundamental chamber'';
there need not be any computable fundamental domain for the action.  One
must thus use the full algorithm, and this can most certainly change the
kind of equation.

\medskip

As an example of how birational maps can change the type of an equation,
consider the case of a nonsymmetric $q$-difference equation with three
polar singularities.  We recall that such equations (with even twisting; as
mentioned, this allows $A(\infty)$ to be a general multiple of the
identity) correspond to sheaves on $\P^1\times \P^1$ with specified
intersection with a union of two bilinear curves meeting in two distinct
points.  The constraint on singularities means that the sheaf meets the
component of $C_\alpha$ corresponding to poles in the specified three
points (and the restriction is a direct sum of structure sheaves of
subschemes of the corresponding degree $3$ scheme).  Now, a bilinear curve
in $\P^1\times \P^1$ has self-intersection $2$, so after blowing up three
points, the strict transform has self-intersection $-1$.  The assumption on
singularities means that this $-1$-curve $e$ is disjoint from the lift of
$M$, and thus there are blowdown structures for which $c_1(M)$ is in the
fundamental chamber and such that $e_m=e$.  Such a blowdown structure blows
down one of the two components of the anticanonical curve, so produces an
integral anticanonical curve on the eventual Hirzebruch surface.

In other words, there is a birational map taking nonsymmetric
$q$-difference equations with three polar singularities (and regular at $0$
and $\infty$, modulo twisting) to symmetric $q$-difference equations (no
longer regular at $0$ and $\infty$).  By looking at what the algorithm does
to move $e$ to $e_m$, we see that this involves reflecting in $s-f$
precisely once.  Before reflecting, the anticanonical curve on $\P^1\times
\P^1$ has two components with classes $s+2f$ and $s$, while after
reflecting it has components of classes $2s+f$ and $f$.

\begin{rem}
Of course, something similar applies if we have more than three polar
singularities or the nonsymmetric equation is singular at $0$ or $\infty$;
the only difference is that the resulting symmetric equation will be more
singular at $0$ and $\infty$.
\end{rem}

On $F_0$, we have a total of 16 possible ways the anticanonical curve can
decompose, each of which corresponds to a different kind of generalized
Fourier transform.  There are 10 such transforms that preserve the type of
equation, and three pairs that change the type.  Of those, one changes
between symmetric and nonsymmetric $q$-difference equations, one is the
ordinary difference analogue, and a final one changes between differential
and ordinary difference equations (a Mellin/z transform).  For the
transforms that preserve type, there are one each for the three integral
types, as well as three transforms on nonsymmetric $q$-difference
equations, two for ordinary difference equations, and two for differential
equations.  Of the latter, the most degenerate is just the Fourier/Laplace
transform, while the other is essentially the transform used in
\cite{KatzNM:1996} (usually called ``middle convolution'' in the later
literature, though this is something of a misnomer).  See Appendix
\ref{chap:fourier} for more details.

\medskip

As we mentioned above, we can also model certain Deligne-Simpson problems
via moduli spaces of sheaves on rational surfaces, and this gives rise to
additional interesting maps of moduli spaces.  One, of course, is the
(essentially trivial) observation that moduli spaces of Fuchsian
differential equations correspond to moduli spaces of solutions to additive
Deligne-Simpson problems; in our terms, we can see this by noting that the
anticanonical curves in the latter case become a double $\P^1$ once we
contract all fiber components.  There is another relation, though, which we
consider in the multiplicative case.  Recall that we modeled four-matrix
Deligne-Simpson problems via sheaves with a quadrangular anticanonical
curve in $\P^1\times \P^1$ (with components of class $f$, $f$, $s$, and
$s$).  This is a somewhat cumbersome model for three-matrix problems, as we
need to take one of the four matrices to be the identity.  The
corresponding component of the anticanonical curve contains a single
singular point; if we blow up this point and blow down both the fiber and
the section containing it, we obtain a sheaf on $\P^2$.  If $g_1$, $g_2$,
$g_3$ are the three matrices with product $1$, the sheaf on $\P^2$ is
modeled by the cokernel of the matrix $x+g_1y+g_1g_2z$, and the conjugacy
classes are determined by the restriction to the anticanonical curve
$xyz=0$.  (The additive variant involves sheaves with specified restriction
to $xy(x+y)=0$.)

Much as in the case of a nonsymmetric difference equation with three poles,
a three-matrix Deligne-Simpson problem in which one matrix has a quadratic
minimal polynomial gives rise to a surface in which the anticanonical curve
contains a $-1$-curve disjoint from the relevant sheaf.  In particular, we
obtain a sheaf on an even Hirzebruch surface by blowing up the two roots of
the minimal polynomial, then blowing down the anticanonical component.  On
that Hirzebruch surface, the anticanonical curve has two components, both
of class $s+f$, and we thus obtain a nonsymmetric $q$-difference equation.
(More precisely, we obtain such an equation {\em after} choosing one of the
two rulings; this is tantamount to choosing one of the two roots of the
minimal polynomial.)

This map can be made precise as follows.  Let $g_0$, $g_1$, $g_\infty$ be
a solution to the Deligne-Simpson problem, with
$(g_1-1)(g_1-\beta^{-1})=0$.  (We have rescaled the chosen root of the
minimal polynomial to $1$.)  Define a matrix $A^+(z)$ with rational
coefficients by
\[
A^+(z) = (1-g_\infty^{-1}z)^{-1}(1-g_0z)
       = (g_\infty-z)^{-1}(g_\infty-g_1^{-1}z).
\]
The matrix $A^+(z)$ fixes any vector fixed by $g_1$, and thus has a
well-defined action on the quotient $\im(g_1-1)$.  If $A(z)$ is the matrix
of this action in some basis, then we find $A(0)=1$, $A(\infty)=\beta$, so
that $A$ represents a twisted $q$-difference equation which is regular at
$0$ and $\infty$.  Moreover, we see that the zeros of $A$ occur at the
eigenvalues of $g_0^{-1}$, and the poles of $A$ occur at the eigenvalues of
$g_\infty$, and thus we obtain a rational map between the two moduli
spaces.  (It is unclear how to make the inverse map explicit, though it
certainly exists, due to the description in terms of sheaves.)  Note that
if $g_\infty$ has a cubic minimal polynomial, then we can proceed further,
turning this $q$-difference equation with three polar singularities into a
symmetric $q$-difference equation.  Similarly, a solution to a three-matrix
additive Deligne-Simpson problem with a quadratic minimal polynomial
produces an ordinary difference equation, which can be further transformed
to a symmetric equation if another minimal polynomial is cubic.

\section{Rigid equations}

Just as the elliptic hypergeometric equation corresponds to a $-2$-curve,
so is rigid, most other hypergeometric difference/differential equations
can be seen to be rigid in the same way.  As an example, we consider a
maximally degenerate case: the Airy function, which satisfies the
non-Fuchsian differential equation
\[
\Ai''(z) = z\Ai(z),
\]
or in matrix form
\[
v'(z) = \begin{pmatrix} 0 & z \\ 1 & 0\end{pmatrix} v(z).
\]
(Note that there is an ambiguity when passing between straight-line and
matrix forms of an equation: the matrix form is only determined up to a
gauge transformation, so (as long as we can avoid apparent singularities)
the sheaf will only be determined up to ``pseudo-twist''.)  As we have
seen, differential equations correspond to sheaves on $F_2$ with
anticanonical curve of the form $y^2=0$; we find that the above matrix
translates to
\[
\begin{pmatrix}
w^2 & yx\\
y & w^3
\end{pmatrix}
:
\sO_{F_2}(-s_{\min}-2f)\oplus \sO_{F_2}(-s_{\min}-3f) \to \sO_{F_2}^2.
\]
In the $q$-difference case, we noted that we can often absorb particularly
well-behaved singularities into a twist; something similar applies here,
and we should perform an elementary transformation centered at the
subscheme with ideal $(y,w^2)$.  That is, blow up this subscheme, minimally
desingularize, then blow down the original fiber and the $-2$-curve.  The
resulting morphism on $\P^1\times \P^1$ is
\[
\begin{pmatrix}
y_0 & y_1 x_1\\
y_1 & y_0 x_0
\end{pmatrix},
\]
where the new coordinates relate to the original coordinates by
\[
x_1/x_0 = x/w,\qquad y_1/y_0 = y/w^2,
\]
and the new anticanonical curve has equation $y_1^2x_0^2$.  The cokernel
has smooth support, so there is no difficulty in resolving its intersection
with the anticanonical curve.  The support $y_0^2x_0=y_1^2x_1$ meets the
anticanonical curve in a $6$-jet, and each of the first five blowups
introduces a new component to the anticanonical curve, with multiplicities
$3$, $4$, $3$, $2$, and $1$ respectively.  The result is a curve of Kodaira
type $\text{III}^*$, except with one reduced fiber removed.  As for the sheaf
itself, we started with support of class $2s+f$ and blew up 6 points, so
the result is a $-2$-curve as expected.  (Had we started from the sheaf on
$F_2$, we would have obtained a sheaf with first Chern class class
$2s+3f-2e_1-2e_2-e_3-e_4-e_5-e_6-e_7-e_8$, which is naturally also a
positive (real) root.)

In general, a rigid second-order equation (of whatever kind) can always be
transformed by a sequence of elementary transformations into one of class
$2s+f-e_1-e_2-e_3-e_4-e_5-e_6$.  (That is, this is the unique class in the
fundamental chamber with respect to $D_m$ and satisfying $D\cdot f=2$.)  We
can thus describe a moduli space of rigid equations, namely the locally
closed substack of ${\cal X}^{\alpha}_6$ on which this class represents a
$-2$-curve.  Since we can readily rule out the existence of $-d$-curves for
$d>2$, it is a (reasonably small) finite problem to determine the different
types of anticanonical surfaces that arise.  There are a total of $3182$
such types, but the action of $W(D_6)$ reduces this to only $41$
equivalence classes, each of which describes a different kind of
hypergeometric equation.  These range from the elliptic hypergeometric
equation (satisfied by the order 1 elliptic beta integral) down to the two
maximally degenerate cases, the Airy equation and the $q$-difference
equation $v(q^2z)=\beta z v(z)$, passing through examples such as the
differential and difference equations satisfied by the hypergeometric
function of type ${}_2F_1$.

In general, given such a classification, it remains an open problem to
determine the corresponding limiting relations, or equivalently to
determine which types are contained in the closures of other types.  In
this case, there is a simplification available: if we choose a fiber
transverse to the anticanonical curve and perform the two additional
blowups needed to separate it from the anticanonical curve, then the
resulting surface has the form $X_8$, with the strict transform of the
fiber being $f-e_7-e_8$.  We can do this for a surface over a dvr (if the
fiber is transverse on the surface over the residue field, it will be
transverse generically as well), and thus see that any limiting relation
between types of rigid second-order equations induces a limit between types
of surfaces with $K^2=0$ and $f-e_7-e_8$ a $-2$-curve disjoint from the
anticanonical curve.  Conversely, given a limit between such surfaces,
blowing down the last two points gives a limit between surfaces of the type
we want to classify.  We can then (in the relaxed case) use an integral
fiber of the anticanonical pencil to label the surfaces by line bundles on
that curve, and thus see that a na\"{i}ve morphism of types is effective
iff the corresponding inclusion of lattices is saturated.  (The lattice a
priori should include $f-e_7-e_8$, but since this is orthogonal to the
remaining lattice, it may safely be ignored.)

We can simplify the latter question somewhat by noting that the stabilizer
of $f-e_7-e_8$ in $W(E_9)$ is isomorphic to $W(E_7)\ltimes \Z^8$, so that
we can enlarge the symmetry group without changing what limiting relations
exist.  This leads to the following diagram of types (Figure
\ref{fig:class2F1}; Figure \ref{fig:num2F1} indicates how the different
types break up into different types of equation), which we may now label by
the type of the anticanonical curve.  (Here $A_5$ splits into two cases, as
there are two different orbits of $A_5$ subsystem inside $E_7$.)

\begin{figure}
\centering
\begin{tikzpicture}[scale = 1]
\begin{scope}
\coordinate (P11s) at (0,0);
\coordinate (P11e) at ($(P11s)-(0.8,0)$);
\coordinate (P10s) at ($(P11e)-(0.6,0)$);
\coordinate (P12e) at ($(P11s)+(0.6,0)$);
\coordinate (P12s) at ($(P12e)+(0.8,0)$);
\coordinate (P13e) at ($(P12s)+(0.6,0)$);
\coordinate (P13s) at ($(P13e)+(0.8,0)$);
\coordinate (P14e) at ($(P13s)+(0.6,0)$);
\coordinate (P14s) at ($(P14e)+(0.8,0)$);
\coordinate (P15e) at ($(P14s)+(0.6,0)$);
\coordinate (P15s) at ($(P15e)+(0.8,0)$);
\coordinate (P16e) at ($(P15s)+(0.6,0)$);
\coordinate (P16s) at ($(P16e)+(0.8,0)$);
\coordinate (P17e) at ($(P16s)+(0.6,0)$);
\coordinate (P17s) at ($(P17e)+(0.8,0)$);
\coordinate (P18e) at ($(P17s)+(0.6,0)$);
\coordinate (P18s) at ($(P18e)+(0.8,0)$);
\coordinate (P19e) at ($(P18s)+(0.6,0)$);
\coordinate (P19s) at ($(P19e)+(0.8,0)$);
\coordinate (P21s) at (0,-1);
\coordinate (P21e) at ($(P21s)-(0.8,0)$);
\coordinate (P20s) at ($(P21e)-(0.6,0)$);
\coordinate (P22e) at ($(P21s)+(0.6,0)$);
\coordinate (P22s) at ($(P22e)+(0.8,0)$);
\coordinate (P23e) at ($(P22s)+(0.6,0)$);
\coordinate (P23s) at ($(P23e)+(0.8,0)$);
\coordinate (P24e) at ($(P23s)+(0.6,0)$);
\coordinate (P24s) at ($(P24e)+(0.8,0)$);
\coordinate (P25e) at ($(P24s)+(0.6,0)$);
\coordinate (P25s) at ($(P25e)+(0.8,0)$);
\coordinate (P26e) at ($(P25s)+(0.6,0)$);
\coordinate (P26s) at ($(P26e)+(0.8,0)$);
\coordinate (P27e) at ($(P26s)+(0.6,0)$);
\coordinate (P27s) at ($(P27e)+(0.8,0)$);
\coordinate (P28e) at ($(P27s)+(0.6,0)$);
\coordinate (P28s) at ($(P28e)+(0.8,0)$);
\coordinate (P29e) at ($(P28s)+(0.6,0)$);
\coordinate (P29s) at ($(P29e)+(0.8,0)$);
\coordinate (P31s) at (0,-2);
\coordinate (P32e) at ($(P31s)+(0.6,0)$);
\coordinate (P32s) at ($(P32e)+(0.8,0)$);
\coordinate (P33e) at ($(P32s)+(0.6,0)$);
\coordinate (P33s) at ($(P33e)+(0.8,0)$);
\coordinate (P34e) at ($(P33s)+(0.6,0)$);
\coordinate (P34s) at ($(P34e)+(0.8,0)$);
\coordinate (P35e) at ($(P34s)+(0.6,0)$);
\coordinate (P35s) at ($(P35e)+(0.8,0)$);
\coordinate (P36e) at ($(P35s)+(0.6,0)$);
\coordinate (P36s) at ($(P36e)+(0.8,0)$);
\coordinate (P37e) at ($(P36s)+(0.6,0)$);
\coordinate (P37s) at ($(P37e)+(0.8,0)$);
\coordinate (P38e) at ($(P37s)+(0.6,0)$);
\coordinate (P38s) at ($(P38e)+(0.8,0)$);
\coordinate (P39e) at ($(P38s)+(0.6,0)$);
\coordinate (P39s) at ($(P39e)+(0.8,0)$);
\coordinate (P41s) at (0,-3);
\coordinate (P42e) at ($(P41s)+(0.6,0)$);
\coordinate (P42s) at ($(P42e)+(0.8,0)$);
\coordinate (P43e) at ($(P42s)+(0.6,0)$);
\coordinate (P43s) at ($(P43e)+(0.8,0)$);
\coordinate (P44e) at ($(P43s)+(0.6,0)$);
\coordinate (P44s) at ($(P44e)+(0.8,0)$);
\coordinate (P45e) at ($(P44s)+(0.6,0)$);
\coordinate (P45s) at ($(P45e)+(0.8,0)$);
\coordinate (P46e) at ($(P45s)+(0.6,0)$);
\coordinate (P46s) at ($(P46e)+(0.8,0)$);
\coordinate (P47e) at ($(P46s)+(0.6,0)$);
\coordinate (P47s) at ($(P47e)+(0.8,0)$);
\coordinate (P48e) at ($(P47s)+(0.6,0)$);
\coordinate (P48s) at ($(P48e)+(0.8,0)$);
\coordinate (P49e) at ($(P48s)+(0.6,0)$);
\coordinate (P49s) at ($(P49e)+(0.8,0)$);
\node at ($(P10s)-(0.4,0)$){$A_0^e$};
\node at ($(P21s)-(0.4,0)$){$A_0^*$};
\node at ($(P22s)-(0.4,0)$){$A_1^*$};
\node at ($(P23s)-(0.4,0)$){$A_2^*$};
\node at ($(P24s)-(0.4,0)$){$A_3$};
\node at ($(P25s)-(0.4,0)$){$A_4$};
\node at ($(P26s)-(0.4,0)$){$A_5$};
\node at ($(P27s)-(0.4,0)$){$A_6$};
\node at ($(P28s)-(0.4,0)$){$A_7$};
\node at ($(P16s)-(0.4,0)$){$A_5'$};
\node at ($(P32s)-(0.4,0)$){$A_0^+$};
\node at ($(P33s)-(0.4,0)$){$A_1^+$};
\node at ($(P34s)-(0.4,0)$){$A_2^+$};
\node at ($(P35s)-(0.4,0)$){$D_4$};
\node at ($(P36s)-(0.4,0)$){$D_5$};
\node at ($(P37s)-(0.4,0)$){$D_6$};
\node at ($(P47s)-(0.4,0)$){$E_6$};
\node at ($(P48s)-(0.4,0)$){$E_7$};
%
%
%
%
\draw [->, thick] (P21s)--(P22e);
\draw [->, thick] (P22s)--(P23e);
\draw [->, thick] (P23s)--(P24e);
\draw [->, thick] (P24s)--(P25e);
\draw [->, thick] (P25s)--(P26e);
\draw [->, thick] (P26s)--(P27e);
\draw [->, thick] (P27s)--(P28e);
\draw [->, thick] (P32s)--(P33e);
\draw [->, thick] (P33s)--(P34e);
\draw [->, thick] (P34s)--(P35e);
\draw [->, thick] (P35s)--(P36e);
\draw [->, thick] (P36s)--(P37e);
\draw [->, thick] (P47s)--(P48e);
\draw [->, thick] ($(P25s)+(0,0.2)$)--($(P16e)-(0,0.2)$); 
\draw [->, thick] ($(P10s)-(0,0.2)$)--($(P21e)+(0,0.2)$); 
\draw [->, thick] ($(P21s)-(0,0.2)$)--($(P32e)+(0,0.2)$); 
\draw [->, thick] ($(P22s)-(0,0.2)$)--($(P33e)+(0,0.2)$); 
\draw [->, thick] ($(P23s)-(0,0.2)$)--($(P34e)+(0,0.2)$); 
\draw [->, thick] ($(P24s)-(0,0.2)$)--($(P35e)+(0,0.2)$); 
\draw [->, thick] ($(P25s)-(0,0.2)$)--($(P36e)+(0,0.2)$); 
\draw [->, thick] ($(P26s)-(0,0.2)$)--($(P37e)+(0,0.2)$); 
\draw [->, thick] ($(P36s)-(0,0.2)$)--($(P47e)+(0,0.2)$); 
\draw [->, thick] ($(P37s)-(0,0.2)$)--($(P48e)+(0,0.2)$); 
\draw [->, thick] ($(P26s)-(0.1,0.35)$)--($(P47e)+(0.1,0.35)$);
\draw [->, thick] ($(P16s)-(0.1,0.35)$)--($(P37e)+(0.1,0.35)$);
\draw [->, thick] ($(P27s)-(0.1,0.35)$)--($(P48e)+(0.1,0.35)$);
\end{scope}
\end{tikzpicture}
\caption{The poset of hypergeometric functions generalizing ${}_2F_1$.}
\label{fig:class2F1}
\end{figure}
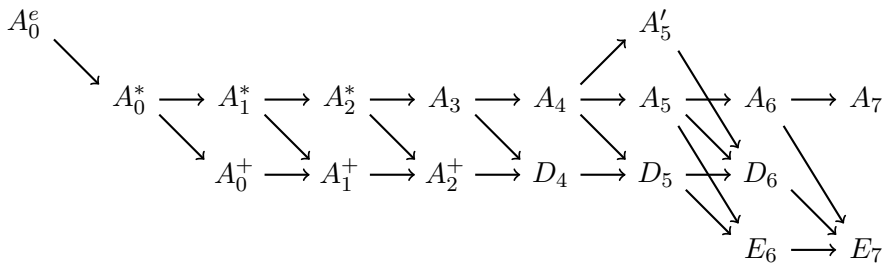

\begin{figure}
  \centering
  \begin{tikzpicture}[scale = 1]
\begin{scope}
\coordinate (P11s) at (0,0);
\coordinate (P11e) at ($(P11s)-(0.8,0)$);
\coordinate (P10s) at ($(P11e)-(0.6,0)$);
\coordinate (P12e) at ($(P11s)+(0.6,0)$);
\coordinate (P12s) at ($(P12e)+(0.8,0)$);
\coordinate (P13e) at ($(P12s)+(0.6,0)$);
\coordinate (P13s) at ($(P13e)+(0.8,0)$);
\coordinate (P14e) at ($(P13s)+(0.6,0)$);
\coordinate (P14s) at ($(P14e)+(0.8,0)$);
\coordinate (P15e) at ($(P14s)+(0.6,0)$);
\coordinate (P15s) at ($(P15e)+(0.8,0)$);
\coordinate (P16e) at ($(P15s)+(0.6,0)$);
\coordinate (P16s) at ($(P16e)+(0.8,0)$);
\coordinate (P17e) at ($(P16s)+(0.6,0)$);
\coordinate (P17s) at ($(P17e)+(0.8,0)$);
\coordinate (P18e) at ($(P17s)+(0.6,0)$);
\coordinate (P18s) at ($(P18e)+(0.8,0)$);
\coordinate (P19e) at ($(P18s)+(0.6,0)$);
\coordinate (P19s) at ($(P19e)+(0.8,0)$);
\coordinate (P21s) at (0,-1);
\coordinate (P21e) at ($(P21s)-(0.8,0)$);
\coordinate (P20s) at ($(P21e)-(0.6,0)$);
\coordinate (P22e) at ($(P21s)+(0.6,0)$);
\coordinate (P22s) at ($(P22e)+(0.8,0)$);
\coordinate (P23e) at ($(P22s)+(0.6,0)$);
\coordinate (P23s) at ($(P23e)+(0.8,0)$);
\coordinate (P24e) at ($(P23s)+(0.6,0)$);
\coordinate (P24s) at ($(P24e)+(0.8,0)$);
\coordinate (P25e) at ($(P24s)+(0.6,0)$);
\coordinate (P25s) at ($(P25e)+(0.8,0)$);
\coordinate (P26e) at ($(P25s)+(0.6,0)$);
\coordinate (P26s) at ($(P26e)+(0.8,0)$);
\coordinate (P27e) at ($(P26s)+(0.6,0)$);
\coordinate (P27s) at ($(P27e)+(0.8,0)$);
\coordinate (P28e) at ($(P27s)+(0.6,0)$);
\coordinate (P28s) at ($(P28e)+(0.8,0)$);
\coordinate (P29e) at ($(P28s)+(0.6,0)$);
\coordinate (P29s) at ($(P29e)+(0.8,0)$);
\coordinate (P31s) at (0,-2);
\coordinate (P32e) at ($(P31s)+(0.6,0)$);
\coordinate (P32s) at ($(P32e)+(0.8,0)$);
\coordinate (P33e) at ($(P32s)+(0.6,0)$);
\coordinate (P33s) at ($(P33e)+(0.8,0)$);
\coordinate (P34e) at ($(P33s)+(0.6,0)$);
\coordinate (P34s) at ($(P34e)+(0.8,0)$);
\coordinate (P35e) at ($(P34s)+(0.6,0)$);
\coordinate (P35s) at ($(P35e)+(0.8,0)$);
\coordinate (P36e) at ($(P35s)+(0.6,0)$);
\coordinate (P36s) at ($(P36e)+(0.8,0)$);
\coordinate (P37e) at ($(P36s)+(0.6,0)$);
\coordinate (P37s) at ($(P37e)+(0.8,0)$);
\coordinate (P38e) at ($(P37s)+(0.6,0)$);
\coordinate (P38s) at ($(P38e)+(0.8,0)$);
\coordinate (P39e) at ($(P38s)+(0.6,0)$);
\coordinate (P39s) at ($(P39e)+(0.8,0)$);
\coordinate (P41s) at (0,-3);
\coordinate (P42e) at ($(P41s)+(0.6,0)$);
\coordinate (P42s) at ($(P42e)+(0.8,0)$);
\coordinate (P43e) at ($(P42s)+(0.6,0)$);
\coordinate (P43s) at ($(P43e)+(0.8,0)$);
\coordinate (P44e) at ($(P43s)+(0.6,0)$);
\coordinate (P44s) at ($(P44e)+(0.8,0)$);
\coordinate (P45e) at ($(P44s)+(0.6,0)$);
\coordinate (P45s) at ($(P45e)+(0.8,0)$);
\coordinate (P46e) at ($(P45s)+(0.6,0)$);
\coordinate (P46s) at ($(P46e)+(0.8,0)$);
\coordinate (P47e) at ($(P46s)+(0.6,0)$);
\coordinate (P47s) at ($(P47e)+(0.8,0)$);
\coordinate (P48e) at ($(P47s)+(0.6,0)$);
\coordinate (P48s) at ($(P48e)+(0.8,0)$);
\coordinate (P49e) at ($(P48s)+(0.6,0)$);
\coordinate (P49s) at ($(P49e)+(0.8,0)$);
\node at ($(P10s)-(0.4,0)$){$1$};
\node at ($(P21s)-(0.4,0)$){$1$};
\node at ($(P22s)-(0.4,0)$){$2$};
\node at ($(P23s)-(0.4,0)$){$2$};
\node at ($(P24s)-(0.4,0)$){$4$};
\node at ($(P25s)-(0.4,0)$){$2$};
\node at ($(P26s)-(0.4,0)$){$4$};
\node at ($(P27s)-(0.4,0)$){$2$};
\node at ($(P28s)-(0.4,0)$){$1$};
\node at ($(P16s)-(0.4,0)$){$2$};
\node at ($(P32s)-(0.4,0)$){$1$};
\node at ($(P33s)-(0.4,0)$){$2$};
\node at ($(P34s)-(0.4,0)$){$2$};
\node at ($(P35s)-(0.4,0)$){$3$};
\node at ($(P36s)-(0.4,0)$){$4$};
\node at ($(P37s)-(0.4,0)$){$3$};
\node at ($(P47s)-(0.4,0)$){$2$};
\node at ($(P48s)-(0.4,0)$){$1$};
%
%
%
%
\draw [->, thick] (P21s)--(P22e);
\draw [->, thick] (P22s)--(P23e);
\draw [->, thick] (P23s)--(P24e);
\draw [->, thick] (P24s)--(P25e);
\draw [->, thick] (P25s)--(P26e);
\draw [->, thick] (P26s)--(P27e);
\draw [->, thick] (P27s)--(P28e);
\draw [->, thick] (P32s)--(P33e);
\draw [->, thick] (P33s)--(P34e);
\draw [->, thick] (P34s)--(P35e);
\draw [->, thick] (P35s)--(P36e);
\draw [->, thick] (P36s)--(P37e);
\draw [->, thick] (P47s)--(P48e);
\draw [->, thick] ($(P25s)+(0,0.2)$)--($(P16e)-(0,0.2)$); 
\draw [->, thick] ($(P10s)-(0,0.2)$)--($(P21e)+(0,0.2)$); 
\draw [->, thick] ($(P21s)-(0,0.2)$)--($(P32e)+(0,0.2)$); 
\draw [->, thick] ($(P22s)-(0,0.2)$)--($(P33e)+(0,0.2)$); 
\draw [->, thick] ($(P23s)-(0,0.2)$)--($(P34e)+(0,0.2)$); 
\draw [->, thick] ($(P24s)-(0,0.2)$)--($(P35e)+(0,0.2)$); 
\draw [->, thick] ($(P25s)-(0,0.2)$)--($(P36e)+(0,0.2)$); 
\draw [->, thick] ($(P26s)-(0,0.2)$)--($(P37e)+(0,0.2)$); 
\draw [->, thick] ($(P36s)-(0,0.2)$)--($(P47e)+(0,0.2)$); 
\draw [->, thick] ($(P37s)-(0,0.2)$)--($(P48e)+(0,0.2)$); 
\draw [->, thick] ($(P26s)-(0.1,0.35)$)--($(P47e)+(0.1,0.35)$);
\draw [->, thick] ($(P16s)-(0.1,0.35)$)--($(P37e)+(0.1,0.35)$);
\draw [->, thick] ($(P27s)-(0.1,0.35)$)--($(P48e)+(0.1,0.35)$);
\end{scope}
\end{tikzpicture}
  \caption{The number of types of equation for each type of hypergeometric
    function}
\label{fig:num2F1}
\end{figure}
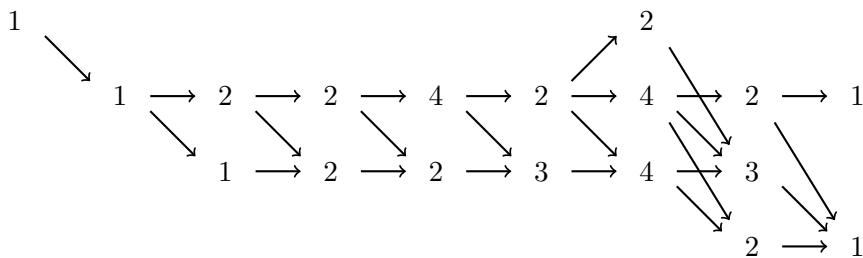

The $E_7$ case is Airy, while the other maximally degenerate case is the
$A_7$ case, which corresponds to the $q$-difference equation
\[
v(qz) = \begin{pmatrix}
  0& 1\\
  -1/z & 0\end{pmatrix}
  v(z),
\]
with solution $(\theta_{q^2}(z),\theta_{q^2}(qz))$.  (Note that although
the corresponding $q^2$-difference equation is a direct sum of first-order
equations, this $q$-difference is actually irreducible!)  The $D_4$ case is
${}_2F_1$.

To see that this diagram corresponds to hypergeometric {\em functions}, as
opposed to hypergeometric {\em equations}, note that the presence of a
$-2$-curve disjoint from the anticanonical curve lets us take the Lax pair
for the corresponding discrete Painlev\'e equation to be reducible.  When
the $-2$-curve is $s-f$, this precisely corresponds to the construction of
the hypergeometric solution via semiclassical biorthogonal functions (or
orthogonal polynomials).  Although our choice of $-2$-curve is $f-e_7-e_8$,
any reflection in a simple root not orthogonal to that $-2$-curve is
generically admissible, and thus we may transform it to $s-f$ instead.

Note that although this diagram tells us that there {\em are} limiting
relations (or, more precisely, {\em will} tell us once we have established
the noncommutative theory!), it does not tell us what that limit actually
looks like.  For instance, consider the equation given by
\[
v(qz)
=
\begin{pmatrix}
  0&1\\-1/qz&1/qz
\end{pmatrix}
v(z),
\]
satisfied by Ismail's $q$-Airy function \cite{IsmailMEH:2005}.  This is one
of the two types of equation of type $A_6$ (the other being
gauge-equivalent to the corresponding $q^2$-difference equation; in
addition to a solution coming from Ismail's $q$-Airy function, this
difference equation also has a solution corresponding to the $q$-Airy
function of \cite{KajiwaraK/MasudaT/NoumiM/OhtaY/YamadaY:2004}), and thus
we see from the above diagram that it degenerates to type $E_7$, a.k.a. the
Airy equation.  But this does not tell us how to find the actual limiting
relation between equations!  If one takes the natural guess that should be
a corresponding limiting relation between solutions, then one can
reconstruct the appropriate scaling limit by looking at the behavior of the
largest zeros for $q$ near 1.  One finds in particular (by numerical
experiment!) that these are close to $1/4$, differing at a scale of
$(1-q)^{2/3}$, suggesting that we should take a scaling limit of the form
$q=1-\epsilon^3$, $z=(1-\epsilon^2 x)/4$.  This limit makes the shift
matrix have eigenvalue $2$ with multiplicity $2$, so we need to do a scalar
gauge transformation to make the limit 1.  (We could avoid the scalar gauge
issue by introducing apparent singularities, but then would need to worry
about how those behave in the limit!)

Thus let $f(z)$ be any rational function with $f(1/4)=1/2$, and gauge by a
solution of the equation $c(qz)=f(z)c(z)$, so that we get a new equation
\[
\hat{v}(qz) = f(z)\begin{pmatrix}0&1\\-1/qz&1/qz\end{pmatrix} \hat{v}(z).
\]
Then defining
\[
w(x) := \begin{pmatrix} -2\epsilon & 0\\2 & -1\end{pmatrix}
  \hat{v}((1-\epsilon^2 x)/4),
\]
we see that
\[
w(x+\epsilon - \epsilon^3 x)
=
\bigl(1+\epsilon \begin{pmatrix} 0&1\\x&0\end{pmatrix} + O(\epsilon^2)\bigr)
w(x)
\]
which in the limit $\epsilon\to 0$ gives us the equation
\[
u'(x) = \begin{pmatrix} 0&1\\x&0\end{pmatrix} u(x)
\]
as desired.

It would be nice to have a more systematic approach to determining such
limit relations, not only for applications to explicit equations, but also
in the hopes that it might give another approach to the problem of
determining when such relations exist!

Thus each arrow in Figure \ref{fig:class2F1} (or, more generally, each
inclusion relation in the poset) gives rise to a nontrivial question of
finding {\em explicit} limiting relations between equations of the given
type.  Of course, once one has such a relation, one then has a purely
analytic question: can one find solutions of the equation for which the
corresponding limit exists (and satisfies the limiting equation)?

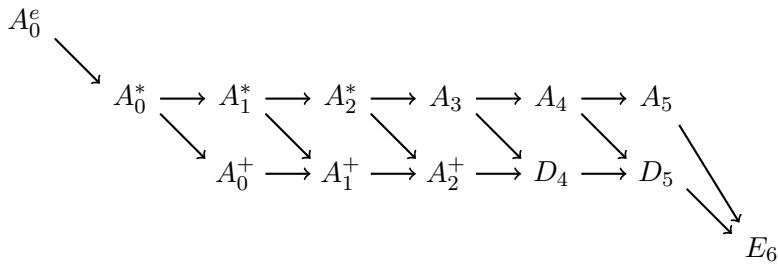
\begin{figure}
  \centering
  \begin{tikzpicture}[scale = 1]
\begin{scope}
\coordinate (P11s) at (0,0);
\coordinate (P11e) at ($(P11s)-(0.8,0)$);
\coordinate (P10s) at ($(P11e)-(0.6,0)$);
\coordinate (P12e) at ($(P11s)+(0.6,0)$);
\coordinate (P12s) at ($(P12e)+(0.8,0)$);
\coordinate (P13e) at ($(P12s)+(0.6,0)$);
\coordinate (P13s) at ($(P13e)+(0.8,0)$);
\coordinate (P14e) at ($(P13s)+(0.6,0)$);
\coordinate (P14s) at ($(P14e)+(0.8,0)$);
\coordinate (P15e) at ($(P14s)+(0.6,0)$);
\coordinate (P15s) at ($(P15e)+(0.8,0)$);
\coordinate (P16e) at ($(P15s)+(0.6,0)$);
\coordinate (P16s) at ($(P16e)+(0.8,0)$);
\coordinate (P17e) at ($(P16s)+(0.6,0)$);
\coordinate (P17s) at ($(P17e)+(0.8,0)$);
\coordinate (P18e) at ($(P17s)+(0.6,0)$);
\coordinate (P18s) at ($(P18e)+(0.8,0)$);
\coordinate (P19e) at ($(P18s)+(0.6,0)$);
\coordinate (P19s) at ($(P19e)+(0.8,0)$);
\coordinate (P21s) at (0,-1);
\coordinate (P21e) at ($(P21s)-(0.8,0)$);
\coordinate (P20s) at ($(P21e)-(0.6,0)$);
\coordinate (P22e) at ($(P21s)+(0.6,0)$);
\coordinate (P22s) at ($(P22e)+(0.8,0)$);
\coordinate (P23e) at ($(P22s)+(0.6,0)$);
\coordinate (P23s) at ($(P23e)+(0.8,0)$);
\coordinate (P24e) at ($(P23s)+(0.6,0)$);
\coordinate (P24s) at ($(P24e)+(0.8,0)$);
\coordinate (P25e) at ($(P24s)+(0.6,0)$);
\coordinate (P25s) at ($(P25e)+(0.8,0)$);
\coordinate (P26e) at ($(P25s)+(0.6,0)$);
\coordinate (P26s) at ($(P26e)+(0.8,0)$);
\coordinate (P27e) at ($(P26s)+(0.6,0)$);
\coordinate (P27s) at ($(P27e)+(0.8,0)$);
\coordinate (P28e) at ($(P27s)+(0.6,0)$);
\coordinate (P28s) at ($(P28e)+(0.8,0)$);
\coordinate (P29e) at ($(P28s)+(0.6,0)$);
\coordinate (P29s) at ($(P29e)+(0.8,0)$);
\coordinate (P31s) at (0,-2);
\coordinate (P32e) at ($(P31s)+(0.6,0)$);
\coordinate (P32s) at ($(P32e)+(0.8,0)$);
\coordinate (P33e) at ($(P32s)+(0.6,0)$);
\coordinate (P33s) at ($(P33e)+(0.8,0)$);
\coordinate (P34e) at ($(P33s)+(0.6,0)$);
\coordinate (P34s) at ($(P34e)+(0.8,0)$);
\coordinate (P35e) at ($(P34s)+(0.6,0)$);
\coordinate (P35s) at ($(P35e)+(0.8,0)$);
\coordinate (P36e) at ($(P35s)+(0.6,0)$);
\coordinate (P36s) at ($(P36e)+(0.8,0)$);
\coordinate (P37e) at ($(P36s)+(0.6,0)$);
\coordinate (P37s) at ($(P37e)+(0.8,0)$);
\coordinate (P38e) at ($(P37s)+(0.6,0)$);
\coordinate (P38s) at ($(P38e)+(0.8,0)$);
\coordinate (P39e) at ($(P38s)+(0.6,0)$);
\coordinate (P39s) at ($(P39e)+(0.8,0)$);
\coordinate (P41s) at (0,-3);
\coordinate (P42e) at ($(P41s)+(0.6,0)$);
\coordinate (P42s) at ($(P42e)+(0.8,0)$);
\coordinate (P43e) at ($(P42s)+(0.6,0)$);
\coordinate (P43s) at ($(P43e)+(0.8,0)$);
\coordinate (P44e) at ($(P43s)+(0.6,0)$);
\coordinate (P44s) at ($(P44e)+(0.8,0)$);
\coordinate (P45e) at ($(P44s)+(0.6,0)$);
\coordinate (P45s) at ($(P45e)+(0.8,0)$);
\coordinate (P46e) at ($(P45s)+(0.6,0)$);
\coordinate (P46s) at ($(P46e)+(0.8,0)$);
\coordinate (P47e) at ($(P46s)+(0.6,0)$);
\coordinate (P47s) at ($(P47e)+(0.8,0)$);
\coordinate (P48e) at ($(P47s)+(0.6,0)$);
\coordinate (P48s) at ($(P48e)+(0.8,0)$);
\coordinate (P49e) at ($(P48s)+(0.6,0)$);
\coordinate (P49s) at ($(P49e)+(0.8,0)$);
\node at ($(P10s)-(0.4,0)$){$A_0^e$};
\node at ($(P21s)-(0.4,0)$){$A_0^*$};
\node at ($(P22s)-(0.4,0)$){$A_1^*$};
\node at ($(P23s)-(0.4,0)$){$A_2^*$};
\node at ($(P24s)-(0.4,0)$){$A_3$};
\node at ($(P25s)-(0.4,0)$){$A_4$};
\node at ($(P26s)-(0.4,0)$){$A_5$};
\node at ($(P32s)-(0.4,0)$){$A_0^+$};
\node at ($(P33s)-(0.4,0)$){$A_1^+$};
\node at ($(P34s)-(0.4,0)$){$A_2^+$};
\node at ($(P35s)-(0.4,0)$){$D_4$};
\node at ($(P36s)-(0.4,0)$){$D_5$};
\node at ($(P47s)-(0.4,0)$){$E_6$};
%
%
%
%
\draw [->, thick] (P21s)--(P22e);
\draw [->, thick] (P22s)--(P23e);
\draw [->, thick] (P23s)--(P24e);
\draw [->, thick] (P24s)--(P25e);
\draw [->, thick] (P25s)--(P26e);
\draw [->, thick] (P32s)--(P33e);
\draw [->, thick] (P33s)--(P34e);
\draw [->, thick] (P34s)--(P35e);
\draw [->, thick] (P35s)--(P36e);
\draw [->, thick] ($(P10s)-(0,0.2)$)--($(P21e)+(0,0.2)$); 
\draw [->, thick] ($(P21s)-(0,0.2)$)--($(P32e)+(0,0.2)$); 
\draw [->, thick] ($(P22s)-(0,0.2)$)--($(P33e)+(0,0.2)$); 
\draw [->, thick] ($(P23s)-(0,0.2)$)--($(P34e)+(0,0.2)$); 
\draw [->, thick] ($(P24s)-(0,0.2)$)--($(P35e)+(0,0.2)$); 
\draw [->, thick] ($(P25s)-(0,0.2)$)--($(P36e)+(0,0.2)$); 
\draw [->, thick] ($(P36s)-(0,0.2)$)--($(P47e)+(0,0.2)$); 
\draw [->, thick] ($(P26s)-(0.1,0.35)$)--($(P47e)+(0.1,0.35)$);
\end{scope}
\end{tikzpicture}
  \caption{The poset of hypergeometric sums of type ${}_2F_1$}
  \label{fig:Askey1}
\end{figure}

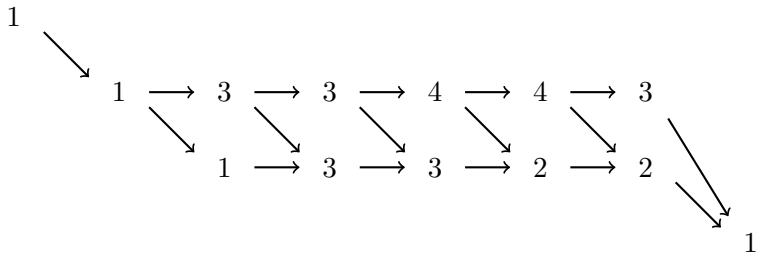
\begin{figure}
  \centering
  \begin{tikzpicture}[scale = 1]
\begin{scope}
\coordinate (P11s) at (0,0);
\coordinate (P11e) at ($(P11s)-(0.8,0)$);
\coordinate (P10s) at ($(P11e)-(0.6,0)$);
\coordinate (P12e) at ($(P11s)+(0.6,0)$);
\coordinate (P12s) at ($(P12e)+(0.8,0)$);
\coordinate (P13e) at ($(P12s)+(0.6,0)$);
\coordinate (P13s) at ($(P13e)+(0.8,0)$);
\coordinate (P14e) at ($(P13s)+(0.6,0)$);
\coordinate (P14s) at ($(P14e)+(0.8,0)$);
\coordinate (P15e) at ($(P14s)+(0.6,0)$);
\coordinate (P15s) at ($(P15e)+(0.8,0)$);
\coordinate (P16e) at ($(P15s)+(0.6,0)$);
\coordinate (P16s) at ($(P16e)+(0.8,0)$);
\coordinate (P17e) at ($(P16s)+(0.6,0)$);
\coordinate (P17s) at ($(P17e)+(0.8,0)$);
\coordinate (P18e) at ($(P17s)+(0.6,0)$);
\coordinate (P18s) at ($(P18e)+(0.8,0)$);
\coordinate (P19e) at ($(P18s)+(0.6,0)$);
\coordinate (P19s) at ($(P19e)+(0.8,0)$);
\coordinate (P21s) at (0,-1);
\coordinate (P21e) at ($(P21s)-(0.8,0)$);
\coordinate (P20s) at ($(P21e)-(0.6,0)$);
\coordinate (P22e) at ($(P21s)+(0.6,0)$);
\coordinate (P22s) at ($(P22e)+(0.8,0)$);
\coordinate (P23e) at ($(P22s)+(0.6,0)$);
\coordinate (P23s) at ($(P23e)+(0.8,0)$);
\coordinate (P24e) at ($(P23s)+(0.6,0)$);
\coordinate (P24s) at ($(P24e)+(0.8,0)$);
\coordinate (P25e) at ($(P24s)+(0.6,0)$);
\coordinate (P25s) at ($(P25e)+(0.8,0)$);
\coordinate (P26e) at ($(P25s)+(0.6,0)$);
\coordinate (P26s) at ($(P26e)+(0.8,0)$);
\coordinate (P27e) at ($(P26s)+(0.6,0)$);
\coordinate (P27s) at ($(P27e)+(0.8,0)$);
\coordinate (P28e) at ($(P27s)+(0.6,0)$);
\coordinate (P28s) at ($(P28e)+(0.8,0)$);
\coordinate (P29e) at ($(P28s)+(0.6,0)$);
\coordinate (P29s) at ($(P29e)+(0.8,0)$);
\coordinate (P31s) at (0,-2);
\coordinate (P32e) at ($(P31s)+(0.6,0)$);
\coordinate (P32s) at ($(P32e)+(0.8,0)$);
\coordinate (P33e) at ($(P32s)+(0.6,0)$);
\coordinate (P33s) at ($(P33e)+(0.8,0)$);
\coordinate (P34e) at ($(P33s)+(0.6,0)$);
\coordinate (P34s) at ($(P34e)+(0.8,0)$);
\coordinate (P35e) at ($(P34s)+(0.6,0)$);
\coordinate (P35s) at ($(P35e)+(0.8,0)$);
\coordinate (P36e) at ($(P35s)+(0.6,0)$);
\coordinate (P36s) at ($(P36e)+(0.8,0)$);
\coordinate (P37e) at ($(P36s)+(0.6,0)$);
\coordinate (P37s) at ($(P37e)+(0.8,0)$);
\coordinate (P38e) at ($(P37s)+(0.6,0)$);
\coordinate (P38s) at ($(P38e)+(0.8,0)$);
\coordinate (P39e) at ($(P38s)+(0.6,0)$);
\coordinate (P39s) at ($(P39e)+(0.8,0)$);
\coordinate (P41s) at (0,-3);
\coordinate (P42e) at ($(P41s)+(0.6,0)$);
\coordinate (P42s) at ($(P42e)+(0.8,0)$);
\coordinate (P43e) at ($(P42s)+(0.6,0)$);
\coordinate (P43s) at ($(P43e)+(0.8,0)$);
\coordinate (P44e) at ($(P43s)+(0.6,0)$);
\coordinate (P44s) at ($(P44e)+(0.8,0)$);
\coordinate (P45e) at ($(P44s)+(0.6,0)$);
\coordinate (P45s) at ($(P45e)+(0.8,0)$);
\coordinate (P46e) at ($(P45s)+(0.6,0)$);
\coordinate (P46s) at ($(P46e)+(0.8,0)$);
\coordinate (P47e) at ($(P46s)+(0.6,0)$);
\coordinate (P47s) at ($(P47e)+(0.8,0)$);
\coordinate (P48e) at ($(P47s)+(0.6,0)$);
\coordinate (P48s) at ($(P48e)+(0.8,0)$);
\coordinate (P49e) at ($(P48s)+(0.6,0)$);
\coordinate (P49s) at ($(P49e)+(0.8,0)$);
\node at ($(P10s)-(0.4,0)$){$1$};
\node at ($(P21s)-(0.4,0)$){$1$};
\node at ($(P22s)-(0.4,0)$){$3$};
\node at ($(P23s)-(0.4,0)$){$3$};
\node at ($(P24s)-(0.4,0)$){$4$};
\node at ($(P25s)-(0.4,0)$){$4$};
\node at ($(P26s)-(0.4,0)$){$3$};
\node at ($(P32s)-(0.4,0)$){$1$};
\node at ($(P33s)-(0.4,0)$){$3$};
\node at ($(P34s)-(0.4,0)$){$3$};
\node at ($(P35s)-(0.4,0)$){$2$};
\node at ($(P36s)-(0.4,0)$){$2$};
\node at ($(P47s)-(0.4,0)$){$1$};
%
%
%
%
\draw [->, thick] (P21s)--(P22e);
\draw [->, thick] (P22s)--(P23e);
\draw [->, thick] (P23s)--(P24e);
\draw [->, thick] (P24s)--(P25e);
\draw [->, thick] (P25s)--(P26e);
\draw [->, thick] (P32s)--(P33e);
\draw [->, thick] (P33s)--(P34e);
\draw [->, thick] (P34s)--(P35e);
\draw [->, thick] (P35s)--(P36e);
\draw [->, thick] ($(P10s)-(0,0.2)$)--($(P21e)+(0,0.2)$); 
\draw [->, thick] ($(P21s)-(0,0.2)$)--($(P32e)+(0,0.2)$); 
\draw [->, thick] ($(P22s)-(0,0.2)$)--($(P33e)+(0,0.2)$); 
\draw [->, thick] ($(P23s)-(0,0.2)$)--($(P34e)+(0,0.2)$); 
\draw [->, thick] ($(P24s)-(0,0.2)$)--($(P35e)+(0,0.2)$); 
\draw [->, thick] ($(P25s)-(0,0.2)$)--($(P36e)+(0,0.2)$); 
\draw [->, thick] ($(P36s)-(0,0.2)$)--($(P47e)+(0,0.2)$); 
\draw [->, thick] ($(P26s)-(0.1,0.35)$)--($(P47e)+(0.1,0.35)$);
\end{scope}
\end{tikzpicture}
  \caption{The numbers of subtypes corresponding to classical
    (bi)orthogonal functions}
  \label{fig:Askey2}
\end{figure}

Just as we can view Figure \ref{fig:class2F1} as classifying semiclassical
biorthogonal functions, there is an analogous classification for {\em
  classical} biorthogonal functions, now corresponding to hypergeometric
equations in the ${}_2F_1$ hierarchy for which the anticanonical components
are orthogonal to $s-f$.  Again we can blow up two additional points to
make $f-e_7-e_8$ effective and disjoint from the anticanonical curve, and
thus project to the classification of deformed elliptic surfaces with an
$A_2$ singularity disjoint from the anticanonical curve, corresponding to
the diagram of subsystems of $E_6$, Figure \ref{fig:Askey1}.  These may
also be interpreted as those hypergeometric functions which are
representable by finite sums.  Figure \ref{fig:Askey2} gives the further
decomposition into classical (bi)orthogonal functions.  (This should, in
fact, be decomposed further, as the treatment via reducible Lax pairs does
not actually fix the parameters corresponding to the poles of the two
orthogonal bases.)

Again, these figures are only giving an algebraic classification; in
particular, the subdiagram corresponding to the Askey scheme has a further
decomposition into different ways the norms can be positive.

\section{More classifications}

Similarly, in the case $2s+2f-e_1-e_2-e_3-e_4-e_5-e_6-e_7-e_8$,
corresponding in the generic case to the elliptic Painlev\'e equation,
which we saw above corresponded to $22$ orbits of $W(E_8)$ (as in Figure 1
above), there are a total of $139981$ types (we omit the figure!).  The
$W(E_8)$ orbits correspond to different types of Painlev\'e equation, but
the Lax pairs themselves are essentially classified by the $W(D_8)$ orbits
(since as we have seen, the map between equations and sheaves is only
defined up to scalar gauge equivalence in any event).  There are $61$ such
orbits, implying that a given (discrete or continuous) Painlev\'e equation
can have qualitatively different second-order Lax pairs.

For instance, the classification indicates that the Painlev\'e VI equation
should have a Lax pair in the form of a second order symmetric difference
equation.  Indeed, there is a $W(D_8)$-orbit of types of surfaces in which
the anticanonical curve decomposes as
\[
(2s+f-e_1-e_2-e_3-e_4-e_5-e_6)
+(f-e_5-e_6)
+(e_5-e_6)
+2(e_6-e_7)
+(e_7-e_8),
\]
which by our discussion of singularities above corresponds to a symmetric
difference equation with four finite singularities and an indecomposable
singularity at infinity.  The corresponding Lax pair was constructed in
\cite{OrmerodCM/RainsEM:2017b}, along with degenerations to Painlev\'e V
and III.

This work actually predated the proof below that such equations should have
continuous isomonodromy deformations, which points out an important fact:
for those readers who are interested in {\em explicit} Lax pairs and
relations between them, it suffices to know that the noncommutative theory
{\em works}, without needing to know the details of {\em how} or {\em why}.
Indeed, once one knows (or believes!) that a Lax pair of a given type
exists, actually writing it down only requires understanding the kind of
equation, the nature of its singularities (along with the generally
nontrivial question of finding a {\em nice} rational parametrization of the
corresponding moduli space), and how the pseudo-twists behave, all of which
are just the obvious deformations of their commutative analogues.

The case $r(2s+2f-e_1-\cdots-e_8)+e_8-e_9$, which corresponds to an
elliptic version of the ``matrix Painlev\'e equation''
\cite{KawakamiH:2015}, has precisely the same classification (apart from
blowing up a point of an anticanonical component of self-intersection
$-1$); there are two apparent additional types for $r=2$, but they have a
component $e_8-e_9$, preventing the linear system from being generically
integral.

The case $2s+3f-e_1-\cdots-e_{10}$, which generically corresponds to the
simplest non-Painlev\'e case of the elliptic Garnier equation
\cite{ellGarnier}, gives rise to a total of $6374578$ types in $84$
$W(D_{10})$-orbits.  Note that since the divisor class is $-K+f$, it is not
preserved by reflection in $s-f$; as a result, each degeneration of the
simplest elliptic Garnier equation corresponds to a unique type of
second-order Lax pair up to scalar gauge equivalence.  (There are, of
course, higher order Lax pairs, e.g., the third-order Lax pair obtained by
reflecting in $s-f$.)

This combinatorial explosion suggests that there should be better ways to
perform these classifications.  In particular, for the Garnier cases
\[
2s+(m+1)f-e_1-\cdots-e_{2m+6},
\]
modulo the action of $W(D_{2m+6})$, we can greatly speed up the calculation
using the classification of singularities.  Indeed, there are relatively
few types of irreducible singularities of rank $1$ or $2$, each of which
requires a certain number of blowups, and thus one may ask for all ways of
combining them to account for a total of $2m+6$ blowups.  This enables the
calculation for $m=3$ as well; combined with the results for $m=\le 2$,
this suggests the following patterns: in addition to the one symmetric
elliptic difference equation, there are $2m+8$ symmetric $q$-difference
equations, $2m+6$ symmetric difference equations, $m^2+8m+14$ nonsymmetric
$q$-difference equations, $2m+7$ nonsymmetric difference equations, while
the number of types of differential equation is given by the number of
partitions of $m+3$ into parts of two colors, with one color of $1$
forbidden.

Note that the $m=1$ Garnier case together with the first nontrivial matrix
Painlev\'e case together account for all cases with $4$-dimensional moduli
space.  Types of differential equation with that property were classified
in \cite{KawakamiH/NakamuraA/SakaiH:2013}, giving the same number of types
as the additive cases of Garnier and matrix Painlev\'e (and the same
predicted number of continuous isomonodromy deformations).  This comes from
the fact that any ordinary difference case can be transformed to a
differential equation using the Weyl group symmetry (possibly with some
additional blowups required), though of course this in general increases
the {\em order} of the equation.  Indeed, given a symmetric difference
equation (in $F_2$, say), doing an elementary transformation at the cusp
and some finite point gives a configuration in $F_0$ for which the
anticanonical curve decomposes as $(2s+f)+(f)$, and thus the ``Fourier''
transform turns this into $(s+2f)+(s)$, corresponding to a nonsymmetric
equation.  Further elementary transformations at finite points can then
give us $(s+f)+(s+f)$, and then doing two elementary transformations at the
point of tangency gives us $(s)+(s)+2(f)$, which transforms to a
differential equation as required.  As usual, of course, one issue with
doing this is that the minimal representation of the isomonodromy problem
with differential Lax pair can be higher-order than the minimal
representation without that constraint, and thus each possibility from the
general theory splits further in the differential theory.  Another issue is
that since different cases in the same poset can be represented by
equations of different orders, the limiting relations are that much harder
to understand.

\smallskip

In the $m=1$ Garnier case, we were able to use the theory of elliptic
surfaces to directly compute the poset of types.  This, of course, breaks
for higher-order Garnier cases, but there is another geometric
interpretation that at least allows one to represent the closures of types
via explicit ideals.  (Actually computing those ideals may require
additional ideas, however!)

The basic idea is that since we are fixing a ruling in the form of a
divisor class $f$, the surface comes with a map $\rho:X\to \P^1$.
Moreover, we have a natural divisor of relative degree 2 (the anticanonical
divisor!)  which induces a map to the $\P^2$ bundle corresponding to the
direct image.  Using our algorithm for computing Betti numbers, we see that
\[
\rho_*\Gamma(\sO_X(C_\alpha))
\cong
\sO_{\P^1}(0,1-m,-m-1),
\]
where we use $\sO_{\P^1}(d_1,d_2,\dots)$ to denote the corresponding direct
sum of line bundles.  Now, the image of $X$ in this $\P^2$ bundle is a
conic bundle, and thus is cut out by an equation of relative degree 2.  We
can determine this by noting
\[
\rho_*\Gamma(\sO_X(2C_\alpha))
\cong
\sO_{\P^1}(0,1-m,-1-m,2-2m,-2m)
\]
This differs from the symmetric power by $\sO_{\P^1}(-2m-2)$, which tells
us that the equation must have this degree (as the map from the symmetric
power must be surjective).  If the leading coefficient were zero, then
$\sO_{\P^1}(-2m-2)$ would survive as a summand, and thus the coefficient is
a unit, allowing us to complete the square in characteristic not 2.  This
gives an equation of the form
\[
u^2 = p_4(x,w)y^2 + p_{m+3}(x,w)y + p_{2m+2}(x,w)
\]
for $m\ge 2$, and similarly for $m=0$ and $1$.  We thus find in general
that $X$ is the minimal desingularization of a double cover of $F_{m-1}$
(with a marked section of self-intersection $1-m$ for $m\le 1$), ramified
along a suitably embedded hyperelliptic curve of genus $m+2$.
(Generically, this is a smooth hyperelliptic curve embedded in the
Hirzebruch surface corresponding to the direct image of a uniquely
effective line bundle of degree 4, agreeing with the dimension $2m+7$ of
the moduli stack of surfaces on which $2s+(m+1)f-e_1-\cdots-e_{2m+6}$ is
effective.)

Each given type corresponds to having a particular type of singularity
along the minimal (or marked for $m\le 1$) section; since the singularities
of the double cover of a smooth surface come from singularities of the
ramification divisor, we see that we can in principle translate each type
into a corresponding family of curves.  This is fairly straightforward if
we use the symmetries to fix the locations of the singularities, but it is
not clear how to efficiently restore symmetry.  (On the other hand, given a
guess for the poset structure, it suffices to prove that the covering
relations work, and a given guess for how the corresponding limit might
work may very well require less restoration of symmetry.)

\medskip

In general, given a particular Lax pair in the literature, one would like
to know how it fits into the present theory, in particular so as to
determine how it generalizes and specializes.  This is mostly
straightforward given our discussions above, but there are a couple of
subtle points, so we give a more explicit discussion below.  We focus on
the cases corresponding to rational surfaces.

Step 1. If necessary, put the equation into matrix form.  Lax pairs are
typically given in matrix form or in straight-line form.  In the latter
case, one can put it into matrix form by taking a solution vector of the
form $(v(z),v(qz),v(q^2z),\dots)$ or analogue thereof.  This will tend to
introduce apparent singularities, however.  (Similarly, if the equation
comes as a (discrete) connection on a nontrivial vector bundle, making the
bundle trivial will introduce apparent singularities.)

Step 2. Eliminate apparent singularities.  Globally, this is done by
looking for irreducible singularities of the equation that differ by
isomonodromy transformations, and doing the corresponding global
transformations.  The one problem with this is that it makes the ambient
vector bundle nontrivial, but we can always introduce apparent
singularities at a suitable point to eliminate this issue.  The eventual
goal is that any two irreducible singularities away from that point will
either be equal or in different orbits under isomonodromy transformations.

There is a somewhat less ad-hoc version of Step 2: we introduced
singularities at an auxiliary point in order to keep the bundle trivial.
This gives an equation with apparent singularities along that orbit, which
we can kill at the expense of introducing apparent singularities along a
different orbit.  This gives us two equations
\[
v(z+q) = A_1(z) v(z), \qquad v(z+q)=A_2(z) v(z)
\]
such that $A_1$ has no apparent singularities away from the orbit of $u_1$,
$A_2$ has no apparent singularities away from the orbit of $u_2$, and $A_1$
and $A_2$ are related by a gauge transformation which is invertibly
holomorphic on the complement of both orbits.  This tells us that the
``true'' connection is on a vector bundle $V$ that has explicit
trivializations on the two complements related by precisely that gauge
transformation.  This actually determines the vector bundle, and thus this
pair of equations with gauge transformation may be viewed as an explicit
description of the global {\em twisted} equation.

Below the symmetric elliptic case, one can always take one of the orbits to
be $\{\infty\}$, and can then get a description valid on a complement of
finitely many finite orbits by working over the local ring at $\infty$.
(In the nonsymmetric $q$-difference and differential cases, one can work
over the complements of $0$ and $\infty$ and get gluing data of a more
familiar kind.)

Step 3. Forget ``q'' to obtain a relaxed equation.  Compute its horizontal
spectral curve and perform suitable blowups to separate its strict
transform from the anticanonical curve.  Here one should ignore any blowups
at points where we introduced apparent singularities in Step 2.

Note that here we must be over a trivial bundle ($A$ must be an
endomorphism for its characteristic polynomial to make sense), but given a
global description as above, the blowups for the spectral curve of $A_1$
will indeed correspond to the singularities of $A_1$ away from the orbit of
$u_1$, and similarly for $A_2$.  In particular, on the common complement,
$A_1$ and $A_2$ differ by gauging in an invertibly holomorphic matrix, and
thus have the same singularities.

Step 4. Blow down any $-1$-curves of the form $f-e_i$ that are disjoint from the
spectral curve (i.e., permute $e_i$ to $e_1$, perform an elementary
transformation, then permute $e_1$ to $e_m$) to get a minimal
representation as a twisted equation.

Once one has exhibited the equation as a sheaf on a surface of a particular
type, one should sanity-check that the moduli space has the correct
dimension.  This can fail, for instance, if the sheaf has some additional
structure, e.g., if the equation is gauge-equivalent to a symmetric
difference equation but one performed the above steps as if it were a
nonsymmetric equation.  In practice, one can often guess such
symmetries from the induced symmetries of the singularities, when it is not
already apparent in the description of the Lax pair.

If such a hidden symmetry takes the form of the equation having a smaller
structure group than $\GL_n$ (e.g., the group of similitudes for a
quadratic or symplectic form), the above theory mostly does not apply, with
the main exception being that similitude groups correspond to sheaves which
are related to their duals by an anti-Poisson automorphism of the surface.
It would be nice to have a version of the theory that works for more
general structure groups, though of course the Fourier transform causes
problems in that regard, since it already changes the structure group in
the $\GL$ case.  This already applies to cases in which the structure group
is $\GL$ but relative to a different representation; even something as
simple as the exterior square can give an equation with more complicated
singularities, and there is no reason to expect its Fourier transform to
itself be an exterior square.  However, as long as one is willing to fix a
ruling and a particular representation of the structure group, it should be
possible to describe the moduli space of equations with that group as a
subspace of the moduli space of sheaves, and (following
\cite{HurtubiseJC/MarkmanE:2002b}) one can even hope for it to be a Poisson
subspace.

In any event, once one has a model for the given Lax pair as a deformation
of a sheaf on a given surface, one can then classify all surfaces with a
sheaf of that Chern class disjoint from the anticanonical curve, putting
the pair into a hierarchy.  One can also apply the generalized Fourier
transform (using the explicit description in Appendix \ref{chap:fourier},
or simply swapping the rulings in the relaxed model) along with the action
of $W(D_m)$ to find lower order Lax pairs for the same isomonodromy
systems.  Note that this may be worth doing even in cases with unknown
symmetries, and indeed this can sometimes make the symmetries more
apparent.  (We discuss an example in Section \ref{sec:morphisms_noncomm}
below in which a Lax pair for Painlev\'e VI with structure group $\so_8$
(given in \cite{NoumiM/YamadaY:2002}) arises as the Mellin transform of the
nonsymmetric equation obtained by forgetting the symmetry of a symmetric
difference Lax pair, and thus can be explained by the present theory
despite the structure group not being $\GL$.)


\part{Noncommutative geometry}
\chapter{Noncommutative quasi-ruled surfaces}
\label{chap:daha}

\section{An algebra coming from involutions}

Since the relaxed/autonomous/isospectral version of the theory of
isomonodromy problems is related to sheaves on {\em Poisson} surfaces, it
is reasonable to guess (as we have above) that the nonrelaxed version
should correspond to sheaves on {\em noncommutative} surfaces.  (This is
reinforced by the observation that the algebra of differential operators on
a curve is naturally a noncommutative deformation of an $\A^1$-bundle over
the curve.)  In particular, before considering blowups, we should first
consider the noncommutative versions of ruled surfaces themselves, and
indeed the surfaces we study below will all be related to such surfaces via
a sequence of blowups and blowdowns.

The basic construction is due to Van den Bergh \cite{VandenBerghM:2012},
but it will be helpful to give an alternate construction.  In particular,
since we expect that there should be a relation between sheaves on such
surfaces and difference or differential equations, we would like to have a
construction in which that relation is visible.

We should first address what we {\em mean} by a noncommutative scheme.
Much like ``integrable system'', this actually does not have a fully
satisfactory definition.  The usual Zariski-topological ideas of scheme
theory break in the noncommutative setting, as there tend to be too few
points; instead, inspired by the fact that a commutative scheme can be
recovered from its category of coherent sheaves, we actually {\em identify}
noncommutative schemes with their categories of coherent sheaves.  That is,
a noncommutative scheme $X$ is an abelian category $\coh X$ that ``looks
like'' the category of coherent sheaves on a scheme; in particular, it has
a marked object $\sO_X$ corresponding to the structure sheaf.  (With this
in mind, an abelian category with a marked object is sometimes called a
``quasi-scheme''.)  In particular, many of the properties of the category
of coherent sheaves of a smooth projective surface carry over to the
noncommutative surfaces we study (they are Gorenstein, the $\Ext$ groups
between sheaves vanishes outside degrees $0,1,2$, there is an analogue of
Cohen-Macaulay duality, etc.).

Note that just as a commutative projective surface can be constructed as
the Proj of many different graded algebras, its category of sheaves can be
obtained not only from those algebras but from various Morita-equivalent
graded algebras.  For instance, if $X$ is a commutative projective scheme
with ample bundle $\sO_X(1)$, then for any vector bundle $V$ on $X$, the
(noncommutative) graded algebra $\bigoplus_d \Hom(V,V(d))$ has the same
category of graded modules as $\bigoplus_d \Gamma(\sO_X(d))$, and the same
subcategory of ``torsion'' modules (modules in which every element is
killed in sufficiently high degree, possibly depending on the element), and
thus gives rise to the same quotient category (which when $V$ has rank 1 is
just Serre's construction of $\coh X$).  Our first construction of
noncommutative surfaces will be similarly Morita equivalent to Van den
Bergh's construction.  (This will not only make the connection to
difference or differential operators clearer, but will make it easier to
show that certain instances have large centers.)

With this in mind, we begin by considering a special case of that
commutative construction.  (Our surfaces will be a sort of free product of
two instances of this construction.)  Let $Y/S$ be a scheme and $\pi:X\to
Y$ a finite flat morphism of degree 2 such that $\pi_*\sO_X$ is locally
free.  Then for any invertible sheaf ${\cal L}$ on $X$, $\pi_*{\cal L}$ is
a vector bundle of rank 2 on $Y$, and thus we obtain an $\sO_Y$-algebra
\[
{\cal A}_{{\cal L},\pi}:=\sEnd_Y(\pi_*{\cal L}).
\]
This is Morita equivalent to $\sO_Y$, with an explicit equivalence
$\sO_Y\text{-mod}\to {\cal A}_{{\cal L},\pi}\text{-mod}$ given by
\[
M\mapsto M\otimes_{\sO_Y} \pi_*{\cal L},
\]
and the inverse given by
\[
M\mapsto \Hom_{{\cal A}_{{\cal L},\pi}}(\pi_*{\cal L},M).
\]
Moreover, since this algebra contains $\pi_*\sO_X$, we may interpret it as
a $\pi_*\sO_X$-bimodule, which we abbreviate to ``$\sO_X$-bimodule'' (as
justified below).

\begin{lem}
  There are canonical isomorphisms
  ${\cal A}_{{\cal L}\otimes_{\sO_X} \pi^*{\cal L}',\pi}\cong {\cal
    A}_{{\cal L},\pi}$
  for any line bundle ${\cal L}'$ on $Y$, as well as a canonical involution
  ${\cal A}_{{\cal L},\pi}^{\text{op}}\cong {\cal A}_{{\cal L},\pi}$.
\end{lem}

\begin{proof}
  For the first claim, we have
  \[
  \sEnd_Y(\pi_*({\cal L}\otimes_{\sO_X} \pi^*{\cal L}'))
  \cong
  \sEnd_Y(\pi_*({\cal L})\otimes_{\sO_Y}{\cal L}')
  \cong
  \sEnd_Y(\pi_*({\cal L})),
  \]
  while for the second claim, we observe that for any rank 2 vector bundle
  $V$ on $Y$, we have isomorphisms
  \[
  \sEnd_Y(V)^{\text{op}}
  \cong
  \sEnd_Y(\sHom(V,\sO_Y))
  \cong
  \sEnd_Y(V\otimes \det(V)^{-1})
  \cong
  \sEnd_Y(V),
  \]
  which in the case of the trivial bundle is the (adjoint) involution
  \[
  \begin{pmatrix} a&b\\c&d\end{pmatrix}
  \mapsto
  \begin{pmatrix} d&-b\\-c&a\end{pmatrix},
  \]
  so is an involution in general.  
\end{proof}

\begin{rem}
  It is worth noting that although the {\em algebra} ${\cal A}_{{\cal
      L},\pi}$ only depends on ${\cal L}$ modulo $\Pic(Y)$, the explicit
  Morita equivalence depends on ${\cal L}$.
\end{rem}

\begin{prop}
  There is a natural isomorphism $\pi^!\cong
  \pi^*({-}\otimes_{\sO_Y}\det(\pi_*\sO_X)^{-1})$
\end{prop}

\begin{proof}
  This reduces to showing that there is a natural isomorphism
  \[
  \sHom_Y(\pi_*\sO_X,N)\cong \pi_*\sHom_X(\sO_X,\pi^*(N\otimes_{\sO_Y}
  \det(\pi_*\sO_X)^{-1})),
  \]
  which in turn reduces to the case $N=\sO_Y$, where it becomes
  \[
  \sHom_Y(\pi_*\sO_X,\sO_Y)\cong \pi_*\sO_X\otimes
  \det(\pi_*\sO_X)^{-1}.
  \]
\end{proof}

\begin{rem}
  When $Y$ is smooth and $X$ is integral, this is essentially
  \cite[Prop.~0.1.3]{CossecFR/DolgachevIV:1989}.
\end{rem}

There is an alternate form for the inverse Morita equivalence.

\begin{prop}
  The functors $M\mapsto \pi_*{\cal L}\otimes_{\sO_Y} M$ and $M\mapsto
  \pi_*({\cal L}^{-1}\otimes \pi^!\sO_Y)\otimes_{{\cal A}_{{\cal L},\pi}}
  M$ are inverse equivalences between $\sO_Y\text{-mod}$ and ${\cal
    A}_{{\cal L},\pi}\text{-mod}$.
\end{prop}
 
\begin{proof}
  This reduces to showing that
  \[
  \pi_*{\cal L}
  \otimes_{\sO_Y}
  \pi_*({\cal L}^{-1}\otimes \pi^!\sO_Y)
  \cong
  {\cal A}_{{\cal L},\pi},
  \]
  which follows from the (bimodule!) isomorphism
  \[
  \pi_*({\cal L}^{-1}\otimes \pi^!\sO_Y) \cong \sHom_Y(\pi_*{\cal L},\sO_Y).
  \]
\end{proof}

The adjoint involution can be written the form $x\mapsto \Tr(x)-x$, and
thus in particular restricts to an involution $s:\pi_*\sO_X\to \pi_*\sO_X$
of the same form.  (This involution is nontrivial except when $\ch(k)=2$
and either $X$ is nonreduced or $\pi$ is inseparable.)  Given a coherent
sheaf $M$ on $X$, let $M s$ denote the $\sO_X$-bimodule with multiplication
$r_1 m r_2=r_1 s(r_2) m$; we will also conflate $\sO_X$-modules on $X$ with
the corresponding bimodules.  (Recall that for the moment, an
$\sO_X$-bimodule means a $\pi_*\sO_X$-bimodule on $Y$.)  And as we have
already mentioned, ${\cal A}_{{\cal L},\pi}$ has a natural $\sO_X$-bimodule
structure.

\begin{prop}\label{prop:aha_filtration}
  There is a natural short exact sequence of $\sO_X$-bimodules
  \[
  0\to \sO_X\to {\cal A}_{{\cal L},\pi}\to {\cal L}\otimes_{\sO_X} (\pi^!\sO_Y
  s)\otimes_{\sO_X} {\cal L}^{-1}\to 0
  \]
\end{prop}

\begin{proof}
  Since ${\cal A}_{{\cal L},\pi}\cong {\cal L}\otimes_{\sO_X} {\cal
    A}_{\sO_X,\pi}\otimes_{\sO_X} {\cal L}^{-1}$, it suffices to prove this
  in the case ${\cal L}=\sO_X$.  If we forget about the right module
  structure, then not only do we have such an exact sequence, but it
  splits.  Indeed, the natural inclusion $\sO_X\to {\cal A}_{\sO_X,\pi}$ is
  naturally split by the evaluate-at-1 map, so it remains to identify the
  other direct summand (the endomorphisms that annihilate 1) and show that
  the right action is triangular.

  There is a natural short exact sequence
  \[
  0\to \sO_Y\to \pi_*\sO_X\to \det(\pi_*\sO_X)\to 0;
  \]
  that the quotient is invertible follows as in
  \cite[\S0.1]{CossecFR/DolgachevIV:1989}, and can then be identified by
  comparing determinants.  An endomorphism that annihilates 1 annihilates
  $\sO_Y$ and is thus determined by its action on $\det(\pi_*\sO_X)$, so
  that we have the left $\sO_X$-module decomposition
  \[
  {\cal A}_{\sO_X,\pi}\cong \sO_X\oplus \sHom_Y(\det(\pi_*\sO_X),\pi_*\sO_X)
  \cong \sO_X\oplus \sHom_Y(\pi_*\sO_X,\sO_Y)
  \cong \sO_X\oplus \pi^!\sO_Y.
  \]
  Since the adjoint involution acts as $x\mapsto \Tr(x)-x$, it is clearly
  triangular with respect to this decomposition, and thus takes the
  diagonal left action of $\sO_X$ to a triangular right action.  Moreover,
  the induced action on the quotient is itself induced by the adjoint, or
  in other words by $s$.
\end{proof}

Note that the other (left) direct summand of ${\cal A}_{\sO_X,\pi}$ consists
(locally) of $s$-derivations, i.e. $\sO_Y$-linear endomorphisms that
satisfy
\[
\phi(r_1r_2) = s(r_2)\phi(r_1)+r_1\phi(r_2).
\]
If $\pi_*\sO_X$ is actually free over $\sO_Y$, with generator $\xi$, then
the morphism $\nu:a+b\xi\mapsto b$ is such an $s$-derivation, and
any other $s$-derivation has the form $f\mapsto g \nu(f)$.
If $\xi-s(\xi)$ is not a zero divisor, then we can write this in the form
\[
\nu(f)=\frac{f-s(f)}{\xi-s(\xi)}.
\]
The $\sO_Y$-module it generates is self-adjoint (the adjoint has eigenvalue
$-1$) and may be characterized as the space of {\em nilpotent}
endomorphisms that annihilate $\sO_Y$, which makes sense even when $\sO_X$
is only locally free over $\sO_Y$.  This module is locally free of rank 1,
and is easily seen to be isomorphic to $\pi^!\sO_Y\cong
\pi^*\det(\pi_*\sO_X)^{-1}$.

\medskip

Although it is most natural to view ${\cal A}_{{\cal L},\pi}$ as a sheaf of
$\sO_Y$-modules, we will need a slightly different perspective coming from
the bimodule structure.  Indeed, the $\pi_*\sO_X$-bimodule structure on
${\cal A}_{{\cal L},\pi}$ induces a $\pi_*\sO_X\otimes_{\sO_Y}
\pi_*\sO_X$-module structure, which in turn allows us to interpret it as a
coherent sheaf on $X\times_Y X=\Spec(\pi_*\sO_X\otimes_{\sO_Y}\pi_*\sO_X)$.
The embedding $X\times_Y X\to X\times X$ then makes it a coherent sheaf on
$X\times X$ supported on the union of the diagonal and the graph of $s$, so
finite over either projection.  In other words, ${\cal A}_{{\cal L},\pi}$
is a {\em sheaf bimodule} over $X\times X$ (in the sense of
\cite{VandenBerghM:1996,ArtinM/VandenBerghM:1990}).  Moreover, the algebra
structure on ${\cal A}_{{\cal L},\pi}$ is compatible with this bimodule
structure, in that the multiplication induces a morphism
\[
{\cal A}_{{\cal L},\pi}\otimes_X {\cal A}_{{\cal L},\pi}
:=
\pi_{2*}(\pi_3^*{\cal A}_{{\cal L},\pi}\otimes_{\sO_{X\times X\times X}}
\pi_1^*{\cal A}_{{\cal L},\pi})
\to
{\cal A}_{{\cal L},\pi},
\]
and thus ${\cal A}_{{\cal L},\pi}$ is a {\em sheaf algebra}.  (Indeed, this
is immediate from the corresponding fact for the fiber product.)

The Morita equivalences are also naturally expressed in terms of sheaf
bimodules; the sheaf $(1\times \pi)_*{\cal L}$ is a sheaf bimodule on
$X\times Y$, and the action of ${\cal A}_{{\cal L},\pi}$ induces a morphism
\[
{\cal A}_{{\cal L},\pi}\otimes_X (1\times \pi)_*{\cal L}
=
\pi_{2*}(\pi_3^*{\cal A}_{{\cal L},\pi}\otimes_{\sO_{X\times X\times Y}}
\pi_1^*(1\times\pi)_*{\cal L})
\to
(1\times \pi)_*{\cal L}
\]
making $(1\times\pi)_*{\cal L}$ a $({\cal A}_{{\cal L},\pi},\sO_Y)$-sheaf
bimodule.  Similarly, there is a $(\sO_Y,{\cal A}_{{\cal L},\pi})$-sheaf
bimodule given by
\[
(\pi\times 1)_*({\cal L}\otimes \pi^*\det\pi_*{\cal L}^{-1}),
\]
and these induce Morita equivalences.

One caution is that the natural adjoint involution does not in general
respect the sheaf bimodule structure, since it acts as $s$ on $\sO_X$.
Luckily, there is a variant that works more generally, at the cost of
changing ${\cal L}$.

\begin{prop}
  There is a natural sheaf algebra isomorphism ${\cal A}_{{\cal
      L},\pi}^{\text{op}}\cong {\cal A}_{{\cal L}^{-1},\pi}$.
\end{prop}

\begin{proof}
  Consider the composition
  \[
  \pi_*{\cal L}\otimes_{\sO_Y} \pi_*({\cal L}^{-1})
  \to
  \pi_*\sO_X
  \to
  \det(\pi_*\sO_X),
  \]
  where the first map is multiplication and the second is the quotient by
  the subbundle $\sO_Y$.  The corresponding pairing respects multiplication
  by $\sO_X$, and thus will induce a sheaf algebra isomorphism as required
  so long as it is perfect, i.e., the induced map
  \[
  \pi_*{\cal L}\to \sHom_{\sO_Y}(\pi_*({\cal L}^{-1}),\det(\pi_*\sO_X))
  \]
  is an isomorphism.  It suffices to check this locally, so that we may
  assume ${\cal L}=\sO_X$.  In that case, the (now symmetric) pairing
  vanishes on $\sO_Y\otimes \sO_Y$, and thus induces a pairing
  $\sO_Y\otimes \det(\pi_*\sO_X)\to \det(\pi_*\sO_X)$; since this is an
  isomorphism, the original pairing is perfect.
\end{proof}

\begin{rem}
  If ${\cal L}=\sO_X$, then $r$ and $s(r)$ always pair to $0$, and thus
  this involution is the composition of the adjoint with conjugation by
  $s$.  If also $\pi_*\sO_X$ is free, with basis $(1,\xi)$, then the
  pairing is given by $\begin{pmatrix} 0 & 1\\1 & \Tr(\xi)\end{pmatrix}$.
\end{rem}

Since this isomorphism is also an adjoint (albeit with respect to a
different pairing), and much more important for our purposes, we will refer
to the original involution on ${\cal A}_{{\cal L},\pi}$ as the
``intrinsic'' adjoint, reserving the unadorned term for this involution.

\section{A sheaf algebra coming from two involutions}

Now, suppose we are given two finite flat degree 2 morphisms $\pi_i:X\to
Y_i$, $i\in \{0,1\}$, as well as two invertible sheaves ${\cal L}_i$ on
$X$.  This gives rise to a pair ${\cal A}_{{\cal L}_0,\pi_0}$, ${\cal
  A}_{{\cal L}_1,\pi_1}$ of sheaf algebras on $X$, and we define the
(quasicoherent) sheaf algebra ${\cal H}_{{\cal L}_0,{\cal
    L}_1,\pi_0,\pi_1}$ on $X$ to be their pushforward over $\sO_X$.  That
is, ${\cal H}_{{\cal L}_0,{\cal L}_1,\pi_0,\pi_1}$ is generated by ${\cal
  A}_{{\cal L}_0,\pi_0}$ and ${\cal A}_{{\cal L}_1,\pi_1}$ subject only to
the relation that the common subalgebras $\sO_X$ should be identified.
Since the respective adjoints act trivially on $\sO_X$, they combine to
form an automorphism
\[
  {\cal H}_{{\cal L}_0,{\cal L}_1,\pi_0,\pi_1}^{\text{op}}
  \cong
  {\cal H}_{{\cal L}_0^{-1},{\cal L}_1^{-1},\pi_0,\pi_1}.
\]  
Note that when $Y_0\cong Y_1\cong \P^1$ and $X$ is a smooth genus 1 curve,
the resulting sheaf algebra is an instance of the ``type $C_1$ elliptic
double affine Hecke algebra'' construction of \cite{elldaha} (with the
Morita-equivalent algebra ${\cal S}$ being the spherical algebra in those
terms).  This suggests that there should be a multivariate version of this
more general construction, possibly for general root systems, but at the
very least for type $C$.

To work with this sheaf algebra, it will be convenient to work locally.
Localization is somewhat trickier with sheaf algebras than with sheaves of
algebras, since the restriction to an open subset is usually not a sheaf
algebra.  In this case, however, there is not too much difficulty.  If an
open subset is invariant under the action of the two involutions $s_0$,
$s_1$, then the restriction to that open subset will again be a sheaf
algebra.  This is still too restrictive (if $s_0s_1$ has infinite order, we
cannot expect to have a covering by affine opens of that form), but a
slight generalization will suffice in most cases.  Define a {\em
  localization} of $X$ to be a nonempty intersection of a nonempty (and
possibly infinite) collection of affine opens.  This is a fiber product of
affine morphisms, so inherits a scheme structure, and is affine over each
open set in the original collection, so affine, with coordinate ring given
by the limit of $\Gamma(U;\sO_X)$ over all open subsets $U$ containing the
intersection.  The significance of this notion is that although there are
not in general any nontrivial $\langle s_0,s_1\rangle$-invariant affine
opens in $X$, there are typically many invariant localizations.  Moreover,
${\cal H}_{{\cal L}_0,{\cal L}_1,\pi_0,\pi_1}$ induces a sheaf algebra on any
invariant localization, which will be an honest algebra containing the
structure sheaf of the localization.

If we not only have nontrivial invariant localizations but have a locally
finite covering by such localizations with consistent algebras on each
localization, then the theory of fpqc descent gives rise to an induced
sheaf algebra.  (This is discussed in detail in \cite[\S 5]{elldaha}.)  We
assume that not only does such a covering exist, but that there is such a
covering such that each bundle ${\cal L}_0$, ${\cal L}_1$, $\pi_{0*}\sO_X$,
$\pi_{1*}\sO_X$ becomes trivial.  Note that if $Y_0$ and $Y_1$ are smooth
curves over a field (which is the case we care about for the present
application), the existence of such a configuration is straightforward.
Indeed, any actual configuration of curves and bundles is the pullback of a
configuration over a field which is not algebraically closed, and any
bundle can be trivialized by removing finitely many points which are not
defined over that field (and thus such that no point of the orbit is in
that field).  It follows that any point of $X$ is contained in an invariant
localization that trivializes the bundles.  Moreover, each such
localization only omits finitely many orbits, and thus there is a finite
subcovering.

Thus let $R$ be a ring which is free of rank 2 over its subalgebras $S_0$
and $S_1$, and consider the algebra $H$ generated over $R$ by
$\End_{S_0}(R)$ and $\End_{S_1}(R)$.  We can give a more explicit
presentation of this algebra as follows.  The subspace of $\End_{S_0}(R)$
consisting of nilpotent elements annihilating $S_0$ is a free $S_0$-module,
so that we may choose a generator $\nu_0$ of this module, and similarly for
$\nu_1$.  Then $H$ has the presentation
\[
H = R\langle N_0,N_1\rangle/(N_0^2,N_0 r-s_0(r) N_0-\nu_0(r),N_1^2,N_1
r-s_1(r) N_1-\nu_1(r)).
\]
Note that since ${\cal L}_0$ and ${\cal L}_1$ are trivial, the adjoint
becomes an involution of this algebra, acting trivially on $R$, $N_0$, and
$N_1$; one has
\[
N_0 r-s_0(r) N_0
=
N_0 r+r N_0-\Tr_0(r)N_0
=
N_0 r+r N_0-N_0\Tr_0(r)
=
r N_0 - N_0 s_0(r),
\]
so that the relations are preserved, and thus this is a contravariant
automorphism, and is clearly of order 2.

Let $D_\infty$ denote the infinite dihedral group, with generating
involutions denoted by $s_0$, $s_1$, acting in the obvious (not necessarily
faithful) way on on $R$.  Each element of $D_\infty$ is represented by a
unique reduced word $s_{i_1}s_{i_2}\cdots$ in which $i_j\ne i_{j+1}$ for
all $j$.  (Thus there are two such words of each positive length,
corresponding to the two possible values of $i_1$.)  Given any such element
$w$, let $N_w$ denote the corresponding product $N_{i_1}N_{i_2}\cdots$.
There is a natural partial ordering on $D_\infty$ given by $w<w'$ iff
$\ell(w)<\ell(w')$, where $\ell(w)$ denotes the length of the reduced
word.\footnote{This is the Bruhat ordering on the infinite dihedral group,
  viewed as the free Coxeter group on two generators.}

\begin{prop}
  As a left $R$-module, one has
  \[
  H\cong \bigoplus_{w\in D_\infty} R N_w,
  \]
  and thus the partial ordering on $D_\infty$ induces a left $R$-module
  filtration of $H$.  This is in fact a bimodule filtration, and the
  subquotient corresponding to $w$ is the bimodule $R w$.
\end{prop}

\begin{proof}
  Define a left action of $H$ on the free module $\bigoplus_{w\in
    D_\infty}R e_w$ by extending the $R$ action by
  \[
  N_i\cdot (r e_w) =
  \begin{cases}
    s_i(r) e_{s_iw} + \nu_0(r)e_w & \ell(s_i w)>\ell(w),\\
    \nu_0(r)e_w & \ell(s_i w)<\ell(w).
  \end{cases}
  \]
  This is easily seen to satisfy the relations, so indeed gives a module
  structure, and since $r N_w\cdot e_1 = r e_w$, may be identified with the
  regular representation, establishing the desired isomorphism.  It remains
  to show that right multiplication by $R$ respects the filtration and acts
  correctly on the associated graded: i.e., that for $r\in R$, $N_w r-w(r)
  N_w\in \bigoplus_{w'<w} R N_{w'}$.  Moving $r$ to the left uses the
  relations $N_0 r=s_0(r) N_0+\nu_0(r)$ and $N_1 r = s_1(r)N_1+\nu_1(r)$,
  from which the result follows by induction in the length.
\end{proof}

\begin{rem}
  The adjoint immediately implies that $H$ is also free as a right
  $R$-module.  Moreover, the adjoint respects the filtration by $D_\infty$,
  modulo the (order-preserving) map $w\mapsto w^{-1}$.
\end{rem}

For the global version, we need to understand how the gluing interacts with
the filtration; in other words, we need to show that the corresponding
automorphisms of $H$ are triangular and understand their diagonal
coefficients.  The automorphisms are given by $N_0\mapsto u_0^{-1} N_0
u_0$, $N_1\mapsto u_1^{-1} N_1 u_1$ for $u_0,u_1\in R^*$, which clearly
respect the relations.  (And, of course, the adjoint involution inverts
both $u_0$ and $u_1$ as expected.)  Since
\[
u_0^{-1} N_0 u_0 = s_0(u_0)u_0^{-1} N_0 + u_0^{-1}\nu_0(u_0),
\]
we conclude that the action of this automorphism on $H$ is indeed
triangular with respect to the filtration, and on the subquotient
corresponding to $w=s_0s_1s_0\cdots$ acts as left multiplication by
\[
s_0(u_0)u_0^{-1}\times s_0(s_1(u_1)) s_0(u_1)^{-1}\times s_0(s_1(s_0(u_0)))
s_0(s_1(u_0))^{-1}\times \cdots
\]
(Note that for gluing purposes, we also need to take into account the fact
that $N_0$ is only determined up to multiplication by $S_0^*$; this gives
rise to factors $\pi^!\sO_{Y_0}$ in the resulting line bundles, and
similarly for $N_1$.)

With this in mind, let ${\cal N}_w$ denote the line bundle defined
inductively by ${\cal N}_1 = \sO_X$ and
\[
  {\cal N}_{s_i w} = s_i^*{\cal L}_i^{-1}\otimes {\cal L}_i\otimes
  \pi^!_i\sO_{Y_i}\otimes s_i^* {\cal N}_w
\]
whenever $\ell(s_iw)>\ell(w)$; note that more generally one has
\[
  {\cal N}_{ww'}\cong {\cal N}_w\otimes (w^{-1})^*{\cal N}_{w'}.
\]  
Similarly, let ${\cal H}_w$ be the sheaf subbimodule defined inductively by
${\cal H}_1 = \sO_X$ and
\[
{\cal H}_{s_i w} = {\cal A}_{{\cal L}_i,\pi_i}{\cal H}_w
\]
for $\ell(s_i w)>\ell(w)$.  In the affine case, this is precisely the
submodule corresponding to the interval under $w$ under the Bruhat
filtration.

\begin{thm}
  For any $w\in D_\infty$, ${\cal H}_w$ is locally free as a left
  $\sO_X$-module, of rank $|\{w':w'\le w\}|$, and in the corresponding
  bimodule filtration of ${\cal H}$, the subquotient corresponding to $w$
  is ${\cal N}_w w$.  Moreover, if $\ell(ww')=\ell(w)+\ell(w')$, then
  ${\cal H}_{ww'}={\cal H}_w{\cal H}_{w'}$.
\end{thm}

\begin{proof}
  Everything except the precise identification of the line bundle in the
  subquotient follows immediately from the corresponding statement in the
  affine case, while the identification of the line bundle follows by an
  easy induction from Proposition \ref{prop:aha_filtration}.
\end{proof}

Although this filtration is particularly convenient for calculations, we
will need some slightly coarser filtrations for the final construction.
Define subbimodules $\overline{\cal H}_{ij}$ for $i,j\in \Z$ with $j\ge i$ as
follows:
\[
  \overline{\cal H}_{i,i+2l} = {\cal H}_{s_{i}(s_{i+1}s_{i})^l},\qquad
  \overline{\cal H}_{i,i+2l+1} = {\cal H}_{(s_{i+1}s_{i})^{l+1}}
\]
where by convention $s_i:=s_{i\bmod 2}$.
(We add 1 to the length here to reflect the fact that ${\cal H}_{s_0}$ and
${\cal H}_{s_1}$ are subalgebras.)  We similarly denote the affine version
by $\overline{H}_{ij}$, and note that the adjoint identifies $\overline{H}_{ij}$ with
$\overline{H}_{-j,-i}$.

\begin{prop}
  For $i\le j\le k$ one has $\overline{\cal H}_{jk}\overline{\cal H}_{ij}\subset
  \overline{\cal H}_{ik}$.
\end{prop}

In other words, the sheaf bimodules $\overline{\cal H}_{j,k}$ fit together
to form a (positively graded) sheaf $\Z$-algebra $\overline{\cal H}$.
(I.e., an enhanced category with objects $\Z$ and Hom spaces given by sheaf
bimodules which are 0 unless the degree is nonnegative.)  This sheaf
$\Z$-algebra is manifestly invariant under shifting the degrees by 2, and
thus one could also consider the corresponding graded sheaf algebra, the
Rees algebra of ${\cal H}$ with respect to the filtration $\{\overline{\cal
  H}_{0,2l}:l\in \N\}$, without changing the category of coherent sheaves.

In the $\Z$-algebra form, each object has an endomorphism ring of the form
${\cal A}_{{\cal L},\pi}$, and thus there is a Morita-equivalent sheaf
$\Z$-algebra in which $\overline{\cal S}_{ij}$ is a sheaf bimodule on
$Y_i\times Y_j$.  (Here and below, we extend the indices on $Y_i$, $N_i$,
etc.~to be periodic of period 2.)  Locally, this has the following
description.  For $i\in \Z$, consider the (cyclic) $H$-module $H/\langle
N_i\rangle$, with generator $e_i$.  Then $\overline S_{ij}$ is the subspace
of $\Hom(H e_j,H e_i)^\text{op}$ such that the image of $e_j$ is in the
image of $\overline H_{ij}$; equivalently (since it is determined by the
image of $e_j$), it is the subspace of $\overline H_{ij} e_i$ which is
annihilated by $N_j$, and the composition $\overline S_{ij}\times \overline
S_{jk}\to \overline S_{ik}$ may be computed by taking $x * y$ to be the
image of $y$ after replacing $e_j$ by $x$.  For gluing purposes, in
addition to $N_i\mapsto u_i^{-1} N_i u_i$, we must take $e_i\mapsto
u_i^{-1} e_i$.  (Note that in the description of $\overline S_{ij}$ as a
subspace of $\overline H_{ij} e_i$, we must left multiply by $u_j$ after
applying the automorphism to $\overline H_{ij}e_i$.)

The following slightly modified local description is useful.

\begin{lem}
  We have $\overline S_{ij}=N_j\overline{H}_{ij}e_i$.
\end{lem}

\begin{proof}
  Since $\overline{H}_{ij}e_i$ is a projective $\End_{S_j}(R)$-module, this
  reduces to the fact that $\ker(N_j)=\im(N_j)$ inside $R$ (the unique
  indecomposable projective module).
\end{proof}

The relations in $\overline{H}_{ij}$ make it straightforward to give a direct
sum decomposition of $\overline{S}_{ij}$.  Indeed, $\overline{H}_{ij}e_i$ is clearly
the span of
\[
R N_j N_{j-1}\cdots N_{i+1} e_i,
R N_{j-1}\cdots N_{i+1} e_i,
\cdots
\]
since every other element of the standard basis annihilates $e_i$.
Since $N_j \overline{H}_{jj}=N_j R$, we can write any element of
$N_j \overline{H}_{ij} e_i=N_j\overline{H}_{i(j-1)} e_i$ as a sum
\[
N_j s_j(c_j) N_{j-1}\cdots N_{i+1} e_i
+N_{j-2} s_j(c_{j-2}) N_{j-1}\cdots N_{i+1} e_i
+\cdots
\]
with $c_j\in R$, except that when $j-i$ is even, the last term has the form
$c_i e_i$ with $c_i\in S_i$.  This, of course, glues to give an analogue
for $\overline{\cal S}_{ij}$ of the Bruhat filtration of $\overline{\cal H}_{ij}$.

\begin{lem}\label{lem:bruhat_for_S}
  One has $\overline{\cal S}_{ii}\cong \sO_{Y_i}$, while for $j>i$ one has a
  short exact sequence
  \[
  0\to \overline{\cal S}_{i(j-2)}\to \overline{\cal S}_{ij} \to (\pi_i\times
  \pi_j)_*( {\cal L}_j^{-1}\otimes{\cal N}_{s_j\cdot s_{i+1}}s_j\cdots
  s_{i+1}\otimes {\cal L}_i ) \to 0
  \]
  of $(\sO_{Y_i},\sO_{Y_j})$-bimodules.
\end{lem}

\begin{rem}
  Note that when taking the direct image under $\pi_i\times \pi_j$ that one
  must take into account the twisting by $s_j\cdots s_{i+1}$.  In addition,
  one can move the twist by $s_j$ to the left and absorb it into $\pi_j$
  to get an alternate description as the bimodule induced by
  \[
  {\cal L}_j^{-1}\otimes \pi^!_j\sO_{Y_j}\otimes {\cal N}_{s_{j-1}\cdots
    s_{i+1}}s_{j-1}\cdots s_{i+1}\otimes {\cal L}_i.
  \]
\end{rem}

\begin{lem}
  The sheaf $\Z$-algebra $\overline{\cal S}$ is a quadratic algebra; that is, it
  is generated by the $\Hom$ bimodules of degree 1 and the relations are
  generated by relations in degree 2.
\end{lem}

\begin{proof}
  It suffices to show this locally, i.e., for the $\Z$-algebra
  $\overline{S}$.  For generation in degree 1, we need to show that
  \[
  N_{j+1} \overline{H}_{j(j+1)} N_j \overline{H}_{0j} e_0
  =
  N_{j+1} \overline{H}_{0(j+1)} e_0,
  \]
  but this is an immediate consequence of the fact that $\overline{H}_{jj} N_j
  \overline{H}_{jj} = \overline{H}_{jj}$.

  It thus remains only to show that the only relations are in degree 2.
  For $i\in \{0,1\}$, let $\xi_i\in R$ be such that $\nu_i(\xi_i)=1$.  Then
  as a (bimodule) $\Z$-algebra, $\overline S$ is generated by elements $N_1
  e_0, N_1 \xi_1 e_0\in \overline S_{2i(2i+1)}$ and $N_0 e_1$, $N_0 \xi_0
  e_1\in \overline S_{(2i-1)2i}$ subject to the bimodule relations (i.e.,
  that right-multiplying a generator by an element of the relevant $S_i$
  equals the appropriate left-$S_{i+1}$-linear combination of generators)
  and the quadratic relations
  \begin{align}
  (N_0 e_1)(N_1 \xi_1 e_0)-(N_0 s_1(\xi_1) e_1)(N_1 e_0) &= 0\notag\\
  (N_1 e_0)(N_0 \xi_0 e_1)-(N_1 s_0(\xi_0) e_0)(N_0 e_1) &= 0.\notag
  \end{align}
  (This follows from $N_1 \xi_1-s_1(\xi_1)N_1=\nu_1(\xi_1)=1$ and $N_0
  e_0=0$.)  Using these quadratic relations, any monomial in the generators
  that has $(N_0 e_1)$ followed by $(N_1 \xi_1 e_0)$ or $(N_1 e_0)$
  followed by $(N_0 \xi_0 e_1)$ can be expressed in terms of monomials
  which are strictly smaller in an appropriate ordering: take $(N_1
  e_0)<(N_1 \xi_1 e_0)$ and $(N_0 e_1)<(N_0 \xi_0 e_1)$ and order monomials
  reverse lexicographically (i.e., compare the last factor, then the second
  to last, etc.).  We thus conclude that in the $\Z$-algebra presented in
  this way, each $\Hom$ bimodule has a left spanning set consisting of
  those monomials such that all appearances of $\xi_i$ are as far to the
  left as possible.  In $\Hom(i,j)$ for $j\ge i$, this spanning set
  contains $j-i+1$ elements; since the Bruhat filtration shows that
  $\overline S_{ij}$ is free of rank $j-i+1$, there can be no further
  relations.
\end{proof}

\begin{rem}
  In particular, to compute the adjoint, it suffices to compute it in
  degree 1, where we find that it swaps $(N_0 r e_1)$ and $(N_1 r e_0)$.
  Globally, we need to twist slightly to make the $\Hom$ bimodules agree,
  and thus either take ${\cal L}_i\mapsto {\cal L}_i^{-1}\otimes
  \pi_i^*\det(\pi_{i*}\sO_{Y_i})$ or ${\cal L}_i^{-1}\otimes
  \pi_i^*\omega_{Y_i}$.  (These are not the same but differ only by a
  common factor of $\omega_X$, which has no effect on the sheaf
  $\Z$-algebra.)
\end{rem}

Locally, the quadratic relation has the following less basis-dependent
description: we have a natural (surjective!) map
\[
\overline{H}_{11}\to \overline{S}_{02}
\]
given by $x\mapsto N_0 x e_0$, and thus an induced map
\[
R\to \overline{H}_{11}\to \overline{S}_{02}.
\]
The quadratic relation starting at $0$ is then the composition of this map
with the inclusion $S_0\to R$.  Much the same holds globally, except that
$\overline{H}_{11}$ gets twisted so that we instead have a surjection
\[
\sHom_{Y_1}(\pi_{1*}{\cal L}_1,\pi_{1*}({\cal L}_1\otimes \pi^!_0\sO_{Y_0}))
\to
\overline{\cal S}_{02}
\]   
and thus an induced map
\[
\pi_{0*}\pi^!_0\sO_{Y_0}\to \overline{\cal S}_{02},
\]
and the relation is the composition with the natural map
\[
\det(\pi_{0*}\sO_X)^{-1}\to \pi_{0*}\pi^!_0\sO_{Y_0}.
\]
Note that since the composition is a relation in $\overline{\cal S}_{02}$, the
cokernel maps to $\overline{\cal S}_{02}$, giving a natural global section of
$\overline{\cal S}_{02}$, locally given by the element $e_0$.

Since the $\Z$-algebra is generated in degree 1, we record the
corresponding special case of Lemma \ref{lem:bruhat_for_S}.

\begin{cor}
  One has the bimodule isomorphism
  \[
  \overline{\cal S}_{01}
  \cong
  (\pi_0\times \pi_1)_*({\cal L}_1^{-1}\otimes \pi^!_1\sO_{Y_1}\otimes {\cal L}_0)
  \]
\end{cor}

Let $\qcoh \overline{\cal S}$ denote the category of quasicoherent sheaves over
$\overline{\cal S}$; that is, the quotient of the category of $\overline{\cal
  S}$-modules by the subcategory generated by right-bounded modules.

\begin{thm}
  There is a natural equivalence between $\qcoh \overline{\cal S}$ and the
  noncommutative $\P^1$-bundle \cite{VandenBerghM:2012} corresponding to
  the sheaf bimodule $(\pi_0\times\pi_1)_*({\cal L}_0\otimes {\cal
    L}_1^{-1})$.
\end{thm}

\begin{proof}
  Both categories are constructed from quadratic sheaf $\Z$-algebras, so it
  will suffice to show that the algebras are twists of each other; that is,
  that we can associate a line bundle ${\cal L}'_i$ on $Y_i$ in each degree
  such that ${\cal L}'_i\otimes_{Y_i}\overline{\cal S}_{ij}\otimes_{Y_j}
  {\cal L}'_j$ agrees with the algebra constructed in
  \cite{VandenBerghM:2012}.  Since
  \[
  (\pi_0\times \pi_1)_*({\cal L}_1^{-1}\otimes \pi^!_1\sO_{Y_1}\otimes {\cal L}_0)
  \cong
  (\pi_0\times \pi_1)_*({\cal L}_1^{-1}\otimes {\cal L}_0)
  \otimes_{Y_1} \det(\pi_{1*}\sO_{Y_1})^{-1},
  \]
  this holds for $\overline{\cal S}_{01}$, and $\overline{\cal S}_{-10}$ is the
  appropriate adjoint bimodule, so is also correct.  Both constructions are
  invariant under shifting the object group, and thus it remains only to
  verify that the quadratic relations agree, and this is straightforward.
\end{proof}

This is most useful in the case that $Y_0$, $Y_1$ are smooth
quasiprojective curves over a field $k$, as in that case {\em every}
noncommutative $\P^1$-bundle arises in this way.

\begin{prop}
  Let $k$ be a field and $C_0/k$, $C_1/k$ smooth quasiprojective curves.
  For any sheaf bimodule ${\cal E}$ on $C_0\times C_1$ such that
  $\pi_{0*}{\cal E}$ and $\pi_{1*}{\cal E}$ are both locally free of rank
  2, there is a curve $\whQ$ with finite flat morphisms
  $\phi_i:\hat{Q}\to C_i$ of degree 2 and an invertible sheaf ${\cal L}$ on
  $\whQ$ such that ${\cal E}\cong (\phi_0\times \phi_1)_*{\cal L}$.
\end{prop}

\begin{proof}
  Suppose first that ${\cal E}$ is not a rank 2 vector bundle supported on
  the graph of an isomorphism, and let ${\cal Q}:=\sHom_{C_0\times
    C_1}({\cal E},{\cal E})$.  Then $\pi_{0*}{\cal Q}$ may be interpreted
  as the subalgebra of $\sEnd_{C_0}(\pi_{0*}{\cal E})$ consisting locally
  of elements that commute with the action of $\sO_{C_1}$.  The generic
  fiber has rank 2, and the cokernel is torsion-free, so flat, and thus
  $\pi_{0*}{\cal Q}$ is itself flat.  In particular $\whQ:=\Spec{\cal
    Q}$ is a double cover of $C_0$ (and thus of $C_1$, by symmetry) as
  required.  Since $\pi_{0*}{\cal Q}$ is saturated in
  $\sEnd_{C_0}(\pi_{0*}{\cal E})$, it follows that $\pi_{0*}{\cal E}$ is
  locally a cyclic module over $\pi_{0*}{\cal Q}$, and thus ${\cal E}$ is
  the image of an invertible sheaf as required.

  Now suppose ${\cal E}$ is a rank 2 vector bundle supported on the graph
  of an isomorphism, or WLOG that $C_0=C_1=C$ and ${\cal E}$ is a vector
  bundle on the diagonal.  Let $\whQ$ be any curve in the corresponding
  projective bundle over $C$ (a commutative ruled surface!) that does not
  contain any fiber and meets the generic fiber twice.  Then there is a
  line bundle ${\cal L}_0$ on $C$ such that $\pi_*\sO_{\whQ}\cong {\cal
    E}\otimes {\cal L}_0$, and thus ${\cal E}\cong \pi_*(\pi^*{\cal
    L}_0^{-1})$.
\end{proof}

\begin{rem}
  The notation $\whQ$ here reflects the fact that this curve is only the
  ``horizontal'' part of a more natural curve $Q$.
\end{rem}

Note that although both constructions give the same results geometrically,
they behave quite differently in families.  This appears not just in the
commutative case (where the $\whQ$-based construction depends on a
suitable choice of curve in the surface), but in the noncommutative case as
well.  The problem is that all we can say about the algebra ${\cal Q}$ in
general is that its cokernel in the relevant endomorphism ring is
torsion-free.  Over a curve, this is not a problem, but once we are dealing
with a family, it is quite possible for the quotient to fail to be flat.
Indeed, if $\overline{Q}$ is a nodal biquadratic curve in $\P^1\times
\P^1$, then there is a natural flat family of bimodules in $\P^1\times
\P^1$ parametrized by $\overline{Q}$, such that the fiber over each point
is the image of the corresponding ideal sheaf in $\sO_{\overline{Q}}$.
Over the smooth locus of $\overline{Q}$, the bimodule is the image of an
invertible sheaf on $\overline{Q}$, so has associated algebra
$\sO_{\overline{Q}}$, but over the singular point, the associated algebra
is instead the structure sheaf of the normalization.  But this cannot be
flat: $\sO_{\overline{Q}}$ has Euler characteristic 0, while the structure
sheaf of the normalization has Euler characteristic 1.  This also gives
rise to a sheaf bimodule on $(\P^1\times\overline{Q})\times (\P^1\times
\overline{Q})$ (or on the corresponding normalization!) such that the
corresponding noncommutative $\P^1$-bundle does not arise from the above
construction.  Luckily, our present interest is only in the curve case, and
many of the results we wish to show need only be checked on geometric
fibers.  (For instance, if we did not already know flatness from
\cite{VandenBerghM:2012}, we could easily prove it when $C_0$, $C_1$ are
projective by using the Bruhat filtration to check that the Euler
characteristic is the same on every geometric fiber.)

One nice aspect of the curve case is that there is a natural representation
of $\overline{\cal S}_{ij}$ in terms of (twisted) difference or differential
operators.  It suffices to give such an interpretation in degree 1 (the
quadratic relation will turn out to be automatically satisfied).  We may
write $\overline{\cal S}_{ij}$ as the tensor product
\[
(1\times \pi_1)_*({\cal L}_1^{-1}\otimes \pi^!_1\sO_{Y_1})
\times
(\pi_0\times 1)_*{\cal L}_0\otimes_{X}
.
\]
Here the second factor may be interpreted as the sheaf bimodule of
$\sO_{Y_0}$-linear maps $\sO_{Y_0}\to {\cal L}_0$, while the first factor
may similarly be viewed as the sheaf bimodule of $\sO_{Y_1}$-linear maps
${\cal L}_1\to \sO_{Y_1}$.  So if ${\cal L}_0={\cal L}_1$, then we may
interpret $\overline{\cal S}_{01}$ as the sheaf bimodule of compositions
\[
\sO_{Y_0}\to {\cal L}_0\to \sO_{Y_1}
\]
with the first factor $\sO_{Y_0}$-linear and the second factor
$\sO_{Y_1}$-linear.  In particular, if we take ${\cal L}_0=\sO_{X}$, then
the first factor can clearly be taken to be the natural inclusion, and thus
the map $\sO_{Y_0}\to \sO_{Y_1}$ is the restriction to $\sO_{Y_0}$ of a
general $\sO_{Y_1}$-linear map $\sO_{X}\to \sO_{Y_1}$.

Thus suppose $C_0/k$, $C_1/k$ are projective curves over a field $k$, or
localizations thereof, and let $\pi_i:\whQ\to C_i$ be a pair of double
covers.  There are three main cases to consider: both maps $\pi_i$ are
separable, both maps are inseparable, or precisely one of the two maps is
inseparable (and thus $k$ is a field of characteristic 2).  (In the first
two cases, we assume that the two maps cannot be identified by an
isomorphism $C_0\cong C_1$, and thus when both maps are inseparable,
$\whQ$ must be nonreduced.)

If $\whQ$ is reduced and both maps are separable, then the typical
$\sO_{C_1}$-linear map $\sO_{\whQ}\to \sO_{C_1}$ locally has the form
$f\mapsto c_0 f + s_1(c_0) s_1(f)$ with $c_0$ in the inverse different,
i.e., the sheaf of elements of $\sO_{\whQ}\otimes k(C_1)$ such that any
element of the corresponding fractional ideal has integral trace.  (As a
line bundle, the inverse different is isomorphic to $\pi_1^!\sO_{C_1}$, and
gives an expression of the latter in terms of an effective Cartier
divisor.)  The restriction to $\sO_{C_0}$ can be expressed in the same
form, but we can also write it as $f\mapsto \alpha f + s_1(\alpha)
s_1(s_0(f))$, and similarly for $\overline{\cal S}_{12}$.  In each case, the
operator involves either $s_1\circ s_0$ or its inverse, so that $\overline{\cal
  S}_{0d}$ maps to an operator of the form
\[
f\mapsto 
\sum_{-\lfloor d/2\rfloor\le i\le \lceil d/2\rceil} c_i (s_1\circ s_0)^i(f)
\]
with coefficients satisfying $c_{(d\bmod 2)-i}=s_d(c_i)$ (as well as
conditions along the ramification divisor which can be fairly complicated
in general, especially where $\whQ$ is singular).  In the case $\whQ$
is elliptic and $C_0,C_1$ rational, this recovers the untwisted version of
the symmetric elliptic difference operators of \cite{generic}.

If $\whQ$ is reduced but not integral, then each component induces an
isomorphism between $C_0$ and $C_1$, so we may as well use one component to
identify them with a single curve $C$, and then the other component induces
a nontrivial automorphism $\alpha\in \Aut(C)$.  In that case, the
description of $\overline{\cal S}_{01}$ simplifies to operators of the form
\[
f\mapsto c_0 f + c_1 \alpha(f)
\]
where $c_0$, $c_1\in \sO_C$ must agree to appropriate order wherever the
corresponding branches of $\whQ$ meet (necessarily at fixed points of
$\alpha$).  More generally, $\overline{\cal S}_{0d}$ maps to operators
\[
f
\mapsto 
\sum_{-\lfloor d/2\rfloor\le i\le \lceil d/2\rceil} c_i \alpha^i(f).
\]

In the nonreduced case, the reduced subscheme of $\whQ$ induces an
isomorphism $C_0\cong C_1$, so that again we may as well identify both
curves with some fixed curve $C$.  Then the coordinate ring of $\whQ$
over $C_0$ is
\[
\sO_{\whQ} = \sO_C\oplus \epsilon {\cal L}
\]
for some line bundle ${\cal L}$ on $C$.  The embedding $\sO_{C_1}\to
\sO_{\whQ}$ then has the form $f\mapsto f+\epsilon \partial f$ where
$\partial$ is a derivation $\sO_{C_1}\to {\cal L}$.  It follows that ${\cal
  L}$ contains $\omega_C$, and thus there is some divisor $D_c$
(essentially a conductor) such that ${\cal L}\cong \omega_C(D_c)$ and
$\partial f=df$.  The general $\sO_{C_1}$-linear map $\sO_{\whQ}\to
\sO_{C_1}$ then takes the form
\[
f+g \epsilon\mapsto \alpha f + \beta(g-df)
\]
with $\alpha\in \sO_{\whQ}$ and $\beta\in \omega_{C_0}(D_c)^{-1}$.
Restricting this to $\sO_{C_0}$ gives
\[
f\mapsto \alpha f - \beta df,
\]
or in other words the general first-order differential operator such that
the coefficient of differentiation vanishes along $D_c$.  Similarly, the
restriction to $\sO_{C_1}$ of the general $\sO_{C_0}$-linear map
$\sO_{\whQ}\to \sO_{C_0}$ takes the form $f\mapsto \alpha f + \beta
df$.  We thus see more generally that $\overline{\cal S}$ has a
representation inside the sheaf of differential operators on $C$ such that
the image of ${\cal S}_{ij}$ locally consists of operators of order $j-i$
such that the coefficient of $(d/dz)^l$ vanishes along $lD_c$.

Finally, if only one of the two maps is separable (our usual odd
characteristic 2 case), then we end up with a hybrid operator
representation.  Supposing $\pi_1$ is the separable morphism, then
$\overline{\cal S}_{01}$ has the same interpretation as in the purely
separable case ($f\mapsto \alpha f + s_1(\alpha f)$ for suitable $\alpha$).
On the other hand, $\overline{\cal S}_{12}$ consists of the restriction to
$\sO_{C_1}$ of linear maps $\sO_{\whQ}\to \sO_{C_0}$.  Since $\pi_0$ is
inseparable, for any local section $f$ of $\sO_{\whQ}$, $df$ may be
interpreted as a local section on $C_0$ of $\omega_{C_0}\otimes
\det(\pi_{0*}\sO_{\whQ})^{-1}$, so that if $g$ is any local section of
$\pi_0^*(\omega_{C_0}^{-1}\otimes \det(\pi_{0*}\sO_{\whQ}))$, then
$f\mapsto d(fg)$ is a linear map $\sO_{\whQ}\to \sO_{C_0}$, and any such
map arises in this way.  We thus see that $\overline{\cal S}_{12}$ has a
natural interpretation as differential operators.  It is worth noting that
the resulting hybrid operator representation is far from faithful; indeed,
it already has a nontrivial kernel in degree 2.  (We will see below that in
this case $\overline{\cal S}$ has the same category of coherent sheaves as
a maximal order in a quaternion algebra on a commutative ruled surface.)

When ${\cal L}_0\otimes {\cal L}_1^{-1}$ is nontrivial, the above
interpretation must of course be twisted accordingly.  To understand the
nature of such twisting, it suffices to consider the corresponding
automorphisms in the affine case.  Although this in principle involves a
pair of units, only their ratio actually appears, and one thus obtains the
automorphism $N_1 r e_0\mapsto N_1 u r e_0$ and $N_0 r e_1\mapsto N_0
u^{-1} r e_1$ of $\overline S$ for some unit $u\in R^*$.  To interpret this,
note that $N_1 e_0$ is an operator from $S_0$ to $S_1$ that annihilates
$1$, and similarly for $N_0 e_1$, while $N_1 \xi_1 e_0$ and $N_0 \xi_0 e_1$
similarly map $1$ to $1$.  Let $F_{0,u}$ be a formal solution (in $S_0$) of
the equation $N_1 u e_0\cdot F_{0,u}=0$, and let $F_{1,u}$ be the formal
image $N_1 u \xi_1 e_0 \cdot F_{0,u}$.  We then find that
\[
(N_0 u^{-1} e_1)\cdot F_{1,u}
=
(N_0 u^{-1} N_1 u \xi_1 e_0)\cdot F_{0,u}
=
(N_0 u^{-1} s_1(\xi_1) N_1 u e_0)\cdot F_{0,u}
=
0,
\]
so that $F_{1,u}$ is a formal solution to $N_0 u^{-1} e_1\cdot F_{1,u}=0$.
Similarly,
\[
(N_0 u^{-1} \xi_0 e_1)\cdot F_{1,u} = e_0\cdot F_{0,u},
\]
so that one may reasonably identify $F_{2,u}$ with $F_{0,u}$.  We also find
that for $a\in S_0$,
\begin{align}
(N_1 u a e_0)\cdot F_{0,u} &= \nu_1(a) F_{1,u},\notag\\
(N_1 u \xi_1 a e_0)\cdot F_{0,u} &= \nu_1(a\xi_1) F_{1,u}
\end{align}
corresponding to the formal identities
\begin{align}
F_{1,u}^{-1}(N_1 u e_0)F_{0,u} &= N_1 e_0.\notag\\
F_{1,u}^{-1}(N_1 u \xi_1 e_0)F_{0,u} &= N_1 \xi_1 e_0.
\end{align}
In other words, the automorphism corresponding to $u$ may be interpreted as
gauging by the system of formal symbols $F_{i,u}$.

In the differential case, the formal equation satisfied by $F_{0,u}$ has
the form $(D-v)F_{0,u}=0$ for suitable $v$; that is, $F_{0,u}$ is a formal
symbol with logarithmic derivative $v$, and one further finds that
$F_{1,u}=(u\bmod \epsilon)F_{0,u}$.  Similarly, in the nonsymmetric
difference case, $F_{0,u}$ is a formal solution to $\alpha(F_{0,u}) = v
F_{0,u}$, while in the symmetric difference case, it is a formal solution
to the equations $s_0(F_{0,u})=F_{0,u}$ and $s_1(uF_{0,u})=uF_{0,u}$.  Note
that in suitably analytic settings, one can represent $F_{0,u}$ and
$F_{1,u}$ by honest functions, e.g., as the exponential of the integral of
the appropriate meromorphic differential, or (when $s_0s_1$ has infinite
order) as a suitable infinite product.

Any local trivialization of ${\cal L}_0\otimes {\cal L}_1^{-1}$ leads to an
representation of $\overline{\cal S}$ as sheaves of {\em meromorphic}
operators of the appropriate kind, and the above calculation shows that any
other such representation will be related by a suitable scalar gauge
transformation (by the solution of a first-order equation).  This is
particularly useful in the rational case; when $C_0\cong C_1\cong \P^1$ and
the image of $\whQ$ in $\P^1\times \P^1$ is singular, we can always
localize to the complement of a singular point, which will make all curves
affine and all line bundles trivial.  This covers every case for which
$C_0\cong C_1\cong \P^1$ except the case when $\whQ$ is smooth genus 1, for
which one can use theta functions or a certain elliptic gauge on $F_1$, see
\cite{generic}.  (Note that although theta functions are generally thought
of as analytic objects, when suitably normalized, they extend as canonical
sections of line bundles on the moduli stack of genus 1 curves with marked
points, see \cite[\S2]{elldaha}, and thus the explicit operators of
\cite{generic} make sense in general.)

\medskip

One useful consequence of the operator interpretation is that it makes it
straightforward to show that $\overline{\cal S}_{ij}$ is a domain.  We in fact
have the following.

\begin{prop}
  Let $Y_0$, $Y_1$ be integral schemes, and let ${\cal E}$ be a sheaf
  bimodule of birank $(2,2)$ on $Y_0\times Y_1$.  Then the corresponding
  sheaf $\Z$-algebra is a domain.
\end{prop}

\begin{proof}
  Enlarging the sheaf $\Z$-algebra can only introduce zero divisors, so we
  may as well restrict to the generic points of $Y_0$ and $Y_1$.  This
  makes ${\cal E}$ invertible on its support, and thus the resulting
  $\Z$-algebra is of the form $\overline{S}_{ij}$ with
  $S_0=\overline{S}_{00}$, $S_1=\overline{S}_{11}$ fields.

  If $R$ is a field, then given any pair of morphisms composing to 0, one
  of their leading coefficients with respect to the Bruhat filtration must
  vanish, since the leading coefficient of the product is the (automorphism
  twisted) product of leading coefficients.  But in that case we can factor
  out the degree 2 element $e_0$ or $e_1$ on the appropriate side to obtain
  a smaller such pair, eventually yielding a contradiction.

  If $R$ is a sum of two fields, then the proof of the difference operator
  interpretation gives a faithful homomorphism from $\overline{S}_{ij}$ to
  the twisted group algebra $S_0[\langle \alpha\rangle]$.  But this is a
  domain (by the same leading coefficient argument).

  Finally, if $R$ is nonreduced, then $\overline{S}_{ij}$ consists of
  polynomials in some fixed derivation (or, rather, an element acting as
  such a derivation) of $S_0$, and again leading coefficients produce the
  domain property.
\end{proof}

\section{The curve of points}

As mentioned, the curve $\whQ$ is only a piece of a somewhat larger
curve contained in the noncommutative surface.  That $\whQ$ itself
embeds in the surface follows by observing that the system of elements
$(e_i)\in \overline{S}_{i(i+2)}$ is ``central'' in a suitable sense;
indeed, if $D\in \overline{S}_{ij}$ then $(e_j)D = D(e_i)$, where we use
the natural identification between $\overline{S}_{ij}$ and
$\overline{S}_{(i+2)(j+2)}$.  The quotient by the corresponding ideal takes
$\overline{S}_{ij}$ to the quotient
$\overline{S}_{ij}/\overline{S}_{i(j-2)}$, which as computed above
globalizes to the image of an invertible sheaf on $\whQ$.  We thus
find that the quotient $\Z$-algebra is a twisted version (\`a la
\cite{ArtinM/VandenBerghM:1990}) of the homogeneous coordinate ring of
$\whQ$, and thus the corresponding category of sheaves is just
$\qcoh \whQ$.  Moreover, this not only embeds $\whQ$ in the
noncommutative surface, but it does so {\em as a divisor} in the sense of
\cite{VandenBerghM:1998}, which in particular means that any point of
$\whQ$ is a suitable candidate for the blowing up construction given
therein.  (Note that everything stated here applies more generally to embed
$X$ as a divisor in the corresponding noncommutative $\P^1$-bundle.)  To be
precise, there is an endofunctor ${-}(-\whQ)$ of $\qcoh
\overline{S}_{ij}$, which simply shifts indices by $2$, and a natural
transformation ${-}(-\whQ)\to \text{id}$ which is 0 on a subcategory
isomorphic to $\qcoh\whQ$, and has cokernel mapping to that
subcategory.

The larger curve is also embedded as a divisor, and arises by asking for
the largest quotient $\Z$-algebra that satisfies the same twisted
commutativity relation as $\whQ$ (or $X$).  Passing to the local case,
we find that we should consider the two-sided ideal generated by
quasi-commutators of the form
\[
(N_0 x e_1)(N_1 s_1(y) e_0) - (N_0 y e_1)(N_1 s_1(x) e_0)
\quad\text{or}\quad
(N_1 x e_0)(N_0 s_0(y) e_1) - (N_1 y e_0)(N_0 s_0(x) e_1).
\]
We may simplify
\[
(N_0 x e_1)(N_1 s_1(y) e_0) - (N_0 y e_1)(N_1 s_1(x) e_0)
=
N_0 (x N_1 s_1(y)-y N_1 s_1(x)) e_0,
\]
and observe that since $N_1$ is nilpotent, so is negated by the intrinsic
adjoint in $\End_{S_1}(R)$, the terms $x N_1 s_1(y)$ and $-y N_1 s_1(x)$
are intrinsic adjoints, and thus their sum is in $S_1$.  (Moreover, taking
$x=1$ gives $\nu_1(y)$, and thus any element of $S_1$ arises in this way.)
We thus see that the two-sided ideal of quasi-commutators is the same as
that generated by $N_i S_{i+1}e_i = \nu_i(S_{i+1}) e_i\subset
\overline{S}_{i(i+2)}$.  In particular, it is indeed contained in the ideal
corresponding to $\whQ$.  Moreover, it is essentially given by a conductor.

\begin{lem}
  We have $\nu_0(S_1)R=\nu_1(S_0)R=\{x:x\in R\mid xR\subset S_0S_1\}$.
\end{lem}

\begin{proof}
  For $f\in S_1$, $g\in R$, we have $\nu_0(f)g = \nu_0(f
  s_0(g))+f\nu_0(g)\in S_0S_1$, so that
  \[
  \nu_0(S_1)R\subset \{x:x\in R\mid xR\subset S_0S_1\}.
  \]
  Conversely, if $cR\subset S_0S_1$, and $\xi_0$ is such that
  $\nu_0(\xi_0)=1$, then both $c$ and $c s_0(\xi_0)$ are in $S_0S_1$,
  so that $\nu_0(c)$ and $\nu_0(c s_0(\xi_0))$ are in
  $\nu_0(S_0S_1)=S_0\nu_0(S_1)$, and thus
  \[
  c = \nu_0(c s_0(\xi_0))+\xi_0\nu_0(c)\in \nu_0(S_1)R
  \]
  as required.
\end{proof}

\begin{rem}
  Note that if $S_0$, $S_1$ are Dedekind domains with distinct images
  inside $R$, then both $R$ and $S_0S_1$ are orders inside the
  normalization, and $\{x:x\in R\mid xR\subset S_0S_1\}$ is the conductor
  of $S_0S_1$ as a suborder of $R$.  (In the commutative case, this ideal
  is of course $0$.)  This can be described in more global terms by noting
  that, as an ideal in $S_0$, the conductor is $\nu_0(S_0S_1)$, or in other
  words the image of the suborder under the natural map $R\mapsto R/S_0$.
  We thus see that globally (for curves) the corresponding ideal sheaf is
  the image of the natural map $\det(\pi_{0*}\sO_{\overline{Q}})\otimes
  \det(\pi_{0*}\sO_{\whQ})^{-1}\to \sO_{C_0}$, where $\overline{Q}$
  is the image of $\whQ$ in $C_0\times C_1$ and the map is induced
  by the natural inclusion $\sO_{\overline{Q}}\to \sO_{\whQ}$, which
  is the identity on $\sO_{C_0}$.  This measures the failure of ${\cal E}$
  to be invertible along $\overline{Q}$.
\end{rem}

In the curve case (and assuming $\pi_0$ and $\pi_1$ are not related by an
isomorphism $C_0\cong C_1$), the ideals $S_0\nu_0(S_1)$ and $S_1\nu_1(S_0)$
globalize to line bundles contained in $\sO_{C_0}$, $\sO_{C_1}$
respectively, with the property that both bundles pull back to isomorphic
bundles on $\whQ$.  Normally, twisting ${\cal L}_0$ or ${\cal L}_1$ by
pulled back line bundles gives an equivalence between different surfaces,
but in this case the twisting has no effect on $\overline{\cal S}_{i(i+1)}$ and
thus gives rise to an autoequivalence.  Composing this with the shift-by-2
autoequivalence gives an autoequivalence equipped with a natural
transformation to the identity such that the cokernel is quasicommutative.
In particular, we find that the result is a twisted version of the
homogeneous coordinate ring of a commutative scheme.  Of course, if $\pi_0$
and $\pi_1$ {\em are} related by an isomorphism, then the ideal is trivial,
and $\overline{\cal S}$ is itself the twisted homogeneous coordinate ring of a
commutative scheme (in this case a ruled surface).  We denote the
corresponding autoequivalence of $\qcoh \overline{\cal S}$ by ${-}(-Q)$, and
note that the natural transformation ${-}(-Q)\to \text{id}$ factors through
${-}(-\whQ)\to\text{id}$.

To identify this scheme, suppose $k\subset S_0\cap S_1$ is some subring,
and let $M\subset R$ be an $S_0S_1$-submodule such that $M/R$ is locally
free of rank 1 over $k$.  Base changing to make this free, we then see that
there are ideals $I_i\subset S_i$ such that $M\supset I_i R$, and then
for $y\in I_1$,
\[
(N_1 R e_0)(N_0 M e_1)
\ni
(N_1 e_0)(N_0 \xi y e_1)-(N_1 s_0(\xi) e_0)(N_0 y e_1)
=
N_1 y e_1
=
\nu_1(y) e_1,
\]
so that $(N_1 R e_0)(N_0 M e_1)\supset \nu_1(I_1)e_1$.  Any element of
$S_1$ can be written as $a+y$ for $a\in k$, $y\in I_1$, and thus
$\nu_1(I_1)=\nu_1(S_1)$, so that the left ideal generated by $N_0 M e_1$
contains the left ideal generated by $\nu_1(S_1)$, and thus the two-sided
ideal so generated.  In particular, for $r\in R$, $m\in M$, we have
\[
(N_1 r e_0)(N_0 m e_1)\equiv (N_1 s_0(m) e_0)(N_0 s_0(r) e_1)
\]
modulo that ideal, and thus
\[
(N_1 R e_0)(N_0 M e_1)
=
(N_1 s_1(M) e_0)(N_0 R e_1).
\]
We thus obtain two conclusions: first that any such $M$ induces a
representation of $\overline{S}_{ij}$ which is locally free of rank 1 over $k$
in every degree, and second the degree 2 elements of the form
$\nu_1(S_0)e_1$ and $\nu_0(S_1)e_0$ act as 0 on any such representation.
In other words (and globalizing), the corresponding commutative scheme is
precisely the moduli space of such representations, which in turn may be
identified with $\Quot(\overline{\cal S}_{01},1)$ (the moduli space globalizing the moduli space of
submodules $M$).

These identifications were already known for general noncommutative
$\P^1$-bundles \cite{VandenBerghM:2012,NymanA:2004}; the one new
observation is that in the surface case (i.e., when $Y_0$, $Y_1$ are
curves) the resulting curve (denoted $Q$) is embedded as a divisor in the
noncommutative surface.  In other words, {\em any} point of a
noncommutative surface obtained from this construction is eligible to be
blown up \`a la \cite{VandenBerghM:1998}.  To see that $Q$ is a curve, note
that since a $S_0S_1$-module is also an $S_0$-module, there is a natural
morphism $Q=\Quot(\overline{\cal S}_{01},1)\to \Quot(\pi_{0*}\overline{\cal
  S}_{01},1)$, embedding $Q$ as a subscheme of a $\P^1$ bundle over $C_0$.
Moreover, that subscheme is locally cut out by the equation $v\wedge
\phi(v)=0$, where $\phi\in \End_{\sO_{C_0}}(\pi_{0*}\overline{\cal
  S}_{01})$ is the image of a local generator of $\sO_{\overline{Q}}$ over
$\sO_{C_0}$.  Since
\[
v\wedge \phi(v) = \gamma (v\wedge \xi v)
\]
where $\xi$ is a local generator of $\sO_{\whQ}$ over $\sO_{C_0}$ and
$\gamma$ is a local generator of the conductor, we see that $Q$ splits as a
divisor isomorphic to $\whQ$ and transverse to every fiber, plus a
vertical divisor given by the pullback of the conductor.  Note in
particular that this implies that the arithmetic genera of the three curves
are related by $g(Q)=g(\whQ)+c=g(\overline{Q})$, where $c$ is the degree of
the conductor.  Indeed, the first equation follows from intersection
theory, while the second follows by comparing differents.

We should also note that in the 1-dimensional case, the pullback of any
$S_0$-ideal containing the conductor descends to an $S_1$-ideal containing
the conductor, and thus gives rise to a corresponding commutative curve
embedded as a divisor.  This is mainly relevant in the commutative case
(when the conductor is 0), as it means that we still have natural models of
$\whQ$ and various curves $Q$ inside the surface.

\section{Ruled and quasi-ruled surfaces}

We would like to have a classification of sheaf bimodules associated to
noncommutative surfaces, but this in general appears to be quite difficult
to control; indeed, any smooth curve with a pair of order 2 automorphisms
gives rise to at least one such bimodule for every invertible sheaf (and
possibly more, if the composition of the automorphisms has fixed points).
Although this is largely intractable, it turns out that nearly every such
configuration is uninteresting from a noncommutative geometry perspective,
as the corresponding noncommutative $\P^1$-bundle can be described as a
finite algebra over a corresponding commutative ruled surface.

With this in mind, we define a noncommutative {\em quasi-ruled} surface to
be any surface obtained from the above construction in the case that $C_0$,
$C_1$ are smooth projective curves, so that we may single out a subclass of
noncommutative {\em ruled} surfaces, which (by Proposition
\ref{prop:qr_is_semicomm} and Theorem \ref{thm:semicomm} below) include all
of the interesting cases and are much easier to classify.  We say that the
quasi-ruled surface corresponding to $\whQ$, $\pi_i$, ${\cal L}_i$ is a
(geometrically) ruled surface if (possibly after a field extension) there
is an isomorphism $\phi:C_0\cong C_1$ such that the divisor $(\pi_0\times
\pi_1)(\whQ)$ is algebraically equivalent to twice the graph of $\phi$.
This of course includes the commutative case $\pi_1=\pi_0\circ\phi$, as
well as the case that $\whQ$ is nonreduced (since then the reduced
subscheme is the graph of such a $\phi$).  Moreover, any component of the
moduli stack of quasi-ruled surfaces that contains a ruled surface contains
either a commutative surface or a surface on which $\whQ$ is nonreduced,
and it is easy to see that the latter can be degenerated to the commutative
case.  So roughly speaking the noncommutative ruled surfaces are those
cases which are deformations of commutative ruled surfaces (though there is
one exotic type of such a deformation which is not ruled, see below; this
corresponds to the quasi-standard Poisson surfaces of Chapter
\ref{chap:classify}).

The ruled surfaces are also nearly characterized by the following fact.

\begin{prop}\label{prop:genus_1_almost_implies_ruled}
  The curve $Q$ associated to a quasi-ruled surface has arithmetic genus
  $\ge 1$, with equality precisely when either the surface is ruled or $Q$
  is smooth of genus 1 and $\pi_0$ and $\pi_1$ are distinct $2$-isogenies.
\end{prop}

\begin{proof}
  We have already shown that $g(Q)=g(\overline{Q})$, so it suffices to consider
  the latter.  In particular, since the double diagonal in $C\times C$ has
  arithmetic genus 1, we certainly have equality in the ruled surface case.
  Similarly, in the $2$-isogeny case, $Q$ is embedded in the product of
  quotients, and thus $\overline{Q}\cong Q$ has arithmetic genus 1.

  For the other direction, we may assume $\overline{Q}$ reduced.  Let
  $\widetilde{Q}$ be the normalization of $\overline{Q}$, and observe that
  $g(\widetilde{Q})\le g(\overline{Q})$.  Since $\widetilde{Q}$ may not be connected,
  its arithmetic genus could be negative, but since it has at most two
  components, the arithmetic genus is at least $-1$.  In any event, if
  $g(\widetilde{Q})\le 0$, then some component of $\widetilde{Q}$ has genus 0, and
  $C_0$ and $C_1$ are images of that component, so must themselves have
  genus 0.  But in that case, we (geometrically) have $C_0\times C_1\cong
  \P^1\times \P^1$ with $\overline{Q}$ of bidegree $(2,2)$, so are in the ruled
  case.

  We thus reduce to the case that $g(\widetilde{Q})=1$ and no component has
  genus 0.  Moreover, $\widetilde{Q}$ must embed in $C_0\times C_1$, since
  otherwise $\overline{Q}$ has strictly larger arithmetic genus than
  $\widetilde{Q}$.  If $\widetilde{Q}$ is reducible, then both components
  must be smooth genus 1 curves isomorphic to both $C_0$ and $C_1$.
  Without loss of generality, we may use one component to identify $C_0$
  and $C_1$, so that $\widetilde{Q}$ is the union of the diagonal and the
  graph of some automorphism of $C$.  If that automorphism had fixed
  points, then the two graphs would intersect, and thus the automorphism
  must be a translation; since the group of translations is connected,
  $\overline{Q}$ is algebraically equivalent to the double diagonal, so the
  surface is ruled.  If $\widetilde{Q}$ is integral, then any degree 2
  morphism is either a $2$-isogeny or a map to $\P^1$.  If both maps are to
  $\P^1$, then we are again in the bidegree $(2,2)$ case, while if both
  maps are $2$-isogenies then they are distinct (since otherwise
  $\overline{Q}$ would be nonreduced).  In the remaining case, $\bar{Q}$ is
  a bidegree $(2,2)$ curve on a surface $E\times \P^1$, so has arithmetic
  genus 3 by Hirzebruch-Riemann-Roch.
\end{proof}

\begin{rems}
  Unlike $\whQ$, the curve $Q$ {\em is} flat for any flat family of rank
  $(2,2)$ sheaf bimodules, and thus in particular its arithmetic genus is
  constant on any component of the moduli stack.  As a result, we find that
  any component on which the genus is $>1$ cannot possibly contain any
  commutative fibers.  Since the components corresponding to ruled surfaces
  certainly do contain such fibers, the only question involves the
  $2$-isogeny case.  When $\Pic^0(C)[2]$ is \'etale, the requirement that
  the $2$-isogenies are distinct is preserved under arbitrary deformation,
  and thus any commutative fiber must have characteristic 2.  On the other
  hand, the modular curve parametrizing pairs of $2$-isogenies
  (equivalently the modular curve classifying ``cyclic'' $4$-isogenies)
  splits into three components in characteristic 2, and on one of those
  components, the two $2$-isogenies agree.  In other words, there exists a
  curve over a $2$-adic dvr such that the generic fiber admits a pair of
  distinct $2$-isogenies which become the same on the special fiber.
  Taking that curve to be $\whQ$ gives a $2$-adic point of the given
  component of the moduli stack such that the special fiber is indeed
  commutative.  In the untwisted case, the special fiber is
  $\P(\pi_*\sO_{\whQ})$, and $\pi_*\sO_{\whQ}$ is the unique indecomposable
  bundle of rank 2 and trivial determinant.  The deformation induces a
  Poisson structure on that ruled surface which vanishes precisely on the
  image of $\whQ$.  This explains the exotic characteristic 2 Poisson
  surfaces of Chapter \ref{chap:classify}, and in particular shows that
  those Poisson structures do, indeed, arise as semiclassical limits of
  noncommutative surfaces.
\end{rems}

\begin{rems}
  In the ruled and $2$-isogeny cases, since every connected component of
  $\overline{Q}$ has arithmetic genus 1, there is a noncanonical isomorphism
  $\omega_{\overline{Q}}\cong \sO_{\overline{Q}}$.  Since $\overline{Q}$ is a double cover
  of $C_0$, we can also write $\omega_{\overline{Q}}\cong \pi_0^!\omega_{C_0}
  \cong \pi_0^*(\omega_{C_0}\otimes \det(\pi_{0*}\sO_{\overline{Q}})^{-1})$, and
  similarly for $C_1$.  In the ruled cases, one in fact has
  $\omega_{C_i}\cong \det(\pi_{i*}\sO_{\overline{Q}})$ for $i\in \{0,1\}$;
  indeed one can check in each case that $\Pic(C_i)\to \Pic(\overline{Q})$ is
  injective.  In the $2$-isogeny case, this fails for at least one of the
  $C_i$.  Indeed, in that case, $\det(\pi_{i*}\sO_{\overline{Q}})$ corresponds
  to the nontrivial point of the kernel of the dual $2$-isogeny; since
  $\overline{Q}$ has at most one $2$-isogeny with inseparable dual, we have
  $\det(\pi_{i*}\sO_{\overline{Q}})\not\cong \sO_{C_i}\cong \omega_{C_i}$ for at
  least one $i$.  We will see below that the isomorphism $\omega_{C_i}\cong
  \det(\pi_{0*}\sO_{\overline{Q}})^{-1}$ is equivalent to the divisor $Q$ of
  points being anticanonical (i.e., such that ${-}(-Q)[2]$ is a Serre
  functor), giving us another characterization of ruled surfaces.
\end{rems}

The key to showing that non-ruled quasi-ruled surfaces are nearly
commutative is the following result.

\begin{prop}\label{prop:qr_is_semicomm}
  If $\whQ$, $\pi_0$, $\pi_1$ corresponds to a noncommutative
  quasi-ruled surface which is not ruled, then the automorphism $s_0s_1$
  has finite order.
\end{prop}

\begin{proof}
  Suppose otherwise.  Since we have excluded the ruled surface case, either
  $g(\overline{Q})>1$ or $\overline{Q}$ is a smooth genus 1 curve and
  $s_0$, $s_1$ are translations by distinct $2$-torsion points (one of
  which may be the identity in characteristic 2).  In the latter case,
  $s_0s_1$ is again translation by a $2$-torsion point, so has order 2.  We
  may thus assume $g(\overline{Q})>1$.  Since $s_0s_1$ acts on the
  normalization $\widetilde{Q}$ preserving the components, and smooth
  projective curves of genus $>1$ have finite automorphism groups, the
  claim is automatic unless the (isomorphic via $s_0$) components have
  genus $\le 1$; since a genus 0 component would force the configuration to
  correspond to a ruled surface, the components in fact have genus $1$.
  But $s_0s_1$ not only acts on $\widetilde{Q}$, but preserves the finite
  set of preimages of singular points of $\overline{Q}$ (which exist, since
  otherwise it would have arithmetic genus 1), so that some power of
  $s_0s_1$ fixes all preimages and a further power fixes $\overline{Q}$.
\end{proof}

In general, if $\pi_0$ and $\pi_1$ are separable, and $s_0s_1$ has finite
order, then the \'etale algebra $k(\whQ)^{\langle s_0s_1\rangle}$ is a
quadratic extension of the field $k(\whQ)^{\langle s_0,s_1\rangle}$.  We
let $C'$ be the smooth curve with the latter function field, and define
$\whQ'$ to be the unique curve with the former generic \'etale algebra such
that the conductor of $\whQ'$ in its normalization is the norm of the
corresponding conductor for $\whQ$.  We will see in Chapter
\ref{chap:semicomm} that the center of the untwisted quasi-ruled surface
corresponding to $\pi_i:\whQ\to C_i$ is the (commutative) untwisted
quasi-ruled surface corresponding to $\pi:\whQ'\to C'$, in the sense that
the category of coherent sheaves on the noncommutative surface is the
category of coherent modules of a sheaf of algebras on the commutative
surface.  The twisted situation is slightly trickier to describe, but we
note that any line bundle on $\whQ$ can be represented by a Cartier divisor
supported on the smooth locus of $\overline{Q}$, and thus by taking images
induces a Cartier divisor on the corresponding locus of $\whQ'$.  This
induces a well-defined map from the possible twists of the noncommutative
surface to the possible twists of the commutative surface, compatible with
the sheaf of algebras structure.

The hybrid case is similar; if $\pi_0$ is inseparable, then we take $C'$ to
be the Frobenius image of $C_1$, so that $C_0$ is a separable double
cover of $C'$, and we take $\whQ'$ to be the corresponding singular
model of $C_0$, such that the conductor is the norm of the conductor of
$\whQ$.

Finally, if $\whQ$ is nonreduced of characteristic $p$ (so that WLOG
$C_0=C_1=C$ and $\overline{Q}$ is the double diagonal), we take $C'$ to be
the image of $C$ under Frobenius.  If $\whQ=\overline{Q}$, then we take
$\whQ'$ to be the double diagonal in $C'$, and in general construct $\whQ'$
so that its conductor relative to the double diagonal is the norm of that
of $\whQ$.  Nontrivial twists (i.e., not in $\Pic(C)$) only exist when the
conductor is trivial, in which case they are given (modulo $\Pic(C)$) by a
class in $H^1(\omega_C)$.  To describe the resulting map $H^1(\omega_C)\to
H^1(\omega_{C'})$, it suffices to give a local ($\F_p$-linear) map
$\Omega_C\to \Omega_{C'}$ inducing it.

The key idea is the following observation.

\begin{prop}\label{prop:twist_differential}
  Let $u\in k(C)$ be a function which is a uniformizer at some point.  Then
  in the algebra of differential operators on $k(C)$, one has the relation
  \[
  \Bigl(\frac{d}{du} + f\Bigr)^p
  =
  \Bigl(\frac{d}{du}\Bigr)^p + f^p + \frac{d^{p-1}f}{du^{p-1}}
  \]
  for any function $f$.
\end{prop}

\begin{proof}
  Since $\frac{d}{du}\mapsto \frac{d}{du}+f$ induces an automorphism of the
  algebra of differential operators and $(d/du)^p$ is central, the operator
  $(\frac{d}{du}+f)^p$ must also be central, and thus a (monic)
  $k(C')$-linear combination of $(d/du)^p$ and $1$.  Applying this to a
  composition of such automorphisms, we find that the map
  \[
  f\mapsto 
  \Bigl(\frac{d}{du} + f\Bigr)^p
  -
  \Bigl(\frac{d}{du}\Bigr)^p
  \]
  from $k(C)$ to $k(C')$ is $\F_p$-linear.

  This map is given by a polynomial in $f$ and its iterated derivatives
  which is homogeneous of degree $p$ with respect to the grading in which
  $\deg(d^lf/du^l)=l+1$, corresponding to rescaling both $u$ and $f$.
  Moreover, $\F_p$-linearity implies that every monomial that appears must
  have degree a power of $p$ with respect to the grading by rescaling $f$
  alone.  It follows that
  \[
  \Bigl(\frac{d}{du}+f\Bigr)^p = \Bigl(\frac{d}{du}\Bigr)^p+a f^p + b
  \frac{d^{p-1}f}{du^{p-1}}
  \]
  for suitable constants $a,b\in \F_p$.  Each of those monomials arises
  from exactly one of the $2^p$ terms on the left, and thus one easily
  verifies that $a=b=1$ as required.
\end{proof}

\begin{lem}\label{lem:quasi-norm-differential}
  For any meromorphic differential $\omega$ on $C$, the differential
  \[
  \tau(\omega)
  :=
  \Bigl[\Bigl(\frac{d+\omega}{du}\Bigr)^p-\Bigl(\frac{d}{du}\Bigr)^p\Bigr]
  d(u^p)
  =
  \Bigl[\Bigl(\frac{\omega}{du}\Bigr)^p + \frac{d^{p-1}}{du^{p-1}}\frac{\omega}{du}\Bigr]d(u^p)
  \]
  on $C'$ is independent of $u$, so induces a canonical $\F_p$-linear sheaf
  map $\tau:\Omega_C\to \Omega_{C'}$.
\end{lem}

\begin{proof}
  We first claim that for $f\in k(C)$, we have the invariant description
  $\tau(df) = d(f^p)$.  Indeed, plugging into the definition of $\tau$ and
  using the fact that $p$-th derivatives vanish gives us the expression
  $\tau(df) = \Bigl(\frac{df}{du}\Bigr)^p d(u^p)$, which agrees with
  $d(f^p)$ when $f$ is a power of $u$ and thus by $\F_p$-linearity for all
  $f\in k((u))$, and in particular on the subspace $k(C)$.

  Similarly, if $g\in k(C')$, so that $dg=0$, we claim that $\tau(g f^{-1}
  df) = (g^p-g) f^{-p} d(f^p)$ for any $f\in k(C)^*$.  It suffices to
  verify this in the completion, and thus by continuity when $f$ is a
  polynomial, and thus, since both sides are multiplicative homomorphisms
  in $f$, for $f=u-u_0$.  But in that case we explicitly compute
  \[
  \frac{\tau(g (u-u_0)^{-1}d(u-u_0))}{d((u-u_0)^p)}
  =
  g^p (u-u_0)^{-p} + g \frac{d^{p-1}(u-u_0)^{-1}}{du^{p-1}}
  =
  (g^p-g)(u-u_0)^{-p}.
  \]

  The topological span of the two $\F_p$-subspaces on which these invariant
  descriptions are defined is $k((u))du$ (indeed, each differential $\alpha
  u^k du$ can be expressed in one of the two forms), and thus the two
  formulas uniquely determine $\tau$, making it invariant.  Moreover, if
  $\omega$ is holomorphic at $x$, then we may use a uniformizer at $x$ to
  see that that $\tau(\omega)$ is also holomorphic, so that it gives a
  sheaf map as required.
\end{proof}

\begin{cor}\label{cor:definition_of_Nmprime}
  There is a multiplicative monoid map $\Nm':k(C)\oplus
  \epsilon\Omega_{k(C)}\to k(C')\oplus \epsilon \Omega_{k(C')}$ given by
  \[
  \Nm'(f+\epsilon \omega)\mapsto f^p + f^p\tau(\omega/f) \epsilon
  \]
  when $f\ne 0$ and with $\Nm'(\epsilon g dh) = g^p d(h^p)\epsilon$.
  Moreover, $\Nm'(f)=\Nm'(f+\epsilon df)=f^p$.
\end{cor}

\begin{proof}
  The definition of $\tau$ via a uniformizer $u$ gives an expression for
  $\Nm'(f+\epsilon g du)$ which is polynomial in the derivatives of $f$ and
  $g$, and agrees with the given description when $f=0$.  Furthermore, the
  values of $\Nm'(f)$ and $\Nm'(f+\epsilon df)$ follow from the definition of
  $\tau$.  To show that $\Nm'$ is multiplicative, it will suffice to show it
  for $f$ nonzero, i.e., that it gives a homomorphism on
  $(\sO_{\overline{Q}}\otimes k(C))^*$.  Since
  \[
  \Nm'(h(f+\epsilon \omega)) = h^p \Nm'(f+\epsilon \omega) =
  \Nm'(h)\Nm'(f+\epsilon\omega),
  \]
  we may further reduce to the case $f=1$, where it follows from the fact
  that $\tau$ is an additive homomorphism.
\end{proof}

By considering how $\tau$ acts when applied to a differential with poles,
we find that the restriction of $\Nm'$ to $\sO_{\whQ}$ takes values in
$\sO_{\whQ'}$, and thus in particular gives a well-defined homomorphism
$H^1(\sO_{\whQ}^*)\to H^1(\sO_{\whQ'}^*)$ between the respective Picard
groups.

\medskip

In each case, we have associated a double cover of a smooth curve and a
line bundle on that double cover, and thus a commutative ruled surface
corresponding to the direct image of that line bundle.  We will see below
that the original quasi-ruled (possibly ruled) surface is a maximal order
in a certain division ring over the function field of this commutative
ruled surface.

\section{Classification of ruled surfaces}

As we mentioned above, although quasi-ruled surfaces are difficult to fully
classify, the ruled case is fairly straightforward.  We should first note
that the isomorphism $C_0\cong C_1$ making the sheaf bimodule algebraically
equivalent to the double diagonal is not in general unique: we can compose
it with anything in the identity component of $\Aut(C_1)$ without actually
changing the noncommutative surface (which depends only on the category of
representations and the choice of structure sheaf, i.e., the pullback of
$\sO_{C_0}$).  We can similarly twist by the pullback of any line bundle on
$C_1$ without any effect.  (Twisting by a line bundle on $C_0$ gives a
Morita-equivalent noncommutative surface.)  Thus when classifying ruled
surfaces, we should feel free to identify $C_0\cong C_1$ with $C$, but must
take these symmetries into account.

We assume we are over an algebraically closed field; although the
classification could be done over a general field with only slightly more
work, the result would be misleading, since not every noncommutative
surface over $k$ which is ruled over $\overline{k}$ is ruled over $k$.  (This is
true even in the commutative case: not only must one contend with conic
bundles in general, but the fact that $\P^1\times \P^1$ has two geometric
rulings, neither of which need be defined over the ground field.  These
issues can be dealt with, but the simplest approach is via yet another
construction to be discussed below.)

In the nonreduced case (e.g., if $g(C_0)\ge 2$), we have already seen that
$\whQ$ has the form $\Spec(\sO_{C_0}\oplus \omega_{C_0}(D)\epsilon)$ for some
effective divisor $D$, with $\epsilon^2=0$, and it remains only to
understand the possible twists, i.e., $\Pic(\whQ)/\Pic(C_1)$.  (The
automorphism freedom is fixed by taking $\overline{Q}$ to be the double
diagonal.)  We claim that in fact $\Pic(\whQ)/\Pic(C_1)=1$ unless $D=0$,
when it is isomorphic to $\G_a(k)$.  Restricting to the underlying reduced
curve gives a splitting of $\Pic(C_1)\to \Pic(\whQ)$, so it suffices to
understand the case that ${\cal L}$ is trivial on the reduced curve.  In
that case, the corresponding cocycle takes values in the group $1+\epsilon
\omega_{C_0}(-D)$, or equivalently in the sheaf $\omega_{C_0}(-D)$, and
thus $\Pic(\whQ)/\Pic(C_1)$ is canonically isomorphic to
$H^1(\omega_{C_0}(D))\cong H^0(\sO_C(-D))$.  In other words, if $D>0$,
there is no scope for twisting, while for $D=0$ there is a one-parameter
family of twists (corresponding to the usual notion of twisting for
differential operators).

In the reduced but nonintegral case, we split into three cases.  The first
is that $g(C_0)=g(C_1)=1$.  In that case, we note that any divisor which is
algebraically equivalent to the double diagonal is a union of two
translates of the diagonal, and thus we may use the $\Aut(C_1)^0$ freedom
to make one of the two components be the diagonal.  Then the other
component is disjoint from the diagonal, so that $\whQ=\overline{Q}$ is
smooth.  We can further use the $\Pic(C_1)$ freedom to make the line bundle
trivial on the diagonal.  We thus find that such cases are classified by
$\Pic^0(C)\times \Pic(C)$, with the first coordinate determining the
translation and the second determining the line bundle on the corresponding
component, and the corresponding algebra is a compactification of the
relevant algebra of twisted elliptic difference operators.  (To be precise,
these surfaces are classified by the quotient of $\Pic^0(C)\times \Pic(C)$
by the involution corresponding to swapping the two components.)

In the reduced but nonintegral cases with $g(C_0)=g(C_1)=0$, we similarly
find that the curve is determined by a choice of nontrivial automorphism.
Geometrically, an element of $\PGL_2(k)$ is either diagonalizable,
corresponding to a linear fractional transformation $z\mapsto qz$ (the
multiplicative case), or unipotent, corresponding to $z\mapsto z+\hbar$
(the additive case).  There are a few possible choices for $\whQ$: in
the multiplicative case, there are four, depending on whether $\whQ$ is
singular over $0$ or $\infty$, while in the additive case there are three:
$\whQ=\overline{Q}$, $\whQ=\widetilde{Q}$, and an intermediate case for which
the conductor has degree 1.  There is a nontrivial map
$\Pic(\whQ)/\Pic(C_0)\to \Z$ given by taking the degree on the
nonidentity component, which is an isomorphism unless $\whQ=\overline{Q}$,
when the kernel is $\G_m$ or $\G_a$ as appropriate.

Finally, in the reduced integral case, where necessarily $g(C_0)=g(C_1)=0$,
$\overline{Q}$ is an integral biquadratic curve, which is either smooth, nodal,
or cuspidal.  The smooth case is classified by a smooth genus 1 curve $Q$
along with two classes $\eta_0,\eta_1\in \Pic^2(Q)$ (i.e., the two maps to
$\P^1$) and one in $\Pic(Q)/\langle \eta_1\rangle$ (the line bundle for
twisting), and corresponds to the construction in \cite{generic}, i.e.,
symmetric elliptic difference operators.

In the nodal and cuspidal cases, $\widetilde{Q}\cong \P^1$ is (assuming
both maps are separable) equipped with a pair of involutions, and $\whQ\in
\{\widetilde{Q},\overline{Q}\}$.  If the product of involutions has two
fixed points, then the involutions swap the fixed points and can be made to
look like $z\mapsto 1/z$, $z\mapsto q/z$, giving symmetric $q$-difference
operators (with $\overline{Q}$ a nodal curve).  If the product of
involutions has a single fixed point, then the involutions have the form
$z\mapsto a-z$, $z\mapsto a+\hbar-z$, where we may take $a=0$ except in
characteristic 2.  In characteristic 2, we also have the (hybrid)
inseparable case in which one of the maps is $z\mapsto z^2$ and the other
can be taken to be $z\mapsto z^2+z$.  In each of these cases, there are two
possible twists when $\whQ\ne \overline{Q}$ (since $\Pic(\P^1)$ injects as
an index 2 subgroup of $\Pic(\whQ)$), and otherwise the group of twists is
an extension of $\Z/2\Z$ by $\G_m$ or $\G_a$ as appropriate.

For Hirzebruch surfaces (i.e., ruled surfaces over $\P^1$), there is
another approach to the classification which is particularly convenient for
some purposes.  Recall that $Q$ naturally embeds in a ruled surface over
$C_1=\P^1$; furthermore, it intersects the generic fiber transversely with
multiplicity 2, and has arithmetic genus 1, thus must be anticanonical.
There is thus a natural map from the moduli stack of noncommutative
Hirzebruch surfaces to the moduli stack of (commutative) anticanonical
Hirzebruch surfaces.  Given a noncommutative Hirzebruch surface, there is
also an associated line bundle $q$ of degree 0 on $Q$ given by
$\pi_0^*\sO_{\P^1}(-1)\otimes \pi_1^*\sO_{\P^1}(1)$, giving a map from the
moduli stack of noncommutative Hirzebruch surfaces to the relative $\Pic^0$
(that is, the stack classifying invertible sheaves which are degree 0 on
every component of every fiber) of the universal anticanonical curve
of the moduli stack of anticanonical Hirzebruch surfaces.

We claim that this map to the relative $\Pic^0$ is an isomorphism.  Indeed,
we can recover line bundles $\eta_0=\pi_0^*\sO_{\P^1}(1)$ and
$\eta_1=\eta_0\otimes q$ from the given data.  Both bundles are acyclic of
degree 2: their inverses have degree $-2$ and have (the same) nonpositive
degree on every component, so are ineffective.  Thus both bundles induce
(isomorphism classes of) maps from $Q$ to $\P^1$, and since they have
degree 0 on every vertical component, this restricts to an embedding
$\overline{Q}\to \P^1\times \P^1$.  The sheaf bimodule can then be
recovered as the direct image of the relative $\sO(1)$ of the commutative
Hirzebruch surface.  (Since $\whQ$ meets every fiber twice, there is
a natural map $\sO_Q(1)\to \sO_{\whQ}(1)$ with kernel an extension
of sheaves $\sO_f(-1)$ supported on fibers of the ruling, and thus both
sheaves have the same direct image.)

Note here that $\Pic^0(Q)=\Pic^0(\overline{Q})$ is an elliptic curve,
$\G_m$, or $\G_a$, depending on whether $Q$ is smooth, nodal, or
cuspidal/nonreduced.  Moreover, when $\overline{Q}$ is reduced, $q$
corresponds directly to the shift in the difference operator
interpretation, while for $\overline{Q}$ nonreduced, it corresponds to the
traditional $\hbar$ used to view differential operators as deformations of
commutative polynomials.

There is a minor technical issue ignored above: to recover the data from
the surface, we must in general include the embedding of $Q$ (i.e.,
consider the moduli stack of {\em anticanonical} noncommutative Hirzebruch
surfaces), as when $q$ is the identity, the corresponding surface is
commutative.  We also caution the reader that although every
(anticanonical) commutative Hirzebruch surface over the generic point of a
dvr extends (nonuniquely), this fails in the noncommutative case, since
$\Pic^0(Q)$ need not be proper.  (This is related to the fact that although
$\overline{Q}$ has a limit inside $C_0\times C_1$, that limit may contain
preimages of points of $C_0$ or $C_1$, making it impossible for a limiting
sheaf on $\overline{Q}$ to be a sheaf bimodule.)

Something similar holds for other ruled surfaces; in general, the moduli
stack of (anticanonical) noncommutative ruled surfaces over $C$ is the
relative $\Pic^0(Q)/\Pic^0(C)$ over the moduli stack of anticanonical ruled
surfaces over $C$.  (If $C$ is a genus 1 curve of characteristic 2, we must
also exclude those anticanonical ruled surfaces for which the anticanonical
curve is connected and smooth, as these correspond to the $2$-isogeny
case.)  The map from noncommutative surfaces is the same as for genus 0
(where we note that $\Pic^0(C_0)$ and $\Pic^0(C_1)$ have the same embedding
in $\Pic^0(Q)$), but the inverse is trickier to construct.  However, we can
verify in both the differential and elliptic difference cases that the map
is surjective and that we can reconstruct the sheaf bimodule from this
data.

\chapter{The semicommutative case}
\label{chap:semicomm}

\section{Quasi-ruled surfaces as maximal orders}

As mentioned above, many examples of noncommutative quasi-ruled surfaces
are not ``truly'' noncommutative, in that the corresponding sheaf
$\Z$-algebra has a large center.  Our objective in the present section is
to make this more precise, by constructing a (coherent) sheaf of algebras
on a corresponding commutative ruled surface $Z$ having the same category
of coherent modules.  In fact, we will see that this $\sO_Z$-algebra is a
{\em maximal order} (see \cite{ReinerI:2003} for an expository treatment).
(Recall that an {\em order} on an integral scheme $Z$ is a coherent,
torsion-free, $\sO_Z$-algebra such that the generic fiber is a central
simple algebra, and an order is {\em maximal} if it is not a proper
subalgebra of an order with the same generic fiber.)  This not only
completes our justification for frequently restricting our attention to
ruled surfaces (as maximal orders have been extensively studied via
commutative means), but gives useful additional information about
semicommutative ruled surfaces (the differential case in characteristic
$p$, as well as difference cases in which the translation is torsion).

Our first goal is to establish the following result.

\begin{thm}\label{thm:semicomm}
  Suppose that $X$ is the quasi-ruled surface associated to maps
  $\pi_i:\hat{Q}\to C_i$ (with $C_i$ smooth projective curves) and a line
  bundle ${\cal L}$ on $\hat{Q}$.  Suppose that $k(C_0)$ is algebraic of
  degree $r$ over the intersection $k(C_0)\cap k(C_1)$ inside $k(\hat{Q})$,
  and let $C'$ be the curve with function field $k(C_0)\cap k(C_1)$.  Then
  there is a commutative ruled surface $Z$ over $C'$ and a maximal order
  ${\cal A}$ on $Z$ with generic fiber a division ring of index $r$ such
  that there is a natural equivalence $\coh {\cal A}\cong \coh X$ taking
  the regular representation to the structure sheaf.
\end{thm}

\begin{rem}
  By Proposition \ref{prop:qr_is_semicomm}, the hypothesis is automatically
  satisfied (for some $r$) for any quasi-ruled surface which is not ruled.
  We will further be able to identify the rank 2 vector bundle on $C'$
  corresponding to $Z$, although this identification is surprisingly subtle
  when ${\cal L}$ is nontrivial.
\end{rem}

A key observation is that being a maximal order is a local condition (up to
some mild subtleties in how one identifies the center), so to a large
extent we may restrict our attention to the analogous affine question.

\section{The differential case}

Let us first consider the differential case in characteristic $p$, where as
usual we may assume that $R=S_0[\epsilon]/\epsilon^2$, with $S_1\cong S_0$
embedded as $f\mapsto f+(df/\omega)\epsilon$ for some nonvanishing
meromorphic differential $\omega$.  Let $u$ be a uniformizer at some point
of the curve $\Spec(S_0)$, and observe that the differential $du$ is
holomorphic and nonvanishing in a neighborhood of that point.  We may thus
reduce by further localization to the case that $du$ is nowhere vanishing,
and thus the conductor is principal, generated by $c:=du/\omega\in S_0$.
We then find that $\overline{S}_{01}=\overline{S}_{12}=S_0+S_0 c D$, where
$D$ satisfies the commutation relation $D f - f D = df/du$.  In other
words, $\overline{S}$ is the $\Z$-algebra associated to the subalgebra
$A_c:=S_0\langle t, cDt\rangle$ of the graded algebra $A=S_0\langle
t,Dt\rangle$, where $t$ is central.

\begin{lem}
  If $S_0$ has characteristic $p$, then the $p$-th Veronese of the center
  of $A_c$ is $S_0^p[t^p,c^p D^p t^p]$, and the corresponding sheaf of
  algebras over $\Proj(S_0^p[t^p,c^p D^p t^p])$ is locally free.
\end{lem}

\begin{proof}
  We first need to show that the central element $c^p D^p t^p\in A$ is in
  $A_c$.  But this follows easily from the identity
  \[
  c^a (c D t) = ((c D t) - a c' t) c^a
  \]
  valid over $S_0\langle t,Dt\rangle$, where $c'=dc/du$; in
  particular, we find 
  \[
  c^k D^k t^k =
  c^{k-1} (cDt)D^{k-1}t^{k-1}
  = (c D t-(k-1)c't) c^{k-1}D^{k-1} t^{k-1},
  \]
  so that by induction, $c^kD^kt^k\in A_c$ for all $k$.

  As an operator, $c^pD^pt^p$ annihilates $S_0$, and any homogeneous
  element of degree $d$ in $A_c$ can be expressed as $t^d$ times a
  polynomial of degree $<p$ in $cD$ plus a multiple of $c^p D^p$.  Since
  differential operators of degree $<p$ act faithfully, we conclude that
  $c^p D^p t^p$ generates the kernel of the natural map $A_c\to
  \End_{S_0^p} (S_0)[t]$.  A central element $x\in A_c$ of degree a
  multiple of $p$ maps to a central element of $\End_{S_0^p}(S_0)[t]$, and
  thus has the form $f_0 t^{pd} + (c^p D^p t^p)y$ where $f_0\in S_0^p$ and
  $x\in A_c$.  Since $A_c$ is a domain, $y$ must also be central, and since
  it has lower degree than $x$, the first claim follows by induction.

  For local freeness, it suffices to show that the submodule of $A_c$
  consisting of elements of degree congruent to $-1$ modulo $p$ is free
  over the center of $A_c$.  A degree $dp+(p-1)$ element of $A_c$ has a
  unique expression of the form
  \[
  \sum_{0\le i\le d} x_i t^{p(d-i)}(c^p D^p t^p)^i
  \]
  with each $x_i$ homogeneous of degree $p-1$.  Since the space of degree
  $p-1$ elements is free as a $S_0^p$-module, the claim follows.
\end{proof}

In particular, we have established that $\coh \overline{S}$ is the category of
coherent sheaves on an order ${\cal A}$ over the surface
$Z:=\Proj(S_0^p[t^p,c^p D^p t^p])$, which is a ruled surface over
$\Spec(S_0^p)$.  Note that since $\overline{S}$ is a domain, the generic fiber
of ${\cal A}$ is a division ring.  Over most of $Z$, maximality is
straightforward, and in fact we have the following stronger statement.

\begin{lem}
  On the complement of the divisor $c^p t^p=0$, ${\cal A}$ is an Azumaya
  algebra.
\end{lem}

\begin{proof}
  Without loss of generality, we may assume that $c$ is a unit, and thus
  this reduces to showing that $S_0\langle D\rangle$ is an Azumaya algebra
  over $\Spec(S_0^p[D^p])$.  The fiber over any point (with field $\ell$)
  of $\Spec(S_0^p[D^p])$ has a presentation of the form $\ell\langle
  u,D\rangle$ with relations $u^p=U$, $D^p=\delta$, $Du-uD=1$, with
  $U,\delta\in \ell$ and both $u$ and $D$ commuting with $k$.  Passing to
  the algebraic closure of $\ell$ lets us subtract $p$th roots of $U$ and
  $\delta$ from $D$ and $u$ respectively, and thus reduce to the case
  $U=\delta=0$, in which case the algebra is readily identified with
  $\Mat_p(\overline\ell)$ as required.
\end{proof}

Since ${\cal A}$ is locally free, and thus reflexive, it is maximal iff its
codimension 1 localizations are maximal, so that it remains only to check
maximality along components of $c^p t^p=0$.  Here we may use the criterion
that an order over a dvr is maximal if its radical is principal (as a left
ideal) and its semisimple quotient is central simple.  (If the radical is
principal, it is invertible, and thus the order is hereditary, and the
result follows by \cite[Thm.~2.3]{AuslanderM/GoldmanO:1960}).  There are
two cases to consider: the fibers where $c$ vanishes, and the curve
$t^p=0$.

For points where the conductor vanishes, we may assume $S_0$ is a dvr and
the conductor has the form $c=z^c$ for some uniformizer $z$ and some
positive integer $c$, and consider the completion $\hat{S}$ of $S$ along
the divisor $z=0$, or equivalently the base change of $S|_{t=1}$ to the
extension $k(z^{cp}D^p)[[z^p]]$ of its center.  Consider the two-sided
ideal $\hat{S} z \hat{S}$.  Since $c$ is positive, $z$ is normalizing:
\[
z (z^c D) = (z^c D-z^{c-1})z,
\]
and thus $\hat{S}z \hat{S}=z \hat{S}$ is principal as a left ideal.
Moreover, the $p$-th power ideal $(z \hat{S})^p = z^p \hat{S}$ is contained
in the radical, since $z^p$ generates the maximal ideal of the center.  We
conclude that $\hat{S} z\hat{S}$ is contained in the radical of $\hat{S}$.
The quotient by this two-sided ideal is the field $k(z^c D)$, and thus
$\hat{S} z \hat{S}$ actually is the radical, and the quotient is a field,
so is central simple as required.

At a point where $t^p=0$, we have a completion of the form $K\langle\langle
D^{-1}\rangle\rangle$ where $K$ is the field of fractions of $S_0$ and
$D^{-1}$ satisfies the commutation relation
\[
D^{-1} f = \sum_{0\le j} (-1)^j \frac{d^j f}{du^j} D^{-1-j}
         = \sum_{0\le j<p} (-1)^j \frac{d^j f}{du^j} D^{-1-j}.
\]
Here we find that the radical is the principal one-sided ideal generated by
$D^{-1}$, and the quotient is again a field.

The above calculation carries over to the untwisted global case, showing
that $\coh {\cal S}$ is a maximal order on a ruled surface over the image
$C^{(p)}$ of $C_0=C_1=C$ under Frobenius.  In fact, our computation allows
us to identify the corresponding rank 2 vector bundle as
$\sO_{C^{(p)}}\oplus \omega_{C^{(p)}}(c')$ where $c'$ is the norm of the
conductor, corresponding to a curve $\hat{Q}'$ embedded as the appropriate
double section.  In this case, twisting is relatively easy to deal with,
since the center is preserved by the automorphisms used under gluing, and
we can compute the effect of nontrivial twists using Proposition
\ref{prop:twist_differential}.  (This will be trickier to deal with in
general, when we need a more subtle notion of center.)

Note that when the conductor is nontrivial, we should really consider the
curve $Q$ rather than $\hat{Q}$, which leads us to consider the curve $Q'$
obtained by adding fibers with the appropriate multiplicities to
$\hat{Q}'$.  We can then verify that $Q'$ is an anticanonical curve on the
center.

\section{The difference cases}

A useful observation for the non-differential cases is that since
$\overline{H}$ and $\overline{S}$ are Morita equivalent, they have the same
centers (appropriately defined), and one is maximal iff the other is
maximal.  This means that we can feel free to work with whichever of the
two is most convenient at any given point in the argument.  In particular,
in the difference case, it will be more convenient to work with
$\overline{H}$ for the bulk of the argument.  As in the differential case,
a key step is identifying the kernel of the action on $R$.

We first suppose that $S_0$ and $S_1$ are complete dvrs (and we are in the
difference case, so $R$ is separable over both $S_0$ and $S_1$).  Then $R$
is an order in either a dvr or a sum of two copies of $S_0$.  In the first
case, $R$ inherits a valuation from $S_0$, and we take $\xi$ to be an
element of smallest absolute value that generates $R$ over $S_0=k[[x_0]]$.
Valuation considerations then imply that $\xi$ also generates $R$ over
$S_1=k[[x_1]]$.  In the disconnected case, we similarly suppose that $\xi$
vanishes in one summand, and among such elements has maximal absolute value
in the other summand; this again forces $\xi$ to generate over both $S_0$
and $S_1$.  In either case, we use this element to define the requisite
twisted derivations: for $f,g\in S_0$, $\nu_0(f+g\xi)=g$, and similarly for
$\nu_1$.

In the differential case, the fact that the conductor is normalizing was
crucial.  In the differential case, it turns out that we get a nontrivial
normalizer not only when the conductor is nontrivial but also when there is
suitable ramification.  Associate a positive integer $m$ to the
configuration $(R,S_0,S_1)$ as follows.  If $R\ne S_0S_1$ (which in
particular happens if $S_0=S_1$), let $m=1$; otherwise, if $r>1$ is the
order of $s_0s_1$, we take $m$ to be the $p'$-part of $r$, unless $r$ is a
power of $p$, in which case $m=p$.

\begin{lem}\label{lem:normalizing_element_of_order_m}
  The element $(x_1/x_0)^m\in R^*$.
\end{lem}

\begin{proof}
  It will suffice to show that $\nu_0(x_1^m)\in x_0^m S_0$, since then
  valuation considerations force $x_1^m-\nu_0(x_1^m)\epsilon\in x_0^m S_0$
  as well, so that $(x_1/x_0)^m\in R$, with the unit property following by
  symmetry.  If $R\ne S_0S_1$, so $m=1$, then $x_1$ cannot generate $R$
  over $S_0$, and thus neither can $\nu_0(x_1)\epsilon =
  x_1-(x_1-\nu_0(x_1)\epsilon)\in x_1+S_0$.  But this implies that
  $\nu_0(x_1)\in S_0$ is not a unit, and thus $\nu_0(x_1)\in x_0 S_0$ as
  required.

  In the remaining cases, since $\nu_0$ is a twisted derivation, we may
  write
  \[
  \nu_0(x_1^m) = \nu_0(x_1)\sum_{0\le i<m} x_1^i s_0(x_1)^{m-1-i}.
  \]
  Since $s_0$ preserves the valuation(s) of $x_1$ (it either preserves the
  unique valuation of the normalization of $R$ or swaps the two valuations;
  since $x_1$ is $s_1$-invariant and $s_1$ has the same effect on the
  valuations, its valuation is $s_0$-invariant), we find that
  $s_0(x_1)/x_1\in \widetilde{R}$, and thus we may compute
  \[
  \nu_0(x_1^m) = \nu_0(x_1) x_1^{m-1} \sum_{0\le i<m}
  (s_0(x_1)/x_1)^{m-1-i}.
  \]
  In particular, if
  \[
  \sum_{0\le i<m} (s_0(x_1)/x_1)^{m-1-i}
  \]
  is in the radical of $\widetilde{R}$, then $x_1^{1-m}\nu_0(x_1^m)$ has
  positive valuation, so that $x_0^{1-m}\nu_0(x_1^m)$ has positive
  valuation as required.

  If $\Spec(R)$ is disconnected, then $x_1$ is a uniformizer in each direct
  summand, and thus modulo the radical, $s_0(x_1)/x_1=s_0(s_1(x_1))/x_1$ is
  given by the action of $s_0s_1$ on the two tangent vectors.  Since
  $S_0\ne S_1$, we have $s_0s_1\ne 1$; if it has order a power of $p$, then
  it fixes the tangent vectors, and otherwise it acts as inverse $m$-th
  roots of unity, and in either case the sum vanishes.

  If $\Spec(R)$ is connected and $\pi$ is a uniformizer of $\widetilde{R}$,
  then $s_0s_1(\pi)/\pi$ is congruent to $1$ since both $s_0$ and $s_1$ act
  as $-1$ on the tangent vector.  It follows that the order of $s_0s_1$ is
  a power of $p$, and again the claim follows.
\end{proof}

\begin{cor}
  The elements $x_0^m$ and $x_1^m$ normalize $H$.
\end{cor}

\begin{proof}
  Since $x_1^m$ commutes with $N_1$, it suffices to show that
  $x_1^m N_0\in N_0 H$.  But this follows by writing
  \[
  x_1^m N_0 = (x_1/x_0)^m x_0^m N_0 = (x_1/x_0)^m N_0 x_0^m = (x_1/x_0)^m
  N_0 (x_0/x_1)^m x_1^m.
  \]
\end{proof}

Define a sequence of elements $x_0^{(l)}$, $1\le l$, by $x_0^{(1)}:=x_0$,
$x_0^{(l)}:=s_{(l-1)\bmod 2}(x_0^{(l-1)})$, and similarly for $\xi^{(l)}$.
Note that all of these elements have the same valuation(s) as $x_0$.

\begin{cor}\label{cor:nux0l_normalizes}
  For any integer $1\le l$, we have
  \[
  \nu_{l\bmod 2}(x_0^{(l)}) H = H \nu_{l\bmod 2}(x_0^{(l)}).
  \]
\end{cor}

\begin{proof}
  It will suffice to show that $\nu_{l\bmod 2}(x_0^{(l)})$, when nonzero,
  has valuation(s) a multiple of that of $x_{l\bmod 2}^m$.  This is
  immediate if $m=1$, so we may assume we have trivial conductor; i.e.,
  $R=S_0S_1$.  Since we are in the difference case, we can write
  \[
  \nu_{l\bmod 2}(x_0^{(l)}) = \frac{x_0^{(l)}-s_{l\bmod
      2}(x_0^{(l)})}{\xi-s_{l\bmod 2}(\xi)}
  \]
  Now, $x_0^{(l)}$ is invariant under the involution $\cdots
  s_1s_0s_1\cdots $ of length $2l-1$, and thus we can rewrite the second
  term as $(s_0s_1)^{\pm l}(x_0^{(l)})$ where the sign depends on the
  parity of $l$.  That the ratio has valuation as required then follows
  from the Hasse-Arf theorem.
\end{proof}

To understand the kernel of the action, it will be helpful to consider a
larger algebra of ($k$-linear) operators.  For each integer $l\ge 0$, the
subspace of $H$ with length bounded by $l$ has rank $2l+1$, and thus we
would naively expect that there should be a unique (up to scale) such
operator that annihilates the $2l$ elements $1,\dots,x_0^{l-1}$ and
$\xi,\dots,x_0^{l-1}\xi$.  The nonuniqueness could then be eliminated (at
the cost of making the section of $H$ meromorphic) by insisting that the
image of $x_0^l$ be 1.  This of course fails once $l$ is large enough, as
the image of $H$ in $\End_k(R)$ has finite rank, but either via
experimentation for small $l$ or by computing appropriate determinants, one
finds that the resulting action can be computed explicitly, and extends to
arbitrary $l$.  Thus for each integer $l\ge 0$, let $L_l$ denote the unique
continuous $k$-linear operator on $R$ such that
\begin{align}
  L_l \frac{1}{1-tx_0} &= \prod_{1\le i\le l+1} \frac{t}{1-t x_0^{(i)}}\\
  L_l \frac{\xi}{1-tx_0} &= \xi^{(l)}\prod_{1\le i\le l+1} \frac{t}{1-t x_0^{(i)}}.
\end{align}
This should be interpreted as saying that if we expand both sides as formal
power series in $t$, then each coefficient transforms as stated.  This is
clearly well-defined, since it specifies $L_l$ on a (topological) basis of
$R$.  Moreover, it follows by triangularity that the operators
$L_0,N_0L_0,L_1,N_0L_1,\dots$ form a topological basis of the algebra of
continuous $k$-linear endomorphisms of $R$.

By our original motivation, we expect the operators $L_l$ to be
representable as meromorphic sections of $H$, at least for $l$ sufficiently
small.  Such an expression is an immediate consequence of the following
recurrence.

\begin{lem}
  These operators satisfy the recurrence
  \[
  \nu_l(x_0^{(l)}) L_l
  =
  (N_l-\frac{\nu_l(\xi^{(l)})}{\nu_{l-1}(\xi^{(l-1)})}N_{l-1})
  L_{l-1}.
  \]
\end{lem}

\begin{proof}
  We need to show that applying the right-hand side to $(1-t x_0)^{-1}$ and
  $\xi(1-t x_0)^{-1}$ gives the correct result.  The key observation is
  that
  \[
  \prod_{1\le i\le l} \frac{t}{1-t x_0^{(i)}}
  \]
  is $s_{l-1}$-invariant and thus annihilated by $\nu_{l-1}$, and similarly
  \[
  \prod_{1\le i\le l-1} \frac{t}{1-t x_0^{(i)}}
  \]
  is $s_l$-invariant.  Since $\nu_l$ is an $s_l$-twisted derivation, we
  compute
  \begin{align}
  \nu_l\biggl(\prod_{1\le i\le l} \frac{t}{1-t x_0^{(i)}}\biggr)
  &=
  \nu_l\biggl(\frac{t}{1-t x_0^{(l)}}\biggr)
  \prod_{1\le i\le l-1} \frac{t}{1-t x_0^{(i)}}\notag\\
  &=
  \frac{1}{\xi-s_l(\xi)}
  \biggl(\frac{t}{1-t x_0^{(l)}}-\frac{t}{1-t x_0^{(l+1)}}\biggr)
  \prod_{1\le i\le l-1} \frac{t}{1-t x_0^{(i)}}
  \notag\\
  &=
  \nu_l(x_0^{(l)})
  \prod_{1\le i\le l+1} \frac{t}{1-t x_0^{(i)}},
  \end{align}
  so that both sides indeed have the same action on $(1-tx_0)^{-1}$.
  Since
  \[
  (N_l-\frac{\nu_l(\xi^{(l)})}{\nu_{l-1}(\xi^{(l-1)})}N_{l-1})\xi^{(l-1)} =
  \xi^{(l)}N_l -
  s_{l-1}(\xi^{(l-1)})\frac{\nu_l(\xi^{(l)})}{\nu_{l-1}(\xi^{(l-1)})}N_{l-1},
  \]
  the actions on $\xi (1-tx_0)^{-1}$ also agree.
\end{proof}

\begin{rem}
  Since $\xi^{(l)}$ satisfies the same valuation conditions as $\xi$ (apart
  from possibly swapping the two components when $R$ is not a domain), we
  find that $\nu_l(\xi^{(l)})$ is a unit.
\end{rem}

Although this only gives a meromorphic expression for $L_l$, we can control
the poles.

\begin{cor}
  For any integer $l\ge 0$, there is an element of
  \[
  H_{\le s_{l-1}\cdots s_0}+H_{\le s_l\cdots s_1}
  \]
  with unit leading coefficients that acts as the operator
  \[
  \Bigl(\prod_{1\le i\le l} \nu_i(x_0^{(i)})\Bigr) L_l.
  \]
\end{cor}

\begin{proof}
  We have the operator identity
  \[
  \Bigl(\prod_{1\le i\le l} \nu_i(x_0^{(i)})\Bigr) L_l
  =
  \Bigl(\prod_{1\le i\le l-1} \nu_i(x_0^{(i)})\Bigr)
  (N_l-\frac{\nu_l(\xi^{(l)})}{\nu_{l-1}(\xi^{(l-1)})}N_{l-1})
  L_{l-1}.
  \]
  The coefficient normalizes $H$ by Corollary \ref{cor:nux0l_normalizes},
  and thus we can move it to the right to obtain an expression of the form
  \[
  \Bigl(\prod_{1\le i\le l} \nu_i(x_0^{(i)})\Bigr) L_l
  =
  (\alpha_l N_l + \beta_l N_{l-1} + \gamma_l)
  \Bigl(\prod_{1\le i\le l-1} \nu_i(x_0^{(i)})\Bigr) L_{l-1},
  \]
  with $\alpha_l,\beta_l\in R^*$, $\gamma_l\in R$.  The claim follows
  by induction.
\end{proof}
  
\begin{rem}
  This expression is uniquely determined unless $\nu_i(x_0^{(i)})=0$ for
  some $1\le i<l$.
\end{rem}  

\begin{cor}
  Suppose that $S_0$ has degree $r<\infty$ over $S_0\cap S_1$.  Then
  the image of $H$ in $\End_{S_0\cap S_1}(R)$ is the same as the images of
  the Bruhat intervals corresponding to the two elements of length $r$,
  both of which inject in $\End_{S_0\cap S_1}(R)$.
\end{cor}

\begin{proof}
  Since $r$ is the order of $s_0s_1$, we find that $\nu_i(x_0^{(i)})=0$ iff
  $i$ is a multiple of $r$, and thus there is an element in the union of
  the Bruhat intervals with unit leading coefficients that acts as
  \[
  \Bigl(\prod_{1\le i\le r} \nu_i(x_0^{(i)})\Bigr) L_r = 0.
  \]
  We can use the resulting element of the kernel to reduce any
  word of length $>r$ to a linear combination of shorter words with the
  same image, and thus the Bruhat intervals of length $r$ both span the
  image of $H$.

  For injectivity, we observe that the elements
  \[
  (\prod_{1\le i\le l} \nu_i(x_0^{(i)})) L_l,
  (\prod_{1\le i\le l} \nu_i(x_0^{(i)})) (N_{l+1}-\nu_{l+1}(\xi^{(l+1)})) L_l
  \]
  for $0\le l<r$ form a basis of the relevant Bruhat interval, so it
  suffices to show that their images are linearly independent.  But this
  follows by considering how they act on $1,\xi,x_0,x_0\xi,\dots$.
\end{proof}

\begin{rem}
Although we did the calculations assuming $S_0$ and $S_1$ are complete
local rings, the corresponding result before completion follows
immediately.
\end{rem}

Now, let us consider the nonintegral case; that is, $S_0$ and $S_1$ are
regular with $S_0\cap S_1$ a dvr.  Assuming again that $S_0\cap S_1$ is
complete, the normalization of $R$ contains a number of idempotents, all of
which must in fact be contained in $R$; indeed, we can represent any
idempotent of the normalization as a product of $s_0$- and $s_1$-invariant
idempotents.  Since $S_0\cap S_1$ is local, the dihedral group acts
transitively on the idempotents of $R$, and if we multiply an idempotent
which is not $s_i$-invariant by general element of $S_i$, we obtain a
general element of that summand of the normalization.  In other words, if
$R$ is not integral, then it is normal.  Furthermore, the automorphism
$s_0s_1$ permutes the idempotents; if it acts as an $r'$-cycle, then the
action of $D_{2r}$ on $R$ can be obtained by induction from an action of
$D_{2r/r'}$ on either a complete dvr or a sum of two complete dvrs swapped
by the reflections.  In particular, the nonintegral case is always obtained
by such an induction from a normal instance $R_0$ of the integral case.

We can use this structure to construct a generator of the kernel in
general.  If $R$ has $2r'$ idempotents (so $R_0$ has two summands), then
$H$ is equal to the twisted group algebra (the coefficient of $N_i$ in
$s_i$ is a unit), and thus the kernel is generated by the difference of the
reduced words of length $r$.  Thus suppose that $R$ has $r'$ idempotents.
Then precisely two of those idempotents are invariant under some reflection
(which will be the same reflection if $r'$ is even, and different
otherwise), giving rise to a pair of operators $\iota_i N_{j(i)} \iota_i$
where $\iota_1$, $\iota_2$ are the two idempotents and $j(i)$ is the index
of the reflection fixing $\iota_i$.  Note that each operator acts on the
corresponding summand as an operator of the form $\nu$.  Using this, we
obtain a collection of $2r'$ operators from the two reduced words of length
$r'$ by replacing single reflections by the corresponding operators
$\iota_i N_{j(i)}\iota_i$.  (We can either replace any reflection in this
way, or can only replace half the reflections, but can do so in two
different ways.)  Each resulting operator annihilates all but one summand
of $R$ and maps into the image of that summand under either word of length
$r'$.  We have two operators for each such pair of summands, which
essentially act as the $N_0$ and $N_1$ on the corresponding copy of $R_0$.
We can thus construct a sum of words of length $\le r/r'$ in such elements
with unit (in $R_0$) leading coefficients acting trivially, giving rise to
a corresponding sum of words of length $\le r$ in $H$ inside the kernel.
The resulting expressions no longer have unit leading coefficients, but
their leading coefficients are units inside the appropriate summands, and
thus adding up $r'$ such expressions gives the desired generator of the
kernel.  (Uniqueness in either case follows by base changing to the field
of fractions of $S_0\cap S_1$.)

Thus in any difference case, we have constructed a generator of the kernel
of the action of $H$, which we denote by ${\bf k}$.  As an element of the
twisted group algebra, we have an expression
\[
{\bf k} = f ((\cdots s_1s_0)-(\cdots s_0s_1))
\]
for some element $f$ of the base change to the field of fractions.  Unlike
in the differential case, this element does not even commute with
multiplication by $R$.  This at least partly arises from the fact that
``center'' is a somewhat ill-defined concept in the case of a $\Z$-algebra
like $\overline{H}$.  The point is that in order to interpret a $\Z$-algebra as
coming from a graded algebra, we need to choose an isomorphism between the
algebra and a suitable shift, and such an isomorphism need not be unique!
In our case, although ${\bf k}$ does not {\em commute} with $R$, it does
act as an automorphism of $R$; indeed, if $\omega$ denotes the longest
element of $D_{2r}$, then we find that ${\bf k}x = \omega(x) {\bf k}$ for
all $x\in R$.  We also find
\[
\omega(\nu_i(\omega(f)))
=
\nu_{i+r}(\omega(\xi))^{-1} \nu_{i+r}(f)
\]
and thus this automorphism of $R$ induces an (order 2) automorphism of $H$
by
\[
\omega(N_i) = \nu_{i+r}(\omega(\xi))^{-1} N_{i+r}.
\]
We may thus define the ``quasi-center'' of $\overline{H}$ to consist of
those morphisms $\phi$ of degree $dr$ such that $\phi \psi = (\omega^d\psi)
\phi$.  This agrees with the naive notion of center in degrees a multiple
of $2r$, but is better-behaved for odd multiples of $r$, and globalizes
just as easily, at least in the untwisted case.  It is then straightforward
to determine which elements of the rank 2 module $R(\cdots s_1s_0)+R {\bf
  k}$ are quasi-central.

\begin{lem}
  An element $f(\cdots s_1s_0)+g(\cdots s_0s_1)$ of degree $r$
  is quasi-central iff $f=s_0(g)=s_0s_1(f)$.
\end{lem}

Any quasi-central element of degree $r$ must have this form, and thus to
determine the degree $r$ elements of the quasi-center, it suffices to
understand the corresponding conditions on $f$ and $g$.  It follows from
our above calculations that there is a fractional ideal $I$ over $R$ such
that $f(\cdots s_1s_0)+g(\cdots s_0s_1)$ is in $H$ iff $f,g\in I$ and
$f+g\in R$, and thus we need to understand the elements $f\in R^{\langle
  s_0s_1\rangle}$ such that $f\in I$ and $f+s_0(f)\in R$.

The fractional ideal $I$ can be expressed as a norm (of the inverse of the
conductor of $R$ inside its normalization) times a product of inverse
differents (of the normalization of $R$ relative to the $r$ subrings fixed
by involutions).  The product of inverse differents is itself an inverse
different, due to the following fact.

\begin{prop}
  Let $K$ be the field of fractions of a Dedekind ring, and let $L/K$ be a
  Galois \'etale $K$-algebra with Galois group $D_{2r}$, let
  $F_1,\dots,F_r$ be the subalgebras fixed by the $r$ reflections, and let
  $E$ be the subalgebra fixed by $C_r$.  Then the different of $E/K$ can be
  computed in terms of the differents of $L/F_i$:
  \[
  d_{E/K} = \prod_{1\le i\le r} d_{L/F_i}.
  \]
\end{prop}

\begin{proof}
  This is equivalent to the corresponding claim for discriminants, which in
  turn follows from the conductor-discriminant formula and basic properties
  of the Artin conductor.
\end{proof}

\begin{cor}
  The quasi-central elements of degree $r$ in $H$ are naturally identified
  with the elements of the inverse different of a quadratic extension of
  $S_0\cap S_1$.
\end{cor}

\begin{proof}
  Letting $K,L,\dots$ be the \'etale algebras corresponding to the action
  of $D_{2r}$ on $R$, we find that the coefficient ideal $I$ is the product
  of the inverse different $d_{E/K}^{-1}$ by the inverse of the norm $C'$ of
  the conductor of $R$ inside $O_K$.  We thus find that the quasi-central
  elements may be identified with the space
  \[
  \{x:x\in C^{\prime{-}1} d_{E/K}^{-1}|x+s_0(x)\in O_E\}.
  \]
  But this is the inverse different of the suborder of conductor $C'$ in
  $O_E$, which can be expressed as
  \[
  \{x:x\in O_E|x-s_0(x)\in C' d_{E/K}\}.
  \]
\end{proof}

\begin{cor}
  The rank 2 free module $R(\cdots s_0s_1)+R{\bf k}$ has an $R$-basis consisting
  of quasi-central elements.
\end{cor}

\begin{proof}
  Let $R'$ be the above quadratic order over $O_K$, so that the
  quasi-center may be described as the space of elements
  \[
  s_0(f)(\cdots s_0s_1)+f(\cdots s_1s_0)
  \]
  with $f$ in the inverse different of $R'$.  If $\xi'$ generates $R'$ over
  $O_K$, then this has an $O_K$-basis of the form
  \[
  s_0(f_0)(\cdots s_0s_1)+f_0(\cdots s_1s_0),
  s_0(\xi' f_0)(\cdots s_0s_1)+\xi' f_0(\cdots s_1s_0),
  \]
  where $f_0$ generates the inverse different.  It suffices to show that
  the corresponding $R'$-module contains $(\cdots s_0s_1)$, since then it
  also contains
  \[
  f_0((\cdots s_0s_1)-(\cdots s_1s_0))\in R^*{\bf k}.
  \]
  The $R'$-module clearly contains
  \[
  (s_0(\xi'f_0)-\xi's_0(f_0)) (\cdots s_0s_1);
  \]
  since
  \[
  s_0(\xi'f_0)-\xi's_0(f_0)
  =
  s_0((\xi'-s_0(\xi'))f_0)
  \]
  and $(\xi'-s_0(\xi'))f_0$ is a unit, the result follows.
\end{proof}

\begin{cor}
  The quasi-center of $\overline{H}$ is generated in degree $r$, and
  $\overline{H}$ is locally free over its quasi-center.
\end{cor}

\begin{proof}
  The space of degree $r$ quasi-central elements has an $S_0\cap S_1$-basis
  $Z_0$, $Z_1$ such that $Z_1$ is a unit multiple of ${\bf k}$ and $Z_0$
  acts (as an operator) on $R$ as $\omega$.  (Indeed, we may take $Z_1 = f
  ((\cdots s_0s_1)-(\cdots s_1s_0))$ where $f$ is an anti-invariant element
  of $I$, and then the other basis element in any extension to a basis will
  act as a unit multiple of $\omega$).  If $\phi$ is quasi-central of
  degree $dr$, then $\phi$ acts on $R$ as an $S_0\cap S_1$-multiple of
  $\omega^d$, and thus has the same operator action as some $S_0\cap
  S_1$-multiple of $Z_0^d$.  Subtracting this multiple to make the operator
  action trivial makes $\phi$ a left multiple of $Z_1$, and right dividing
  by $Z_1$ gives a quasi-central element of degree $(d-1)r$, so that the
  quasi-center of $\overline{H}$ is generated by $Z_0$, $Z_1$.

  Similarly, if $\phi$ is a morphism of degree $dr+(r-1)$, then there is a
  unique element $\mu_0$ of degree $r-1$ acting on $R$ as $\omega^d \phi$,
  and subtracting $Z_0^d$ times this element gives $Z_1$ times a morphism
  of degree $(d-1)r+(r-1)$, inducing a unique expression of the form
  \[
  \phi = \sum_{0\le i\le d} Z_1^i Z_0^{d-i} \mu_i.
  \]
  (Here we have used $\omega$ to identify the different $\Hom$ spaces of
  degree $r$ so that it makes sense to take a product of elements $Z_i$,
  and such products commute as expected.)  Local freeness follows
  immediately.
\end{proof}

In particular, ${\overline H}$ indeed corresponds to a (locally free) order
over the appropriate ruled surface $Z$, and to show maximality it remains
only to check maximality in codimension 1.  This is again an Azumaya
algebra over a large open subset of $Z$; indeed, if $R$ is regular and
unramified over $S_0\cap S_1$, then the fiber over any point not ``at
infinity'' (i.e., not on the image of $\hat{Q}$) is a ``dihedral'' algebra,
i.e., the analogue of a cyclic algebra in which the cyclic group is
replaced by a dihedral group.  The same argument as in the cyclic case
tells us that the fiber is central simple of rank $4r^2$.  (This is $4$
times what we had in the differential case since we are working with
${\overline H}$ rather than ${\overline S}$.)

It remains to consider fibers over ramified points and components of
$\hat{Q}$.  The first case reduces to the case that $S_0\cap S_1$ is a
complete dvr.  If both $S_0$ and $S_1$ are integral with generators $x_0$,
$x_1$, then Lemma \ref{lem:normalizing_element_of_order_m} above gives us a
minimal positive integer $m$ such that $(x_1/x_0)^m\in R^*$, and we find as
in the differential case that $H x_0^m H = x_0^m H$ is contained in the
radical of the localization of $H$.  The quotient by this ideal is central
simple (a dihedral algebra of rank $4m^2$), and thus we have maximality as
required.  (Note that if $m=r$, then we still obtain an Azumaya algebra
away from the intersection of $\hat{Q}$ with the given fiber.)  If at least
one of $S_0$ or $S_1$ fails to be integral, then $R$ must be integrally
closed, and one finds that the radical of the localization of $H$ is
generated by the (principal) radical of $S_0$, with quotient central simple
of degree $4c^2$ where $c$ is the number of irreducible components of
$\Spec(R)$.  (So again, if $c=r$, we obtain an Azumaya algebra away from
$\hat{Q}$.)

For the completion along a component of $\hat{Q}$, we may first base change
to the function field $K$ of the base curve.  Then $S_0$ and $S_1$ are
fields, and $R$ is either a field or a sum of two fields, and in each case
the completion along $\hat{Q}$ is the completion along the two-sided
homogeneous ideal $\overline{H}t^{2r}\overline{H} = t^{2r}\overline{H}$ (which may not be
prime).  (The factor of 2 comes from the fact that $t^r$ is not
quasi-central in general!)  We find as before that $\overline{H} t^2 \overline{H} =
t^2\overline{H}$ is contained in the radical of this completion.  The $2r$-th
Veronese of the quotient is precisely $R t^{2r}$, so that the quotient is a
field whenever $R$ is a field, giving maximality in that case.  When $R$ is
a sum of two fields, we note that it suffices to prove the corresponding
maximality for the spherical algebra, which is a complete Ore ring of the
form $S_0\langle\langle \rho\rangle\rangle$ where $\rho$ is an order $r$
automorphism depending on the choice of component of $\hat{Q}$.  The
radical is then the principal ideal generated by $\rho$, and the quotient
is the field $S_0$ as required.

We thus conclude that for a quasi-ruled surface of difference type such
that $s_0s_1$ has order $r$, $\overline{H}$ is (in large degree) isomorphic to
the algebra of sections of a maximal order of rank $4r^2$ over the
appropriate commutative ruled surface, and this continues to hold globally
for an {\em untwisted} surface.  Twisting is somewhat tricky in this
instance, as the quasi-center is not actually invariant under twisting.
The existence of {\em some} global surface is easy enough, as we can simply
take the second Veronese of the quasi-center (i.e., the elements of degree
a multiple of $2r$), as then the notion of quasi-center agrees with the
naive notion of center and is thus invariant under twisting.  This
immediately tells us that there is at least a global {\em conic} bundle
over which the noncommutative surface is a maximal order.  We could in
principle pin this down by carefully looking at how the twisting occurs,
but the precise identification will in any event follow once we have a
sufficient understanding of noncommutative elementary transformations.
Indeed, we will see below that any twisted quasi-ruled surface of
difference type can be related to an untwisted such surface by a sequence
of (noncommutative) elementary transformations in smooth points of
$\hat{Q}$, and thus (by Theorem \ref{thm:blowup_of_order}) its center is
related to the untwisted center by a corresponding sequence of commutative
elementary transformations.

\section{The hybrid case}

In the remaining (hybrid) case, in which one of the two maps (say
$\Spec(R)\to \Spec(S_0)$) is inseparable, we have a similar result, in
which $\overline{S}$ corresponds to a quaternion algebra over an appropriate
ruled surface.  Note that the argument from the difference case still gives
an element ${\bf k}$ generating the kernel, so the only missing step is
showing that the span of $s_1$ and ${\bf k}$ is generated by quasi-central
elements.  This can in principle be done by a direct computation, but we
also note that there is an argument via semicontinuity.  This is slightly
tricky, as we do not actually have suitable semicontinuity locally, but we
can rescue this by embedding the desired local situation in a suitable
global configuration.  It is easy to see that any complete local hybrid
case appears inside a global case in which $C_0$ (and thus $C'$) is
isomorphic to $\P^1$, so that $\hat{Q}$ (and thus $C_1$) is hyperelliptic.
We can embed such cases in a larger family by allowing the degree 2 map
$C_0\to \P^1$ to vary (constructing $\overline{Q}$ via the compositum of field
extensions).  In the untwisted case with $\hat{Q}=\overline{Q}$, the global
version of the span of $s_0s_1$ and ${\bf k}$ is flat over this family,
while the $\sO_{\hat{Q}}$-submodule generated by quasi-central elements is
upper semicontinuous (as a $\sO_{C'}$-module) and generically agrees with
the larger module, so always agrees with the larger module.  It thus
follows that any {\em local} configuration with trivial conductor has the
desired basis of quasi-central elements.  To deal with conductors, we
observe that the general configuration ($C_0\cong \P^1$, $C_1$
hyperelliptic, $s_0s_1=s_1s_0$) makes sense in characteristic 0, and there
are characteristic 0 examples with nontrivial conductor.  A similar
semicontinuity argument gives some characteristic 2 examples with
nontrivial conductor, and we can then get more general characteristic 2
examples with nontrivial conductor by suitable limits (e.g., taking a limit
in which several simple cusps coalesce).

\section{Identifying the center}

One reason we did not spend too much effort determining the explicit center
in the twisted case is that there is an alternate approach to computing the
center once we know that we actually {\em have} a maximal order.  The key
idea is that, most easily on the Azumaya locus, we can identify the center
as a moduli space of $0$-dimensional sheaves of degree $r$ on the
noncommutative surface.  Indeed, any point on the center pulls back to a
degree $r^2$ representation of the maximal order corresponding to
$\overline{S}$, which if the point is in the Azumaya locus is the $r$-th power
of an irreducible representation of degree $r$.  Thus in particular the
fiber of the Azumaya locus over a point of $C'$ will correspond to such a
sheaf which is supported over that point.  Using the semiorthogonal
decomposition (Theorem \ref{thm:semiorth_for_quasiruled} below), such a
sheaf determines and is determined by a morphism of the form
\[
\pi_1^* Z_1\to \pi_0^* Z_0\otimes {\cal L}
\]
where $Z_0$, $Z_1$ are the respective pullbacks of the point sheaf on $C'$.
Thus the space of such morphisms may be naturally identified with the fiber
of $\pi_*{\cal L}$ over the given point of $C'$.  Two such morphisms
determine the same object of the derived category iff they are in the same
orbit under the actions of $\Aut(Z_0)$ and $\Aut(Z_1)$, which act as the
corresponding fibers of $\sO_{C_0}^*$ and $\sO_{C_1}^*$ respectively.

There are some technical issues here in general, since not every morphism
will correspond to an injective morphism of sheaves on the noncommutative
surface, and one must also consider the sheaves only up to $S$-equivalence
to obtain a well-behaved moduli space.  The first issue can be resolved by
considering only those morphisms which are isomorphisms of sheaves on
$\hat{Q}$; this gives an open subset which is still large enough to
determine the surface.  The second issue is more subtle, in general, though
only an issue away from the Azumaya locus.

This is, however, enough to let us understand twisting in the difference
setting.  Consider the locus of $C'$ over which $\hat{Q}$ is unramified and
equal to $\overline{Q}$.  Over this locus, the fibers of $\sO_{\hat{Q}}^*$ (over
which $\pi_*{\cal L}^*$ is a torsor) are algebraic tori, and it is
straightforward to see that two points are in the same coset of
$\sO_{C_0}^*\sO_{C_1}^*$ iff their norms down to $\sO_{\hat{Q}'}^*$ are in
the same coset over $\sO_{C'}^*$.  We thus find more generally that this
component of the moduli space of $0$-dimensional sheaves compactifies
(fiberwise) to a $\P^1$-bundle (which will, in fact, be the GIT quotient):
two nonzero elements of an unramified fiber of $\pi_*{\cal L}$ correspond
to the same point of this $\P^1$-bundle iff their norms in the
corresponding fiber of
$N_{\hat{Q}^{\text{unr}}/\hat{Q}^{\prime\text{unr}}}({\cal L})$ are in the
same $\sO_{C'}^*$-orbit.

This description makes sense even on the ramified fibers, with one
important caveat: we do not know in general that the norm on $k(\hat{Q})$
takes $\sO_{\hat{Q}}$ to $\sO_{\hat{Q}'}$.  This holds trivially whenever
$\hat{Q}$ is smooth, and is easy to verify at simple nodes of $\hat{Q}$,
and thus in particular holds as long as $k(\hat{Q})/k(\hat{Q}')$ is tamely
ramified wherever it is totally ramified.  As long as it holds, we can
define norms of line bundles on $\hat{Q}$ by taking norms of their
$1$-cocycles, and find that the norm indeed gives a polynomial function of
degree $r$ from ${\cal L}$ to its norm, compatible with the norms from
$\sO_{C_0}^*$ and $\sO_{C_1}^*$ to $\sO_{C'}^*$.  In other words, the norm
induces a morphism from the family of quotient stacks
$\sO_{\hat{Q}}/\sO_{C_0}^*\sO_{C_1}^*$ to the $\P^1$-bundle on $C'$
associated to $N_{\hat{Q}/\hat{Q}'}({\cal L})$.  Over the locus where
$k(\hat{Q})$ is unramified over $k(\hat{Q}')$, it is easy to see the
restriction to $\sO_{\hat{Q}}^*/\sO_{C_0}^*\sO_{C_1}^*$ makes the image a
coarse moduli space (i.e., that the ratio of the coefficients of the norm
generates the field of invariants); outside the Azumaya locus, this fails
since the resulting sheaves are no longer irreducible and the
$S$-equivalence classes contain multiple isomorphism classes.

In fact, since any local trivialization of ${\cal L}$ can be refined to one
in which the overlaps are contained in the smooth (or unramified) locus, we
can define $N_{\hat{Q}/\hat{Q}'}({\cal L})$ by taking the norm of the
associated cocycle; since this is consistent with the gluing of the center
over the unramified locus, this gives us the desired description of the
global center.  (This strongly suggests that the norm indeed always takes
$\sO_{\hat{Q}}$ to $\sO_{\hat{Q}'}$.)  We can rephrase this in terms of
divisors: if $D$ is a divisor supported on the smooth locus of $\hat{Q}$,
then
\[
N_{\hat{Q}/\hat{Q}'}(\sO_{\hat{Q}}(D))
\cong
\sO_{\hat{Q}'}(D')
\]
where $D'$ is the image of $D$ in $\hat{Q}'$.  (Note that any line bundle
has such a description, which corresponds to a sequence of (noncommutative)
elementary transformations that untwist ${\cal L}$.)

In the differential case, the situation is more subtle, as now the quotient
is no longer even generically the quotient by a torus.  It is still
relatively straightforward to identify two semi-invariants of degree $p$ on
the fibers of $\sO_{\hat{Q}}$, namely the coefficients of the restriction of
$\Nm'$ (as defined in Corollary \ref{cor:definition_of_Nmprime}) to the fiber.
Furthermore, at least over the locus with trivial conductor, these are the
{\em only} semi-invariants of degree $p$.  To see this, it suffices to
consider the restriction to fibers of $\sO_{\hat{Q}}^*$, or equivalently to
$k[z,\epsilon]/(z^p,\epsilon^2)^*$.  That $\sO_{C_0}^*$ acts as the norm
lets us further restrict to the subgroup $1+(k[z]/z^p)\epsilon$, and then
using $\sO_{C_1}^*$-invariance to the quotient of $(k[z]/z^p) dz$ by the
$\F_p$-span of the logarithmic differentials.  These include the elements
$\frac{ay}{1+byz} dz$ for $a,b\in \F_p^*$, $y\in k$, and thus (by taking
suitable $\F_p$-linear linear combinations, essentially a discrete Fourier
transform relative to the group $\F_p^*$) the elements $(y+y^p z^{p-1})dz$
and $y^l z^{l-1} dz$ for $1\le l<p$.  It follows immediately that the field
of invariants is generated by a unique function of degree $p$, giving the
desired result.  (This fails when the conductor is nontrivial just as in
the difference case.)  Again, it follows immediately that $\Nm'$ gives the
desired map on twisting data, and here there is no difficulty in dealing
with the case of nontrivial conductor, since then $\hat{Q}$ has trivial
Picard group.

In the hybrid case, we suppose that $\hat{Q}\to C_0$ is inseparable and
consider the locus of $C'$ over which $C_0$ is \'etale.  The corresponding
fibers of $\sO_{\hat{Q}}$ have the form $k[z]/z^2\oplus k[z]/z^2$, with
$\sO_{C_1}^*$ acting diagonally and $\sO_{C_0}^*$ acting as multiplication
by $k^*\times k^*$.  We again find that there are (only) two invariants of
degree 2, which on $(a_0+a_1z,b_0+b_1z)$ take the values $a_0b_0$ and
$a_0b_1+a_1b_0$.  In particular, just as in the difference case, this is
the norm down to $k(\hat{Q}')$, so that we can deal with twisting in the
same way.  (Again, we expect this norm to take $\sO_{\hat{Q}}$ to
$\sO_{\hat{Q}'}$ even at singular points of $\hat{Q}$.)

\medskip

One reason why twisting was tricky to deal with is that not every line
bundle on the noncommutative surface will be the pullback of a line bundle
on the center, and we cannot expect there to be a natural choice of such a
pullback of degree $r$ relative to the fibration.  In the untwisted case,
there was no difficulty: the above construction associates sheaves to any
point of a fiber of $\sO_{\hat{Q}}^*$, and thus in particular associates a
{\em canonical} family of sheaves to the identity section.  This gives rise
to a corresponding section of the $\P^1$-bundle over $C'$, which is why in
this case we not only can identify the relevant rank 2 vector bundle on
$C'$, but have a natural choice of map from that vector bundle to
$\sO_{C'}$, such that the corresponding section is disjoint from the image
of $\hat{Q}$.


\section{Blowing up}

Since we are not only interested in quasi-ruled surfaces themselves, but
also those surfaces which are merely birationally quasi-ruled, we also want
to understand how to translate a blowup \`a la van den Bergh into an
operation on maximal orders.  There is a folklore result (subject to some
implicit conditions that exclude some quasi-ruled surfaces in finite
characteristic) saying that such a blowup corresponds to blowing up the
corresponding point of the underlying commutative surface and choosing a
maximal order containing the pulled back algebra.  We give a fully general
proof of this below.  In principle, we would also need to understand
blowing down, but this will be dealt with (at least for blowdowns of
blowups of quasi-ruled surfaces) via Theorem
\ref{thm:birationally_qr_is_rationally_qr} below.

Our objective is to prove the following.

\begin{thm}\label{thm:blowup_of_order}
  Let $Z$ be a (commutative) smooth surface, let ${\cal A}$ be an order on
  $Z$, and let ${\cal I}\subset {\cal A}$ be an invertible two-sided ideal
  such that there is an equivalence $\phi:\coh Y\cong \coh({\cal A}/{\cal
    I})$ for a commutative curve $Y$.  Let $x\in Y$ be a point such that
  $\phi(\sO_x)$ has finite injective dimension as a ${\cal A}$-module and
  ${\cal A}$ is locally free at the point $z\in Z$ on which $\phi(\sO_x)$
  is supported.  Then the van den Bergh blowup \cite{VandenBerghM:1998} of
  $\coh {\cal A}$ at $x$ is equivalent to the category $\coh \widetilde{A}$
  where $\widetilde{A}$ is a maximal element of the poset of orders on the
  blowup $\widetilde{Z}$ of $Z$ at $z$ that agree with ${\cal A}$ on the
  complement of the exceptional divisor.
\end{thm}

\begin{rem}
  In particular, if ${\cal A}$ is a maximal order, then so is
  $\widetilde{\cal A}$.
\end{rem}

\begin{proof}
Since ${\cal A}$ is an honest sheaf of algebras, the construction of the
van den Bergh blowup simplifies considerably, and we find that the blowup
may be expressed as the relative $\Proj$ of a Rees algebra:
\[
\widetilde{X}=\Proj({\cal A}[{\frakm}_x\otimes {\cal I}^{-1} t]),
\]
where ${\frakm}_x$ is the maximal ideal corresponding to $x$ and
$t$ is an auxiliary central variable.  On the complement of $z$, we get
\[
\widetilde{X}\cong \Proj({\cal A}[{\cal I}^{-1} t])\cong \coh {\cal A},
\]
so it will suffice to understand what happens over the local ring at $z$.

For each integer $l$, define a maximal ideal ${\frakm}_l\subset {\cal
  A}$ by
\[
{\frakm}_l = {\cal I}^l \otimes_{\cal A} {\frakm}_x \otimes_{\cal A}
{\cal I}^{-l}.
\]
The quotient ${\cal A}/{\frakm}_l$ is supported on $z$ for all $l$, so we
may view these as maximal ideals of ${\cal A}_z$.  The significance of
these maximal ideals is that
\[
({\frakm}_x {\cal I}^{-1})^l
=
{\frakm}_0{\frakm}_{-1}\cdots {\frakm}_{1-l} {\cal I}^{-l}.
\]
Since ${\cal I}$ permutes the maximal ideals, the sequence ${\frakm}_i$ is
periodic of some period $n$, and this lets us twist the $rn$-th Veronese of
the graded algebra:
\[
\widetilde{X}
\cong \Proj({\cal A}[{\frakm}_0\cdots {\frakm}_{1-rn} {\cal I}^{-n} t]
\cong \Proj({\cal A}[{\frakm}_0\cdots {\frakm}_{1-rn} t]).
\]
The latter turns out to have the better behaved center; in particular, we
claim that for a suitable choice of $r$, its center is precisely
\[
\sO_Z[{\frakm}_z t],
\]
or in other words that
\[
({\frakm}_0\cdots {\frakm}_{1-n})^{lr}\cap \sO_Z
=
{\frakm}_z^l
\]
for all $l\ge 0$.  For this claim, it suffices to consider the completion
$\hat{\cal A}_z\cong {\cal A}\otimes_{\sO_Z} \hat\sO_{Z,z}$.

A somewhat different completion of ${\cal A}$ was studied in
\cite{vanGastelM/VandenBerghM:1997}, namely the completion with respect to
the intersection $\cap_{0\le i<n} {\frakm}_{-i}$.  (To be precise, they
studied a ring constructed from the Serre subcategory generated by the $n$
simple modules of ${\cal A}$, but this is Morita equivalent to the given
completion.)  Luckily, this completion is equal to $\hat{\cal A}_z$, due to
the following fact.

\begin{lem}
  Every maximal ideal of ${\cal A}_z$ is of the form ${\frakm}_i$ for some
  $i$.
\end{lem}

\begin{proof}
  Suppose otherwise, and let ${\frakm}'$ be some other maximal ideal lying
  over $z$.  Idempotents in ${\cal A}_z/\rad {\cal A}_z$ lift over
  $\hat{\cal A}_z$, and thus we may find an idempotent $e\in \hat{\cal
    A}_z$ such that $e\in {\frakm}_i$ for all $i$ but $e\notin {\frakm}'$.
  Since central simple algebras do not have nontrivial two-sided ideals,
  the ideal $\hat{\cal A}_z e \hat{\cal A}_z$ contains an ideal of the form
  ${\cal J}\hat{\cal A}_z$ where ${\cal J}$ is a nonzero ideal of
  $\hat\sO_{Z,z}$.  In particular, we find that
  \[
    {\cal J}\hat{\cal A}_z
    \subset
    \hat{\cal A}_z e \hat{\cal A}_z
    \subset
    {\frakm}_0\cdots {\frakm}_{1-l}
  \]
  for all $l\ge 0$, so that
  \[
    {\frakm}_0\cdots {\frakm}_{1-l}
    \supset
    ({\frakm}_z^l+{\cal J}) \hat{\cal A}_z.
  \]
  It follows in particular that
  \[
  \dim(\hat{\cal A}_z/{\frakm}_0\cdots {\frakm}_{1-l})
  \le
  \rank({\cal A}) \dim(\hat\sO_{Z,z}/{\frakm}_z^l+{\cal J})
  =
  O(l).
  \]
  On the other hand, it was shown in \cite[Cor.~5.2.4]{VandenBerghM:1998}
  that this quotient has length $l(l+1)/2$, giving a contradiction.
\end{proof}

It thus follows from \cite{vanGastelM/VandenBerghM:1997} that the
completion $\hat{\cal A}_z$ is Morita equivalent to an algebra $A$ of the
following form.  If $A=R$ is local (so $\hat{\cal A}_z$ is local), then $R$
has the form $k\langle x,y\rangle/\langle \phi\rangle$ where the relation
$\phi\in {\frakm}_R^2$ and its image in ${\frakm}_R^2/{\frakm}_R^3$ is a
rank 2 element of $({\frakm}_R/{\frakm}_R^2)^{\otimes 2}$, and $Y$ has the
form $\Spec(R/UR)$ where $U\in {\frakm}_R$ is a normalizing element such
that $[R,R]\subset UR$.  (Note that per
\cite[Prop.~7.1]{vanGastelM/VandenBerghM:1997}, the commutator $[x,y]$ is
itself a normalizing element.)  More generally, if $A$ (thus ${\cal A}_z$)
has $n$ maximal ideals, then there is a pair $(R,U)$ as above with $U\notin
{\frakm}_R^2$ such that $A$ is isomorphic to the subalgebra of $\Mat_n(R)$
in which the coefficients above the diagonal are multiples of $U$.  In this
case, the curve $Y$ is cut out by the normalizing element $N$ which has
$1$s just below the diagonal, $U$ in the upper right corner, and $0$
everywhere else.  We note that the center of $A$ consists precisely of the
matrices $z I$ with $z\in Z(R)$, so that that $A$ is free over its center
iff $R$ is free over its center.  The maximal ideals of $A$ are
${\frakm}_1$,\dots,${\frakm}_n$ where ${\frakm}_i$ is cut out by the
condition that $a_{ii}\in {\frakm}_R$, and these satisfy $N {\frakm}_i =
{\frakm}_{i+1} N$, with subscripts interpreted cyclically.

Note that in our case, the algebra $A$ satisfies the additional condition
that it is free as a module over its center.  Moreover, the claim about how
the center of $\hat{\cal A}_z$ meets the products of maximal ideals reduces
to the corresponding claim about the center of $A$, which we refine to
\[
({\frakm}_n\cdots {\frakm}_1)^{lr}
\cap
Z(A)
=
{\frakm}_{Z(A)}^l
\]
where $r$ is the square root of the rank of $R$ over its center.
(Equivalently, $\rank(A) = r^2n^2$.)  An easy induction shows that
\[
{\frakm}_n {\frakm}_{n-1}\cdots {\frakm}_1\subset {\frakm}_R\Mat_n(R),
\]
or in other words that given any sequence of elements of the appropriate
ideals, their images in $\Mat_n(k)$ multiply to 0.  The codimension of this
product is known (again by \cite[Cor.~5.2.4]{VandenBerghM:1998}), so that
one immediately concludes
\[
  {\frakm}_n {\frakm}_{n-1}\cdots {\frakm}_1
  =
  A\cap {\frakm}_R\Mat_n(R),
\]
and similarly
\[
  ({\frakm}_n {\frakm}_{n-1}\cdots {\frakm}_1)^l
  =
  A\cap {\frakm}_R^l\Mat_n(R).
\]
(These identities are implicit in \cite[Prop.~5.2.2]{VandenBerghM:1998}.)
Since the center of $A$ is the center of the diagonal copy of $R$, we
conclude that
\[
({\frakm}_n {\frakm}_{n-1}\cdots {\frakm}_1)^{lr}
\cap
Z(A)
=
{\frakm}_R^{lr}\cap Z(R),
\]
so that it remains only to show that
\[
  {\frakm}_R^{lr}\cap Z(R)
  =
  {\frakm}_{Z(R)}^l.
\]

To finish the argument, we will need a couple of lemmas about orders over
dvrs with sufficiently nice residue field.

\begin{lem}
  Let $R$ be a dvr with residue field $k$ and field of fractions $K$, and
  suppose that $k$ has trivial Brauer group.  Let $A$ be an $R$-order in a
  central simple $K$-algebra of degree $r$.  Then any simple $A$-module has
  dimension at most $r$ as a $k$-space.
\end{lem}

\begin{proof}
  We may as well pass to the completion, as this has no effect on the
  simple modules.  Suppose first that $A$ is a maximal order.  Then
  $A/\rad(A)\cong \Mat_d(l)$ where $l/k$ is a field extension of degree
  $r/d$.  If $A\otimes_R K$ is a division ring, then this follows from
  \cite[Thm.~14.3]{ReinerI:2003} with $d=\sqrt{r}$; this result was only
  stated for the case of finite residue field, but the proof only used
  triviality of $\Br(k)$.  More generally, $A$ is a matrix ring over a
  maximal order in a division ring with center $K$, and thus the result
  follows in general from \cite[Thm.~17.3]{ReinerI:2003}.  Thus in the
  maximal order case, $A$ has a unique simple module isomorphic to $l^d$,
  and thus of dimension $r$ over $k$.

  In general, there exists a maximal order $B$ containing $A$.  Then
  $\rad(B)\cap A\subset \rad(A)$, and $A/(\rad(B)\cap A)\subset
  B/\rad(B)\cong \Mat_d(l)$.  It follows that $A/\rad(A)$ is a sum of
  matrix algebras of the form $\Mat_{d'}(l')$ where $l'/l$ is a field
  extension such that $d'[l':l]\le d$.  Each summand corresponds to a
  simple module with $k$-dimension $d'[l':k]\le d[l:k]=r$ as required.
\end{proof}

\begin{rem}
  To see that $\rad(B)\cap A$ is contained in the radical of $A$, observe
  that this is a two-sided ideal, and thus it suffices to show that $1-z$
  is a unit for all $z\in \rad(B)\cap A$.  But this follows from the fact
  that $\lim_{n\to\infty} z^n=0$ relative to the natural topology in $B$.
\end{rem}

\begin{cor}
  With $R$ and $A$ as above, the left $A$-module $A\otimes_R k$ has length
  at least $r$ (with equality if $A$ is maximal).
\end{cor}

\begin{proof}
  Indeed, $A\otimes_R k$ has dimension $r^2$, while each simple constituent
  contributes at most $r$ to the dimension.
\end{proof}

\begin{prop}
  Let $k$ be an algebraically closed field, and let $R=k\langle\langle x
  ,y\rangle\rangle/\Phi$ with $\Phi\in {\frakm}_R^2$, $\Phi_2$
  nondegenerate.  If $Z(R)$ is a regular $2$-dimensional ring and $R$ is a
  free $Z(R)$-module of rank $r^2$, then for any $l\ge 0$,
  \[
    {\frakm}_R^{lr}\cap Z(R) = {\frakm}_{Z(R)}^l.
  \]
\end{prop}

\begin{proof}
  Since $Z(R)$ is regular, we may write it as $k[[u,v]]$ for suitable
  central elements $u$ and $v$, where WLOG $\deg(u)\le \deg(v)$; moreover,
  if $\deg(u)=\deg(v)$, we may assume the leading terms are linearly
  independent, as otherwise we can subtract a multiple of $u$ from $v$ to
  increase the latter's degree.  The claim will follow immediately if we
  can show that $\deg(u)=\deg(v)=r$, since then their images in $\gr R$ are
  algebraically independent.  (The center of $\gr R$ is easily seen to be
  of the form $k[u',v']$ with $\deg(u')=\deg(v')$ a divisor of $r$, and
  linearly independent homogeneous polynomials of the same degree are
  algebraically independent.)

  The freeness condition implies that $\dim_k(R/{\frakm}_{Z(R)}R)=r^2$,
  while the Hilbert series of ${\frakm}_{Z(R)} R=uR+vR$ (i.e., the Hilbert
  series of the associated graded) is upper bounded (in each coefficient)
  by
  \[
  \frac{z^{\deg(u)}}{(1-z)^2}+
  \frac{z^{\deg(v)}}{(1-z)^2}.
  \]
  Thus the Hilbert series of the quotient is lower bounded by
  \[
  \frac{1-z^{\deg(u)}-z^{\deg(v)}}{(1-z)^2}
  =
  \frac{(1-z^{\deg(u)})(1-z^{\deg(v)})}{(1-z)^2}
  +
  O(z^{\deg(u)+\deg(v)}).
  \]
  The first term has degree $\deg(u)+\deg(v)-2$, and thus the dimension is
  at least the value of the first term at $z=1$, i.e., $\deg(u)\deg(v)$, so
  that $\deg(u)\deg(v)\le r^2$.  

  Now, $u$ generates a prime ideal in $Z(R)$, and the residue field of the
  localization is isomorphic to $k((u))$, so has trivial Brauer group by
  \cite{LangS:1952}.  We may thus apply the Corollary to conclude that
  $R_{(u)}/uR_{(u)}$ has length $\ge r$, so that we may choose a chain
  \[
  R_{(u)}=I_0\supsetneq I_1\supsetneq \cdots\supsetneq I_r=u R_{(u)}
  \]
  of left ideals.  Since $R$ has global dimension $2$
  (\cite{vanGastelM/VandenBerghM:1997}), each module $R\cap I_i$ is free,
  and thus we may choose a sequence of elements $z_i$ such that $R\cap I_i
  = z_i R$, with $z_0=1$, $z_r=u$, so that we have represented our chain as
  a descending chain of principal left ideals.  Since $z_i R\supsetneq
  z_{i+1}R$, it follows that $z_{i+1}\in z_i{\frakm}$ for $1\le i\le n$,
  and thus $u=z_r\in {\frakm}^r$, so that $r\le \deg(u)\le \deg(v)$, which
  together with $\deg(u)\deg(v)\le r^2$ implies $\deg(u)=\deg(v)=r$.
\end{proof}

\begin{rem}
  It follows from the proof that $uR+vR$ has Hilbert
  series $(2z^r-z^{2r})/(1-z)^2$, and thus that $\gr u$ and $\gr v$ have no
  nontrivial syzygies.  In particular, they have trivial $\gcd$ as elements
  of the polynomial ring $Z(\gr R)$.
\end{rem}

We thus find that the van den Bergh blowup indeed corresponds to a sheaf
$\widetilde{\cal A}$ of algebras on $\widetilde{Z}$.  Moreover, that sheaf
is coherent, since if we quotient by the degree 1 elements of the center,
the resulting graded algebra has bounded degree; this reduces immediately
to the completion and thus to $A$ since it is really a question about modules.

It remains only to show that $\widetilde{\cal A}$ is maximal among orders
agreeing with ${\cal A}$ outside the exceptional locus $e$.  This reduces to
showing that it is locally free at any point of $e$ and that its completion
at the corresponding valuation is maximal.  Here we may use what we know
about sheaves on the blowup.

The local freeness condition reduces to showing that $\widetilde{\cal
  A}\otimes_{\sO_{\widetilde{Z}}}\sO_e$ is a torsion-free sheaf on $e$.  The
result is an extension of $\widetilde{\cal A}$-modules of the form
$\sO_e(d)$, so it suffices to prove the result for those modules.  Twisting
by $\sO_{\widetilde{Z}}(-ne)$ for $n\gg 0$ lets us reduce to the case $d<0$,
where it follows from the fact that $\sO_e(d)$ has no global sections.
(Here we meant $\sO_e(d)$ as a sheaf on the blowup, though it follows from
torsion-freeness and the Euler characteristic that it is also $\sO_e(d)$ as
a sheaf on $\widetilde{Z}$\dots)

It remains only to show maximality over the local ring at $e$.  Since
${\cal A}_z$ has global dimension 2 (its simple modules have injective
dimension 2, since this holds for $A$), the same holds for the blowup, and
thus the localization of the blowup has global dimension 1, so is
hereditary.  Moreover, it has a unique simple module (the now isomorphic
sheaves $\sO_e(d)$), and thus is maximal as required, finishing the proof
of Theorem \ref{thm:blowup_of_order}.
\end{proof}

\begin{cor}
  With hypotheses as above, if the orbit of $x\in Y$ has size $r$, then
  $\widetilde{\cal A}$ is an Azumaya algebra on every point of $e$ not
  meeting the blowup of $Y$.
\end{cor}

\begin{proof}
  In this case, $R=Z(R)$ and $A$ is the order in $\Mat_r(R)$ consisting of
  matrices in which the entries above the diagonal are multiples of the
  equation $U$ of $Y$ inside $\Spec(R)$.  The relevant affine patch of the
  blowup is then given by $A[{\frakm}_r\cdots {\frakm}_1 U^{-1}]\cong
  \Mat_r(R[{\frakm}_R/U])$, so that $\widetilde{\cal A}$ is a matrix
  algebra on the complement of $Y$ in the pullback of the local ring at
  $z$.
\end{proof}

\begin{rem}
  This result is particularly nice in the case of ruled surfaces, when
  every orbit has size $1$ or $r$, so that when an iterated blowup of a
  noncommutative ruled surface is a maximal order, it is an Azumaya algebra
  on the complement of the curve of points.
\end{rem}

In the cases of present interest, the invertible ideal corresponding to the
divisor $Y$ is $[{\cal A},{\cal A}]{\cal A}$ (i.e., the ideal sheaf of the
scheme classifying ${\cal A}$-module quotients of ${\cal A}$ with Hilbert
polynomial $1$ over $Z$), and we would like to understand the corresponding
ideal in $\widetilde{\cal A}$.  That is, we want to know how the curve of
points on $\widetilde{X}$ is related to the curve of points on $X$.

\begin{prop}
  Let $Z$ be a (commutative) smooth surface, let ${\cal A}$ be an order on
  $Z$ with invertible commutator ideal, and let ${\frakm}$ be a maximal
  ideal of ${\cal A}$ containing the commutator ideal.  If $\widetilde{\cal
    A}$ is the order obtained by blowing up ${\frakm}$, then the commutator
  ideal of $\widetilde{\cal A}$ is invertible, and in the notation of
  \cite{VandenBerghM:1998} is isomorphic to $\sO(-1)$.
\end{prop}

\begin{proof}
  This is a local question, and thus reduces to showing that for $l\gg 0$,
  the degree $l$ component of the commutator ideal of the graded algebra
  $
  A[({\frakm}_n\cdots {\frakm}_1)^rt]
  $
  is equal to 
  \[
  N {\frakm}_{n-1}\cdots {\frakm}_1({\frakm}_n\cdots {\frakm}_1)^{rl-1}
  t^l,
  \]
  corresponding to the ideal generated by $t$ in the original graded
  algebra before taking the Veronese and twisting.

  For $n>1$, we first observe that
  \[
  A [{\frakm}_n\cdots {\frakm}_1,{\frakm}_n\cdots {\frakm}_1]A
  \subset
  N {\frakm}_{n-1}\cdots {\frakm}_1
  ({\frakm}_n\cdots {\frakm}_1)
  \]
  or in other words the subspace of $A$ in which the coefficients on or
  above the diagonal are in ${\frakm}_RU$ and the coefficients below the
  diagonal are in ${\frakm}_R^2$.  This is a two-sided ideal, so it suffices
  to consider commutators, and since the ideal differs from
  $
  ({\frakm}_n\cdots {\frakm}_1)^2
  $
  only in the diagonal, it suffices to consider the diagonal coefficients
  of commutators.  The only commutators that can possibly contribute are
  those of diagonal matrices, where the claim follows from $[R,R]\subset UR$.

  We can now compute by induction that for $a,b\ge 1$,
  \[
    A [({\frakm}_n\cdots {\frakm}_1)^a,({\frakm}_n\cdots {\frakm}_1)^b]A
    \subset
    N {\frakm}_{n-1}\cdots {\frakm}_1
    ({\frakm}_n\cdots {\frakm}_1)^{a+b-1},
  \]
  while a similar calculation gives
  \[
  A[A,({\frakm}_n\cdots {\frakm}_1)^a]A
  =
  N {\frakm}_{n-1}\cdots {\frakm}_1
  ({\frakm}_n\cdots {\frakm}_1)^{a-1}.
  \]
  It follows that the given ideal agrees with the commutator ideal in
  positive degree, so the saturations agree.

  It remains only to consider the case $n=1$.  Let $r'=p$ if $\gr R$ is
  abelian, and otherwise let $r'$ be the rank of $\gr R$ over its center.
  Then $r$ is a multiple of $r'$, and we may reduce to showing
  \[
  R[{\frakm}^{r'},{\frakm}^{r'}]R \subset {\frakm}^{2r'-1}U
  \]
  and
  \[
  R[R,{\frakm}^{r'}]R={\frakm}^{r'-1}U.
  \]

  If $\gr R$ is abelian, then $(g,h)\mapsto [g,h]U^{-1}$ induces a
  Poisson structure on $\gr R=k[x,y]$, namely the degree $-2$ bracket
  $\{x,y\} = 1$, and the claims become that the bracket vanishes
  when both arguments have degree $p$ and that any homogeneous
  polynomial of degree $p-1$ is in the span of the brackets with
  one argument of degree $p$.  For the first claim, we have
  \[
  \{x^a y^{p-a},x^b y^{p-b}\}
  =
  (a(p-b)-(p-a)b) x^{a+b-1} y^{2p-a-b-1} = 0,
  \]
  and the second claim follows by noting that the brackets with $x$ span
  the polynomials not of degree $p-1$ in $y$, while the brackets with $y$
  span the polynomials not of degree $p-1$ in $x$.

  When $\gr R$ is not abelian, then the calculation reduces to one
  in $\gr R$, so that we reduce to considering $k\langle
  x,y\rangle/(yx-qxy)$ and $k\langle x,y\rangle/(yx-xy-x^2)$.  The
  first ($q$-Weyl) case is straightforward:
  \[
    [x^a y^b,x^c y^d]
    =
    (q^{bc}-q^{ad}) x^{a+c} y^{b+d},
  \]
  which again is 0 if $q^{a+b}=q^{c+d}=1$, and the span
  $
  \langle [x,x^a y^b],[y,x^a y^b]\rangle
  $
  for $a+b=\ord(q)$ has the correct dimension.

  For $k\langle x,y\rangle/(yx-xy-x^2)$, we observe that this has a
  representation in differential operators: $x=t$, $y=t^2 D_t$.
  We in particular easily find
  $
    [y,x^a y^{p-a}] = a x^{a+1} y^{p-a}
  $
  which together with $[x,x^{p-1} y]=-x^{p+1}$ span a space of the correct
  dimension.  It thus remains only to show that
  $
    [x^a y^{p-a},x^b y^{p-b}]=0.
  $
  Doing the change of variable $t\mapsto 1/t$ in the corresponding
  commutator of differential operators and observing that $t^p$ is
  central reduces this to showing that
  $
    [t^a D_t^a,t^b D_t^b] = 0,
  $
  which is true in arbitrary characteristic since $t^a D_t^a$ is a
  polynomial in $t D_t$ for all $a\ge 0$.
\end{proof}

\medskip

In the proof that the blowup of a maximal order is a maximal order, we
skirted around the question of determining the precise structure of the
localization along the exceptional divisor.  Although this is somewhat
tricky to determine in general, it is fairly straightforward in the cases
corresponding to ruled surfaces.  

For instance, in the difference cases of elliptic or multiplicative type,
the division ring is a localization of the ring of (symmetric) difference
operators.  Over the local ring of a codimension 1 point, the residue field
of the maximal order is a separable extension of the residue field of the
center.  Choose an element $u$ whose image generates that extension (of
degree $n$ if the division ring has dimension $n^2$), and observe that if
take the corresponding unramified base change, then this splits the
division ring, and lets us see that for each element $\sigma$ of the Galois
group, there is an $n$-dimensional subspace of the division ring consisting
of elements $v$ satisfying $vu=\sigma(u)v$, which is then a 1-dimensional
space over $k(u)$.  In particular, there exists an element of this form
having valuation 1, so that $v^n$ is a uniformizer of the center.

For the remaining cases, a key observation is that given a discrete
valuation on the center of a division ring, there is a unique extension to
a maximal order over the valuation ring.  In particular, blowing up does
not change the division ring, merely the set of available valuations.  We
claim that the division ring is generated over its center by elements $u$,
$v$ such that $[u,v]=1$.  This is clear in the differential cases, while
the remaining cases (additive difference and the hybrid case) are
birational to the differential case, so have the same division ring.  We
then find the following.

\begin{prop}\label{prop:local_structure_of_order_single_component}
  Let $X$ be a rationally ruled surface over a field of characteristic $p$
  such that $Q$ has a non-nodal singular point.  Then $X=\Spec({\cal A})$
  for a maximal order ${\cal A}$ with generic fiber a division ring, and
  for any associated point $x\in Q$, the localization ${\cal A}_x$ is
  generated by elements $u$, $v$ with $u$ a unit, $v$ of valuation $1$, and
  $[u,v]=v^m$, where $m$ is the multiplicity of $Q$ along $x$.
\end{prop}

\begin{proof}
  We need to understand the extension of the generic fiber of ${\cal A}$ to
  the corresponding valuation ring, and have already observed that the
  generic fiber is generated by elements $u,v$ with $[u,v]=1$, $u^p=U$,
  $v^p=V$ for appropriate elements $U$, $V$ of the center.  If these
  elements are both integral, then $1\in [{\cal A}_x,{\cal A}_x]$, so that
  $x$ was not a component of $Q$.  If $U$ has valuation a multiple of $p$
  and its leading coefficient is a $p$-th power, then we can subtract a
  central element from $u$ to make $U$ smaller, and similarly for $V$.  If
  $U$ and $V$ still both have valuation a multiple of $p$, then we can
  rescale $(u,v)\mapsto (\pi^l u,\pi^{-l} v)$, where $\pi$ is a uniformizer
  of the center, to make $U$ a unit which is not a $p$-th power in the
  residue field.  In particular $U+{\frakm}$ generates the residue field
  over its $p$-th power subfield, so that we may add a suitable polynomial
  in $u$ to $v$ to reduce its valuation.

  In this way, we reduce to the case that at least one of $u$ or $v$ (say
  $v$) has valuation prime to $p$.  For any integers $a$, $l$ with $a$
  prime to $p$, there is an automorphism
  \[
  (u,v)\mapsto (a^{-1}\pi^{-l}v^{1-a}u,\pi^l v^a)
  \]
  of the division ring, and thus by suitable choice of $l$ and $a$, we may
  arrange to have $v$ (and thus $V$) of valuation $1$.  If $U$ has
  valuation $-l$, then $v^l u$ is a unit, and we have
  \[
  [v^l u,v] = v^l.
  \]

  These elements still generate the division ring, and the order they
  generate is easily verified to be maximal, so that they generate ${\cal
    A}_x$ as required.  Moreover, it follows that the 2-sided ideal
  generated by commutators is the principal ideal generated by $v^l$, so
  that the quotient has length $l=m$.
\end{proof}

\begin{rems}
   Note that when the multiplicity is $1$, this looks like $uv=v(u+1)$,
   which is an instance of the separable case, and in particular shows
   directly that the additive difference cases are birational to the
   differential cases.
\end{rems}

\begin{rems}
The situation for quasi-ruled surfaces is more complicated, as the residue
field extension may be an inseparable extension of degree larger than $p$.
\end{rems}

It would be nice to have a similar description for the localizations of
${\cal A}$ at closed points of $Q$, even if only for the corresponding
completion.  For the completion, there is no difficulty at smooth points,
as in the ruled case the relevant automorphism acts faithfully, and thus we
obtain the situation above with $A\subset \Mat_n(R)$ and $R$ commutative.  At
nodes, the situation is only slightly more complicated.  In that case, the
leading term of the relation of $R$ must be $vu-quv$ with $q$ a primitive
$n$-th root of unity, and in any such ring, one can perform changes of
variables to make the relation have the form $vu = quv + u\Phi v$ where
$\Phi$ is in the span of $u^{ni}v^{nj}$ with $i+j>0$.  (We certainly have a
relation of the form $vu=quv+\Phi'$ with $\Phi'$ in the span of $u^i v^j$
with $i+j>0$, and if $\Phi'$ is only of the form $u\Phi v$ to degree $d$,
then we can improve the agreement by adding suitable terms of degree $d$ to
$u$ and $v$.)  Then $(u,v)\mapsto (qu,v)$ and $(u,v)\mapsto (u,qv)$ are
automorphisms, so act on the center.  The center is generated by two
elements of degree $n$, the leading terms of which must be central in the
associated graded of $R$, so are in the span of $u^n$ and $v^n$.  It then
follows that the center is generated by invariant elements, and must
therefore (by Hilbert series considerations) be isomorphic to
$k[[u^n,v^n]]$.  In particular, $\Phi$ is central, and thus
\[
v^n u^n = (q+\Phi)^{n^2} u^n v^n,
\]
so that $\Phi=0$.  In other words, the completion at a node has the form
$k\langle\langle u,v\rangle\rangle/(vu-quv)$ where $q$ is a primitive
$n$-th root of unity.

Thus only the differential (or additive) case remains open.  It is unclear
whether simply having degree $p^2$ over the center is enough to pin down
the completion (possibly together with the known possibilities for the
structure of $R/R[R,R]R$) as it was in the above case of degree prime to
$p$.

\chapter{Noncommutative surfaces as $t$-structures}
\label{chap:derived}

\section{Quasi-ruled surfaces}

Although the above construction of noncommutative ruled surfaces was well
suited to an interpretation via difference/differential operators as well
as to understanding the center, it is not ideal for a number of other
applications, especially since it does not behave well in families.  Van
den Bergh's original construction does not have this issue, but in some
respects turns out to be too explicit; constructing the various
isomorphisms we require (e.g., between the two interpretations of a
noncommutative $\P^1\times \P^1$ as a ruled surface) via that construction
either requires working with sheaves on highly singular curves or requires
that one split into a large number of separate cases that can be dealt with
explicitly.  Luckily, it turns out that there is another way to construct
noncommutative $\P^1$-bundles: there is a relatively simple description of
the corresponding derived category, as well as the associated
$t$-structure.  Although this construction cannot quite replace van den
Bergh's construction (our proof that the category is the derived category
of the heart of the $t$-structure uses van den Bergh's construction), it
gives an alternate approach to constructing isomorphisms, namely as derived
equivalences respecting the $t$-structures.  In each of the cases of
interest, it is easy to construct the derived equivalence corresponding to
the desired isomorphism, and not too difficult to show sufficient
exactness.

The overall approach works for each of the main constructions in the
literature (noncommutative (projective) planes
\cite{ArtinM/TateJ/VandenBerghM:1990,BondalAI/PolishchukAE:1993},
noncommutative $\P^1$-bundles \cite{VandenBerghM:2012}, and blowups
\cite{VandenBerghM:1998}), but we begin by considering the $\P^1$-bundle
case.  The basic idea is that there is a natural semiorthogonal
decomposition of the derived category of a commutative $\P^1$-bundle over a
smooth projective scheme.  Not only does this extend to the noncommutative
setting, but we shall see that one can recover the $t$-structure from the
decomposition.  Although the discussion below works for general
noncommutative $\P^1$-bundles over smooth projective schemes, we consider
only the surface case, and at least initially work over an algebraically
closed field.  (Otherwise, we would need to consider things like conic
bundles!)


Recall that the category $\qcoh \bar{\cal S}$ is defined as the quotient of
the category of $\bar{\cal S}$-modules by the subcategory of bounded such
modules (i.e., which become 0 in sufficiently large degree), with $\coh
\bar{\cal S}$ the subcategory of Noetherian modules.  For each $d$, there
is a natural functor $\rho_d^{-1}$ from $\coh (C_d:=C_{d\bmod 2})$ to the
category of $\bar{\cal S}$-modules taking a coherent sheaf $M$ on $C_d$ to
the representation which in degree $d'$ is $\bar{\cal S}_{dd'}\otimes_{C_d}
M$.  Composition with the quotient morphism gives a functor $\rho_d^*:\coh
C_d\to \coh \bar{\cal S}$.  Since the quotient morphism has a right
adjoint, as does $\rho_d^{-1}$, it follows that $\rho_d^*$ has a right
adjoint $\rho_{d*}$.

\begin{lem}\label{lem:exact_tri_quasi_ruled}
  For any vector bundle $V$ on $C_0$, the quadratic relation gives rise to
  an exact sequence
  \[
  0\to \rho_2^*(V\otimes \det(\pi_{0*}\sO_{\hat{Q}})^{-1})
  \to \rho_1^*(\pi_{1*}(\pi_0^*V\otimes \bar{\cal S}_{01}))
  \to \rho_0^*V
  \to 0.
  \]
\end{lem}

\begin{proof}
  It suffices to check exactness locally, where it reduces to showing
  exactness of
  \[
  0\to \bar{S}_{2n}\to \bar{S}_{1n}\otimes \bar{S}_{01}\to \bar{S}_{0n}\to 0.
  \]
  The second (multiplication) map is surjective since $\bar{S}$ is
  generated in degree 1, while the first map is
  \[
  x\mapsto (x N_0 e_1\otimes N_1 \xi_1 e_0)-(x N_0 s_1(\xi_1) e_1\otimes
  N_1 e_0),
  \]
  which maps into the kernel of the second map.  Since
  \[
  x = (x N_0 e_1)(N_1 \xi_1 \xi_0 e_0)-(x N_0 s_1(\xi_1) e_1)(N_1 \xi_0 e_0),
  \]
  the map from $\bar{S}_{2n}$ to the kernel of the second map splits, and
  is thus an isomorphism by rank considerations.
\end{proof}  

\begin{lem}
   For $M\in D_{\qcoh} \bar{\cal S}$, if $R\rho_{0*}M=R\rho_{1*}M=0$, then $M=0$.
\end{lem}

\begin{proof}
  Suppose $M$ is an object such that $R\rho_{0*}M=R\rho_{1*}M=0$.  From the
  previous Lemma, we deduce that for any vector bundle $V$ on $C_0$,
  $R\Hom(L\rho_2^*V,M)=0$, and thus $R\rho_{2*}M=0$.  It follows by
  induction that $R\rho_{d*}M=0$ for all $d\ge 0$, and thus that the
  cohomology modules of the corresponding complex in $\bar{\cal
    S}\text{-mod}$ are torsion, so that $M$ is 0 in $D_{\qcoh}\bar{\cal
    S}$.
\end{proof}

The following was essentially shown in \cite{MoriI:2007}, subject to
the (luckily unused) assumption that the two base schemes (i.e., $C_0$ and
$C_1$) are equal.  (To be precise, the reference showed the result when $M$
and $N$ are line bundles, but this implies it in general.)

\begin{lem}\cite[Lem.~4.4]{MoriI:2007}
  For $M,N\in D^b_{\coh} C_d$, one has
  $R\Hom(L\rho_d^*M,L\rho_d^*N)=R\Hom(M,N)$, while for $N'\in D^b_{\coh}
  C_{d+1}$ one has $R\Hom(L\rho_d^*M,L\rho_{d+1}^*N')=0$.
\end{lem}

\begin{cor}
  For $M\in D^b_{\coh} C_{d+1}$, one has
  $R\rho_{d*}L\rho_{d+1}^*M=0$.
\end{cor}

\begin{proof}
  Indeed, for any $N\in D^b_{\coh} C_d$, one has
  \[
  R\Hom(N,R\rho_{d*}L\rho_{d+1}^*M)\cong
  R\Hom(L\rho_d^*N,L\rho_{d+1}^*M)=0,
  \]
  and thus in particular
  \[
  R\Hom(R\rho_{d*}L\rho_{d+1}^*M,R\rho_{d*}L\rho_{d+1}^*M)=0,
  \]
  from which the conclusion follows immediately (the identity endomorphism
  is null homotopic!).
\end{proof}

Combining the above, we obtain the following, compare Proposition
\ref{prop:weak_semiorth_on_ruled}.

\begin{thm}\label{thm:semiorth_for_quasiruled}
  The subcategories $L\rho_0^*D^b_{\coh} C_0$ and $L\rho_1^*D^b_{\coh} C_1$
  form a semiorthogonal decomposition of $D^b_{\coh} \bar{\cal S}$: for any
  object $M\in D^b_{\coh} \bar{\cal S}$, there is a unique distinguished
  triangle of the form
  \[
  L\rho_0^* N_0\to M\to L\rho_1^*N_1\to.
  \]
\end{thm}

\begin{proof}
  Adjunction gives a natural morphism $L\rho_0^* R\rho_{0*}M\to M$, which
  extends to a distinguished triangle
  \[
  L\rho_0^*R\rho_{0*}M\to M\to M'\to.
  \]
  Applying $R\rho_{0*}$ and using $R\rho_{0*}L\rho_0^*\cong \text{id}$ implies
  that $R\rho_{0*}M'=0$.  The natural distinguished triangle
  \[
  L\rho_1^*\rho_{1*}M'\to M'\to M''\to
  \]
  then gives an object $M''$ such that $R\rho_{1*}M''=0$ and
  $R\rho_{0*}M''\cong R\rho_{0*}M'=0$, and thus $M''=0$, so that $M'$ is in
  the image of $L\rho_1^*$.

  For uniqueness, observe that for any such triangle, applying $R\rho_{0*}$
  gives $N_0\cong R\rho_{0*}M$.
\end{proof}

\begin{rem}
  As we mentioned, the same proof shows that this holds more generally for
  the noncommutative $\P^1$-bundle associated to any sheaf bimodule on
  Noetherian schemes $Y_0$ and $Y_1$, except that if the base schemes are
  singular, one should instead work with the categories of compact objects
  (which on $Y_i$ are just the perfect complexes).
\end{rem}

The semiorthogonal decomposition yields two immediate corollaries.  The
first is that we can easily compute the Grothendieck group of $\coh
\bar{\cal S}$, since the Grothendieck of a derived category with a
semiorthogonal decomposition is just the sum of the Grothendieck groups of
the two subcategories.

\begin{cor}
  We have $K_0\coh \bar{\cal S}\cong K_0\coh C_0\oplus K_0\coh C_1$.
\end{cor}

Furthermore, the Mukai pairing on $K_0\coh \bar{\cal S}$ is easily
reconstructed from the pairings on $C_0$ and $C_1$ and the induced pairing
\[
(M,N)\mapsto \chi R\Hom(\rho_1^*M,\rho_0^*N)
\]
between the two curves.  The relevant Hom space was again computed in
\cite{MoriI:2007}.

\begin{lem}
  We have $R\Hom(\rho_1^*M,\rho_0^*N)\cong R\Hom(M,N\otimes_{C_0} \bar{\cal
    S}_{01})$.
\end{lem}

\begin{rems}
  This of course not only allows us to compute the Mukai pairing, but tells
  us that it is well-defined, i.e., that the $\Hom$ complexes in
  $D^b_{\coh}\bar{\cal S}$ are finite-dimensional, since they are built out
  of $\Hom$ complexes in $D^b_{\coh}C_0$ and $D^b_{\coh}C_1$ and a complex
  of the above form.
\end{rems}

\begin{rems}
  The functor $N\mapsto \rho_{1*}\rho_0^*N\cong N\otimes_{C_0} \bar{\cal
    S}_{01}$ from $D^b_{\coh} C_0$ to $D^b_{\coh} C_1$ is nothing other
  than the Fourier-Mukai functor with kernel $\bar{\cal S}_{01}$.
\end{rems}

The other corollary is slightly more subtle.  Since we have not just the
one semiorthogonal decomposition, but one for each $d$, and $L\rho_d^*D^b
\coh C_d$ appears in two such decompositions, once on each side, we find
that each such subcategory is admissible.  Since both subcategories have
Serre functors (being derived categories of projective schemes), the same
is true for $D^b_{\coh} \bar{\cal S}$ (per
\cite[Prop.~3.8]{BondalAI/KapranovMM:1990}).  Moreover, we can compute the
Serre functor explicitly.

\begin{prop}\label{prop:Serre_for_quasiruled}
  The triangulated category $D^b_{\coh} \bar{\cal S}$ has a Serre functor $S$
  which is the composition of $M\mapsto M(-\hat{Q})[2]$ with twisting by
  the pair of line bundles $\omega_{C_i}\otimes
  \det(\pi_{i*}\sO_{\hat{Q}})^{-1}$ having the same pullback to $\hat{Q}$.
\end{prop}

\begin{proof}
  It suffices to compute how the Serre functor acts on the two pieces of
  the semiorthogonal decomposition, and thus (since shifting degrees gives
  another quasi-ruled surface) how it acts on an object $L\rho_0^*M$.  We
  thus need an object $SL\rho_0^*M$ such that
  \[
  \Hom(L\rho_0^*N,SL\rho_0^*M)\cong \Hom(L\rho_0^*M,L\rho_0^*N)^*
  \cong \Hom(M,N)^*\cong \Hom(N,SM),
  \]
  where we also denote the Serre functor on $D^b_{\coh} C_0$ by $S$, and
  \[
  \Hom(L\rho_1^*N,SL\rho_0^*M)\cong \Hom(L\rho_0^*M,L\rho_1^*N)^*=0.
  \]
  In other words, we want $R\rho_{0*}SL\rho_0^*M\cong SM$ and
  $R\rho_{1*}SL\rho_0^*M=0$.

  Now, the short exact sequence of Lemma \ref{lem:exact_tri_quasi_ruled}
  extends to a distinguished triangle for perfect complexes in general, of
  the form
  \[
  L\rho_2^*(N\otimes \det(\pi_{0*}\sO_{\hat{Q}})^{-1})
  \to
  L\rho_1^*(N\otimes_{C_0} \bar{\cal S}_{01})
  \to
  L\rho_0^*(N)
  \to.
  \]
  Applying $R\rho_{0*}$ makes the connecting map an isomorphism
  \[
  R\rho_{0*}L\rho_2^*(N\otimes \det(\pi_{0*}\sO_{\hat{Q}})^{-1})
  \cong
  N[-1]
  \]
  while
  \[
  R\rho_{1*}L\rho_2^*(N\otimes \det(\pi_{0*}\sO_{\hat{Q}})^{-1})=0.
  \]
  We thus conclude that
  \[
  SL\rho_0^*M\cong L\rho_2^*(SM\otimes \det(\pi_{0*}\sO_{\hat{Q}})^{-1})[1]
  \cong L\rho_2^*(M\otimes \omega_{C_0}\otimes
  \det(\pi_{0*}\sO_{\hat{Q}})^{-1})[2],
  \]
  agreeing with the desired formula for the Serre functor.
\end{proof}

\begin{rem}
  The analogous result was shown (again with the assumption that the base
  schemes are isomorphic) for general smooth projective base in
  \cite{NymanA:2005}.  The semiorthogonal decomposition yields a
  significant streamlining of the argument, however.
\end{rem}

Since $M(-Q)$ itself differs from $M(-\hat{Q})$ by a pair of twists by
line bundles, we can also write the Serre functor as the composition of
${-}(-Q)[2]$ with twisting by $\omega_{C_i}\otimes
\det(\pi_{i*}\sO_Q)^{-1}$.  As we remarked after Proposition
\ref{prop:genus_1_almost_implies_ruled}, these bundles are trivial iff the
surface is ruled, and thus we obtain the following.

\begin{cor}
  The surface $\bar{\cal S}$ is ruled iff the Serre functor is equivalent
  to ${-}(-Q)[2]$.
\end{cor}

In other words, a quasi-ruled surface is ruled iff the curve of points is
anticanonical.
\medskip

There is one important notational issue to consider here, namely that since
$\hat{Q}$ does not behave well in flat families, neither do the functors
$L\rho_d^*$ in general (though for $d\notin \{0,1\}$, we can rephrase
everything in terms of the sheaf bimodule).  We could fix this by using the
labelling of \cite{VandenBerghM:2012}, but it turns out that there is an
even nicer choice.  If we leave $\rho_0^*$ and $\rho_1^*$ fixed, then we
may twist the remaining functors by line bundles and automorphisms so that
\[
L\rho_{d+2}^* = S L\rho_d^* S^{-1}[-1]
\]
With this relabelling, the distinguished triangle coming from the
semiorthogonal decomposition becomes
\[
L\rho_{d+1}^*R\rho_{(d-1)*}M\to L\rho_d^*R\rho_{d*}M\to M\to.
\]
Indeed, applying $R\rho_{(d-1)*}$ to the semiorthogonal decomposition
reduces this to showing
\[
R\rho_{(d-1)*}L\rho_{d+1}^*[1]\cong \text{id},
\]
which follows (for coherent sheaves, and thus in general) from the
computation
\begin{align}
R\Hom(M,R\rho_{(d-1)*}L\rho_{d+1}^*N[1])
&\cong
R\Hom(L\rho_{d-1}^*M,S L\rho_{d-1}^*S^{-1}N)\notag\\
&\cong
R\Hom(L\rho_{d-1}^*S^{-1}N,L\rho_{d-1}^*M)^*\notag\\
&\cong
R\Hom(S^{-1}N,M)^*\notag\\
&\cong
R\Hom(N,SM)^*\notag\\
&\cong
R\Hom(M,N).
\end{align}
We will use this alternate labelling below.

\medskip

The results of \cite{OrlovD:2016} tell us that we can reverse the above
construction: given a pair of (dg-enhanced, automatic for commutative
projective schemes) triangulated categories and a suitable functor from one
to the other, there is a corresponding ``glued'' category with a
semiorthogonal decomposition such that the $\Hom$s between the two
subcategories are induced by the functor.  This gives an immediate
construction of a triangulated category from a sheaf bimodule: simply glue
the two categories using the associated Fourier-Mukai functor.  In the case
of sheaf bimodules such that both direct images are locally free of rank 2,
the above discussion shows that the resulting category is precisely the
derived category of the corresponding noncommutative $\P^1$-bundle.

As mentioned above, we can actually extend this to give an alternate
construction of noncommutative $\P^1$-bundles.  To determine the
noncommutative $\P^1$-bundle from its derived category, we need to specify
the $t$-structure, i.e., those objects with only nonnegative or only
nonpositive cohomology.  We again consider the surface case, leaving the
modifications for more complicated base schemes to the reader.  In that
case, the Serre functor has the form $\theta[2]$ where $\theta$ is an
autoequivalence of the desired abelian category.  Moreover, since $\theta$
is just a twist of a degree shift of the sheaf algebra, we find that
$\theta^{-1}$ is {\em relatively} ample, in that any nonzero sheaf $M\in
\coh \bar{\cal S}$ has the property that $R\rho_{0*}\theta^{-d}M$ is a
nonzero sheaf for all sufficiently large $d$.  It then follows more
generally that $M\in D^b_{\coh} \bar{\cal S}^{\ge 0}$ iff
$R\rho_{0*}\theta^d M\in D^b_{\coh} C_0^{\ge 0}$ for all $d$.  Indeed,
$\rho_{0*}$ is right exact and $\theta$ is exact, so the condition
certainly holds for all $M\in D^b_{\coh} \bar{\cal S}^{\ge 0}$, while if
$M$ has cohomology in some negative degree, then some power of
$\theta^{-1}$ makes that cohomology have nonzero direct image.  We can then
reconstruct the other half $D^b_{\coh} \bar{\cal S}^{\le 0}$ of the
$t$-structure as the objects with no maps to objects in $D^b_{\coh} \bar{\cal
  S}^{\ge 1}$.

The significance of this is that the Serre functor is intrinsic to the
triangulated category, and thus so is $\theta$.  In other words, the
``geometric'' $t$-structure can be reconstructed given only the
triangulated category and the functor $\rho_{0*}$.  This construction works
fairly generally: given any functor $f:{\cal A}\to {\cal B}$ of
triangulated categories such that ${\cal A}$ has a Serre functor and ${\cal
  B}$ is equipped with a $t$-structure, and any integer $d$, there is a
putative $t$-structure defined by taking ${\cal A}^{\ge 0}$ to be those
objects such that $f(S^k M[-dk])\in {\cal B}^{\ge 0}$ for all $k\in \Z$.
If we already have a $t$-structure in mind for ${\cal A}$, then we say that
$f$ is ``pseudo-canonical (of dimension $d$)'' if the putative
$t$-structure agrees with the specified $t$-structure.  (Note that this
forces $S[-d]$ to be exact, so is probably not quite the right notion in
general, but is good enough for present purposes since we only care about
analogues of smooth surfaces.  In general, one should probably conjugate
$f$ by the respective shifted Serre functors.) In the case of a smooth
commutative projective $d$-fold $X$, this condition holds if either $-K_X$
or $K_X$ is relatively ample.

Recall that a quasi-scheme is {\em Gorenstein} of dimension $d$ if it has a
Serre functor such that $S[-d]$ is an abelian autoequivalence.

\begin{lem}
  Suppose that $X$ and $Y$ are Gorenstein quasi-schemes ($X$ of dimension
  $d$) with $f:\qcoh X\to \qcoh Y$ a left exact functor such that for all
  $M\in \coh X$, there exists $l\in \Z$ such that $f(S^l M[-dl])\ne 0$.
  Then $Rf$ is pseudo-canonical of dimension $d$.
\end{lem}

\begin{proof}
  Since $f$ is left exact, $Rf$ certainly takes $D^b_{\qcoh} X^{\ge 0}$ to
  $D^b_{\qcoh} Y^{\ge 0}$, and the same applies if we compose it with any
  power of $S[-d]$.  So it will suffice to show that if $\tau_{<0}Rf S^l
  M[-dl]=0$ for all $l$ then $\tau_{<0}M=0$.  If not, then let $p$ be the
  smallest integer such that $h^pM\ne 0$.  This has a coherent subsheaf,
  and thus there is an integer $l$ such that $f(S^l(h^p M)[-ld])\ne 0$.
  But the standard spectral sequence then tells us that $h^p(Rf S^l
  M[-ld])\ne 0$, giving a contradiction.
\end{proof}

\begin{rem}
  The converse is also true: if for some nonzero sheaf $M$,
  $h^0(Rf S^l M[-ld])=0$ for all $l\in \Z$, then $M[1]$ is not a
  nonnegative complex, but its image under $RfS^l[-ld]$ is for all $l$.
\end{rem}

\section{Planes and blowups}

It turns out that the other two constructions we need also come with
semiorthogonal decompositions such that one of the associated functors is
pseudo-canonical.  For noncommutative planes, this is already in the
literature.

\begin{prop}\cite{BondalAI/PolishchukAE:1993}
  Any noncommutative plane has a strong exceptional collection consisting
  of three objects $\sO_X(-2)$, $\sO_X(-1)$, $\sO_X$ such that
  $\Hom(\sO_X(-2),\sO_X(-1))$ and $\Hom(\sO_X(-1),\sO_X)$ are
  three-dimensional and
  \[
  \Hom(\sO_X(-2),\sO_X(-1))\otimes \Hom(\sO_X(-1),\sO_X)
  \to
  \Hom(\sO_X(-2),\sO_X)
  \]
  is surjective with three-dimensional kernel.  Moreover, $D_{\qcoh} X$ is
  compactly generated, has Serre functor $M\mapsto M(-3)[2]$, and the
  functor $R\Hom(\sO_X,{-})$ is pseudo-canonical of dimension 2.
\end{prop}

\begin{rems}
  The condition that the resulting category has a well-behaved
  $t$-structure is open in the set of surjections $k^3\otimes k^3\to
  k^6$, and reduces to the condition that there be a cubic curve $C\subset
  \P^2$ and a degree $0$ line bundle $q$ on $C$ such that
  $\Hom(\sO_X(-1),\sO_X)\cong \Gamma(C;q(1))$,
  $\Hom(\sO_X(-2),\sO_X(-1))\cong \Gamma(C;\sO_C(1))$, and the composition
  is given by multiplication.  The scheme of points is isomorphic to $C$,
  which embeds as a divisor via a natural transformation ${-}(-3)\to \text{id}$.
\end{rems}

\begin{rems}
  Again, it follows immediately from this description that $\Hom$ complexes
  in $D^b_{\coh}(X)$ are finite-dimensional.
\end{rems}

It thus remains to consider blowups.  If $C$ is a (commutative) divisor on
the noncommutative surface $X$, a ``weak line bundle relative to $C$'' is a
sheaf on $X$ such that $X|_C$ is a line bundle\footnote{These were called
``line bundles'' in \cite{VandenBerghM:1998}, but we reserve this
terminology for a particular class of weak line bundles to be considered in
Definition \ref{defn:line_bundle} below.}, and we say that $X$ is
``generated by weak line bundles relative to $C$'' if every element of
$\qcoh X$ is a quotient of a direct sum of weak line bundles.  (Presumably
one could weaken this to allow sheaves with locally free restriction.)
Note that this is easily seen to be hold for a quasi-ruled surface with $C$
the curve of points, as any sheaf of the form $\rho_d^*{\cal L}$ is a weak
line bundle.  A similar statement also holds for noncommutative planes,
since the sheaves $\sO_X(d)$ are weak line bundles relative to any curve in
$X$.

\begin{thm}\cite[Thm.~8.4.1]{VandenBerghM:1998}
  Let $X$ be a noncommutative surface with $C$ a commutative divisor and
  $p\in C$ a point, and suppose that $X$ is generated by weak line bundles
  relative to $C$. Let $\widetilde{X}$ be the blowup at $p$, with associated
  functors $\alpha^*$, $\alpha_*$, and let $\sO_e(-1)$ denote the
  corresponding exceptional sheaf.  Then the subcategories $D^b_{\qcoh}
  k\otimes \sO_e(-1)$ and $L\alpha^* D^b_{\qcoh} X$ form a semiorthogonal
  decomposition of $D^b_{\qcoh} \widetilde{X}$.
\end{thm}

\begin{rem}
  To be precise, the associated distinguished triangle is
  \[
  L\alpha^*R\alpha_*M\to M\to R\Hom_k(R\Hom(M,\sO_e(-1)),\sO_e(-1))\to
  \]
  as long as $M$ is coherent.  It follows from the explicit form of Serre
  duality (Corollary \ref{cor:Serre_for_blowup} below) that one can instead
  take
  \[
  L\alpha^*R\alpha_*M\to M\to R\Hom(\sO_e,M[2])\otimes_k \sO_e(-1)\to,
  \]
  and this works even when $M$ is only quasicoherent, see the proof of
  Proposition \ref{prop:blowup_acyclic} below.
\end{rem}

\begin{cor}
  We have $K_0(\widetilde{X})\cong K_0(X)\oplus \Z$.
\end{cor}

To control the Serre functor, we will need to know that the two
subcategories are admissible.  Define a new functor $\alpha^!$ by
$\alpha^!M = \alpha^*(M(C))(-1)$.

\begin{lem}\label{lem:alpha_shriek}
  The derived functor $L\alpha^!$ is right adjoint to $R\alpha_*$.
\end{lem}

\begin{proof}
  It suffices to show that for any object $M\in D^b_{\qcoh} X$,
  \[
  R\alpha_*L\alpha^!M\cong M\qquad\text{and}\qquad
  R\Hom(\sO_e(-1),L\alpha^!M)=0.\notag
  \]
  Indeed, the second claim lets us reduce to checking the adjunction on the
  image of $L\alpha^*$, which in turn reduces by the adjunction between
  $\alpha^*$ and $\alpha_*$ to checking the first claim.
  
  By \cite[Prop.~8.3.1]{VandenBerghM:1998}, we have
  $R\alpha_*((\alpha^*E)(-1))\cong E(-C)$ for any sheaf $E$ which is a
  direct sum of weak line bundles, and since these generate, the same holds
  for any object in $D^b_{\qcoh} X$.  Applying this to $M(C)$ gives the first
  result.  It also follows that
  \[
  R\alpha_*((\alpha^*M)(-1))(C)\cong M\cong R\alpha_*L\alpha^*M
  \]
  and thus $R\Hom(\sO_e,L\alpha^*M)=0$ by \cite[Lem.~8.3.3(2)]{VandenBerghM:1998}.
  We then find
  \[
  R\Hom(\sO_e(-1),L\alpha^!M)\cong R\Hom(\sO_e,L\alpha^*M(C))=0
  \]
  as required.
\end{proof}

\begin{rem}
  It follows that there is another semiorthogonal decomposition, with
  distinguished triangle
  \[
  \sO_e(-1)\otimes_k R\Hom(\sO_e(-1),M)\to M\to L\alpha^!R\alpha_*M\to.
  \]
\end{rem}

This lets us compute the gluing functor, which again can be expressed in
terms of an adjoint pair of functors.

\begin{prop}
  We have
  \[
  R\Hom_{\widetilde{X}}(V\otimes \sO_e(-1),L\alpha^*M)
  \cong
  R\Hom_X(V\otimes \sO_{\tau p}[-1],M)
  \cong
  R\Hom_k(V,R\Hom_X(\sO_{\tau p}[-1],M)).
  \]
\end{prop}

\begin{proof}
  We first observe that
  \[
  R\Hom_{\widetilde{X}}(\sO_e(-1),L\alpha^*M)
  \cong
  R\Hom_{\widetilde{X}}(\sO_e(-2),L\alpha^!(M(-C)))
  \cong
  R\Hom_X(R\alpha_*\sO_e(-2),M(-C)),
  \]
  so that we need to compute $R\alpha_*\sO_e(-2)$.  By
  \cite[Prop.~8.3.2]{VandenBerghM:1998}, we have a distinguished triangle
  \[
  L\alpha^* \sO_{\tau p}\to \sO_e\to \sO_e(-1)[2]\to
  \label{eq:dtri_e_em1}
  \]
  which twists to give
  \[
  L\alpha^!(\sO_{\tau p}(-C))\to \sO_e(-1)\to \sO_e(-2)[2]\to.
  \]
  Applying $R\alpha_*$ then gives
  \[
  R\alpha_*\sO_e(-2)\cong \sO_{\tau p}(-C)[-1],
  \]
  so that
  \[
  R\Hom_{\widetilde{X}}(\sO_e(-1),L\alpha^*M)\cong R\Hom_X(\sO_{\tau p}[-1],M),
  \]
  from which the claim follows.
\end{proof}

\begin{rem}
  It follows that if $D^b_{\coh}(X)$ is proper, then so is
  $D^b_{\coh}(\tilde{X})$.
\end{rem}

In particular, $D_{\qcoh}(\widetilde{X})$ is compactly generated (so has a
Serre functor) iff $D_{\qcoh}(X)$ is compactly generated.  Indeed, if $G$
is a compact generator of $D_{\qcoh}(X)$, then $L\alpha^*G\oplus \sO_e(-1)$
certainly generates $D_{\qcoh}(X)$, and the Proposition tells us that it is
compact.  Moreover (\cite{BondalAI/KapranovMM:1990}), given such a
semiorthogonal decomposition, one can compute the Serre functor on the
ambient category from that on the subcategories (and vice versa, if
desired).  Since the Serre functor on $\perf(k)$ is just the identity, this
is not too hard in our case.

\begin{cor}\label{cor:Serre_for_blowup}
  If $X$ is compactly generated (i.e., $D_{\qcoh}(X)$ is compactly
  generated), then the Serre functor on $\widetilde{X}$ satisfies $S
  L\alpha^*M\cong L\alpha^! S M$ and $S \sO_e\cong \sO_e(-1)[2]$.
\end{cor}

\begin{proof}
  The first claim follows immediately from Lemma \ref{lem:alpha_shriek},
  while the second claim follows from \eqref{eq:dtri_e_em1}, since that
  distinguished triangle implies that ${-}\otimes \sO_e(-1)[2]$ is
  similarly right adjoint to the projection to ${-}\otimes \sO_e$ in the
  semiorthogonal decomposition coming from $R\Hom(\sO_e,L\alpha^* M)=0$.
\end{proof}

Note that since the Serre functor is intrinsic, it commutes with any
autoequivalence of $\perf \widetilde{X}$.

\smallskip

For the pseudo-canonical condition, we need a lemma about how the various
functors related to a divisor interact.

\begin{lem}
  Let ${\cal A}$ be an abelian category, let $G:{\cal A}\to {\cal A}$ be an
  autoequivalence, and let $\tau:G\to \text{id}$ be a natural
  transformation.  If ${\cal A}$ is generated by objects on which $\tau$ is
  surjective or cogenerated by objects on which $\tau$ is injective, then
  $\tau G = G\tau$.
\end{lem}

\begin{proof}
  The ``cogenerated'' case is essentially just
  \cite[Lem.~8.3]{vanGastelM/VandenBerghM:1997} (using the cogenerators in
  place of injective objects), and the ``generated'' case is dual.
\end{proof}

\begin{rem}
  In particular, if $X$ is generated by weak line bundles along $C$, then
  this applies to the autoequivalence ${-}(-C)$ and the natural
  transformation ${-}(-C)\to\text{id}$.  We may thus safely refer to ``the''
  natural map $M(aC)\to M(bC)$ for any $a<b$, since all ways of
  constructing such a map out of the natural transformation will agree.
\end{rem}

\begin{lem}
  Let $C$ be a divisor on the noncommutative surface $X$ and suppose that
  $X$ is generated by weak line bundles along $C$ and has a Serre functor.
  Let $i_*:\qcoh C\to \qcoh X$ be the inclusion, with adjoints $i^*$,
  $i^!$.  Then there is a natural isomorphism $Ri^!\cong Li^*({-}(C))[-1]$.
  Furthermore, for $M\in D_{\qcoh} C$, $S (i_* M)(C)\cong i_* SM[1]$.
\end{lem}

\begin{proof}
  By \cite[Lem.~8.1]{vanGastelM/VandenBerghM:1997}, there is a 4-term exact
  sequence
  \[
  0\to i_*i^!M\to M\to M(C)\to i_*i^*(M(C))\to 0.
  \]
  Since there are enough acyclic objects for both functors (injectives for
  $i^!$ and weak line bundles for $i^*$), we find that the functors fit
  into functorial distinguished triangles
  \begin{align}
  {-}(-C)\to {}&\text{id}\to Li^*\to\\
  Ri^!\to {}&\text{id}\to {-}(-C)\to,
  \end{align}
  from which the first claim follows.  

  Serre duality immediately gives
  \[
  Ri^! \cong S Li^* S^{-1},
  \]
  and thus comparing the two expressions of $Ri^!$ gives a natural isomorphism
  \[
  Li^*(S^{-1}M(-C))\cong S^{-1} Li^* M[-1].
  \]
  Taking right adjoints gives $S(i_*M)(C)\cong i_* SM[1]$ as required.
\end{proof}

\begin{rem}
  The fact that the twist by $C$ is an autoequivalence identifying the two
  adjoints of $i_*$ establishes that $i_*$ is a {\em spherical functor} in
  the sense of \cite{AnnoR/LogvinenkoT:2017}.
\end{rem}

\begin{prop}
  If $X$ and $C$ are Gorenstein, then so is $\widetilde{X}$, and the functor
  $R\alpha_*:D^b_{\qcoh} \widetilde{X}\to D^b_{\qcoh} X$ is pseudo-canonical.
\end{prop}

\begin{proof}
  Express the Serre functor on $X$ by $\theta[2]$ with $\theta$ an exact
  functor.  The autoequivalence $M\mapsto \theta M(C)$ takes the sheaf
  $\sO_{\tau p}=i_*\sO_{\tau p}$ to $i_*(\sO_{\tau p}\otimes \omega_C)\cong
  \sO_{\tau p}$, so respects the gluing, and thus induces an
  autoequivalence of $D^b_{\qcoh} \widetilde{X}$ acting trivially on
  $\sO_e(-1)$.  Composing this autoequivalence with $M\mapsto M(-1)$ gives
  an autoequivalence acting by
  \[
  L\alpha^* M \to (L\alpha^* \theta M(C))(-1)
  \cong S(L\alpha^! M(C))(-1)[-2]
  \cong S L\alpha^*M[-2]
  \]
  and
  \[
  \sO_e(-1)\to \sO_e(-2)\cong S\sO_e(-1)[-2],
  \]
  which must therefore actually agree with the shifted Serre functor on
  $\widetilde{X}$ (which we also denote $\theta$).  We then find
  \[
  R\alpha_* \theta^{-d} M
  \cong
  \theta^{-d} (R\alpha_* M(d))(-dC).
  \]
  Since $\theta$ and $M\mapsto M(C)$ are exact on $X$, we conclude that
  $\theta^{-1}$ is relatively ample for $R\alpha_*$ iff $M\mapsto M(1)$ is
  relatively ample for $R\alpha_*$, and the latter holds by construction.
\end{proof}

The construction of van den Bergh depends on a choice of divisor containing
$p$, but it follows immediately from the semiorthogonal decomposition and
the pseudo-canonical property that the result does not actually depend on
the divisor.  For quasi-ruled surfaces, there is a natural choice of such
divisor, since the moduli scheme of points is embedded as a divisor.  It is
thus worth noting that this is preserved under blowup.

\begin{prop}\label{prop:nice_divisor_on_blowup}
  Suppose $X$ is a noncommutative surface and $C$ is a curve embedded as a
  divisor, and let $\widetilde{X}$ be the blowup of $X$ in a point $p\in C$.
  Then there is a curve $C^+$ embedded in $\widetilde{X}$ as a divisor via a
  natural transformation ${-}(-1)\to \text{id}$, which is Gorenstein
  whenever $C$ is Gorenstein.  If $C$ is anticanonical on $X$, so is $C^+$
  on $\widetilde{X}$.
\end{prop}

\begin{proof}
  The graded bimodule algebra associated to the blowup is such that the
  homogeneous component of degree $i$ is a sub-bimodule of the homogeneous
  component of degree $i+1$, which gives the desired natural
  transformation.  We need to show that the quotient is a commutative
  curve.  If $\tau(p)\ne p$, then the quotient is just the curve
  $\widetilde{C}$, and thus the claim holds, and $C^+\cong C$, so the
  Gorenstein claim is immediate.

  If $\tau(p)=p$, then we can still factor the natural transformation as
  \[
  {-}(-1)\to {-}(-\widetilde{C})\to \text{id},
  \]
  so that the (possibly noncommutative) scheme $C^+$ is isomorphic to $C$
  away from $p$.  We thus reduce to showing that if we quotient by any
  power of the bimodule ideal associated to $p$, then the result is
  commutative, and since every homogeneous component of the result is an
  extension of $\sO_p$, we may pass (via
  \cite{vanGastelM/VandenBerghM:1997}) to a local calculation in a ring of
  the form $R=k\langle\langle x,y\rangle\rangle/(uy-vx)$ where $u,v$
  generate the maximal ideal ${\frakm}$, in such a way that the point
  sheaves of $C$ correspond to maximal ideals of a commutative quotient
  $R/U$ where $U$ is normalizing.  The blowup is then locally given by the
  Proj of the Rees algebra $\bigoplus_i {\frakm}^i U^{-i} t^i$ where $t$
  is a central variable of degree 1.  (Note that since conjugation by $U$
  is an invertible automorphism of $R$, it makes sense to adjoin an inverse
  of $U$ to $R$.)  We need to show that the quotient by $t$ is a
  commutative curve.

  Since multiplying the degree $0$ element $U$ by any element of degree $1$
  gives a multiple of $t$ (in fact, an element of ${\frakm}t$), we find
  that $U$ is contained in the saturation of the ideal generated by $t$.
  Thus the objective is to show that the scheme of points is the $\Proj$ of
  \[
  R/\langle U\rangle \oplus \bigoplus_{i\ge 1} ({\frakm}^i/{\frakm}^{i-1} U)
  U^{-i}t^i.
  \]
  The corresponding $\Z$-algebra are unchanged if we remove the factors
  $U^{-1}$ from these algebras, so that we may as well consider instead the
  algebra
  \[
  R/\langle U\rangle \oplus \bigoplus_{i\ge 1} ({\frakm}^i/{\frakm}^{i-1}
  U)t^i.
  \]
  This is the $R/\langle U\rangle$-algebra generated by ${\frakm}/\langle
  U\rangle t$ subject to the single relation $(ut)(yt)=(vt)(xt)$, and thus
  to show that its $\Proj$ is the category of coherent sheaves of a
  commutative scheme, it will suffice to find an $R/\langle
  U\rangle$-linear automorphism of the degree 1 submodule that twists this
  into a commutator.

  In other words, we need an $R/\langle U\rangle$-module automorphism of
  ${\frakm}/\langle U\rangle$ that takes $x$ to $u$ and $y$ to $v$.  By
  \cite[Prop.~7.1]{vanGastelM/VandenBerghM:1997}, there is a short exact
  sequence
  \[
  \begin{CD}
    0@>>> R @> (v,-u) >> R^2 @> (x,y) >> {\frakm}\to 0
  \end{CD}
  \]
  of $R$-modules, and thus there is an $R$-linear morphism
  $\phi^+:{\frakm}\to {\frakm}/\langle U\rangle$ taking $x$ to $u$ and $y$
  to $v$ iff the corresponding map on $R^2$ annihilates $(v,-u)$.  Since
  this has image $vu-uv\in \langle U\rangle$, such a map indeed exists and
  is unique.  Moreover, the map on $R^2$ takes $(y,x)$ to $yu-xv\in
  uy-vx+\langle U\rangle=\langle U\rangle$, so that $\phi^+$ annihilates
  $\langle U\rangle$, so induces a unique $R$-linear morphism
  $\phi:{\frakm}/\langle U\rangle\to {\frakm}/\langle U\rangle$.  The same
  reasoning shows that there is a unique map of right $R$-modules taking
  $(u,v)$ to $(x,y)$, which since $R/\langle U\rangle$ is commutative must
  be the inverse of $\phi$.

  We may thus twist by $\phi$ to obtain a description of this scheme as the
  $\Proj$ of the free commutative graded $R/\langle U\rangle$-algebra
  generated by the maximal ideal of $R/\langle U\rangle$.  Since $R/\langle
  U\rangle\cong \hat\sO_C$, we may consider this in entirely commutative
  terms, and find that the local scheme of points is obtained by blowing up
  the origin of $\Spec k[[x,y]]$, taking the total transform of
  $\hat\sO_C$, and removing a single copy of the exceptional divisor.  In
  particular, the new scheme is indeed a Gorenstein curve when $C$ is
  Gorenstein.

  For the final claim, we note that the endofunctor $M\mapsto \theta
  M(C^+)$ of $D_{\qcoh}(X)$ is the endofunctor induced via the
  semiorthogonal decomposition by the pair $(M\mapsto \theta
  M(C),\text{id})$, and is thus trivial iff $C$ is anticanonical.
\end{proof}

\begin{rems}
  If one can embed $C$ as a Cartier divisor in a smooth projective surface
  $X'$, then it follows from the above local description that $C^+$ is
  naturally isomorphic to the effective Cartier divisor $\pi^*C-e$, where
  $\pi:\widetilde{X}'\to X'$ is the blowup at $\tau p$.  Note that $C$ being
  Gorenstein implies only that such an embedding exists locally, but we
  have already seen that it exists globally for quasi-ruled surfaces, while
  for noncommutative planes, the curve of points embeds in a commutative
  plane.
\end{rems}

\begin{rems}
  We would like to have a statement along the lines that if $C$ is the
  moduli space of point sheaves on $X$, then $C^+$ is the moduli space of
  point sheaves on $\widetilde{X}$.  Unfortunately, there is no intrinsic
  notion of a ``point sheaf'' on a noncommutative surface, so this is
  difficult to state in general.  For the surfaces of present interest
  (iterated blowups of quasi-ruled surfaces and projective planes), we do
  have such a notion, and it follows easily from the above calculation that
  $C^+$ is the moduli space of such sheaves on $\widetilde{X}$ iff $C$ is
  the moduli space of point sheaves on $X$.  (This involves only local
  calculations near the exceptional curve.)
\end{rems}


We also need to know that generation by weak line bundles is inherited, so
that further blowups still have the semiorthogonal decomposition.

\begin{lem}
  If $X$ is generated by weak line bundles relative to $C$,
  then $\widetilde{X}$ is generated by weak line bundles relative to $C^+$.
\end{lem}

\begin{proof}
  If $M$ is a weak line bundle on $X$ relative to $C$, then $\alpha^*M$ is
  a weak line bundle relative to the new curve $C^+$, and $\theta$ takes
  weak line bundles to weak line bundles.  Since $\theta^{-1}$ is
  relatively ample for $\alpha_*$, $\qcoh \widetilde{X}$ is generated by sheaves
  of the form $\theta^{-d}\alpha^*M$ where $M$ ranges over any set of
  generators of $\qcoh X$.
\end{proof}

For the next result, we assume not just that $C$ is Gorenstein, but that it
is a surface curve, i.e., that it embeds as a Cartier divisor in a
commutative smooth projective surface.

\begin{lem}\label{lem:blowup_generates_correct_curve}
  Let $X$ be a compactly generated noncommutative surface and $i:C\to X$
  an embedding of a surface curve as a divisor.  Let $\widetilde{X}$ be a
  blowup of $X$ in some point of $C$.  If the image of $\perf(X)$ generates
  $\perf(C)$, then the image of $\perf(\widetilde{X})$ generates $\perf(C^+)$.
\end{lem}

\begin{proof}
  Let $X'$ be a smooth projective surface containing $C$.  The line bundle
  $\sO_C(d)$ is certainly in the image of $\perf(X')$ for any $d$, and
  since $\sO_C(1)$ is ample, these bundles generate $\perf(C)$.  By the
  same argument, $\perf(C^+)$ is generated by the image of
  $\perf(\widetilde{X}')$.  By the semiorthogonal decomposition,
  $\perf(\widetilde{X}')$ is generated by $L\alpha^*\perf(X')$ and
  $\sO_e(-1)$, and thus $\perf(C^+)$ is generated by the image of
  $L\alpha^*\perf(X')$ and the image of $\sO_e(-1)$.  The result follows by
  noting that the Karoubian subcategory generated by the image of
  $L\alpha^*\perf(X')$ in $\perf(C^+)$ is the same as the image of
  $\perf(C)$ in $\perf(C^+)$, and the image of $\sO_e(-1)$ in $\perf(C^+)$
  is the same for either $\widetilde{X}$ and $\widetilde{X}'$.
\end{proof}

In particular, both noncommutative planes and noncommutative quasi-ruled
surfaces have the property that any point can be blown up, and this
property is inherited by such blowups.  This motivates the following two
definitions

\begin{defn}
  A {\em noncommutative rationally (quasi-)ruled surface} is an iterated
  blowup of a noncommutative (quasi-)ruled surface, while a {\em
    noncommutative rational surface} is an iterated blowup of either a
  noncommutative plane or a noncommutative ruled surface over a smooth
  curve of genus 0.  (For technical reasons, when the surface is
  commutative, we include the curve of points as part of the data.)
\end{defn} 

One could similarly define a noncommutative {\em birationally} quasi-ruled
surface to include any noncommutative rationally quasi-ruled surface and
any surface such that an iterated blowup is birationally quasi-ruled, but
we will see in Theorem \ref{thm:birationally_qr_is_rationally_qr} below
that this adds no additional generality, and similarly for the rational
case.

\section{Commuting blowups}

Thus in each case, we could have constructed the desired abelian category
by first gluing the appropriate commutative derived categories and using
the pseudo-canonical property to construct the $t$-structure.  Of course,
it is nontrivial to show that this is a $t$-structure and even less trivial
to show that the resulting triangulated category is the derived category of
its heart.  So although these constructions are relatively simple to state,
they are difficult enough to control to make this a less than ideal
definition.  It is still useful to think about the categories this way,
however, as it leads to rather different approaches to constructing
morphisms.

In particular, there are a number of isomorphisms we would like to
establish, analogous to standard constructions in commutative birational
geometry: blowups in (sufficiently) distinct points should commute, a
blowup of a quasi-ruled surface should be a blowup of a different
quasi-ruled surface (i.e., elementary transformations), a blowup of a
noncommutative plane should be a noncommutative Hirzebruch surface, and a
noncommutative $\P^1\times \P^1$ should be a noncommutative Hirzebruch
surface in two different ways.  In each case, we will show this by first
constructing a corresponding {\em derived} equivalence, and then showing
that the derived equivalence respects the $t$-structure.  The latter will
in turn reduce to showing that the derived equivalence respects some
pseudo-canonical functor.

For commuting blowups, we note that if $X$ is a noncommutative Gorenstein
surface with successive blowups $\alpha_{1*}:X_1\to X$ and
$\alpha_{2*}:X_2\to X_1$, then there are three natural subcategories
$\sO_{e_2}(-1)\otimes D^b_{\coh} k$, $\alpha_2^*\sO_{e_1}(-1)\otimes
D^b_{\coh} k$ and $\alpha_2^*\alpha_1^*D^b_{\coh} X$ of $D^b_{\coh}
X_2$, and these moreover form a three-step semiorthogonal
decomposition.  Note that
\[
R\Hom(\sO_{e_2}(-1),\alpha_2^*\sO_{e_1}(-1))
\cong
R\Hom(\sO_{e_1}(-1),\sO_{p_2})^*,
\]
and thus if $R\Hom(\sO_{e_1}(-1),\sO_{p_2})=0$, we can swap the first two
subcategories and still have a semiorthogonal decomposition.  We find that
the new semiorthogonal decomposition is again the decomposition associated
to a two-step blowup, but now with $p_1$ and $p_2$ swapped.  (Note that if
$p_1$ is a singular point of $Q$, then the condition implies that $p_2$ is
not on the component $e_1$ of $Q_1=Q^+$, and thus it makes sense to swap
$p_1$ and $p_2$.)  The resulting derived equivalence fixes the functor
$\alpha_{1*}\alpha_{2*}$, and thus when this functor is pseudo-canonical,
the derived equivalence respects the $t$-structure.

To see that this is pseudo-canonical (in fact relatively Fano in a suitable
sense) under reasonable conditions, it will be helpful to have a way of
checking in the case of a single blowup whether a sheaf is acyclic and
globally generated for $\alpha$.  We assume the original surface is
Gorenstein (which holds in all cases of interest in the present
work), but with care one should be able to replace $\theta$ by suitable
functors coming from the relevant divisors of points.

\begin{prop}\label{prop:blowup_acyclic}
  Let $X$ be a noncommutative surface and $\alpha:\widetilde{X}\to
  X$ the van den Bergh blowup of $X$ in the point $x$.  Then the object
  $M\in \qcoh \widetilde{X}$ is acyclic and globally generated for $\alpha$ iff
  $\Ext^2(\sO_e,M)=0$.
\end{prop}

\begin{proof}
  We first recall that the distinguished triangle associated to the
  semiorthogonal decomposition takes the form
  \[
  L\alpha^*R\alpha_*M\to M\to FM\otimes_k \sO_e(-1)\to
  \]
  for some functor $F:D_{\qcoh} \widetilde{X}\to D_{\qcoh} k$, and we claimed
  above that $FM=R\Hom(\sO_e,M[2])$.  To see this, note that since
  $R\Hom(\sO_e,L\alpha^*N)=0$, we have
  \[
  R\Hom(\sO_e,M[2])\cong FM\otimes_k R\Hom(\sO_e,\sO_e(-1)[2])
  \cong FM.
  \]
  Substituting into the distinguished triangle and taking the corresponding
  long exact sequence gives an exact sequence ending
  \[
      h^0(L\alpha^*R\alpha_*M)
  \to M
  \to \Ext^2(\sO_e,M)\otimes_k \sO_e(-1)
  \to \alpha^* R^1\alpha_*M
  \to 0.
  \]
  In particular, if $\Ext^2(\sO_e,M)=0$, then $\alpha^* R^1\alpha_*M=0$
  and thus $R^1\alpha_*M=0$ so $M$ is acyclic for $\alpha_*$.

  In general, if $M$ is acyclic for $\alpha_*$, then the tail of the exact
  sequence becomes
  \[
  \alpha^*\alpha_*M\to M \to \Ext^2(\sO_e,M)\otimes_k \sO_e(-1)\to 0
  \]
  and thus $\alpha^*\alpha_*M\to M$ is surjective iff
  $\Ext^2(\sO_e,M)=0$.
\end{proof}

\begin{rems}
  Of course, if $M$ is coherent and $X$ has a Serre functor, then we may
  use Serre duality to rephrase the condition as $\Hom(M,\sO_e(-1))=0$.
\end{rems}

\begin{rems}
  The distinguished triangle from the $\sO_e(-1),L\alpha^!$ decomposition
  similarly tells us that
  \[
  \alpha^!R^1\alpha_*M\cong \Ext^2(\sO_e(-1),M)\otimes\sO_e(-1),
  \]
  so that $M$ is acyclic for $\alpha_*$ iff $\Ext^2(\sO_e(-1),M)=0$.  This
  can be viewed as a relative Castelnuovo-Mumford regularity statement,
  which in this case is an equivalence: relative to $\alpha$, $M$ is
  acyclic and globally generated iff $\theta M$ is acyclic.
\end{rems}

\begin{prop}
  If $M\in D^b_{\qcoh} X$ is such that $R\Hom(\sO_{\tau p},M)=0$, then
  $L\alpha^*M\cong L\alpha^!M$.
\end{prop}

\begin{proof}
  We have a distinguished triangle
  \[
  L\alpha^* M\to L\alpha^!M\to R\Hom(\sO_e,L\alpha^!M[2])\otimes \sO_e(-1)\to,
  \]
  and thus it suffices to have $R\Hom(\sO_e,L\alpha^!M)=0$.  But then
  \[
  R\Hom(\sO_e,L\alpha^!M)
  \cong
  R\Hom(L\alpha^*\sO_{\tau p},L\alpha^!M)
  \cong
  R\Hom(\sO_{\tau p},M),
  \]
  so that the desired result follows.
\end{proof}

\begin{rem}
  Of course, if $X$ is Gorenstein and $M$ is coherent, then the hypothesis
  is equivalent to $R\Hom(M,\sO_p)=0$.
\end{rem}
  
\begin{prop}
  Let $X_0$ be a noncommutative Gorenstein surface, and let $\alpha:X_1\to
  X_0$, $\beta:X_2\to X_1$ be a sequence of van den Bergh blowups.
  If the images of the corresponding points $x_1$, $x_2$ in $X_0$ are not
  in the same orbit under $\tau^{\Z}$, then $\alpha_*\beta_*$ is
  pseudo-canonical.
\end{prop}

\begin{proof}
  Let $M\in \coh X_2$ be a nonzero sheaf.  We need to show that for some
  $l$, $\alpha_*\beta_* \theta^{-l} M\ne 0$.  It will of course suffice to
  find $l$ such that $\theta^{-l} M$ is globally generated relative to the
  composed morphism.  With this in mind, let $b(M)$ be the function that
  assigns to each $M$ the smallest integer such that
  $\Ext^2(\sO_{e_2},\theta^{-b}M)=0$, and if $b(M)\le 0$, let $a(M)$ be the
  smallest integer such that $\Ext^2(\sO_{e_1},\theta^{-a}\beta_*M)=0$.  If
  $b(M)\le 0$ and $a(M)\le 0$, then we find that $M$ is acyclic and
  globally generated for $\beta_*$, while $\beta_*M$ is acyclic and
  globally generated for $\alpha$, which implies that $M$ is acyclic and
  globally generated for $\alpha\circ\beta$ as required.  It will thus
  suffice to show that $b(\theta^{-1}M)=b(M)-1$ and if $b(M)\le 0$, then
  $a(\theta^{-1}M)=a(M)-1$.

  The claim for $b(M)$ is trivial, so it remains to consider the claim for
  $a(M)$.  We have
  \[
  R\Hom(\sO_{e_1},\theta^{-a}\beta_*\theta^{-1}M)
  \cong
  R\Hom(\sO_{e_1},\theta^{-a}R\beta_*\theta^{-1}M)
  \cong
  R\Hom(L\beta^!\theta^a \sO_{e_1}(-1),M)
  \]
  If $R\Hom(\theta^a \sO_{e_1}(-1),\sO_{x_2})=0$, then we have
  \begin{align}
  R\Hom(L\beta^!\theta^a \sO_{e_1}(-1),M)
  &\cong
  R\Hom(L\beta^*\theta^a \sO_{e_1}(-1),M)\notag\\
  &\cong
  R\Hom(\theta^a \sO_{e_1}(-1),R\beta_* M)\notag\\
  &\cong
  R\Hom(\sO_{e_1},\theta^{-a-1} R\beta_* M),
  \end{align}
  and thus the claim follows.

  It thus remains only to show that
  $R\Hom(\theta^a\sO_{e_1}(-1),\sO_{x_2})=0$, which we may rewrite using
  \[
  R\Hom(\theta^a \sO_{e_1}(-1),\sO_{x_2})
  \cong
  R\Hom(\sO_{e_1}(-1),\sO_{\tau^{-a}x_2}).
  \]
  By assumption, $e_1$ is either not a component of $C_1^+$ or does not
  contain $x_2$, and thus we may choose the divisor $C_2$ containing $x_2$
  so that $e_1$ is not a component of $C_2$.  We then have
  \[
  R\Hom(\sO_{e_1}(-1),\sO_{\tau^{-a}x_2})
  \cong
  R\Hom_{C_2}(\sO_{e_1}(-1)|_{C_2},\sO_{\tau^{-a}x_2}).
  \]
  Now, $\sO_{e_1}(-1)|_{C_2}$ either vanishes (if $x_1$ and $x_2$ came from
  different components of the original curve of points) or equals $x_1$,
  in which case
  \[
  R\Hom(\sO_{e_1}(-1),\sO_{\tau^{-a}x_2})
  \cong
  R\Hom_{C_2}(\sO_{x_1},\sO_{\tau^{-a}x_2})
  \]
  vanishes unless $x_2=\tau^a x_1$.
\end{proof}

We then have the following immediate consequence.

\begin{thm}
  Let $X$ be a Gorenstein surface containing a divisor $C$, and let $x_1$,
  $x_2$ be a pair of points of $C$ which are not in the same orbit under
  $\tau$.  Then the blowup of $X$ in $x_1$ then $x_2$ (as a point in
  $\widetilde{C}$) is isomorphic to the blowup of $X$ in $x_2$ then $x_1$.
\end{thm}

\smallskip

It will be useful to have an understanding of what happens when $x_1$ and
the image of $x_2$ {\em are} in the same orbit.  If $x_1$ is a fixed point,
then $e_1$ will be an actual divisor in the first blowup $X_1$, and we can
use its strict transform to control things.  So suppose that $x_1$ is not a
fixed point, but rather that it has an orbit of size $r\in [2,\infty]$, and
thus that the curve of points on $X_1$ locally agrees with $C$, and
similarly (since we are assuming $x_2$ in the same orbit) for the two-fold
blowup.  Applying a suitable equivalence of the form ${-}(de_2)$ then lets
us assume that $x_2$ is equal to $x_1=:x$.  Then there is a natural
morphism $\sO_{e_2}(-1)\to \sO_{e_1}(-1)$, and we define
$\sO_{e_1-e_2}(-1)$ to be its cokernel.  Note that since $\sO_{e_1}(-1)$
and $\sO_{e_2}(-1)$ both restrict to $\sO_x$ on $C$, the cone of this
morphism has trivial restriction to $C$, and thus is invariant under the
autoequivalence ${-}(C)$.  Since this autoequivalence is relatively ample
for $\beta_*$ and the direct image of the cone is the sheaf
$\sO_{e_1}(-1)$, we conclude that the morphism is in fact injective.  Note
also that we have $\theta\sO_{e_1-e_2}(-1)\cong \sO_{e_1-e_2}(-1)(-C)\cong
\sO_{e_1-e_2}(-1)$.  This object is readily verified to be spherical, and
thus gives rise to an inverse pair of spherical twists.  (This of course
also makes sense when $r=1$, taking the appropriate line bundle on the
strict transform of $e_1$.)

The size of the orbit plays a role via the following fact.

\begin{lem}
  Let $\widetilde{X}$ be the blowup of $X$ in a point $x$ with orbit of size
  $r\in [1,\infty]$.  Then for all $0\le d<r$, one has
  \[
  \Ext^1_{\widetilde{X}}(\sO_{e_1}(-1),\sO_{e_1}(d-1))=
  \Ext^2_{\widetilde{X}}(\sO_{e_1}(-1),\sO_{e_1}(d-1))=0
  \]
  and $\dim\Hom_{\widetilde{X}}(\sO_{e_1}(-1),\sO_{e_1}(d-1))=1$.
\end{lem}

\begin{proof}
  This holds for $d=0$ since $\sO_{e_1}(-1)$ is exceptional.  We proceed by
  induction in $d$.  Using the short exact sequence
  \[
  0\to \sO_{e_1}(d-1)\to \sO_{e_1}(d)\to \sO_{x_d}\to 0
  \]
  (where $x_d$ is an appropriate point of $\widetilde{C}$), we find that
  \[
  R\Hom(\sO_{e_1}(-1),\sO_{e_1}(d-1))\cong
  R\Hom(\sO_{e_1}(-1),\sO_{e_1}(d))
  \]
  unless
  $R\Hom(\sO_{e_1}(-1),\sO_{x_d})\ne 0$, or equivalently unless
  $x_d=x_{-1}$. Since $x_d$ ranges over consecutive points in the orbit of
  $x$, this fails only when $d+1$ is a multiple of $r$.
\end{proof}

Define a sheaf $\sO_{X_2}(ae_1+be_2)$ by
\[
\sO_{X_2}(ae_1+be_2):=
\theta^b \beta^* \theta^{a-b} \alpha^* \theta^{-a} \sO_X.
\]
(This could of course be defined using ${-}(-C)$ in place of $\theta$; the
use of $\theta$ is only for notational convenience.)  Our key result
for use below can be viewed as computing certain special cases of the
spherical twists of such sheaves.

\begin{prop}\label{prop:weak_reflection}
  If $0\le b-a\le r$, then there is a short exact sequence
  \[
  0\to \sO_{X_2}(-be_1-ae_2)\to \sO_{X_2}(-ae_1-be_2)\to
  \sO_{e_1-e_2}(-1)^{b-a}\to 0
  \]
\end{prop}

\begin{proof}
  We show this by induction in $b-a$, noting that the base case $a=b$ is
  trivially true.  Since the functor ${-}(e_1+e_2)$ preserves
  $\sO_{e_1-e_2}(-1)$ (it agrees with ${-}(-C)$ and $\theta$ on $D^b_{\coh}
  X^{\perp}$), we may as well assume that $a=0$.  We thus need to show that
  if the claim holds for some $b<r$, then it holds for $b+1$.

  There is certainly a natural injection $\sO_{X_2}(-(b+1)e_1)\to
  \sO_{X_2}(-be_1)$, and composing with the already constructed map
  $\sO_{X_2}(-be_1)\to \sO_{X_2}(-be_2)$ gives an injection such that the
  cokernel $M_b$ is an extension of $\sO_{e_1-e_2}(-1)^b$ by
  $\beta^*\sO_{e_1}(b)$.  But
  \begin{align}
  \Ext^1(\sO_{e_1-e_2}(-1),\beta^*\sO_{e_1}(b))
  &\cong
  \Ext^1(\theta^{-1}\sO_{e_1-e_2}(-1),\beta^*\sO_{e_1}(b))\notag\\
  &\cong
  \Ext^1(\sO_{e_1-e_2}(-1),\beta^!\theta\sO_{e_1}(b))\notag\\
  &\cong
  \Ext^1_{X_1}(\sO_{e_1}(-1),\sO_{e_1}(b-1))\notag\\
  &=0,
  \end{align}
  and thus this extension splits.  Since
  \[
  \Hom(\sO_{e_1-e_2}(-1),\sO_{e_2}(b))\cong
  \Hom(\sO_{e_1-e_2}(-1),\sO_{e_2}(-1))=0
  \]
  and
  \[
  \Hom(\beta^*\sO_{e_1}(b),\sO_{e_2}(b))
  \cong
  \Hom(\sO_{e_1}(b),\beta_*\sO_{e_2}(b))
  \cong
  k
  \]
  (using the fact that $b<r$, so $\beta_*\sO_{e_2}(b)$ is the direct sum of
  $x_0$,\dots,$x_b$, only one of which has a map from $\sO_{e_1}(b)$), we
  see that there is a unique morphism from $M_b$ to $\sO_{e_2}(b)$, and
  thus the map $\sO_{X_2}(-(b+1)e_1)\to \sO_{X_2}(-be_2)$ factors through
  $\sO_{X_2}(-(b+1)e_2)$.  The cokernel of the resulting map is the sum of
  $\sO_{e_1-e_2}(-1)^b$ and the cokernel of the map $\beta^*\sO_{e_1}(b)\to
  \sO_{e_2}(b)$.  Since $b<r$, we have
  \[
  \beta^*\sO_{e_1}(b)\cong \theta^{-1}\beta^!\theta \sO_{e_1}(b)
  \cong \theta^{-1}\beta^*\sO_{e_1}(b-1),
  \]
  and thus this cokernel is $\theta^{-1}$ of the cokernel for $b-1$, and
  thus by induction is $\theta^{-b-1}\sO_{e_1-e_2}(-1)\cong
  \sO_{e_1-e_2}(-1)$.
\end{proof}

\begin{rems}
  This construction of a morphism $\sO_{X_2}(-be_1-ae_2)\to
  \sO_{X_2}(-ae_1-be_2)$ is somewhat reminiscent of the construction of
  (multivariate) operators of degree $d(s-f)$ in \cite[Prop.~8.7]{elldaha}.
  This is not a coincidence, as if we blow up a point on a noncommutative
  $\P^1\times \P^1$, the result is isomorphic to a two-point blowup of a
  noncommutative plane, and the commutation of blowups symmetry on the
  latter is the same as the exchange of rulings symmetry on the former.
\end{rems}

\begin{rems}
  When $X$ is rational or rationally quasi-ruled, we can replace $\sO_X$ by
  any line bundle (see Definition \ref{defn:line_bundle} below) by applying
  the corresponding twist autoequivalence before blowing up.  If
  $r=\infty$, we find that every line bundle on $X_2$ fits into such an
  exact sequence (and thus, modulo checking that the other $\Ext$ groups
  vanish, we can compute one of its two spherical twists).  If $r$ is
  finite, this fails, but every line bundle fits into a sequence along the
  above lines with $\sO_{e_1-e_2}(-1)$ replaced by some
  $\sO_{e_1-e_2}(-1-dr)$.
\end{rems}

\section{More semiorthogonal decompositions}

For the remaining cases, showing the pseudo-canonical property will require
some additional facts about sheaves, so we will postpone that to a future
section and consider only the construction of the derived equivalences.
The $\P^1\times \P^1$ case can be dealt with similarly, but the other two
cases (elementary transformations and $F_1$ as a blowup of $\P^2$) require
more complicated modifications of the semiorthogonal decomposition.  To
deal with these cases, it will be helpful to consider different
semiorthogonal decompositions.  In the rational case, the object $\sO_X$ is
exceptional, and thus induces a semiorthogonal decomposition
$(\sO_X^\perp,\langle \sO_X\rangle)$, while for blowups of quasi-ruled
surfaces, the image of $D^b_{\coh} C_0$ is admissible, so induces a
decomposition.  It turns out that in each case, the orthogonal subcategory
is also essentially commutative, in that it appears in the same way as a
subcategory of the derived category of a commutative surface.  Moreover,
the gluing data can also be expressed in commutative terms.

A key observation is that when $X$ has an anticanonical curve $C$, the gluing
data can be expressed as an $R\Hom$ in $D^b_{\coh} C$.  Let $i_*:D^b_{\coh}
C\to D^b_{\coh} X$ be the embedding, with left adjoint $Li^*$ and
associated twist ${-}(-C)$.

\begin{lem}
  Let $X$ be a noncommutative surface with an anticanonical divisor $C$.
  If $M,N\in D^b_{\coh} X$ are such that $R\Hom(M,N)=0$, then $R\Hom(N,M)\cong
  R\Hom_C(Li^* N,Li^* M)$.
\end{lem}

\begin{proof}
  By Serre duality, we have $R\Hom_X(N,M(-C)[2])\cong R\Hom_X(M,N)=0$, so that
  \[
  R\Hom_X(N,M)\cong R\Hom_X(N,i_*Li^*M)\cong R\Hom_C(Li^*N,Li^*M)\notag
  \]
  as required.
\end{proof}

We need a slight variation of this in order to deal with iterated blowups
of quasi-ruled surfaces.

\begin{lem}
  Let $X$ be a noncommutative surface with a Serre functor and a curve $C$
  embedded as a divisor, and suppose that ${\cal A}$ is a full subcategory
  of $D^b_{\coh} X$ preserved by $M\mapsto S M(C)$.  Then for any $A\in {\cal
    A}$, $B\in {\cal A}^\perp$, $B'\in {}^\perp{\cal A}$
  \begin{align}
    R\Hom(B,A)&\cong R\Hom_C(Li^*B,Li^*A)\\
    R\Hom(A,B')&\cong R\Hom_C(Li^*A,Li^*B').
  \end{align}
\end{lem}

\begin{proof}
  For the first claim, it again suffices to show that $R\Hom(B,A(-C))=0$.
  By Serre duality, this is isomorphic to
  $R\Hom(A,SB(C))=R\Hom(S^{-1}A(-C),B)=0$ since $S^{-1}A(-C)\in {\cal A}$.
  Similarly, $R\Hom(A,B'(-C))\cong R\Hom(B',SA(C))=0$.
\end{proof}

\begin{rem}
  Note that if $M\mapsto SM(C)$ preserves ${\cal A}$, then it also
  preserves ${\cal A}^\perp$ and ${}^\perp{\cal A}$.
\end{rem}

For a noncommutative rational surface, the curve $Q$ of points is
anticanonical, and thus $M\mapsto SM(Q)$ preserves $\langle \sO_X\rangle$.
For an iterated blowup of a quasi-ruled surface, this can fail, but the
functor always preserves the category $L\rho_0^*D^b_{\coh} C_0$.

\begin{thm}
  Suppose the noncommutative surface $X$ is an iterated blowup of either a
  noncommutative plane or a noncommutative Hirzebruch surface, with
  anticanonical curve $Q$.  Then there is a commutative surface $X'$, an
  embedding $Q\subset X'$ as an anticanonical curve, and a point $q\in
  \Pic^0(Q)$ such that there is an equivalence $\kappa:\sO_X^\perp\cong
  \sO_{X'}^\perp$ satisfying
  \[
  R\Hom_X(M,\sO_X)\cong R\Hom_{X'}(\kappa(M),q)\cong R\Hom_Q(\kappa(M)|_Q,q).
  \]
\end{thm}

\begin{proof}
  The construction of $X$ as an iterated blowup gives rise to an
  exceptional collection
  \[
  \sO_{e_m}(-1),\alpha_m^*\sO_{e_{m-1}}(-1),\cdots,
  \]
  ending with an exceptional collection for the original surface, in turn
  ending with $\sO_X$.  The forward maps in this exceptional collection
  (from which we can reconstruct the full derived category via gluing) can
  be computed via the restrictions to $Q$.  What we need to prove is that
  there is a line bundle $q$ in the identity component of $\Pic(Q)$ such
  that if we tensor the restrictions with $q$ and then replace the last
  term by $\sO_Q$, then the result is the sequence corresponding to an
  exceptional collection in some commutative $X'$ with anticanonical curve $Q$.

  For $X$ a noncommutative plane, the exceptional collection is
  $\sO_X(-2)$, $\sO_X(-1)$, $\sO_X$, with restrictions ${\cal L}_{-6}$,
  ${\cal L}_{-3}$, $\sO_Q$ with ${\cal L}_{-6}$, ${\cal L}_{-3}$ line
  bundles of the given degree such that ${\cal L}_{-3}^{-1}$, ${\cal
    L}_{-6}^{-1}\otimes {\cal L}_{-3}$ are effective, acyclic, and
  numerically equivalent.  The commutative case is that ${\cal L}_{-6}\cong
  {\cal L}_{-3}^2$, and thus if we tensor everything with $q$ and replace
  the last term with $\sO_Q$, the result will corresponding to a
  commutative surface as long as $q\otimes {\cal L}_{-6}\cong q^2\otimes
  {\cal L}_{-3}^2$, or in other words when $q\cong {\cal L}_{-6}\otimes
  {\cal L}_{-3}^{-2}$.  The degree condition ensures that $q$ has degree 0
  on every component, so is in $\Pic^0(Q)$ as required.  Note that the
  commutative surface is the $\P^2$ in which $Q$ is embedded via the line
  bundle $q^{-1}\otimes {\cal L}_{-3}^{-1}$.

  For a noncommutative Hirzebruch surface, we note that $D^b_{\coh} \P^1$
  has an exceptional collection $\sO_{\P^1}(-1),\sO_{\P^1}$, and thus we
  can refine the semiorthogonal decomposition of Theorem
  \ref{thm:semiorth_for_quasiruled} to an exceptional collection.
  Restricting to $Q$ gives a sequence
  \[
    {\cal L}'_{-2},{\cal L}',{\cal L}_{-2},\sO_X
  \]
  which is commutative whenever ${\cal L}'_{-2}\cong {\cal L}'\otimes
  {\cal L}_{-2}$.  Moreover, intersection theory forces $q:={\cal
    L}^{\prime{-}1}\otimes {\cal L}_{-2}^{-1}\otimes {\cal L}'_{-2}$ to be
  degree 0 on every component, and thus be in $\Pic^0(Q)$ as required.  The
  resulting commutative Hirzebruch surface is the one obtained by using
  ${\cal L}_{-2}^{-1}\otimes q^{-1}$ to determine a degree 2 morphism $\pi:Q\to
  \P^1$ and then taking the $\P^1$-bundle corresponding to $\pi_*({\cal
    L}^{\prime{-}1}\otimes q^{-1})$.

  Finally, if we blow up a point, the effect on the sequence of objects in
  $D^b_{\coh} Q$ is to (derived) pull back all of the existing objects to
  the new anticanonical curve $Q'$ and then prepend the appropriate object
  to the beginning of the sequence.  If $p$ is a smooth point of $Q$, this
  is simply a point sheaf, while in general it is given by the (perfect!)
  complex $\sO_{Q'}\to \sO_{Q'}(e)$, which can be computed inside the
  commutative surface.  This complex is trivial on the generic point of the
  blowup $\widetilde{Q}$, and is thus invariant under twisting by $q$,
  while for the other sheaves we note that tensoring by $q$ then pulling
  back is the same as pulling back then tensoring by the pullback of $q$
  (which is again of degree 0 on every component).  We thus find that the
  effect of blowing up $X$ is the same as that of blowing up the
  corresponding point of $X'$ and then pulling back $q$.
\end{proof}

\begin{rems}
  Given any triple $(X',Q,q)$ where $X'$ is a smooth projective rational
  surface, $Q\subset X$ is an anticanonical curve, and $q$ is a line bundle
  in the identity component of $\Pic(Q)$, there is a corresponding
  (dg-enhanced) triangulated category obtained by using $q$ to glue
  $\sO_{X'}^\perp$ to $\sO_{X'}$, such that taking $q=\sO_Q$ recovers
  $\sO_{X'}$.  Of course, to obtain the noncommutative rational surfaces
  themselves, one needs to determine the $t$-structure, which turns out not
  to be possible without imposing some additional structure on $X'$.  This
  can, however, be done inductively given a choice of blowdown structure on
  $X'$ (i.e., a sequence of monoidal transformations blowing it down to
  either a Hirzebruch surface or $\P^2$; see Section
  \ref{sec:blowdown_structures_comm}), by insisting that the corresponding
  functors be pseudo-canonical.  As usual, it appears to be hard to prove
  independently that this gives a $t$-structure, let alone that it is the
  derived category of its heart.  It follows, of course, from the above
  considerations, and thus we can use this to construct a family of
  noncommutative surfaces parametrized by the relative $\Pic^0$ of the
  universal anticanonical curve of the moduli stack of anticanonical
  rational surfaces with blowdown structure.
\end{rems}

\begin{rems}
  In the fully noncommutative case, the gluing data is again in the form of
  a pair of adjoint functors, since $q\in \sO_{X'}^\perp$.  On the other hand,
  when $q=\sO_Q$, the above description is no longer of that form.  This
  can be rectified, however, since for $M\in \sO_{X'}^\perp$, one has
  \[
  R\Hom_{X'}(M,\omega_{X'})\cong R\Hom_{X'}(\sO_X,M[2])^*=0,
  \]
  and thus
  \[
  R\Hom_{X'}(M,\sO_{X'})\cong R\Hom_{X'}(M,\cone(\Tr)[-1])
  \]
  where $\Tr$ is the trace map $\sO_{X'}\to \omega_{X'}[2]$.  Since
  $R\Hom(\sO_{X'},\Tr)$ is the identity and $\cone(\Tr)[-1]\in
  \sO_{X'}^\perp$, this induces the desired functor.
\end{rems}

\begin{rems}
  It follows from \cite[Thm. 8.5]{KuznetsovA:2009} that the Hochschild
  cohomology of $\sO_{X'}^\perp$ is the hypercohomology of $\sO_{X'}\oplus
  T_{X'}$, so that (since $\sO_{X'}$ is acyclic) the deformation theory of
  $\sO_X^\perp$ is given by the commutative deformation theory of $X'$.
  Similarly, the deformations of the gluing functor (i.e., the deformations
  of $X$ with $\sO_{X}^\perp$ fixed) are controlled by the deformations
  of the image of $q$ in $\sO_{X'}^\perp$.  Thus if $q\not\cong \sO_Q$,
  then the deformation theory is controlled by the hypercohomology of
  $\sO_Q\oplus \sO_Q(Q)$, while for $X=X'$, the deformation theory is
  controlled by the hypercohomology of $\wedge^2 T_{X'}\cong
  \omega_{X'}^{-1}$.  This implies at the very least that any first-order
  deformation of $\perf(X)$ is a noncommutative rational surface of the
  above form, and the agreement of obstruction spaces suggests that
  this should extend to arbitrary formal deformations.
\end{rems}

In particular, any isomorphism $X'\cong Y'$ of rational surfaces with
blowdown structures induces a derived equivalence between the corresponding
noncommutative surfaces, giving the desired derived equivalences
corresponding to $\P^1\times \P^1\cong \P^1\times \P^1$ and the
representation of $F_1$ as a blowup of $\P^2$.  In each case, the resulting
equivalence will be an abelian equivalence as long as the functor $R\Gamma$
is pseudo-canonical, since this functor is clearly preserved by the equivalence.

\begin{prop}
  The image of the restriction map $Li^*:\sO_X^\perp\to \perf(Q)$ generates
  $\perf(Q)$.
\end{prop}

\begin{proof}
  Without loss of generality, $X$ is a commutative rational surface.  As in
  Lemma \ref{lem:blowup_generates_correct_curve}, we find that
  $Li^*\perf(X)$ generates $\perf(Q)$, so it will suffice to show that
  $\sO_Q=Li^*\sO_X$ is in the Karoubian triangulated subcategory
  generated by $Li^*\sO_X^\perp$.  Define an element $M\in \sO_X^\perp$
  by the distinguished triangle
  \[
  M\to \sO_X\to \theta\sO_X[2]\to
  \]
  coming from the natural element of $\Ext^2(\sO_X,\theta \sO_X)$.
  Since $\Ext^2$ vanishes for line bundles on $Q$, this map is annihilated
  by $Li^*$, producing a splitting:
  \[
  Li^*M\cong \theta \sO_Q[1]\oplus \sO_Q.
  \]
  But this is what we needed to show.
\end{proof}

\smallskip

For quasi-ruled surfaces, we note the following.

\begin{lem}
  Let $X$ be a quasi-ruled surface over $C_0,C_1$, and let $Q$ be the
  corresponding curve of points.  Then there are morphisms $\pi_i:Q\to C_i$
  such that $Li^*L\rho_0^*\cong L\pi_0^*$ and
  \[
  Li^*L\rho_1^*\cong {\cal L}_s^{-1}\otimes L\pi_1^*.
  \]
  for some line bundle ${\cal L}_s$ having degree $1$ on every vertical
  component.
\end{lem}

\begin{proof}
  Let $\bar\pi_i:\bar{Q}\to C_i$ be the natural degree 2 morphisms, let
  $\pi:Q\to\bar{Q}$ be the morphism contracting the vertical fibers, and
  define $\pi_i=\pi\circ \bar\pi_i$.  Then for any line bundle ${\cal L}$
  on $C_i$,
  \[
  Li^*\rho_i^*{\cal L}\cong Li^*\rho_i^*\sO_{C_i}\otimes \pi_i^*
  {\cal L},
  \]
  since twisting by ${\cal L}$ gives another quasi-ruled surface.  The
  claim for $i=0$ follows immediately, and for $i=1$ reduces to showing
  that $i^*\rho_1^*\sO_{C_i}$ is an invertible sheaf having degree $-1$
  on every vertical component.  Writing this object as ${\cal L}_s^{-1}$,
  we find
  \[
  R\Hom_X(\rho_1^*M,\rho_0^*N)
  \cong
  R\Hom_Q(\pi_1^*M,{\cal L}_s\otimes \pi_0^*N)
  \cong
  R\Hom_{C_1}(M,N\otimes_{C_0} R(\pi_0\times\pi_1)_*{\cal L}_s),
  \]
  so that $R(\pi_0\times\pi_1)_*{\cal L}_s$ is the sheaf bimodule from
  which we constructed the quasi-ruled surface.  Since $(\pi_0\times
  \pi_1)_*$ has at most $1$-dimensional fibers over its image, it has
  cohomological dimension $\le 1$, and thus the standard spectral sequence
  implies that ${\cal L}_s$ must itself be a sheaf.  We similarly find that
  the adjoint sheaf bimodule is given by
  \[
  R(\pi_0\times\pi_1)_*({\cal L}_s^{-1}\otimes \pi_0^!\sO_{C_0})
  \]
  so that ${\cal L}_s^{-1}$ is also a sheaf.  Thus ${\cal L}_s$ is
  reflexive, and since it has rank 1 away from the vertical fibers, must
  have rank 1.  If ${\cal L}_s$ had degree $\ge 2$ on any vertical fiber, it
  would have a subsheaf isomorphic to the structure sheaf of that fiber,
  and thus its direct image would fail to be pure $1$-dimensional.  On the
  other hand, ${\cal L}_s^{-1}$ must have negative degree on every vertical
  fiber, since the pullback of a point sheaf agrees with that point sheaf
  in cohomological degree $\ge 0$, and thus the corresponding map
  $i^*\pi_1^*\sO_{p_1}\to i^*\pi_0^*\sO_{p_0}$ must be injective
  with cokernel $\sO_p$.
\end{proof}

\begin{prop}\label{prop:semiorth_of_quasi-ruled}
  Let $X$ be an iterated blowup of a noncommutative quasi-ruled surface
  over $(C_0,C_1)$, with curve of points $Q$, and let $L\rho_0^*:D^b_{\coh}
  C_0\to D^b_{\coh} X$ be the corresponding embedding.  Then there is a
  commutative rationally ruled surface $\rho:X'\to C_1$ such that there is
  an equivalence
  \[
  \kappa:(\rho_0^*D^b_{\coh} C_0)^\perp\cong (\rho^*D^b_{\coh} C_1)^\perp
  \]
  satisfying $\kappa(M)|_Q\cong i^* M$, and thus
  \[
  R\Hom(M,\rho_0^*N)\cong R\Hom_Q((\kappa M)|_Q,\pi_0^*N).
  \]
\end{prop}

\begin{proof}
  If $X$ is a noncommutative quasi-ruled surface, then we can take
  $X'$ to be $\P(\pi_{1*}{\cal L}_s)$, so that $Q$ embeds in $X'$
  in such a way that $\sO_{X'}(s)|_Q\cong {\cal L}_s$.  Blowups then work
  as in the rational case.
\end{proof}

\begin{rem}
  Note that when $X$ is rational, the commutative surface obtained in this
  way will not in general be isomorphic to that associated to
  $\sO_X^\perp$, although this just involves a minor twist: the above
  equivalence extends to an isomorphism $\rho_0^*\sO_{\P^1}(1)^\perp\cong
  \rho^*\sO_{\P^1}(1)^\perp$ that still respects the restriction to $Q$,
  and there is an abelian equivalence taking $\rho_0^*\sO_{\P^1}(1)$ to the
  structure sheaf of a slightly different noncommutative rational surface.
\end{rem}

Again, this gives derived equivalences corresponding to any isomorphism of
commutative ruled surfaces, giving the desired equivalences corresponding to
elementary transformations.  Since this respects the projection
$\rho_{0*}\alpha_*$ to $\im\rho_0^*$, it will be an abelian equivalence
whenever that projection is pseudo-canonical.

There is a similar statement associated to iterated blowups of a
noncommutative surface: if $\pi:X\to Y$ is such a blowup such that the
curve of points $Q\subset Y$ is a surface curve, then there is a
commutative blowup $\pi':X'\to Y'$ with an isomorphism $(\pi^*D^b_{\coh}
Y)^\perp \cong (\pi^*D^b_{\coh} Y')^\perp$ respecting restriction to $Q$.
Indeed, we need simply embed $Q$ in a smooth projective surface $Y'$ and
then perform the corresponding sequence of blowups.  In particular, the
results in Section \ref{sec:invisible_comm} on the category of
sheaves/objects annihilated by $R\pi_*$ carry over immediately to the
noncommutative case; we have seen that the derived categories agree, as do
the inductive descriptions of the $t$-structures.  (In fact, it suffices
for $Q$ to be Gorenstein, as it still locally embeds in a surface, and the
categories only depend on how $Q$ behaves near the points being blown up.)

One consequence is that any noncommutative rational surface over an
algebraically closed field is derived equivalent to a ring.  Indeed, this
is clearly true for commutative planes and for Hirzebruch surfaces, as in
either case we have a full exceptional collection of line bundles with no
higher morphisms (either $(\sO_{\P^2}(-2),\sO_{\P^2}(-1),\sO_{\P^2})$ or
$(\sO_X(-s-(d+1)f),\sO_X(-s-df),\sO_X(-f),\sO_X)$ for $d\gg 0$), and thus
the endomorphism dg-algebra of the corresponding direct sum of line bundles
is an ordinary ring.  More generally, writing a rational surface as an
iterated blowup of a Hirzebruch surface gives a semiorthogonal
decomposition in which the other component is the derived category of
invisible sheaves, so generated by $\bigoplus_c P_c$ where $P_c$ are the
projective objects in that abelian category.  We then find that
\[
\bigoplus_c P_c[-1]\oplus \pi^*\sO_X(-s-(d+1)f)\oplus
\pi^*\sO_X(-s-df)\oplus \pi^*\sO_X(-f)\oplus \pi^*\sO_X
\]
generates the derived category and again has endomorphism dga an ordinary
ring.  More generally, for a noncommutative rational surface, we still have
a generator of $\sO_X^\perp$ of the form
\[
\bigoplus_c P_c[-1]\oplus \pi^*\sO_X(-s-(d+1)f)\oplus
\pi^*\sO_X(-s-df)\oplus \pi^*\sO_X(-f),
\]
and need merely verify that the $R\Hom$ from this to $\sO_X$ is supported
in degree 0.

Note that higher genus (quasi-)ruled surfaces are {\em not} derived
equivalent to rings.  Indeed, since Hochschild homology is additive for
semiorthogonal decompositions, the Hochschild homology of a higher genus
(quasi-)ruled surface has as a summand the Hochschild homology of $C_0$,
and thus ${\mathbb{H}\mathbb{H}}_{-1}(X)$ is nonzero, while the Hochschild
homology of a ring is nonnegative.  Also, the above construction can fail
for rational surfaces over a non-closed field, in particular for surfaces
without a rational ruling (such as the generic del Pezzo surface of degree
1).

For ruled surfaces, the two maps $\pi_i:Q\to C$ are in the same component
of the moduli space of such maps, and thus one would like to understand
that moduli space.  Of course, composing with an automorphism of $C$ has no
effect on the dg-category, so one should really consider the corresponding
quotient.

\begin{lem}\label{lem:Pic0Q/Pic0C}
  Let $C/S$ be a smooth projecctive curve of positive genus, and $(X,Q)$ a
  commutative anticanonical rationally ruled surface over $C$.  Then the
  component of $\Map(Q,C)/\Aut(C)$ containing the induced map $\pi:Q\to C$
  is naturally isomorphic to $\Pic^0(Q)/\Pic^0(C)$.
\end{lem}

\begin{proof}
  The map $\pi$ has a Stein factorization
  \[
  Q\to \Spec_C(\pi_*\sO_Q)\to C
  \]
  and thus we may deform $\pi$ by deforming the map $\bar\pi$ from
  $\bar{Q}:=\Spec_C(\pi_*\sO_Q)$ to $C$.  From the (canonical up to
  scalars) exact sequence
  \[
  0\to \omega_X\to \sO_X\to \sO_Q\to 0,
  \]
  one obtains an exact sequence
  \[
  0\to \sO_C\to \pi_*\sO_Q\to \omega_C\to 0
  \]
  with $R^1\pi_*\sO_Q=0$.  In particular, $\bar{Q}$ is a double cover of
  $C$.  Moreover, since $R\pi_*\sO_Q\cong R\bar\pi_*\sO_{\bar{Q}}$, the
  pull-back map $\Pic(\bar{Q})\to \Pic(Q)$ induces an isomorphism on tangent
  spaces, and is thus easily seen to be injective with discrete cokernel.
  The map $Q\to\bar{Q}$ is the image under the map to weighted projective
  space corresponding to the pullback of any degree 1 bundle on $C$, and it
  it easy to see via the case analysis below that any algebraically
  equivalent bundle on $\bar{Q}$ arises as the pullback of a bundle on $C$
  through a deformation of $\bar\pi$.  It follows that $Q\to \bar{Q}$ is
  independent of $\pi$ and thus $\Map(Q,C)\cong \Map(\bar{Q},C)$.

  For $g\ge 2$, there is a canonical section $C\to \bar{Q}$ (the minimal
  section on $X$ induces a section $C\to Q$ and thus a section $C\to
  \bar{Q}$), and thus we may identify $\Map(\bar{Q},C)/\Aut(C)$ with the
  family of maps that respect the section.  The algebra $\pi_*\sO_Q$ is the
  square-zero extension of $\sO_C$ by $\omega_C$ (the section $C\to
  \bar{Q}$ splits the short exact sequence giving $\pi_*\sO_Q$), and thus
  maps $C\to\bar{Q}$ respecting the section are in natural correspondence
  with derivations taking values in $\omega_C$, so are parametrized by
  $\G_a\cong \Pic^0(Q)/\Pic^0(C)$.  More explicitly, given $\hbar\in \G_a$,
  the corresponding algebra map $\Gamma(U;\sO_C)\to \Gamma(U;\pi_*\sO_Q)$
  is given by
  \[
  f\mapsto f + \hbar df \epsilon.
  \]
  
  When $g(C)=1$, there is no longer a canonical section, but $X$ either has
  a section of negative degree, inducing a section $C\to Q$, or is a ruled
  surface with no such section, forcing $X\cong E\times \P^1$, and in
  either case $\pi_*\sO_Q\cong \sO_C\oplus \omega_C\cong \sO_C^2$.  Thus
  the sections of $\bar{Q}\to C$ are classified by a double cover $T/S$,
  and we may identify $\bar{Q}\cong C\times T$ in such a way that $\pi$ is
  the projection to $C$.  It follows that $\Map(\bar{Q},C)\cong
  \Map(T,\Map(C,C))$, with $\pi$ contained in the component
  $\Map(T,\Aut^0(C))\cong \Map(T,\Pic^0(C))\cong \Pic^0(Q)$.  Moreover, the
  action of $\Aut^0(C)\cong \Pic^0(C)$ on this component is nothing other
  than the usual group law, and thus again
  $\Map_\pi(T,\Pic^0(C))/\Aut^0(C)\cong \Pic^0(Q)/\Pic^0(C)$.
\end{proof}

\begin{rem}
  Note that in both cases, the deformation can be obtained by composing
  with an automorphism of $\bar{Q}$; this is obvious in the genus 1 case,
  while for $g\ge 2$ one uses the automorphism $f+\epsilon g\mapsto
  f+\epsilon (g+\hbar df)$.
\end{rem}

\begin{rem}
  It is straightforward from the explicit description of the deformation to
  see that any such deformation of $\pi_0$ produces a sheaf bimodule and
  thus (inducting along blowups) a noncommutative rationally ruled surface.
\end{rem}

\begin{rem}
  The hypothesis that $g(C)>0$ is not quite necessary; it is no longer true
  that $\Map(Q,C)\cong \Map(\bar{Q},C)$, but it {\em is} true by the same
  argument that $\Map(\bar{Q},C)\to \Map(Q,C)$ is a closed and open
  embedding.
\end{rem}

There is a more uniform description of a map from the component
$\Map_\pi(Q,C)$ to $\Pic^0(Q)/\Pic^0(C)$.  As discussed in the proof, given
a point $L\in \Pic^1(C)$ and a deformation $\phi$ in the given component
$\Map_\pi(Q,C)$, there is an induced element of $\Pic^0(Q)$ given by
$\phi^*L\otimes \pi^*L^{-1}$.  Changing $L$ changes the class in
$\Pic^0(Q)$ by an element of $\Pic^0(C)$ (both $\pi^*:\Pic^0(C)\to
\Pic^0(Q)$ and $\phi^*:\Pic^0(C)\to \Pic^0(Q)$ having the same image), and
thus we get a well-defined morphism from $\Map_\pi(Q,C)\to
\Pic^0(Q)/\Pic^0(C)$ that contracts $\Aut(C)$-orbits.

For $g\ge 2$, if we start with a line bundle with transition function $u$,
then the the corresponding point in $\Pic^0(Q)$ is the pullback of the line
bundle on $\bar{Q}$ with transition function $1+(\hbar du/u)\epsilon$.
This is nothing other than the usual section $\G_a\to \Pic^0(\bar{Q})$, so
agrees with the earlier identification with $\G_a$.  (In particular, this
establishes our earlier claim that every bundle in the relevant component
of $\Pic(Q)$ arises as a pullback.)  Similarly, for $g=1$, our deformations
arise by translation in the group $\Pic^0(Q)$, and it is easy to see that
the induced map $\Map_\pi(\bar{Q},C)\cong \Pic^0(Q)$ agrees with the map
from the proof.  (Again, it follows that everything in the relevant
component of $\Pic(Q)$ arises as a pullback through a deformation of $\pi$,
completes the proof that $\Map_\pi(Q,C)\cong \Map_\pi(\bar{Q},C)$.)

\begin{rem}
  The map from deformations to $\Pic^0(Q)/\Pic^0(C)$ has a natural
  interpretation in terms of the corresponding deformations of
  dg-categories: it is the image of the class of a point under a natural
  map $K_0^{\text{num}}(X)\to \Pic(Q)/\Pic^0(C)$, so by Proposition
  \ref{prop:order_iff_q_torsion} below is torsion iff $X$ is a maximal
  order and trivial iff $X$ is commutative.
\end{rem}

\chapter{Derived equivalences}
\label{chap:derived_eqs}

\section{Duality}

Of course, we can also use these semiorthogonal decompositions to construct
derived equivalences which are not abelian equivalences.  The simplest
instance of such a derived equivalence is the following duality.  Note that
in the following proposition, only the category $D^b_{\coh} X_q$ is
determined by the data, as we have not chosen a blowdown structure on $X'$.
For quasi-ruled surfaces, such a duality can be obtained from the adjoint
involution discussed above, and the duality we give is, in fact, closely
related to that involution.

\begin{prop}
  Let $(X',Q)$ be an anticanonical rational surface, and for each $q\in
  \Pic^0(Q)$, let $X_q$ denote the corresponding (derived) noncommutative
  rational surface.  Then there is a family of equivalences $\ad_q:(D^b_{\coh}
  X_q)^{\text{op}}\cong D^b_{\coh} X_{q^{-1}}$ such that
  $\ad_{q^{-1}}\ad_q\cong \text{id}$.
\end{prop}

\begin{proof}
  On $X'$ we have the natural (Cohen-Macaulay) duality functor
  $M^D:=R\sHom(M,\omega_{X'})$.  This takes $\sO_{X'}$ to $\omega_{X'}$,
  and if $M\in \sO_{X'}^\perp$, then
  \[
  R\Hom_{X'}(\sO_{X'},M^D)\cong R\Hom_{X'}(M^D,\omega_{X'}[2])^*
  \cong R\Hom(\sO_{X'},M[2])^* = 0,
  \]
  so that ${-}^D$ restricts to a duality on $\sO_{X'}^\perp$.  We claim that
  there is a duality $\ad_q:(D^b_{\coh} X_q)^{\text{op}}\to D^b_{\coh}
  X_{q^{-1}}$ such that $\ad_q\sO_{X_q}\cong \theta \sO_{X_{q^{-1}}}$ and
  \[
  \ad_q \kappa_q^{-1}M = \kappa_{q^{-1}}^{-1} M^D.
  \]
  This transforms the semiorthogonal decomposition as
  $(\sO_{X_q}^\perp,\sO_{X_q})\to (\theta
  \sO_{X_{q^{-1}}},\sO_{X_{q^{-1}}}^\perp)$, so it remains only to show
  that it respects the gluing functor.  In other words, we need to show
  that
  \[
  R\Hom(\ad_q\sO_{X_q},\ad_qM)\cong R\Hom(M,\sO_{X_q})
  \]
  for any $M\in \sO_{X_q}^{\perp}$.  Since
  \[
  R\Hom_{X_{q^{-1}}}(\theta \sO_{X_{q^{-1}}},\kappa_{q^{-1}}^{-1} M^D)
  \cong
  R\Hom_{X_{q^{-1}}}(\cone(\Tr_{q^{-1}})[-2],\kappa_{q^{-1}}^{-1} M^D),
  \]
  where $\Tr_{q^{-1}}:\sO_{X_{q^{-1}}}\to \theta \sO_{X_{q^{-1}}}[2]$
  is the natural trace map associated to the Serre functor, and
  $\cone(\Tr_{q^{-1}})\in \sO_{X_{q^{-1}}}^\perp$, we find that
  \begin{align}
  R\Hom_{X_{q^{-1}}}(\theta \sO_{X_{q^{-1}}},\kappa_{q^{-1}}^{-1} M^D)
  &\cong
  R\Hom_{X'}(\kappa_{q^{-1}}\cone(\Tr_{q^{-1}})[-2],M^D)\notag\\
  &\cong
  R\Hom_{X'}(M,(\kappa_{q^{-1}}\cone(\Tr_{q^{-1}})[-2])^D),
  \end{align}
  so that it remains only to compute
  $\kappa_{q^{-1}}\cone(\Tr_{q^{-1}})[-2]$.
  For any $M\in \sO_{X'}^\perp$, we have
  \[
  R\Hom_{X'}(\kappa_{q^{-1}}^{-1}M,\cone(\Tr_{q^{-1}})[-2])
  \cong
  R\Hom_{X'}(\kappa_{q^{-1}}^{-1}M,\sO_{X_{q^{-1}}}[-1])
  \cong
  R\Hom_{X'}(M,q^{-1}[-1]),
  \]
  so that $\kappa_{q^{-1}}\cone(\Tr_{q^{-1}})[-2]\cong q^{-1}[-1]$
  (adjusted accordingly when $q=\sO_Q$) and dualizing gives $q$ as
  required.

  Since $\kappa_q \ad_{q^{-1}}\ad_q\kappa_q^{-1}={-}^{DD}\cong \text{id}$,
  to show that this is an involution it remains only to check that
  $\ad_{q^{-1}}\theta \sO_{X_{q^{-1}}}\cong \sO_{X_q}$.  But this follows
  from the fact the Serre functor is intrinsic, and thus any contravariant
  equivalence satisfies $S \ad_q S\cong \ad_q$.
\end{proof}

This duality respects the blowup and ruled surface decompositions.

\begin{lem}
  Let $\alpha:\tilde{X}'\to X'$ be a monoidal transformation centered at a
  point of $Q$, and let $\tilde{X}_{q}$, $X_{q}$ be the corresponding
  noncommutative families, with $\tilde{X}_q$ the blowup of $X_q$.  Then
  $\ad_q\sO_e(-1)\cong \sO_e(-1)[-1]$, $\ad_q L\alpha^*M\cong
  L\alpha^!\ad_q M$, and $\ad_q R\alpha_* M\cong \alpha_* \ad_q M$.
\end{lem}

\begin{proof}
  Since $\kappa_q \sO_e(-1)\cong \sO_e(-1)$, the action of the duality on
  $\sO_e(-1)$ reduces to the commutative case.  For the action on
  $L\alpha^*M$, we first note that if $M\in \sO_{X_q}^\perp$, then
  $L\alpha^*M\in \sO_{\tilde{X}_q}^\perp$, and thus
  $\ad_q L\alpha^*M$ can be computed inside $\tilde{X}'$, giving
  \[
  \kappa_{q^{-1}} \ad_q L\alpha^* M
  \cong
  \ad_q L\alpha^* \kappa_q M
  \cong
  L\alpha^! \ad_q\kappa_q M
  \cong
  L\alpha^! \kappa_{q^{-1}} \ad_q M
  \cong
  \kappa_{q^{-1}} L\alpha^! \ad_q M.
  \]
  Since
  \[
  \ad_q L\alpha^*\sO_{X_q}\cong \ad_q \sO_{\tilde{X}_q}\cong \theta
  \sO_{\tilde{X}_{q^{-1}}} \cong L\alpha^!\theta \sO_{X_{q^{-1}}} \cong
  L\alpha^! \ad_q \sO_{X_q},
  \]
  the action on $L\alpha^*M$ is as described.  Finally, we find
  \begin{align}
  R\Hom(N,R\alpha_* \ad_q M)
  &\cong
  R\Hom(L\alpha^* N,\ad_q M)\notag\\
  &\cong
  R\Hom(M,\ad_q L\alpha^* N)\notag\\
  &\cong
  R\Hom(M,L\alpha^! \ad_q N)\notag\\
  &\cong
  R\Hom(R\alpha_* M,\ad_q N)\notag\\
  &\cong
  R\Hom(N,\ad_q R\alpha_* M)
  \end{align}
  so that $R\alpha_* \ad_q M\cong \ad_q R\alpha_* M$ as required.
\end{proof}

\begin{lem}
  Let $(X',Q')$ be an anticanonical Hirzebruch surface with $\rho:X'\to
  \P^1$.  Then on $X_q$ one has $\ad_q \rho_1^*M\cong \rho_1^* M^D$ and
  $\ad_q \rho_0^*M\cong \theta \rho_0^* \theta^{-1} M^D$.
\end{lem}

\begin{proof}
  Since $\im\rho_1^*\subset \sO_{X_q}^\perp$, we can compute the duality
  inside $\im\rho_1^*$ on $X'$, where it becomes
  \[
  \ad_q (\rho^*M(-s))\cong \rho^*(M^D)(-s),
  \]
  easily checked on line bundles.  Similarly, a commutative computation gives
  $\ad_q \rho_0^*\sO_{\P^1}(-1)\cong \theta \rho_0^*\theta^{-1} \sO_{\P^1}(-1)$,
  and the claim for $\sO_{\P^1}$ follows from the value of $\ad_q$ on $\theta
  \sO_X$.
\end{proof}

This suggests that there should be a similar duality on a general iterated
blowup of a quasi-ruled surface.  Although this could in principle be
constructed using the $\rho_0^*D^b_{\coh} C_0$-based semiorthogonal
decomposition, the calculation involves some somewhat tricky adjoint
computations, and thus we use an inductive approach instead.

\begin{prop}
  Let $C_0$, $C_1$ be smooth projective curves, let ${\cal E}$ be a sheaf
  bimodule of birank $(2,2)$ on $C_0\times C_1$, and let $X$ be the
  corresponding noncommutative $\P^1$-bundle.  Similarly, let $X^{\ad}$ be
  the noncommutative $\P^1$-bundle associated to the sheaf bimodule
  $\sExt^1_{C_0\times C_1}({\cal E},\omega_{C_0\times C_1})$.  Then
  $(D^b_{\coh} X)^{\text{op}}\cong D^b_{\coh} X^{\ad}$
\end{prop}

\begin{proof}
  Per the rational case, we want a duality defined in terms of the duality
  on $C_0$, $C_1$ by
  \[
  \ad \rho_1^*M\cong \rho_1^*\sHom(M,\omega_{C_1})\qquad\text{and}\qquad
  \ad \theta \rho_0^*M\cong \rho_0^*\sHom(M,\sO_{C_0}).
  \]
  Since $(\im\theta \rho_0^*,\im\rho_1^*)$ gives a semiorthogonal
  decomposition, this certainly defines a contravariant equivalence to {\em
    some} glued triangulated category $X^{\ad}$, so we just need to show that
  the gluing map corresponds to the given sheaf bimodule.  We compute
  \begin{align}
  R\Hom_{X^{\ad}}(\rho_1^* M,\rho_0^*N)
  &\cong
  R\Hom_{X^{\ad}}(\ad \rho_1^*\sHom(M,\omega_{C_1}),\ad \theta
  \rho_0^*\sHom(N,\sO_{C_0}))\notag\\
  &\cong
  R\Hom_X(\theta \rho_0^*\sHom(N,\sO_{C_0}),\rho_1^*\sHom(M,\omega_{C_1})).
  \end{align}
  Using the functorial distinguished triangle
  \[
  \theta \rho_0^*\theta^{-1}\to \rho_1^* ({-}\otimes_{\sO_{C_0}}{\cal
    E})\to \rho_0^*\to,
  \]
  we find that
  \begin{align}
  R\Hom_X(\theta \rho_0^*\sHom(N,\sO_{C_0}),{}&\rho_1^*\sHom(M,\omega_{C_1}))\notag\\
  &\cong
  R\Hom_X(\rho_1^* (\sHom(N,\sO_{C_0})\otimes_{\sO_{C_0}}{\cal E}),
  \rho_1^* \sHom(M,\omega_{C_1}))\notag\\
  &\cong
  R\Hom_{C_1}((\sHom(N,\sO_{C_0})\otimes_{\sO_{C_0}}{\cal
    E}),\sHom(M,\omega_{C_1})\notag\\
  &\cong
  R\Hom_{C_1}(M,\sHom(\sHom(N,\sO_{C_0})\otimes_{\sO_{C_0}}{\cal
    E},\omega_{C_1})).
  \end{align}
  Using the definition of the tensor product of sheaf bimodules, the
  functor being applied to $N$ has the form
  \begin{align}
  N
  \mapsto
  \sHom(\pi_{1*}(\pi_0^* \sHom(N,\sO_{C_0})\otimes {\cal E}),\omega_{C_1})
  &\cong
  \pi_{1*}\sExt^1(\pi_0^* \sHom(N,\sO_{C_0})\otimes {\cal E},\omega_{C_0\times C_1})\notag\\
  &\cong
  \pi_{1*}(\pi_0^* N\otimes \sExt^1({\cal E},\omega_{C_0\times C_1}))\notag\\
  &=
  N\otimes_{\sO_{C_0}} \sExt^1({\cal E},\omega_{C_0\times C_1})
  \end{align}
  as required.
\end{proof}
  
\begin{rem}
  This duality is of course related to the adjoint involution defined
  above in terms of the sheaf $\Z$-algebra: any module over the
  $\Z$-algebra has a dual which is a module over the opposite algebra, and
  an appropriate degree shift and twist by a line bundle recovers the
  duality.  In particular, when applied to a $1$-dimensional sheaf
  corresponding to a difference or differential equation, this is just a
  (cohomological) shift of the usual duality (i.e., coming from the dual
  representation of $\GL_n$ or $\mathfrak{gl}_n$ as appropriate).  (This
  agrees with what we conjectured in Section \ref{sec:elldiff} based on
  commutative considerations.)
\end{rem}

This duality also respects the curve of points, in the following way.

\begin{prop}
  Let $X$ be a quasi-ruled surface with curve of points $Q$, and let
  $X^{\ad}$ be its dual.  Then the curve of points of $X^{\ad}$ can be
  identified with $Q$ in such a way that $\ad[1]$ restricts to the
  Cohen-Macaulay duality on $D^b_{\coh} Q$.
\end{prop}

\begin{proof}
  Any $S$-point $p$ of $X$ arises from a short exact sequence
  \[
  0\to \rho_1^*\sO_{p_1}\to \rho_0^*\sO_{p_0}\to \sO_p\to 0.
  \]
  Applying $\ad$ and computing the relevant commutative duals gives
  a distinguished triangle
  \[
  \ad \sO_p\to \theta \rho_0^*\sO_{p_0}[-1]\to \rho_1^*\sO_{p_1}[-1]\to.
  \]
  It follows from this that $(\ad \sO_p)[2]$ is an $S$-point sheaf on
  $X^d$, and thus that the corresponding functors (so moduli
  schemes) are isomorphic.  Moreover, one has
  \[
  \ad \sO_p\cong \sO_p[-2]\cong \sHom_Q(\sO_p,\omega_Q)[-1]
  \]
  so that $\ad[1]$ agrees with the Cohen-Macaulay duality on point sheaves.
  It then suffices to verify that it behaves correctly on the structure
  sheaf.  We have the short exact sequence
  \[
  0\to \sO_X(-Q)\to \sO_X\to \sO_Q\to 0
  \]
  and the adjoint gives
  \[
  0\to \theta \sO_{X^{\ad}}\to \theta \sO_{X^{\ad}}(Q)\to \ad\sO_Q[1]\to 0,
  \]
  But $\theta \sO_X(Q)|^\dL_Q\cong \omega_Q$ by Proposition
  \ref{prop:Serre_for_quasiruled}.
\end{proof}

This claim about the restriction is crucial to let us make blowups work.
  
\begin{lem}
  Let $X$ be a noncommutative surface with curve of points $Q$ such that
  $X$ is generated by weak line bundles along $Q$, and let $\tilde{X}$ be
  the blowup of $X$ at $p$.  If $X^{\ad}$ is another noncommutative surface
  with the same curve of points and there is a duality $\ad:(D^b_{\coh}
  X)^{\text{op}}\to D^b_{\coh} X^{\ad}$ such that $\ad(M(-Q))\cong (\ad M)(Q)$
  and $\ad[1]$ restricts to the Cohen-Macaulay duality on $Q$, then there
  is a duality $\widetilde\ad:(D^b_{\coh} \tilde{X})^{\text{op}}\to D^b_{\coh}
  \tilde{X}^{\ad}$, where $\tilde{X}^{\ad}$ is the blowup of $X^{\ad}$ at $p$,
  and $\widetilde\ad[1]$ satisfies the same conditions with respect to the new
  curve of points.
\end{lem}

\begin{proof}
  Following the rational case, we aim to define $\widetilde\ad$ so that
  \[
  \widetilde\ad \sO_e(-1)\cong \sO_e(-1)[-1]\qquad\text{and}\qquad
  \widetilde\ad L\alpha^! M\cong L\alpha^* \ad M.
  \]
  Consistency of the gluing reduces to computing $\ad$ on the appropriate
  point sheaf, which follows by a computation in $Q$.

  It remains only to show that $\widetilde\ad[1]$ interacts correctly with
  the new curve of points $Q^+$.  That $\widetilde\ad(M(-1))\cong
  (\widetilde\ad M)(1)$ follows from the description of $M\mapsto \theta
  M(1)$ in terms of $\theta M(-Q)$ and the fact that $\widetilde\ad$ and
  $\ad$ both invert $\theta$ while $\ad$ inverts ${-}(-Q)$.  Any line
  bundle on $Q$ pulls back via $\alpha^*$ to a line bundle on $Q^+$, and
  one has
  \[
  (\ad \alpha^* {\cal L})[1]
  \cong
  \alpha^! (\ad {\cal L})[1]
  \cong
  \alpha^! \sHom({\cal L},\omega_Q)
  \cong
  \alpha^* \sHom({\cal L},\omega_Q) \otimes \sO_{Q^+}(e).
  \]
  This depends only on the commutative morphism $\alpha:Q^+\to Q$, and may
  thus be computed inside a commutative surface, where it reduces to the
  corresponding fact for the Cohen-Macaulay dual.

  It follows more generally that for any $d\in \Z$, one has
  \[
  \ad (\alpha^*{\cal L}(d))[1]
  \cong
  \sHom(\alpha^*{\cal L},\omega_{Q^+})(-d)
  \cong
  \sHom(\alpha^*{\cal L}(d),\omega_{Q^+})
  \]
  But the line bundles on $Q^+$ of this form are ample, and thus any object
  $M$ in $D^b_{\coh} Q^+$ is quasi-isomorphic to a right-bounded complex in
  which the terms are sums of such line bundles.  This complex is acyclic
  for $\ad[1]$, so $\ad M[1]$ is quasi-isomorphic to the left-bounded
  complex obtained by applying $\ad[1]$ term-by-term.  It follows that $\ad
  M[1]\cong \sHom_{Q^+}(M,\omega_{Q^+})$ as required.
\end{proof}

This gives rise to an inductive construction of a duality $\ad$ on any
iterated blowup of a quasi-ruled surface or noncommutative plane, which
agrees with the previous definition in the rational surface case.

\section{Line bundles}

Although this duality cannot be expected to be exact (any more than the
Cohen-Macaulay duality it deforms), we can in fact control fairly well how
it interacts with the $t$-structure.  The key idea is that there is a
particularly nice class of generators of $\qcoh X$ that is preserved by
duality.

\begin{defn}\label{defn:line_bundle}
  A {\em line bundle} on a rational or rationally quasi-ruled surface is a
  member of the smallest class of sheaves on such surfaces such that (a)
  any sheaf isomorphic to a line bundle is a line bundle, (b) on a
  noncommutative plane, $\sO_X$, $\sO_X(-1)$ and $\sO_X(-2)$ are line
  bundles, (c) on a quasi-ruled surface, $\rho_d^*L$ is a line bundle for
  any line bundle $L\in C_d$, and (d) the image of a line bundle under
  $\theta$ or $\alpha_i^*$ is a line bundle.
\end{defn}

\begin{rem}
  This differs from other notions of line bundles in the literature; in
  addition to the weak notion considered above (the restriction to $Q$ is a
  line bundle), one might also consider those weak line bundles which are
  reflexive under $h^0\ad$ (this was considered in
  \cite{NevinsTA/StaffordJT:2007} for noncommutative planes).  The above
  notion has the merit that the corresponding family of moduli spaces is
  flat; in particular, in a sufficiently large family of rationally ruled
  surfaces, the line bundles are deformations of line bundles on the
  commutative fibers of the family.  Although it in principle depends on
  the way in which we described the surface as an iterated blowup, this is
  not the case, see Corollary \ref{cor:lbs_are_intrinsic} below.  It is
  unclear whether there is a corresponding notion for vector bundles (i.e.,
  as a subclass of reflexive coherent sheaves).
\end{rem}

One key fact about line bundles is that there is a notion of twisting by
line bundles.  To be precise, for any line bundle $L$ on $X$, there is
another rational (resp. rationally quasi-ruled) surface $X'$ such that
there is a Morita equivalence $\coh X\cong \coh X'$ taking $L$ to
$\sO_{X'}$, taking point sheaves to point sheaves, and taking line
bundles to line bundles.  This is trivial for noncommutative planes, is
straightforward for quasi-ruled surfaces, and follows by an easy
induction in general.  (For a blowup, we simply take the composition of a
power of $\theta$ with the result of applying such an equivalence to the
original surface, making sure to keep track of the point being blown up.)
When $X$ is commutative, this is agrees with the usual notion of twisting
by a line bundle, but in general this twist is not an autoequivalence.

Line bundles give an alternate description of the $t$-structure.

\begin{lem}
  An object $M\in D_{\qcoh}(X)$ is in $D_{\qcoh}(X)^{\ge 0}$ iff
  $\Hom(L[p],M)=0$ whenever $p$ is a negative integer and $L$ is a line
  bundle.
\end{lem}

\begin{proof}
  This is automatic for a plane, while on a quasi-ruled surface, this
  follows via adjunction from the corresponding fact for line bundles on
  commutative curves; for rationally quasi-ruled surfaces, it follows by an
  easy induction.
\end{proof}

We also have the following important fact, which in particular shows that
line bundles generate $\coh X$, so $\qcoh X$.

\begin{prop}
  If $M\in \coh(X)$, then there is a line bundle $L$ such that
  $\Ext^p(L,M)=0$ for $p>0$ and $L\otimes \Hom(L,M)\to M$ is surjective.
\end{prop}

\begin{proof}
  This is certainly true by construction for planes and quasi-ruled
  surfaces and for blowups follows by an easy induction from the fact that
  for $l\gg 0$, $\alpha_{m*}\theta^{-l}M$ is a sheaf and
  $\alpha_m^*\alpha_{m*}\theta^{-l}M\to M$ is globally generated.  Indeed,
  if $L'$ is the line bundle satisfying the conclusion for the sheaf
  $\alpha_{m*}\theta^{-l}M$ on $X_{m-1}$, then we may take $L=\theta^l
  \alpha_m^* L'$.
\end{proof}

\begin{rems}
  We will prove a number of refinements of this fact below.  One easy
  refinement is that we need not use every line bundle here; in particular,
  if we fix points $x_0\in C_0$, $x_1\in C_1$, we need only use those
  bundles coming from $\sO_{C_d}(nx_d)$ for $d,n\in \Z$.  (This refinement
  also applies to the description of the $t$-structure.)
\end{rems}
 
\begin{rems}
  More generally, given a finite collection of coherent sheaves, we can
  find a single line bundle that acyclically globally generates all of
  them: simply apply the Proposition to the direct sum.
\end{rems}

Since $\ad \theta\cong \theta^{-1}\ad$ and $\ad \alpha^*\cong \alpha^!\ad$,
it is easy to verify that for any line bundle $L$, $\ad L$ is also a line
bundle.  We may thus view the above derived duality as the derived functor
of a contravariant functor $\ad$ of abelian categories, and will thus in
the sequel denote it as $R\ad$, particularly when applying it to sheaves.

\begin{cor}
  The duality $\ad$ is left exact of homological dimension $\le 2$.
\end{cor}

\begin{proof}
  Since line bundles generate $\coh X$ and $\ad$ takes line bundles to line
  bundles, it is left exact.  Now, for any line bundle $L$ and coherent
  sheaf $M$, we have
  \[
  R\Hom(L,R\ad M)\cong R\Hom(M,\ad L)\cong R\Hom(\ad L,\theta M[2])^*,
  \]
  and thus $h^p R\Hom(L,R\ad M)=0$ for $p\notin \{0,1,2\}$.  For each $p'$,
  we may choose $L$ so that the first $p'+1$ cohomology sheaves of $R\ad M$
  are acyclically globally generated by $L$.  But then the spectral
  sequence gives $h^i R\Hom(L,R\ad M) = \Hom(L,R^i\ad M)$ for $0\le i<p'$,
  and thus $R^i\ad M=0$ for $2\le i<p'$.  Since this is true for all $p'$,
  we conclude that $\ad$ has homological dimension 2 as required.
\end{proof}

\section{Derived equivalences of deformed elliptic surfaces}

The other family of true derived equivalences we consider are related to
derived autoequivalences of commutative elliptic surfaces.  If $C$ is a
smooth curve and $X/C$ is the minimal proper regular model of an elliptic
curve over $k(C)$, then there is an associated action of $\SL_2(\Z)$ on
$D^b_{\coh} X$.  If $X$ has only smooth fibers, this is just a relative
version of the action of $\SL_2(\Z)$ on the derived category of an elliptic
curve, and in general a similar construction works (albeit with some issues
of nonuniqueness when $X$ has reducible fibers).

Since the deformation theory of an abelian category depends only on the
associated derived category, the autoequivalences of the derived category
act on the space of infinitesimal deformations.  It is therefore not
unreasonable to expect this to extend to an action on algebraic
deformations, such that two deformations in the same orbit should be
derived equivalent.  This is of course purely heuristic, but works for both
types of elliptic surface that have a rational ruling, and thus admit
noncommutative deformations via the above construction.

A rationally ruled surface over a curve of genus $>1$ does not contain any
curves of genus 1 (a vertical curve has genus 0, while a horizontal curve
maps to the base of the ruling, so has genus $>1$).  If the base curve has
genus 1, then every fiber of the elliptic surface contains a smooth genus 1
curve isogenous to the base, and thus the minimal proper regular model must
be a constant family, so that the only rationally ruled case is $E\times
\P^1$.  In the genus 0 case, we have a rational elliptic surface, and it is
well-known that such a surface, if relatively minimal, must have $K^2=0$,
with the elliptic fibration given by the anticanonical pencil $|{-}K|$.
These cases then deform to (a) noncommutative ruled surfaces over $E$ in
the same component of the moduli stack as $E\times \P^1$, and (b)
noncommutative rational surfaces with $\theta\sO_Q\in \Pic^0(Q)$.

For the genus 1 case (which is simpler, since the fibers are all smooth),
we recall that a sheaf bimodule on $E\times E$ corresponding to a
noncommutative ruled surface necessarily has reducible (or nonreduced)
support, and thus may be viewed as an extension of two rank $(1,1)$ sheaf
bimodules (i.e., line bundles on bidegree $(1,1)$ curves in $E\times E$),
both algebraically equivalent to the diagonal.  The moduli space of such
rank $(1,1)$ sheaf bimodules can itself be identified with $E\times
\Pic(E)$: the bimodule is given by a line bundle in $\Pic^0(E)\cong E$
supported on the graph of a translation.  So as to obtain a surface
sufficiently similar to $E\times \P^1$, we insist that the line bundle be
degree 0, so that the moduli space is $E\times E$.  This moduli space has a
universal family, which in turn induces a Fourier-Mukai functor $D^b_{\coh}
(E\times E)\to D^b_{\coh} (E\times E)$ taking point sheaves to their
corresponding sheaf bimodules.  This is actually an equivalence, allowing
us to identify the relevant component of the moduli stack of rank $(2,2)$
sheaf bimodules with the moduli stack $(E\times E)^{\langle 2\rangle}$ of
0-dimensional sheaves on $E\times E$ of Euler characteristic 2.  (This
stack is a mild thickening of the 2-point Hilbert scheme of $E\times E$,
gluing in sheaves of the form $\sO_p^2$ as a divisor of stacky points.)

We thus obtain a family of noncommutative ruled surfaces over this stack.
Any automorphism of the scheme $E\times E$ has an induced action on the
stack, and we see that the translation subgroup does not actually change
the noncommutative ruled surface.  (The effect on the sheaf bimodule is to
modify the identification between the two copies of $E$ by a translation,
and then twist by the pullback of a line bundle on $E$.)  Thus we should
really view this as a family of noncommutative ruled surfaces over the
quotient $(E\times E)^{\langle 2\rangle}/E\times E$.  There still remains a
group $\Aut(E\times E)$, and the subgroup $U(E\times E)$ (automorphisms
which are unitary with respect to the involution induced by the Poincar\'e
bundle), which is generically $\SL_2(\Z)$, acts as derived equivalences on
the fiber $E\times \P^1$ over the point corresponding to $\sO_p^2$.

\begin{prop}
  The action of $U(E\times E)$ on $(E\times E)^{\langle
    2\rangle}/E\times E$ lifts to an action by equivalences on the derived
  categories of the corresponding noncommutative ruled surfaces.
\end{prop}

\begin{proof}
  By \cite{PolishchukA:1996}, there is an action of $U(E\times E)$ on
  $D^b_{\coh} E$, and thus on both subcategories in the semiorthogonal
  decomposition.  One thus obtains a derived equivalence to a glued
  category in which the gluing functor has been conjugated by the given
  element of $\Aut(D^b_{\coh} E)$.  This reduces to computing the action on
  the sheaf bimodules corresponding to $E\times E$, which again follows by
  \cite{PolishchukA:1996}.
\end{proof}

\begin{rem}
  More generally, if $A$ is an abelian variety, the same argument gives a
  family of rank $(2,2)$ sheaf bimodules on $A\times A$ parametrized by the
  stack analogue $(A\times A^\vee)^{\langle 2\rangle}$ of the $2$-point
  Hilbert scheme, and the action of $U(A\times A^\vee)$ on $D^b_{\coh} A$
  induces derived equivalences between the corresponding noncommutative
  $\P^1$-bundles, deforming the action on $A\times \P^1$.
\end{rem}

\medskip

In the rational case, the elliptic surfaces we consider are obtained by
blowing up the base points of a pencil of cubic curves on $\P^2$.  (We say
``base points'' rather than ``base locus'', as the latter will in general
be singular; to obtain an elliptic surface, we should repeatedly blow up a
base point of the pencil until the resulting anticanonical pencil is base
point free.)  If we instead blow up only 8 of the base points, the result
is a (possibly degenerate) del Pezzo surface $Y$ of degree 1, and the
elliptic surface is obtained by blowing up the unique base point of the
anticanonical pencil on $Y$.  If we fix an anticanonical curve $Q$ on $Y$,
then $Q$ blows up to a fiber of the elliptic surface.  Since $Q$ comes
equipped with a choice of smooth point (the base point of the anticanonical
pencil), we can identify the corresponding component of the smooth locus of
$Q$ with $\Pic^0(Q)$.  This, then, lets us construct a family of
deformations parametrized by $\Pic^0(Q)^2$: we first blow up a point $z\in
\Pic^0(Q)$ viewed as the ``identity'' component of the smooth locus of $Q$,
and then take the noncommutative deformation corresponding to $q\in
\Pic^0(Q)$.  (These are thus nothing other than noncommutative deformations
of the surfaces arising in Sakai's construction \cite{SakaiH:2001} of
discrete Painlev\'e equations.)  As a family of noncommutative surfaces,
this depends on a choice of blowdown structure on $Y$, but the
corresponding derived categories do not, as we have a semiorthogonal
decomposition
\[
(\langle \sO_e(-1)\rangle,\sO_Y^\perp,\langle \sO_X\rangle)
\]
in which the gluing data is determined by $z$ and $q$, with
$\sO_e(-1)|^\dL_Q\cong \sO_z$.  

If we consider the family of all such surfaces (i.e., a $2$-dimensional
fiber bundle over the stack of degenerate del Pezzo surfaces of degree 1
with blowdown structure), we expect that any derived equivalence of the
generic elliptic fiber should extend to the full family.  The relevant
group has the form $\Z^2.\Lambda_{E_8}^2\rtimes \SL_2(\Z)$, where the
center $\Z^2$ is generated by the shift and the Serre functor, one copy of
$\Lambda_{E_8}$ is translation by the Mordell-Weil group, and the other is
the quotient by the canonical bundle of the group of line bundles on $X$
that have degree 0 on $Q$.  The action of $\SL_2(\Z)$ intertwines these two
actions, and one of its generators is the twist by $\sO_X(e)$ where $e$ is
the built-in section (the exceptional curve of the final blowup).  We have
already seen that twists by line bundles extend to the full family (with a
nontrivial $q$-dependent action on the base of the family), and thus need
only one more generator to generate the full group.  One should note here
that the ``correct'' group is really $\Z^2.\Lambda_{E_8}^2\rtimes
W(E_8)\times \SL_2(\Z)$, since we should include the action of $W(E_8)$ by
changing the blowdown structure of the del Pezzo surface.  Also, though
every element of this group extends, the group itself does not; even for
$W(E_8)$, the extension to fibers with $-2$-curves involve a spherical
twist, and there are choices involved which cannot be made in a uniformly
consistent fashion.

\medskip

To get the other generator, we will need to consider how modifications of the
semiorthogonal decomposition interact with the restriction to $Q$.
Given objects $M$ and $N$, let $T_M(N)$ be the cone
\[
M\otimes R\Hom(M,N)\to N\to T_M(N)\to.
\]
If $M$ is exceptional, this is precisely the operation appearing in
modifications of semiorthogonal decompositions involving $\langle M\rangle$
(more precisely, those modifications not changing $M$), while if $M$ is
spherical, this is an autoequivalence (the spherical twist).  (Since we are
really working with a dg-enhanced triangulated category, we can compute the
cone functor $T_M$ inside the dg-category, and thus sidestep the fact that
cones of natural transformations of triangulated functors need not exist.)

\begin{lem}
  Let $X$ be a noncommutative surface with a Serre functor and a curve $C$
  embedded as a divisor, and let $E\in D^b_{\coh} X$ be an exceptional object
  such that $S E \cong E(-C)[d]$ for some $d\in \Z$.  Then for $M\in
  {}^\perp E$, one has
  \[
  T_E(M)|^{\dL}_C\cong T_{E|^{\dL}_C}(M|^{\dL}_C).
  \]
\end{lem}

\begin{proof}
  We recall that $R\Hom(E,M)\cong R\Hom_C(E|^\dL_C,M|^\dL_C)$, and thus
  one can restrict the distinguished triangle defining $T_E(M)|^\dL_C$ to $C$
  to obtain
  \[
  E|^\dL_C\otimes R\Hom_C(E|^\dL_C,M|^\dL_C)\to M|^\dL_C\to T_E(M)|^\dL_C\to,
  \]
  which is precisely the distinguished triangle defining $T_{E|^\dL_C}(M|^\dL_C)$.
\end{proof}

\begin{rem}
  If $X$ is generated by weak line bundles along $C$, then one can compute
  the functor $T_{E|^\dL_C}$ on $D^b_{\coh} C$ by computing it on sheaves
  $M|^\dL_C$.  In particular, the fact that we can undo the modification forces
  $T_{E|^\dL_C}$ to be an autoequivalence, and thus makes $E|^\dL_C$ a spherical
  object.  Something similar holds for more general subcategories in a
  semiorthogonal decomposition: as long as the corresponding
  restriction-to-$C$ functor has adjoints, it will be spherical, and the
  modification acts as the corresponding twist.  This is primarily of
  interest for quasi-ruled surfaces, where the corresponding spherical
  twists are essentially given by the involutions corresponding to double
  covers $Q\to C_i$.
\end{rem}  

Let $X_{z,q}$ be the family defined above associated to some fixed
degenerate del Pezzo surface $Y$.  Note that by construction, any object
in $D^b_{\coh} X_{z,q}$ restricts to a {\em perfect} object in $D^b_{\coh} Q$,
and thus there is a well-defined map $K_0(X_{z,q})\to \Pic(Q)$ given by
$[M]\mapsto \det(M|^{\dL}_Q)$.  This gives us the following identification.

\begin{prop}
  We have $\det(\sO_Q|^{\dL}_Q)\cong z$, $\det(\sO_x|^\dL_Q)\cong q$.
\end{prop}

\begin{proof}
  For $\det(\sO_x|^\dL_Q)$, we first observe that all point sheaves are
  equivalent in $K_0(X_{z,q})$; indeed, a numerically trivial class in
  $K_0(X_{z,q})$ is trivial (since the Mukai pairing is independent of $q$
  and this is true for commutative rational surfaces), and any two point
  sheaves are algebraically equivalent.  We may thus assume that $x$ does
  not lie on any of the exceptional curves, and thus reduce to a
  calculation in $\P^2$ or a Hirzebruch surface, where we can compute it
  via numerical considerations.  Similarly, for $\det(\sO_Q|^\dL_Q)$, we
  can replace $\sO_Q$ by any sheaf $w\in \Pic^0(Q)$, and in particular may
  assume that $w$ has no global sections, so lies in $\sO_X^\perp$.  Since
  $w|^\dL_Q$ has rank 0, we find that $\det(w|^\dL_Q)$ is independent of
  $q$, so equals $z$ by a commutative calculation.
\end{proof}

Now, starting with the semiorthogonal decomposition $(\sO_X^\perp,\sO_X)$
of $D^b_{\coh} X_{z,q}$, let us first modify it to $(\sO_X,{}^\perp\sO_X)$ and
then rotate to obtain the semiorthogonal decomposition
$({}^\perp\sO_X,\theta^{-1}\sO_X)$.  Let $\nu:\sO_{X_{z,0}}^\perp\to
{}^\perp\sO_X$ be the corresponding inclusion functor, and note that
\[
\nu(M)|^\dL_Q\cong T_q^{-1}(M|^\dL_Q)\otimes q^{-1}.
\]
Moreover, we have
\[
\theta^{-1}\sO_X|^\dL_Q\cong
\sO_X(Q)|^\dL_Q
\cong
\det(\sO_Q(Q))
\cong
\det(\sO_Q(Q)|^\dL_Q)
\cong
z,
\]
since $\sO_X(Q)|^\dL_Q$ is a line bundle, $\det(\sO_Q)\cong \sO_Q$, and
$\sO_Q(Q)\in \Pic^0(Q)$.  We have thus changed the gluing data (as functors
to $D^b_{\coh} Q$) from $(M\mapsto M|^\dL_Q,q)$ to $(M\mapsto
T_q^{-1}(M|^\dL_Q),qz)$.  Composing both functors with the autoequivalence
$T_q$ gives the equivalent data $(M\mapsto M|^\dL_Q,T_q(qz))$.  But
$T_q(qz)\cong qz$: this is true if $z\not\cong \sO_X$ since then
$R\Hom(\sO_q,\sO_{qz})=0$, and if $z\cong \sO_X$ since $T_q(q)\cong
q\otimes \Ext^1(q,q)\cong q$.  In other words, the resulting gluing data
has the form $(M\mapsto M|^\dL_Q,qz)$, which is precisely the gluing data
associated to $X_{z,qz}$.  Since $\theta$ is an autoequivalence, we obtain
the following.

\begin{thm}\label{thm:weird_langlands}
  There is a derived equivalence $\Phi_{\sO_X}:D^b_{\coh} X_{z,q}\cong
  D^b_{\coh} X_{z,qz}$ acting on $D^b_{\coh} Q$ as the spherical twist
  $T_q$.  Moreover, this equivalence fixes $\sO_X$ and takes $\sO_e(-1)$ to
  $\theta\sO_X(e)$.
\end{thm}

\begin{proof}
  Everything except for the action on $\sO_e(-1)$ is immediate from the
  above discussion, which also shows that for $M\in \sO_X^\perp$, there is
  a distinguished triangle
  \[
  \theta^{-1}\Phi(M)\to M\to R\Hom(R\Hom(M,\sO_X),\sO_X)\to
  \]
  from which it follows that $\theta^{-1}\Phi_{\sO_X}(\sO_e(-1))$ is the unique
  non-split extension of $\sO_e(-1)$ by $\sO_X$.
\end{proof}

We of course also have an equivalence $\Phi_{\sO_e(-1)}:M\mapsto
\theta^{-1} M(e)$ taking $D^b_{\coh} X_{z,q}$ to $D^b_{\coh} X_{qz,q}$,
fixing $\sO_e(-1)$ and taking $\sO_X$ to $\theta^{-1}\sO_X(e)$.  As in the
case $Q$ smooth, these generate a faithful action of $\SL_2(\Z)$ on
$K_0^{\text{perf}}(Q)$, and thus $\Phi$ provides the desired extensions of
derived autoequivalences of elliptic surfaces.

\begin{cor}\label{cor:weird_langlands}
  For any element $\begin{pmatrix} a&b\\c&d\end{pmatrix}\in \SL_2(\Z)$,
    the noncommutative surfaces $X_{z,q}$ and $X_{z^aq^b,z^cq^d}$ are
    derived equivalent.
\end{cor}

\begin{rems}
  Here we should caution the reader that this does not correspond to an
  action of $\SL_2(\Z)$, as the relations of $\SL_2(\Z)$ correspond to
  nontrivial autoequivalences (acting on the Grothendieck group in the same
  way as suitable powers of $\theta$).  If we allow contravariant
  equivalences, we can of course include the duality $R\ad$, extending
  $\SL_2(\Z)$ to $\GL_2(\Z)$.  Note that since $R\ad$ is contravariant, it
  inverts the Serre functors, and thus each of $\Phi_{\sO_X}$
  or $\Phi_{\sO_{e_8}(-1)}$ conjugates to its inverse under $R\ad$.
  We should also note that the abelian equivalence
  $X_{z,q}\cong X_{1/z,1/q}$ does {\em not} correspond to an element of the
  above family, since the corresponding autoequivalence of $D^b\coh C$ is
  the action of a hyperelliptic involution.
\end{rems}

\begin{rems}
  It is useful to record how the above operations act on $K_0(X_{z,q})\cong
  K_0(X_{\sO_Q,\sO_Q})$, in terms of the numerical invariants (Definition
  \ref{defn:num_invs} below).  Note first that for any root $\alpha$ of
  $E_8\subset E_9$, the two generating equivalences preserve the class with
  numerical invariants $(0,\alpha,0)$, as this class is in both ${\cal
    K}_{\sO_{e_8}(-1)}$ and ${\cal K}_{\sO_X}$ and is invariant under
  $\theta$.  We thus need only determine how they act on the
  $4$-dimensional orthogonal complement in $K_0(X_{z,q})$ of that
  $8$-dimensional space, or in other words on classes of the form
  $(r,se_8+tQ,u)$.  We find that the action on such classes of the derived
  equivalence $D^b\coh X_{z,q}\to D^b\coh X_{qz,q}$ is
  \begin{align}
    (1,0,0)&\mapsto (1,-e_8-Q,0)\notag\\
    (0,e_8,0)&\mapsto (0,e_8,0)\notag\\
    (0,Q,0)&\mapsto (0,Q,-1)\notag\\
    (0,0,1)&\mapsto (0,0,1)\notag
  \end{align}
  while the action of $D^b\coh X_{z,q}\to D^b\coh X_{z,q/z}$ is
  \begin{align}
    (1,0,0)&\mapsto (1,-Q,0)\notag\\
    (0,e_8,0)&\mapsto (1,e_8-Q,0)\notag\\
    (0,Q,0)&\mapsto (0,Q,0)\notag\\
    (0,0,1)&\mapsto (0,Q,1).\notag
  \end{align}
  The images of $(0,Q,0)$ and $(0,0,1)$ in each case follow immediately
  from the action on parameters, and the remaining degree of freedom (given
  that we must respect $\chi(\_,\_)$) is determined by the known fixed
  element $(0,e_8,0)$, resp. $(1,0,1)$.  For the general equivalence
  $D^b\coh X_{z,q}\to D^b\coh X_{z^a q^b,z^c q^d}$, we cannot be quite as
  precise (due to the aforementioned central extension); all we can say is
  that
  \begin{align}
    (1,0,0)&\mapsto (d,-be_8+\frac{dh-b+c}{2} Q,\frac{-bh-b-a+d}{2})\\
    (0,e_8,0)&\mapsto (-c,ae_8+\frac{-ch+c+a-d}{2} Q,\frac{ah+b-c}{2})\\
    (0,Q,0)&\mapsto (0,dQ,-b)\\
    (0,0,1)&\mapsto (0,-cQ,a)
  \end{align}
  for some $h\in ab+bc+cd+2\Z$ (corresponding to the application of a power
  of $\theta$).  Note that this gives a non-split extension of $\SL_2(\Z)$
  by $\Z$, split (by taking $h=0$) over the index 2 (congruence) subgroup
  on which $ab+bc+cd\in 2\Z$.
\end{rems}

In particular, there is a derived equivalence $D^b_{\coh} X_{z,1}\cong
D^b_{\coh} X_{1,z}$ for any $X$.  By Theorem
\ref{thm:Painleve_moduli_spaces} below, $X_{z,1}$ may be interpreted (in a
countably infinite number of different ways) as a moduli space of
difference or differential equations, while $X_{1,z}$ is a noncommutative
deformation of the corresponding moduli space $X_{1,1}$ of ``Higgs
bundles'' (or other relaxed equations).  In the differential case, this is thus
precisely the sort of derived equivalence arising in geometric Langlands
theory, and is indeed a nontrivial special case thereof (with nonabelian
structure group $\GL_2$ and allowing the equation to have singularities,
not necessarily Fuchsian, with total complexity corresponding to four
Fuchsian singularities).  This suggests that our $\SL_2(\Z)$ of derived
equivalences is the shadow of a more general theory which would extend the
geometric Langlands correspondence to difference equations.

\medskip

We single out one particular composition of the above functors:
\[
\Psi:=\theta^{-1}\Phi_{\sO_{e}(-1)}\Phi_{\sO_X}\Phi_{\sO_{e}(-1)},
\]
or in other words
\[
\Psi:M\mapsto \Phi_{\sO_X}\theta M(-e))(-e)
\]
taking $D^b\coh X_{z,q}\to D^b\coh X_{q,1/z}$.

\begin{lem}\label{lem:Psi_values}
  We have $\Psi(\sO_e)\cong\sO_X$ and $\Psi(\sO_X)\cong\sO_{e}(-1)[-1]$.
\end{lem}

\begin{proof}
  Using the fact that $\theta$ commutes with any derived equivalence, we
  may rewrite the second claim using
  \[
  \Psi\cong \Phi_{\sO_{e}(-1)}\Phi_{\sO_X}\theta^{-1}\Phi_{\sO_{e}(-1)},
  \]
  and thus reduce to showing
  \begin{align}
    \Phi_{\sO_X}(\sO_{e}) &\cong \sO_X(e),\\
    \Phi_{\sO_X}(\sO_X(-e)) &\cong \sO_{e}(-1)[-1].
  \end{align}
  These are equivalent via the distinguished triangle
  \[
  \Phi_{\sO_X}(\sO_X(-e))\to \Phi_{\sO_X}(\sO_X)\to
  \Phi_{\sO_X}(\sO_{e})\to
  \]
  with $\Phi_{\sO_X}(\sO_X)\cong \sO_X$, so we need merely compute
  $\Phi_{\sO_X}(\sO_X(-e))$.  Since $\sO_X(-e)\in \sO_X^\perp$, this may be
  computed via the distinguished triangle
  \[
  \theta^{-1} \Phi_{\sO_X}(\sO_X(-e))\to
  \sO_X(-e)\to \Hom_k(R\Hom(\sO_X(-e),\sO_X),\sO_X) \to
  \]
  which since $R\Hom(\sO_X(-e),\sO_X)\cong k$ reduces to
  \[
  \theta^{-1} \Phi_{\sO_X}(\sO_X(-e))\cong \sO_e[-1],
  \]
  giving the desired formula.
\end{proof}

\begin{cor}
  Let $Y$ be a (possibly degenerate) commutative del Pezzo surface of
  degree 1 with anticanonical curve $Q$.  For $q\in \Pic^0(Q)$, let $Y_q$
  denote the corresponding noncommutative deformation, and let $X=X_{q,1}$
  denote the blowup in the point corresponding to $q$ in the identity
  component of $Q$.  Then $D^b_{\coh}(Y_q)\cong \sO_{X_{q,1}}^\perp$.
\end{cor}

\begin{proof}
  If we blow up the identity point of $Q$, we obtain $X_{1,q}$, and we can
  embed $D^b_{\coh}(Y_q)$ in $X_{1,q}$ via $\alpha_8^!$.  This
  identifies $D^b_{\coh}(Y_q)$ with $\sO_{e_8}^\perp$, and applying $\Psi$
  identifies this with $\sO_{X_{q,1}}^\perp$.
\end{proof}

\begin{rem}
  It follows more generally that the derived category of any noncommutative
  del Pezzo surface embeds in $\sO_X^\perp$ for some commutative surface
  $X$ with $K_X^2=0$.  (The results of \cite{OrlovD:2016} imply that any
  noncommutative surface has a derived embedding in {\em some} smooth
  projective variety, but the construction leads to varieties of rather
  high dimension.)
\end{rem}

\begin{cor}
  One has a natural isomorphism
  \[
  \Phi_{\sO_{e}(-1)} \Phi_{\sO_X} \Phi_{\sO_{e}(-1)}
  \cong
  \Phi_{\sO_X} \Phi_{\sO_{e}(-1)} \Phi_{\sO_X}.
  \]
\end{cor}

\begin{proof}
  Consider the two compositions
  \[
  \Psi\circ \Phi_{\sO_X}
  \qquad\text{and}\qquad
  \Phi_{\sO_{e}(-1)}\circ\Psi.
  \]
  Since $\Psi(\sO_X)=\sO_{e}(-1)[-1]$, we find that $\Psi$ takes
  $\sO_X^\perp$to $\sO_e(-1)^\perp$.  Since each of $\Phi_{\sO_X}$ and
  $\Phi_{\sO_{e}(-1)}$ acts as the internal Serre functor on the
  appropriate subcategory, and Serre functors are canonical, we conclude
  that the two compositions are naturally equivalent.  Applying $\theta
  \Phi_{\sO_{e}(-1)}^{-1}$ to both sides and simplifying gives the desired
  result.
\end{proof}

\begin{rem}
  The action on $K_0(X_{z,w})$ suggests that $\theta^{-2}\Psi^4[2]$ should
  be the identity.  If $2Q+e$ is ample, then this is true:
  $\theta^{-2}\Psi^4[2]$ preserves $\sO_X$ and commutes with
  $\theta^{-1}\Phi_{\sO_e(-1)}$ and thus preserves $\sO_X(r(2Q+e))$ for all
  $r$.  Since it preserves a set of generators, it must be trivial.  On the
  other hand, it is easy to see that if $M$ is any nonzero sheaf with
  $M\cong \theta M$, $R\Gamma(M)=R\Hom(\sO_{e}(-1),M)=0$, then $\Psi(M)=M$,
  so $\theta^{-2}\Psi^4 M[2]\cong M[2]\not\cong M$.  Such a sheaf exists
  whenever there is an effective root $\alpha$ of $E_8$.  Presumably one
  finds in such cases that $\theta^{-2}\Psi^4[2]$ is a composition of
  spherical twists.
\end{rem}

\medskip

One application of Corollary \ref{cor:weird_langlands} is the following
extension of Theorem \ref{thm:Picr} to the case of nonreduced special fiber
(and general $d$).  (Another, more significant, application is Theorem
\ref{thm:Painleve_moduli_spaces} below.)

\begin{cor}\label{cor:Picr}
  Let $z$ be an $l$-torsion point of $\Pic^0(Q)$ and let $d$ be an integer.
  Then the minimal proper regular model of the relative $\Pic^d$ of the
  (quasi-)elliptic fibration $X_{z,1}$ is isomorphic to the center of
  $X_{z^a,z^{l/\gcd(l,d)}}$, where $a$ is such that $ad\equiv \gcd(l,d)
  \pmod l$.
\end{cor}

\begin{proof}
  Set $b=(ad-\gcd(l,d))/l$, so that
  \[
  \begin{pmatrix}
    a&b\\
    l/\gcd(l,d)&d/\gcd(l,d)
  \end{pmatrix}
  \in
  \SL_2(\Z).
  \]
  Then Theorem \ref{thm:weird_langlands} induces a derived equivalence
  between $X:=X_{z,1}$ and $Y:=X_{z^a,z^{l/\gcd(l,d)}}$.  Moreover,
  since derived equivalences respect Serre functors, it follows
  that the (quasi-)elliptic fibration of $X$ induces a pencil of
  natural transformations $\theta^l\to\text{id}$, and equivalences
  between the triangulated categories of perfect objects on the
  respective fibers.  In particular, a $0$-dimensional sheaf on $Y$
  of degree $\gcd(l,d)$ and disjoint from $Q$ maps to a perfect
  complex of rank $1$ and degree $d$ on the corresponding fiber of
  $X$.  If the fiber on $X$ is smooth, then this perfect complex
  is itself a shifted simple sheaf, so necessarily a line bundle;
  since this is a derived equivalence, it follows that line bundles
  of degree $d$ on smooth fibers of $X$ map to $0$-dimensional
  sheaves on $Y$ of degree $\gcd(l,d)$.  Such a sheaf is supported
  on a unique point of the center $Z$ of $Y$, so that smooth fibers
  of $Z$ may indeed be identified with $\Pic^d$ of the corresponding
  fiber of $Z$.  (We should, of course, work over a family of such
  sheaves on $X$, but there is no difficulty.)  Since $Z$ is a
  minimal proper regular elliptic surface, it is the minimal proper
  regular model of the relative $\Pic^d$ of $X$ as required.

  If $X$ is quasi-elliptic, so there are no smooth fibers, then
  $\Pic^0(Q)\cong \G_a$, and thus we must have $l=\ch(k)=p$, with $p=2$ or
  $p=3$ since the surface is quasi-elliptic.  As $d$ only matters modulo
  $l$ and we can use $R^1\ad$ to negate $d$, we see that we need only
  consider $d=1$ or $d=0$.  The claim is obvious for $d=1$ (the relative
  $\Pic^1$ of $X$ is always $X$!), while for $d=0$, it remains only to show
  that the minimal proper regular model of the relative Jacobian is the
  center of $X_{1,z}$.  But in that case the derived equivalence is simple
  enough that we can directly check that it generically takes line bundles
  on fibers to sheaves.  Indeed, it suffices to consider the image of the
  structure sheaf of a fiber, which reduces (since the derived equivalence
  respects restriction to a fiber) to showing that the derived equivalence
  takes $\sO_X$ to $\sO_{e}(-1)$.  Apart from a shift, this is nothing
  other than Lemma \ref{lem:Psi_values}.
\end{proof}

\begin{rem}
  If $Q$ is reduced, then the center of $Y$ may be obtained by taking the
  images of the various parameters under the isogeny with kernel generated
  by $z$.  In general, the center can in principle be made explicit via our
  understanding of blowups of maximal orders: we need simply keep track of
  the morphism of anticanonical curves as we blow up points.
\end{rem}

\chapter{Families of surfaces}
\label{chap:families}

\section{Split families}

For a number of purposes, we will need to consider not just individual
noncommutative surfaces but also families; e.g., even for results purely
over fields, some of the results below require us to work over dvrs.  The
most satisfying theory applies in the case of a Noetherian base, although
again for technical reasons we will need to consider more general bases.

Here the main technical issue in the direct approach is that the theory of
blowups has only been developed over a field.  Indeed, there are a number
of places in van den Bergh's arguments where he does a case-by-case
analysis depending on the multiplicity of the point in the curve $Y$ or on
the size of its orbit under $\theta$, and neither of these is well-behaved
in families.  And, of course, there are further issues even over fields,
since one might wish to blow up a Galois cycle of points.  (One might also
consider the $p^l$-fold blowup of a point defined over an inseparable
extension, but even in the commutative case, there's no way to resolve the
singularity without a field extension.)  Further issues over fields already
arising in the commutative case are that geometrically ruled surfaces might
only be conic bundles and geometrically rational surfaces might not have a
rational ruling; one must also consider Brauer-Severi schemes of relative
dimension $2$.

Most of these issues can be dealt with via descent, so we can at least
initially consider only $\P^1$-bundles, blowups in single points, etc.
Thus for the moment, let $R$ be a general commutative ring\footnote{We will
  eventually be gluing over the \'etale topology, so there is no need to
  leave the affine case yet.}.  A {\em split} noncommutative plane over $R$
will be a category $\qcoh X$ with choice of structure sheaf $\sO_X$,
defined in the following way.  Let $Q\subset \P^2_R$ be a cubic curve, and
let $q$ be an invertible sheaf on $Q$ which on every geometric fiber is in
the identity component of $\Pic^0(Q)$.  Then we can define a $\Z$-algebra
over $R$ as in \cite{BondalAI/PolishchukAE:1993}, and can take $\qcoh X$ to
be the appropriate quotient of the category of modules over that
$\Z$-algebra.

One difficulty that arises is showing that the result is actually flat.
The main issue is that since the $\Z$-algebra is defined via a
presentation, there is no a priori reason why the $\Z$-algebra should be
flat.  We would also like to have the semiorthogonal decomposition, which
reduces to establishing an exact sequence of the form
\[
0\to \sO_X(d-3)\otimes_R L_d\to \sO_X(d-2)\otimes_R V_d\to
\sO_X(d-1)\otimes_R \Gamma(q^{d+1}(1))
\to \sO_X(d)\to 0
\]
where $V_d$, $L_d$ are determined by exact sequences
\[
0\to V_d\to \Gamma(q^{d}(1))\otimes_R \Gamma(q^{d+1}(1))
\to \Gamma(q^{2d+1}(2))\to 0
\]
and
\begin{align}
0\to L_d&{}\to \Gamma(q^{d-1}(1))\otimes_R  \Gamma(q^{d}(1))\otimes_R
\Gamma(q^{d+1}(1))\notag\\
&{}\to
\Gamma(q^{2d-1}(2))\otimes_R \Gamma(q^{d+1}(1))
\oplus
\Gamma(q^{d-1}(1))\otimes_R \Gamma(q^{2d+1}(2)).
\end{align}

If $R$ is Noetherian, we can establish both facts by a common induction.
Indeed, the $R$-module $\Hom(\sO_X(-n),\sO_X(d))$ is defined inductively as
the cokernel of the map
\[
\Hom(\sO_X(-n),\sO_X(d-2))\otimes_R V_d\to
\Hom(\sO_X(-n),\sO_X(d-1))\otimes_R \Gamma(q^{d+1}(1))
\]
(with $\Hom(\sO_X(-n),\sO_X(d))=0$ for $d<-n$ and
$\Hom(\sO_X(-n),\sO_X(-n))=R$), and thus flatness and exactness reduce to
showing that the complex
\begin{align}
\Hom(\sO_X(-n),\sO_X(d-3))\otimes_R L_d
&{}\to
\Hom(\sO_X(-n),\sO_X(d-2))\otimes_R V_d\notag\\
&{}\to
\Hom(\sO_X(-n),\sO_X(d-1))\otimes_R \Gamma(q^{d+1}(1))
\end{align}
represents a flat $R$-module.  By induction, this is a complex of finitely
generated flat $R$-modules, and thus represents a (finitely generated) flat
$R$-module iff it represents a vector space over every closed point.  We
may thus reduce to the case of a field, which we can even assume is
algebraically closed (since the desired claim descends through fpqc base
change), and thus reduce to the known case.  Similarly, the requisite
vanishing of Ext spaces reduces to considerations of complexes of f.g. flat
$R$-modules, so again can be checked on geometric fibers.  In this way, we
find that $\sO_X(-2),\sO_X(-1),\sO_X$ is still a full exceptional
collection, with gluing functors as expected.

More generally, this property is inherited by any base change, and as
already mentioned descends through fpqc, and thus smooth base change.  Thus
if there is a smooth cover $R'$ of $R$ such that $(C_{R'},q_{R'})$ is a
base change of a family over a locally Noetherian family, then we still
obtain a flat family of categories, and the derived category still has the
expected exceptional collection.  But {\em any} family has this property!
Indeed, cubic curves in a fixed $\P^2$ are classified by a projective fine
moduli space (isomorphic to $\P^9$, to be precise), and thus pairs $(C,q)$
are classified by a smooth group scheme over that moduli space.  The latter
is not quite a fine moduli space (there is a Brauer obstruction to a point
of $\Pic^0$ being represented by a line bundle), but this can be resolved
by an \'etale base change.  We thus see that any family of pairs $(C,q)$ is
\'etale-locally given by a base change from the universal family over the
Noetherian moduli stack of such pairs, and thus inherits flatness and the
exceptional collection.

We can also recover the embedding of $Q$ in $X$ by considerations from
geometric fibers.  The existence of a map $\Hom_X(\sO_X(d'),\sO_X(d))\to
\Hom_Q(q^{d'(d'+3)/2}(d'),q^{d(d+3)/2}(d))$ is immediate from the
construction, and when $d'=d-3$ is surjective with invertible kernel (since
on geometric fibers it is surjective with 1-dimensional kernel).  This
establishes the requisite natural transformation ${-}(-Q)\to\text{id}$, and
we easily see that $\Hom_X(L,L'|_Q)$ behaves as expected when $L,L'$ are
line bundles, so recovers a twisted homogeneous coordinate ring of $Q$.

For inductive purposes, we would like to know that the resulting category
$\qcoh X$ is locally Noetherian.  Unfortunately, the proof of
\cite{ArtinM/TateJ/VandenBerghM:1991} for the field case does not quite
carry over, as it uses valuations on the field in a significant way (in
particular, \cite[Lem.~8.9]{ArtinM/TateJ/VandenBerghM:1991}).  Luckily,
those parts of the argument are in service of showing that the quotient of
the graded algebra by the normalizing element of degree $3$ is Noetherian;
we now recognize that quotient as a twisted homogeneous coordinate ring
\cite{ArtinM/VandenBerghM:1990}, and thus it is Noetherian since its
$\Proj$ is Noetherian.  We thus conclude that any ``split'' noncommutative
plane over a Noetherian base ring is locally Noetherian.  (This is a
stronger form of the ``strongly Noetherian'' condition: any Noetherian base
change is still a family over a Noetherian base, so remains Noetherian.)

\smallskip

The argument for quasi-ruled surfaces is analogous.  Given a pair of
projective smooth curves $C_0/R$, $C_1/R$ with geometrically integral
fibers and a sheaf ${\cal E}$ on $C_0\times_R C_1$ which is a sheaf
bimodule on every geometric fiber, we can use the construction of van den
Bergh to get a global sheaf $\Z$-algebra, and flatness, the semiorthogonal
decomposition, and the embedding of $Q$ as a divisor follow over a
Noetherian base using the analogue of the exact sequence of Lemma
\ref{lem:exact_tri_quasi_ruled}.  (Here when checking that the $R\Hom$
complexes behave correctly, we need only check sheaves of the form
$\rho_d^* L$ and can further arrange for the line bundle on the source side
to be arbitrarily negative, since it suffices to check generators; this
lets us reduce to cases in which there is a single cohomology sheaf.)
Again, the corresponding moduli stack is itself locally Noetherian, so the
claims follow in general.  Moreover, over a Noetherian base, the category
is again locally Noetherian; in this case, the argument from the literature
\cite{ChanD/NymanA:2013} can be adapted with no difficulty (since they were
already considering constant families over general Noetherian
$k$-algebras).  For technical reasons, we will define a ``split'' family of
quasi-ruled surfaces to be a surface obtained from the above data in which
we have also marked points $x_i\in C_i(R)$.  (This ensures that every
component of the moduli space of line bundles on $X$ has a point.  Of
course, this is only relevant data for quasi-ruled surfaces themselves; on
blowups, we can simply take the images of the first point to be blown up!)

\smallskip

The case of blowups is somewhat trickier, as the relevant graded bimodule
algebra is no longer given by a presentation.  However, we can show by a
similar induction (again for $R$ Noetherian) that for any line bundle $L$
on $X$, $\alpha_*(L(d))$ is $R$-flat, and thus since line bundles generate
$\coh X$, the bimodule algebra is again $R$-flat.  (For the inductive
argument, we need to represent $\alpha_*(L(d))$ by a complex of finite sums
of line bundles, apply the appropriate functor to obtain $R\alpha_*L(d+1)$,
and observe by checking geometric fibers that the cone of the natural
transformation $(R\alpha_*L(d+1))(-Q)\to \alpha_*(L(d))$ is an $R$-flat
coherent sheaf, and thus so is $R\alpha_*L(d+1)$.)  The semiorthogonal
decomposition is similar, although now we need analogues of the two triangles
\[
L\alpha^* L\to L\alpha^! L\to \sO_e(-1)\to
\]
and
\[
\sO_e(-1)[1]\to L\alpha^* \sO_{\tau p}\to \sO_e\to
\]
in order to establish that $L\alpha^*\coh X$ and $\sO_e(-1)$ generate the
derived category.  Here the sheaf $\sO_e(-1)$ is essentially defined by the
special case $L=\sO_X$ of the first triangle (with $\sO_e:=\sO_e(-1)(1)$),
and is an $R$-flat sheaf since it is a sheaf on every geometric fiber.
Moreover, semicontinuity shows that $\Hom(\alpha^!L/\alpha^*L,\sO_e(-1))$
is a line bundle on $R$ (presumably the fiber of $L$ over the appropriate
point of $Q$), and we may verify that
\[
(\alpha^!L/\alpha^*L)\otimes_R \Hom(\alpha^!L/\alpha^*L,\sO_e(-1))\to
\sO_e(-1)
\]
is an isomorphism by checking geometric fibers.  The second triangle is
somewhat trickier, but we may rewrite it (up to line bundles on $R$) as
\[
\sO_e(-1)[1]\to L\alpha^*R\alpha_* \sO_e\to \sO_e\to
\]
and represent the latter map as a chain map of complexes of finitely
generated $R$-flat modules, and thus again reduce to checking things
geometrically.  And once more the relevant moduli stack is locally
Noetherian (see Proposition \ref{prop:moduli_stack_of_surfaces_is_small}
below), so the general case follows from the Noetherian case.

Inducting along blowups requires that we show that the blowup is again
locally Noetherian as long as the base ring is Noetherian.  Here the
argument of \cite{VandenBerghM:1998} is more difficult to adapt, as it
involves the quotient by a certain invertible ideal that depends on the
multiplicity of the point being blown up in the ambient curve, and this is
not flat.  Luckily, there is a different (and in many ways more natural)
ideal we can use instead, which being flat allows the argument to carry
over.

\begin{prop}
  Let $X/R$ be a split family of noncommutative rational or rationally
  quasi-ruled surfaces over a Noetherian ring $R$.  Then $\qcoh(X)$ is
  locally Noetherian.
\end{prop}

\begin{proof}
  Per the above discussion, it remains only to show that this is inherited
  by blowups.  We use the notation of \cite{VandenBerghM:1998}, and note
  that the basic properties we need are easily deduced from the fact that
  they hold on geometric points.  Thus $\qcoh(\widetilde{X})$ has been
  represented as $\Proj({\cal D})$, where ${\cal D}$ is a suitable ``graded
  bimodule algebra on $X$'', and the main missing ingredient in the
  argument is an invertible ideal in ${\cal D}$ {\em defined over $S$} such
  that the quotient is Noetherian.  Consider the invertible ideal ${\cal
    D}(-1)$.  If we can show that the $\Proj$ of the quotient is a
  commutative curve over $S$, then the $\Proj$ will be Noetherian.  On any
  geometric fiber, this follows from Proposition
  \ref{prop:nice_divisor_on_blowup}, and thus is true globally.  To lift to
  the graded quotient itself, we need simply observe that since $S$ is
  Noetherian, there is a maximum degree in which the graded algebra on any
  fiber differs from its saturation.  The saturation is certainly
  Noetherian, and every homogeneous component where there is a difference
  is Noetherian, and thus the original graded algebra is Noetherian.
\end{proof}

\begin{rems}
  In fact, one can check more directly that the quotient by ${\cal D}(-1)$
  is equal to its saturation in all but degree 0, where the saturation has
  $o_Y$ in place of $o_X$.  Furthermore, the saturation (which since $Y$ is
  commutative may be rephrased in terms of graded coherent sheaves on
  $Y\times Y$) is essentially a relative version of a twisted homogeneous
  coordinate ring (based on the fact that the autoequivalence
  $M\mapsto M(Y^+)=M(1)$ is relatively ample for $Y^+\to Y$).
\end{rems}

\begin{rems}
  In general, the ideal van den Bergh uses is one of the sequence of
  invertible ideals ${\cal D}(-m)((m-1)Y)$.  This is essentially ascending,
  as there is only one degree where ${\cal D}(-m)((m-1)Y)\not\subset {\cal
    D}(-m-1)(mY)$.  It thus describes a descending chain of closed
  subschemes, the limit of which is $\tilde{Y}$.  (For general $m$, the
  quotient is obtained by removing $\min(m-1,\mu-1)$ copies of the
  exceptional curve from $Y^+$.)  Only the case $m=1$ is flat, however.
\end{rems}

\begin{rems}
  Considering the commonalities between the proofs in the different cases,
  it seems likely that it should be possible to prove a quasi-scheme is
  Noetherian whenever it has a morphism to a Noetherian quasi-scheme, a
  relatively ample autoequivalence $\tau$, and a natural transformation
  $T:\text{id}\to\tau$ such that the category of sheaves with $T_M=0$ is
  also Noetherian.  (It may also be necessary to have $T_{\tau M}=\tau
  T_M$.)  In the three cases, $\tau(M)=M(Q)$, $T$ is the obvious natural
  transformation, and the first morphism is $\Gamma$, $\rho_*$, or the
  blow-down respectively.
\end{rems}

\section{Perfect objects over families}

Since in each case we retain the semiorthogonal decomposition over a
general base, we can give an alternate description of these categories as
$t$-structures on derived categories of suitable dg-algebras.  Consider
the quasi-ruled case.  For each $i$, the category $\perf(C_i)=D^b_{\coh}
C_i$ is generated by $G_i:=\sO_{C_i}\oplus \sO_{C_i}(-x_i)$.  This lets us
represent $D_{\qcoh}(C_i)$ as (the homotopy category of) the category of
right modules over the finite type dg-algebra $R\End(G_i)$ (which we think
of as acting on $G_i$ on the left).  Then $\rho_0^*G_0\oplus \rho_1^*G_1$
is a compact generator of the quasi-ruled surface, and has endomorphism
dg-algebra of the form
\[
\begin{pmatrix}
  R\End_{C_0}(G_0) & R\Hom_{C_0\times C_1}(\pi_1^*G_1,\pi_0^*G_0\otimes {\cal E})\\
  0 & R\End_{C_1}(G_1)
\end{pmatrix}.
\]
In particular, the derived category of modules over this dg-algebra agrees
with $D_{\qcoh}(X)$, with the best analogue of $D^b_{\coh}(X)$ being the
category of {\em perfect} modules over the dg-algebra.  (Over a general
base, ``coherent'' is not very well-behaved.)

We may deal with blowups similarly; if we have already constructed a sheaf
of dg-algebras $A_{m-1}$ corresponding to a family of surfaces $X_{m-1}$,
then any object $\sO_{x_m}$ in $\perf(A_{m-1})$ gives rise to a new sheaf
of dg-algebras of the form
\[
A_m=
\begin{pmatrix}
  A_{m-1} & R\Hom(\sO_{x_m},R)[-1] \\
      0   &          R
\end{pmatrix},
\]
with associated family $X_m$.  (This corresponds to the gluing functor
$R\Hom_{X_m}(\sO_{e_m}(-1),L\alpha_m^*M)\cong R\Hom_X(M,\sO_{x_m}[1])^*$
for $M\in \perf(X_{m-1})$.)  If every geometric fiber of $\sO_{x_m}$ is the
structure sheaf of a point on a divisor on the corresponding fiber of
$X_{m-1}$, then this dg-algebra is again Morita-equivalent to the blown up
surface.

We can again recover the $t$-structure without reference to the original
abelian category.  Indeed, the definition of line bundles extends easily
from the field case (the only change being that anything \'etale locally
isomorphic to a line bundle should be considered a line bundle), and we see
that line bundles generate the abelian category, and further that $M\in
D_{\qcoh}(X)^{\ge 0}$ iff $R\Hom(L,M)\in D_{\qcoh}(R)^{\ge 0}$ for all line
bundles $L$ (or, say, only those coming from line bundles
$\sO_{C_i}(nx_i)$).  Equivalently, $D_{\qcoh}(X)^{\le 0}$ is the smallest
``cocomplete preaisle'' containing the given set of line bundles; that is,
it is the smallest subcategory containing $L[p]$ for all $p<0$ and closed
under extensions and arbitrary direct sums.  By
\cite[Thm.~A.1]{TarrioLA/LopezAJ/SalorioMJS:2003}, a cocomplete preaisle
generated by a set of compact generators is an ``aisle'' (the negative part
of an actual $t$-structure), and thus in particular this holds for the
aisle generated by line bundles.  (Of course, we already knew this, but
this argument does not require that we have a direct construction of the
abelian category.)

One useful consequence (which again could be shown using the construction
of the abelian category, but is particularly simple using the dg-algebra
description) is that if $M,N\in \perf(X)$, then $R\Hom(M,N)\in \perf(R)$,
and thus satisfies semicontinuity.  (This, of course, is the natural
generalization of the fact that over a field, $\Hom$ complexes are
finite-dimensional.)

\begin{prop}
  Let $\perf(X)$ be the dg-category of perfect complexes on a flat family
  of rational or rationally quasi-ruled surfaces over $R$.  Then for any
  $M,N\in \perf(X)$ and any integer $p$, the map
  \[
  s\mapsto \dim\Ext^p_{X_s}(M\otimes^L k(s),N\otimes^L k(s))
  \]
  is an upper semicontinuous function of $s\in \Spec(R)$.  Moreover, if
  $\Ext^q(M,N)=0$ for $q>p$ and $\dim\Ext^p_{X_s}(M\otimes^L
  k(s),N\otimes^L k(s))$ is constant, then $\Ext^p(M,N)$ is a locally free
  $R$-module.
\end{prop}

\begin{proof}
  We need simply note that
  \[
  \Ext^p_{X_s}(M\otimes^L k(s),N\otimes^L k(s))
  \cong
  \Ext^p_X(M,N\otimes^L k(s))
  \cong
  h^p(R\Hom(M,N)\otimes^L k(s)),
  \]
  and thus reduce to standard semicontinuity results in $\perf(R)$.
\end{proof}

Of course, this is only useful when combined with an understanding of the
perfect objects.  An object $M\in D_{\qcoh}(X)$ is perfect iff its
projections to the various pieces of the semiorthogonal decomposition are
perfect.  Each of these lives on a scheme which is smooth and projective
over $R$, and thus we may apply the following.

\begin{prop}
  Let $Y$ be a smooth projective scheme over $R$ of relative dimension $n$.
  An object $M\in D_{\qcoh}(Y)$ is perfect iff there are integers $[a,b]$
  such that for any point $s\in \Spec(R)$, $M\otimes^L_R k(s)$ has
  coherent cohomology supported on the interval $[a,b]$.
\end{prop}

\begin{proof}
  An object is perfect iff $M\otimes^L_{\sO_Y} k(y)$ has coherent
  cohomology of uniformly bounded amplitude as $y$ ranges over points of
  $Y$.  If $\pi:Y\to \Spec(R)$ is the structure map, then we may write
  \[
  M\otimes^L_{\sO_Y} k(y)
  \cong
  M\otimes^L_R k(\pi(y))
  \otimes^L_{\sO_{Y_{\pi(y)}}} k(y).
  \]
  Since $Y_s$ is a projective scheme over a field, $k(y)$ has $\Tor$
  dimension at most $n$.  Thus if $M\otimes^L_R k(\pi(y))$ has coherent
  cohomology supported on $[a,b]$, then $M\otimes^L_{\sO_Y} k(y)$ has
  coherent cohomology supported on $[a-n,b]$, giving one implication.  For
  the converse, note that for any point $s\in \Spec(R)$, if
  $h^p(M\otimes^L_R k(s))$ is the nonzero cohomology sheaf of lowest
  degree and $y$ is a generic point of a component of its support in $Y_s$,
  then $h^{p-\text{codim}(y)}(M\otimes^L_{\sO_Y} k(y))\ne 0$.
\end{proof}

In our case, every fiber has a $t$-structure, and we have the following.

\begin{prop}
  Let $X/R$ be a family of noncommutative rational or rationally
  quasi-ruled surfaces as above.  Then an object $M\in D_{\qcoh}(X)$ is
  perfect iff there are integers $[a,b]$ such that $M\otimes^L_R k(s)\in
  D_{\qcoh}(X_s)$ has coherent cohomology supported on the interval $[a,b]$
  for all $s\in \Spec(R)$.
\end{prop}

\begin{proof}
  On each fiber and for each piece of the semiorthogonal decomposition, the
  projection and the inclusion functors have finite cohomological
  amplitude, and this amplitude is uniformly bounded on $\Spec(R)$.  Thus
  $M\otimes^L_R k(s)$ has uniformly bounded amplitude iff its projections
  have uniformly bounded amplitude.
\end{proof}

\begin{cor}
  If $M\in D_{\qcoh}(X)$ is a flat family of coherent sheaves (i.e.,
  $M\otimes^L_R k(s)$ is a coherent sheaf for all $s\in \Spec(R)$), then
  $M\in \perf(X)$.
\end{cor}

\begin{rem}
  In particular, flat families of coherent sheaves satisfy semicontinuity.
  When $R$ is Noetherian, $\qcoh(X)$ is locally Noetherian, and thus we can
  define $\coh(X)$ to be the subcategory of Noetherian objects, and
  automatically have $\perf(X)\subset D^b_{\coh}(X)$.  Over a more general
  base, ``coherent'' is both tricky to define and of unclear use.
\end{rem}

There is of course no real difficulty in gluing the above construction in
the Zariski topology.  We thus make the following definition.  Note that
the data we used to construct a family of noncommutative planes came with a
choice of curve $Q$, and we make a similar choice as necessary for
quasi-ruled surface, which we include as part of the data.  On fibers which
are not commutative, we can recover $Q$ from the surface, so it does not
give extra data, but commutative fibers come with a choice of anticanonical
curve.

\begin{defn}
  A {\em split family of noncommutative rationally quasi-ruled surfaces} over a
  scheme $S$ is an $\sO_S$-linear category $\qcoh(X)$ (or $D_{\qcoh}(X)$
  with a $t$-structure) which is Zariski locally of the above form.
\end{defn}

\begin{defn}
  A {\em split family of noncommutative rationally ruled surfaces} is a split
  family of noncommutative rationally quasi-ruled surfaces equipped with an
  isomorphism $C_0\cong C_1$ identifying $x_0$ with $x_1$ and making
  $\overline{Q}$ algebraically equivalent to the double diagonal; in addition,
  $Q$ may not be geometrically integral on any fiber such that $g(C_0)=1$.
\end{defn}

\begin{defn}
  A {\em split family of noncommutative rational surfaces} is either an 
  iterated blowup of a family of noncommutative planes as constructed
  above, or a split family of noncommutative rationally ruled surfaces such
  that $g(C_0)=0$.
\end{defn}

\begin{rem}
  Since we do not deal with any other kind of noncommutative surface, we
  will feel free to abbreviate this to ``split family of noncommutative
  surfaces'', only adding ``rational'', ``ruled'', or ``quasi-ruled'' when
  we are imposing further restrictions.
\end{rem}

\begin{rem}
  As we have already mentioned, we will show below that a blowup of a
  noncommutative plane is a noncommutative Hirzebruch surface, and thus a
  (split) noncommutative rational surface is either a noncommutative
  rationally ruled surface over $\P^1$ or a noncommutative plane.
\end{rem}
  
\section{General families}

As we have already discussed, the above construction is clearly too
restrictive to be the right notion of family.  An obvious further approach
is to consider surfaces that can locally be put in the above form relative
to some suitable topology.  To ensure we make the correct choice of
topology here, we should first determine the correct choice in the
commutative case.

We first observe that $C_0$ and $C_1$ have sections \'etale locally
(indeed, {\em any} smooth surjective morphism has a section \'etale
locally), so that $x_0$ and $x_1$ can be chosen \'etale locally.  Moreover,
they only played a role in constructing a set of generating line bundles on
the surface, and thus the $t$-structure, but the resulting $t$-structure
was independent of the choices, so descends.  For the same reason, a flat
family of smooth surfaces geometrically isomorphic to $\P^2$ again has an
\'etale local section, and is thus \'etale locally isomorphic to $\P^2$,
and the $t$-structure defined using that isomorphism descends.

A somewhat more subtle question has to do with the fact that a
geometrically ruled surface need not be a projective bundle.  Note that
although we are assuming a particular choice of ruling, there is no loss of
generality when $g(C)>0$, since then the canonical map $X\to \Alb^1(X)$ to
the Albanese torsor factors through $C$ on any geometric fiber, so that we
actually have a canonical ruling.

\begin{prop}
  Let $C/S$ be a smooth proper curve, and let $\rho:X\to C$ be a smooth
  proper curve which over geometric points of $C$ is isomorphic to $\P^1$.
  Then any geometric point $s\in S$ has an \'etale neighborhood over which
  $X$ is isomorphic to a projective bundle.
\end{prop}

\begin{proof}
  First note that $C$ \'etale locally has a section, so may be assumed
  projective, and $-K_X$ is relatively ample over $C$, so that $X$ is then
  also projective.  The anticanonical embedding of $X_s/C_s$ is a conic,
  and thus by Tsen's Theorem has a section $\sigma:C_s\to X_s$.  This
  induces a line bundle $\sO_{X_s}(\sigma(C_s))$ such that $X_s\cong
  \P(\rho_*\sO_{X_s}(\sigma(C_s))$.  Let $\Pic^\sigma_{X/S}$ denote the
  component of the Picard scheme $\Pic_{X/S}$ containing the isomorphism
  class of this line bundle.  This is a torsor over $\Pic^0_{C/S}$, thus
  smooth and surjective, so has a section over an \'etale neighborhood $S'$
  of s in S, and the isomorphism class of line bundles has a representative
  ${\cal L}$ defined over an \'etale neighborhood $T$ of $s$ in $S'$.  The
  sheaf of graded algebras over $C$ corresponding to ${\cal L}$ is
  naturally isomorphic to the symmetric algebra of $\rho_*{\cal L}$, so
  that $X_T\cong \P_C(\rho_*{\cal L})$ as required.
\end{proof}

\begin{rem}
  It is tempting to simply specify the component $\Pic^\sigma(X/S)$ by
  numerical considerations (i.e., specifying the Hilbert polynomials of
  ${\cal L}$ and ${\cal L}(f)$).  The issue is that the resulting component
  could be (and, half the time, is!) empty; we need a splitting over some
  geometric point to pick out a good component.
\end{rem}

The hardest case is when there is no longer a unique way to interpret the
commutative surface as an iterated blowup of a ruled surface.  This is
likely to be especially difficult to deal with for rational surfaces, as in
that case there may in fact be infinitely many such interpretations, and
thus the moduli stack of anticanonical rational surfaces is not algebraic!
(If one imposes a blowdown structure, this is no longer an issue, as we saw
in Theorem \ref{thm:antican_moduli_lci}.)  To get a feel for the issue,
consider the case of a {\em projective} family over a connected base with
generic geometric fiber $\P^1\times \P^1$.  If the family has geometric
fibers of the form $F_d$ for $d>0$, then the ample bundle must have
bidegree $(d,d')$ with $d>d'$, and this distinguishes the two rulings on
the generic fiber, so that there is a unique choice of global ruling.  If
all fibers are $\P^1\times \P^1$, then the line bundle of bidegree $(1,1)$
extends \'etale locally, and gives an embedding of the surface as a smooth
quadric surface.  Over a non-closed field, a typical such surface has no
ruling, and thus the only rational line bundles are multiples of the
bidegree $(1,1)$ bundle.  However, the family has a natural double cover
over which it has a ruling, which in odd characteristic simply takes the
square root of the determinant.  In particular, by removing the singular
surfaces (which correspond to smooth surfaces on which the given bundle is
not ample), we have made this cover \'etale.  (If we had not removed the
singular surfaces, the cover would be ramified, but the family would not be
smooth, and does not have a global desingularization.)  Over the cover, the
surface is isomorphic to a product of two conics, and thus at worst two
further \'etale covers make it globally $\P^1\times \P^1$.

More generally, as long as the family is projective, the ample divisor
$D_a$ cuts out a finite set of blowdown structures on the geometric fibers
of $X$, namely those for which $D_a$ is in the fundamental chamber (as
defined in Section \ref{sec:nef_comm}; see also Chapter \ref{chap:effnef}
below).  Moreover, the resulting set of blowdown structures is a torsor
over an appropriate finite Coxeter group, and thus a blowdown structure
exists \'etale locally.

\medskip

With the above discussion in mind, we define a family of noncommutative
rational or rationally quasi-ruled surfaces (again abbreviated to ``family
of noncommutative surfaces'') over $S$ to be a family $D_{\qcoh}(X)$ of
dg-categories over $S$ such that (a) the dg-category is \'etale locally a
split family of noncommutative surfaces, and (b) the associated
$t$-structures are compatible, in that the two ways of pulling back to the
fiber square of the cover agree.

Note that the locally Noetherian condition easily descends through
\'etale covers, and thus we have the following, over a Noetherian base.

\begin{thm}\label{thm:noetherian_base_implies_noetherian}
  Let $X/S$ be a family of noncommutative rational or rationally
  quasi-ruled surfaces over a Noetherian base.  Then $\qcoh(X)$ is locally
  Noetherian.
\end{thm}
  
In principle, we could construct such a category by specifying the specific
data $(C_0,C_1,\dots)$ on each piece of the \'etale cover, and giving
suitable gluing data.  (Indeed, by fpqc descent, quasicoherent sheaves can
always be glued along \'etale covers.)  For rational or rationally ruled
surfaces, this is not too difficult: we have an associated family of
commutative surfaces, so as long as that family is projective, there is an
\'etale cover $T/S$ over which there is a global choice of blowdown
structure.  Since the associated tuple $(C_0,C_1,\dots)$ over $T$ satisfies
the natural compatibility conditions, the result indeed glues to give a
dg-category with $t$-structure over $S$.

It turns out that the dg-category itself is easy enough to construct
directly.  Let us consider first the case that $X$ is rational.  Then the
exceptional object $\sO_X$ of $\perf(X\times_S T)$ certainly descends, and
thus induces a semiorthogonal decomposition over $S$, suggesting the
following construction.  Given a family $Y/S$ of commutative rational
surfaces (i.e., $Y$ is smooth and proper over $S$ of relative dimension
$2$, and its fibers are geometrically rational), not only does $\perf(Y)$
have a generator over $S$, but so does the subcategory $\sO_Y^\perp$;
indeed, if $G$ generates $Y$, then the cone $G'$ of $\sO_Y\otimes
R\Hom(\sO_Y,G)\to G$ generates $\sO_Y^\perp$.  If $Q$ is an anticanonical
divisor on $Y$ and $q\in \Pic^0(Q)$ (that is, $q$ is a line bundle on $Q$
which has degree 0 on every component of every geometric fiber of $Q$),
then we may consider the sheaf of triangular dg-algebras
\[
\begin{pmatrix}
  \sO_S & R\Hom_Q(G'|^{\bf L}_Q,q)\\
     0  & R\End(G')
\end{pmatrix}
\]
This {\em almost} gives a fully general construction of $D_{\qcoh}(X)$ for
noncommutative rational surfaces.  The only remaining issue is that all we
actually obtain from the \'etale local description is a point of
$\Pic_{(Q/S)(\text{\'et})}(S)$ which is in the identity component of every
fiber; in particular, it defines a line bundle on the base change to some
\'etale cover such that the two pullbacks are isomorphic, but the
isomorphisms need not be compatible.  (Basically, the problem is that $q$
is normally given as the determinant of the restriction of a point on the
noncommutative surface, but $Q$ may not have rational points, so we only
know $q$ up to isomorphism!)  This is not too hard to fix as long as $Q$ is
projective (e.g., if $Y$ is projective); the idea is that although $q$ may
not descend, the objects $q\otimes R\Gamma(q^{-1}(d))\in \perf(Q_{T})$ and
$R\End(R\Gamma(q^{-1}(d)))$ descend for any $d$.  Taking $d\gg 0$ gives a
well-defined sheaf of dg-algebras
\[
\begin{pmatrix}
  \sO_S & R\Hom(G'|^{\bf L}_Q,q\otimes \Gamma(q^{-1}(d)))\\
  0 & R\End(G')\otimes \End(\Gamma(q^{-1}(d))).
\end{pmatrix}
\]
This is clearly Morita-equivalent to the original dg-algebra over $T$, and
thus the category of (perfect) modules over this dg-algebra may indeed be
viewed as a family of (derived categories of) noncommutative rational
surfaces, in that every geometric fiber is such a surface.  Moreover, any
category $\perf(X)$ which is \'etale locally equivalent to the category of
perfect objects on an iterated blowup of a noncommutative Hirzebruch
surface (or noncommutative plane) can be obtained via the above
construction.

So we may instead construct families of noncommutative rational surfaces by
specifying a triple $(Y,Q,q)$ (with $q\in \Pic(Q/S)(S)$) along with a
splitting over some \'etale cover $T/S$ such that both pullbacks of the
$t$-structure agree.  If $T$ is Noetherian, then we will see below that
$\qcoh(X)$ is equivalent to the $\Proj$ of a sheaf of $\Z$-algebras.

For ruled surfaces which are not rational, Proposition
\ref{prop:f_is_unique_if_quasi-ruled} below tells us that the functor
$\alpha_m^*\cdots \alpha_1^*\rho_0^*$ descends (i.e., noncommutative ruled
surfaces of genus $>0$ also have canonical rulings), and thus we may use
the corresponding semiorthogonal decomposition to construct $\perf(X)$ from
a commutative family $Y$ of geometrically rationally ruled surfaces over a
curve $C$ of genus $>0$, an anticanonical divisor $Q$, and a morphism $Q\to
C$ which contracts vertical fibers and is in the same component of the
moduli stack of such morphisms as the morphism coming from the ruling on
$Y$.  (By Lemma \ref{lem:Pic0Q/Pic0C}, such morphisms are classified,
modulo the action of $\Aut^0(C)$, by $\Pic^0(Q)/\Pic^0(C)$.)

We should note that even when $Y$ is geometrically ruled, the subcategory
$\perf(C)^\perp$ is not in general equivalent to $\perf(C)$.  Indeed, if
$Y/C$ is a nonsplit conic bundle, then $\perf(C)^\perp$ is the category of
perfect objects in the derived category of the corresponding (quaternion)
Azumaya algebra over $C$.

We ignore families of surfaces which are rationally quasi-ruled but not
rationally ruled, as the fibers themselves are not particularly interesting
as noncommutative surfaces (and there is no intrinsic difficulty in
constructing them in any case; Proposition
\ref{prop:f_is_unique_if_quasi-ruled} implies that the ruling
is again canonical).  Note that we
cannot quite dismiss these as being maximal orders on families of
rationally ruled commutative surfaces, as the degree over the center can
fail to be constant: see the discussion above of quasi-ruled surfaces of
``2-isogeny'' type.

\section{The moduli stack of noncommutative surfaces}

The one issue with the above notion of family, though more natural than the
split case, is that the corresponding moduli problem is not well-behaved.
We can resolve this by putting back part of the splitting data.  A {\em
  blowdown structure} on a family $X/S$ of noncommutative ruled surfaces is
an isomorphism class of sequences
\[
X=X_m\to X_{m-1}\to\cdots\to X_0\to C
\]
of morphisms of families, in which all but the last map is a blowdown, and
$X_0\to C$ is a ruling.  This is essentially combinatorial data; as we will
see, the blowdown structures on a noncommutative ruled surface over an
algebraically closed field form a discrete set.  Note however, that in the
rational case, this set can be infinite, and in general, the set is far
from flat (in a typical family, the size changes on dense subsets!).

\begin{prop}\label{prop:moduli_stack_of_surfaces_is_small}
  The moduli stack of noncommutative rationally ruled surfaces with
  blowdown structures is an algebraic stack locally of finite type over
  $\Z$, and the components classifying $m$-fold blowups of ruled surfaces
  of genus $g$ have dimension $m+3$ if $g=0$, $m+1$ if $g=1$, and $m+2g-2$
  otherwise.
\end{prop}

\begin{proof}
  We first note that if $X/S$ is a family of noncommutative rationally
  ruled surfaces with a blowdown structure, then the splittings of $X$
  compatible with the blowdown structure are essentially classified by a
  smooth scheme over $S$.  More precisely, the splittings are classified by
  $C\times_S \Pic(C)$, with $C$ recording the marked point and $\Pic(C)$
  capturing the fact that the sheaf bimodule is only determined modulo
  $\pi_1^*\Pic(C)$.  If we insist that the Euler characteristic of the
  bimodule be $0$ or $-1$, then this eliminates all but one component.
  Thus it will suffice to construct the corresponding stack for split
  surfaces.
  
  For ruled surfaces, we first note that the curve $C$ is classified by
  the union over $g$ of the classical moduli stacks ${\cal M}_g$ of
  smooth genus $g$ curves with a marked point.

  Relative to the ample divisor $\sO_C(x)\boxtimes \sO_C(x)$ on $C\times
  C$, a sheaf bimodule of rank 2 has Hilbert polynomial $2t+\chi$ for some
  $\chi$.  The condition for a sheaf on $C\times C$ with that Hilbert
  polynomial to have Chern class algebraically equivalent to the double
  diagonal is open and closed, while having no subsheaf supported on a fiber
  of $\pi_1$ or $\pi_2$ is open.  We thus find that the stack ${\cal
    M}_{sb}$ of sheaf bimodules corresponding to ruled surfaces is smoothly
  covered by open substacks of
  \[
  \Quot_{C\times C/{\cal M}_{g,1}}(\sO_C(-dx)\boxtimes \sO_C(-dx),2t+\chi).
  \]

  For a family without commutative fibers, we can recover the anticanonical
  curve as the Quot scheme $\Quot_{C\times C/{\cal M}_{sb}}({\cal E},1)$.
  Over commutative fibers, this is not a curve, but we can resolve this as
  follows.  The condition that ${\cal E}$ be a vector bundle of rank 2 on a
  curve is a closed condition, and thus we may blow the corresponding locus
  up to obtain a new stack ${\cal M}_0$.  The $\Quot$ scheme classifying
  the anticanonical curve is still badly behaved over the exceptional
  locus, but now the bad behavior comes in the form of an extra component
  that we can remove to obtain a family of curves of genus 1.  We thus find
  that ${\cal M}_0$ is the desired moduli stack of noncommutative ruled
  surfaces.  (We easily check that when $g=1$, the anticanonical curve is
  either two disjoint copies of $C$ or a double copy of $C$ with some
  number of fibers, so in either case is not integral.)

  The universal anticanonical curve over ${\cal M}_0$ embeds in the
  commutative ruled surface
  \[
  X':=\Quot_{C/{\cal M}_0}(\pi_{1*}{\cal E},1),
  \]
  which lets us construct the moduli stack in general.  For each $m>0$, ${\cal
    M}_m$ is the universal anticanonical curve over ${\cal M}_{m-1}$, and
  the universal anticanonical curve on ${\cal M}_m$ is obtained by blowing
  up the corresponding point of the universal commutative ruled surface
  and then removing one copy of the exceptional curve from the pullback of
  $Q$.

  To compute the dimension, we first note that each blowup adds $1$ to the
  dimension, so it suffices to consider the ruled surface case, and for
  $g\ge 1$ since the calculation of Theorem \ref{thm:antican_moduli_lci}
  applies for $g=0$ with only minor modifications.  We may also compute the
  dimension for $C$ fixed and add $3g-3$.  For $g\ge 2$, the sheaf bimodule
  is supported on the double diagonal, and is either invertible (a
  $0$-dimensional stack since the double diagonal has trivial dualizing
  sheaf), a torsion-free sheaf determined by where it fails to be
  invertible, or a vector bundle.  The invertible case has dimension $0$,
  but we must the mod out by the action of the Picard stack of $C$
  (twisting by $\rho_1^*{\cal L}$), so the correct dimension is $1-g$,
  giving the total $2g-2$ as stated.

  For the non-invertible cases, pick a point where ${\cal E}$ is not
  invertible (there being only finitely many such points on non-commutative
  geometric fibers).  Then $Q$ has a component lying over that point, and
  thus we obtain a 1-parameter family of blowups in a point of that
  component and not on any other component.  Each of those blowups has a
  unique blowdown (using Theorem \ref{thm:elem_xform_works}) to a surface
  with one fewer point of non-invertibility.  We can recover the original
  surface and its blowup from the resulting $Q'$ along with the point that
  got blown up, and thus we find that the two pieces of the moduli stack
  have the same dimension.

  In the genus 1 case, the nonreduced support case picks up an additional
  dimension (which reduced curve it is supported on) which is then
  cancelled out by the fact that $C$ has infinitesimal automorphisms, and
  we need to mod out by the choices of isomorphism $C_1\cong C_0$.  A
  bimodule with reduced support is supported on the union of two graphs of
  translations ($2$ dimensions), is a line bundle on each ($2-2$ more
  dimensions), modulo a line bundle on $C_1$ ($1-1$ fewer dimensions) and
  automorphisms of $C_1$ ($-1$ dimensions).
\end{proof}

\begin{rems}
  In the rational case, we could also construct this as the relative
  $\Pic^0(Q)$ (the scheme, not the stack!) over the moduli stack of
  commutative rational surfaces (with blowdown structure), which also deals
  with the case of noncommutative planes.  (Note that $\Pic^0(Q)$ is cut
  out from $\Pic(Q)$ by the open conditions $\dim H^1(q^2),\dim
  H^1(q^{-2})\le 1$; if $q$ had negative degree along some component, then
  $q^2$ would have degree $\le -2$ along that component, and the resulting
  quotient would have too large an $H^1$.  In particular, $\Pic^0(Q)$ is
  indeed a flat group scheme over the moduli stack.)  It follows from
  Theorem \ref{thm:antican_moduli_lci} that the stack of noncommutative
  rational surfaces is a local complete intersection.  Something similar is
  true in higher genus, with fibers $\Pic^0(Q)/\Pic^0(C)\cong
  \Pic^0(\Spec_C(\rho_*\sO_Q))/\Pic^0(C)$, using Lemma
  \ref{lem:Pic0Q/Pic0C} to deform the map $Q\to C$ and construct the
  dg-category.  (There is a mild technical issue for $g=1$, in that the
  deformed map $Q\to C$ is only defined over the smooth cover $\Pic^0(Q)$,
  and thus one must first construct the family of dg-categories over the
  cover and then descend.)  In either case, we find that the moduli stack
  of noncommutative surfaces retracts to the substack of anticanonical
  commutative surfaces via a smooth morphism.
\end{rems}

\begin{rems}\label{rem:derived_moduli_stack_of_surfaces}
  One can also compute the expected dimension of this moduli stack, and
  find that it is {\em highly} obstructed for $g\ge 2$.  (That there are
  some obstructions should not be surprising: for $g\ge 2$, the moduli
  stack has nonreduced components!)  This is most easily seen via the
  moduli stack of anticanonical commutative surfaces, for which the
  expected dimension is $6g-6$ (the dimension of the stack of
  $\PGL_2$-bundles on smooth curves of genus $g$) plus $2m$ (a choice of
  $m$ blowups) plus $8-9g-m$ (the expected dimension of the anticanonical
  linear system), for a total of $2-3g+m$.  Since the moduli stack of
  noncommutative surfaces is a $1$-dimensional fiber bundle over the moduli
  stack of anticanonical surfaces, we end up with a total expected
  dimension of $3-3g+m$.  This is correct if $g=0$, but for $g=1$, the true
  dimension is greater by 1, and for $g=2$, the true dimension is greater
  by $5g-5$, so that we must have at least that many obstructions.  This
  suggests (along with the fact that we are already working with
  dg-algebras) that one should consider the {\em derived} moduli stack of
  noncommutative surfaces.  This is reasonably straightforward to
  construct: the construction of Proposition
  \ref{prop:stack_of_all_surfaces} readily generalizes to show that the
  moduli stack of rationally ruled surfaces of genus $g$ is smooth of
  dimension $6g-6+2m$, and thus we obtain a derived analogue of the moduli
  stack of anticanonical surfaces by taking the derived linear system
  $|{-}K_X|$ over this stack.  (This is an open substack of the quotient by
  $\G_m$ of the derived stack of morphisms $\omega_X\to \sO_X$.)  It
  remains only to show that the noncommutative deformations are classified
  by a geometric derived stack.  For planes, this is immediate (the derived
  structure on the stack of anticanonical planes is trivial!), while for
  surfaces ruled over $\P^1$, we may use the isomorphism $\Pic^0(Q)\cong
  \Pic^0(\Spec_{\P^1}(\pi_*\sO_Q))$ to obtain $\Pic^0(Q)$ as the pullback
  from an underived family (the family of anticanonical curves on $F_2$
  disjoint from the $-2$-curve).  In higher genus, we need to construct
  $\Map_\pi(Q,C)\cong \Map_\pi(\bar{Q},C)$.  This is open in
  $\Map(\bar{Q},C)$ (since it is open in the underived case) and thus it
  suffices to show that $\Map(\bar{Q},C)$ is relatively geometric and
  locally of finite presentation.  This ought to be directly constructible
  (since $\bar{Q}$ is still a double cover of $C$, so can be represented by
  a more or less explicit sheaf of dg-algebras on $C$), but in any event
  can be shown to be geometric and locally of finite presentation using
  \cite[Thm. 5.1.1]{Halpern-LeistnerD/PreygelA:2019}.
\end{rems}

In Section \ref{sec:types_comm}, we considered a decomposition of the
moduli stack of rational surfaces based on the structure of the
anticanonical curve (i.e., the divisor classes and multiplicities of the
components along with whether $\Pic^0(Q)$ is elliptic, multiplicative, or
additive), and showed that in each case, the corresponding substack was
irreducible and smooth.  (Indeed, there is an \'etale covering making it a
smooth scheme, in which we specify an ordering within each set of
components of $Q$ with the same divisor class and multiplicity.)  This of
course carries over immediately to the moduli stack of noncommutative
rational surfaces, since $\Pic^0(Q)$ is also irreducible and smooth.

The basic idea also works for noncommutative rationally ruled surfaces of
higher genus.  There is an easy induction reducing to $X_0$: each time we
blow up a point, the new combinatorics either determines the multiplicity
with which each component contained that point or gives a smooth open curve
over which the point varies.  For $X_0$ itself, if $g=1$ and $Q$ is smooth,
then we need simply specify the pair of line bundles and translations, and
otherwise we are in the differential case, so need simply specify the
torsion-free sheaf on the double diagonal.  Such a sheaf is either
invertible, so by the constraint on the Euler characteristic is classified
by $\Pic^0(Q)\cong \G_a$, or strictly torsion-free, in which case it is
uniquely determined by the divisor on which it fails to be torsion-free, or
equivalently the points of $C_0$ (with multiplicities) over which the
components of $Q$ of class of $f$ lie.  In addition, we note that the
degree of the Fitting scheme of the sheaf bimodule is semicontinuous, so
for $g>1$ each subfamily is locally closed, and the same follows for $g=1$
once we observe that the support being reduced is an open condition.

One immediate consequence is that characteristic 0 points are dense in the
moduli stack of noncommutative surfaces; i.e, that any surface over a field
of characteristic $p$ can be lifted (possibly after a separable field
extension) to characteristic $0$, which may even be taken to have the same
type.  Another is that for any family $X/S$, there is an induced
decomposition of $S$ into locally closed subsets, which must be finite if
the base is Noetherian.

A natural further question is how the closures of these substacks interact.
As we saw in the commutative case, there are already examples showing in
the rational case that this is not a stratification in general over $\Z$,
though it is still open whether this phenomenon can occur if one excludes
points of finite characteristic.  Unfortunately, despite some questions
being easier for ruled surfaces with $g>0$, this appears unlikely to be one
of them.

\chapter{Dimension and numerical invariants of sheaves}
\label{chap:dimens}

\section{The numerical Grothendieck group}

In \cite{ChanD/NymanA:2013}, a definition was given for a ``non-commutative
smooth proper $d$-fold'', and it was shown there that ruled surfaces (more
precisely, quasi-ruled surfaces in which the two curves are isomorphic)
satisfy their definition (with $d=2$, of course).  Our objective is to show
that this continues to hold for iterated blowups of arbitrary quasi-ruled
surfaces.  In this section, our focus is on those of their conditions
related to dimensions of sheaves, though in the process we will also be
considering the Grothendieck groups in greater detail.  This will give us
the information required to finish our considerations of birational
geometry; we will then prove the remaining Chan-Nyman axiom (that $\Quot$
schemes behave well) in Chapter \ref{chap:quot} below.

Thus suppose $X=X_m$ is an $m$-fold blowup of a quasi-ruled surface $X_0$
(over an algebraically closed field), with intermediate blowups
$X_1$,\dots,$X_{m-1}$ and associated morphisms $\alpha_1$,\dots,$\alpha_m$.
We wish to show that $X$ satisfies the Chan-Nyman axioms.  (Note that they
implicitly include ``strongly Noetherian'' as an unnumbered axiom, but this
of course follows as a special case of Theorem
\ref{thm:noetherian_base_implies_noetherian}.)

The first two of their axioms are automatic: the Gorenstein condition
simply states that there is a Serre functor such that $S[-2]$ is an
autoequivalence, which follows from the above descriptions of the derived
categories and $t$-structures; this in turn immediately implies that
$\Ext^p$ vanishes for $p\notin \{0,1,2\}$, giving their ``smooth, proper of
dimension $2$'' condition.  Thus the first case that requires some work is
the existence of a suitable dimension function (satisfying their axioms
3,4,6 as well as continuity and finite partitivity).

To define this, we recall that the description of the derived category
immediately gives us a description of the Grothendieck group of $\coh X$:
it is the direct sum of $K_0(C_0)$, $K_0(C_1)$, and a copy of $\Z$ for each
point being blown up.  We define the {\em rank} of a class in the
Grothendieck group to be the rank of the image in $K_0(C_0)$ minus the rank
of the image in $K_0(C_1)$.  This then immediately defines a rank for any
coherent sheaf (or object in $D^b_{\coh} X$) as the rank of its class in
the Grothendieck group.

\begin{lem}
  If $m>0$, then for any class $[N]\in K_0(X_{m-1})$,
  $\rank(\alpha_m^*[N])=\rank([N])$, while for any $[M]\in K_0(X)$,
  $\rank(\alpha_{m*}[M])=\rank([M])$.
\end{lem}

\begin{proof}
  The first claim is immediate from the definition, while the second claim
  follows from the first together with the fact that
  $\alpha_m^*\alpha_{m*}[M]-[M]$ is a multiple of the rank 0 class
  $[\sO_{e_m}(-1)]$.
\end{proof}

Applying the construction of Proposition \ref{prop:nice_divisor_on_blowup}
inductively gives an embedding of a commutative curve $Q$ as a divisor,
with morphism $\iota:Q\to X$.  Note that there is a natural map $L\iota^*$
from $D^b_{\coh} X$ to the derived category $\perf(Q)$ of perfect complexes
on $Q$, and thus an induced map from $K_0(X)$ to $K_0^{\perf}(Q)$.  The
latter of course also has a well-defined rank function, which is almost
everywhere on $Q$ given by the alternating sums of the ranks of the
cohomology sheaves.  In particular, we may compute the rank by considering
the restriction to $\hat{Q}$ (the scheme-theoretic union of component(s) of
$Q$ which are not components of fibers).

\begin{lem}
  For any class $[M]\in K_0(X)$, $\rank([M])=\rank(\iota^*[M])$.
\end{lem}

\begin{proof}
  It suffices to check this on each summand of $K_0(X)$.  Each of the
  summands other than $K_0(C_0)$ and $K_0(C_1)$ is generated by a sheaf
  $\sO_{e_i}(-1)$ of rank 0 that meets $\hat{Q}$ transversely, and thus has
    rank 0 restriction as required.  Similarly, $K_0(C_0)$ and $K_0(C_1)$
    are generated by line bundles, and in either case the restriction to
    $\hat{Q}$ is again a line bundle, so has rank 1 as required.
\end{proof}

\begin{cor}
  For any sheaf $M\in \coh Q$, $\rank(\iota_*M)=0$.
\end{cor}

\begin{proof}
  We have $\rank(\iota_*M)=\rank(L\iota^*\iota_*M)$,
  and the cohomology sheaves of $L\iota^*\iota_*M$ are $M$ and the twist of
  $M$ by a suitable invertible sheaf, so almost everywhere have the same
  rank.
\end{proof}

\begin{cor}
  For any class $[M]\in K_0(X)$, $\rank(\theta[M])=\rank([M])$.
\end{cor}

\begin{proof}
  We have $[M]-\theta[M]=\iota^*\iota_*[M]$.
\end{proof}

\begin{prop}
  For any object $M\in \coh X$, $\rank(M)\ge 0$.
\end{prop}

\begin{proof}
  Since $\alpha_{m*}$ preserves the rank, we may reduce to the case $m=0$,
  so that $X=X_0$ is a quasi-ruled surface.  Consider the function
  $r(l,[M]):=\rank \rho_{l*}[M]$ on $\Z\times K_0(X_0)$.  This is clearly
  linear in $[M]$ (since rank and the action of $\rho_{l*}$ on the
  Grothendieck group are both linear), while in $l$ it satisfies (as an
  immediate consequence of Lemma \ref{lem:exact_tri_quasi_ruled})
  \[
  r(l+1,[M])+r(l-1,[M])=2r(l,[M]),
  \]
  so that $r(l,[M])-r(0,[M])$ is bilinear.   Since
  \[
  \rank([M])=r(0,[M])-r(-1,[M]),
  \]
  we find that $r(l,[M]) = r(0,[M])+l\rank([M])$.
  For $M$ a sheaf, $\rho_{0*}\theta^{-l}M$ is a sheaf for $l\gg 0$, and
  thus
  $\rank(\rho_{0*}\theta^{-l}M)=r(2l,M)\ge 0$
  for $l\gg 0$, implying $\rank(M)\ge 0$ as required.  (Here we use the
  fact that $R\rho_{0*}\theta^{-l}M$ differs from $R\rho_{(2l)*}M$ via a
  twist by a line bundle.)
\end{proof}

The rank of course only depends on the {\em algebraic} class in $K_0(X)$,
or in other words on the quotient $K_0^{\num}(X)$ of $K_0(X)$ by the
identity component $\Pic^0(C_0)\oplus \Pic^0(C_1)$ of $K_0(C_0)\oplus
K_0(C_1)$.  (The notation is justified by the fact that numeric and
algebraic equivalence turn out to agree in this case, as we will see.)  Let
$[\pt]\in K_0^{\num}(X)$ denote the class of a point of $Q$; this is
well-defined since the corresponding map $Q\to K_0(X)$ is a map from a
connected scheme, so lies inside a single coset of the identity component.
This can also be described in terms of the semiorthogonal decomposition: it
is the class of a point in $K_0^{\num}(C_0)$ minus the class of a point in
$K_0^{\num}(C_1)$, and is 0 in every other component.  This is dual to the
definition of rank, and indeed one finds that $\rank([M])=\chi([M],[\pt])$,
where $\chi$ denotes the Mukai pairing (which on a pair of complexes is the
alternating sum of $\Ext$ dimensions, and is bilinear on $K_0(X)$).  We of
course have $\rank([\pt])=0$, and thus we may use $\rank$ and $[\pt]$ to
define a filtration of the Grothendieck group, and in turn define the
dimension of a coherent sheaf.

\begin{defn}
  Given a nonzero sheaf $M\in \coh(X)$, the dimension $\dim(M)$ of $M$ is
  defined to be $2$ if $\rank(M)>0$, $0$ if $[M]\propto [\pt]$, and $1$
  otherwise.  We say that $M$ is {\em pure} $d$-dimensional if $\dim(M)=d$
  while any nonzero proper subsheaf of $M$ has dimension $<d$.  A sheaf is
  {\em torsion-free} if it is either $0$ or pure $2$-dimensional.
\end{defn}

To show that this is truly a dimension function, we need to show that a
subsheaf or quotient of a $d$-dimensional sheaf has dimension at most $d$.
(Chan and Nyman also ask for compatibility with the Serre functor, but this
follows immediately from $\rank(\theta[M])=\rank([M])$.)  For $d=2$ this
is trivial, while for $d=1$ it reduces to the fact that the rank is
additive and nonnegative.  For $d=0$, this will require some additional
ideas.  We first note that the nonnegativity of the rank has an analogue
for $0$-dimensional sheaves.

\begin{prop}
  If $M\in \coh X$ is a sheaf of numeric class $d[\pt]$, then $d\ge 0$, with
  equality iff $M=0$.
\end{prop}

\begin{proof}
  If $m>0$, then $\alpha_{m*}\theta^{-l}M$ is a sheaf for $l\gg 0$, and we
  readily verify that it also has class $d[\pt]$.  Thus by induction $d\ge
  0$, and if $d=0$, then $\alpha_{m*}\theta^{-l}M=0$ for $l\gg 0$; since
  $\alpha_m^*\alpha_{m*}\theta^{-l}M\to \theta^{-l}M$ is surjective for
  $l\gg 0$, we conclude that $M=0$ as required.  Similarly, for $m=0$,
  $\rho_{l*}M$ has class $d[\pt]$ and is a sheaf for $l\gg 0$, so that we
  reduce to the corresponding claims in $C_0$ and $C_1$.
\end{proof}

\begin{cor}
  If $M$, $N$ are coherent sheaves with $[M]=[N]$ in $K_0^{\num}(X)$, then
  any injective or surjective map $M\to N$ is an isomorphism.
\end{cor}

\begin{proof}
  The quotient resp. kernel would be a sheaf with trivial class in
  $K_0^{\num}(X)$, and thus 0 by the Proposition.
\end{proof}

\begin{defn}
  The {\em N\'eron-Severi lattice} $\NS(X)$ is the subquotient
  $\ker(\rank)/\Z[\pt]$ of $K_0^{\num}(X)$.  Given two classes $D_1,D_2\in
  \NS(X)$, their {\em intersection number} $D_1\cdot D_2$ is given by
  $-\chi(D_1,D_2)$.  The {\em canonical class} $K_X\in \NS(X)$ is the class
  $[\theta \sO_X]-[\sO_X]$
\end{defn}

\begin{rem}
  Since $\chi(D,[\pt])=\chi([\pt],D)=0$ iff $\rank(D)=0$, this pairing is
  indeed well-defined.
\end{rem}

By mild abuse of notation, we will refer to the classes in $\NS(X)$ as
``divisor'' classes.  This might more properly be reserved for the
extension of $\NS(X)$ obtained as the quotient of the subgroup of
$1$-dimensional classes in $K_0(X)$ by the subgroup of classes of
$0$-dimensional sheaves, but this quotient is badly behaved when $X$ is not
a maximal order, since the subgroup of such classes need not be closed.

\begin{defn}
  The (numeric first) {\em Chern class} $c_1:K_0(X)\to \NS(X)$ is defined
  by $c_1([M])=[M]-\rank([M])[\sO_X]$.
\end{defn}

In particular, a $1$-dimensional sheaf is $0$-dimensional iff its Chern
class is 0.

\begin{prop}
  The N\'eron-Severi lattice of a noncommutative rationally quasi-ruled
  surface has rank $m+2$, and the intersection pairing is a nondegenerate
  symmetric bilinear form of signature $(+,-,-,\dots,-)$.
\end{prop}

\begin{proof}
  The semiorthogonal decomposition makes it straightforward to determine
  the Fourier-Mukai form on $K_0^{\num}(X_m)$ from that of
  $K_0^{\num}(X_{m-1})$, and we in particular find that the induced form on
  $\NS(X_m)$ satisfies
  \[
  (\alpha_m^*D_1+r_1 e_m)\cdot (\alpha_m^*D_2+r_2 e_m)
  =
  (D_1\cdot D_2) -r_1r_2,
  \]
  so that we may reduce to the case $m=0$.  (Here we have defined $e_m$ to
  be the Chern class of the exceptional sheaf $\sO_{e_m}(-1)$.)

  For $m=0$, we have $K_0^{\num}(C_i)\cong \Z^2$ (determined by rank and
  degree), and thus $K_0^{\num}(X_0)\cong \Z^4$ and $\NS(X_0)\cong \Z^2$.
  We may then define $f$ to be the (numeric) class of a point of $C_0$
  (noting that the numeric class of a point in $C_1$ gives the same class
  in $\NS(X)$, since they differ by $[\pt]$), and let $s$ be any class
  which is rank $1$ in both $C_0$ and $C_1$.  We then readily compute that
  $s\cdot f=f\cdot s=1$, giving the required symmetry.  Moreover, we have
  $f^2=0$, making the signature $(+,-)$ as required.
\end{proof}

We will often refer to the basis $s,f,e_1,\dots,e_m$ of $\NS(X)$
arising from this proof.  The class $s$ is only defined modulo $f$, but
since $(s+df)^2=s^2+2d$, we can fix this by insisting that $s^2\in
\{-1,0\}$ to obtain a canonical basis of $\NS(X_m)$.  We should
caution that this basis is only canonical relative to the given
representation of $X_m$ as an iterated blowup of a quasi-ruled surface;
e.g., if $X_0$ is a sufficiently nondegenerate deformation of $\P^1\times
\P^1$, then it can be represented in two ways as a ruled surface, and the
two resulting bases differ by interchanging $s$ and $f$.

\begin{cor}
  The radical of the Fourier-Mukai pairing is $\Pic^0(C_0)\times
  \Pic^0(C_1)$, and thus numerically equivalent classes in $K_0(X)$ are
  algebraically equivalent.
\end{cor}

\begin{proof}
  Since the Fourier-Mukai pairing is constant in flat families, the
  identity component of $K_0(X)$ is certainly contained in the radical, and
  we have just shown that the induced pairing on $[\pt]^\perp/[\pt]$ is
  nondegenerate, so that the radical of the induced pairing on
  $K_0^{\num}(X)$ is contained in $\Z[\pt]$.  Since
  $\chi([\sO_X],[\pt])=1$, it follows that the pairing on $K_0^{\num}(X)$
  is nondegenerate.
\end{proof}

We have the following by an easy induction, where we recall the notion of
line bundle from Definition \ref{defn:line_bundle}.

\begin{cor}
  For any divisor class $D$, there is a line bundle $L$ with $c_1(L)=D$,
  and any two line bundles with the same Chern class have the same class in
  $K_0^{\num}(X)$.
\end{cor}

\begin{rem}
  The reader should bear in mind that line bundles do not form a group in
  any reasonable sense!  Indeed, the map from the set of isomorphism
  classes of line bundles to $\NS(X)$ (which, of course, {\em is} a
  group) may not even have constant fibers: half of the fibers are
  $\Pic^0(C_0)$-torsors, and half are $\Pic^0(C_1)$-torsors, depending on
  the parity of $D\cdot f$.
\end{rem}

Since $\sO_X$ has Chern class $0$, and $[\pt]$ can be expressed as a linear
combination of classes of line bundles, we conclude the following.

\begin{cor}
  $K_0^{\num}(X)$ is generated by classes of line bundles.
\end{cor}

For a split family of noncommutative surfaces, the above considerations
imply that the family $K_0^{\num}(X_s)$ of abelian groups is locally
constant as $s$ varies, letting us make the following statement.

\begin{cor}
  If $X/S$ is a split family of noncommutative surfaces, then for any $M\in
  \perf(X)$, the class of $M|^{\bf L}_s$ in $K_0^{\num}(X_s)$ is locally
  constant on $S$.
\end{cor}

\begin{proof}
  The class in $K_0^{\num}(X_s)$ of an object $M$ is uniquely determined by
  the linear functional $\chi({-},M)$ coming from the Fourier-Mukai pairing,
  and thus on the values $\chi(L,M)$.  But $\chi(L|^{\bf L}_s,M|^{\bf
    L}_s)$ is locally constant since $R\Hom(L,M)$ is perfect.
\end{proof}

This implies the ``no shrunken flat deformations'' axiom of
\cite{ChanD/NymanA:2013}: if $M$ and $N$ are distinct fibers of a flat family of
coherent sheaves over a connected base, then any injective or surjective
morphism between them is an isomorphism, since the cokernel or kernel is
then a sheaf with trivial class in $K_0^{\num}(X)$, and thus 0.


\begin{cor}
  A class in $K_0^{\num}(X)$ is uniquely determined by its rank, Chern
  class and Euler characteristic $\chi(M):=\chi(\sO_X,M)$.  The
  Fourier-Mukai pairing is given in these terms by
  \begin{align}
    \chi(M,N) 
    = {}&-\rank(M)\rank(N)\chi(\sO_X)\notag\\
  &+\rank(M)\chi(N)+\rank(N)\chi(M)\notag\\
  &-c_1(M)\cdot (c_1(N)-\rank(N)K_X).
  \end{align}
\end{cor}

\begin{proof}
  We first observe that Serre duality implies
  \[
  \chi([\sO_X],\theta[\sO_X])
  =
  \chi([\sO_X],[\sO_X])
  \]
  and
  \begin{align}
  \chi(M,[\sO_X])
  &=
  \chi(M,[\sO_X]-\theta[\sO_X])
  +
  \chi(M,\theta[\sO_X])\notag\\
  &=
  \chi(M-\rank(M)[\sO_X],[\sO_X]-\theta[\sO_X])
  +
  \chi([\sO_X],M)\notag\\
  &=
  c_1(M)\cdot K_X + \chi(M),
  \end{align}
  while the definition of the intersection pairing gives
  \[
  -c_1(M)\cdot c_1(N) = \chi(M-\rank(M)[\sO_X],N-\rank(N)[\sO_X]).
  \]
  Expanding the right-hand side via bilinearity and simplifying gives the
  desired result.
\end{proof}

This leads to the following definition, as a convenient way to represent
classes in $K_0(X)$.

\begin{defn}\label{defn:num_invs}
  The {\em numerical invariants} of a class $[M]\in K_0(X)$ are its rank,
  Chern class and Euler characteristic, generally given as the triple
  $(\rank([M]),c_1([M]),\chi([M]))$.
\end{defn}

In the above expression, $\chi(\sO_X)$ is easy to compute, as $\sO_X$ lives
entirely inside the component $D^b_{\coh} C_0$ of the semiorthogonal
decomposition; we thus conclude that
$\chi(\sO_X)=\chi(\sO_{C_0})=1-g(C_0)$.  More generally, if $L$ is a line
bundle, one has $\chi(L,L)=1-g(C_{c_1(L)\cdot f})$, from which one can
solve for $\chi(L)$:
  \[
  \chi(L) = 1 - \frac{g(C_0)+g(C_{c_1(L)\cdot f})}{2} + \frac{c_1(L)\cdot
    (c_1(L)-K_X)}{2}.
  \]
  This, of course, agrees with the standard Riemann-Roch formula for the
  Euler characteristic in the case of a ruled surface (as it must: $\chi(L)$ is
  locally constant as we vary the surface, so for ruled surfaces may be
  computed on the commutative fiber).

For $K_X$, the situation is slightly more complicated, but we have the
following.

\begin{prop}
  In terms of the standard basis of $\NS(X)$, we have
  \[
  K_X = \begin{cases}
    -2s-(2-g(C_0)-g(C_1))f+e_1+\cdots+e_m, & s^2 = 0\\
    -2s-(3-g(C_0)-g(C_1))f+e_1+\cdots+e_m, & s^2 = -1
  \end{cases}.
  \]
\end{prop}

\begin{proof}
  We readily reduce to the case $m=0$, where the description of the action
  of $\theta$ tells us that
  \[
    [\theta\sO_X]+[\rho_0^*\omega_{C_0}] = \rho_1^*V
  \]
  for some rank 2 vector bundle $V$ on $C_1$.  It follows that $K_X =
  -2s+df$ for some $d$, and since $K_X^2$ is linear in $d$, we reduce to
  showing $K_X^2 = 4(2-g(C_0)-g(C_1))$.  Since $\chi(\rho_1^*V,\rho_1^*V)
  =\chi(\End(V)) = 4(1-g(C_1))$ by Riemann-Roch, we find
  \[
  \chi([\theta\sO_X]+[\rho_0^*\omega_{C_0}],[\theta\sO_X]+[\rho_0^*\omega_{C_0}])
  =
  4-4g(C_1).
  \]
  The left-hand side expands as a sum of 4 terms, three of which simplify
  via Serre duality to a calculation inside $D^b_{\coh} C_0$ and the fourth of
  which can be simplified using the general expression for $\chi(M,N)$.  We
  thus find
  \[
  4-4g(C_1)
  =
  K_X^2 - c_1(\rho_0^*\omega_{C_0})\cdot K_X.
  \]
  Since $c_1(\rho_0^*\omega_{C_0}) = (2g(C_0)-2)f$, this gives $K_X^2 =
  4(2-g(C_0)-g(C_1))$ as required.
\end{proof}

\begin{rem}
  We similarly find $[\sO_Q]=-K_X + (g(Q)-1)f$.
\end{rem}

\begin{cor}
  The action of $\theta$ on $K_0^{\num}(X)$ is given by
  \begin{align}
    \rank(\theta M) &= \rank(M)\\
    c_1(\theta M) &= c_1(M)+\rank(M)K_X\\
    \chi(\theta M) &= \chi(M) + c_1(M)\cdot K_X.
  \end{align}
\end{cor}

\begin{proof}
  For each $M$, $\chi(M,N)=\chi(N,\theta M)$ by Serre duality.  Both sides
  are linear functionals on $K_0^{\num}(X)$, and comparing coefficients
  gives the desired expressions.
\end{proof}

The action of the duality $R\ad$ is also straightforward to compute, since
we know how it acts on line bundles.

\begin{prop}
  The action of $R\ad$ on $K_0^{\num}(X)$ is given by
  \begin{align}
    \rank(R\ad M) &= \rank(M)\\
    c_1(R\ad M) &= -c_1(M)+\rank(M)K_X\\
    \chi(R\ad M) &= \chi(M).
  \end{align}
\end{prop}  

\begin{proof}
  Indeed, this is linear and acts correctly on line bundles.
\end{proof}

\begin{prop}
  One has $\rank(\rho_{d*}M) = (d+1)\rank(M)+c_1(M)\cdot f$.
\end{prop}

\begin{proof}
  We have $\rank(\rho_{d*}M) =
  \chi([\sO_X(-ds-(l+1)f)],M)-\chi([\sO_X(-ds-lf)],M)$ for some $l$, and can
  compute the numerical invariants of the two line bundles; in particular,
  we see that the result does not depend on $l$.
\end{proof}

\begin{rem}
  For $\chi(\rho_{d*}M)$, the calculation is more delicate, and in
  particular depends on the precise normalization of the sheaf bimodule;
  we have $[\rho_d^*\sO_{C_d}]=[\sO_X(-ds-lf)]$ for some $l$, but now $l$
  matters, though only when $\rank(\rho_{d*}M)\ne 0$.
\end{rem}

For simplicity, we only give the next result in the ruled case.

\begin{prop}
  Let $X$ be a rationally ruled surface, and let $L$ be a line bundle
  from $X$ by twisting by a line bundle, so that one has a functor
  ${-}(L):\coh X\to \coh X'$ taking $\sO_X$ to $\sO_{X'}$.  Then
  \begin{align}
    \rank(M(L)) &= \rank(M)\\
    c1(M(L)) &= c_1(M) + \rank(M)c_1(L)\\
    \chi(M(L)) &= \chi(M) + c_1(M)\cdot c_1(L) + \rank(M)(c_1(L)\cdot
    (c_1(L)-K_{X'})/2
  \end{align}
\end{prop}

When $X$ is a maximal order, we are also interested in how the direct image
and pullback functors relate the two Grothendieck groups.  The nicest
answer regards the N\'eron-Severi groups.

\begin{prop}\label{prop:NS_of_center}
  Let $X$ be a rationally quasi-ruled surface, and suppose that $X\cong
  \Spec{\cal A}$ where ${\cal A}$ is a maximal order of rank $r^2$ over the
  commutative surface $Z$, with associated morphism $\pi:X\to Z$.  Then $X$
  and $Z$ have the same parity, and in terms of the standard bases of
  $\NS(X)$ and $\NS(Z)$, the maps $\pi^*$ and $\pi_*$ are both
  multiplication by $r$.
\end{prop}

\begin{proof}
  Since $\pi^*\pi_*M\cong M\otimes_{\sO_Z} {\cal A}$ and ${\cal A}$ is a
  vector bundle of rank $r^2$, we see that $c_1(\pi^*\pi_*M)=r^2 c_1(M)$.
  Adjunction then gives $\pi_*D_1\cdot \pi_*D_2 = r^2 D_1\cdot D_2$, and
  the same holds for $\pi^*$ since their product is multiplication by $r^2$.

  We next address the question of parities, which is of course really about
  quasi-ruled surfaces.  An elementary transformation flips both parities,
  and thus we may assume that we are in the untwisted case with
  $\hat{Q}=\bar{Q}$.  As we discussed when considering point sheaves, in
  this case there is a natural morphism $\rho_1^*\sO_{C_1}\to
  \rho_0^*\sO_{C_0}$ the cokernel of which is a sheaf of class in $s+\Z f$
  disjoint from $Q$.  The direct image of this sheaf is a vector bundle of
  rank $r$ on its support, a section of $Z$.  Since we have shown that the
  direct is a similitude relative to the intersection form, we conclude
  that the classes in $s+\Z f$ on the respective surfaces have the same
  self-intersection, and thus the surfaces have the same parity, letting us
  identify their N\'eron-Severi lattices via the standard bases.

  It remains to show that $\pi_*$ and $\pi^*$ are multiplication by $r$
  relative to these bases.  For $m\ge 1$, we note that $\pi_*\sO_{e_m}(-1)$
  is supported on the corresponding exceptional curve of $Z$, and thus has
  Chern class a multiple of $e_m$, which must be $r e_m$ by the isometry
  property.  Similarly, for any point $x\in C'$, we have
  $\pi^*\rho^*\sO_x\cong \rho_0^* \phi^*\sO_x$, where $\phi$ is the
  projection $C_0\to C'$, and thus $\pi^*(f) = rf$.  Finally, $\pi_*s$ is
  orthogonal to each $e_i$ and has intersection $r$ with $f$, so has the
  form $\pi_*s = rs+df$, with $d=0$ then forced by $(\pi_*s)^2=r^2$.
\end{proof}

In fact, we can give the full map between the Grothendieck groups, though
this is more complicated.

\begin{prop}
  If the rationally quasi-ruled surface $X$ is a maximal order over $Z$,
  then using the standard bases to identify their N\'eron-Severi lattices,
  one has
  \begin{align}
  \rank(\pi_*M) &= r^2\rank(M)\notag\\
  c_1(\pi_*M) &= r c_1(M)+\rank(M) \frac{r^2 K_Z - r K_X}{2}\notag\\
  \chi(\pi_*M) &= \chi(M)
  \end{align}
  and
  \begin{align}
  \rank(\pi^*M) &= \rank(M)\notag\\
  c_1(\pi^*M) &= r c_1(M)\notag\\
  \chi(\pi^*M) &= r^2\chi(M)+\frac{c_1(M)\cdot (r^2K_Z-rK_X)}{2}+\rank(M)(\chi(\sO_X)-r^2\chi(\sO_Z)).
  \end{align}
\end{prop}

\begin{proof}
  Since $\pi^*\sO_Z\cong \sO_X$, we have $\chi(\pi_*M)=\chi(M)$, and since
  $\pi_*\sO_X={\cal A}$, $\rank(\pi_*\sO_X) = r^2$.  Moreover, our Chern
  class computation tells us how $\pi_*$ acts on the Chern class of
  $1$-dimensional sheaves.  This almost determines how $\pi_*$ acts on a
  basis of $K_0^{\num}(X)$, except that we still need to determine
  $c_1(\pi_*\sO_X)$.  Since $\pi_*\theta \sO_X\cong \sHom(\pi_*\sO_X,\omega_Z)$,
  we find that
  \[
  c_1(\pi_*\theta \sO_X) = r^2 K_Z - c_1(\pi_*\sO_X)
  \]
  and thus
  \[
  2c_1(\pi_*\sO_X) = r^2 K_Z - r K_X.
  \]
  The action of $\pi^*$ follows by adjunction.
\end{proof}  

\section{Effective divisor classes}

In the commutative setting, the notion of an effective divisor is of course
quite crucial.  In the commutative setting, a divisor on a smooth surface
is effective if it can be represented by a curve; although curves
themselves do not make sense in the noncommutative case (apart from
components of $Q$), we are still led to the following definition.

\begin{defn}
  A divisor class is {\em effective} if it is the Chern class of a
  $1$-dimensional sheaf.
\end{defn}

\begin{rems}
Note that the effective classes form a monoid, since we can always take the
direct sum of the corresponding $1$-dimensional sheaves.
\end{rems}

\begin{rems}
  In the commutative case, one often uses the equivalence between a divisor
  being effective and the corresponding line bundle having a global
  section.  As stated, this does not hold in the noncommutative setting;
  there are cases (e.g., $e_1-e_2$ when one blew up two points in the same
  orbit) of effective divisors such that no line bundle with the given
  Chern class has a global section.  The situation is somewhat better if
  one allows both line bundles to vary, fixing the difference of Chern
  classes, but even then it is unclear whether the resulting classes form a
  monoid.  (The monoid they generate is, however, correct, at least for
  rational or rationally ruled surfaces.)
\end{rems}

\begin{rems}
  This should not be confused with the notion of effective divisor from
  \cite{JorgensenP:2000}.
\end{rems}

If a $0$-dimensional sheaf had a $1$-dimensional subsheaf, then both the
subsheaf and the quotient would have nonzero effective Chern classes, and
thus to rule this out, we need to show that the effective monoid intersects
its antipode only in 0.

\begin{prop}\label{prop:neg_of_eff_not_eff}
  If $D\in \NS(X)$ is such that both $D$ and $-D$ are effective,
  then $D=0$.
\end{prop}

\begin{proof}
  Let $M$ be a $1$-dimensional sheaf with Chern class $D$, and let $N$ be a
  $1$-dimensional sheaf with Chern class $-D$.  If $m>0$, then for $l\gg
  0$, $\alpha_{m*}\theta^{-l}M$ and $\alpha_{m*}\theta^{-l}N$ are both
  sheaves, and their Chern classes add to 0 since $\alpha_{m*}$ induces a
  well-defined homomorphism $\NS(X_m)\to \NS(X_{m-1})$ with
  kernel $\Z e_m$.  Thus by induction, we have $D=d e_m$ for some $d\in
  \Z$.  But then $[\alpha_{m*}\theta^{-l}M]$ and
  $[\alpha_{m*}\theta^{-l}N]$ are multiples of $[\pt]$ for all $l$, with the
  coefficient depending linearly on $l$; since they must be nonnegative for
  $l\gg 0$, both linear terms must be nonnegative, implying that $d,-d\ge
  0$ as required.

  For $m=0$, suppose $D=ds+d'f$ (relative to the standard basis discussed
  above).  Then $\rank(\rho_{l*}M)=d=-\rank(\rho_{l*}N)$, and since both
  are sheaves for $l\gg 0$, we must have $d=0$.  But then $[\rho_{l*}M]$
  and $[\rho_{l*}N]$ are proportional to the class of a point for all $l$,
  again depending linearly on $l$, with coefficients $\pm d'$, so that
  $d'=0$ as well.
\end{proof}

It follows immediately that our dimension function is {\em exact} in the
sense of \cite{ChanD/NymanA:2013}: given an exact sequence
\[
0\to M'\to M\to M''\to 0,
\]
we have $\dim(M)=\max(\dim(M'),\dim(M''))$.  (We have also implicitly shown
that the shifted Serre functor preserves dimension.)

\section{Irreducibility of surfaces}

The next two results essentially say that $X$ is irreducible.

\begin{lem}
  If $L_1$, $L_2$, $L_3$ are line bundles on $X_m$ and $\phi_1\in
  \Hom(L_1,L_2)$, $\phi_2\in \Hom(L_2,L_3)$ are morphisms such that
  $\phi_2\circ \phi_1=0$, then $\phi_1=0$ or $\phi_2=0$.
\end{lem}

\begin{proof}
  For $m>0$, we may write $L_i = \theta^{l_i}\alpha_m^* \theta^{-l_i}L'_i$
  for line bundles on $X_m$, and we find that there are
  composition-respecting injections
  \[
  \Hom(L_1,L_2)\subset \Hom(L'_1,L'_2)\qquad\text{and}\qquad
  \Hom(L_2,L_3)\subset \Hom(L'_2,L'_3).
  \]
  (Indeed, by the construction of the blowup, $\Hom(L_1,L_2)$ is
  essentially defined to be the subspace of $\Hom(L'_1,L'_2)$ satisfying an
  appropriate condition of the form ``vanishes to multiplicity $l_1-l_2$'')
  We thus reduce by induction to the case $m=0$, where it is simply the
  fact that the $\Z$-algebra $\bar{\cal S}$ corresponding to a quasi-ruled
  surface is a domain.
\end{proof}

\begin{prop}
  Any nonzero morphism between line bundles is injective.
\end{prop}

\begin{proof}
  Suppose otherwise, so that $\phi_2:L_2\to L_3$ is a morphism of line
  bundles with nonzero kernel.  Then there is a line bundle $L_1$ such that
  $L_1\otimes \Hom(L_1,\ker(\phi_2))\to \ker(\phi_2)$ is surjective, and
  thus $\Hom(L_1,\ker(\phi_2))\ne 0$.  A nonzero homomorphism $L_1\to
  \ker\phi_2$ induces a nonzero homomorphism $\phi_1:L_1\to L_2$ such that
  $\phi_2\circ \phi_1=0$, giving a contradiction.
\end{proof}

\begin{cor}
  Line bundles are torsion-free.
\end{cor}

\begin{proof}
  If $M\subset L$ is a rank 0 subsheaf of a line bundle $L$, then as
  before, there exists a nonzero morphism $L'\to M$ for some line bundle
  $L'$ and thus the composition $L'\to L$ has image a nonzero subsheaf of
  $M$.  But this composition is nonzero, so injective, and thus its image
  has rank 1, giving a contradiction.
\end{proof}

Serre duality immediately gives the following.

\begin{cor}
  If $L$ is a line bundle and $M$ a $\le 1$-dimensional sheaf, then
  $\Ext^2(L,M)=0$.
\end{cor}

We would like more generally to know that the ``cohomology'' of a coherent
sheaf $M$ (i.e., $R\Gamma(M):=R\Hom(\sO_X,M)$) vanishes in degree
$>\dim(M)$.  For $2$-dimensional sheaves, this follows trivially from the
global bound on $\Ext$ groups, while for $1$-dimensional sheaves, this is a
special case of the corollary.  For $0$-dimensional sheaves, it will be
handy to use the following dichotomy.  We note that although the
restriction map $K_0(X)\to K_0^{\perf}(Q)$ is not well-defined on
$K_0^{\num}(X)$, the rank is certainly well-defined, as is the class of the
determinant in $\Pic(Q)/\Pic^0(C_0)\times \Pic^0(C_1)$.  In particular,
there is a well-defined class $q$ in the latter group obtained as the
determinant of the restriction of $[\pt]$.  Note that if $q^r\not\sim
\sO_Q$, then a sheaf of class $r[\pt]$ has nontrivial restriction to $Q$;
in particular, if $q\not\sim \sO_Q$, then any sheaf of class $[\pt]$ is a
point of $Q$.

\begin{prop}\label{prop:order_iff_q_torsion}
  The class $q$ is torsion iff $X$ is a maximal order on a smooth
  commutative projective surface $Z$.
\end{prop}

\begin{proof}
  If $X$ is a maximal order in a central simple algebra of degree $r$, then
  there are $0$-dimensional sheaves on $X$ of class $r[\pt]$ which are
  disjoint from $Q$, and thus the determinant of their restriction is
  trivial, implying $rq\sim \sO_Q$.

  Now, suppose that $q$ is torsion.  The conclusion is inherited under
  blowing up, so we may assume $m=0$, and since the conclusion is automatic
  for surfaces which are quasi-ruled but not ruled, we may assume $X=X_0$
  is a ruled surface, and not of the hybrid type (since the conclusion also
  holds there) or commutative.  In the differential case, we have
  $\Pic^0(Q)/\Pic^0(C_0)\times \Pic^0(C_1)\cong \G_a$ and find that $q=1$,
  which is torsion only in finite characteristic, where the conclusion
  holds.  In the difference cases, we can express $q$ in terms of the
  action of $s_0s_1$, and find that the latter is torsion iff $q$ is
  torsion.
\end{proof}

\begin{rem}
  The same calculation shows that for a rationally ruled surface, $q$ is
  trivial iff $X$ is commutative.
\end{rem}

\begin{prop}
  If $M$ is $0$-dimensional of class $d[\pt]$, then for any line bundle
  $L$, $\Ext^p(L,M)=0$ for $p>0$ and $\dim\Hom(L,M)=d$.
\end{prop}

\begin{proof}
  If $X$ is a maximal order, then the claim follows immediately from the
  corresponding fact on $Z$, so we suppose that $q$ is not torsion.

  Since $\Ext^2(L,M)=0$, we find that
  $\dim\Hom(L,M)-\dim\Ext^1(L,M)=\chi(L,M)=d$, and thus any line bundle has
  a nonzero morphism to $M$.  The image and cokernel (if nonzero) of such a
  morphism will again be $0$-dimensional, and thus vanishing of
  $\Ext^1(L,M)$ will follow by induction in $d$.  We may thus assume that
  $M$ has no proper sub- or quotient sheaf.  Since $\det(M|^{\dL}_Q)\ne 0$, we
  conclude that $M|^{\dL}_Q\ne 0$, and thus that $M$ is actually supported on
  $Q$.  In other words, we may write $M=\iota_* N$ for some $0$-dimensional
  sheaf $N$ on $Q$, and find
  \[
  R\Hom(L,M)\cong R\Hom(\iota^*L,N),
  \]
  and vanishing follows since $\iota^*L$ is an invertible sheaf on $Q$.
\end{proof}

\begin{rem}
  By Serre duality, it follows that $\Ext^1(M,L)=0$.
\end{rem}

\begin{cor}
  If $M$ is $\le d$-dimensional, then $R^p\ad M=0$ for $p>2$ and $p<2-d$.
\end{cor}

\begin{proof}
  We have already shown that $\ad$ is left exact of homological dimension
  $2$.  Since this is bounded, there is a line bundle $L$ on $X^{\ad}$ such
  that each cohomology sheaf is acyclic and globally generated relative to
  $L$.  On the other hand, we have
  \[
  \Hom(L,R^p\ad M)
  \cong
  h^p(R\Hom(L,R\ad M))\cong \Ext^p(M,\ad L)\cong \Ext^{2-p}(\ad L,\theta M)^*
  \]
  and thus $R^p\ad M$ vanishes if $2-p>\dim(\theta M)=\dim(M)$.
\end{proof}

\begin{prop}
  For any coherent sheaf $M$, $\ad M$ is torsion-free.
\end{prop}

\begin{proof}
  We have a surjection of the form $L^n\to M$ and thus by left exactness
  $\ad M$ is a subsheaf of the torsion-free sheaf $\ad L^n$.
\end{proof}

Call a coherent sheaf $M$ ``reflexive'' if $R\ad M$ is a sheaf.

\begin{cor}
  A coherent sheaf $M$ is reflexive iff there exists a sheaf $N$ such that
  $M\cong \ad N$.
\end{cor}

\begin{proof}
  If $M$ is reflexive, then $M\cong R\ad \ad M=\ad(\ad M)$, so it remains to
  show that $\ad N$ is always reflexive.  But this follows from the
  spectral sequence $R^p\ad R^q\ad N\Rightarrow N$: if $R^p\ad \ad N\ne 0$ for
  $p\in \{1,2\}$, then this would contribute to the positive cohomology of
  the limit of the spectral sequence.
\end{proof}

\begin{cor}
  If $M$ is reflexive, then $M|^L_Q$ is a reflexive sheaf.
\end{cor}

\begin{proof}
  The known interaction between duality and restriction to $Q$ tells us
  that for any reflexive sheaf $M$, $R\ad(M)|^L_Q\cong
  R\sHom(M|^L_Q,\omega_Q)$.  The first object is nonpositive (since $R\ad
  M$ is a sheaf), while the second object is clearly nonnegative, and thus
  both objects are sheaves on $Q$.  The standard spectral sequence for the
  second object immediately implies that $M|^L_Q$ is a reflexive sheaf.
\end{proof}

\begin{cor}
  For any coherent sheaf $M$, $R^p\ad M$ has dimension $\le 2-p$.
\end{cor}

\begin{proof}
  We similarly conclude using the negative part of the spectral sequence
  that $R^p\ad R^2\ad M=0$ for $p<2$, and that for any line bundle $L$,
  $R^2\ad R^2\ad M$ is acyclic for $R\Hom(L,{-})$.  It follows that $R^2\ad
  R^2\ad M$ is $0$-dimensional, and thus (dualizing yet again) that $R^2\ad
  M$ is $0$-dimensional.

  For $R^1\ad M$, vanishing of the negative cohomology only directly tells
  us that $\ad R^1\ad M$ injects in $R^2\ad R^2\ad M$, but since $\ad
  R^1\ad M$ is torsion-free, its injection in a $0$-dimensional sheaf
  forces it to vanish.  But then
  \[
  \rank(R^1\ad M) = \rank(R\ad R^1\ad M) = \rank(R^2\ad R^1\ad
  M)-\rank(R^1\ad R^1\ad M),
  \]
  and since $R^2\ad R^1\ad M$ is $0$-dimensional by the previous paragraph,
  we find that $\rank(R^1\ad M)+\rank(R^1\ad R^1\ad M)=0$.  Since sheaves
  have nonnegative rank, it follows that $R^1\ad M$ has rank 0, so is $\le
  1$-dimensional as required.
\end{proof}

\begin{cor}
  A sheaf $M$ has $R^2\ad M=0$ iff it has no $0$-dimensional subsheaf.
\end{cor}

\begin{proof}
  For any $0$-dimensional sheaf $N$, we have $R\ad N\cong R^2\ad N[-2]$ and thus
  \[
  \Hom(N,M)\cong h^0(R\Hom(R\ad M,R^2\ad N[-2]))\cong \Hom(R^2\ad M,R^2\ad N),
  \]
  so that there is a nonzero morphism $N\to M$ iff there is a nonzero
  morphism $R^2\ad M\to R^2\ad N$.  In other words, if $N$ is a nonzero
  $0$-dimensional subsheaf of $M$, then $R^2\ad N$ is a nonzero quotient of
  $R^2\ad M$, while if $N$ is a nonzero quotient of $R^2\ad M$, then
  $R^2\ad N$ is a nonzero $0$-dimensional subsheaf of $M$.
\end{proof}

\begin{rem}
  We may think of this as an analogue of Proposition
  \ref{prop:hd1_is_pointless}, with ``homological dimension $1$'' replaced
  by ``$\ad$-dimension $1$''.
\end{rem}

\begin{cor}
  If $M$ is pure 1-dimensional, then $R^1\ad R^1\ad M\cong M$.
\end{cor}

This gives an easy proof of the ``finitely partitive'' property considered
in \cite{ChanD/NymanA:2013}.

\begin{prop}
  Let $M=M_0\supset M_1\supset\cdots $ be a descending chain of coherent
  sheaves such that $M_i/M_{i+1}$ has dimension $d=\dim(M)$ for all $i\ge
  0$.  Then the chain must be finite.
\end{prop}

\begin{proof}
  If $d=2$, then this follows from the fact that each subquotient has
  positive rank, so the chain terminates in at most $\rank(M)$ steps, while
  if $d=0$, we similarly must terminate in at most $\chi(M)$ steps.  For
  $d=1$, we first replace each $M_i$, $i>0$ by a slightly larger sheaf
  $\widetilde{M}_i$ defined by letting $\widetilde{M}_i/M_i$ be the maximal
  $0$-dimensional subsheaf of $\widetilde{M}_{i-1}/M_i$.  This makes each
  subquotient a pure $1$-dimensional sheaf without affecting the length of
  the chain.  Furthermore, if $N$ is the maximal $0$-dimensional subsheaf
  of $M$, then $N$ is contained in each $\widetilde{M}_i$, and thus we may
  quotient the chain by $N$ without affecting the dimensions of the
  subquotients.  We have thus reduced to the case that $M$ and every
  subquotient is pure $1$-dimensional.  But then applying $R^1\ad$ gives
  an {\em ascending} chain
  \[
  0\subset R^1\ad(M_0/M_1)\subset R^1\ad(M_0/M_2)\subset\cdots\subset R^1\ad M_0
  \]
  which terminates since $\coh X$ is Noetherian.
\end{proof}

\begin{rem}
  Of course, once we have an understanding (and have established the
  existence!) of ample divisors, we will be able to sharpen the
  $1$-dimensional considerably; the length of the chain is bounded by
  $D_a\cdot c_1(M)$ for any ample divisor $D_a$.
\end{rem}

\section{Planes}

The above results also hold with minor changes in the case of an iterated
blowup of a noncommutative plane.  The inductive steps are all the same,
with the only difference being the base case.  For a noncommutative plane,
the semiorthogonal decomposition shows
\[
K^0(X)\cong \Z[\sO_X]\oplus \Z[\sO_X(-1)]\oplus \Z[\sO_X(-2)],
\]
from which (together with constancy on flat deformations) it follows that
$\chi([\sO_X(-d)],[M])$ is a quadratic polynomial in $d$ and $[M]$ is
uniquely determined by this polynomial.  Moreover, this polynomial is
nonnegative for $d\gg 0$ (since $\theta^{-1}$ is ample), with equality iff
$M=0$.  The rank of $M$ is given by the coefficient of $d^2/2$ in this
polynomial, while $[pt]$ is represented by the polynomial $1$, and
$NS(X)\cong \Z$ is generated by the class $h$ (with $h^2=1$) corresponding
to the polynomials with leading term $d$, with $K=-3h$.

That morphisms between line bundles are injective can be readily
established by choosing an integral component $C$ of $Q$ and observing that
if a product of morphisms of line bundles vanishes, then so does the
product of their restrictions to $C$, and thus one of the two morphisms
vanishes on $C$; the result then follows by induction on the overall
difference of Chern classes.

This is all we need for the proof of Theorem
\ref{thm:blowup_of_plane_is_F1}, which states that a one-point blowup of a
noncommutative plane is a noncommutative Hirzebruch surface.  (We need to
know that $H^0$ of a nonzero $0$-dimensional sheaf is nonzero, but this
follows from vanishing of $H^2$ and positivity of the Euler
characteristic.)  Given that, we can then reduce the remaining claims to
statements on such surfaces.  In particular, that a noncommutative plane is
a maximal order iff $q$ is torsion (and commutative iff $q$ is trivial)
follows by choosing points in distinct orbits, and observing that the
maximal order property follows from the fact that both blowups are maximal
orders.

\chapter{Birational geometry}
\label{chap:birat_noncomm}

\section{Elementary transformations}

We can now finish the proof that the usual isomorphisms in the birational
geometry of commutative surfaces carry over to the noncommutative setting.
We recall that in each case we have constructed a derived equivalence and a
pair of left exact functors related by the derived equivalence, and thus to
show that the derived equivalence induces an abelian equivalence, it
suffices to show that the functors $f$ satisfy the property that for any
nonzero sheaf $M$, $f_*\theta^l M\ne 0$ for some $l\in \Z$.

For elementary transformations, the functor is $\rho_{0*}\alpha_*$,
so suppose $M$ is a nonzero sheaf such that $\rho_{0*}\alpha_*\theta^l M=0$
for all $l\in Z$.  Since
\[
-\rank(R^1(\rho_{0*}\alpha_*)\theta^l M)
=
\rank(R\rho_{0*}R\alpha_*\theta^l M)
=
\rank(R\rho_{0*}R\alpha_*M)-l\rank(M),
\]
we must have $\rank(M)=0$, so that $M$ is a (pure) $1$-dimensional sheaf.
(Here we have used the fact that $\rho_{0*}\alpha_*$ has cohomological
dimension 1.)  We then have
\[
f\cdot c_1(M)
=
-\rank(R^1(\rho_{0*}\alpha_*)M).
\]
If this were not 0, then there would be line bundles on $C_0$ such that
$\Ext^2(\alpha^*\rho_0^*L,M)\ne 0$, which is impossible for a
$1$-dimensional sheaf, and thus $f\cdot c_1(M)=0$, so that $c_1(M)\in \Z
f+\Z e_1$.  Moreover, since now $h^1(R\rho_{0*}R\alpha_*\theta^l M)$ is
$0$-dimensional, we find that $R\Hom(\sO_{C_0},R\rho_{0*}R\alpha_*\theta^l
M)$ is supported in degree $1$, so its Euler characteristic is negative.
It follows that $c_1(M)\cdot K=0$, and thus $c_1(M)=a(f-2e_1)$ for some
nonzero $a\in \Z$.

Choose such an $M$ that minimizes $|a|$.  We have $\chi(M,M)=4a^2\ge 4$,
and thus $\dim\Hom(M,M)+\dim\Ext^2(M,M)\ge 4$.  If $\dim\Hom(M,M)>1$, then
the endomorphism ring of $M$ contains zero divisors, and thus there are
nonzero endomorphisms with nontrivial kernel.  The kernel and cokernel are
both subsheaves of $M$, so have nonzero Chern class proportional to
$f-2e_1$, which by minimality of $|a|$ forces them to have opposite signs,
contradicting Proposition \ref{prop:neg_of_eff_not_eff}.  We thus conclude
that $\dim\Hom(M,M)=1$ and thus $\dim\Ext^2(M,M)\ge 3$.  But this is
equivalent to $\dim\Hom(M,\theta M)\ge 3$, and the same argument tells us that
any such morphism must be injective.  Since $M$ and $\theta M$ have the
same class in the numeric Grothendieck group, such a morphism must be an
isomorphism, and we get the desired contradiction.  We have thus shown the
following.

\begin{thm}\label{thm:elem_xform_works}
  Let $X$ be a quasi-ruled surface and $x\in Q$ a point.  Then the blowup
  of $X$ at $x$ is also a blowup of the quasi-ruled surface obtained by
  replacing the sheaf bimodule ${\cal E}$ by the kernel of the point in
  $\Quot({\cal E},1)$ associated to $x$.
\end{thm}

\begin{rem}
  The construction of the derived equivalence via the associated
  commutative ruled surface shows that this operation is an involution.
  The above description may appear to violate this (both steps make ${\cal
    E}$ smaller!), but this is an illusion coming from the fact that ${\cal
    E}$ is not uniquely determined by the quasi-ruled surface.  Indeed,
  what one finds is that the sheaf bimodule resulting after two steps is
  isomorphic to the tensor product of ${\cal E}$ with the line bundle
  $\pi_1^*{\cal L}(-\pi_1(x))$, and thus determines the same quasi-scheme
  $(X,\sO_X)$.
\end{rem}  

We also note the following useful fact we essentially proved in the course
of the above argument.

\begin{prop}\label{prop:f_is_nef}
  If $D\in \NS(X)$ is effective, then $D\cdot f\ge 0$.
\end{prop}

\begin{proof}
  Indeed, if $M$ is a pure $1$-dimensional sheaf, then
  \[
  f\cdot c_1(M)
  =
  \rank(R^0\rho_{0*}\alpha_{1*}\cdots\cdots \alpha_{m*}M)
  -
  \rank(R^1\rho_{0*}\alpha_{1*}\cdots\cdots \alpha_{m*}M),
  \]
  and the argument above shows
  $
  R^1\rho_{0*}\alpha_{1*}\cdots\cdots \alpha_{m*}M=0$.
\end{proof}

Since $K\cdot f=-2$, we immediately conclude the following.

\begin{cor}
  The class $K$ is ineffective.
\end{cor}

This has the following useful consequence.

\begin{cor}
  If $M$, $N$ are torsion-free sheaves of rank 1, then at least one of
  $\Hom(M,N)$ and $\Ext^2(M,N)$ vanishes.
\end{cor}

\begin{proof}
  A nonzero morphism $M\to N$ must be injective, since otherwise the image
  would be a $\le 1$-dimensional subsheaf of $N$.  The cokernel is thus a
  $1$-dimensional sheaf of Chern class $c_1(N)-c_1(M)$, which must
  therefore be effective whenever $\Hom(M,N)\ne 0$.  Since
  $\Ext^2(M,N)\cong \Hom(N,\theta M)^*$, we similarly find that if
  $\Ext^2(M,N)\ne 0$ then $c_1(M)-c_1(N)+K$ is effective.  But these
  divisor classes cannot be simultaneously effective, since their sum $K$
  is ineffective.
\end{proof}

\section{Isomorphisms of Hirzebruch surfaces}

For the remaining (rational) cases, the functor we consider is just the
usual global sections functor $\Gamma(M):=\Hom(\sO_X,M)$.  Here we have
some additional complications coming from the fact that this functor has
cohomological dimension 2.  Luckily, this is only a possibility when $M$ is
$2$-dimensional, and we can rule that case out.

\begin{lem}
  Let $M$ be a sheaf on a noncommutative rational or rationally quasi-ruled
  surface $X$.  Then $\Hom(M,\theta^l \sO_X)=0$ for $l\gg 0$.
\end{lem}

\begin{proof}
  It suffices to show this when $M$ is a line bundle, since in general
  there is a surjection to $M$ from a power of a line bundle.  On a
  noncommutative plane, the line bundles have the form $\sO_X(d)$, and we
  have $\Hom(\sO_X(d),\theta^l \sO_X)=0$ for $l>-d/3$.  Otherwise, if
  $\Hom(L,\theta^l\sO_X)\ne 0$, then $c_1(L)+l K_X$ is effective, but
  $(c_1(L)+lK_X)\cdot f=c_1(L)\cdot f-2l$ is eventually negative.
\end{proof}

\begin{rem}
  By Serre duality, this implies that $H^2(\theta^{-l}M)=0$ for $l\gg 0$.
\end{rem}

This implies that if $M$ is a sheaf on a noncommutative Hirzebruch surface
or a one-point blowup of a noncommutative plane such that $H^0(\theta^l
M)=0$ for all $l$, then $M$ must be $1$-dimensional.  Indeed, if
$\rank(M)>0$, then then we have $\chi(\theta^{-l}M) = 4\rank(M) l^2 +
O(l)$, and thus $\chi(\theta^{-l}M)>0$ for $l\gg 0$, but this is impossible
since $H^2(\theta^{-l}M)=0$ for $l\gg 0$.

We then find as in the elementary transformation case that $c_1(M)\cdot
K=0$, which uniquely determines $c_1(M)$.  For $F_1$ and the one-point
blowup of $\P^2$, we have $\chi(M,M)\ge 8$, and obtain a
contradiction as before.

\begin{thm}\label{thm:blowup_of_plane_is_F1}
  Let $X$ be a noncommutative Hirzebruch surface such that the
  corresponding commutative surface is $F_1$.  Then $X$ is the one-point
  blowup of a noncommutative plane.  Conversely, any one-point blowup of a
  noncommutative plane is the noncommutative Hirzebruch surface associated
  to a point of $\Pic^0(Q)$ with $Q$ an anticanonical curve in $F_1$.
\end{thm}

The $\P^1\times \P^1$ case is somewhat more subtle, as the inequality we
obtain is only that $\chi(M,M)\ge 2$, which is in fact possible.  However,
this does at least pin down that $c_1(M)=\pm (s-f)$, and $-s+f$ can be
ruled out since $(-s+f)\cdot f<0$.  If $M$ is not transverse to $Q$ (i.e.,
if $M|^{\bf L}_Q$ is not a sheaf), then this gives a nonzero subsheaf of
$\theta M$ supported on $Q$, which must have the same Chern class, and thus
some component of $Q$ is a $-2$-curve on the corresponding commutative
surface (and we may take $M=\sO_{Q'}(-1)$ where $Q'$ is that component).
Otherwise, since $c_1(M)\cdot Q=0$, we must have $M|^{\bf L}_Q=0$ and thus
$\det(M|^{\bf L}_Q)\cong \sO_Q$.  In particular, the image of $s-f$ in
$\Pic^0(Q)$ must be a power of $q\in \Pic^0(Q)$, and we conclude the
following.

\begin{thm}
  Let $X$ be the noncommutative Hirzebruch surface associated to a triple
  $(\P^1\times \P^1,Q,q)$ with $Q\subset \P^1\times \P^1$ an anticanonical
  curve and $q\in \Pic^0(Q)$.  If the line bundle $\sO_{\P^1\times
    \P^1}(1,-1)|_Q$ is not a power of $q$, then swapping the two rulings
  gives an isomorphic noncommutative Hirzebruch surface.
\end{thm}

Note that when the ``not a power of $q$'' hypothesis fails, then the
commutative surface associated to some twist of $X$ will be $F_2$.  There
we have the following corollary of Proposition \ref{prop:weak_reflection}.

\begin{prop}
  Suppose that $X$ is the noncommutative Hirzebruch surface associated to a
  triple $(F_2,Q,q)$ with $Q\subset F_2$ an anticanonical curve disjoint
  from the $-2$-curve and $q\in \Pic^0(Q)$ of order $r\in [1,\infty]$.
  Then there is an irreducible $1$-dimensional sheaf $\sO_{s-f}(-1)$ of
  Chern class $s-f$ and Euler characteristic 0, and if $0\le b-a\le r$,
  there is a short exact sequence
  \[
  0\to \sO_X(-bs-af)\to \sO_X(-as-bf)\to \sO_{s-f}(-1)^{b-a}\to 0.
  \]
\end{prop}

\begin{proof}
  Let $x\in Q$ be any point.  If we blow up $X$ in $x$, then an elementary
  transformation lets us interpret $\widetilde{X}$ as a blowup of a
  noncommutative $F_1$, and thus as a two-fold blowup of a noncommutative
  $\P^2$, in which the same point is blown up each time.  The exceptional
  classes on the latter are $e'_1=s-e_1$ and $e'_2=f-e_1$, and thus there
  is a sheaf $\sO_{e'_1-e'_2}(-1)$, the direct image of which is the
  desired sheaf $\sO_{s-f}(-1)$ on $X$.  Moreover, there is a short exact
  sequence
  \[
  0\to \sO_{\widetilde{X}}(-be'_1-a'e_2)\to \sO_{\widetilde{X}}(-ae'_1-be'_2)
  \to \sO_{e'_1-e'_2}(-1)^{b-a}\to 0
  \]
  or equivalently
  \[
  0\to \sO_{\widetilde{X}}(-bs-af+(a+b)e_1)
  \to \sO_{\widetilde{X}}(-as-bf+(a+b)e_1)
  \to \sO_{e'_1-e'_2}(-1)^{b-a}\to 0
  \]
  The result follows by taking direct images down to $X$.
\end{proof}

\begin{rem}
  If $Q$ is not disjoint from the $-2$-curve, then it contains the
  $-2$-curve, so that one may define $\sO_{s-f}(-1)$ as the appropriate
  line bundle on that component of $Q$, and one again obtains a short exact
  sequence for $0\le b-a\le 1$.
\end{rem}

\section{Blowing down}

The other thing we want to show vis-\`a-vis birational geometry is that
``rationally quasi-ruled'' and ``birationally quasi-ruled'' are essentially
the same, or in other words the following.

\begin{thm}\label{thm:birationally_qr_is_rationally_qr}
  Suppose that $X$ is a noncommutative surface such that some iterated
  blowup of $X$ is an iterated blowup of a quasi-ruled surface with
  associated curves $(C_0,C_1)$.  Then either $X$ is an iterated blowup of a
  quasi-ruled surface over the same pair of curves, or $C_0\cong C_1\cong
  \P^1$ and $X$ is a noncommutative plane.
\end{thm}

By an easy induction, this reduces to showing that if a one-fold blowup of
$X$ is an iterated blowup of a quasi-ruled surface or noncommutative
plane, then the same is true for $X$, and thus to understanding when we
can blow down a $-1$-curve on a rationally quasi-ruled surface.

Thus let $\alpha:\widetilde{X}\to X$ be a van den Bergh blowup with
exceptional sheaf $\sO_e(-1)$, such that $\widetilde{X}$ is rational or
rationally quasi-ruled.  We first need to deal with the technical issue
that $X$ may not have a well-defined dimension function, and thus it is not
immediately obvious that the exceptional sheaf is $1$-dimensional relative
to the dimension function on $\tilde{X}$.  This is not too difficult to
show: $X$ inherits a Serre functor from $\widetilde{X}$, and if
$M_l:=\alpha^*(\sO_X(lC))(-l)$ (where $C$ is the divisor on $X$ containing
the point being blown up), then we have short exact sequences
\[
0\to M_l\to M_{l+1}\to \sO_e(-l-1)\to 0,
\]
since $\sO_X(lC)$ is a weak line bundle along $C$.  An easy induction from
$\theta^{-1}\sO_e(-1)\cong \sO_e$ shows that $\sO_e(-l-1)\cong
\theta^{-l}\sO_e(-1)$ and thus has the same rank, implying that
$\rank(M_l)=\rank(\sO_e(-1))l+O(1)$.  Since ranks of coherent sheaves on
$\tilde{X}$ are nonnegative, we conclude that $\rank(\sO_e(-1))=0$ as
desired.

Let $e$ denote the Chern class of $\sO_e(-1)$.  As one might expect, this
behaves numerically like a $-1$-curve on a commutative surface.  Since
$\sO_e(-1)$ is exceptional, $-e^2 = \chi(\sO_e(-1),\sO_e(-1))=1$.  (In
particular, $\widetilde{X}$ cannot be a noncommutative plane, since
$\NS(\widetilde{X})\cong \NS(\P^2)$ has no classes of self-intersection $-1$.)
In addition, since
\[
1 = \chi(\sO_{\tau p}) =
\chi(\sO_e)-\chi(\sO_e(-1))=\chi(\sO_e)-\chi(\theta\sO_e),
\]
we have $e\cdot K=-1$.  (Note that we can compute $\chi(\sO_{\tau p})$
inside $\coh(\widetilde{C})$, and when $e$ is a component of
$\widetilde{C}$, we can do the same for $\chi(\sO_e(-l))$.)  Moreover,
since $R\Gamma(\sO_e(-1))=0$, $\sO_e(-1)$ must be pure $1$-dimensional.

Since $\rank(\sO_e(d))=0$, $\theta\sO_e(d)$ and $\sO_e(d)$ have the same
Chern class, so that $c_1(\sO_e(d))=e$.  It then follows that
$\sO_e/\sO_e(-1)$ has numerical class $\chi(\sO_e/\sO_e(-1))[\pt]=[\pt]$.
Thus this sheaf is a point, and the same holds (by algebraic equivalence)
for any point $x\in\widetilde{C}$.  In particular, for any $x\in C$,
$[L\alpha^*\sO_x]=[\pt]$.


Suppose $M$ is a quotient of $\sO_e(-1)$ by a proper subsheaf $N$.  Then
the long exact sequence gives $R^1\alpha_*M=0$ and $\alpha_*M\cong
R^1\alpha_*N$, so that $\alpha_*M$ is an extension of point sheaves on $C$,
thus $[L\alpha^*\alpha_*M]\propto [\pt]$ and $[M]$ is in the span of
$[\sO_e(-1)]$ and $[\pt]$.  Since $c_1(M)+c_1(N)=e$ is a sum of effective
divisors, the only possibilities are $c_1(M)=e$ or $c_1(M)=0$, and the
former case can be ruled out since $\sO_e(-1)$ is pure $1$-dimensional.  We
thus conclude that $\sO_e(-1)$ is irreducible in a suitable sense: any
proper quotient is $0$-dimensional.

When $X_0$ (the surface of which $\tilde{X}$ is an iterated blowup)
is strictly quasi-ruled, there is an additional constraint on $e$,
namely that $e\cdot Q=1$.  (This constraint is equivalent to $e\cdot
K=-1$ if $X_0$ is ruled, and is actually redundant unless
$g(C_0)+g(C_1)=1$!)  Since $K+Q$ is a nonnegative multiple of $f$,
Proposition \ref{prop:f_is_nef} tells us that $e\cdot Q\ge e\cdot
(-K)=1$, so it suffices to show that if equality fails, then
$e\cdot f=0$.  Suppose $e\cdot Q>1$.  Then $e$ meets $Q$ in at least
2 points (possibly by being a component of $Q$), and thus one finds
that $\dim\Hom(\sO_e(-1),\sO_e)>1$.  The cokernel of any such
homomorphism is a point sheaf, and thus there is an induced family
of point sheaves parametrized by $\P^1$.  Taking images to $X_0$
gives a family of point sheaves on $X_0$ and thus a morphism from
$\P^1$ to $\Quot(\sO_{X_0},1)$ (which is well-defined since $X_0$
is a maximal order).  Since $X_0$ is strictly quasi-ruled, its curve
of points is the union of an irreducible curve $\hat{Q}$ of positive
genus and a collection of fibers, and thus the morphism must map
to a point or a fiber, implying in either case that $e\cdot f=0$
(which we can compute as an intersection in the center).

Thus in general, given an iterated blowup $X$ of a quasi-ruled surface, we
need to classify the elements $e\in \NS(X)$ such that $e^2=e\cdot K=-e\cdot
Q=-1$ and there exists an irreducible $1$-dimensional sheaf $\sO_e(-1)$
with Chern class $e$ and vanishing $R\Gamma$.  Call such a class a ``formal
$-1$-curve'' on $X$; our objective is to show that any such class can be
blown down, and the result is as described in the Theorem.

Note first that if $m=0$ with $s^2=-1$ ($e^2=-1$ cannot happen in the even
parity case, since then any self-intersection is even), then the equations
$e^2=e\cdot K=-1$ imply $g(C_0)+g(C_1)\in \{0,2\}$, as $e\cdot f$ satisfies
a quadratic equation with discriminant $9-4g(C_0)-4g(C_1)$.  If
$g(C_0)+g(C_1)=2$, the unique solution of the equations is $e=-s$, which is
ineffective since then $e\cdot f<0$, while if $g(C_0)=g(C_1)=0$, then the
unique integral solution is $e=s$, which can be blown down (to a
noncommutative plane) precisely when there is an irreducible sheaf
$\sO_e(-1)$ (i.e., when the corresponding commutative surface is $F_1$). So
we may reduce to the case $m>0$.

A key observation is that we have shown that $X$ is not in general uniquely
representable as an iterated blowup of a quasi-ruled surface.  Changing the
representation of course has no effect on $X$, but {\em does} change the
basis we are using for $\NS(X)$.  Thus to show that formal $-1$-curves can
be blown down and the result is again rationally quasi-ruled (or rational),
it suffices to show that in {\em some} representation of $X$ as an iterated
blowup, the expression for $e$ in the resulting basis is $e_m$.  (Compare
the proof of Theorem \ref{thm:two_orbits_of_blowdowns_comm}.)

\begin{prop}\label{prop:formal_-1_can_be_blown_down}
  Let $X_m$ be an iterated blowup of a quasi-ruled surface $X_0$, with
  $m>0$, and let $e$ be a formal $-1$-curve on $X_m$.  Then there is a
  sequence of interchanges of commuting blowups, elementary
  transformations, and (if $X$ is rational) swappings of rulings on
  $\P^1\times \P^1$ such that in the resulting basis of $\NS(X_m)$,
  $e$ is represented by $e_m$.
\end{prop}

\begin{proof}
  Consider the representation of $e$ in the associated basis of $\NS(X_m)$.
  If $e\cdot (e_i-e_{i+1})<0$ for some $i$, then there are two
  possibilities: either we may interchange the two blowups, and thus obtain
  a basis in which $e\cdot (e_i-e_{i+1})>0$ (essentially reflecting in
  $e_i-e_{i+1}$), or the blowups do not commute, in which case there is a
  $1$-dimensional sheaf $M$ on $X_{i+1}$ of class $e_i-e_{i+1}$, which is
  either a component of $Q$ or has irreducible pullback to $X_m$.  But then
  the fact that the intersection number is negative forces
  $\Hom(\sO_e(-1),\alpha_m^*\cdots \alpha_{i+2}^*M)\ne 0$ or
  $\Hom(\theta^{-1}\alpha_m^*\cdots \alpha_{i+2}^*M,\sO_e(-1))\ne 0$,
  neither of which can happen by irreducibility of $\sO_e(-1)$.  (An easy
  induction shows that any irreducible constituent of the pullback of $M$
  has Chern class of the form $e_j-\sum_{k\in S} e_k$ where $j\ge i$ and
  $|S|$ consists of integers greater than $j$.)  Similarly, if $D\cdot
  e_1\ge D\cdot f/2$, then we may perform an elementary transformation to
  subtract $D\cdot e_1$ from $D\cdot f$.

  Conjugation by an elementary transformation preserves the finite Weyl
  group $W(D_m)$, and thus this process will eventually terminate in a
  basis such that $e\cdot e_m\le e\cdot e_{m-1}\cdots\le e\cdot e_1\le
  (e\cdot f)/2$.  If $e\cdot e_m<0$, then irreducibility gives $e=e_m$ and
  we are done, so we may assume $e\cdot e_m\ge 0$, and thus we have an
  expression
  \[
  e = ds+d'f-\lambda_1e_1-\cdots-\lambda_me_m
  \]
  in which $\lambda$ is an ordered partition ($\lambda_1\ge \lambda_2\ge\cdots$) with largest part at most $d/2$.

  We now split into three cases, depending on the ``genus''
  $g:=(g(C_0)+g(C_1))/2$ appearing in the expression for $K_X$.
  If $g\ge 1$, then we may write
  \[
  ds+d'f = -(d/2)K_{X_0} + ((2g-2)d+n/2)f
  \]
  for some integer $n$ (with parity depending on the parity of $2g$ and
  $s^2$), in terms of which we have
  \[
  (ds+d'f)^2 = (2g-2)d^2+dn
  \qquad\text{and}\qquad
  (ds+d'f)\cdot K_{X_0} = -n.
  \]
  Since $e\cdot K = -1$, we find that $\lambda$ is a partition of $n-1$
  into at most $m$ parts all at most $d/2$ (and thus in particular $n\ge
  1$).  It follows that
  \[
  \sum_i \lambda_i^2 \le d(n-1)/2,
  \]
  and thus
  \[
  e^2 = (2g-2)d^2+dn-\sum_i \lambda_i^2
  \ge (2g-2)d^2+d(n+1)/2
  \ge 0,
  \]
  giving a contradiction.

  If $g=1/2$, then $X_0$ is strictly quasi-ruled, and not of the
  $4$-isogeny type, so that $Q=-K+lf$ for some $l>0$.  Since $e\cdot
  (Q+K)=0$, we conclude that $e\cdot f=0$, and thus $d=0$, forcing
  $\lambda=0$ and $e^2=0$.

  In the $g=0$ case (i.e., rational surfaces), we must consider an
  additional reflection, namely the reflection in $s-f$ when $X_0$ is an
  even Hirzebruch surface or in $s-e_1$ when $X_0$ is an odd Hirzebruch
  surface.  Again, we can perform this reflection precisely when the
  corresponding class is not effective, and thus irreducibility of
  $\sO_e(-1)$ again ensures that we can perform the reflection whenever we
  would want to.  Moreover, if we have already put $e$ into standard form
  relative to interchanges of blowups and elementary transformations, then
  the reflection in $s-f$ or $s-e_1$ will decrease $e\cdot f$, and since
  $e$ is effective (thus $e\cdot f\ge 0$ relative to {\em any} blowdown
  structure), the resulting procedure (alternating between the two
  reductions) will necessarily terminate.  (Note that it is no longer the
  case that there is a finite Coxeter group forcing termination!)

  For convenience, we then perform an elementary transformation if needed
  to make $X_0$ an even Hirzebruch surface, so obtain an expression of the
  form
  \[
  e = af + b(s+f)+c(s+f-e_1)+d(2s+2f-e_1-e_2) - \sum_i \lambda_i e_{i+2}
  \]
  in which $a,b,c,d\ge 0$ and $\lambda$ is a partition with $\lambda_1\le d$.
  (Of course, if $m=1$, we must have $\lambda=d=0$ and find
  $e^2=2ab+2ac+2b^2+4bc+c^2\ge 0$) We have $|\lambda| = 2a+4b+3c+6d-1$ and
  thus
  \begin{align}
  e^2 &= 2ab+2ac+4ad+2b^2+4bc+8bd+c^2+6cd+6d^2
  - \sum_i \lambda_i^2\notag\\
  &\ge 2ab+2ac+4ad+2b^2+4bc+8bd+c^2+6cd+6d^2
  - \lambda_1|\lambda|\notag\\
  &\ge 2ab+2ac+2ad+2b^2+4bc+4bd+c^2+3cd+d\notag\\
  &\ge 0,
  \end{align}
  again giving a contradiction.
\end{proof}  

\begin{rems}
  It is worth noting that this only applies to the van den Bergh blowup; in
  the case of a maximal order, the underlying commutative surface $Z$ may
  very well have $-1$-curves that cannot be blown down via the above
  Proposition.  The issue is that the local structure of the resulting
  order on the blown down surface at the base point may be ``singular'' in
  a suitable sense.  One example of this comes from biquadratic extensions
  of $\P^1$.  Any invertible sheaf on such an extension produces a
  quasi-ruled surface corresponding to a quaternion order on a Hirzebruch
  surface, and when the Euler characteristic of the invertible sheaf has
  the correct parity (depending on the various genera), the Hirzebruch
  surface is, at least generically, $F_1$.  In characteristic 0, one can
  represent this as the Clifford algebra associated to a quadratic form,
  and find that the blown down algebra is the Clifford algebra associated
  to a quadratic form on a vector bundle on $\P^2$ such that the form
  vanishes identically at the base point of the blowup.  In particular, the
  base change of the order to the complete local ring at that point is
  local but not of global dimension 2.  Thus the blow down in this case
  more closely resembles the contraction of a $-2$-curve than the
  contraction of a $-1$-curve.
\end{rems}

\begin{rems}
  When $g>0$ and thus the relevant group is finite, we may conclude that in
  the original basis for $\NS(X)$, we have $e=e_i$ or $e=f-e_i$ for
  some $1\le i\le m$.  It is worth noting that these classes are always
  effective and satisfy all of the numerical conditions to be a formal
  $-1$-curve; the only condition they can violate is irreducibility, which
  happens iff $e_i-e_j$ or $f-e_i-e_j$ (respectively) is effective for some
  $j>i$.
\end{rems}

The above discussion shows that a typical rationally quasi-ruled surface
has many blowdown structures, but also gives us a fair amount of control.
Specifying a blowdown structure is essentially equivalent to specifying a
basis of $\NS(X)$ in standard form, and thus any two blowdown structures
specify a change of basis matrix which respects both $K$ and the
intersection pairing.  The atomic birational transformations all act in
relatively simple ways on the basis; in particular, the commutation of
blowups and the exchange of rulings on $\P^1\times \P^1$ both act as
reflections relative to the intersection pairing, with corresponding roots
$e_i-e_{i+1}$ and $s-f$ respectively.  Moreover, although an elementary
transformation does not have such a description (since it changes the
parity and thus the intersection matrix), the conjugate of the reflection
in $e_1-e_2$ by an elementary transformation is the reflection in
$f-e_1-e_2$.  (Similarly, if $X_0$ is an odd Hirzebruch surface that blows
down to $\P^2$, then reflection in $s-e_1$ acts as commutation on the
two-point blowup of $\P^2$.)  In each case, the condition for the image
under the reflection to still be a blowdown structure is precisely that the
corresponding root is ineffective.  With this in mind, we call the divisor
classes $f-e_1-e_2, e_1-e_2,\dots, e_{m-1}-e_m$ ``simple roots'', along
with $s-f$ when $X_0$ is an even Hirzebruch surface and $s-e_1$ when $X_0$
is an odd Hirzebruch surface.  Note that these are the simple roots of a
Coxeter group of type $D_m$ for $g>0$ and $E_{m+1}$ for $g=0$.  This,
together with the reduction of Proposition
\ref{prop:formal_-1_can_be_blown_down}, motivates the following definition,
following the commutative rational case.

\begin{defn}
  Let $X$ be a noncommutative rationally quasi-ruled surface with $m\ge 1$.
  Then the {\em fundamental chamber} in $\NS(X)$ (relative to a given
  blowdown structure) is the cone of divisors having nonnegative
  intersection with every simple root and with $e_m$, as well as $f-e_1$
  when $m=1$.
\end{defn}

\begin{prop}
  Let $X$ be a noncommutative rationally quasi-ruled surface.  Then any two
  blowdown structures on $X$ with $X_0$ of the same parity are related by a
  sequence of reflections in ineffective simple roots.
\end{prop}

\begin{proof}
  Let $s,f,e_1,\dots,e_m$ and $s',f',e'_1,\dots,e'_m$ be the bases
  corresponding to the two blowdown structures.  For $m\ge 1$, we have
  shown that there is a sequence of simple reflections and elementary
  transformations taking the second blowdown structure to one in which
  $e'_m=e_m$.  For $m\ge 2$, this is unchanged by an elementary
  transformation, and since the conjugate of a simple reflection by an
  elementary transformation is a simple reflection, we may arrange to
  eliminate all of the elementary transformations.  We thus see that the
  second blowdown structure is related by a sequence of reflections in
  ineffective simple roots to one coming from a blowdown structure on
  $X_{m-1}$ of the same parity, so that the claim follows by induction.

  For $m=1$, a formal $-1$-curve determines the parity of the blowdown
  structure, and the only case in which there is more than one of the
  correct parity is when $X_0$ is an odd Hirzebruch surface and both $s$
  and $e_1$ are formal $-1$-curves.  But in that case $s$ is a formal
  $-1$-curve iff $s-e_1$ is ineffective, and thus again they are related by
  a reflection.  Similarly, for $m=0$, the ruling is uniquely determined
  unless $X_0$ is an even Hirzebruch surface on which $s-f$ is ineffective,
  in which case the corresponding reflection swaps the rulings.
\end{proof}

The fact that a commutative ruled surface of genus $>0$ has a canonical
ruling has the following noncommutative analogue.

\begin{prop}\label{prop:f_is_unique_if_quasi-ruled}
  Let $X$ be a noncommutative rationally quasi-ruled surface of genus $>0$.
  Then the functor $R\rho_{0*}R\alpha_{1*}\cdots R\alpha_{m*}$ is independent
  of the choice of blowdown structure.
\end{prop}

\begin{proof}
  Reflection in $e_i-e_{i+1}$ preserves the three factors
  \[
  R\rho_{0*}R\alpha_{1*}\cdots R\alpha_{(i-1)*},
  R\alpha_{i*}R\alpha_{(i+1)*},
  R\alpha_{(i+2)*}\cdots R\alpha_{m*},
  \]
  so preserves the composition, while an elementary transformation
  preserves the two factors
  \[
  R\rho_{0*}R\alpha_{1*},
  R\alpha_{2*}\cdots R\alpha_{m*}.
  \]
  Since any two blowdown structures are related by a sequence of simple
  reflections and elementary transformations, the claim follows.
\end{proof}

\begin{rem}
  An alternate approach is to observe that $f$ is the unique effective
  class such that $f^2=0$, $f\cdot K=-2$, and that $C_0$ is the moduli
  space of irreducible $1$-dimensional sheaves of Chern class $f$ and Euler
  characteristic 1, with the morphism to $C_0$ following as in
  \cite{ChanD/NymanA:2013}.
\end{rem}

\section{The minimal lift}

If $X$ is equipped with a blowdown structure that passes through the
surface $Y=X_j$, then we have a ``morphism'' from $X$ to $Y$ in the form
of a pair of adjoint functors $\pi_*:\coh(X)\to \coh(Y)$,
$\pi^*:\coh(Y)\to \coh(X)$, as well as a functor $\pi^!:=\theta
\pi^*\theta^{-1}$, with $L\pi^!$ right adjoint to $R\pi_*$.

Say that a sheaf $M$ is {\em pointless} if it has no nonzero
$0$-dimensional subsheaf, or equivalently if $R^2\ad M=0$.

\begin{prop}
  Pointless sheaves are $\pi^*$-acyclic, and the pullback of a pointless
  sheaf is pointless.
\end{prop}

\begin{proof}
  We may immediately reduce to the case that $\pi=\alpha$ is a monoidal
  transform.  Since $R\alpha_*L\alpha^*M\cong M$, the corresponding
  spectral sequence implies $\alpha_*L_1\alpha^*M=R^1\alpha_*\alpha^*M=0$
  and gives a short exact sequence
  \[
  0\to R^1\alpha_*L_1\alpha^*M\to M\to \alpha_*\alpha^*M\to 0
  \]
  Since $R^1\alpha_*$ of a sheaf is always $0$-dimensional
  (\cite[Prop.~6.5.2]{VandenBerghM:1998}), we conclude that $R\alpha_*L_1\alpha^*M=0$ and
  thus $L_1\alpha^*M\cong \sO_e(-1)^d$ for some $d$.  But we also have
  \[
  \Hom(\sO_e(-1),L_1\alpha^*M)
  \cong
  \Hom(\sO_e,\theta L_1\alpha^*M)
  \cong
  \Hom(\sO_e,L_1\alpha^!\theta M)
  \cong
  \Hom(\sO_p,\theta M)=0,
  \]
  and thus $d=0$.  Thus $L\alpha^* M$ is indeed a sheaf, and
  $\alpha_*\alpha^*M\cong M$.  The direct image of a $0$-dimensional
  subsheaf of $\alpha^*M$ would be a $0$-dimensional subsheaf of $M$, and
  thus $\alpha^*M$ is pointless as required.
\end{proof}

As in the commutative case (Chapter \ref{chap:*!}), adjunction gives a
natural map $\pi^*M\to \pi^!M$ for any pointless sheaf, and thus we may
define the minimal lift $\pi^{*!}M:=\im(\pi^*M\to \pi^!M)$.  Most of the
results of Chapter \ref{chap:*!} carry over with little difficulty, as we
generally tried to avoid using anything about supports of sheaves other
than the dimension.  (The primary exception is those results like
Proposition \ref{prop:*!_of_line_bundle_on_curve} that include a support
condition in the hypotheses; also, one must replace ``homological
dimension'' throughout by ``$\ad$ dimension''.)  We will thus restrict our
attention to mentioning the few changes that need to be made.

In Section \ref{sec:minimal_lift}, there are only two results where we need
to make changes.  For Lemma \ref{lem:pi*_dim_1} (that $R\pi_*$ has
cohomological dimension 1), this reduces to the monoidal case (shown in
\cite{VandenBerghM:1998}) together with the fact that $R^1\alpha_*M$ is
$0$-dimensional, so acyclic for further blowdowns.  The subcase of
Proposition \ref{prop:*!_of_line_bundle_on_curve} stating that
$\pi^{*!}\sO_{Q_Y}\cong \sO_{Q_X}$ also reduces to the monoidal case, where
it is straightforward, while the main case of the Proposition makes no
sense in the noncommutative context.

We have already seen from the relevant semiorthogonal decomposition that
the results of Section \ref{sec:invisible_comm} on the category of
invisible sheaves carry over (even for those results where the proof uses
supports of sheaves!).  The exceptions that do not carry over are
Proposition \ref{prop:invisible_image} (as the hypothesis does not make
sense in the noncommutative context) and Corollary \ref{cor:invis_on_e_pi}
(since the conclusion makes no sense).

For Section \ref{sec:*!_when_anticanon}, the exact sequence of Proposition
\ref{prop:four_term_for_monoidal} becomes
\[
0\to E_1\to \alpha^*M\to \alpha^!M\to E_2\to 0,
\]
where
\begin{align}
  E_1&\cong\Ext^1(\sO_{?p},M)\otimes_k \sO_e(-1),\\
  E_2&\cong\Hom_k(\Hom(M,\sO_{?p}),\sO_e(-1));\\
\end{align}
with corresponding changes in Corollary \ref{cor:monoidal_pseudotwist}.

Theorem \ref{thm:elem_pseudo_twists} (that elementary pseudo-twists
preserve pointless sheaves) is problematic, as the argument used specific
curves on the surface.  However, there is an alternate approach.  Indeed,
the short exact sequences of Corollary
\ref{cor:pseudo_twist_is_isospectral} hold, at least as distinguished
triangles, and thus the claim that $\pi_*(\pi^{*!}M(e_f))$ is pointless
reduces to showing that in the sequence
\[
0\to M\to \pi_*(\pi^{*!}M(e_f))\to \sO_p^{r_1}\to 0,
\]
there is no partial splitting $\sO_p\to \pi_*(\pi^{*!}M(e_f))$.  But this
holds in the monoidal case, and implies the result in general.  Similarly,
in the distinguished triangle of the form
\[
R\pi_*(\pi^{*!}M(-e_f))\to M\to \sO_p^{r_2}\to
\]
we need to show that $M\to \sO_p^{r_2}$ is surjective, which again reduces
to the monoidal case.

In any event, the key result is Lemma \ref{lem:resolving_pseudotwist},
which carries over (omitting ``Poisson surface'', since we already have an
anticanonical curve) with no difficulty, and allows us (up to
pseudo-twists) to reduce many questions about $1$-dimensional sheaves to
the case that $M\cong \theta M$.

There is also a result that holds only in the noncommutative setting.

\begin{prop}\label{prop:pseudo_twists_coalesce}
  Suppose $q$ is non-torsion, $M$ is a pure $1$-dimensional
  sheaves on $X$ transverse to $Q$, and $N\subset M$ is a subsheaf with
  $c_1(N)=c_1(M)$.  Then $M$ and $N$ have a common pseudo-twist.
\end{prop}

\begin{proof}
  If $M$ and $N$ are disjoint from $Q$, then $[M]=[N]+d[\pt]$ and thus
  $\det([M]|_Q)\cong \det([N]|_Q)\otimes q^d$, implying since $q$ is
  non-torsion that $d=0$.  Otherwise, let $x$ be a point in the support of
  $M|_Q$, and let $\tilde{X}$ be the blowup of $X$ in $x$.  Applying the
  snake lemma to
  \[
  \begin{CD}
    0@>>> \theta \alpha_* \theta^{-1}\alpha^{*!}N @>>> N @>>>
    \Hom_k(\Hom(N,\sO_x),\sO_x)@>>> 0\\
    @. @VVV @VVV @VVV @.\\
    0@>>> \theta \alpha_* \theta^{-1}\alpha^{*!}M @>>> M @>>>
    \Hom_k(\Hom(M,\sO_x),\sO_x)@>>> 0
  \end{CD}
  \]
  tells us that the left vertical map is injective, and thus the two
  downward pseudo-twists again satisfy the hypotheses.  Using this, we may
  ensure that any given orbit of points in $Q$ contains only one point in
  the support of $M|_Q\oplus N|_Q$ (otherwise, do a downwards pseudo-twist
  at the ``largest'' point in the orbit).  If the orbit consists of a
  single element (necessarily a singular point of $Q$), then we do further
  downwards pseudo-twists to ensure that the minimal lifts remain
  transverse from the anticanonical curve.

  Since then $\Hom(M,\sO_{\tau^d x})=0$ for $d\ne 0$, an easy induction
  tells us that $\Hom(\alpha^*M,\sO_e(d))=0$ for all $d$  Since
  $\alpha^{*!}M$ is a quotient of $\alpha^*M$, we also have
  $\Hom(\alpha^{*!}M,\sO_e(d))=0$ for all $d$.  Similarly, for $d\ne 0$,
  \[
  \Ext^1(\sO_{\tau^{d+1}x},N)\cong \Ext^1(N,\sO_{\tau^d x})^*
  \cong \Ext^1(N|_Q,\sO_{q^d x})^*
  \cong \Hom(\sO_{q^d x},N|_Q)
  \]
  and thus $\Hom(\sO_e(d),\alpha^!N)=0$ for all $d$, implying
  $\Hom(\sO_e(d),\alpha^{*!}N)=0$.

  Now, consider the morphism
  \[
  \phi:\alpha^{*!}N\to \alpha^{*!}M.
  \]
  The direct image is injective with $0$-dimensional cokernel, and thus
  comparing the two spectral sequences for the direct image implies that
  \[
  \alpha_*\ker(\phi)=R^1\alpha_*\coker(\phi)=0,
  \]
  while $R^1\alpha_*\ker(\phi)$,$\alpha_*\coker(\phi)$ are $0$-dimensional.
  It follows immediately that the Chern classes $c_1(\ker(\phi))$ and
  $c_1(\coker(\phi))$ are both multiples of $e$.  Since $\alpha^{*!}N$ has
  no such subsheaf and $\alpha^{*!}M$ has no such quotient, it follows that
  $\phi$ is itself injective with $0$-dimensional cokernel, so we may
  induct.  Since $\alpha^{*!}M$ remains transverse to the anticanonical
  curve and $\chi(\alpha^{*!}M|_Q)<\chi(M|_Q)$, the claim follows by
  induction.
\end{proof}

\begin{rems}
  More generally, given two subsheaves $N_1$, $N_2$ of $M$ with
  $0$-dimensional quotient, applying the above to $N_1\oplus
  N_2\subset M^2$ gives a sequence of pseudotwists making $N_1\oplus N_2$
  and $M^2$ agree, and thus in particular making $N_1$ and $N_2$ agree.
\end{rems}

\begin{rems}
  One can view this as a stronger version of the statement that
  $0$-dimensional sheaves are supported on $Q$ when $q$ is non-torsion.
  Indeed, when $q$ is torsion, the claim fails precisely because $M/N$
  could have support away from $Q$.
\end{rems}

\chapter{The nef cone and ample divisors}
\label{chap:effnef}

\section{The effective monoid and the nef cone}

Since the surfaces we are deforming are projective, not just proper, we
would like to know that they have some analogue of ample divisors, and
ideally be able to control the set of such divisors.  It is unreasonable to
hope that ample divisors will correspond to graded algebras (let alone
quotients of noncommutative projective spaces), but we can still hope for
some version of Serre vanishing and global generation.  By Proposition
\ref{prop:NS_of_center}, this is at least in principle easy when $X$ is a
maximal order (though the lack of an anticanonical curve in the general
quasi-ruled case makes it difficult to make everything explicit) as then we
may identify the effective and nef cones of $X$ and its center, so that
anything in the interior of the nef cone is ample.

We will thus restrict our attention to the ruled (or rational) surface
case.  Note here that we are focused on {\em noncommutative} surfaces, and
thus any commutative fiber will come with a choice of anticanonical curve.
(It is also worth noting that our arguments do not actually work for the
excluded case of a (characteristic 2) commutative ruled surface of genus 1
with integral anticanonical curve!)

In the commutative setting, the ample divisors on a surface are best
understood as the integer interior points of the cone of nef
divisors, for which we certainly have a well-behaved analogue in the
noncommutative setting: a divisor class is nef iff it has nonnegative
intersection with every effective divisor class.

Of course, to understand the nef cone, we will first need to fully
understand the effective monoid.  To describe this, we will need one
more set of effective classes.

\begin{defn}
  A formal $-2$-curve is a $\theta$-invariant class of self-intersection
  $-2$ which is the Chern class of an irreducible sheaf.
\end{defn}

\begin{prop}
  If $\alpha$ is a formal $-2$-curve, then there is a blowdown structure in
  which $\alpha$ is a simple root.
\end{prop}

\begin{proof}
  This is by the same argument as in Proposition
  \ref{prop:formal_-1_can_be_blown_down}, with only minor changes to the
  final inequalities.
\end{proof}

\begin{thm}
  Let $X$ be a rationally ruled surface with $m\ge 1$.  Then the effective
  monoid of $X$ is generated by components of $Q$, formal $-1$-curves, and
  formal $-2$-curves.
\end{thm}

\begin{proof}
  Let $\mu$ denote the monoid so generated, which is clearly contained in
  the effective monoid.  We will proceed in two steps: first showing that
  the effective monoid is generated by $\mu$ and its dual, then showing
  that $\mu$ contains its dual.

  For the first claim, let $M$ be a pure $1$-dimensional sheaf.
  If $c_1(M)\cdot Q_i<0$ for some component $Q_i$ of $Q$, then
  $\chi(M|^{\bf L}_{Q_i})<0$ and thus $h^{-1}(M|^{\bf L}_{Q_i})\ne 0$.
  Since this gives a subsheaf of $M(-Q_i)$, it must have positive rank, and
  thus there is a subsheaf of $M$ of Chern class $Q_i$.  Similarly, if
  $D\cdot e<0$ for some formal $-1$-curve $e$, then $\chi(\sO_e(-1),M)>0$,
  so that either $\Hom(\sO_e(-1),M)\ne 0$ or $\Hom(M,\theta \sO_e(-1))\ne
  0$, giving either a subsheaf of class $e$ or a quotient sheaf of class
  $e$.  Finally, if $\alpha$ is a formal $-2$-curve with $D\cdot \alpha<0$,
  then $\chi(M_\alpha,M)>0$ for any irreducible sheaf $M_\alpha$ of the
  appropriate Chern class, and thus again there is a sub- or quotient sheaf
  of Chern class $\alpha$.

  We may thus obtain a descending sequence of subquotients of $M$ by
  repeatedly choosing sub- or quotient sheaves of Chern class a component
  of $Q$ or a formal $-1$- or $-2$-curve.  The corresponding pair of
  chains of subsheaves of $M$, one ascending, one descending, both
  terminate, and thus this process terminates.  The only way it can
  terminate is with a subquotient of $M$ having nonnegative intersection
  with all generators of $\mu$, giving an expression for $c_1(M)$ as the
  sum of an element of $\mu$ and an element of the dual monoid as required.

  To proceed further, we note first that any effective simple root (for any
  blowdown structure!) is either a formal $-2$-curve or a sum of components
  of $Q$ (with negative self-intersection) and formal $-1$-curves, so is in
  $\mu$.  Moreover, if $e$ is in the $W(D_m)$ or $W(E_{m+1})$-orbit of
  $e_m$, then when applying the procedure of Proposition
  \ref{prop:formal_-1_can_be_blown_down} to move $e$ back to $e_m$, each
  reflection in an effective simple root subtracts a positive multiple of
  that root, and thus $e$ is a sum of effective roots and a formal
  $-1$-curve, so again in $\mu$.  In particular, both $e_1$ and $f-e_1$ are
  in $\mu$, so that $f$ is in $\mu$.  Furthermore, the pullback of a class
  in $\mu(X_{m-1})$ is in $\mu(X_m)$.

  Thus if $D$ is in the dual of $\mu$ and $\alpha$ is a simple root such
  that $D\cdot \alpha<0$, then $\alpha$ is ineffective and we may reflect
  in $\alpha$, and iterating this process will eventually terminate since
  there can be only a bounded number of steps in a row which do not
  decrease $D\cdot f$ and $D\cdot f$ cannot become negative.  We may thus
  assume that $D$ is in the fundamental chamber.

  If $g=0$, we conclude that the dual of the monoid is contained in the
  union over even blowdown structures of the simplicial monoids generated
  by
  \[
  f,s+f,s+f-e_1,2s+2f-e_1-e_2,\dots,2s+2f-e_1-\cdots-e_m.
  \]
  (This is just the dual of the monoid generated by simple roots and
  $e_m$.)  For each $i\ge 2$, the class $2s+2f-e_1-\cdots-e_i$ is certainly
  a sum of components of $Q$ on $X_i$, and thus (being a pullback) is
  in $\mu$.  We have already shown that $f$ is in $\mu$, and since $s-e_1$
  is in the orbit of $f-e_1$, it is also in $\mu$, so that
  $(s-e_1)+(f-e_1)+(e_1)$ and $(s-e_1)+(f-e_1)+2(e_1)$ are in $\mu$.
  It follows that any element of the dual monoid is in $\mu$ as required.

  If $g>0$, we note that $X_0$ (which we again take to be even) has at
  least one component of $Q$ of class $s-df$ for $d\ge 0$, which is again
  in $\mu$ (as a pullback).  We thus conclude that $D\cdot (s-df)\ge 0$,
  telling us that $D$ is an element of the simplicial monoid with
  generators
  \[
  f,s+df,s+df-e_1,2s+df-e_1-e_2,\dots,2s+df-e_1-\cdots-e_m.
  \]
  It is again easy to see that each of these classes is in $\mu$,
  using
  \[
  Q = -K_X=2s+(2-2g)f-e_1-\cdots-e_m
  \]
  for all but the first three generators and the expressions
  \begin{align}
    f &= (f-e_1)+(e_1)\\
    s+df &= (s-df) + d(f-e_1)+d(e_1)\\
    s+df-e_1 &= (s-df)+d(f-e_1)+(d-1)(e_1).
  \end{align}
  If $d=0$, this last case fails, but then $Q$ on $X_0$ has no vertical
  components, and thus the point being blown up must lie on one of the at
  most two components of class $s$, so that $s-e_1$ is in $\mu$.
\end{proof}

\begin{rems}
  This fails if $m=0$, but is easy enough to work around since a class on
  $X_0$ is effective iff its pullback is effective.  We find that the
  effective monoid of $X_0$ is generated by $s'$ and $f$ where $s'$ is the
  effective class of minimal self-intersection such that $s'\cdot f=1$, and
  the only case in which $s'$ is not a formal $-1$- or $-2$-curve or a
  component of $Q$ is when $X_0$ is noncommutative $\P^1\times \P^1$ and
  $s'=s$.  Similarly, the effective monoid of a noncommutative plane is
  the same as the nef monoid, and is generated by the class
  $h:=c_1(\sO_X(1))$.
\end{rems}

\begin{rems}
  It is worth noting that effectiveness as we have defined it is a {\em
    numerical} condition; there will be {\em some} $1$-dimensional sheaf
  with the given Chern class, but we cannot expect to have any particular
  control over the continuous part of its class in $K_0(X)$.  Thus, for
  instance, if $g>1$, then $f$ is effective, but not every line bundle of
  Chern class $f$ has a global section.  (This can also be an issue for
  effective simple roots in general.)
\end{rems}

\begin{cor}
  Any nef divisor class is effective.
\end{cor}

\begin{proof}
  In the course of the above proof, we showed that the dual of the
  effective monoid is contained in the effective monoid.
\end{proof}

It is important to note that some of these generators are redundant.

\begin{prop}
  For $m\ge 1$, the effective monoid is generated by formal $-1$- and
  $-2$-curves and those components $Q_i$ of $Q$ such that $Q_i\cdot Q\le
  0$, unless $g=0$, $m=7$, and $Q$ is integral.
\end{prop}

\begin{proof}
  We need to show that if $Q_i$ is a component of $Q$ such that $Q_i\cdot
  Q>0$, then $Q_i$ is in the stated span.  If $Q_i$ is a vertical component
  with positive intersection with $Q$, then it is either exceptional or the
  strict transform of a fiber, and thus is either $f$ or a $-1$-curve, and
  is redundant in either case.

  If $g>1$ or in the differential case for $g=1$, any horizontal component
  has negative intersection with $Q$.  In the remaining (elliptic
  difference) case with $g=1$, there are two disjoint components,
  $s-df-\sum_{i\notin I} e_i$ and $s+df-\sum_{i\in I} e_i$, and only the
  latter can have positive intersection with $Q$.  We moreover must have
  $|I|<2d$.  The difference of the two components may then be written as
  \[
  (2d-|I|)(f) + \sum_{i\in I} (f-e_i) + \sum_{i\notin I} (e_i),
  \]
  each term of which is a sum of formal $-1$- and $-2$-curves.
  
  We are thus left with the rational case.  When $Q$ is integral, we have
  $Q^2=8-m$, and thus need only establish that $Q$ is redundant for $1\le
  m\le 6$.  For $m=1$, this is easy: $2s+2f-e_1=2(f-e_1)+2(s-e_1)+e_1$,
  while for $2\le m\le 6$, there is a commutative rational surface with
  anticanonical curve an $8-m$-gon of $-1$-curves, again giving an
  expression of $-K_X$ as a sum of elements in the orbit of $e_m$.

  When $Q$ is not integral, its components are rational, and $Q_i\cdot
  Q>0$ implies $Q_i^2\ge -1$.  If $Q_i^2=-1$, then it is a formal
  $-1$-curve, so we may assume $Q_i^2\ge 0$.  As in the proof of
  Proposition \ref{prop:formal_-1_can_be_blown_down}, we may assume that
  $Q_i$ is in the fundamental chamber.  We then find (by the analogous
  inequality) that $Q_i\cdot Q-2=Q_i^2\ge Q_i\cdot Q$ unless $Q_i\cdot
  e_2=0$, so that we reduce to the case $m=1$.  The generators of the
  relevant monoid are contained in the monoid generated by $e_1$, $f-e_1$,
  $s-e_1$, so the claim follows.
\end{proof}

\begin{rem}
  As in the commutative case, when $g=0$, $m=7$ and $Q$ is irreducible, the
  generator $Q$ is {\em nearly} redundant, since $2Q$ {\em can} be
  expressed as a sum of $-1$-curves on a suitable commutative surface.
\end{rem}

Our understanding of the effective monoid also gives us the following
stronger version of the finitely partitive property.

\begin{cor}
  For any divisor class $D$, there are only finitely many $D'$ such
  that both $D'$ and $D-D'$ are effective.
\end{cor}

\begin{proof}
  Since any formal $-1$- or $-2$-curve can be expressed as a nonnegative
  linear combination of simple roots and $e_m$, we find that the effective
  monoid is contained in a finitely generated monoid (generated by
  components of $Q$, simple roots, and $e_m$).  Moreover, all of the
  generators are lexicographically positive (relative to the ordering $s$,
  $f$, $e_1$,\dots, $e_m$ of the standard basis), and thus any linear
  dependence between the generators has two coefficients of opposite sign.
  It follows that the dual cone has dimension $m+2$, and thus its interior
  contains integer elements.  An element of the interior induces a grading
  on the larger monoid with respect to which it is finite and thus there
  are only finitely many elements of the monoid of degree strictly smaller
  than $D$.
\end{proof}

\begin{rem}
It is worth noting that this result continues to hold even if we allow the
surface to vary over a Noetherian base, and simply ask for $D'$ and $D-D'$
to be effective on {\em some} surface in the family.  Indeed, there are
still only finitely many possible divisor classes of components of $Q$
(this reduces to a statement about families of commutative surfaces),
giving a finitely generated monoid containing the effective monoid of any
fiber and having dual of full dimension.
\end{rem}

As for commutative rational surfaces, we can adapt the above arguments to
give an essentially combinatorial algorithm for determining whether a class
is effective (resp. nef).  Indeed, if we {\em start} by reducing to the
fundamental chamber, each step either just changes the blowdown structure
(and possibly reduces $D\cdot f$) or exhibits a formal $-2$-curve or
component of $Q$ having negative intersection with $D$, which we may
subtract without changing whether $D$ is effective.  Thus after a finite
number of steps, we will either have put $D$ into the fundamental chamber
or made $D\cdot f$ negative (in which case it was not effective).  If the
result has nonnegative intersection with every component of $Q$, then it is
effective, and otherwise, we may subtract the offending component and
rereduce to the fundamental chamber as necessary.  Each such step subtracts
a nontrivial effective divisor from $D$, and thus again we can only do this
finitely many times before eventually reaching a divisor having negative
intersection with some element of the interior of the dual of the cone
generated by simple roots, $e_m$ and components of $Q$.

To determine if $D$ is nef, the algorithm is even simpler: any time we
would have wanted to subtract an effective divisor class from $D$, we
simply terminate saying that $D$ was not nef.

\section{Ample divisors}

Having mostly settled the question of what the effective and nef divisor
classes look like, we now turn to understanding {\em ample} divisors.
Since our divisor classes are only numerical, the two things we would like
an ample divisor $D_a$ to satisfy are (a) that for any sheaf $M$, there is
an integer $B$ such that for any line bundle $L$ of Chern class $-bD_a$
with $b\ge B$, $\Ext^p(L,M)=0$ for $p>0$, and (b) that there is a similar
bound guaranteeing that $L\otimes \Hom(L,M)\to M$ is surjective.  (If both
hold for a given pair $(L,M)$, we say that $M$ is {\em acyclically globally
  generated} by $L$.)  We will not entirely characterize the divisor
classes satisfying these conditions (there are some technical issues coming
from the fact that there can be $1$-dimensional sheaves disjoint from $Q$
without proper subsheaves), but will show that anything in the interior of
the nef cone will do.  (It is also easy to see that (a) and (b) fail for
divisors which are not nef, so the only issues arise on the boundary.)

A key observation is that global generation for general sheaves follows if
we can establish global generation for line bundles; there is also a more
subtle reduction to line bundles for acyclic global generation (see below).
Since the only thing we know about $D_a$ is certain numerical inequalities,
we are led to look for a statement along the following lines: if $L$ and
$L'$ are line bundles such that $c_1(L)-c_1(L')$ satisfies certain
inequalities, then $L$ is acyclically globally generated by $L'$.  The
presence of an anticanonical curve turns out to be extremely helpful in
this regard; the condition will essentially be that $c_1(L)-c_1(L')$ is nef
and that the ratio of their restrictions to $Q$ is acyclic and globally
generated.  The latter condition is somewhat different depending on whether
$X$ is rational or has genus $\ge 1$.  Note that in either case $Q$ may be
embedded as an anticanonical curve on a commutative ruled surface.

\begin{lem}
  Let $Q$ be an anticanonical curve on a (commutative) rational surface and
  let $L$ be a line bundle having nonnegative degree on every component of
  $Q$.  Then $L$ is acyclic if $\deg(L)\ge 1$ and globally generated if
  $\deg(L)\ge 2$.
\end{lem}

\begin{proof}
  When the ambient surface is $\P^2$ or $\P^1\times \P^1$, this was shown
  in \cite[\S7]{ArtinM/TateJ/VandenBerghM:1991}.  In fact, the only way in
  which the ambient surface was used in those proofs was in showing two
  numerical connectivity properties of $Q$: first, that if $Q=A+B$ with
  $A,B$ nonzero effective divisors on the ambient surface, then $A\cdot
  B>0$, and second that if $Q=A+B+C$ with $A,B$ nonzero effective and $C$ a
  smooth irreducible component of $Q$, then $A\cdot B>0$.  The first is
  just Lemma \ref{lem:Ca_is_num_conn}, which further showed that $A\cdot B$
  is even.  Thus in the latter situation, we find that $A\cdot (B+C)$ and
  $B\cdot (A+C)$ are both positive even integers, and thus their sum
  \[
  2(A\cdot B) + (A+B)\cdot C\ge 4.
  \]
  Since $C$ is a smooth rational curve, $(A+B)\cdot C=2$, and thus $A\cdot
  B\ge 1$ as required.
\end{proof}

\begin{prop}
  Let $X$ be a noncommutative rational surface, and let $L$, $L'$ be line
  bundles on $X$ such that $D=c_1(L)-c_1(L')$ is nef.  If $D\cdot Q\ge 1$,
  then $\Ext^p(L',L)=0$ for $p>0$, and if $D\cdot Q\ge 2$, then
  $\Hom(L',L)\otimes L'\to L$ is surjective.
\end{prop}

\begin{proof}
  Since the nef cone is invariant under twisting $X$ by a line bundle, we
  may as well assume that $L'\cong \sO_X$, so that we need to show that for
  any nef divisor $D$ and line bundle $L$ of class $D$, $L$ is acyclic if
  $D\cdot Q\ge 1$ and globally generated if $D\cdot Q\ge 2$.  If $m\ge 1$,
  choose a blowdown structure on $X$ such that $D$ is in the fundamental
  chamber.  In particular, $D$ is in the monoid generated by
  \[
  f,s+f,s+f-e_1,2s+2f-e_1-e_2,\dots,2s+2f-e_1-\cdots-e_m.
  \]
  Acyclicity and global generation are inherited by pullbacks, so if $D$ is
  a pullback (i.e., if $D\cdot e_m=0$) then we may reduce to $X_{m-1}$ and
  induct.  This also applies when $m=1$, except that now the reduction
  applies if either $D\cdot e_1$ or $D\cdot (f-e_1)$ is 0, and we may need
  to perform an elementary transformation before blowing down.  Similarly,
  for the odd $m=0$ case, if $D\cdot s=0$, then $s$ must be a $-1$-curve
  and we can reduce to $m=-1$.

  Thus suppose that $m\ge 1$ and $D\cdot e_m>0$ (with $D\cdot (f-e_1)>0$ if
  $m=1$).  To show that $L$ is acyclic (resp. and globally generated), it
  will suffice to show that $L(-Q)$ and $L|_Q$ are acyclic (resp. and
  globally generated).  The latter holds since $L|_Q$ has nonnegative
  degree on every component (since $L$ is nef) and has total degree $D\cdot
  Q$.  For the former, we need to show that $D-Q$ is nef and either
  vanishes or satisfies the same inequality.  The result is certainly still
  in the fundamental chamber, so that for nefness we need only show that it
  has nonnegative intersection with the components of $Q$ having
  nonnegative intersection with $Q$.  If $Q$ is reducible and $Q_i$ is such
  a component, then $(D-Q)\cdot Q_i=D\cdot Q_i-(Q\cdot Q_i)\ge D\cdot Q_i$
  as required.  It thus remains only to consider when $(D-Q)\cdot Q\ge 1$
  (or $2$, as appropriate).  If $Q^2=8-m\le 0$, there is again no
  difficulty, while if $m<7$, then any generator of the standard monoid has
  intersection $\ge 2$ with $Q$, and $\sO_X$ is certainly acyclic and
  globally generated.

  The only case in which $D-Q$ fails to satisfy the requisite inequality is
  when $m=7$ and $D=2Q$.  In that case, we still have acyclicity, and need
  only establish global generation.  One global section of $\sO_X(2Q)$ is
  easy to understand: simply take the composition of the maps $\sO_X\to
  \sO_X(Q)$ and $\sO_X(Q)\to \sO_X(2Q)$ coming from the anticanonical
  natural transformation.  The cokernel of this global section is an
  extension
  \[
  0\to \sO_X(Q)|_Q\to \sO_X(2Q)/\sO_X\to \sO_X(2Q)|_Q\to 0.
  \]
  Since $\sO_X$ is acyclic and globally generated, $\sO_X(2Q)$ is globally
  generated iff $\sO_X(2Q)/\sO_X$ is globally generated, and since this is
  an extension of sheaves on $Q$ with acyclic kernel, it is globally
  generated iff its restriction to $Q$ is globally generated.  But
  $\sO_X(2Q)|_Q$ is globally generated by degree considerations.

  It remains only to consider the base cases $m=0$ and $m=-1$.  Here, the
  relevant inequalities force $D=0$ or $D\cdot Q\ge 2$, so we always want
  acyclic global generation.  For $m=-1$, any effective line bundle is
  acyclically globally generated, essentially by construction.  For $m=0$,
  the effective monoid is generated by $f$ and either $s$ or the most
  negative horizontal component of $Q$, so if $D\cdot f\ge 2$, then $D-Q$
  will again be nef and either be 0 or have sufficiently large degree on
  $Q$.  For $D\cdot f=1$, acyclicity and global generation reduce to the
  corresponding properties for ${\cal E}$, so again reduce to what happens
  on $Q$, while for $D\cdot f=0$, so $D\propto f$, we find that $\sO_X(D)$
  is the pullback of an acyclic and globally generated sheaf on $\P^1$.
\end{proof}

The argument for higher genus surfaces is analogous, the main difference
being in the specific numeric condition along $Q$.

\begin{lem}
  Let $Q$ be an anticanonical curve on a commutative rationally ruled
  surface of genus $g\ge 1$ (with $Q$ not integral), and let $L$ be a line
  bundle having nonnegative degree on every component of $Q$.  Then $L$ is
  acyclic if it has degree at least $2g-1$ on the horizontal component(s)
  of $Q$ and globally generated if it has degree at least $2g$ on the
  horizontal component(s).
\end{lem}

\begin{proof}
  Since $0$-dimensional sheaves are acyclic and globally generated and both
  ``acyclic'' and ``acyclic and globally generated'' are closed under
  extensions, we may feel free to replace $L$ by any line bundle it
  contains that still satisfies the numerical conditions.  In particular,
  we may as well assume that $L$ has degree $0$ on every vertical component
  of $Q$.  But then we may as well contract that component before checking
  the desired conditions, and we may achieve this effect by blowing down to
  a ruled surface and performing a suitable sequence of elementary
  transformations.  We may thus reduce to the case that $Q$ has only
  horizontal components.  If $Q$ is reduced, then $g=1$, $Q$ is smooth, and
  $L$ is clearly acyclic (resp. globally generated) on each component.
  If $Q$ is nonreduced, then we have a short exact sequence
  \[
  0\to L'\to L\to L|_{Q^{\red}}\to 0
  \]
  where $L'$ is a line bundle on $Q^{\red}\cong C$.  That $L|_{Q^{\red}}$
  is acyclic (resp. acyclic and globally generated) follows by standard
  facts about line bundles on smooth curves, and $L'$ has degree at least
  that of $L|_{Q^{\red}}$, so is also acyclic.
\end{proof}

\begin{rem}
  Since a component of $Q$ is vertical iff $Q_i\cdot f=0$ and horizontal
  iff $Q_i\cdot f=1$, we may rephrase the numerical condition for
  $\sO_X(D)|_Q$ as saying that $D-(2g-1)f$ (resp. $D-2gf$) has nonnegative
  intersection with every component of $Q$.  Since any formal $-1$- or
  $-2$-curve on a higher genus rationally ruled surface is orthogonal to
  $f$, a divisor $D$ is nef with $\sO_X(D)|_Q$ satisfy the relevant
  condition iff $D-(2g-1)f$ (resp. $D-2gf$) is nef.
\end{rem}

\begin{prop}
  Let $X$ be a noncommutative rationally ruled surface of genus $g\ge 1$,
  and let $L$, $L'$ be line bundles on $X$.  If $c_1(L)-c_1(L')-(2g-1)f$ is
  nef, then $\Ext^p(L',L)=0$ for $p>0$, and if $c_1(L)-c_1(L')-2gf$ is nef,
  then $L$ is acyclically globally generated by $L'$.
\end{prop}

\begin{proof}
  We may again assume that $L'=\sO_X$ and choose a blowdown structure such
  that $D=c_1(L)-(2g-1)f$ (resp. $c_1(L)-2gf$) is in the fundamental
  chamber with $D\cdot e_m\ge 1$ if $m\ge 1$ (and $D\cdot (f-e_1)\ge 1$ if
  $m=1$).  For $m\ge 1$, there is no difficulty since if $D$ is nef with
  $D\cdot e_m>0$, then $D-Q$ is again nef, and thus we readily reduce to
  $m=0$.  In that case, we again find that $D-Q$ remains nef as long as
  $D\cdot f\ge 2$, and thus reduce to the cases $D\cdot f=1$ and $D\cdot
  f=0$.  The latter is immediate (since then $L$ is the pullback of an
  acyclic or acyclic and globally generated line bundle on $C$), while the
  former reduces to the corresponding claim for ${\cal E}$ and thus for
  $Q$.
\end{proof}

\begin{thm}
  Let $X$ be a noncommutative rational or rationally ruled surface, and let
  $D_a$ be an integral element of the interior of the nef cone.  Then for
  any sheaf $M$ on $X$, there is an integer $B$ such that for any line
  bundle $L$ of Chern class $-bD_a$ with $b\ge B$, $M$ is acyclically
  globally generated by $L$.
\end{thm}

\begin{proof}
  There is certainly {\em some} line bundle $L_0$ that acyclically globally
  generates $M$, and a line bundle $L_1$ that acyclically globally
  generates the kernel, etc., giving us a resolution
  \[
  \cdots\to L_1^{n_1}\to L_0^{n_0}\to M.
  \]
  For any other line bundle $L$, we may use this resolution to compute
  $R\Hom(L,M)$, and find that the only possible contributions to
  $\Ext^p(L,M)$ for $p\ge 0$ come from $\Ext^2(L,L_1)$ and $\Ext^p(L,L_0)$
  for $p\in \{1,2\}$.  Thus if $c_1(L_1)-c_1(L)$ and $c_1(L_0)-c_1(L)$ are
  nef divisors satisfying the appropriate additional inequality, then
  $\Ext^p(L,M)=0$ for $p>0$.  Similarly, if $L_0$ is globally generated by
  $L$, then so is $M$.  Since $D_a$ is in the interior of the nef cone, it
  has positive intersection with every effective class, and thus there is
  some $B$ such that for $b\ge B$, $bD_a+c_1(L_0)$ and $bD_a+c_1(L_1)$ are
  both nef and have intersection at least $2$ with $Q$ or $2g$ with the
  horizontal components of $Q$ as appropriate.
\end{proof}

\begin{rem}
  Since we have shown that the nef cone contains a subcone of full
  dimension, there are indeed ample divisors.
\end{rem}

\begin{cor}
  If we choose for each $b\in \Z$ a line bundle of Chern class $bD_a$, then
  $X$ may be recovered as the $\Proj$ of the corresponding full subcategory
  of $X$, viewed as a $\Z$-algebra.
\end{cor}

Now, suppose that $X/S$ is a family of noncommutative surfaces over some
Noetherian base $S$.  Let $Y/S$ be the corresponding family of commutative
surfaces, and let $D$ be an ample divisor on that family.  Since $D^2>0$,
it follows that $\langle D,Q\rangle^\perp$ is negative definite, and thus
that it meets the root system $E_{m+1}$ or $D_m$ (as appropriate) in a
finite subsystem.  If we first base change so that $Y$ has a blowdown
structure, then change the blowdown structure so that $D$ is in the
fundamental chamber, then $\langle D,Q\rangle^\perp$ meets the ambient root
system in a finite {\em parabolic} subsystem, and it follows that the set
of such blowdown structures is a torsor over the corresponding finite
Coxeter group.  (None of these roots can be effective on $Y$, since $D$ is
ample.)  It follows that there is a finite \'etale (and Galois) cover of
$S$ over which $Y$ has a blowdown structure.  Over this larger family, we
may choose a new divisor $D'$ which is not only ample on every fiber but in
the interior of the fundamental chamber.  Although this divisor class
itself will typically not descend to $X$, we note that since the
$t$-structure is invariant under the action of monodromy, the monodromy
preserves the ample cone.  Thus the trace (in the obvious sense) of $D'$
will still be ample, and now descends to a class in $\NS(X)$ over $S$.
Moreover, any blowdown structure on $Y$ that puts $\Tr(D')$ into the
fundamental chamber will give rise to the same $t$-structure on $X$.  The
significance of this is that $\Tr(D')$ makes sense even without first
assuming that $\perf(X)$ has a rational $t$-structure, and thus we may use
it to define such a $t$-structure whenever it exists.

In fact, we may use such a divisor to construct a model of $X$ as the
$\Proj$ of a sheaf of $\Z$-algebras on $S$.  Let $D_a\in \NS(X)$ be a
divisor which is ample on every fiber.  The main difficulty is that we need
to choose a line bundle of each Chern class which is a multiple of $D_a$,
and there may not be any such line bundle over $S$.  The main issues are
that the corresponding $\Pic^0_{C/S}$-torsor could fail to have points over
$S$, and that a point over $S$ need not correspond to an actual line
bundle.  For the first issue, we note that the $\Pic^0_{C/S}$-torsor lies in
an appropriate component of $\Pic_{Q/S}$, and is in the image of the natural
map
\[
\Gamma(S;\Pic_{Q/S}/\Pic^0_{C/S})\to H^1(S;\Pic^0_{C/S}).
\]
The element in $\Gamma(S;\Pic_{Q/S}/\Pic^0_{C/S})$ can be computed as the
image under the determinant-of-restriction map $K_0^{\num}(X)\to
\Pic_{Q/S}/\Pic^0_{C/S}$ of the numeric class of a line bundle of Chern
class $bD_a$.  Since the class of such a line bundle depends quadratically
on $b$, we obtain a quadratic map $\Z\to \Gamma(S;\Pic_{Q/S}/\Pic^0_{C/S})$
and thus a quadratic map $\Z\to H^1(S;\Pic^0_{C/S})$.  Moreover, since each
component of $\Pic_{Q/S}$ is quasi-projective, so is each
$\Pic^0_{C/S}$-torsor, and thus the $\Pic^0_{C/S}$-torsors are projective.
This implies that the resulting classes in $H^1(S;\Pic^0_{C/S})$ are
torsion, and thus that we can make the map to $H^1(S;\Pic^0_{C/S})$ trivial
by replacing $D_a$ by a sufficiently divisible positive multiple.  Moreover,
we find that not only does the map $\Z\to
\Gamma(S;\Pic_{Q/S}/\Pic^0_{C/S})$ map to the image of $\Pic_{Q/S}(S)$, but
we may choose a lifting of the map which is still quadratic.  (Indeed, we
may simply choose arbitrary preimages for the images of $\pm 1$, and take
the image of $0$ to be trivial.)  The obstruction to being able to lift
this map from $\Pic_{Q/S}(S)$ to $\Pic(Q)$ is a quadratic map $\Z\to
H^2(S;G_m)$, and again the image is torsion since $Q$ is projective.  So
multiplying $D_a$ by another sufficiently divisible positive integer makes
the Brauer obstruction vanish and allows us to lift.  We thus obtain a
system of line bundles ${\cal L}_b$, $b\in \Z$ satisfying ${\cal L}_b\cong
{\cal L}_1^{b(b+1)/2}\otimes {\cal L}_{-1}^{b(b-1)/2}$ and such that for
each $b\in \Z$, there is a line bundle $L_b$ on a suitable base change of
$X$ which has Chern class $bD_a$ and satisfies $L_b|_Q\cong {\cal L}_b$.
Moreover, by including this isomorphism, we rigidify $L_b$ and guarantee
that it descends to $S$.  (Indeed, $L_b$ is unique up to automorphisms, and
its automorphisms are scalars, so change the isomorphism to ${\cal L}_b$.)
The full subcategory of $\qcoh X$ with objects $L_b$, $b\in \Z$ is a sheaf
of $\Z$-algebras on $S$ and every fiber of $X$ agrees with the $\Proj$ of
the corresponding fiber of this sheaf of $\Z$-algebras.

By mild (but very suggestive) abuse of notation, we will denote the line
bundles $L_b$ so constructed by $\sO_X(bD_a)$.  The reader should be
cautioned that these sheaves do not descend to the moduli stack of
noncommutative surfaces, as their construction required some additional
choices.

\section{Computing Betti numbers}

In the commutative case, one normally says that a line bundle is effective
iff it has a global section.  This does not quite agree as stated with our
definition above, but leads one to ask not only when a line bundle has
global sections, but how many global sections it has.  Although it is
unclear how to answer this question for general ruled surfaces (it contains
the analogous problem for smooth projective curves as a special case!), we
can give a complete answer in the rational case, and in fact can give an
essentially combinatorial algorithm to compute the dimensions of Ext groups
between any two line bundles.  (In other words, we have a noncommutative
analogue of the algorithm for computing dimensions of global sections of
line bundles on anticanonical rational surfaces given at the end of Section
\ref{sec:effdivs}.)

Let $D_1,D_2\in \Pic^{\num}(X)$ be a pair of divisor classes, and suppose
that we wish to compute $\dim \Ext^i(\sO_X(D_1),\sO_X(D_2))$ for $i\in
\{0,1,2\}$.  We know the alternating sum of these dimensions, and thus it
suffices to compute $\dim\Hom(\sO_X(D_1),\sO_X(D_2))$ and
$\dim\Ext^2(\sO_X(D_1),\sO_X(D_2))=\dim\Hom(\sO_X(D_2),\theta
\sO_X(D_1))$.  If $D_2-D_1$ is not effective, then
$\dim\Hom(\sO_X(D_1),\sO_X(D_2))=0$ while if $K_X+D_1-D_2$ is not
effective, then $\dim\Ext^2(\sO_X(D_1),\sO_X(D_2))=0$; since $K_X$ is
ineffective, at least one of these must be the case.  We thus reduce to the
question of computing $\dim\Hom(\sO_X(D_1),\sO_X(D_2))$ when $D_2-D_1$
is effective.  This in turn reduces to computing $\dim\Gamma(\sO_{X'}(D_2-D_1))$
where $X'$ is a suitable twist of $X$, so that we may as well take $D_1=0$,
$D_2=D$.

As in the algorithms above, a key step is to reduce to the case that $D$ is
nef and in the fundamental chamber.  If there is a component $Q_1$ of $Q$
such that $D\cdot Q_1<0$, then $\sO_{Q_1}(D)$ has no global
sections, and thus
\[
\Gamma(\sO_X(D-Q_1))\cong \Gamma(\sO_X(D)),
\]
allowing us to replace $D$ by $D-Q_1$.  We may then iterate this until $D$
has nonnegative inner product with every component of $Q$.  (Since $D$ is
assumed effective, this process terminates just as in the algorithm for
testing whether $D$ is effective.)

The same argument lets us subtract $e_m$ from $D$ whenever $D\cdot e_m<0$,
and thus it remains to consider the case that $D$ has negative intersection
with a simple root.  If the simple root is ineffective, then we may apply
the corresponding reflection, but this of course fails when the simple root
is effective.  Suppose that the simple root is $e_i-e_{i+1}$.  We have
already disposed of the possibility that $e_i-e_{i+1}$ is an anticanonical
component on $X_{i+1}$, and thus the only way it can be effective is when
the two points being blown up are smooth points of the anticanonical curve
in the same orbit.  We may then for simplicity twist by $e_{i+1}$ to ensure
that the points are the same, at the cost of changing the required
computation to $\dim\Hom(\sO_X(le_{i+1}),\sO_X(D+le_{i+1}))$ for some $l$.
If $q$ is non-torsion, then there is a unique such $l$, while if $q$ is
torsion of order $r$, we choose $l$ so that $-r\le (D+le_{i+1})\cdot
(e_i-e_{i+1})<0$.

If $(D+le_{i+1})\cdot (e_i-e_{i+1})\ge 0$, so that $q$ is non-torsion and
$l=(le_{i+1})\cdot (e_i-e_{i+1})>0$, then we may apply Proposition
\ref{prop:weak_reflection} to suitable twists to obtain short exact
sequences
\[
0\to \sO_X(le_{i+1}) \to \sO_X(le_i) \to \sO_{e_i-e_{i+1}}(-1)^l \to 0
\]
and
\[
0\to \sO_X(D+le_{i+1})
 \to \sO_X(s_i(D)+le_i)
 \to \sO_{e_i-e_{i+1}}(-1)^{(D+le_{i+1})\cdot (e_i-e_{i+1})}
 \to 0,
\]
where $s_i(D)$ is the image of $D$ under the reflection.  We may thus
compute the spherical twists of $\sO_X(le_{i+1})$ and $\sO_X(D+le_{i+1})$,
allowing us to conclude that
\[
\Hom(\sO_X(le_{i+1}),\sO_X(D+le_{i+1}))
\cong
\Hom(\sO_X(le_i),\sO_X(s_i(D)+le_i)),
\]
so that $\sO_{X'}(D)$ and $\sO_{s_i(X')}(s_i(D))$ have isomorphic spaces of
global sections, where $X'$ is the appropriate twist of $X$.  Similarly, if
$l\le 0$, then either $q$ is non-torsion or $-r<l\le 0$, so that we may use
Proposition \ref{prop:weak_reflection} to compute the opposite spherical
twists and again conclude that
\[
\Hom(\sO_X(le_{i+1}),\sO_X(D+le_{i+1}))
\cong
\Hom(\sO_X(le_i),\sO_X(s_i(D)+le_i)).
\]

The remaining case is that $l>0$ but $(D+le_{i+1})\cdot (e_i-e_{i+1})<0$.
In that case, we still have a short exact sequence
\[
0\to \sO_X(s_i(D)+le_i)
 \to \sO_X(D+le_{i+1})
 \to \sO_{e_i-e_{i+1}}(-1)^{-(D+le_{i+1})\cdot (e_i-e_{i+1})}
 \to 0,
\]
but now have
\[
\Hom(\sO_X(le_{i+1}),\sO_{e_i-e_{i+1}}(-1))=0,
\]
so that
\[
\Hom(\sO_X(le_{i+1}),\sO_X(s_i(D)+le_i))
\cong
\Hom(\sO_X(le_{i+1}),\sO_X(D+le_{i+1})).
\]
Setting $D'=s_i(D)+le_i-le_{i+1}$, we find that $D-D'$ is a positive
multiple of $e_i-e_{i+1}$, so that this reduces $-D\cdot (e_i-e_{i+1})$.

We thus find that in each case such that $D\cdot (e_i-e_{i+1})<0$ with
$e_i-e_{i+1}$ effective, we have an isomorphism $\Gamma(\sO_X(D))\cong
\Gamma(\sO_{X'}(D'))$ where $D-D'$ is a positive multiple of $e_i-e_{i+1}$.
Similar calculations also hold for the simple roots $f-e_1-e_2$ (related to
$e_1-e_2$ by an elementary transformation) and $s-f$.  Thus in general if
$D$ has negative intersection with a simple root, we may reduce the problem
to the corresponding calculation on a related surface with smaller
(relative to the effective cone) $D$, sandwiched
between the original $D$ and its image under the reflection.  In
particular, we can perform only finitely many such reductions before ending
up in the fundamental chamber.

It remains to determine what happens when $D$ is nef and in the fundamental
chamber; we may also assume that $D\cdot e_m>0$ as otherwise we may reduce
to the blown down surface.  If $D=0$ or $D\cdot Q>0$, then $\sO_X(D)$ is
acyclic, and thus we have $\dim\Gamma(\sO_X(D))=
\chi(\sO_X(D))=1+(D\cdot (D+Q))/2$.  It remains only to consider what
happens when $D\ne 0$ with $D\cdot Q=0$.  We again have a short exact
sequence
\[
0\to \sO_X(D-Q)\to \sO_X(D)\to \sO_X(D)|_Q\to 0.
\]
If $\sO_X(D)|_Q$ is nontrivial (so has no global sections) or $\sO_X(D-Q)$
is acyclic, then this induces a short exact sequence
\[
0\to \Gamma(\sO_X(D-Q))\to \Gamma(\sO_X(D))\to \Gamma(\sO_X(D)|_Q)\to 0,
\]
again allowing us to reduce the dimension calculation.

At this point, the only way $D$ can avoid all of our reductions is to have
$m=8$, $Q^2=0$ and $D\propto Q$, such that any rational component of $Q$
has self-intersection $-2$.  We thus need to compute $\Gamma(\sO_X(lQ))$
for $l>0$, where $Q$ is of this form and $\sO_Q(lQ)\cong \sO_Q$.  If $l$ is
minimal with this property, then $\sO_X((l-1)Q)$ is acyclic, so that we may
compute
\[
\dim\Gamma(\sO_X(lQ))=\chi(\sO_X((l-1)Q))+1=2.
\]
Since these global sections are generically nonzero on $Q$, but there is a
1-dimensional subspace vanishing on $Q$, they are algebraically
independent, and thus we conclude that
\[
\dim\Gamma(\sO_X(alQ))\ge a+1
\]
for any integer $a\ge 0$.  The long exact sequence tells us that
\[
\dim\Gamma(\sO_X((a+1)lQ))\le 1 + \dim\Gamma(\sO_X(((a+1)l-1)Q))
= 1 + \dim\Gamma(\sO_X(alQ)),
\]
from which it follows that $\dim\Gamma(\sO_X(alQ))=a+1$, finishing off our
calculation of dimensions of global sections.

\section{Representations as graded algebras}

Although the construction of $X$ via the $\Z$-algebra associated to an
ample divisor is in principle explicit, it is still not the nicest model
one might hope for.  In the ideal situation, we could represent $X$ as the
$\Proj$ of an actual graded algebra.  This is too much to hope for in
complete generality, as it requires us $\coh X$ to have an ample
autoequivalence, and this is not the case in general.  There is one case in
which can guarantee such an autoequivalence, namely when $X$ is Fano, i.e.,
when $Q$ is ample.  More generally, as long as $Q$ is nef with $Q^2\ge 1$, we
can hope to obtain a graded algebra deforming the anticanonical coordinate
ring of a degenerate del Pezzo surface.  (Such rings may be of particular
interest as models of noncommutative surfaces with rational singular
points.)

This leads us to the question of what form the presentation of the
anticanonical graded algebra $\bigoplus_i \Gamma(\theta^{-i}\sO_X)$ takes.
A key observation is that there is a central element of degree $1$
corresponding to the anticanonical natural transformation, and the quotient
by the corresponding principal ideal is the graded algebra
\[
\bigoplus_i \Gamma(\theta^{-i}\sO_X|_Q),
\]
which is essentially a twisted homogeneous coordinate ring of $Q$, except
that if $Q$ has components of self-intersection $-2$, its $\Proj$ will
contract those components.  Indeed, for any $i$, $\theta^{-i}\sO_X|_Q$ has
degree 0 on any such component, and thus we may as well work over the
contracted curve.  Moreover, we may think of the anticanonical ring itself
as a filtered deformation of this twisted homogeneous coordinate ring
(i.e., a filtered ring with that associated graded); in particular,
apart from the additional generator of degree 1 and the relations that it
be central, the structure of the presentation is the same.

The case of degree $1$ del Pezzo surfaces is particularly nice, as in that
case the contracted curve is just the Weierstrass model of $Q$.  In general,
let $L_i$ be the line bundle on $Q$ given by $\theta^{-i}\sO_X|_Q$.

\begin{lem}
  If $Q$ is nef with $Q^2=1$ and $q$ is nontrivial, then the anticanonical
  ring is a filtered deformation of a ring generated by elements of degree
  $1$ and $2$, with relations of degree $4$, $5$, and $6$.
\end{lem}

\begin{proof}
  This is really just a statement about the twisted homogeneous coordinate
  ring of a Weierstrass curve relative to a degree 1 line bundle and a
  translation $q$.  We first observe that there is an exact sequence
  \[
  0\to 
  L_{-6}\to
  L_{-5}\oplus L_{-4}\to
  L_{-2}\oplus L_{-1}\to \sO_Q\to 0
  \]
  of sheaves on $Q$, which may be obtained as the Yoneda product of
  \[
  0\to L_{-6}\to L_{-5}\oplus L_{-4}\to L_{-2}\otimes L_{-1}\to 0
  \]
  and
  \[
  0\to L_{-2}\otimes L_{-1}\to L_{-2}\oplus L_{-1}\to \sO_Q\to 0
  \]
  Here the maps $L_{-1}\to \sO_Q$ and $L_{-6}\to L_{-5}$ come from the
  degree $1$ generator $v_1$ (which is unique modulo scalars) and similarly
  for the maps $L_{-2}\to \sO_Q$ and $L_{-6}\to L_{-4}$, and the
  identification of $L_{-2}\otimes L_{-1}$ is via determinant
  considerations.  It follows in particular that the algebra is generated
  by the elements of degree $1$ and $2$ coming from this exact sequence.

  The correspondence between sheaves and graded modules over the twisted
  homogeneous coordinate ring is by applying $\theta^{-i}$ and taking
  global sections.  The only degrees in which the corresponding sequence of
  graded modules is not exact are those for which some terms have
  nontrivial $H^1$ (and some terms have global sections), and thus only
  $0\le i\le 6$ need be considered.  For $0<i$, surjectivity at $\sO_Q$ is
  straightforward (given $q$ nontrivial), and implies exactness for $i\in
  \{1,2,3\}$.  We dually have exactness for $i\in \{3,4,5\}$, so that only
  $i=0$ and $i=6$ remain.  For $i=0$, we clearly have cohomology $0,0,0,k$,
  while for $i=6$ the fact that the spectral sequence converges to $0$ lets
  us read off the cohomology of global sections from the cohomology of
  $H^1$, and thus we obtain $0,0,k,0$.  But this implies that the relations
  are generated by the relations of degree $4$ and $5$ coming from the
  complex, except that there is one missing relation of degree $6$.
\end{proof}

\begin{rem}
  In the commutative case, we need an additional generator and relation of
  degree $3$, and the only relation which is not a commutator is the
  relation of degree $6$.
\end{rem}

For degree 2 del Pezzo surfaces, the situation is slightly more complicated
since the contracted curve may be reducible.  It is always a double cover
of $\P^1$, however, and this makes things tractable.

\begin{lem}
  If $Q$ is nef with $Q^2=2$ and $q$ nontrivial, then the anticanonical
  ring is a filtered deformation of a ring generated by two elements of
  degree $1$ with relations of degree $3$, $3$, and $4$.
\end{lem}

\begin{proof}
  Analogous, except that now the relevant exact sequence has the form
  \[
  0\to L_{-4}\to L_{-3}^2\to L_{-1}^2\to \sO_Q\to 0,
  \]
  again obtained as a Yoneda product.
\end{proof}

\begin{rem}
  Again, we need an additional generator and relation of degree $2$ in the
  commutative case, with the only nontrivial relation that of degree $4$.
\end{rem}

For $Q^2\ge 3$, the algebra is generated in degree $1$, as an easy
consequence of the following result.

\begin{lem}
  Let $Q$ be a nef anticanonical curve on a rational surface, and suppose
  ${\cal L}$, ${\cal L}'$ are invertible sheaves on $Q$ such that ${\cal
    L}$ and ${\cal L}'\otimes {\cal L}^{-1}$ have nonnegative degree on
  every component and ${\cal L}$ has total degree at least 3.  Then
  $\Gamma({\cal L})\otimes \Gamma({\cal L}')\to \Gamma({\cal L}\otimes
  {\cal L}')$ is surjective.
\end{lem}

\begin{proof}
  We have already seen that ${\cal L}$ is globally generated.  It follows
  that for any $0$-dimensional subscheme of $Q$, there is a section of
  ${\cal L}$ (over the algebraic closure) surjective along that
  subscheme. In particular, applying this to a $0$-dimensional subscheme
  hitting every component gives us a section $f_1$ with $0$-dimensional
  cokernel, and applying it to $\supp(\coker(f_1))$ gives a section $f_2$
  such that ${\cal L}$ is generated by $f_1$ and $f_2$.

  Using these sections, we may construct a short exact sequence
  \[
  0\to {\cal L}'\otimes {\cal L}^{-1}
  \to^{(f_2,-f_1)}
  {\cal L}'\oplus {\cal L}'
  \to^{(f_1,f_2)}
  {\cal L}'\otimes {\cal L}
  \to
  0.
  \]
  If ${\cal L}'\otimes {\cal L}^{-1}$ is acyclic, then the induced map
  \[
  \Gamma({\cal L}')\oplus \Gamma({\cal L}')
  \to
  \Gamma({\cal L}\otimes {\cal L}')
  \]
  is surjective, implying the desired result.

  Otherwise, the nonnegativity assumption implies ${\cal L}'\cong {\cal
    L}$.  In that case, let $f_3$ be a section such that both $f_1,f_3$ and
  $f_2,f_3$ generate ${\cal L}$, and consider the short exact sequence
  \[
  0\to V
  \to {\cal L}^3
  \to {\cal L}^{\otimes 2}
  \to 0,
  \]
  where $V\cong \sO_Q^3/{\cal L}^{-1}$.  Since
  $(f_1,f_2,f_3):\Gamma(\sO_Q)\to \Gamma({\cal L})$ is injective,
  the dual map $H^1({\cal L}^{-1})\to H^1(\sO_Q^3)$ is surjective,
  and thus $H^1(V)=0$, giving the desired surjectivity of $\Gamma({\cal
    L}^3)\to \Gamma({\cal L}^{\otimes 2})$.
\end{proof}

\begin{lem}
  If $Q$ is nef with $Q^2=3$ and $q$ nontrivial, then the anticanonical
  ring is a filtered deformation of a ring generated by three elements of
  degree $1$ with three quadratic relations and one cubic relation.
\end{lem}

\begin{proof}
  Similar but now with
  \[
  0\to L_{-3}\to L_{-2}^3\to L_{-1}^3\to \sO_Q\to 0,
  \]
  which can be obtained by observing that $Q$ embeds in a noncommutative
  $\P^2$ in a way that respects the automorphism.
\end{proof}

Something similar should hold in higher degrees, and one expects the
algebra to be quadratic for $Q^2\ge 4$.  (When $Q$ is smooth, this is
shown, along with a number of generalizations, in \cite{filtered}.)

In special cases, one expects additional graded algebras; indeed, if $Q$ is
reducible, then we get a graded algebra corresponding to a singular form of
the surface from any linear combination of components of $Q$ which is nef
and has positive self-intersection.  There are also cases where the
autoequivalence is less canonical; e.g., for a noncommutative $\P^2$,
twisting by $\sO_X(1)$ is not an autoequivalence in general, but does act
trivially on the open substack of the moduli stack where the anticanonical
curve is integral and not an additive curve of characteristic 3.  (That is,
on that locus, the twisted surface is isomorphic to the original surface,
though not canonically.)  The corresponding scheme of isomorphisms is the
$\Pic^0(Q)[3]$-torsor corresponding to $q\in \Pic^0(Q)/\Pic^0(Q)^3$, giving
a $\Pic^0(Q)[3]$-torsor of graded algebra representations.  (This, of
course, gives rise to the classical construction of noncommutative planes
via graded algebras \cite{ArtinM/TateJ/VandenBerghM:1991}.)  In general,
for rational surfaces, a divisor class gives rise to an autoequivalence of
$X$ whenever there is an isomorphism of the corresponding commutative
surfaces that identify the two embeddings of $q$ as sheaves on the
respective anticanonical curves as well as the images of the nef cone of
$X$.  (The first condition gives a derived autoequivalence, and the second
ensures that it preserves the $t$-structure.  The latter can be weakened to
simply asking that the chosen divisor lie in the intersection of the images
of the nef cone of $X$ under the action of the group generated by the
autoequivalence on $\Pic(X)$.)

\chapter{Quot and Hilbert schemes}
\label{chap:quot}

\section{Families of sheaves}

There remains one of Chan and Nyman's axioms we have yet to show: the
``halal Hilbert scheme'' condition\footnote{It is unclear why a condition
on $\Quot$ schemes was given this name.}, which states that for any
coherent sheaf $M$, the functor $\Quot(M)$ classifying flat quotients of
$M$ is representable by a separated scheme, locally of finite type, which
is a countable union of projective schemes.  Since $\Quot(M)$ classifies
families of sheaves, this suggests that we should ask for this condition to
hold when $M$ itself is a family of sheaves, or better yet a sheaf on a
family of surfaces.  Furthermore, the commutative and maximal order cases
strongly suggest that there should be a better version of the $\Quot$
scheme condition: if $M$ is a $S$-flat family of sheaves on $X/S$ (with $S$
Noetherian), then the subfunctor of $\Quot$ classifying quotients of
Hilbert polynomial $P$ (relative to some fixed ample divisor class) should
be (and, as we will show, is) projective.

In order to carry out the usual commutative construction of the $\Quot$
functor, we need some boundedness results, as well as an analogue of the
flattening stratification.  A key issue is that we need to be able to force
not just the original family, but various base changes, to be acyclically
generated by line bundles $\sO_X(-bD_a)$ for uniformly bounded
$b$.  In particular, our previous results give this over a field, but we
both need to control how the bound varies with the fiber and need to deal
with the fact that acyclicity and global generation are not a priori
inherited from fibers when the sheaf is not flat over the base.

Luckily, this last fact turns out not to be an issue.  We will need to
start by showing this for curves.

\begin{lem}
  Let $C/S$ be a smooth proper curve over a Noetherian base $S$.  Then a
  sheaf $E\in \coh(C)$ is acyclic iff every fiber is acyclic.
\end{lem}

\begin{proof}
For any point $s$, the complex $R\Gamma(E)\otimes^L k(s)\cong
R\Gamma(E\otimes^L k(s))$ has vanishing cohomology in all degrees $>1$,
and thus $H^1(E)\otimes k(s)\cong H^1(E\otimes k(s))$.
\end{proof}

For global generation, a globally generated sheaf has globally generated
fibers, but the converse need not hold.  Luckily, if we also ask for
acyclicity, this problem goes away.  The key idea is the following.  In the
proofs of the next several results, we abbreviate ``acyclic and globally
generated'' by ``a.g.g.''.

\begin{lem}
  Let $C/S$ be a smooth proper curve over a Noetherian base $S$, and let
  $\Delta\subset C\times_S C$ be the diagonal, with projections $\pi_1$,
  $\pi_2$.  Then a sheaf $M\in \coh(C)$ is acyclic and globally generated
  iff $\pi_2^*M(-\Delta)$ is $\pi_{1*}$-acyclic.
\end{lem}

\begin{proof}
  We first note that $M$ is globally generated iff
  \[
  \pi_{1*}\pi_2^*M\to \pi_{1*}(\pi_2^*M|_{\Delta})
  \]
  is surjective, as this is nothing other than the map
  $\Gamma(M)\otimes_S\sO_C\to M$.  Let $N$ be the kernel of the
  (surjective) natural map $\pi_2^*M\to \pi_2^*M|_{\Delta}$.  Then we have
  an exact sequence
  \[
  0\to \pi_{1*}N\to \pi_{1*}\pi_2^*M\to \pi_{1*}(\pi_2^*M|_{\Delta}) \to
  R^1\pi_{1*}N\to R^1\pi_{1*}\pi_2^*M\to 0,
  \]
  from which it follows immediately that $M$ is a.g.g.~iff $R^1\pi_{1*}N=0$.
  We then note the short exact sequence
  \[
  0\to \Tor_1(\pi_2^*M,\sO_{\Delta})\to \pi_2^*M(-\Delta)\to N\to 0.
  \]
  Since the kernel is supported on $\Delta$, it is $\pi_{1*}$-acyclic, and
  thus $N$ is $\pi_{1*}$-acyclic iff $\pi_2^*M(-\Delta)$ is
  $\pi_{1*}$-acyclic.
\end{proof}

\begin{cor}
  Under the same hypotheses, the set of points $s\in S$ such that $M_s$ is
  acyclic and globally generated is open in $S$.
\end{cor}

\begin{proof}
  Indeed, it is the complement of the image of the support of
  $R^1\pi_{1*}\pi_2^*M(-\Delta)$.
\end{proof}

\begin{cor}
  A sheaf $M$ on $C/S$ is acyclic and globally generated iff every fiber of
  $M$ is acyclic and globally generated, and then every Noetherian base
  change of $M$ is acyclic and globally generated.
\end{cor}

\begin{proof}
  We may interpret $R^1\pi_{1*}\pi_2^*M(-\Delta)$ as the relative $H^1$ of
  a sheaf on the base change of $C$ to $C$.  Thus
  $R^1\pi_{1*}\pi_2^*M(-\Delta)=0$ iff for every point $x\in C$,
  $H^1(M(-x))=0$.  This is equivalent to asking that $H^1(M_s(-x))=0$ for
  every point $s\in S$ and $x\in C_s$, which is in turn equivalent to
  asking that every fiber be a.g.g.  Acyclic global generation on fibers is
  clearly inherited by base changes, so the remaining claim follows.
\end{proof}
  
\begin{prop}
  Let $X/S$ be a family of noncommutative surfaces over a Noetherian base,
  and let $L$ be a line bundle on $X$.  A sheaf $M\in \coh(X)$ is
  acyclically generated by $L$ iff every fiber of $M$ is acyclically
  generated by $L$, and then every Noetherian base change of $M$ is
  acyclically generated by $L$.
\end{prop}

\begin{proof}
  Acyclic generation descends along \'etale morphisms, so we may
  assume $X$ split.  We can then use twist functors to reduce to the case
  $L=\sO_X$.  If $m>0$, then $M$ is a.g.g.~iff it is a.g.g.~relative to
  $\alpha_m$ and $\alpha_{m*}M$ is a.g.g.  But the argument of Proposition
  \ref{prop:blowup_acyclic} carries over to an arbitrary family and shows
  that $M$ is a.g.g.~relative to $\alpha_m$ iff $\Ext^2(\sO_{e_m},M)=0$.  The
  same argument as in the curve case (since this is the highest possible
  degree of a nonvanishing cohomology sheaf) shows that this is equivalent
  to the same condition on fibers, and thus the result holds by induction.

  If $X$ is a ruled surface, we similarly find that $M$ is a.g.g.~iff it is
  a.g.g.~relative to $\rho_0$ and $\rho_{0*}M$ is a.g.g.  We can again use the
  semiorthogonal decomposition to find that $M$ is a.g.g.~relative to $\rho_0$
  iff $R^1\rho_{-1*}M=0$.  Since $\rho_{-1*}$ has cohomological dimension
  1, we once more find that this is equivalent to the condition on fibers,
  and thus the result reduces to the curve case.

  It remains to consider the case that $X$ is a noncommutative plane.  Take
  the (fppf) base change of $X$ to $Q$, and let $\tilde{X}$ be the blowup
  of this family in the tautological section of $Q$.  The sheaf $\alpha^*M$
  is a.g.g.~iff $\alpha_*\alpha^*M$ is a.g.g.~and $R^1\alpha_*\alpha^*M=0$.
  The spectral sequence associated to $R\alpha_*L\alpha^*M\cong M$ tells us
  that $R^1\alpha_*\alpha^*M$ always vanishes, and that there is a short
  exact sequence
  \[
    0\to R^1\alpha_*L_1\alpha^*M\to M\to \alpha_*\alpha^*M\to 0.
  \]
  The kernel is an extension of sheaves on $Q$ with support finite over
  $S$, and thus is a.g.g.; it follows that $M$ is a.g.g.~iff
  $\alpha_*\alpha^*M$ is a.g.g.

  Thus $M$ is a.g.g iff $\alpha^*M$ is a.g.g.~iff $(\alpha^*M)_q$ is
  a.g.g.~for all $q\in Q$ iff $(\alpha^*M)_s$ is a.g.g.~for all $s\in S$.
  Since $\alpha^*$ is right exact, $(\alpha^*M)_s\cong \alpha^*(M_s)$ and
  thus this is a.g.g.~iff $M_s$ is a.g.g.
\end{proof}

Here and below, we assume that $D_a$ is not only ample, but $D_a\cdot Q\ge
2$ or $D-2gf$ is nef as appropriate.  (This can always be arranged by
replacing $D_a$ by a sufficiently large multiple of $D_a$, without changing
the validity of the final result.  E.g., in the following result, we need
simply combine the bounds for finitely many twists of $M$ by multiples of
the original ample divisor.)

\begin{cor}\label{cor:noetherian_is_bounded}
  Let $X/S$ be a family of noncommutative surfaces over a Noetherian base.
  Then for any $M\in \coh(X)$, there is an integer $b_0$ such that for
  $b\ge b_0$, every Noetherian base change of $M$ is acyclically generated
  by $\sO_X(-bD_a)$.
\end{cor}

\begin{proof}
  It suffices to find a bound that works for $M$ itself.  For each $b$, let
  $M_b$ be the image of the natural map
  \[
  \sO_X(-bD_a)\otimes_{\sO_S}\Hom(\sO_X(-bD_a),M)\to M.
  \]
  Since $\sO_X(-bD_a)$ is acyclically generated by
  $\sO_X(-(b+1)D_a)$, this gives an ascending chain of subsheaves of $M$,
  which therefore terminates and gives a surjection
  \[
  \sO_X(-b_1D_a)\otimes_{\sO_S} N_1\to M
  \]
  for some $N_1\in \coh(S)$.  The kernel is again a coherent sheaf, so we
  can extend this to a resolution
  \[
  \sO_X(-b_0D_a)\otimes_{\sO_S}N_0\to \sO_X(-b_1D_a)\otimes_{\sO_S} N_1\to M\to 0
  \]
  with $b_0\ge b_1$.  Thus $M$ is generated by $\sO_X(-b D_a)$ as long as
  $b\ge b_1$, and is acyclic for $\sO_X(-bD_a)$ as long as $b\ge b_0$.  But
  then the same statement holds on an arbitrary Noetherian base change.
\end{proof}

\begin{rem}
  One consequence is that the support of $M$ over $S$ is closed (and can be
  given a natural closed subscheme structure), as it may be identified with
  the support of the coherent $\sO_S$-module $\Hom(\sO_X(-bD_a),M)$.  This,
  together with the fact that the cohomology sheaves of perfect complexes
  are coherent (since $X/S$ is Noetherian), implies that a perfect complex
  being a sheaf is an open condition (the complement of the unions of the
  supports of the other cohomology sheaves, only finitely many of which can
  be nonzero).
\end{rem}

In particular, this lets us prove the following criterion for flatness.

\begin{lem}
  Let $X/S$ be a family of noncommutative surfaces over a Noetherian base.
  A coherent sheaf $M\in \coh(X)$ is $S$-flat iff for all $b\gg 0$,
  $R\Hom(\sO_X(-bD_a),M)$ is an $S$-flat sheaf.
\end{lem}

\begin{proof}
  Let $b_0$ be such that $M$ is acyclically generated by $\sO_X(-bD_a)$ for
  all $b\ge b_0$.  Then the complex $R\Hom(\sO_X(-bD_a),M)$ is certainly a
  sheaf, and is clearly flat whenever $M$ is flat.  Thus suppose that
  $R\Hom(\sO_X(-bD_a),M)$ is a flat sheaf on $S$ for all $b\ge b_0$.  For
  any closed point $s\in S$, we can find $b'\ge b_0$ such that both
  $\Tor_1(M,k(s))$ and $\Tor_2(M,k(s))$ are acyclically generated
  by $\sO_X(-b'D_a)$.  But then we find that
  \begin{align}
  0
  &=
  \dim(h^{-1} R\Hom(\sO_X(-b'D_a),M)\otimes^L k(s))\notag\\
  &=
  \dim(h^{-1} R\Hom(\sO_X(-b'D_a),M\otimes^L k(s)))\notag\\
  &\ge
  \dim\Hom(\sO_X(-b'D_a),\Tor_1(M,k(s))),
  \end{align}
  so that $\Tor_1(M,k(s))=0$.
\end{proof}

We define the Hilbert polynomial of a coherent sheaf over a field in the
obvious way, namely as the polynomial $P_M(b)=\chi(\sO_X(-bD_a),M)$.  This
is clearly locally constant in flat families (since it respects field
extensions and on a geometric fiber depends only on the class of $M$ in
$K_0^{\num}(X)$).  The usual argument then gives the following.

\begin{cor}
  If $S$ is integral and Noetherian, then a coherent sheaf $M\in \coh X$ is
  $S$-flat iff all of its fibers have the same Hilbert polynomial.
\end{cor}

\begin{proof}
  For $b\gg 0$, so that every fiber of $M$ is acyclically
  generated by $\sO_X(-bD_a)$, we have
  $\dim\Hom(\sO_X(-bD_a),M_s)=\chi(\sO_X(-bD_a),M_s)$, and thus by
  Grauert, $\Hom(\sO_X(-bD_a),M)$ is flat.
\end{proof}

This lets us construct a weak form of flattening stratification.  For any
sheaf $M\in \coh(X)$, suppose that $M$ is acyclically generated by
$\sO_X(-bD_a)$ for all $b\ge b_0$.  Then for each such $b$, the sheaf
$\Hom(\sO_X(-bD_a),M)$ has a universal flattening stratification, and
taking the fiber product over $b\ge b_0$ gives a covering of $S$ by a
countable union of localizations of closed subschemes.  (This fiber product
exists and is affine since each factor is affine.)  This satisfies the
universal property that any $f:T\to S$ such that $f^*M$ is $T$-flat (a
``flattening morphism for $M$'') factors through this fiber product.  Say
that $M$ has ``finite flattening stratification'' if the fiber product is a
finite union of locally closed subschemes.

\begin{lem}  Let $X/S$ be a family of noncommutative surfaces over a
  Noetherian base $S$, and $M\in \coh(X)$.  Suppose that for any integral
  closed subscheme $S'\subset S$, there is an open neighborhood of $k(S')$
  in $S'$ on which $M$ is flat.  Then $M$ has finite flattening
  stratification.
\end{lem}

\begin{proof}
  Let $s$ be the generic point of some component of $S$.  Then by
  hypothesis, there is an open neighborhood of $k(s)$ in its reduced
  closure on which $M$ is flat, and thus in particular the Hilbert
  polynomial of $M$ is constant.  It follows by Noetherian induction that
  $S$ is the disjoint union of finitely many locally closed subsets on
  which the Hilbert polynomial of $M$ is constant.  The corresponding
  disjoint union of reduced subschemes is a flattening morphism, so factors
  through the universal flattening morphism.  Each set in the partition
  maps to at most one component of the flattening morphism, which therefore
  has only finitely many components, each of which must therefore be
  locally closed.
\end{proof}

This is somewhat difficult to prove in general, but a large family of
special cases where generic flatness holds was given in
\cite{ArtinM/ZhangJJ:2001} (a categorical version of the main result of
\cite{ArtinM/SmallLW/ZhangJJ:1999}), implying the following.

\begin{cor}\cite[Thm.~C5.1]{ArtinM/ZhangJJ:2001}
  If $S$ is of finite type over a field or $\Z$, then any sheaf on $X$ has
  finite flattening stratification.
\end{cor}

\begin{rem}
  In fact, it follows from the reference that this holds whenever $S$ is
  ``admissible'', as defined in \cite{ArtinM/SmallLW/ZhangJJ:1999}, but the
  finite type case turns out to be enough to let us bootstrap to a general
  Noetherian base, see Corollary \ref{cor:can_always_flatten} below.
\end{rem}

\section{The $\Quot$ scheme}

We now turn to the $\Quot$ scheme.  Given a family $X/S$ of noncommutative
surfaces over a Noetherian base and a sheaf $M\in \coh(X)$, we define a
functor $\Quot_{X/S}(M,P)$ on the category of schemes over $S$ by taking
$\Quot_{X/S}(M,P)(T)$ to be the set of isomorphism classes of short exact
sequences
\[
0\to I\to M_T\to N\to 0
\]
such that $N$ is a $T$-flat coherent sheaf all fibers of which have Hilbert
polynomial $P$.  (Recall that a ``$T$-flat coherent sheaf'' is a $T$-flat
element of $\perf(X_T)\cap \qcoh(X_T)$, so makes sense even when $T$ is not
Noetherian, and the usual argument tells us that $\Quot$ respects base
change.)  This is of course a contravariant functor of $T$, and also
satisfies a partial functoriality in $M$.  Indeed, given any surjection
$M'\to M$, we obtain a morphism $\Quot_{X/S}(M,P)\to \Quot_{X/S}(M',P)$ by
pullback of short exact sequences (in particular keeping the same quotient
sheaf).  This is a closed embedding, and thus to show that
$\Quot_{X/S}(M,P)$ is projective, it suffices to show this for some sheaf
surjecting on $M$.  We may thus reduce to considering
$\Quot(\sO_X(-bD_a)\otimes_{\sO_S} M,P)$ for $M\in \coh(S)$, and then by
twisting to the case $b=0$.  (Note that although we have only defined
twisting functors for split families, the twist functor descends as long as
the twisting line bundle $L$ is defined over $S$.)

In order to carry out the standard commutative argument, we need to show
that we can make a global choice of $b$ such that for any $T$ and any short
exact sequence in $\Quot_{X/S}(\sO_X\otimes_{\sO_S}M,P)(T)$, $\sO_X(-bD_a)$
acyclically generates $I$.  Since this reduces to a question about fibers
(modulo some technical issues coming from the fact that $T$ is not
Noetherian), a bound that works when $T$ is a geometric point of $S$ will
work in general.  In particular, if we can give a bound that works for a
surface over an algebraically closed field and depends only on the
combinatorics of the surface, this will produce a bound for $S$.

Thus let $X$ be a surface over an algebraically closed field $k$, and
consider a quotient of $\sO_X\otimes V$.  A complete flag in $V$ induces a
filtration of $\sO_X\otimes V$, and thus a filtration of the short exact
sequence.  This suggests that we should able to to reduce this problem to
the case $V=k$.

Thus consider a short exact sequence
\[
0\to I\to \sO_X\to M\to 0
\]
over a geometric point of $S$.  Since $\sO_X$ is torsion-free of rank 1,
$I$ either vanishes or is itself torsion-free of rank 1.

\begin{lem}
  Any nonzero homomorphism between torsion-free sheaves of rank 1 is
  injective.
\end{lem}

\begin{proof}
  The image is a nonzero subsheaf of a torsion-free sheaf, so must have
  positive rank, and thus the kernel is a rank 0 subsheaf of a torsion-free
  sheaf, so is 0.
\end{proof}

\begin{cor}
  If $I$, $I'$ are torsion-free sheaves of rank 1, then $\Hom(I,I')=0$
  unless $c_1(I')-c_1(I)$ is effective.
\end{cor}

We will want to understand the moduli space of such sheaves below, so one
natural question is the dimension of that space.

\begin{prop}
  Suppose that $I$ is a torsion-free sheaf of rank 1 on a noncommutative
  rational surface.  Then
  \[
  \chi(I,I)\le 1,
  \]
  with equality iff $I$ is a line bundle.
\end{prop}

\begin{proof}
Since $c_1(\theta I)-c_1(I)=K$ is ineffective, we find that
$\Ext^2(I,I)=\Hom(I,\theta I)=0$, and since any endomorphism of $I$ is
injective, we must have $\dim \Hom(I,I)=1$.  The inequality follows
immediately.  If equality holds, then since $\chi(I,I)$ is invariant under
twisting, we may as well assume that $c_1(I)=0$, at which point
Riemann-Roch gives $\chi(I)=1$.  We thus conclude that $I$ is numerically
equivalent to $\sO_X$ and $\Hom(\sO_X,I)\ne 0$, forcing $I\cong \sO_X$.
\end{proof}

In particular, the line bundles on a noncommutative rational surface are
precisely the torsion-free sheaves of rank 1 with $\chi(I,I)=1$, and thus
the notion is independent of the blowdown structure.  As this is of
independent interest, we include quasi-ruled surfaces in the next result.

\begin{prop}
  Suppose that $I$ is a torsion-free sheaf of rank 1 on a noncommutative
  rationally quasi-ruled surface over curves of genus $g_0$ and $g_1$.  Then
  \[
  \chi(I,I)\le 1-g_{c_1(I)\cdot f},
  \]
  with equality iff $I$ is a line bundle.
\end{prop}

\begin{proof}
  Again, we may as well assume that $c_1(I)=0$, and may use Riemann-Roch to
  restate the inequality as $\chi(I)\le 1-g_0$.  Let $\rho$ denote the
  composition $\rho_0\circ \alpha_1\circ\cdots \alpha_m$, and note that
  $\rho_*$ has cohomological dimension 1.  A homomorphism from $I$ to a
  line bundle must be injective, and thus if $\Hom(I,L)\ne 0$, then
  $c_1(L)$ is effective.  It follows that for any line bundle ${\cal L}$ on
  $C_0$, we have
  \[
  \Ext^2(\rho^*{\cal L},I)
  \cong
  \Hom(I,\theta \rho^*{\cal L})^*
  =
  0,
  \]
  since $(K+df)\cdot f=-2<0$ for any $d$.  It follows that
  \[
  \Ext^2({\cal L},R\rho_*I)=0
  \]
  for any ${\cal L}$ and thus that $\rank(R^1\rho_*I)=0$.  We then
  compute $\rank(\rho_*I)=\rank(I)=1$ and
  \[
  \chi(I) = \chi(\rho_*I)-\chi(R^1\rho_*I)\le \chi(\rho_*I).
  \]
  In particular, if $\chi(I)\ge 1-g_0$, then $\chi(\rho_*I)\ge 1-g_0$,
  and thus there is a degree 0 line bundle ${\cal L}$ on $C_0$ such that
  $\Hom({\cal L},\rho_*I)\ne 0$, making $\Hom(\rho^*{\cal L},I)\ne 0$.
  Since $\rho^*{\cal L}$ is again torsion-free of rank 1, we again find
  that such a morphism must be injective.  Since the quotient has trivial
  Chern class, it must be $0$-dimensional, and since $\Ext^1$ from a
  $0$-dimensional sheaf to a line bundle always vanishes, the corresponding
  short exact sequence splits.  Since $I$ is torsion-free, the quotient
  must be 0 and thus $I\cong\rho^*{\cal L}$.
\end{proof}

In the rationally ruled (or rational) case, we may use Riemann-Roch to
rewrite $\chi(I,I)$ in terms of $\chi(I)$ and obtain the following.

\begin{cor}\label{cor:line_bundle_bound}
  If $I$ is a torsion-free sheaf of rank 1 on a noncommutative rational or
  rationally ruled surface of genus $g$, then
  \[
  \chi(I)\le 1-g + c_1(I)\cdot (c_1(I)-K_X)/2,
  \]
  with equality iff $I$ is a line bundle.
\end{cor}

\begin{rem}
  An analogous formula of course exists in the quasi-ruled case, but is
  significantly messier.
\end{rem}

\begin{cor}\label{cor:lbs_are_intrinsic}
  Let $X$ be a noncommutative rational or rationally quasi-ruled surface
  and suppose that $L$ is a line bundle with respect to some blowdown
  structure on $X$.  Then $L$ is a line bundle with respect to {\em every}
  blowdown structure on $X$.
\end{cor}

\begin{proof}
  Certainly, $L$ is a torsion-free sheaf of rank 1 and saturates the
  relevant bound.  In the rational or rationally ruled case, the bound is
  manifestly independent of the blowdown structure, so the result is
  immediate.  For a strictly quasi-ruled surface, we need merely observe
  that the class $f$ is the same in any blowdown structure, so the bound
  remains invariant.
\end{proof}

We now return to considering only rational and rationally ruled surfaces.
It turns out that our bound on the Euler characteristic of torsion-free
sheaves of rank 1 also gives a bound on certain globally generated sheaves.

\begin{cor}\label{cor:Ox_quotient_bound}
  If $M$ is a quotient of $\sO_X$, then either $M\cong \sO_X$ or
  $\rank(M)=0$, $c_1(M)$ is effective, and
  \[
  \chi(M)\ge -c_1(M)\cdot (c_1(M)+K_X)/2.
  \]
\end{cor}

\begin{proof}
  Let $I$ be the kernel of the map $\sO_X\to M$, and observe that either
  $I=0$ or $I$ is torsion-free of rank 1, with
  $c_1(I)=-c_1(M)$ and $\chi(I)=\chi(\sO_X)-\chi(M)=1-g-\chi(M)$.
\end{proof}

\begin{cor}
  If $M$ is globally generated, then $c_1(M)$ is effective.
\end{cor}

\begin{proof}
  If $\sO_X\otimes \Hom(\sO_X,M)\to M$ is surjective, then each subquotient
  of the filtration induced by a complete flag in $\Hom(\sO_X,M)$ has
  effective Chern class, since it is either $\sO_X$ or $1$-dimensional.
\end{proof}

\begin{cor}\label{cor:finite_hps_for_glob_generated}
  For any integer $n$ and quadratic polynomial $P$, there is a finite
  subset of $K_0^{\num}(X)$, depending only on $D_a$ and the combinatorics
  of $Q$, that contains the class of any quotient of $\sO_X^n$ with Hilbert
  polynomial $P$.
\end{cor}

\begin{proof}
  The Hilbert polynomial determines the rank, the Euler characteristic, and
  $c_1(M)\cdot D_a$.  Since $c_1(M)$ is effective, we reduce to showing
  that there are finitely many classes that are effective on {\em some}
  surface with the given combinatorics and have given intersection with
  $D_a$.  The effective monoid of $X$ is generated by components of $Q$ and
  formal $-1$- and $-2$-curves on $X$, so it will suffice to find a
  finitely generated monoid containing the latter such that every generator
  has positive intersection with $D_a$.  We assume $m\ge 1$ (the cases
  $m\le 0$ are straightforward) so that we may choose a blowdown structure
  such that $D_a$ is in the fundamental chamber.

  Let $W(D_a)$ be the Weyl group generated by the reflections in the simple
  roots orthogonal to $D_a$.  Since $D_a^2>0$ so $D_a^\perp$ is negative
  definite, this is a finite Weyl group.  We then claim that any formal
  $-1$- or $-2$-curve is contained in the monoid generated by
  the $W(D_a)$-orbits of the simple roots not orthogonal to $D_a$ along
  with the $W(D_a)$-orbit of $e_m$ (and $f-e_1$ for $m=1$).  But this
  follows by the usual argument, once we observe that since $D_a$ is ample,
  a simple root orthogonal to $D_a$ cannot be effective.

  Since $D_a$ is $W(D_a)$-invariant, has positive intersection with $e_m$
  (and $f-e_1$ if needed), and nonnegative intersection with any simple
  root, it follows that this finitely generated monoid is indeed finitely
  graded by the intersection with $D_a$, as is the monoid obtained by
  throwing in the components of $Q$.
\end{proof}

Returning to a split family $X$ over a Noetherian base, let
$\Hilb^{n,+}(X)$ denote the functor such that $\Hilb^{n,+}(X)(T)$ is the
set of $T$-flat coherent sheaves such that every geometric fiber is torsion
free, with $\rank(I)=1$, $c_1(I)=0$ and $\chi(I)=1-g-n$.  (We will
eventually show (Theorem \ref{thm:hilb_is_projective}) that this is
represented by a smooth projective scheme over $S$.)  Note that these
conditions on a sheaf are independent of the blowdown structure, and thus
this functor can be defined for a non-split family.

We need to show that this functor is covered by a Noetherian scheme;
i.e., that there exists a family of sheaves $I$ over some base change of
$X$ such that every geometric point of $\Hilb^{n,+}(X)$ is a fiber of $I$.

\begin{lem}\label{lem:hilb_as_minimal_lift}
  For $m>0$, let $I\in \Hilb^{n,+}(X_m)(k)$ for a field $k$.  Then
  $\Hom(I,\sO_{e_m}(-d))=0$ for all $d>n$.
\end{lem}

\begin{proof}
  If $\Hom(I,\sO_{e_m}(-d))\ne 0$, then the image has the form
  $\sO_{e_m}(-d')$ for $d'>d$ so we may restrict to the case of a
  surjective morphism.  In that case, the kernel $I'\subset I$ is a
  torsion-free sheaf with $\rank(I')=1$, $c_1(I')=-e_m$,
  $\chi(I')=1-g-n-(-d+1)>-g$.  But then $I'(e_m)$ has rank 1, trivial Chern
  class and Euler characteristic $>1-g$, contradicting Corollary
  \ref{cor:line_bundle_bound}.
\end{proof}

\begin{prop}
  Let $X$ be a split family over a Noetherian base.  Then
  $\Hilb^{n,+}(X)$ is empty for $n<0$, while for $n>0$, it is covered by
  a Noetherian scheme.
\end{prop}

\begin{proof}
  For noncommutative planes, an explicit construction was given in
  \cite{NevinsTA/StaffordJT:2007} showing that the functor
  $\Hilb^{n,+}(X_{-1})$ is not only bounded but projective.

  For $m>0$, we have seen that every fiber of $I$ satisfies
  $\Hom(I,\sO_{e_m}(-n-1))=0$ and thus $\Hom(I(-ne_m),\sO_{e_m}(-1))=0$.
  Since we also have $\Hom(\sO_{e_m}(-1),I(-ne_m))=0$ since $I$ is
  torsion-free, we conclude that $I(-ne_m)\cong
  \alpha_m^{*!}\alpha_{m*}I(-ne_m)$.  But $\alpha_{m*}I(-ne_m)$ is easily seen
  to be a point of $\Hilb^{n(n+3)/2,+}(X_{m-1})$, letting us identify
  $\Hilb^{n,+}(X_m)$ with a closed subfunctor of
  $\Hilb^{n(n+3)/2,+}(X_{m-1})$.
    
  For $m=0$, the semiorthogonal decomposition gives a five-term exact
  sequence
  \[
  0\to \rho_1^*M_1\to \rho_0^*M_0\to I\to \rho_1^*N_1\to \rho_0^*N_0\to 0
  \label{eq:five_term_for_Hilb_of_ruled}
  \]
  for sheaves $M_i$, $N_i$ on $C$.  Chern class considerations tell us that
  $\Ext^2(\rho_0^*{\cal L},I)=\Ext^2(\rho_{-1}^*{\cal L},I)=0$ for any line
  bundle ${\cal L}$ on $C$, and thus that $\rank(N_0)=\rank(N_1)=0$, so
  that $\rank(M_1)=0$ and $\rank(M_0)=1$.  If $M_0$ had a nonzero
  $0$-dimensional subsheaf, then the corresponding image in $I$ would have
  nontrivial class in $K_0^{\num}(X)$ and rank 0, which is impossible since
  $I$ is torsion-free.  We thus conclude that $M_0$ is a line bundle and
  $M_1=0$.  Moreover, $\chi(M_0)-\chi(N_0)$ and $\chi(N_1)$ are determined
  by $[I]$, and $\deg(M_0)$ is bounded between $-n$ and $-1$.  Thus for
  each $n$, the moduli stack of potential triples $(M_0,N_1,N_0)$ is
  Noetherian, and remains so when we include the map $\rho_1^*N_1\to
  \rho_0^*N_0$ (on which we impose the open condition of surjectivity) and
  then the $\Ext^1$ from the kernel to $M_0$ and the torsion-free condition
  on the extension.  The result is a Noetherian family of sheaves
  surjecting on $\Hilb^{n,+}(X_0)$ as required.
\end{proof}

\begin{cor}
  For any quadratic polynomial $P(t)$, the functor $\Quot_{X/S}(\sO_X,P)$
  is covered by a Noetherian scheme.
\end{cor}

\begin{proof}
  The kernel in any corresponding short exact sequence is a torsion-free
  sheaf of rank 1, and thus up to a twist is just a point in
  $\Hilb^{n,+}(X)$ for suitable $n$.  There are only finitely many possible
  choices for the associated Chern class, and the short exact sequence is
  determined by a choice of injective map from the torsion-free sheaf to
  $\sO_X$.
\end{proof}

\begin{prop}
  For any coherent sheaf $M\in \coh(S)$ and any quadratic polynomial
  $P(t)$, the functor $\Quot_{X/S}(\sO_X\otimes_S M,P)$ is covered by a
  Noetherian scheme.
\end{prop}

\begin{proof}
  Since $S$ is Noetherian, it suffices to prove this on each piece of an
  affine covering of $S$, and on each such covering we can cover $M$ by a
  free module, so that it suffices to consider the case that $M\cong
  \sO_S^n$.

  Consider a subsheaf of $\sO_X^n$ with fixed class $[I_n]\in
  K_0^{\num}(X)$ having the given Hilbert polynomial.  (By Corollary
  \ref{cor:finite_hps_for_glob_generated}, there are only finitely many
  possible values of $[I_n]$.)  The filtration associated to a complete
  flag gives subquotients each of which is either $0$ or torsion-free of
  rank 1. If we ignore the Euler characteristics for the moment, then there
  are only finitely many possibilities for the rank and Chern class of the
  last subquotient, since the last subquotient is either $0$ or has Chern
  class $-D$ such that both $D$ and $-c_1([I])-D$ are effective.  We thus
  find by induction that there are only finitely many choices for the
  entire sequence of ranks and Chern classes of subquotients.  For each
  such choice, every subquotient has an upper bound on its Euler
  characteristic, and since the total Euler characteristic is fixed, there
  are only finitely many possible choices for the sequence of Euler
  characteristics.  In particular, we conclude that there are indeed only
  finitely many possible choices for the numerical invariants of the last
  subquotient, and thus only finitely many possible choices for its Hilbert
  polynomial.

  Fix such a choice $P_n$, and consider the morphism
  $\Quot_{X/S}(\sO_X^n,P)\to \Quot_{X/S}(\sO_X^{n-1},P-P_n)\times
  \Quot_{X/S}(\sO_X,P_n)$.  The fiber over a point corresponding to
  subsheaves $I_{n-1}\subset \sO_X^{n-1}$ and $I\subset \sO_X$ is
  represented by the affine space $\Hom(I,\sO_X^{n-1}/I_{n-1})$.  Since
  this is Noetherian and $\Quot_{X/S}(\sO_X,P-P_n)$,
  $\Quot_{X/S}(\sO_X,P_n)$ are covered by Noetherian schemes, the claim
  follows.
\end{proof}

Applying Corollary \ref{cor:noetherian_is_bounded} to this covering family
gives the following.

\begin{cor}
  Let $X/S$ be a family of noncommutative surfaces over a Noetherian base,
  and let $M$ be a coherent sheaf on $S$.  Then for any quadratic
  polynomial $p$, there is a bound $b_0$ such that for all $T/S$ and $I\in
  \Quot_{X/S}(\sO_X\otimes_{\sO_S} M,P)(T)$, $I$ is acyclically
  generated by $\sO_X(-bD_a)$ for all $b\ge b_0$.
\end{cor}

\begin{proof}
  Since $T$ need not be Noetherian, our previous arguments do not suffice
  to show that $I$ inherits acyclic generation from its fibers.  Acyclicity
  is not a difficulty, since we can represent $I$ by the complex
  \[
  \sO_X\otimes_{\sO_S} M\otimes_{\sO_S} \sO_T\to N
  \]
  with $N$ flat coherent (thus certainly inheriting acyclicity from fibers)
  and $M$ coherent, and thus
  \[
  \Ext^1(\sO_X(-bD_a),I)\otimes k(t)\cong
  \Ext^1(\sO_X(-bD_a),I\otimes k(t))=0
  \]
  for all points $t\in T$, implying $\Ext^1(\sO_X(-bD_a),I)=0$, since
  $\Ext^1(\sO_X(-bD_a),I)$ is a quotient of the (flat) coherent sheaf
  $\Hom(\sO_X(-bD_a),N)$.
  
  For global generation, it suffices to show that for each $b\ge b_0$,
  there is some $b_1$ such that
  \[
  \Hom(\sO_X(-b'D_a),\sO_X(-bD_a))\otimes \Hom(\sO_X(-bD_a),I)
  \to
  \Hom(\sO_X(-b'D_a),I)
  \]
  is surjective for all $b'\ge b_1$.  By the snake lemma, the cokernel is a
  quotient of the locally free coherent sheaf
  \[
  \ker(
  \Hom(\sO_X(-b'D_a),\sO_X(-bD_a))\otimes \Hom(\sO_X(-bD_a),N)
  \to
  \Hom(\sO_X(-b'D_a),N))
  ,
  \]
  so that again it suffices to check vanishing on fibers.  The existence of
  a global $b_1$ that works for every fiber follows using the Noetherian
  covering family.
\end{proof}

As is common in these constructions, the direct approach will only give
quasiprojectivity, so we need to check the valuative condition.

\begin{lem}
  The functor $\Quot_{X/S}$ is proper.
\end{lem}

\begin{proof}
  Let $X/R$ be a family of noncommutative surfaces, where $R$ is a discrete
  (since the base scheme is Noetherian) valuation ring with quotient field
  $K$ and field of fractions $k$.  For $M\in \coh(X)$, consider an exact
  sequence
  \[
  0\to I\to M_K\to N\to 0
  \]
  in $\coh(X_K)$.  We need to show that there is a unique short exact
  sequence
  \[
  0\to I'\to M\to N'\to 0
  \]
  such that $I'_K=I$ and $N'$ is $R$-flat, or equivalently $N'$ is
  $R$-torsion-free (since a module over a dvr is flat iff it is
  torsion-free).  Let $I'$ be the maximal submodule of $M$ such that
  $I'_K=I$.  Then $M/I'$ cannot have torsion, since we could move any
  torsion to the subsheaf without changing the generic fiber, while any
  proper subsheaf $I''\subsetneq I'$ has $I'/I''$ torsion since they have
  the same generic fiber, and thus makes $N''$ fail to be torsion-free.
  It follows that $I'$ gives rise to the desired unique extension.
\end{proof}

\begin{thm}
  Let $X/S$ be a family of noncommutative surfaces over a Noetherian base.
  Then for any $M\in \coh(X)$, $\Quot_{X/S}(M,P)$ is a projective scheme
  over $S$.
\end{thm}

\begin{proof}
  As already discussed, it suffices to prove this for
  $\Quot_{X/S}(\sO_X\otimes_{\sO_S} M,P)$ with $M\in \coh(S)$.  We thus
  find in general that there is a constant $b_0$ such that for any $T$-point
  of $Q:=\Quot_{X/S}(\sO_X\otimes_{\sO_S}M,P)$, the kernel in the
  corresponding short exact sequence is acyclically generated by
  $\sO_X(-b_0D_a)$.  We thus obtain a resolution
  \[
  \sO_X(-b_0D_a)\otimes_{\sO_T} \Hom(\sO_X(-b_0D_a),I)
  \to
  \sO_X\otimes_{\sO_T} M_T
  \to
  N_T
  \to
  0,
  \]
  giving a $T$-point of the Grassmannian
  \[
  G:=\Grass(\Hom(\sO_X(-b_0D_a),\sO_X)\otimes_{\sO_S} M,P(b_0))
  \]
  classifying rank $P(b_0)=\rank(\Hom(\sO_X(-b_0D_a),N_T))$ quotients of
  $\Hom(\sO_X(-b_0D_a),\sO_X)\otimes_{\sO_S} M$.
  
  In the other direction, any point of the Grassmannian determines a
  quotient sheaf of $\sO_X\otimes_{\sO_T}M_T$, which gives a point of $Q$
  iff it has the correct Hilbert polynomial.  It follows that $Q$ is a
  component of the universal flattening stratification of the universal
  quotient sheaf over the Grassmannian.  If we knew that sheaf had finite
  flattening stratification, we could then conclude that $Q$ was locally
  closed in $G$.  Since $Q$ is proper and $G$ is projective, this would
  imply $Q$ projective as required.  In particular, this implies that the
  theorem holds whenever $S$ (and thus $G$) is of finite type.

  For more general base rings, it clearly suffices to show that $Q$ is
  closed in $G$.  In particular, if we can show that $Q$ is fppf locally
  projective, then it follows that the map $Q\to G$ is fppf locally closed,
  so closed, and thus $Q$ must be projective.  As the question is now fppf
  local, we may now make some simplifying assumptions: since $X$ is
  \'etale locally split, we may assume that $X$ is split, and since
  $M$ is Zariski locally a quotient of a free sheaf, we may assume this is
  so.  Since the corresponding map of $\Quot$ functors is closed, we may
  assume $M\cong \sO_X^n$.
  
  Let ${\cal M}$ denote the moduli stack of split surfaces, and note that
  $X/S$ induces a map $S\to {\cal M}$.  Since $S$ is Noetherian and ${\cal
    M}$ is locally of finite type, there is a morphism $U\to {\cal M}$ with
  $U$ of finite type such that the fiber product is smooth and
  surjective over $S$.  Let $T$ be that fiber product, and observe that one
  has
  \[
  \Quot_{X/S}(\sO_X^n,P)\times_S T
  \cong
  \Quot_{X_T/T}(\sO_X^n,P)
  \cong
  \Quot_{X_U/U}(\sO_X^n,P)\times_U T.
  \]
  Since $U$ is of finite type, the corresponding $\Quot$ scheme is
  projective, and thus $\Quot_{X/S}(\sO_X^n,P)$ is fppf locally projective
  as required.
\end{proof}

\begin{rem}
  Of course, in the last part of the argument, what we are really proving
  is that the functor $\Quot_{X/{\cal M}}(\sO_X^n,P)$, appropriately
  defined, is fppf locally projective.
\end{rem}

Although much of the above proof involved working around the potential lack
of finite flattening, this is in some sense unnecessary, as finite
flattening always holds.  Unfortunately, the only proof we have of this
fact uses projectivity of the $\Quot$ scheme!

\begin{cor}\label{cor:can_always_flatten}
  Let $X/S$ be a family of noncommutative surfaces over a Noetherian base.
  Then any sheaf $M\in \coh(X)$ has finite flattening stratification.
\end{cor}

\begin{proof}
  It suffices to show that when $S$ is integral, there is an open
  neighborhood of $k(S)$ on which $M$ is flat, or equivalently on which the
  Hilbert polynomial is constant.  Let $P$ be the Hilbert polynomial of
  $M_{k(S)}$, and consider the Quot scheme $Q=\Quot_{X/S}(M,P)$, with
  associated short exact sequence
  \[
  0\to I\to M_Q\to N\to 0.
  \]
  We first note that the fiber over $k(S)$ consists of the single short
  exact sequence
  \[
  0\to 0\to M_{k(S)}\to M_{k(S)}\to 0,
  \]
  as the Hilbert polynomial constraint forces the map $M\to N$ to become
  an isomorphism.  In particular, since the image of $\pi:Q\to S$ contains
  $k(S)$, the morphism is surjective.  Since $N$ is flat, the restriction
  of the short exact sequence to any point $q\in Q$ will still be a short
  exact sequence, and thus we find that $I_q\ne 0$ precisely when the
  Hilbert polynomial of $M_{\pi(q)}$ is different from $p$.  In other
  words, the closed subset $\pi(\supp(I))$ of $S$ agrees with the set of
  points where the Hilbert polynomial jumps.  It follows that the set of
  points with the same Hilbert polynomial as the generic point is open as
  required.
\end{proof}

\begin{rem}
  This argument in fact shows that the class in the Grothendieck group is
  upper semicontinuous: it is constant on locally closed subsets, and the
  value at a point on the boundary of such a set differs by the class of a
  sheaf over that point.
\end{rem}

\section{The Hilbert scheme of points}

We now turn to a closer analysis of the functors $\Hilb^{n,+}(X/S)$ for
families over a Noetherian base.  Note that since any torsion-free sheaf of
rank 1 is stable, we automatically find (via the standard argument of
Langton) that the functor is proper, with smoothness following by
deformation theory, since $\Ext^2(I,I)=0$.  We then find that the dimension
is everywhere $1-\chi(I,I)=2n+g$.

\begin{thm}\label{thm:hilb_is_projective}
  For all $n\ge 0$, the functor $\Hilb^{n,+}(X/S)$ is represented by a
  smooth projective scheme of dimension $2n+g$ over $S$ with geometrically
  integral fibers.
\end{thm}

\begin{proof}
  Since $\Hilb^{n,+}(X/S)$ is covered by a Noetherian family, there is a
  line bundle $L=\sO_X(-bD_a)$ such that any point $I\in \Hilb^{n,+}(X)$ is
  acyclically generated by $L$.  As above, a complete flag in $\Hom(L,I)$
  determines a filtration of $I$ and there are only finitely many
  possibilities for the sequence of subquotients; it follows that there are
  only finitely many possibilities for the class in $\Pic^{\num}(X)$ of a
  subsheaf of $I$ which is generated by $L$.  It follows as in
  \cite[Lems.~4.4.5,4.4.6]{HuybrechtsD/LehnM:1997} that there is a
  linearization of the action of $\GL(\Hom(L,I))$ on the $\Quot$
  scheme such that every point corresponding to $I$ is GIT stable.  (Note
  that since $I$ has rank 1, any subsheaf still satisfies the relevant
  bound, even if it has more global sections than expected.)  But then
  $\Hilb^{n,+}(X/S)$ is a subscheme of the GIT quotient, and is therefore
  quasiprojective, so (by properness) projective as required.

  For irreducibility, it suffices to consider the moduli problem over a
  smooth covering of the moduli stack of surfaces; in particular, we may
  assume that the family contains commutative fibers.  When $X$ is
  commutative, the numerical conditions on $I$ force it to be the tensor
  product of a line bundle in $\Pic^0(X)$ by the ideal sheaf of an
  $n$-point subscheme, and thus that fiber is connected, and thus
  (following \cite[Prop.~8.6]{NevinsTA/StaffordJT:2007}) every fiber is
  connected.  Since the fibers are smooth, they are geometrically integral
  as required.
\end{proof}

\begin{rem}
  An alternate proof of irreducibility is given as a special case of Theorem
  \ref{thm:moduli_irred1} below.
\end{rem}

In the rational case, this is a family of deformations of the Hilbert
scheme of the corresponding commutative surface, except when $n=1$, when it
{\em is} the corresponding commutative surface.

\begin{prop}
  Let $X$ be a commutative rational surface with anticanonical curve $Q$,
  and let $X_q$ be the noncommutative deformation associated to $q\in
  \Pic^0(Q)$.  Then $\Hilb^{1,+}(X_q)\cong X$.
\end{prop}

\begin{proof}
  Let $I\in \Hilb^{1,+}(X_q)$.  Then $\Hom(\sO_{X_q},I)=0$ since such a map
  would be injective with $0$-dimensional cokernel of negative Euler
  characteristic, while $\Hom(I,\theta\sO_{X_q})=0$ since such a map would be
  injective with cokernel of Chern class $-Q$.  Since $\chi(\sO_{X_q},I)=0$, we
  conclude that $\Ext^1(\sO_{X_q},I)=0$ and thus $I\in \sO_{X_q}^\perp$.
  We may thus use the functor $\kappa_q$ to move this to $\sO_X^\perp$,
  and the claim is that this induces an isomorphism $\Hilb^{1,+}(X_q)\cong
  \Hilb^1(X)\cong X$.

  We first consider the case $m=-1$.  In this case, we similarly find that
  $\Hom(\sO_{X_q}(\pm 1),I)=\Ext^2(\sO_{X_q}(\pm 1),I)=0$, and thus $I$ has a
  canonical resolution of the form
  \[
  0\to \sO_{X_q}(-2)\to \sO_{X_q}(-1)^2\to I\to 0,
  \]
  while for any $2$-dimensional subspace of
  $\Hom(\sO_{X_q}(-2),\sO_{X_q}(-1))$, the corresponding quotient will be
  torsion-free with the correct invariants.  The action of $\kappa_q$,
  $\kappa_q^{-1}$ on such a presentation is essentially trivial, and thus
  the claim follows.

  Similarly, for $m=0$, we have a short exact sequence of the form
  \[
  0\to \sO_{X_q}(-f)\to I\to \sO_f(-1)\to 0
  \]
  where $\sO_f(-1)$ denotes some pure 1-dimensional sheaf of Euler
  characteristic 0 and Chern class $f$; conversely, every nontrivial
  extension of such a sheaf by $\sO_{X_q}(-f)$ gives rise to a point of
  $\Hilb^{1,+}(X_q)$.  Since $\kappa_q\sO_{X_q}(-f)\cong \sO_X(-f)$
  and $\kappa_q\sO_f(-1)$ has the form $\sO_f(-1)$, the claim follows.

  Finally, for $m>0$, there are two cases.  If $\Hom(I,\sO_{e_m}(-1))=0$,
  then $I$ is $\alpha_{m*}$-acyclic and $\alpha_{m*}I$ is torsion-free with
  $I\cong \alpha_m^*\alpha_{m*}I$, and thus the claim follows by induction.
  Otherwise, the map must be surjective by Lemma
  \ref{lem:hilb_as_minimal_lift}, and thus we have a non-split short exact
  sequence
  \[
  0\to \sO_{X_q}(-e_m)\to I\to \sO_{e_m}(-1)\to 0,
  \]
  and again the claim follows by observing that
  $\kappa_q(\sO_{X_q}(-e_m))\cong \sO_X(-e_m)$ and
  $\kappa_q(\sO_{e_m}(-1))\cong \sO_{e_m}(-1)$.
\end{proof}

\begin{rem}
  With this in mind, we will feel free to omit the ``$+$'' in the rational
  case.
\end{rem}

To recover something closer to the usual Hilbert scheme in higher genus, we
observe that there is a natural morphism $\Hilb^{n,+}(X)\to \Pic^0(Q)$
given by $I\mapsto \det(I|^{\bf L}_Q)$.  Since the numerical class of $I$
is fixed, the image is contained in a particular coset of $\Pic^0(C)$.  Any
line bundle in $\Pic^0(C)$ induces an autoequivalence of $\coh X$, and we
find that this autoequivalence simply twists the determinant of the
restriction to $Q$ by that same line bundle.  We thus find that the fibers
of the determinant-of-restriction map are all isomorphic, and we thus we
have an isomorphism of the form $\Hilb^{n,+}(X)\cong \Hilb^n(X)\times
\Pic^0(C)$ where $\Hilb^n(X)$ is any fiber of the
determinant-of-restriction map.  (One caveat is that this presumes the
image of the determinant-of-restriction map has a point!  This is automatic
if $g\ge 2$, when there is a canonical retraction $\Pic^0(Q)\to \Pic^0(C)$,
while for $g=1$, it at least holds \'etale locally, and the transition maps
simply involve translation in $\Pic^0(C)$, so that there is a canonical
isomorphism of the form $\Hilb^n(X)\times T$ where $T$ is the relevant
coset of $\Pic^0(C)$.)  We conclude that $\Hilb^n(X)$ is smooth,
irreducible, and projective of dimension $2n$.  When $X$ is commutative,
the relevant coset of $\Pic^0(C)$ is the trivial coset, and thus there is a
canonical choice of fiber, and we find that the result agrees with the
usual Hilbert scheme of $n$ points.  In general, $\Hilb^1(X)$ is a Poisson
surface, and thus gives a canonical map from the moduli stack of
noncommutative surfaces to the moduli stack of commutative anticanonical
surfaces.  This map is more natural than the map implicit in Proposition
\ref{prop:semiorth_of_quasi-ruled}, but harder to control; thus we do not
attempt to prove that it induces an isomorphism with the relative
$\Pic^0(Q)/\Pic^0(C)$.  (By the construction in the next paragraph, it is
at the very least rationally ruled over $C$.)

In the commutative case, there is also a map from the Hilbert scheme of $n$
points to $\Sym^n(C)$.  This too has an analogue in the noncommutative
setting.  We saw in the case of a ruled surface that $I$ always satisfies
$\rho_{-1*}I=0$ and $R^1\rho_{-1*}I$ is 0-dimensional of degree $n$ (these
are the sheaves $M_1$ and $N_1$ in equation
\eqref{eq:five_term_for_Hilb_of_ruled}), inducing a map $\Hilb^{n,+}(X)\to
\Sym^n(C)$ which is preserved by the action of $\Pic^0(C)$, and thus
factors through $\Hilb^n(X)$.  In general, one finds that for $I\in
\Hilb^{n,+}(X_m)$, there is an integer $0\le d\le n$ such that
$\alpha_{m*}I$ is in $\Hilb^{d,+}(X_{m-1})$ and $R^1\alpha_{m*}I$ is a
$0$-dimensional sheaf of degree $n-d$.  It follows by induction that
\[
h^1(R\rho_{-1*}R\alpha_{1*}\cdots R\alpha_{m*}I)
\]
is again a 0-dimensional sheaf of degree $n$, so the map to $\Sym^n(C)$
survives.

We also have the following regarding the direct image in $C$.

\begin{prop}\label{prop:hilb_rho*_is_line_bundle}
  There is a nonempty open subscheme of $\Hilb^{n,+}(X)$ on which the
  direct image of $I$ to $C_0$ is a line bundle of degree $-n$.
\end{prop}

\begin{proof}
  Let $y_1,\dots,y_n$ be $n$ points in $Q$ such that their images in $C_0$
  are distinct, and let $Z$ be the structure sheaf of the corresponding
  $n$-point subscheme.  Then we may define a point $I\in \Hilb^{n,+}(X)$ by
  the short exact sequence
  \[
  0\to I_Z\to \sO_X\to Z\to 0,
  \]
  and taking direct images to $C$ gives a four-term sequence
  \[
  0\to \rho_*I_Z\to \sO_{C_0}
   \to \rho_*Z
   \to R^1\rho_*I_Z
   to 0.
  \]
  The constraint that the images be distinct ensures that the map
  $\sO_{C_0}\to \rho_*Z$ is surjective, and thus $R^1\rho_*I_Z=0$.  By
  semicontinuity, this implies $R^1\rho_*I=0$ holds on a nonempty open
  subset of $\Hilb^{n,+}(X)$, which in turn implies that $\rho_*I$ is a
  line bundle of degree $-n$ as required.
\end{proof}

\begin{thm}
  Let $X$ be a noncommutative rationally ruled surface over a curve $C$
  over an algebraically closed field $k$.  Then $\Hilb^n(X)$ is birational
  to the $n$-th symmetric power of $C\times \P^1$.
\end{thm}

\begin{proof}
  We first observe that if $\alpha:\tilde{X}\to X$ is a one-point blowup,
  and $I|^L_Q$ is a line bundle, then $\alpha^*I$ is torsion-free and thus
  represents a class in $\Hilb^{n,+}(\tilde{X})$.  Since $I|^L_Q$ being a
  line bundle is generic behavior, we get a rational map from
  $\Hilb^{n,+}(X)$ to $\Hilb^{n,+}(\tilde{X})$, which is birational to its
  image (since $\alpha_*\alpha^*I\cong I$) and thus birational by
  irreducibility.  Moreover, since twisting by $\Pic^0(C)$ does not affect
  this condition, we have an induced birational map $\Hilb^n(X)\to
  \Hilb^n(\tilde{X})$.  It follows that it suffices to prove the claim when
  $X$ is a ruled surface.

  In that case, we generically have a short exact sequence of the form
  \[
  0\to \rho^*_0{\cal L}\to I\to \rho_1^*Z\to 0,
  \]
  where ${\cal L}$ is a line bundle of degree $-n$, $Z$ is a
  $0$-dimensional sheaf of degree $n$, and ${\cal L}$ is determined from
  $Z$ by the chosen value of $\det I|^L_Q$.  (In particular, we could
  take ${\cal L}\cong \det(Z)^{-1}$.)  For each point $x\in C$,
  $\Ext^1(\rho_1^*\sO_x,\rho_0^*{\cal L})\cong k^2$, so the actual
  non-split extension sheaves are classified by $\P^1$; applying this to
  the $n$ distinct points in the support of $Z$ tells us that the fiber
  over $Z$ in $\Hilb^n(X)$ is $(\P^1)^n$, and thus that $\Hilb^n(X)$ is
  birational to $(\P^1)^n\times \Sym^n(C)$.  Since this applies equally
  well in the commutative case, the claim follows.
\end{proof}

\begin{rem}
  In the above proof, we saw that there was a particular choice of point in
  the relevant coset of $\Pic^0(C)$ in $\Pic(Q)$.  This is somewhat
  misleading, however, as it depended on having a well-defined functor
  $\rho_1^*$, while this is only in general determined up to
  $\Aut(C)\ltimes \Pic(C)$, and thus can fail to be rational over
  non-closed fields or in more general families.
\end{rem}

\section{The Hilbert scheme and duality}
\label{sec:duality_of_Hilb}

One surprising feature of the noncommutative case is that there is
typically a duality on the Hilbert scheme, more precisely the following.

\begin{prop}
  Let $n\ge 0$ be a nonnegative integer, and suppose that $|\langle
  q\rangle|>n$.  Then the functor $I\mapsto \theta^{-1}\ad I$ induces
  a birational map $\Hilb^{n,+}(X_q)\cong \Hilb^{n,+}(X_{q^{-1}})$.
\end{prop}

\begin{proof}
  We claim more precisely that if $I\in \Hilb^{n,+}(X_q)$ with $I|_Q$ a
  line bundle, then $R\ad I$ is a torsion-free sheaf, which immediately
  implies that $\theta^{-1}\ad I\in \Hilb^{n,+}(X_{q^{-1}})$.

  We first recall that $R^2\ad I$ is $\le 0$-dimensional and $R^1\ad I$ is
  $\le 1$-dimensional.  If either bound were tight, then the isomorphism
  $R\ad R\ad I\cong I$ would produce a $<2$-dimensional subsheaf of $I$,
  contradicting the requirement that $I$ be torsion-free.  In particular,
  we see that $R^2\ad I=0$ and $R^1\ad I$ is $\le 0$-dimensional.
  Similarly, if $\ad I$ were not torsion-free, then $R\ad R\ad I$ could
  not be a sheaf.  Since
  \[
  (R\ad I)^\dL_Q\cong R\sHom_Q(I|^\dL_Q,\omega_Q)
  \]
  is a line bundle, $R^1\ad I$ must be disjoint from $Q$, and thus
  $\chi(R^1\ad I)$ must be a multiple of $|\langle q\rangle|$.  But this
  would cause $\theta^{-1}\ad I$ to violate Corollary
  \ref{cor:line_bundle_bound} unless $R^1\ad I=0$ as required.
\end{proof}

\begin{rems}
  In fact, there is always a birational map $\Hilb^{n,+}(X_{q})\cong
  \Hilb^{n,+}(X_{q^{-1}})$; if $q\cong \sO_Q$, this is trivial (though the
  map does not take the above form!), while otherwise one can verify using
  the above argument that the locus where $R\ad I$ is not a sheaf has
  positive codimension.
\end{rems}

\begin{rems}
  It seems likely that this birational map is an isomorphism in general.
  Indeed, we can always embed $I$ in its reflexive hull $I^+:=\ad\ad I$,
  giving a short exact sequence
  \[
  0\to I\to I^+\to Z\to 0
  \]
  with $Z$ $0$-dimensional, from which we see that $\ad I\cong \ad I^+$ and
  $R^1\ad I\cong R^2\ad Z$.  We thus find that
  \[
  \Hom(\theta^{-1}\ad I,R^1\ad I)
  \cong
  \Hom(\ad \theta I^+,R^2\ad Z)
  \cong
  \Ext^2(Z,\theta I^+)
  \cong
  \Hom(I^+,Z)^*
  \]
  where the second isomorphism uses the fact that $R\ad$ is a contravariant
  equivalence.  This also tells us that both Hom complexes are supported in
  degree 0, and thus $\dim\Hom(\theta^{-1}\ad I,R^1\ad I)=\chi(Z)$.
  In particular, we expect that $\Aut(Z)$ has a dense orbit in this $\Hom$
  space such that taking the kernel of any morphism in that orbit gives the
  same element of $\Hilb^{n,+}(X_{q^{-1}})$.
\end{rems}

In the commutative setting, the proof that certain moduli spaces of
1-dimensional sheaves are rational involved a birational map to the
symmetric power of the surface.  In the noncommutative case, the symmetric
power is of course a noncommutative scheme, but since the symmetric power
of a commutative surface is birational to the Hilbert scheme of points,
this suggests we should look for birational maps in that context.  In other
words, we should look for a map between points of the Hilbert scheme and
1-dimensional sheaves of Euler characteristic 1.  Recall that in the
commutative case, a 1-dimensional sheaf of Euler characteristic 1 is
generically a line bundle of degree $g$ on an integral curve, which in turn
generically has a unique effective representative, corresponding to a
$0$-dimensional subscheme of the curve and thus of the surface.  The
difficulty in translating this directly to the noncommutative setting is
that the $0$-dimensional subscheme appears as its structure sheaf rather
than as its ideal sheaf, but we can fix this by duality: indeed, the
$2$-term complex $\sO_X\to M$ is dual to an ideal sheaf (under
$\theta^{-1}R\ad$).

We thus begin by considering the following dual scenario.  Let $X$ be a
noncommutative rational surface, and let $M$ be a pure $1$-dimensional
sheaf on $X$, disjoint from $Q$, of Euler characteristic $-1$ and Chern
class $c_1(M)=D$.  Suppose moreover (for simplicity) that $M$ is {\em
  irreducible}, in the sense that any nontrivial subsheaf also has $c_1=D$.
We find
\[
\chi(M,\theta \sO_X) = \chi(\sO_X,M) = -1
\]
and thus ``generically'' we expect that $R\Hom(M,\theta \sO_X)=k[-1]$.  (We
will see below that under fairly general conditions, the moduli space of
$D_a$-stable sheaves of class $[M]$ is irreducible and has a nonempty open
subset on which the sheaf is irreducible and only has a single cohomology
space.)  When this happens, we obtain a corresponding (non-split) universal
extension
\[
0\to \theta \sO_X\to F\to M\to 0,
\]
with $F\in {}^{\perp}\theta \sO_X\cong \sO_X^\perp$.  Any torsion subsheaf
of $F$ would be a subsheaf of $M$, and cannot have the same Chern class
since that would induce a non-split extension of the $0$-dimensional
quotient by $\theta \sO_X$, so the irreducibility hypothesis ensures that
$F$ is torsion-free.  Since it also has rank $1$, we can write it as
$\theta I(D)$ where $I$ is a point of $\Hilb^{D^2/2+1}(X(-D))$, giving us
a rational map from the space of sheaves to $\Hilb^{D^2/2+1}(X(-D))$.

To recover $M$, we a priori need to include the map $\theta \sO_X\to F$ in
the data.  But since $F\in \sO_X^\perp$, we have
\[
R\Hom(\theta \sO_X,F)\cong R\Hom_Q((\theta \sO_X)|^{\dL}_Q,F|^{\dL}_Q),
\]
and may thus recover the map $\theta \sO_X\to F$ from the isomorphism
\[
F|_Q\cong (\theta \sO_X)|_Q.
\]
We thus see that the rational map is birational to its image, and thus
dimension considerations (plus irreducibility) tell us that it is itself
birational.  That is, for generic sheaves $I$ in the appropriate Hilbert
space (and assuming sheaves of class $[M]$ disjoint from $Q$ exist),
$\theta I(D)$ is in $\sO_X^\perp$ and satisfies $I(D)|_Q\cong \sO_Q$, so
may be used to construct a sheaf $M$.

When $|\langle q\rangle|>D^2/2+1$, the generic sheaf $I$ is reflexive, and
thus we may dualize the above construction, which leads to the following.
Given irreducible $M$ with $\chi(M)=1$, it is generically acyclic and
globally generated, and the kernel $F\in \sO_X^\perp$ of the unique global
section can be written in the form $I(-D)$ for $I\in
\Hilb^{D^2/2+1}(X(D))$, again defining a birational map.

The fact that $M$ is generated by a single global section is of course
quite unexpected from the perspective of commutative algebraic geometry,
but it less so from the perspective of difference equations, where it
corresponds to the fact that one can frequently recover the coefficients of
a solution to a matrix equation in terms of the shifts (or derivatives) of
a single entry in the solution vector.  In other words, the ability to
generate an equation from a single global section corresponds to the
existence of ``straight-line'' forms of equations.  We can see this in the
present setting as follows.  First, note that for $I\in \Hilb^n(X)$,
$I(nf)$ is generically acyclic, as an immediate consequence of Proposition
\ref{prop:hilb_rho*_is_line_bundle}.  Since $\chi(\sO_X(-nf),I)=1$, we
generically have a natural map $\sO_X(-nf)\to I$.  Applying this to the
sheaf $F$ gives a composition
\[
\sO_X(-D-nf)\to F\to \sO_X
\]
expressing $M$ as a quotient of the cokernel of a morphism $\sO_X(-D-nf)\to
\sO_X$ (which we recall can be represented as a (twisted) difference(tial)
operator).  The corresponding subsheaf of the cokernel has Chern class
$nf$, and thus (per the discussion in Section \ref{sec:sheaves_from_eq_nc}
below) at most adds apparent singularities to the equation corresponding to
$M$.  In other words, a ``typical'' 1-dimensional sheaf of Euler
characteristic 1 corresponds to an equation in straight-line form with $n$
apparent singularities, where $n$ is half the dimension of the
corresponding moduli space.

An analogous construction works in the higher genus case, with the main
difference being that we must replace $\sO_X(-nf)$ by $\rho^*\rho_*I$, and
thus we still obtain a straight-line form with $\dim/2$ apparent singularities.

We can also use this idea to prove the claim that irreducible sheaves
generically have a single cohomology group, and thus the above
constructions give birational maps to the Hilbert scheme in the $\chi=\pm
1$ cases.

\begin{lem}\label{lem:irreds_generically_acyclic}
  Let $X$ be a noncommutative rational surface with $q$ non-torsion.  The
  generic sheaf in any irreducible component of the moduli space of
  irreducible sheaves disjoint from $Q$ satisfies either $H^0(M)=0$ or
  $H^1(M)=0$.
\end{lem}

\begin{proof}
  The moduli space in question is representable by an intersection of open
  subsets of a quasi-projective scheme, since any irreducible sheaf is
  stable and stable moduli spaces of rank 0 sheaves are quasi-projective by
  Theorem \ref{thm:semistable_moduli_are_projective} below.

  Since $R^1\ad$ swaps the two cohomology groups, it suffices to show that
  when $\chi(M)=l\le 0$, $M$ generically has no global sections.  Since
  deformation theory (plus smoothness of symplectic leaves) tells that
  irreducible sheaves vary in a $D^2+2$-dimensional family, it suffices to
  show that the irreducible sheaves with $n$ global sections vary over a
  smaller family.  Since $q$ is non-torsion, any global section of $M$ is
  surjective, and thus if $\dim\Hom(\sO_X,M)=n$, we have a $\P^{n-1}$ of
  surjective maps $\sO_X\to M$, giving a $\P^{n-1}$-bundle over the family
  of irreducible sheaves with $n$ global sections.  The kernel of such a
  map has the form $I(-D)$ for $I\in \Hilb^{D^2/2+l}$, and one has a short
  exact sequence
  \[
  0\to \Hom(\sO_X,\sO_X)\to \Hom(I(-D),\sO_X)\to \Ext^1(M,\sO_X)\to 0
  \]
  with $\dim \Ext^1(M,\sO_X)=\dim \Ext^1(\sO_X,M) = n-l$.  It follows that
  $\dim\Hom(I(-D),\sO_X) = n-l$, and thus the moduli space of maps
  $\sO_X\to M$ can also be viewed as a $\P^{n-l}$-bundle over a subscheme
  of $\Hilb^{D^2/2+l}$.  This gives an upper bound of $D^2+n+l$ on the
  dimension of that moduli space, and subtracting $n-1$ gives an upper
  bound of $D^2+1+l<D^2+2$ on the dimension of the subscheme where
  $\dim\Hom(\sO_X,M)=n$.
\end{proof}

\begin{rem}
  In particular, if $\chi(M)\le 0$, the subscheme of {\em irreducible}
  sheaves with global sections has codimension at least $1-\chi(M)$.
\end{rem}

\chapter{Moduli of sheaves and complexes}
\label{chap:moduli_sheaves_noncomm}

\section{Moduli of simple objects}
\label{sec:moduli_simple_noncomm}

Now that we have established that our surfaces behave reasonably, the most
natural next object of study is the corresponding moduli spaces of sheaves.
The construction via the derived category suggests that a more natural
first problem would be to study moduli spaces in $\perf(X)$.  The most
general result along those lines is the following (specializing a result of
\cite{ToenB/VaquieM:2007} for fairly general dg-categories).

\begin{prop}\cite{ToenB/VaquieM:2007}
  For any family $X/S$ of noncommutative surfaces, there is a derived stack
  ${\cal M}^X$ such that for all derived schemes $T/S$, ${\cal
    M}^X(T)$ is the $\infty$-groupoid of objects in
  $\perf(X_T)$.  This stack is locally geometric and locally of finite
  presentation, and the tangent complex to ${\cal M}^X$ at a geometric
  point $M$ is $R\sEnd(M)[1]$.
\end{prop}

\begin{proof}
  Per the reference, this holds if $\perf(X)$ is {\em saturated}, i.e., if
  it is locally perfect, has a compact generator, is triangulated, and is
  smooth.  The first three properties follow immediately from the
  description via a dg-algebra in $\perf(S)$.  Smoothness follows from the
  fact that gluing preserves smoothness \cite{LuntsVA/SchnurerOM:2014}, along
  with the fact that $\perf(S)$ and $\perf(C)$ are smooth over $S$.
\end{proof}

\begin{rem}
  Of course, the most natural object to study here is the derived stack
  classifying objects of the universal $\perf(X)$ over the moduli stack of
  surfaces.  Since the moduli stack of surfaces is singular, this suggests
  that one should instead work over a suitable derived version of the
  moduli stack of sheaves.  (In particular, the construction sketched in
  Remark \ref{rem:derived_moduli_stack_of_surfaces} to Proposition
  \ref{prop:moduli_stack_of_surfaces_is_small} gives a geometric derived
  stack which is locally of finite presentation.)
\end{rem}

This, of course, is quite far from a stack in the usual sense, and even the
truncation is only a higher Artin stack.  However, we can use this to
obtain algebraic stacks in the usual sense.  Following \cite{LieblichM:2006},
call an object in a geometric fiber ``gluable'' if $\Ext^i(M,M)=0$ for
$i<0$.  By semicontinuity, this is an open condition.  (More precisely, it
is open for each $i$, and $R\sEnd(M)$ is perfect, so locally bounded in
cohomological degree.)  Thus we may define a functor ${\cal M}^X_{ug}$
such that (for $T/S$ an ordinary scheme) ${\cal M}^X_{ug}(T)$ is the
set of quasi-isomorphism classes of objects in $\perf(X_T)$ such that every
geometric fiber is gluable, and this functor will be an open substack of
the truncation of ${\cal M}^X$.

\begin{cor}
  The functor ${\cal M}^X_{ug}$ can be represented by an algebraic stack
  locally of finite presentation.
\end{cor}

\begin{rem}
  Note that the substack of ${\cal M}^X$ classifying flat families of
  coherent sheaves is contained in the derived stack of which ${\cal
    M}^X_{ug}$ is a truncation, since any coherent sheaf on $X/k$ is
  gluable.  We thus find that coherent sheaves on $X$ are classified by an
  algebraic stack locally of finite presentation.  In fact, this is true
  for {\em any} family of $t$-structures on $D_{\qcoh}(X)$.
\end{rem}

Similarly, one calls an object over a geometric fiber ``simple'' if the
natural map
\[
k\to \tau_{\le 0} R\End(M)
\]
is a quasi-isomorphism.  Again, this condition is open, so defines an open
substack ${\cal M}^X_{s}$ on which $\sO_S\to \tau_{\le 0}R\End(M)$ is a
quasi-isomorphism.  Here there is a further reduction one can make: since
the family of automorphism groups of any family with simple fibers is just
$\G_m$ (so flat and locally of finite presentation), ${\cal M}^X_{s}$ is a
gerbe over an algebraic space $\Spl_{X/S}$.  The corresponding functor is
easy enough to define: an element of $\Spl_{X/S}(T)$ is an equivalence
class of simple objects in $\perf(X_U)$ with $U/T$ \'etale, such that the
two base changes to $U\times_T U$ are quasi-isomorphic (but not necessarily
canonically so).  Objects $M$ and $M'$ are equivalent if their base changes
to the fiber product \'etale cover are quasi-isomorphic.

\begin{cor}
  The functor $\Spl_{X/S}$ is represented by a algebraic space locally of
  finite presentation.
\end{cor}

\begin{rem}
  Note that a stable sheaf (relative to any stability condition) is
  automatically simple.
\end{rem}

We have seen in Section \ref{sec:comm_poisson_moduli} that in the
commutative case, the moduli space of simple {\em sheaves} has a Poisson
structure.  Such a Poisson structure is in particular a biderivation, and
the construction of the biderivation carries over to $\Spl_{X/S}$.  To wit,
the cotangent sheaf on $\Spl_{X/S}$ is given by $\sExt^1(M,\theta M)$ where
$M$ is the universal family (the cotangent complex on the corresponding
derived stack is $R\sHom(M,\theta M)[1]$, and the $\sExt^2(M,\theta M)$
term is killed by the structure morphism of the gerbe), and one has a
bivector defined by
\[
\sExt^1(M,\theta M)\otimes \sExt^1(M,\theta M)
\to
\sExt^1(M,M)\otimes \sExt^1(M,\theta M)
\to
\sExt^2(M,\theta M)
\to
\sO_{\Spl_{X/S}},
\]
where the map $\theta M\to M$ comes from a choice of isomorphism $\theta
M\cong M(-Q)$.  (In other words, the map $\theta M\to M$ is a Poisson
structure on $X$ that vanishes on $Q$.  This is unique modulo scalars, so
such Poisson structures on noncommutative surfaces form a line bundle over
the moduli stack of noncommutative surfaces.  The associated biderivation
scales linearly with the choice of Poisson structure on the surface, so is
associated to the same line bundle.) In order to define a Poisson
structure, this bivector needs to be alternating and satisfy the Jacobi
identity.  Since it is defined via the Yoneda product, it is easy to see
that it is antisymmetric, so (apart from characteristic 2), only the Jacobi
identity remains.  Again, derived techniques apply, with the one caveat
that the derived analogues of Poisson structures have only so far been
defined in characteristic 0.

\begin{thm}\cite{nc-lagrangian}
  Let $X/S$ be a family of noncommutative surfaces over a base of
  characteristic 0 (i.e., $S$ is a scheme over $\Q$).  Then
  the above bivector defines a Poisson structure on $\Spl_{X/S}$.
  Moreover, the Poisson structure restricts to a nondegenerate Poisson
  structure on any geometric fiber of $|^\dL_Q$.
\end{thm}

\begin{rem}
  The argument there is in two steps: first, showing that the map ${\cal
    M}^X\to {\cal M}^Q$ has a ``$0$-shifted Lagrangian structure'', a
  derived analogue of the statement that ${\cal M}^X$ has a Poisson
  structure with symplectic leaves classified by ${\cal M}^Q$.  This
  structure immediately restricts to the open substack classifying simple
  objects, and is preserved by the action of $B\G_m$ on $\perf(X)$ and
  $\perf(Q)$, so induces an analogous structure on $\Spl_{X/S}\to [{\cal
      M}^Q/B\G_m]$ over the 2-stack $BB\G_m$.  This is then shown to induce
  an ordinary Poisson structure on $\Spl_{X/S}$ as described.  A similar
  argument gives Poisson structure on the simple locus in the moduli space
  of objects on $X$ with an $m$-step filtration of their restriction to
  $Q$, generalizing the parabolic bundles of \cite{BottacinF:2000}.
\end{rem}

By Lemma \ref{lem:closure_is_Poisson}, we deduce the following.

\begin{cor}
  Let $X/S$ be the universal family of noncommutative surfaces.  Then the
  above bivector defines a Poisson structure on every irreducible component
  of $\Spl_{X/S}$ that contains points of characteristic 0.
\end{cor}

Note that if $M$ is a simple object of some $X_k$ with $\Ext^p(M,M)=0$ for
$p>1$, so that $\Spl_{X/S}$ is formally smooth at $M$, then $M$
automatically lifts to characteristic 0 (since we have seen that $X_k$
itself lifts), and thus the claim holds on the component containing $M$.
In particular, the claim holds for the Zariski closure of the locus of
simple sheaves with $\Ext^2(M,M)=0$, or equivalently $\Hom(M,M(-Q))=0$,
thus for any simple sheaf transverse to $Q$.  We can then use the trick of
\cite{HurtubiseJC/MarkmanE:2002b} (which, when formulated as in Section
\ref{sec:moduli_hd1} works equally well in the noncommutative case; the
only difference is that one ends up with reflexive sheaves rather than
locally free sheaves, but those are still unobstructed) to extend this
result to arbitrary simple sheaves.

\begin{cor}
  Let $X/S$ be a family of noncommutative surfaces.  Then the moduli space
  of simple sheaves on $X$ has a Poisson structure given by the above
  bivector.
\end{cor}

Although this argument is (modulo the heavy machinery) considerably cleaner
than our original arguments that the bivector is Poisson (based on reducing
to the case of reflexive sheaves on maximal orders), it is still presumably
not the ``right'' proof of the Poisson structure.  The ideal way to exhibit
a Poisson structure on an algebraic space is to produce an explicit
noncommutative deformation.  In light of the discussion following Theorem
\ref{thm:weird_langlands}, we expect that the symplectic structure on the
moduli space of 1-dimensional sheaves disjoint from $Q$ is the
semiclassical limit of a noncommutative deformation of the moduli space
(with parameter $\Pic^0(Q)/\Pic^0(C)$), and that there should be derived
equivalences between these deformations extending the derived
autoequivalences of the abelian fibration arising in the fully commutative
case.

If $F:\coh X\to \coh X'$ is an equivalence of categories, then there is a
canonical natural isomorphism $F\theta \cong \theta F$, and thus we can
transport the anticanonical natural transformation $\theta\to \text{id}$
through $F$.  In particular, we obtain a notion of a {\em Poisson}
equivalence, in which the natural transformations $\theta F\to F$ and
$F\theta\to F$ agree modulo the canonical natural isomorphism.  We then
find that a Poisson equivalence of categories induces a Poisson isomorphism
between the moduli spaces of simple objects.  (Something similar applies
for the derived moduli stacks, though one should caution that the analogous
statement for {\em derived} equivalences are more delicate.  Also, we
should really be saying ``bivector-preserving'', since we have not shown
the Jacobi identity holds on the entire finite characteristic locus.)  When
$X=X'$, we note that the anticanonical natural transformation is uniquely
determined by the restriction-to-$Q$ functor and a linear functional on
$H^1(\sO_Q)/H^1(\sO_C)$ (the image of the map $H^1(\sO_Q)\to
H^2(\omega_X)$), and thus an autoequivalence of $X$ is Poisson iff the
induced automorphism of $Q$ acts trivially on $H^1(\sO_Q)/H^1(\sO_C)$.
Note that there will typically be no nontrivial automorphisms of $Q$ fixing
$q$, and thus the map being Poisson will be automatic.

Now, given a finite group $G$, we may consider an action of $G$ on $\coh X$
in the following weak sense: to each $g\in G$, we associate a Poisson
autoequivalence $F_g$ of $\coh X$, in such a way that $F_e=\text{id}$ and
$F_gF_h$ is isomorphic to $F_{gh}$.  Then $G$ acts on $\Spl_X$ preserving
the Poisson structure, and the fixed subspace of this action will inherit a
Poisson structure.  (Note that without additional consistency conditions on
the natural isomorphisms, this does {\em not} give an action of $G$ on the
derived moduli stack or its truncation.)

In the commutative case, such an action of $G$ on $\coh X$ is determined by
a combination of its action on point sheaves and a class in
$H^1(G;\Pic(X))$, determining its action on the structure sheaf.  Then a
simple sheaf $M$ corresponds to a $G$-fixed point of $\Spl_X$ iff for each
$g\in G$, there is an isomorphism $M\cong F_gM$.  If we fix a system of
such isomorphisms (which should be the identity for $g=e$), then we have an
induced isomorphism $F_{gh}M\cong M\cong F_h M\cong F_gF_hM$ for each $g$,
$h$, and those isomorphisms will satisfy the obvious consistency
conditions.  In particular, since we can write $F_g M\cong
(g^{-1})^*M\otimes {\cal L}_g$ with ${\cal L}_{gh}\cong {\cal L}_g\otimes
(g^{-1})^*{\cal L}_h$, we find that there are induced choices for the
latter morphisms making them consistent.  In other words, any $G$-fixed
point of $\Spl_X$ promotes the class in $H^1(G;\Pic(X))$ to an equivariant
gerbe (with trivial underlying gerbe), and this makes the corresponding
sheaf $M$ a twisted $G$-equivariant sheaf.  (This choice is not unique, as
we can twist the isomorphisms $M\cong F_g M$ by any class in $H^1(G;k^*)$.)
When the twisting is trivial, these are sheaves on the orbifold quotient
$[X/G]$, which is itself (in characteristic prime to $|G|$) derived
equivalent to a commutative projective surface via the derived McKay
correspondence \cite{BridgelandT/KingA/ReidM:2001}.  (More precisely, it is
derived equivalent to the minimal desingularization of the scheme-theoretic
quotient $X/G$.)

Something similar holds in the noncommutative setting: the only natural
automorphisms of the identity functor are scalars, and thus the obstruction
to making the natural isomorphisms $F_gF_h\cong F_{gh}$ compatible is a
class in $H^3(G;k^*)$, which vanishes as long as there is {\em any} fixed
point in $\Spl_X$.  If that class vanishes, then the different compatible
choices form a torsor over $H^2(G;k^*)$, and each such choice gives a
disjoint (possibly empty) subset of $(\Spl_X)^G$.  Each fixed point then
corresponds to a collection of (twisted) $G$-equivariant sheaves, which
themselves form a torsor over $H^1(G;k^*)$.  We thus see that these fixed
subspaces are closely related to the moduli space of simple (twisted)
$G$-equivariant sheaves, though in addition to the $H^1(G;k^*)$ action, we
must also bear in mind that a simple equivariant sheaf need not have simple
underlying sheaf, and thus the fixed subspaces are at best quotients of
open subsets of the ``correct'' moduli spaces.  (This issue goes away if we
consider the stack version, where we must fix a consistent family of
natural isomorphisms, but should then get the full moduli stack of
$G$-equivariant objects, suitably defined.)  Presumably there is also an
analogue of the derived McKay correspondence in this case as well.

We can also obtain Poisson automorphisms of $\Spl_X$ associated to {\em
  contravariant} equivalences.  To construct such automorphisms, we first
note that for any $d\in \Z$, we have a contravariant derived equivalence
$M\mapsto (R\ad M)[d]$.  This in general changes the surface, by inverting
$q$, but we can sometimes fix this by composing with an abelian equivalence
$\coh \ad X\cong \coh X$, which as above will tend to give a Poisson
autoequivalence.  The one caveat is that $M\mapsto M[1]$ actually negates
the Poisson structure, and thus one must be careful about signs.  (There is
a simple trick that lets one reduce such sign questions to the commutative
case: the sign is a function on the moduli stack of noncommutative
surfaces, which is connected over the moduli stack of anticanonical
commutative surfaces.)

The direct image $R\pi_*$ associated to a birational morphism respects the
Poisson structure on the derived moduli stack (per \cite{nc-lagrangian},
this actually applies to the projection onto any component of a
semiorthogonal decomposition), and thus in particular so do $L\pi^*$ and
$L\pi^!$, with both facts extending to the moduli stack of simple objects,
with the usual exception that for $R\pi_*$, one must restrict to the open
substack where $R\pi_*$ is simple!  The analogous statements for simple
sheaves are then straightforward and include the minimal lift $\pi^{*!}$.

\section{Smoothness of symplectic leaves}
\label{sec:leaves_are_smooth}

Although Poisson structures tend to be less useful in finite characteristic
than in characteristic 0, one particularly nice feature of the commutative
result is the identification of symplectic leaves (especially since in the
pure 1-dimensional case, the symplectic leaves are essentially just moduli
spaces of difference equations with fixed singularities).  Although these
are identified above in characteristic 0, there is one missing ingredient:
the ``leaves'' have nondegenerate Poisson structures, but this only makes
them symplectic when they are smooth!

We thus wish to know that the symplectic leaves are smooth in general
(which in characteristic 2 implies the pairing is alternating, as the
self-pairing gives the obstruction to extending a first-order deformation).
The first hope would be that the corresponding obstructions might vanish.
We can compute this using the derived moduli stack ${\cal M}_{X/k}$.
Indeed, the symplectic leaves (rather, certain gerbes over the symplectic
leaves) are truncations of open substacks of fibers of the obvious map
${\cal M}_{X/k}\to {\cal M}_{Q/k}$.  We can thus compute the tangent
complex $T_M$ to the derived symplectic leaf containing $M$ as the cone:
\[
T_M\to R\Hom(M,M)[1]\to \tau_{\ge 1}R\Hom_Q(M|^{\bf L}_Q,M|^{\bf L}_Q)[1]\to.
\]
(Indeed, $R\Hom(M,M)[1]$ is the tangent complex in ${\cal M}_{X/k}$, and
the other term is the normal complex of the given point, viewed as a higher
Artin stack.)  This makes it easy to compute the homology of $T_M$:
\[
h^p(T_M) = \begin{cases}
  \Ext^{p+1}(M,M) & p<0\\
  \im(\Ext^1(M,\theta M)\to \Ext^1(M,M)) & p=0\\
  \Ext^{p+1}(M,\theta M) & p>0
\end{cases}
\]
Indeed (unsurprisingly), the tangent complex of the symplectic leaf is
self-dual.  Since $M$ is simple, $T_M$ has amplitude $[-1,1]$, but we have
$h^{-1}T_M\cong k$ and thus by duality $h^1T_M\cong k$, giving an
obstruction.

The dimension of the tangent space at $M$ of its symplectic leaf is
$\dim(h^1(T_M))=\chi(T_M)+2$, and the description of $T_M$ as a cone tells
us that
\[
 -\chi(T_M) = \chi(M,M) - \chi(\tau_{\ge 1}R\Hom_Q(M|^{\bf L}_Q,M|^{\bf
    L}_Q)).
\]
(Compare the ``index of rigidity'' of Section \ref{sec:rigidity_comm}.)
The second term is clearly constant on the symplectic leaf, while the first
term depends only on $[M]\in K_0^{\num}(X)$, and is thus at least locally
constant.  In other words, the tangent spaces to symplectic leaves have
locally constant dimension, and thus each component of a symplectic leaf is
either smooth or nowhere reduced.  (The latter can actually happen for the
commutative Poisson surfaces of characteristic 2 that are not fibers of the
moduli stack of noncommutative rational or rationally ruled surfaces.
Though here ``symplectic leaf'' is particularly misleading, as we do not
expect the bivector to give a Poisson structure.)

Luckily, there are two cases in which we can prove smoothness, which
together are large enough that we have been unable to construct an object
to which neither case applies.

\begin{prop}\label{prop:leaves_are_smooth_nc1}
  Let $X/k$ be a noncommutative rational or rationally ruled surface over an
  algebraically closed field, and let $N\in \perf(Q)$ be any nonzero object.
  Then the fiber over $N$ in $\Spl_{X/k}$ is smooth.
\end{prop}

\begin{proof}
  We need to show that for any local Artin algebra $A$ with residue field
  $k$, any small extension $A^+$ of $A$, and any $M\in \perf(X_A)$ with
  $M|^{\bf L}_Q\cong N\otimes_k A$, there is an extension $M^+$ of $M$ to
  $A^+$ with $M^+|^{\bf L}_Q\cong N\otimes_k A^+$.

  We first consider the problem of extending $M$, ignoring for the moment
  the constraint on its restriction to $Q$.  The obstruction to extending
  $M$ is given by a class in $\Ext^2(M_k,M_k)$, and the image of this class
  in $\Ext^2_Q(N,N)$ is the obstruction to extending $M|^{\bf L}_Q$.  Since
  $M|^{\bf L}_Q$ is a trivial deformation, there {\em is} no such
  obstruction, and thus the obstruction to extending $M$ is in the image of
  $\Ext^2(M_k,\theta M_k)$.  The map $\Ext^2(M_k,\theta M_k)\to
  \Ext^2(M_k,M_k)$ vanishes, since it is dual to $\Hom(M_k,\theta M_k)\to
  \Hom(M_k,M_k)$, which vanishes since $\Hom(M_k,M_k)\to \Hom_Q(N,N)$ is
  injective.  It follows that the obstruction to extending $M$ is trivial,
  and thus it has extensions, classified by a torsor $T$ over
  $\Ext^1(M_k,M_k)$.

  It remains to show that at least one of those extensions restricts to
  $N\otimes_k A^+$.  The map taking an extension of $M$ to its restriction,
  viewed as a deformation of $N$, is a torsor map $T\to \Ext^1_Q(N,N)$
  compatible with the homomorphism $\Ext^1(M_k,M_k)\to \Ext^1_Q(N,N)$, and
  thus its image is a coset of the image of $\Ext^1(M_k,M_k)$.  We thus
  find that the image of $T$ contains $0$ iff the composition
  \[
  T\to \Ext^1_Q(N,N)\to \Ext^2(M_k,\theta M_k)
  \]
  is 0.  The composition
  \[
  \Ext^1_Q(N,N)\to \Ext^2(M_k,\theta M_k)\cong k
  \]
  is dual to the morphism
  \[
  k\to \Hom_Q(N,N)
  \]
  taking $1$ to $\id$, and thus has the alternate expression
  \[
  \begin{CD}
  \Ext^1_Q(N,N) @>\Tr>> \Ext^1_Q(\sO_Q,\sO_Q) \to \Ext^2(\sO_X,\theta
  \sO_X)\cong k.
  \end{CD}
  \]
  The trace of a deformation of complexes on a commutative scheme is the
  deformation of its determinant, so that we finally reduce to showing that
  $\det(M^+|^{\bf L}_Q)\otimes \det(N\otimes_k A^+)^{-1}$ is in the image
  of $\Pic^0(C)(A^+)$ for any $M^+$.  But this follows from the fact that
  determinant-of-restriction gives a well-defined map from $K_0^{\num}(X)$
  to $\Pic(Q)/\Pic^0(C)$.
\end{proof}

\begin{rem}
  We used much the same argument for vector bundles on commutative surfaces
  in Lemma \ref{lem:leaves_of_vect_are_smooth}; the main new
  complication is that $\Ext^2(M,M)$ is no longer guaranteed to vanish.  Of
  course, since it does not actually contribute to the obstructions of the
  symplectic leaf, it is unsurprising that we can show that we never see
  nontrivial classes in $\Ext^2(M,M)$.
\end{rem}

When $M|^{\bf L}_Q=0$, we have $\Ext^2(M,\theta M)\cong\Ext^2(M,M)$, and
thus we can no longer get started in the above argument.  The key idea for
dealing with (nearly all of) these remaining cases is that in any
deformation theory problem, the moduli space of formal deformations is a
quotient of the completion at the origin of the affine space of
infinitesimal deformations, and the obstruction space gives a bound on the
number of elements required to generate the ideal.  In particular, since
our obstruction space is $1$-dimensional, the ideal must be principal
(hopefully 0).  Directly showing that the ideal is 0 can be difficult (as
it is unexpected behavior), but when the ideal is nonzero, controlling the
ideal reduces to showing that certain deformations do {\em not} extend,
giving us an alternate angle of attack.

Of course, since our objective is to prove the ideal is 0, this may not
seem particularly helpful.  However, because we are now asking for $M|^{\bf
  L}_Q$ to be 0, we can embed the symplectic leaf in a family of moduli
problems by allowing the surface to vary.  The argument for $M|^{\bf
  L}_Q\ne 0$ suggests that we should look at $\det M|^{\bf L}_Q$, or rather
its image in $\Pic^0(Q)/\Pic^0(C)$.  This can be computed in two ways.  On
the one hand, it is $\det(0)=\sO_Q \Pic^0(C)$ by our assumption on the
restriction of $M$.  On the other hand, we can also compute $\det(M)$ via
its class in $K_0^{\num}(X)$.  Thus if we allow $X$ to vary in such a way
that $\det([M])$ varies, then the resulting moduli problem {\em will} be
visibly obstructed, but the obstruction will necessarily be in the
direction of the deformation of surfaces.
  
\begin{lem}
  Let $X/k$ be a noncommutative surface over an algebraically closed field,
  and let $[M]\in K_0^{\num}(X)$ be a nonzero class such that $[M]|_Q$ is
  numerically trivial.  Then there is a dvr $R$ with residue field $k$ (and
  field of fractions $K$) and an extension $X'/R$ (as a split family) of
  $X_l$ to $R$ such that $\det([M]_K|_{Q_K})\notin \Pic^0(C_K)$.
\end{lem}

\begin{proof}
  Of course, if $\det([M]|_Q)\notin \Pic^0(C)$, then we may simply take the
  trivial extension.  Otherwise, we may consider the (smooth!) locally
  closed substack of the moduli stack of surfaces determined by the
  combinatorial type of $Q$ (and, in the rational case, whether $\Pic^0(Q)$
  is elliptic, multiplicative, or additive).  It suffices to show that
  $\det([M]|_Q)\in \Pic^0(C)$ cuts out a nontrivial divisor in this
  substack, since then any point of the divisor can be extended to a dvr
  not contained in the divisor.  In particular, we may restrict our
  attention to the corresponding stack in characteristic 0.
  
  If $c_1([M])\cdot e_m\ne 0$, then $e_m$ cannot be a component of $Q$, so
  that $x_m$ is a smooth point of the anticanonical curve of $X_{m-1}$, and
  varying $x_m$ varies the determinant.  If $c_1([M])\cdot e_m=0$, then we
  may reduce to the blowdown, and thus reduce to ruled surfaces and planes.
  If $c_1([M])\ne 0$, then $X$ must be a ruled surface, and then $c_1([M])$
  is fixed up to a scalar multiple by the condition $c_1([M])\cdot Q=0$.
  We must also have $\hat{Q}=\overline{Q}$, and by varying the invertible
  sheaf (again possibly including a lift to characteristic 0) on $\bar{Q}$,
  we may again vary the determinant.  Finally, if $c_1([M])=0$, then the
  only way $\det([M]|_Q)$ can be 0 is if $X$ has center of order dividing
  $\chi([M])\ne 0$, which is cut out by a nontrivial codimension 1
  condition on $q$.
\end{proof}

One caveat is that in the above argument, we may need to allow the surface
to ``vary'' over a mixed characteristic dvr, which means we are departing
somewhat from classical deformation theory.  In addition, there are
representability issues in general which we can sidestep, but only by using
the fact that our moduli problem is represented by a nice derived algebraic
space.

\begin{prop}
  Let $S$ be an (underived) regular algebraic space, and let $X/S$ be a
  derived algebraic space locally of finite presentation.  Then for any
  point $x\in X$ lying over $s\in S$, one has
  \[
  \dim_x(X)
  \ge
  \dim_s(S)
  +
  \dim h^0({\mathbb L}_{X/S}\otimes^\dL k(x))
  -
  \dim h^1({\mathbb L}_{X/S}\otimes^\dL k(x)).
  \]
  If equality holds, then $X$ is a local complete intersection at $x$.
\end{prop}

\begin{proof}
  Since $t_0(X)\to X$ is $1$-connected, and thus $h^p{\mathbb
    L}_{t_0(X)/X}=0$ for $p\ge -1$, replacing $X$ by its truncation can
  only improve the bound, so that we may assume that $X$ is an underived
  algebraic space.  Since \'etale maps respect both dimension and cotangent
  complexes, we may further reduce to the case that both spaces are affine
  schemes.  Moreover, since ${\mathbb L}_{\A^n_S/S}\cong \sO_{\A^n_S}^n$,
  any extension of the map to a map $X\to \A^n_S$ produces the same lower
  bound on $\dim_x(X)$.  In particular, for $n\gg 0$, there is such an
  extension which is a closed immersion, letting us reduce to that case.
  We may then of course replace $S$ by the local ring at $s$.

  It thus remains only to show that if $R$ is a regular local ring and
  $I\subset {\mathfrak m}$ is an ideal, then
  \[
  \dim(R/I)
  \ge
  \dim(R)
  -
  \dim_{R/{\mathfrak m}} (h^1{\mathbb L}_{(R/I)/R}\otimes_{(R/I)} (R/{\mathfrak m})).
  \]
  But $h^1{\mathbb L}_{(R/I)/R}\cong I/I^2$, so this becomes
  \[
  \dim(R/I)
  \ge
  \dim(R)
  -
  \dim_{R/{\mathfrak m}} I\otimes_R (R/{\mathfrak m}).
  \]
  The preimages of any basis of $I\otimes_R (R/{\mathfrak m})$ generate $I$,
  and thus $I$ has a generating set of size
  \[
  \dim_{R/{\mathfrak m}} I\otimes_R (R/{\mathfrak m}),
  \]
  immediately implying the dimension bound.
\end{proof}

\begin{rems}
  Note that although truncating $X$ can only improve the bound, it also
  tends to make the bound significantly harder to compute.  In particular,
  the base change of the cotangent complex is the cotangent complex of the
  {\em derived} base change, which is far easier to understand for a nice
  derived algebraic space.  In particular, for the untruncated version of
  $\Spl_{X/S}$, the cotangent complex of the derived base change is nothing
  other than the shifted dual of the endomorphism complex.
\end{rems}

\begin{rems}
  One of course also has an upper bound $\dim_x(X)\le
  \dim_s(S)+\dim h^0({\mathbb L}_{X/S}\otimes k(x))$, with equality
  implying that $X/S$ is smooth at $x$, but we only use this in the
  well-known case that $S$ is a field.
\end{rems}

One useful observation is that $M|^{\bf L}_Q=0$ is an open condition on any
Noetherian substack of $\Spl_{X/k}$ (the complement of the supports of the
cohomology sheaves, only finitely many of which can be nonzero on the
substack).

\begin{prop}\label{prop:leaves_are_smooth_nc2}
  Let $X/k$ be a noncommutative surface over an algebraically closed field.
  The fiber of $\Spl_{X/k}$ over $0$ is smooth, except possibly on those
  components with trivial class in the Grothendieck group.
\end{prop}

\begin{proof}
  Choose $M\in \Spl_{X/k}$ with $[M]\ne 0$ and $M|^{\bf L}_Q=0$, let $X'/R$
  be the extension of $X$ given by the Lemma, and consider the localization
  ${\cal S}$ of $\Spl_{X'/R}$ at $M$.  Since $\dim(R)=1$, we conclude that
  \[
  \dim({\cal S})\ge 1+\dim\Ext^1(M,M)-\dim\Ext^2(M,M)=\dim\Ext^1(M,M).
  \]
  Since $\det([M]|_Q)$ is nontrivial over $K$ and ${-}|^{\bf L}_Q=0$ on
  ${\cal S}$, we conclude that ${\cal S}(K)=0$, and thus that ${\cal S}$ is
  annihilated by some power of the maximal ideal of $R$.  It follows that
  $\dim({\cal S})=\dim({\cal S}_k)$, and thus $\dim({\cal S}_k)\ge
  \dim\Ext^1(M,M)$.  Since $\Ext^1(M,M)$ is the tangent space at $M$,
  ${\cal S}_k$ is smooth at $M$ as required.
\end{proof}

\begin{rem}
  It is not clear if it is even possible for a simple object to have
  trivial class in the Grothendieck group (let alone also have trivial
  restriction to $Q$).  Note in particular that such an object satisfies
  $\chi(M,M)=\chi(0,0)=0$, and thus any such connected component of
  $\Spl_{X/k}$ is either a smooth surface or a nowhere reduced curve.  Note
  that to show smoothness in finite characteristic, it would suffice to
  show smoothness in characteristic $0$: if $M$ lifts to characteristic
  $0$, then semicontinuity of fiber dimension forces its component to be
  smooth, and if it does not lift, the deformation theory argument shows
  that its component is smooth!
\end{rem}

The fact that symplectic leaves are smooth lets us produce analogues of the
results of Section \ref{sec:rigidity_comm} on {\em rigid} sheaves; in
particular, it tells us that a rigid sheaf is necessarily infinitesimally
rigid.  We have already mentioned that the index of rigidity is nothing
other than the (negative) of the Euler characteristic of the tangent
complex to the derived symplectic leaf, and thus a simple object is
infinitesimally rigid iff its index of rigidity is $2$.  Moreover, if $r$
is a simple root in the system of effective roots, then any sheaf of the
form $\sO_r(d)$ is rigid.  The converse is more delicate than in the
commutative case (the argument used supports!), though one finds similarly
that if $M$ cannot be separated from $Q$ by minimal lifts, then $M$ cannot
be rigid.  We thus reduce to the case $M$ disjoint from $Q$, where
infinitesimal rigidity forces $c_1(M)^2=-2$.  In particular, either $M$ is
irreducible, or it has an irreducible quotient of the form $\sO_r(d)$ with
$r$ a root, and one finds by Euler characteristic considerations that there
is a sheaf numerically equivalent to $M$ in which $\sO_r(d)$ is a subsheaf
and having one fewer map to $\sO_r(d)$.  (I.e., there is a non-split
extension in the other direction, since the intersection number is
positive.)  We thus see that it remains true in the noncommutative case
that rigid sheaves correspond to simple effective roots.

\section{Moduli of semistable sheaves}
\label{sec:semistable_noncomm}

We now turn to the more classical types of moduli spaces, classifying {\em
  stable} (or semistable) sheaves.  The usual argument of Langton
\cite{LangtonSG:1975} shows that the stable moduli space (an open subspace
of $\Spl_{X/S}$) is separated and the semistable moduli space is proper,
but of course one expects in general that they should be
(quasi-)projective.  The standard construction involves inequalities proved
by induction involving general hyperplane sections, and thus is difficult
to extend in general.  (Though of course for rank 1 sheaves, we have
already shown this above!)  This can be finessed for $1$-dimensional
sheaves, however, by replacing Hilbert polynomials of hyperplane sections
by analogous polynomials coming from $\Ext$s from the structure sheaves of
the hyperplane sections, which we may then replace by globally generated
$1$-dimensional sheaves; the resulting bound is somewhat weaker than the
standard bound in the commutative case, but is strong enough to make the
arguments work.

\begin{lem}\label{lem:old_11_42}
  Let $X/k$ be a noncommutative rational or rationally ruled surface, and
  let $M\in \coh X$ be globally generated.  Then for any line bundle $L$
  such that $c_1(M)+c_1(L)-Q$ is ineffective and $\Ext^p(L,\sO_X)=0$ for
  $p>0$, we have $\Ext^p(L,M)=0$ for $p>0$.
\end{lem}

\begin{proof}
  Since $M$ is globally generated, we have a surjection $\sO_X^n\to M$ for
  some $n$.  If $n=1$ (the case $n=0$ being trivial), $\rank(M)=1$, this is
  immediate from $M=\sO_X$.  If $n=1$, $\rank(M)=0$, then let $I$ be the
  kernel of the given global section.  Applying $R\Hom(L,{-})$ to the
  corresponding short exact sequence, we find $\Ext^2(L,M)=0$ and
  \[
  \Ext^1(L,M)\cong \Ext^2(L,I)\cong \Hom(I,\theta L)^*
  \]
  But $I$ is torsion-free of rank 1, so a nonzero morphism $I\to \theta L$
  must be injective, but $c_1(\theta L)-c_1(I)=c_1(M)+c_1(L)-Q$ is
  ineffective.

  For $n>1$, choose one global section, and let $M_1$ be its image.  Then
  $c_1(M_1)-c_1(L)-Q\le c_1(M)-c_1(L)-Q$ is ineffective, so
  $\Ext^p(L,M_1)=0$ for $p>0$ from the $n=1$ case, while
  $\Ext^p(L,M/M_1)=0$ for $p>0$ by induction since $c_1(M/M_1)-c_1(L)-Q\le
  c_1(M)-c_1(L)-Q$ is ineffective.
\end{proof}  

\begin{rem}
Of course, the condition that $\Ext^p(L,\sO_X)=0$ holds for rational
surfaces if $-c_1(L)$ is nef and $-c_1(L)\cdot Q>0$, and for rationally ruled
surfaces of higher genus if $-c_1(L)-(2g-1)f$ is nef.
\end{rem}

Call an ample divisor class $D_a$ ``strongly ample'' if $D_a-2gf$ is nef,
and, if $g=0$, $D_a\cdot Q\ge 2$.  Clearly, every ample divisor class has a
positive multiple which is strongly ample, and replacing it by the latter
leaves the associated stability condition unchanged.

\begin{lem}
  Let $D_a$ be a strongly ample divisor class, and let $M$ be a
  $D_a$-semistable $1$-dimensional sheaf with $\chi(M)>0$ and $c_1(M)\cdot
  D_a=r>0$.  Then $M(rcD_a)$ is acyclic and globally generated, with
  $c=\lfloor \frac{r}{D_a^2}+1\rfloor$.
\end{lem}

\begin{proof}
  For $n\ge 0$, let $M_n$ denote the image of the natural map
  \[
  \sO_X(-ncD_a)\otimes_k \Hom(\sO_X(-ncD_a),M)\to M.
  \]
  Since $M_n(ncD_a)$ is globally generated by construction, $M_n((n+1)cD_a)$
  is again globally generated, and is acyclic by Lemma \ref{lem:old_11_42}.
  It follows immediately that $M_{n+1}$ contains the maximal sheaf
  $M_n\subseteq \tilde{M}_n\subseteq M$ such that $\tilde{M}_n/M_n$ is
  $0$-dimensional.  Moreover, if $M\ne \tilde{M}_n$, then
  \[
  \chi(\sO_X(-(n+1)cD_a),M/\tilde{M}_n)>0,
  \]
  so there is a map to the quotient, implying that $M_{n+1}$ is strictly
  larger than $\tilde{M}_n$.  It follows that
  \[
  c_1(M_{n+1})\cdot D_a \ge c_1(\tilde{M}_n)\cdot D_a=c_1(M_n)\cdot D_a,
  \]
  with equality iff $c_1(M_n)\cdot D_a=c_1(M)\cdot D_a$.  Since
  $c_1(M_0)\ne 0$, we conclude that $c_1(M_n)\cdot D_a\ge \min(D\cdot
  D_a,n+1)$, and thus $c_1(M_{D\cdot D_a-1})=c_1(M)$.  It follows that
  $M/M_{D\cdot D_a-1}$ is $0$-dimensional.  Since an extension of a
  $0$-dimensional sheaf by a globally generated acyclic sheaf is globally
  generated and acyclic, we conclude that $M(c(D\cdot D_a)D_a)$ is acyclic
  and globally generated as required.
\end{proof}

\begin{cor}\label{cor:funny_bound}
  Let $D_a$ be a strongly ample divisor class.  If $M$ is a
  $D_a$-semistable 1-dimensional sheaf on $X$ with $c_1(M)\cdot D_a=r>0$
  and $\chi(M)>r^2(r+1)$, then $M$ is acyclic and globally generated.
\end{cor}

\begin{proof}
  Let $c=\lfloor \frac{r}{D_a^2}+1\rfloor$, so that $(cD_a)^2>c_1(M)\cdot
  (cD_a)$.  Thus to show that $M$ is acyclic and globally generated, it
  suffices to show that $\chi(M(-rcD_a))>0$, or equivalently that
  \[
  \chi(M)>r^2c.
  \]
  Since $c\le r+1$, the claim follows.
\end{proof}

\begin{lem}
  If $M$ is a pure $1$-dimensional sheaf of Chern class $df$ on a
  noncommutative quasi-ruled surface, then there is a filtration $M_i$ of
  $M$ such that each subquotient $M_{i+1}/M_i$ is a pure 1-dimensional
  sheaf of Chern class $f$ and such that
  \[
  \chi(M_1)\ge \chi(M_2/M_1)\ge\cdots\ge\chi(M/M_{d-1}).
  \]
\end{lem}

\begin{proof}
  Since $M$ is pure $1$-dimensional, there is some twist of $M$ that has no
  global sections; the claim is invariant under twisting, so we may assume
  $H^0(M(-s))=0$.  (Twisting by $-s$ is normally only defined up to the
  pullback of a line bundle on $C$, but $\rho_{-1*}M$ is $0$-dimensional,
  so invariant under such twists, and thus whether $H^0(M(-s))=0$ is
  independent of that choice.)  Since $\chi(M(ns))>0$ for $n\gg 0$, there
  is some smallest $n_1$ such that $M(n_1s)$ has a global section, and
  again we may assume $n_1=0$.  Now, let $N$ be the image of the natural
  map
  \[
  \sO_X\otimes_k \Hom(\sO_X,M)\to M.
  \]
  Since $c_1(N)$ and $df-c_1(N)$ are effective and orthogonal to the nef
  divisor class $f$, we see that $c_1(N)=d_1f$ for some $1\le d_1\le d$.
  We claim that $\chi(N)\ge d_1$.  Indeed, the subquotients of the
  filtration of $N$ coming from any saturated filtration of $\Hom(\sO_X,M)$
  must again have Chern class proportional to $f$, and each subquotient is
  a quotient of $\sO_X$, so we can bound its Euler characteristic using
  Corollary \ref{cor:Ox_quotient_bound}.  Since $H^0(N(-s))=0$, we also
  have $\chi(N(-s))\le 0$, implying $\chi(N)=d_1$, making the bound tight.
  It follows that any saturated filtration of $\Hom(\sO_X,M)$ induces a
  filtration of $N$ with subquotients of the form $\sO_X/L$ with $L$ a line
  bundle of class $-d'f$.  Each such subquotient is the pullback of a
  $0$-dimensional sheaf on $\P^1$, and thus admits a refined filtration by
  pullbacks of points.

  In particular, we conclude that $N$ is acyclic, and thus $M/N$ has no
  global sections, so is pure.  By induction, $M/N$ has a filtration as
  desired, and $H^0(M/N)=0$ implies that the Euler characteristic of the
  bottom subsheaf is nonpositive, and thus strictly less than the Euler
  characteristics of the subquotients of the filtration of $N$.  Gluing the
  two filtrations gives the desired result.
\end{proof}

\begin{cor}
  Let $X$ be a noncommutative rational or rationally ruled surface over an
  algebraically closed field $k$, and let $D$ be a divisor class such that
  $D-2gf$ is nef and, if $g=0$, $D=0$ or $D\cdot Q\ge 2$.  Then for any
  pure $1$-dimensional sheaf $M$ and any line bundle $L$ of class $-D$,
  there is a homomorphism $\phi:L\to \sO_X$ such that
  $\Hom(\coker\phi,M)=0$.
\end{cor}

\begin{proof}
  We first note that for any filtration of $M$, if $\coker\phi$ maps to
  $M$, then it maps to some subquotient, and thus we may assume without
  loss of generality that $M$ is irreducible.  In particular, the image of
  a nonzero map $\coker\phi\to M$ is a globally generated sheaf with Chern
  class $c_1(M)$, so has Euler characteristic at least $-c_1(M)\cdot
  (c_1(M)+K_X)/2$.  Now, for any point $x\in Q$, $(\coker\phi|_Q)\otimes
  k(x)=0$ for $\phi$ in the complement of a proper subspace of
  $\Hom(L,\sO_X)$, and thus we see that $\Hom(\coker\phi,M|_Q)$ is 0
  outside a finite union of proper subspaces of $\Hom(L,\sO_X)$.  It
  follows that, away from those subspaces, any map $\coker\phi\to M$
  factors through the image of $M(-Q)$ in $M$.  Iterating, we see that for
  $\phi$ outside a (possibly larger) finite union of proper subspaces, any
  map $\coker\phi\to M$ factors through $M(-nQ)$.  Since the Euler
  characteristic of the image in $M$ is uniformly bounded below, we get a
  contradiction unless the image of $M(-nQ)$ eventually stabilizes, which
  can only happen if the image of $\coker\phi$ is disjoint from $Q$.

  In other words, we may assume WLOG that $M$ is irreducible and disjoint
  from $Q$ (in particular ruling out noncommutative planes and odd
  Hirzebruch surfaces).  For any $l\ge 0$ such that $D-lQ$ satisfies the
  hypotheses, we may consider those maps $\phi$ that factor through
  $\theta^{-l}L$, and thus reduce to $D-lQ$.  In particular, if $m\ge 2$
  and we choose a blowdown structure putting $D$ in the fundamental
  chamber, then we may take $l=D\cdot e_m$ and thus reduce to the case that
  $D\cdot e_m=0$; similarly, for $m=1$, we may take $l=\min(D\cdot
  e_1,D\cdot (f-e_1))$ and reduce similarly (up to an elementary
  transformation).  Disjointness from $Q$ implies that $M$ is
  $\alpha_{m*}$-acyclic with pure $1$-dimensional direct image, and thus we
  may reduce to the corresponding statement on $X_{m-1}$.
  
  Since we have already ruled out noncommutative planes and odd ruled
  surfaces, this reduces to the case of even ruled surfaces, for which
  either $g=0$, $Q$ is integral, or $g\ge 0$ and $Q$ is the sum of two
  (possibly identical) curves of class $s-(g-1)f$.  We can then perform a
  similar reduction, asking for $\phi$ to vanish along an integral
  component of $Q$, each time reducing $D\cdot f$ by 2 in the first case
  and 1 in the second case.  Since $f$ is nef, $(D-c_1(M))\cdot f\ge 0$,
  and $c_1(M)\propto s+(g-1)f$, the only possibility after fully reducing
  $D\cdot f$ is $g=0$, $D\cdot f=1$, $c_1(M)=s-f$, and $Q$ integral.  In
  that case, the kernel of a map $\coker\phi\to M$ has Chern class a
  multiple of $f$, so is an extension of sheaves of Chern class $f$,
  implying that the points in the support of $\coker\phi|_Q$ can be paired
  so that their products are in a fixed coset of $q^\Z$.  But the support
  of $\coker\phi|_Q$ is nothing other than the zero locus of a section of
  a line bundle of degree $\ge 4$ on the integral genus 1 curve $Q$, and
  thus this condition fails generically.
\end{proof}

These results imply an analogue of the Le Potier-Simpson bound.

\begin{lem}
  Let $D_a$ be a strongly ample divisor class.  If $M$ is a semistable
  $1$-dimensional sheaf with $c_1(M)\cdot D_a=r$, then $h^0(M)\le
  \max((r+1)^3+\chi(M),0)$.
\end{lem}

\begin{proof}
  By the analogue of Cohen-Macaulay duality, we may instead prove that
  \[
  h^1(M)\le \max((r+1)^3-\chi(M),0).
  \]
  If $\chi(M)>r^2(r+1)$, this follows from Corollary \ref{cor:funny_bound}.
  Otherwise, let $c=\lceil (r+1)^2-\chi(M)/r\rceil$, so that for any line
  bundle $L$ of class $-cD_a$, $\chi(L,M)\ge r(r+1)^2>r^2(r+1)$ and thus
  $R\Hom(L,M)$ is acyclic.  Let $\phi:L\to \sO_X$ be a map such that
  $\Hom(\coker\phi,M)=0$, guaranteed to exist by the previous Lemma.
  Acyclicity of $R\Hom(L,M)$ ensures that $\Ext^2(\coker\phi,M)\cong
  \Ext^2(\sO_X,M)=0$ and that $\Ext^1(\coker\phi,M)\to \Ext^1(\sO_X,M)$ is
  surjective, so that
  \[
  h^1(M)\le \dim\Ext^1(\coker\phi,M)
  = -\chi(\coker\phi,M)
  = cr \le (r+1)^3-\chi(M)
  \]
  as required.
\end{proof}

Although weaker than the bound one can derive in the commutative case
(which is only quadratic in $r$), this is still good enough to let us adapt
the commutative argument.  In particular, the argument of
\cite[Cor.~3.3.8]{HuybrechtsD/LehnM:1997} for pure $1$-dimensional sheaves
carries over directly, with the only change being that the constant $C$
should be replaced by $r^2+3r+5$.  We can then follow the argument of
Thm.~4.4.1 op.~cit.~to conclude that a subsheaf $M'$ such that
\[
\frac{\chi(M')}{c_1(M')\cdot D_a}
<
\frac{\chi(M)}{c_1(M)\cdot D_a}
-
\frac{(c_1(M)\cdot D_a)^3-1}{c_1(M)\cdot D_a}
-2
\]
is not destabilizing on the Quot scheme representing $M$ as a quotient of
$\sO_X(-mD_a)$ for (uniformly) sufficiently large $m$.  There are only
finitely many possible values for the remaining numerical invariants; since
the family of sheaves $M$ is bounded, so is the corresponding family of
subsheaves with those numerical invariants.  We may thus find a twist such
that all of those subsheaves are acyclic, and thus every subsheaf has a
sufficiently good estimate for the dimension of its space of global
sections.  It follows that for $m$ uniformly sufficiently large, every
(semi)stable sheaf $M$ is GIT-(semi)stable as a point of the corresponding
Quot scheme, and thus the semistable moduli space is quasiprojective.
Since we have already shown that it is proper, projectivity follows.

This establishes the difficult part of the following result.

\begin{thm}\label{thm:semistable_moduli_are_projective}
  Let $D$ be an effective divisor class and let $D_a$ be an ample divisor
  class.  Then the moduli problem of classifying $S$-equivalence classes of
  $D_a$-semistable sheaves $M$ with $\rank(M)=0$, $c_1(M)=D$ and
  $\chi(M)=l$ is represented by a projective scheme $M_{D,l}$.  The stable
  locus of $M_{D,l}$ is a Poisson subscheme, which is smooth if $D\cdot
  Q>0$, as is the Zariski closure (in the stable locus) of the sublocus
  with $M|_Q=0$ if $D\cdot Q=0$.  If $\Z[M]\subset K_0^{\num}(X)$ is a
  saturated sublattice, then the stable locus admits a universal family.
\end{thm}

\begin{proof}
  Since stable sheaves are simple, the Poisson structure is induced from
  that on $\Spl_{X/S}$.  If $D\cdot Q>0$, then stability implies
  $\Hom(M,\theta M)=0$ (the slope decreases) and thus deformations are
  unobstructed, implying smoothness.  Similarly, if $D\cdot Q=0$,
  then $M$ and $\theta M$ are stable sheaves of the same slope, and thus
  any map between them is an isomorphism.  Semicontinuity then implies that
  $\dim \Ext^2(M,M)=\dim \Hom(M,\theta M)=1$ on the Zariski closure of the
  given symplectic leaf.  In particular, the dimension of the tangent space
  remains unchanged on the Zariski closure, and thus the Zariski closure
  inherits smoothness from the symplectic leaf.

  It remains to show that the stable locus admits a universal family.  It
  certainly has one \'etale locally, since this is tautologically true for
  $\Spl_{X/S}$.  Since $\Z[M]\subset K_0^{\num}(X)$ is saturated, there is
  a class $[N]$ such that $\chi([N],[M])=1$.  Let $N^\cdot$ be a complex of
  class $[N]$.  If $M$ is a universal family over some \'etale cover of the
  stable moduli space, then $M'=M\otimes \det R\Hom(N^{\cdot},M)^{-1}$ is
  again universal.  The twisted descent data for $M$ induces {\em
    untwisted} descent data for $M'$, which therefore descends to an actual
  family of sheaves, giving the desired universal family.
\end{proof}

\begin{rems}
  If $\chi(M)\ne 0$ and $\Z[M]$ is saturated, we can always modify our
  ample line bundle $D_a$ to some new ample line bundle $D'_a$ such that
  every $D_a$-stable sheaf is $D'_a$-stable, as is every $D'_a$-semistable
  sheaf.  To do this, simply observe that there are only finitely many
  possible divisor classes of subsheaves of a sheaf with Chern class
  $c_1(M)$, and thus a discrete set of possible $D_a$-slopes of subsheaves.
  In particular, there is a lower bound on a nonzero difference between the
  slope of a subsheaf and the $D_a$-slope of $M$.  Thus if we could replace
  $D_a$ by $D_a+\epsilon D_0$ for $\epsilon$ sufficiently small, the slopes
  on either side of $\mu(M)$ would remain on the same side after the
  deformation.  Replacing $D_a$ by $D'_a=nD_a+D_0$ for $n\gg 0$ thus has a
  similar effect (effectively using $D_0$ to break ties between sheaves of
  the same slope).  In particular, this ensures that $D_a$-stable sheaves
  are still $D'_a$-stable.  Similarly, there are only finitely many
  possibilities for the projection
  \[
  c_1(N)-\frac{c_1(N)\cdot D_a}{c_1(M)\cdot D_a} c_1(M)
  \]
  for $N$ ranging over subsheaves, and choosing $D_0$ to have nonzero
  intersection with all such vectors ensures that the only way $M$ could be
  strictly $D'_a$-semistable would be if $c_1(N)$ were proportional to
  $c_1(M)$, at which point the primitivity condition forces $N$ to have a
  different $D_a$-slope.  This gives a smooth projective Poisson moduli
  space with a natural morphism to the original semistable moduli space,
  which if $D\cdot Q=0$ is a symplectic resolution.  (We can also produce
  such a resolution when $\chi(M)=0$ by starting with $M(D_a)$ instead.)
\end{rems}

\begin{rems}
  In general, there is an obstruction in the Brauer group of the stable
  moduli space to the existence of a universal family.  (If $[M]$ is
  divisible by $l$ and $[M]/l$ is primitive, then this obstruction has
  order dividing $l$.)  There are examples, even in the commutative case,
  for which the obstruction is nontrivial, e.g., Corollary \ref{cor:Picr}
  or Proposition \ref{prop:irreds_on_fibers_exist} below.
\end{rems}

\begin{rems}
  One significance of having a universal family is that it gives rise to
  Lax pairs.  To be precise, the universal family induces a sheaf over the
  function field of the moduli space, which in turn gives rise to a
  difference/differential equation.  If we twist that sheaf by a line
  bundle, then the new equation is related by a gauge transformation (see
  Section \ref{sec:sheaves_from_eq_nc}), but is also a point on a
  corresponding moduli space, allowing us to identify the function fields
  of the moduli spaces.  Moreover, since the twist of the original
  universal sheaf is isomorphic to the new universal sheaf, simplicity lets
  us make the change of basis unique (up to scalars!).  If we define the
  universal family globally on some open subset of parameter space, then
  this gives a gauge transformation between the universal family and its
  pullback, and this is precisely what we want from a Lax pair.  (As usual,
  there are issues with scalar gauge freedoms, so that the Lax equation is
  only defined up to a scalar.  Again, there are no difficulties if we work
  projectively!)  Without a universal family, one can either work
  projectively (killing the obstruction) or work on a suitable smooth cover
  of the moduli space, though in the latter case the difficulty is lifting
  the isomorphisms of moduli spaces coming from twisting to isomorphisms of
  the covers.  One such possibility (compare \cite{ellGarnier}) that works
  for twists by $D$ with $D\cdot f=0$ is to take the smooth cover to be the
  $\GL_n$-bundle coming by choosing a trivialization of $\rho_*M$ at some
  suitable auxiliary point not in the orbit of any singularities.  One cost
  of taking the latter cover is that the family is no longer over a
  projective base, and thus it will be more difficult to prove algebraic
  integrability as in Section \ref{sec:integrable}.  (It also becomes more
  difficult to parametrize the base of the family and write down the
  universal equation, both of which are needed if one wishes an {\em
    explicit} Lax pair!)
\end{rems}

\section{Painlev\'e cases}
\label{sec:painleve_moduli}

Since the smallest possible nontrivial symplectic leaf is 2-dimensional,
the simplest interesting moduli space will also be 2-dimensional.  If the
primitivity condition is satisfied, so we can arrange for the stable moduli
space to be projective, then it is in fact a smooth projective Poisson
surface, and thus a surface of the type considered above.  Something
stronger is true, in fact.

First, note that if the moduli space of simple 1-dimensional sheaves with
$c_1=D$, $\chi=l$ contains a nontrivial $2$-dimensional symplectic leaf,
then we must have $D\cdot Q=0$ (for disjointness) and $D^2=0$ (since this
is the Euler characteristic of the tangent complex).  We may also assume
$D$ nef (the generic sheaf on any component of the moduli space has a
constant sub- or quotient sheaf corresponding to any $-d$-curve having
negative intersection with $D$), at which point putting it in the
fundamental chamber and blowing down any unnecessary $-1$-curves leaves us
with two cases: either $X$ is rational or $X$ is genus $1$, in either case
with $Q^2=0$ and $D\propto Q$.  These are of course precisely the cases for
which we constructed derived equivalences above.

The key observation is then that when the moduli space is well-behaved, we
can use it to construct a derived equivalence.  (Note that this will in
general not be a derived equivalence as constructed above, as $X$ itself
can have a large group of derived autoequivalences, in particular including
spherical twists by any sheaf of the form $\sO_C(d)$ with $C^2=-2$, $C\cdot
Q=0$.)  We focus on the rational case, as the genus 1 case should be
analogous and leads to less interesting moduli spaces.

We recall the notation used above in Theorem \ref{thm:weird_langlands}:
we fix a (possibly degenerate) commutative del Pezzo surface $Y$ with
anticanonical curve $Q$, and let $X_{z,q}$ be the 2-parameter family of
noncommutative surfaces obtained from $Y$ by blowing up a point and taking
the noncommutative deformation, with parameters $z\cong \det([\sO_Q]|_Q)$ and
$q\cong \det([\pt]|_Q)$.

\begin{lem}\label{lem:Painleve_derived_equivalence}
  Let $X_{z,q}$ be a surface as above, let $r$, $d$ be relatively prime
  integers such that $z^r q^d\cong \sO_Q$, and let $D_a$ be an ample
  divisor class such that any $D_a$-semistable $1$-dimensional sheaf on
  $X_{z,q}$ with $c_1(M)=rQ$, $\chi(M)=d$ is $D_a$-stable.  Then any
  component of the corresponding stable moduli space that contains sheaves
  disjoint from $Q$ is a smooth rational projective surface $Y$ with
  $K_Y^2=0$, and is derived equivalent to $X_{z,q}$ by a functor taking
  structure sheaves of points to $D_a$-stable $1$-dimensional sheaves with
  the given invariants.
\end{lem}

\begin{proof}
  Let $Y$ be such a component, and note that since it is the Zariski
  closure of a 2-dimensional symplectic leaf, it is smooth and projective.
  We first show that it is derived equivalent to $X_{z,q}$.
  
  The fact that $r$ and $d$ are relatively prime ensures that there is a
  universal family $M$ on the base change of $X_{z,q}$ to $Y$.  We may use
  this to construct a Fourier-Mukai functor
  \[
  \Phi_M:D^b\coh Y\to D^b\coh X_{z,q}:
  \]
  given a complex $N^\cdot$ on $Y$, we take the derived tensor product with
  $M$ to obtain a sheaf on the base change of $X_{z,q}$, and then take
  derived direct images in the two pieces of the natural semiorthogonal
  decomposition.  We claim that $\Phi_M$ is actually an equivalence.  Since
  point sheaves are a spanning class for $D^b\coh Y$, we may apply
  \cite[Thm. 2.4]{BridgelandT/KingA/ReidM:2001} to reduce this to showing
  the following:
  \begin{itemize}
    \item[(1)] $\Phi_M$ has both a left and a right adjoint.
    \item[(2)] For all $x\in Y$, $\Phi_M(\sO_x\otimes \omega_Y)\cong \theta
      \Phi_M(\sO_x)$.
    \item[(3)] For any two points $x,y\in Y$ the
      natural morphism
      \[
      R\Hom_Y(\sO_x,\sO_y)\to R\Hom_{X_{z,q}}(\Phi_M(\sO_x),\Phi_M(\sO_y))
      \]
      is a quasi-isomorphism.
    \item[(4)] The category $D^b \coh X_{z,q}$ is indecomposable.
  \end{itemize}
  For (1), we note that the right adjoint applied to $N^{\cdot}\in
  D^b_{\coh}X_{z,q}$ is just $R\sHom_Y(M,N^{\cdot}_Y)$, where $N_Y$ is the
  corresponding complex on the base change; this follows easily from the
  description via the semiorthogonal decomposition and the analogous fact
  for commutative Fourier-Mukai functors.  We can then obtain the required
  left adjoint by conjugating the right adjoint by the respective Serre
  functors.

  For (2) and (3), note first that $\Phi_M(\sO_x)$ is just the
  corresponding fiber of the universal family, and is thus a stable sheaf
  $M_x$ on $X_{z,q}$ with the given invariants.  For (2), we find that
  $\dim\Hom_{X_{z,q}}(M_x,\theta M_x)\ge 1$ by semicontinuity, and any such
  map induces an isomorphism by stability.  For (3) with $x\ne y$, we have
  two nonisomorphic sheaves of the given kind, and stability prevents there
  being a morphism in either direction, so that Euler characteristic
  considerations imply
  \[
  R\Hom_{X_{z,q}}(M_x,M_y)=0=R\Hom_Y(\sO_x,\sO_y).
  \]
  
  For $x=y$, we note first that the map $\Hom(\sO_x,\sO_x)\to
  \Hom(M_x,M_x)$ takes the identity to the identity;
  since both sides are spanned by the identity, this is an isomorphism.
  Next, the map $\Ext^1(\sO_x,\sO_x)\to
  \Ext^1(M_x,M_x)$ is just the usual (Kodaira-Spencer)
  identification of the tangent space to the moduli space with the space of
  infinitesimal deformations.  Finally, the map $R\Hom(\sO_x,\sO_x)\to
  R\Hom(M_x,M_x)$ respects the Yoneda pairing, and
  Serre duality implies that the image of $\Ext^1(\sO_x,\sO_x)\otimes
  \Ext^1(\sO_x,\sO_x)$ spans $\Ext^2(\sO_x,\sO_x)$. Since
  $\dim\Ext^2(\sO_x,\sO_x)=\dim\Ext^2(M_x,M_x)=1$, this establishes (3) in
  this case as well.

  Finally, for (4) (i.e., that $D^b_{\coh} X_{z,q}$ is not a direct sum),
  note that any two line bundles with Chern class differing by an ample
  divisor are in the same component, and thus there is a component
  containing every line bundle.  Since line bundles generate, this
  component must be all of $D^b_{\coh} X_{z,q}$ as required.

  It follows that $Y$ is derived equivalent to $X_{z,q}$.  In particular,
  the natural transformation $\theta\to \text{id}$ on $X_{z,q}$ induces a
  corresponding natural transformation on $Y$, or equivalently an
  anticanonical map $\omega_Y\to \sO_Y$, which cannot be the identity since
  it is not the identity on $X_{z,q}$.  It follows that $Y$ is nontrivially
  Poisson, and is thus birationally ruled.  Moreover, the derived
  equivalence induces an isomorphism of Grothendieck groups, so that
  $K_0(Y)\cong \Z^{10}$.  (Note here that this is $K_0(Y)$, not just
  $K_0^{\text{num}}(Y)$.)  Since the Grothendieck group of a higher genus
  ruled surface over an algebraically closed field is not finitely
  generated, we conclude that $Y$ must be rational, and thus (since
  $K_0(Y)$ has rank 10) that $K_Y^2=0$.
\end{proof}

\begin{rem}
  Note that the derived equivalence is nonunique, since twisting a
  universal family by a line bundle yields another universal family.
\end{rem}

The proof actually tells us something stronger about $Y$, namely that its
anticanonical curve can be identified with a moduli space of sheaves on
$Q$, and is itself derived equivalent to $Q$.  When $Q$ is smooth, this is
enough to tell us that $Y$ also has anticanonical curve $Q$ (when $Q$ is
smooth, $D^b_{\coh} Q\cong D^b_{\coh} Q'$ forces $Q\cong Q'$), but in general
we will need to work somewhat harder.  There is, however, one particularly
tractable special case.

\begin{lem}\label{lem:Painleve_elliptic}
  Let $X$ be a rational Jacobian elliptic surface with generically smooth
  fiber of class $Q$, let $r$, $d$ be relatively prime integers, and let
  $D_a$ be an ample divisor class on $X$ such that any $D_a$-semistable
  $1$-dimensional sheaf with $c_1=rQ$, $\chi=d$ is $D_a$-stable.  Then the
  $D_a$-semistable moduli space of such sheaves is isomorphic to $X$, and
  the associated derived equivalence restricts to a derived autoequivalence
  of every fiber.
\end{lem}

\begin{proof}
  Let $Y$ be the given moduli space.  We first observe that $Y$ is smooth.
  Indeed, any simple sheaf with the given invariants is supported on a
  single fiber of $X$, and thus is disjoint from {\em some} anticanonical
  curve on $X$, at which point Lemma \ref{lem:Painleve_derived_equivalence}
  tells us that the component of $Y$ containing the corresponding point
  must be a smooth projective rational surface derived equivalent to $X$.
  The pencil of natural transformations $\theta\to\text{id}$ can be
  transported through these derived equivalences, and thus $Y$ itself has
  an anticanonical pencil, so every component is a relatively minimal
  elliptic surface.  Moreover, the derived equivalence from $Y$ to $X$
  restricts to a derived equivalence between corresponding anticanonical
  curves.  For any smooth fiber $C$ of $X$, \cite{AtiyahMF:1957} tells us
  that vector bundles of rank $r$ and degree $d$ are determined by their
  determinant, and thus the corresponding anticanonical curve of $Y$ is
  $\Pic^d(C)\cong C$.  (In particular, this implies that $Y$ itself must be
  irreducible, since it has irreducible fibers.)  But this implies that $Y$
  is the minimal proper regular model of $X$ itself, so is canonically
  isomorphic to $X$.
\end{proof}

This suggests that when $X$ has a unique anticanonical curve, that the
moduli space should have isomorphic anticanonical curve.  The key idea is
that the anticanonical curve of the semistable moduli space is itself a
semistable moduli space, and we can obtain the {\em same} moduli space from
a suitable elliptic surface.

\begin{lem}\label{lem:Painleve_embed_in_elliptic}
  Let $X/k$ be a noncommutative rational surface over an algebraically
  closed field, with anticanonical curve $Q$ such that $Q^2=0$ and no
  component of $Q$ has self-intersection $<-2$, and let $r,d$ be relatively
  prime integers with $r>0$.  Also let $D_a$ be an ample divisor class such
  that any $D_a$-semistable 1-dimensional sheaf with $c_1(M)=rQ$,
  $\chi(M)=d$ is $D_a$-stable.  Then there is a commutative Jacobian
  elliptic surface $X'/k$ with a fiber isomorphic to $Q$ and an ample
  divisor class $D'_a$ inducing the same stability condition on sheaves
  supported on $Q$ with the given invariants.
\end{lem}

\begin{proof}
  Embedding $Q$ as a fiber of some $X'$ is straightforward; indeed, in
  characteristic not 2 or 3, we may simply obtain $X'$ from $X$ by passing
  to the corresponding commutative surface, blowing down a $-1$-curve, then
  blowing up the base point of the anticanonical linear system, and
  observing that the resulting genus 1 fibration has smooth components
  since the characteristic is not 2 or 3.  Alternatively (and in arbitrary
  characteristic), one can use Tate's algorithm to write down general
  rational genus 1 fibrations with any given Kodaira type at $\infty$, and
  find that the only cases that (after Tate's change of coordinates) impose
  closed constraints on the coefficients of the fiber at $0$ are $I_7$,
  $I_8$, $I_9$ and $I_4^*$.  Any fibration with a fiber of the first three
  types must also have smooth fibers, and similarly for type $I_4^*$ except
  in characteristic 2, where we note that the elliptic curve
  \[
  E/\F_2(t):y^2+xy=x^3+tx^2+x
  \]
  has fiber $I_4^*$ at $t=\infty$.
 
  Having chosen $X'$, let $D'$ be a divisor on $X'$ that has the same
  restriction to $Q$ as a positive integer multiple of $D_a$.  (In the
  direct construction of $X'$, we may of course take $D'=D_a$, but note
  that this need not be ample on $X'$.)  If we apply the algorithm to put
  $D'$ in the fundamental chamber, then any time we would wish to reflect
  in an effective root, that root cannot be a component of $Q$, and thus
  performing that reflection will have no effect on the induced stability
  condition.  We may thus assume that $D'$ has nonnegative intersection
  with every simple root, and then by adding a suitable multiple of $Q$,
  that it has positive intersection with $e_8$.  If it is not ample, then
  this is because it has intersection $0$ with some component of a fiber
  other than $Q$; if we multiply $D'$ by 3 and then add all components of
  other fibers that it {\em does} intersect, this again has no effect on
  the stability condition, but strictly decreases the set of bad
  components.  We thus need do this only finitely many times before
  eventually getting an ample divisor as desired, which we can then perturb
  as needed to eliminate any remaining strictly semistable sheaves of the
  desired invariants.
\end{proof}

In constructing the derived equivalence, we needed the moduli space to
contain points disjoint from $Q$.  This is a nontrivial condition in
general, but luckily we can use our earlier derived equivalences to find
suitable sheaves.

\begin{lem}\label{lem:Painleve_stable_exists}
  Let $X_{z,q}$ be a surface as above, and $r,d$ relatively prime integers
  such that $z^r q^d\cong \sO_Q$.  Then for any ample divisor class $D_a$,
  there is a $D_a$-stable $1$-dimensional sheaf $M$ on $X_{z,q}$ with
  $c_1(M)=rQ$, $\chi(M)=d$, and $M|^\dL_Q=0$.
\end{lem}

\begin{proof}
  Let $R$ be a dvr with residue field $k$ and field of fractions $K$ with
  $\ch(K)=0$, and let $Y^+/R$ be an extension of the del Pezzo surface $Y$
  to an anticanonical del Pezzo surface over $R$ in such a way that the
  generic fiber of $Q^+$ is smooth.  Then we may lift $q$ and $z$ (via
  their single degree of freedom) to $Y^+$ to obtain a corresponding
  extension $X^+_{z,q}$ of $X_{z,q}$. (Moreover, the ample line bundle on
  the special fiber lifts to an ample line bundle on $X^+_{z,q}$.)

  Now, let $ar-db=1$, so that Corollary \ref{cor:weird_langlands} gives a
  derived equivalence from $X^+_{z^a q^b,1}$ to $X^+_{z,q}$ acting on $K_0$
  by taking $[pt]$ to the desired class $[M]$.  We claim that there exists
  a point of the generic fiber of $X^+_{z^a q^b,1}$ such that the image has
  the desired properties (up to an even shift).  Being a $D_a$-stable sheaf
  and being disjoint from $Q^+$ are both intersections of open conditions,
  and thus it suffices to show that each of the two conditions holds on
  {\em some} point of $X^+_{z^a q^b,1}$.  For disjointness from $Q^+$, this
  is immediate: the open condition is precisely that the point not be on
  $Q^+$.  The stable sheaf condition is more subtle, but since $Q^+$ is
  smooth, we see that the image of any point of $Q^+$ is a simple object of
  $D^b_{\coh}(Q^+)$ with the desired invariants, and is thus an even shift
  of a simple sheaf.  But a simple sheaf on a smooth genus 1 curve is
  irreducible, and thus stable for any ample divisor class.  We in
  particular see that we have a 2-dimensional family of such sheaves on
  $X^+_{z,q}$, and thus have at least a 2-dimensional family of semistable
  sheaves on $X_{z,q}$.

  To see that we have {\em stable} sheaves disjoint from $Q$, consider a
  reducible semistable sheaf with the given invariants and disjoint from
  $Q$.  A subsheaf of Chern class $D$ will itself be reducible unless
  $D^2\in \{0,-2\}$ (as otherwise the sheaf will have nontrivial
  endomorphisms).  There are only finitely many candidate divisor classes
  with $D^2=-2$ and $D$, $rQ-D$ effective, and only finitely many possible
  Euler characteristics for each $D$, and thus only finitely many possible
  sheaves.  In each of those finitely many cases, having the given subsheaf
  is a closed condition, and is a nontrivial condition for the simple
  reason that it fails for sheaves supported on $Q$.  We thus see that any
  subsheaf of the {\em generic} semistable sheaf has Chern class
  proportional to $Q$.  Since $\gcd(r,d)=1$, such a subsheaf cannot have
  the same slope, and thus the generic semistable sheaf is stable as
  required.
\end{proof}

Putting everything together, we find the following.

\begin{thm}\label{thm:Painleve_moduli_spaces}
  With $Y$, $Q$ fixed, let $w\in \Pic^0(Q)$.  Let $r,d$ be relatively prime
  integers with $r>0$, and let $D_a$ be an ample divisor class such that
  any semistable $1$-dimensional sheaf on $X_{w^{-d},w^r}$ with
  $c_1(M)=rQ$, $\chi(M)=d$ is stable.  Then the corresponding semistable
  moduli space is isomorphic to $X_{w,1}$.
\end{thm}

\begin{proof}
  Let $Z$ be a component of the semistable moduli space containing points
  disjoint from $Q$; such a component exists by Lemma
  \ref{lem:Painleve_stable_exists}.  Moreover, the subscheme of the
  semistable moduli space corresponding to sheaves supported on $Q$ is
  isomorphic to $Q$ by Lemmas \ref{lem:Painleve_embed_in_elliptic} and
  \ref{lem:Painleve_elliptic}, so is in particular connected, and thus $Z$
  meets every component of the semistable moduli space; it follows by
  smoothness that $Z$ is the entire semistable moduli space.

  By Lemma \ref{lem:Painleve_derived_equivalence}, we have a derived equivalence
  \[
  D^b_{\coh} Z\cong D^b_{\coh}(X_{w^{-d},w^r})\cong D^b_{\coh}(X_{w,1}),
  \]
  where the second derived equivalence comes from Corollary
  \ref{cor:weird_langlands}.  In other words, the smooth (commutative)
  projective surfaces $Z$ and $X_{w,1}$ are derived equivalent.  Moreover,
  we can choose the second equivalence in such a way that the composition
  preserves rank and Euler characteristic.  If $w$ is not torsion, this
  already induces the desired isomorphism, by
  \cite[Thm.~3.2]{KawamataY:2002}.  Otherwise, $X_{w,1}$ is a genus 1
  fibration (with one fiber of the form $\ord(w)Q$), and the derived
  equivalence respects this fibration.  It follows that it takes smooth
  points of integral fibers to shifts of smooth points of integral fibers,
  and thus the support of the corresponding Fourier-Mukai kernel is
  $2$-dimensional, so that the proof of \cite{KawamataY:2002} still
  applies.
\end{proof}

The primary significance of this results lies in the fact that $X_{w,1}$ is
a surface of the type considered in \cite{SakaiH:2001}; in other words,
$X_{w,1}$ is the space of initial conditions of a (possibly discrete)
Painlev\'e equation; see the discussion in Section \ref{sec:Painleve_Lax}
below.

We can generalize Lemma \ref{lem:Painleve_stable_exists} somewhat.

\begin{prop}\label{prop:irreds_on_fibers_exist}
  Let $r,d$ be relatively prime integers with $r>0$, and let $z$, $q$ be
  such that $z^r q^d$ is torsion, of exact order $l$.  Then the moduli
  space of stable 1-dimensional sheaves on $X_{z,q}$ with Chern class $lrQ$
  and Euler characteristic $ld$ is nonempty.  Moreover, if $z$ has exact
  order $lr$ in $\Pic^0(Q)/\langle q\rangle$, then there exist irreducible
  such sheaves.
\end{prop}

\begin{proof}
  The proof of Lemma \ref{lem:Painleve_stable_exists} shows that there are
  semistable sheaves with the given invariants and disjoint from $Q$ such
  that any subsheaf has Chern class proportional to $Q$.  For such a sheaf
  to be strictly semistable, there would need to be a subsheaf with
  invariants $c_1=mrQ$, $\chi=md$ with $1\le m<l$, but the determinant of
  the restriction of such a sheaf to $Q$ is $(z^r q^d)^m\ne \sO_Q$.  The
  existence of irreducible sheaves when $z$ has exact order $lr$ in
  $\Pic^0(Q)/\langle q\rangle$, follows similarly; now the determinant
  condition prohibits {\em any} subsheaf.
\end{proof}

\begin{rems}
  If $z^r q^d$ has order strictly dividing $l$, then the generic semistable
  sheaf on $X^+$ is a direct sum, so is not simple, and thus semicontinuity
  shows that this is inherited by all semistable sheaves, including on the
  special fiber.  So the condition on $\ord(z^r q^d)$ is necessary and
  sufficient for the existence of stable sheaves.  If $\ord(z^r q^d)=l$ but
  $z$ has order strictly dividing $lr$ in $\Pic^0(Q)/\langle q\rangle$,
  then $q$ must be torsion, and considering the direct image on the center
  tells us that the sheaf cannot be irreducible.
\end{rems}

\begin{rems}
  There should be an extension of Theorem \ref{thm:Painleve_moduli_spaces}
  to this case, in the following form: if the ample divisor class $D_a$ is
  chosen to eliminate strictly semistable sheaves disjoint from $Q$, then
  (taking $ad-br=1$) there should be a derived equivalence $X_{z^a q^b,z^r
    q^d}\cong X_{z,q}$ taking semistable $0$-dimensional sheaves of degree
  $l$ to $D_a$-semistable sheaves with $c_1=lrQ$, $\chi=ld$.  The main
  difficulty, as in Theorem \ref{thm:Picr}, is that we cannot avoid
  strictly semistable sheaves supported on $Q$, so cannot simply use the
  semistable moduli space directly.  Note that since the proof tells us that
  the standard derived equivalence works generically, it also tells us that
  there is an obstruction to the existence of a universal family,
  given by the central simple algebra of which $X_{z^a q^b,z^rq^d}$ is an
  order.
\end{rems}

\section{Existence of irreducible sheaves}
\label{sec:irreds_exist}

One thing which is not quite settled by the above discussion is when the
moduli space of simple (or stable) sheaves (with fixed class in
$K_0^{\num}(X)$, say) is nonempty.  This is most interesting in the case of
$1$-dimensional sheaves, as these are the ones that correspond to
differential or difference equations, and in particular in the case of
$1$-dimensional sheaves disjoint from $Q$.  For such sheaves, we could also
ask a harder question, namely whether there are any {\em irreducible}
sheaves, i.e., such that any nonzero subsheaf has the same Chern class.
Such a sheaf corresponds to an equation which is irreducible in a quite
strong sense: there is no gauge transformation over $k(C)$ that makes the
equation become block triangular, or its singularities become simpler.
(The singularities becoming simpler corresponds to having a sub- or
quotient sheaf of Chern class $f-e_i-e_j$ or $e_i-e_j$.)

For simple sheaves in general, the above smoothness results give the
following reduction.

\begin{prop}
  Let $X/R$ be a split noncommutative rationally ruled surface over the dvr
  $R$ with field of fractions $K$, and let $[M]\in K_0^{\num}(X)$ be a
  class such that $\rank([M])\ne 0$ or $c_1([M])\cdot Q\ne 0$.  If there is
  a separable extension $l/k$ and a simple sheaf $M_l$ of this class on
  $X_l$, then there is a flat family of sheaves on the base change of $X$
  to an \'etale extension of $R$ such that the special fiber is isomorphic
  to $M_l$ and the generic fiber is simple.
\end{prop}

\begin{proof}
  Since any separable extension $l/k$ lifts to an \'etale extension of $R$,
  we may assume $l=k$.  The conditions on $[M]$ ensure that the moduli
  space of simple sheaves is unobstructed, and thus the point corresponding
  to $M_k$ extends to an \'etale neighborhood, and is in turn \'etale
  locally represented by a sheaf.
\end{proof}

This has a variant in the disjoint-from-$Q$ case.

\begin{prop}
  Let $X/R$ be a split noncommutative rationally ruled surface over the dvr
  $R$ with field of fractions $K$, and let $[M]\in K_0^{\num}(X/K)$ be a
  1-dimensional class such that $\det([M]|^{\bf L}_Q)\in \Pic^0(C_K)$.  If
  there is a separable extension $l/k$ and a simple sheaf $M_l$ of this
  class on $X_l$ with $M_l|_Q=0$, then there is a flat family of sheaves on
  the base change of $X$ to an \'etale extension of $R$ such that the
  special fiber is isomorphic to $M_l$ and the generic fiber is simple and
  disjoint from $Q$.
\end{prop}

\begin{proof}
  Although the moduli space of sheaves is obstructed, the determinant
  condition ensures that the obstruction is ``orthogonal'' to $R$, and thus
  we still have the desired \'etale local extensions.
\end{proof}

\begin{prop}
  Let $X/R$ be a noncommutative rationally ruled surface over the dvr $R$,
  and let $M$ be an $R$-flat coherent sheaf such that $M_k$ is an
  irreducible $1$-dimensional sheaf.  Then so is $M_K$.
\end{prop}

\begin{proof}
  Suppose $M_K$ is not irreducible, and let $N_K$ be an irreducible pure
  $1$-dimensional quotient of $M_K$.  Then $N_K$ extends to $R$ (albeit not
  necessarily as a quotient of $M$), and we may further arrange for that
  extension $N$ to have pure $1$-dimensional fibers.  (Indeed, by
  properness, we may arrange for $N_k$ to be semistable!)  Then
  semicontinuity tells us that $\dim \Hom(M_k,N_k)\ge \dim\Hom(M_K,N_K)>0$,
  so that there is a nonzero homomorphism from $M_k$ to $N_k$.  But the
  kernel and image of this homomorphism have nonzero Chern class,
  contradicting irreducibility of $M_k$.
\end{proof}

In particular, if we start with a geometrically irreducible (thus simple)
sheaf on $k$ and use the previous Propositions to lift it to a family of
simple sheaves on an \'etale cover of $R$, the generic fiber of that family
will also be geometrically irreducible.

In general, if we are given a surface $X/K$ and a class in $K_0^{\num}(X)$
and want to know whether simple sheaves (or simple sheaves disjoint from
$Q$) of the given class exist (over $\bar{K}$, but then by smoothness over
$K^{\text{sep}}$), these results suggest the following strategy.  (This is
inspired by the strategy used in \cite[\S8]{ArtinM/TateJ/VandenBerghM:1990}
to prove noncommutative planes are Noetherian.)  First, since the moduli
stack of surfaces is locally of finite type over $\Z$, the point
corresponding to $X$ factors through a field of finite transcendence degree
over its prime field, and thus $X$ itself is the base change of a surface
defined over such a field.  We thus reduce to the case that $K$ itself has
finite transcendence degree.  Any irreducible Cartier divisor on the
closure of $X$ in the moduli stack induces a valuation on $K$ and an
extension of $X$ to an \'etale extension of the corresponding valuation
ring.  By the Propositions, if the desired simple sheaves exist on the
special fiber of this family, then they exist on some separable cover of
$K$, and thus it will suffice in general to settle the question in the case
that $K$ is finite!  (It will be helpful to make one small refinement: we
may choose the smooth neighborhood in the moduli stack in such a way that
every surface in the closure of $X$ has an anticanonical curve with the
same combinatorial structure, and in the rational case remains smooth or
nodal when appropriate.)

The advantage of having $K$ be finite is that since
$\Pic^0(Q_K)/\Pic^0(C_K)$ is finite, $q$ will necessarily have finite
order, and thus any surface over a finite field is a maximal order.  This
opens up the ability to use commutative techniques in constructing the
desired sheaves.  This is particularly powerful in the case of irreducible
$1$-dimensional sheaves, since we can then consider the support of the
sheaf.

We thus consider the following problem.  For any Chern class $D\in \NS(X)$
and Euler characteristic $\chi\in \Z$, we want to know if there are
irreducible sheaves on $X/K^{\text{sep}}$ with the given invariants (and,
if $D\cdot Q=0$, disjoint from $Q$).  The usual reductions clearly apply:
if $D$ has negative intersection with a component of $Q$, then it must be
that component to be irreducible, so we may as well assume that all such
intersections are nonnegative.  We may then attempt to reduce $D$ to the
fundamental chamber.  If this fails because $D\cdot f$ becomes negative or
because we attempt to reflect in an effective simple root, this shows that
a sheaf of class $D$ cannot be irreducible, unless $D$ is actually equal to
the offending simple root.  (In that case, we know that an irreducible
sheaf exists whenever $D$ and $\chi$ satisfy the appropriate determinant
condition, and the result will be disjoint from $Q$.)

We thus reduce to the case that $D$ is in the fundamental chamber, letting
us ignore the possibility that the specializations we apply to $K$ may
introduce additional $-2$-curves.  (We have already assumed that they do
not change the structure of $Q$.)  When $X$ is given by a maximal order
${\cal A}$ of degree $l$, we may identify the N\'eron-Severi groups of $X$
and its center $Z$, and find by Proposition \ref{prop:NS_of_center} that
the support of a sheaf of Chern class $D$ on $X$ is a divisor of class $D$
on $Z$.  Conversely, given an {\em integral} curve $C_D$ on the center
which has this class and is transverse to the anticanonical curve, the
fiber of $X$ over $K(C_D)$ will be a central simple algebra.  The Brauer
group of $K^{\text{sep}}(C_D)$ is trivial (since it has transcendence
degree 1 over a separably closed field), and thus there is some separable
extension $L/K$ such that ${\cal A}\otimes L(C_D)\cong \Mat_l(L(C_D))$.  We
may WLOG assume that $L=K$ and thus that ${\cal A}\otimes K(C_D)$ has a
simple module of dimension $l$ as a vector space over $K(C_D)$.  If we view
this as an $\sO_{C_D}$-module, we may certainly extend it to a torsion-free
coherent sheaf $M'$ on $C_D$, and can then obtain a coherent ${\cal
  A}|_{C_D}$-module with the same generic fiber as the image of ${\cal
  A}\otimes M'\to M'\otimes K(C_D)$.  This ${\cal A}|_{C_D}$-module is of
course also a coherent ${\cal A}$-module $M$.  Any nonzero subsheaf of $M$
has the same generic fiber over $C_D$, so has $0$-dimensional quotient, and
thus $M$ is an irreducible sheaf of the desired Chern class.  Directly
controlling the Euler characteristic of this sheaf is somewhat tricky, but
is unnecessary.  Indeed, if the support is transverse to $Q$ but not
disjoint, then we can add any integer to the Euler characteristic by
modifying the sheaf at such a point.  If the support is disjoint from $Q$,
we can still add any multiple of $l$ to the Euler characteristic, and this
is enough to give any Euler characteristic satisfying the requisite
determinant condition.

\medskip

In the rational case, we obtain the following for irreducible sheaves
disjoint from $Q$.  Note that since changing the Euler characteristic
multiplies the determinant of restriction by a power of $q$, we have a
natural restriction map $\NS(X)\to \Pic(Q)/\langle q\rangle$.  Let $Q_l$,
$0\le l\le m$ denote the pullback to $X$ of the class of the anticanonical
curve on $X_l$.

\begin{prop}\label{prop:irreds_exist}
  Let $X/k$ be a noncommutative rational surface, and let $D\in \NS(X)$,
  $\chi\in \Z$ be such that (a) $D$ is nef and in the fundamental chamber,
  and (b) if $[M]\in K_0^{\num}(X)$ is the class of rank 0 sheaves with
  Chern class $D$ and Euler characteristic $\chi$, then $\det([M]|^{\bf
    L}_Q)\cong \sO_Q$.  Then there is an
  irreducible sheaf $M$, disjoint from $Q$, of class $[M]$ unless either
  (a) $D=rQ_8$ and $Q_8|_Q\in \Pic(Q)/\langle q\rangle$ has order $r'$
  strictly dividing $r$, or (b) $D=rQ_8+e_8-e_9$ and $Q_8|_Q\in \langle
  q\rangle$.
\end{prop}

\begin{proof}
  If $X$ is a maximal order, then the given conditions ensure that the
  corresponding linear system on the center contains integral divisors
  disjoint from the anticanonical curve (\label{thm:nonintegral}), and thus
  per the discussion above ensures the existence of suitable irreducible
  sheaves on $X$.  If $D\notin \Z Q_8\cup (e_8-e_9)+\Z Q_8$, then this is
  already enough to give the desired existence result on general surfaces
  over fields, by induction along valuations.  To make this work in the
  remaining cases, we need merely observe that we can always choose a
  valuation on a separable extension of $K$ in such a way that the order of
  $Q_8|_Q$ in $\Pic(Q)/\langle q\rangle$ does not change, and thus when we
  reach a finite field will again have the requisite irreducible curves.
\end{proof}

\begin{rem}
  Analogous results of course apply to the case that $D\cdot Q>0$, since we
  again (\cite{HarbourneB:1997}, also Propositions
  \ref{prop:non_integral2}, \ref{prop:non_integral1}) know the possible ways
  such a divisor class on the center can fail to be generically
  irreducible.
\end{rem}

The situation in the two remaining cases is more subtle.  The first case
was dealt with in Proposition \ref{prop:irreds_on_fibers_exist} (and the
remark following), completely determining in which subcases irreducible
sheaves can exist.  In the second case, irreducible sheaves in fact never
exist.  We may assume $q$ non-torsion (since otherwise the support makes
sense and is reducible) and first note that by twisting by a suitable
multiple of $e_7$, we may ensure that there are irreducible sheaves
disjoint from $Q$ with Chern class $Q_8$ and Euler characteristic 0;
twisting by $e_9$ then leaves this condition alone and allows us to set
$\chi=1$, and ensures the existence of a sheaf $\sO_{e_8-e_9}$.  By Lemma
\ref{lem:irreds_generically_acyclic}, the generic irreducible sheaf with
the given invariant is acyclic, so has a unique global section.  The image
of that global section must have the same Chern class, but then must have
the same Euler characteristic by disjointness from $Q$, and thus the sheaf
is globally generated.  Let $I$ be the kernel of the global section, and
note that $I$ is torsion-free of rank $1$, Chern class $-rQ_8-e_8+e_9$, and
Euler characteristic 0 (and in fact $R\Gamma(I)=0$).  Moreover, since
$\theta M\cong M$ is acyclic, $\Hom(I,\sO_X)\cong \Hom(\sO_X,\sO_X)$.  To
obtain a contradiction, it will suffice to show that
$\Hom(I,\sO_X(e_9-e_8))\ne 0$, since then the unique map $I\to \sO_X$
factors through $\sO_X(e_9-e_8)$, implying that $M$ has a nonzero morphism
to $\sO_{e_8-e_9}$.

It is of course equivalent to ask for $\Hom(\theta \sO_X(e_8-e_9),\ad I)\ne
0$.  It will thus (since torsion-free sheaves of rank 1 are reflexive when
$q$ is non-torsion) suffice to show that if $I$ is a torsion-free sheaf of
rank 1, Chern class $rQ_8+e_8-e_9$ and $R\Gamma(I)=0$, then
$\Hom(\sO_X(-Q_8+e_8-e_9),I)\ne 0$.  By semicontinuity, we may reduce to
the generic surface $X$ having sheaves disjoint from $Q$ with Chern class
$Q_8$ and Euler characteristic 0 (but with no sheaf of the form
$\sO_{e_8-e_9}(d)$).  Since both $\sO_X(-Q_8+e_8-e_9)$ and $I$ are in
$\sO_X^\perp$, we may transport the question to the corresponding
commutative surface $X'$, and since the generic surface has $Q$ smooth.  If
$I$ maps to a torsion-free sheaf $I'$ (which by dimension counting and
irreducibility will exhibit generic behavior), then twisting reduces to
showing that if $I'$ is the ideal sheaf of $r$ points in $X'$, then there
is a morphism from $\sO_X(-rQ_8)$ to $I'$, and this is a trivial
consequence of the fact that $X'$ is elliptic.

To see that $I$ maps to a torsion-free sheaf on $X'$, we use the following
two results.

\begin{lem}
  Let $X$ be a noncommutative rational surface with $K_X^2=0$ and smooth
  anticanonical curve.  Then for the generic torsion-free sheaf $M$ with
  numerical invariants $(1,(n-1)Q+e_8,0)$, $R\Gamma(M)=0$ and $\kappa_q M$
  is a torsion-free sheaf.
\end{lem}

\begin{proof}
  We proceed by induction on $n$, noting that when $n=0$,
  $M=\sO_X(-Q+e_8)$, and similarly for $\kappa_q M$.  For $n>0$, since both
  conditions are open, it suffices to exhibit {\em some} torsion-free sheaf
  satisfying the conclusions.  Thus let $M_{n-1}$ be such a sheaf with
  invariants $(1,(n-2)Q+e_8,0)$, let ${\cal L}$ be a generic line bundle of
  degree $0$ on $Q$, and let $M_n$ be the generic extension:
  \[
  0\to M_{n-1}\to M_n\to {\cal L}\to 0.
  \]
  By induction, $\kappa_q M_{n-1}$ is well-defined and a torsion-free
  sheaf, and $\kappa_q{\cal L}={\cal L}\otimes q^{-1}$ unless ${\cal
    L}\cong q$, so again is a sheaf, and thus so is $\kappa_q M_n$.
  Moreover, since (by dimension counting and irreducibility) $\kappa_q
  M_{n-1}$ is a generic torsion-free sheaf with the given invariants, we
  have $\Ext^1(\sO_x,\kappa_q M_{n-1})=0$ for all $x\in Q$, and thus the
  pullback of the extension to any subsheaf of $\kappa_q{\cal L}$ is still
  non-split.  It follows that $\kappa_q M_n$ is torsion-free as required.
\end{proof}

\begin{cor}
  Let $X$ be a noncommutative rational surface with $K_X^2=-1$ and smooth
  anticanonical curve, and suppose that $\det([M]|_Q)\not\cong q$, where
  $[M]$ is the class of a sheaf with numerical invariants
  $(1,nQ+e_8-e_9,0)$.  Then for the generic such sheaf, $M\in \sO_X^\perp$
  and $\kappa_q M$ is a torsion-free sheaf.
\end{cor}

\begin{proof}
  We may argue as above, except that now we take ${\cal L}=\det([M]|_Q)$
  and replace $M_{n-1}$ by its pullback to $X$.
\end{proof}

\smallskip

It should of course be possible to give similar results for higher genus
noncommutative surfaces, with the main obstruction being the requirement to
understand precisely when effective classes in the N\'eron-Severi group of
a {\em commutative} anticanonical rationally ruled surface are generically
irreducible.

\section{Irreducibility of moduli spaces}
\label{sec:moduli_are_irred}

In \cite{MarkmanE:2007}, various instances of stable moduli spaces on
Poisson (or symplectic) surfaces were shown to be irreducible (generalizing
\cite[Thm. 2.1]{EllingsrudG/StrommeSA:1993}).  With some care, the argument
can be adapted to the noncommutative case.  There are additional
technicalities in the case of $1$-dimensional sheaves with $c_1(M)\cdot
Q=0$, so we first consider the simpler case.  Note that the hypothesis that
the stable moduli space is projective is automatic for rank $1$ or (for
suitable choice of ample divisor class that excludes strictly semistable
sheaves) rank $0$.

\begin{lem}
  Let $X$ be a noncommutative rational or rationally ruled surface, let
  $D_a$ be an ample divisor class on $X$, and let $[M]\in K_0^{\num}(X)$ be
  a primitive class such that $(\rank([M]),c_1([M])\cdot Q)>(0,0)$ in
  lexicographic order.  Suppose further that the moduli space ${\cal M}$ of
  $D_a$-stable sheaves of class $[M]$ is nonempty and projective.  Then the
  class of the diagonal in the Chow ring of ${\cal M}\times {\cal M}$ is
  represented by $c_{\dim{\cal M}}(R\Hom_X(M_1,M_2)[1])$, where $M_1$,
  $M_2$ are the two pullbacks of a universal family to $X_{{\cal M}\times
    {\cal M}}$.
\end{lem}

\begin{proof}
  If $N_1$, $N_2$ are two stable sheaves of class $[M]$, then $\theta N_2$
  has strictly smaller slope than $N_1$, and thus $\Hom(N_1,\theta N_2)=0$.
  (In particular, ${\cal M}$ is unobstructed, and therefore smooth.)  It
  follows that if $M_1$, $M_2$ are two families of such sheaves (over $S$
  and $T$, say), then $R\Hom_X(M_1,M_2)$ is a perfect complex of magnitude
  $[0,1]$ on $S\times T$, and thus may be represented by a two-term complex
  of vector bundles.  Moreover, $N_1\cong N_2$ iff $\Hom(N_1,N_2)\ne 0$,
  and thus the pullback of the diagonal of ${\cal M}$ to $S\times T$ is
  precisely the locus where the complex has $\dim H^1=\dim{\cal M}$, or in
  other words (by Porteous) $c_{\dim {\cal M}}(R\Hom_X(M_1,M_2)[1])$.
  Applying this to the pullbacks of the universal family gives the desired
  result.
\end{proof}

\begin{rem}
  We stated the bulk of the argument in terms of a general pair of families
  to make it clear that the argument actually computes the class of the
  diagonal in the moduli {\em stack} of stable sheaves.  In particular, we
  see that with care the projectivity assumption is not needed for the
  result to hold.
\end{rem}

Following \cite[Cor. 10]{MarkmanE:2007}, we fix a sheaf $M_0$ of
(numerical) class $[M]$, and consider the class $c_m(R\Hom_X(M_0,M)[1])$,
where $M$ is the universal family.  On the one hand, this is the pullback
of the diagonal to the given section of $\pi_2:{\cal M}\times {\cal M}\to
{\cal M}$, and thus vanishes except on the component containing $M_0$.  On
the other hand, the base change of a numerically trivial class is still
numerically trivial, and thus the algebraic equivalence class of
$R\Hom(M_0,M)$ in $K_0({\cal M})$ is independent of $M_0$.  We thus
conclude the following (using smoothness to equate connectedness with
irreducibility).

\begin{thm}\label{thm:moduli_irred1}
  Let $X$ be a noncommutative rational or rationally ruled surface, let
  $D_a$ be an ample divisor class on $X$, and let $[M]\in K_0^{\num}(X)$ be
  a primitive class such that $(\rank([M]),c_1([M])\cdot Q)>(0,0)$ in
  lexicographic order.  Suppose further that the moduli space ${\cal M}$ of
  $D_a$-stable sheaves of class $[M]$ is nonempty and projective.  Then it
  is irreducible.
\end{thm}
  
For rank 0 sheaves with $c_1(M)\cdot Q=0$, the $\Hom$ complex can no longer
be written as a two-term complex, but luckily this problem already appears
in the case of sheaves on K3 and abelian surfaces, and the relevant lemma
from \cite{MarkmanE:2002} is directly applicable.

\begin{lem}
  Let $X$ be a noncommutative rational or rationally ruled surface, let
  $D_a$ be an ample divisor class on $X$, and let $[M]\in K_0^{\num}(X)$ be
  a primitive class such that $\rank([M])=c_1([M])\cdot Q=0$ and
  $c_1([M])\ne 0$.  Let ${\cal M}$ be the Zariski closure of the sheaves
  disjoint from $Q$ in the moduli space of $D_a$-stable sheaves of class
  $[M]$.  If ${\cal M}$ is proper and nonempty, then the class of the
  diagonal in the Chow ring of ${\cal M}\times {\cal M}$ is represented by
  $c_{\dim{\cal M}}(R\Hom_X(M_1,M_2)[1])$, where $M_1$, $M_2$ are the two
  pullbacks of a universal family to $X_{{\cal M}\times {\cal M}}$.
\end{lem}

\begin{proof}
  The only difference is that $R\Hom_X(M_1,M_2)$ is now a perfect complex
  of magnitude $[0,2]$, and thus we can only represent it as a {\em
    three}-term complex of vector bundles.  By
  \cite[Lem. 4]{MarkmanE:2002}, it suffices (since $\dim{\cal M}$ is even)
  to prove that the top cohomology sheaves of $R\Hom_X(M_1,M_2)$ and
  $R\Hom_X(M_2,M_1)^*$ are both line bundles on the diagonal.  (The other
  hypotheses are straightforward.)  Since (by semicontinuity) the universal
  family satisfies $\theta M\cong M$ fiberwise, $\Hom(M,\theta M)$ is a
  line bundle on ${\cal M}$, and thus we have
  \begin{align}
  \Hom(R\Hom_X(M_1,M_2),\sO_{{\cal M}\times {\cal M}})
  &\cong
  R\Hom_X(M_2,\theta M_1)[2]\notag\\
  &\cong
  R\Hom_X(M_2,M_1)[2]\otimes \Hom(M_1,\theta M_1),
  \end{align}
  so that the two claims are equivalent.  The claim thus reduces to showing
  that if $N_1$, $N_2$ are two stable sheaves of class $[M]$ with $N_1\cong
  \theta N_1$, $N_2\cong \theta N_2$, then $\Ext^2(N_1,N_2)\ne 0$ iff
  $N_1\cong N_2$, when $\dim\Ext^2(N_1,N_2)=1$.  By Serre duality, this is
  equivalent to the corresponding claim for $\Hom(N_1,N_2)\cong
  \Hom(N_1,\theta N_2)$.
\end{proof}

The proof of irreducibility then carries over directly.

\begin{thm}\label{thm:moduli_irred2}
  Let $X$ be a noncommutative rational or rationally ruled surface, let
  $D_a$ be an ample divisor class on $X$, and let $[M]\in K_0^{\num}(X)$ be
  a primitive class such that $\rank([M])=c_1([M])\cdot Q=0$ and
  $c_1([M])\ne 0$.  Let ${\cal M}$ be the Zariski closure of the sheaves
  disjoint from $Q$ in the moduli space of $D_a$-stable sheaves of class
  $[M]$.  If ${\cal M}$ is nonempty and proper, then it is irreducible.
\end{thm}

\section{Divisor classes and integrability}
\label{sec:integrable}

In the rational case, these results (again following
\cite{EllingsrudG/StrommeSA:1993,MarkmanE:2007}) have immediate
consequences for the Chow and Picard groups of the moduli spaces.

\begin{cor}
  Let $X$ be a noncommutative rational surface, and let $D_a$, $[M]$
  satisfy the hypotheses of Theorem \ref{thm:moduli_irred1} or
  \ref{thm:moduli_irred2}, with ${\cal M}$ nonempty and proper.  Let
  $E_1,\dots,E_{m+4}$ be a full exceptional collection on $X$.  Then the
  Chow ring $A^*({\cal M})$ is generated by the Chern classes of the
  complexes $R\Hom(E_i,M)$ where $M$ is a universal family.  In particular,
  $\Pic({\cal M})$ is generated by the classes of the line bundles $\det
  R\Hom(E_i,M)$.
\end{cor}

\begin{proof}
  Let $E^\perp_{m+4},\dots,E^{\perp}_1$ be the dual exceptional collection,
  so that $R\Hom(E_i,E^\perp_j)\cong k^{\delta_{ij}}$ (this corresponds to
  a full twist in the action of the braid group on semiorthogonal
  decompositions).  Using these, we may expand the class of the universal
  sheaf on ${\cal M}\times X$ in two ways:
  \[
    [M] = \sum_i [R\Hom_{\cal M}(R\Hom(M,E^\perp_i),E_i)]
        = \sum_i [E_i^\perp\otimes_{\cal M} R\Hom(E_i,M)].
  \]
  Using these two expansions and the duality between the exceptional
  collection, we may compute
  \[
    [R\Hom(M_1,M_2)[1]]
    =
    -\sum_i \pi_1^*[R\Hom(M_1,E^\perp_i)]\otimes \pi_2^*[R\Hom(E_i,M_2)]
  \]
  In particular, the Chern classes of $[R\Hom(M_1,M_2)[1]]$ are homogeneous
  polynomials in the Chern classes of
  \[
  \pi_1^*[R\Hom(M_1,E^\perp_i)]\otimes \pi_2^*[R\Hom(E_i,M_2)],
  \]
  which in turn are homogeneous polynomials in the Chern classes of
  the two factors.  In particular, we see that the class of the diagonal
  in ${\cal M}\times {\cal M}$ is a homogeneous polynomial of degree
  $\dim({\cal M})$ in these Chern classes.

  Now, let $z$ be any class in $A^*({\cal M})$.  Since the two projections
  induce isomorphisms with the diagonal, we have the essentially trivial
  identity
  \[
  z = \pi_{2*} (\Delta \cap \pi_1^*z).
  \]
  Expanding $\Delta$ as above, we find that $\Delta\cap \pi_1^*z$ is
  a polynomial in the Chern classes $\pi_2^*[R\Hom(E_i,M_2)]$ with
  coefficients in the image of $\pi_1^*$, and the effect of $\pi_{2*}$ on
  such an expansion is to apply it to the coefficients.  This gives a
  polynomial in the Chern classes $[R\Hom(E_i,M)]$ with integer
  coefficients, and thus proves the desired result.

  In particular, any element in $\Pic^1({\cal M})$ is a homogeneous
  polynomial of degree 1 in the Chern classes, and is thus an integer
  linear combination of the classes $c_1(R\Hom(E_i,M))$.
\end{proof}

There is a less basis-dependent way to describe the claim about the Picard
group of the moduli space.  Given any object $N\in \perf(X)$, we get a line
bundle on ${\cal M}$ by taking $\det R\Hom(N,M)$.  This map from perfect
complexes to $\Pic({\cal M})$ is clearly additive in distinguished
triangles, and thus induces a homomorphism $K_0(X)\to \Pic({\cal M})$.

\begin{cor}
  The homomorphism $K_0(X)\to \Pic({\cal M})$ is surjective.
\end{cor}

One significance of this fact is that the association of line bundles on
the moduli space to classes in $K_0(X)$ extends naturally to the moduli
space of simple sheaves (or, for that matter, simple objects in
$\perf(X)$), with the only caveat being that one must again choose a
universal family.  Note that by simplicity, any two universal families are
themselves related by a line bundle on the moduli space, and thus changing
the universal family has a relatively straightforward effect on this map:
\[
\det R\Hom(N,M\otimes_{\cal M} {\cal L})
\cong
{\cal L}^{\chi([N],[M])}
\otimes
\det R\Hom(N,M).
\]
In particular, given our primitivity assumption (which we needed for a
universal family in any case), there exists a class $[N_0]$ such that
$\chi([N_0],[M])=1$, which allows us to normalize the universal family so
that $\det R\Hom(N_0,M)$ is trivial.  We then find that
\[
\det R\Hom([N],M)
\cong
\det R\Hom([N]-\chi([N])N_0,M)
\]
and thus we may if we like restrict our attention to the sublattice of $K_0(X)$
orthogonal to $M$.

\begin{prop}
  Let $[M]$ be a class supporting sheaves disjoint from $Q$, let $[N]$ be a
  class in $K_0(X)$ and let $C$ be a proper curve in ${\cal M}$.  Then the
  function
  \[
  \det R\Hom([N],M(D))\cap C
  \]
  of $D\in c_1(M)^\perp$ grows at most quadratically in $D$.
\end{prop}

\begin{proof}
  Note that although $M(D)$ may fail to be stable, it will always be
  simple, and thus the line bundle $\det R\Hom([N],M(D))$ still makes sense
  as a divisor class on the moduli space of simple sheaves, and thus can be
  restricted to $C$.  On the moduli space of simple sheaves, we have
  $R\Hom([N],M(D))\cong R\Hom([N](-D),M)$, and thus the result follows from
  linearity in $K_0(X)$ together with the fact that $[N](-D)$ is at most
  quadratic in $D$.  (More precisely, the quadratic term is
  $\rank([N])D^2[\pt]/2$.)
\end{proof}

\begin{rem}
  From the proof, we see that the quadratic term in $\det
  R\Hom([N],M(D))\cap C$ is proportional to $\det R\Hom([\pt],M)\cap C$,
  and one in fact finds that $\det R\Hom([\pt],M)$ is the anticanonical
  class.  This follows from the fact that for $x\in Q$, $R\Hom(\sO_x,M)=0$
  unless $M$ meets $Q$ in $x$, when it is generically $1$-dimensional in
  degrees $1$ and $2$, and thus $\det R\Hom(\sO_x,M)$ is represented by the
  divisor where $M|^L_Q$ is nonzero at $x$.  But this condition is
  independent of $x$, and is nothing other than the divisor where the
  Poisson structure degenerates.  We thus conclude that if $C$ is entirely
  contained in the symplectic locus of the moduli space, then this degree
  function grows at most linearly.  It would thus be interesting to
  understand proper subspaces of the symplectic locus in general.  An
  obvious source of such subspaces is the projective space classifying
  extensions of two fixed sheaves, or more generally filtered sheaves with
  fixed subquotients; are there examples not contained in such families?
\end{rem}

This corresponds to a form of algebraic integrability (as in
\cite{BellonMP/VialletC-M:1999}).  Indeed, given a family of rational
functions on the universal moduli space (over the substack of the moduli
stack of surfaces where the moduli space is nonempty and proper), we may
(possibly after localization of the base of the family) write the global
function as a ratio of sections of a fixed family of line bundles
$R\Hom([N],M)$, and thus if we start with a family of sheaves in a fixed
surface (i.e., a map $C\to {\cal M}$ from a curve), twist by $D$, then
evaluate the function, the degree of the resulting function in $k(C)$ is at
most $\det R\Hom([N],M(D))\cap C$.  (This makes sense as long as
$M_{k(C)}(D)$ is stable, thus in particular if the $M_{k(C)}$ is
irreducible.)  It follows in particular that if we twist by $mD$ and let
$m$ tend to infinity, the degree is $O(m^2)$, and thus in a suitably
relative sense, the rational map ${-}(D)$ has algebraic entropy 0, so is
integrable. (In the $2$-dimensional case, this is essentially the argument
of \cite{TakenawaT:2001}, and one sees by a similar argument that the
relative algebra entropy of a map is the log of its largest eigenvalue on
the relative N\'eron-Severi lattice.)  This can be made non-relative for a
fixed smooth genus 1 curve, for the simple reason that in that case we can
also control the action on the N\'eron-Severi lattice of the base of the
family ($\Pic^2(C)^2\times C^m\times \Pic^0(C)$), and again have quadratic
growth.  A priori, this could yield quartic growth for the total space, but
since we have globally defined exceptional collections, we can define a
global generating set of the relative $\Pic$, and thus can write the Picard
group of the total space as a quotient of a product of $K_0(X)$ and the
Picard group of the base of the family.  Presumably a similar result
continues to hold for non-generic types, but the fact that the base of the
family is no longer proper breaks the current argument.

Something similar presumably holds in the higher genus cases; the main
difficulty is that since we do not have a full exceptional collection,
there are additional contributions to the Chow ring, and thus $\Pic({\cal
  M})$ is harder to describe.  For instance, we have a natural map from
${\cal M}$ to $\Pic(C)$ given by $\det\rho_*M$, and thus any divisor class
on the relevant component of $\Pic(C)$ induces a divisor class on ${\cal
  M}$.  For $g(C)>1$, the N\'eron-Severi lattice of $J(C)$ can vary, and
thus for some $C$, there are divisor classes on ${\cal M}$ that do not
extend to neighboring curves, and thus cannot be of the form given above.

Of course, for the above integrability argument to work, we do not actually
need {\em every} divisor class on the moduli space to have the form $\det
R\Hom(N,M)$; all we actually need is that there are {\em ample} divisor
classes of this form (since we can always write our test function as a
ratio of sections of powers of that ample bundle).  This seems much more
likely; in particular, this is true on $J(C)$ itself (since the theta
divisor on $\Pic^{g-1}(C)$ has the desired form).

\chapter{Differential and difference equations}
\label{chap:diffeq_noncomm}

\section{Sheaves as equations}
\label{sec:sheaves_from_eq_nc}

The construction of noncommutative ruled surfaces via differential and
difference (or hybrid) operators means that we can think of sheaves on such
surfaces as analogous to $D$-modules, with the caveat that the category of
$D$-modules is a deformation of a {\em quasi-}projective surface (the
cotangent bundle of the curve).  In particular, we expect that we should be
able to identify a class of sheaves corresponding to differential (or
difference) equations, or more precisely to meromorphic (possibly discrete)
connections on vector bundles.  (We set aside the hybrid and commutative
cases, as the action by operators is sufficiently far from faithful to make
it difficult to interpret sheaves in this way.  Though of course the
commutative case can be interpreted in terms of Higgs bundles or their
discrete analogues, as we saw in Chapter \ref{chap:sheaves_from_eq_comm}.)

This is simplest to achieve in the untwisted case with $\bar{Q}=Q$.  Let
$M$ be a pure $1$-dimensional sheaf on such a ruled surface.  Then the
semiorthogonal decomposition gives a long exact sequence
\[
0\to \rho_1^*\rho_{{-}1*}M
\to \rho_0^*\rho_{0*}M\to M
\to \rho_1^*R^1\rho_{{-}1*}M
\to \rho_0^*R^1\rho_{0*}M
\to 0
\]
Suppose that $M$ is $\rho_{{-}1*}$-acyclic (so $\rho_{0*}$-acyclic), and
that $W=\rho_{{-}1*}M$ and $V=\rho_{0*}M$ are vector bundles.  Then $M$ is
uniquely determined by the corresponding element of
\[
\Hom(\rho_1^*W,\rho_0^*V)
\cong
\Hom_{\bar{Q}}(\pi_1^*W,\pi_0^*V).
\]

In the differential case, $\bar{Q}$ is the double diagonal in $C\times C$,
and thus there is a natural map
\[
\Hom_{\bar{Q}}(\pi_1^*W,\pi_0^*V)
\to
\Hom_C(W,V)
\to
\Hom_{\bar{Q}}(\pi_1^*W,\pi_1^*V).
\]
If the map from $W$ to $V$ is injective (i.e., if $M$ is transverse to
$Q$), then it is generically invertible, and thus we may compose with its
inverse to get a meromorphic map from $\pi_1^*V$ to $\pi_0^*V$ such that
the corresponding meromorphic endomorphism of $V$ is the identity.  But
this is precisely a meromorphic connection on $V$!  Conversely, given a
meromorphic connection on $V$, we can choose a subsheaf $W\subset V$ with
torsion quotient such that the induced map $\pi_1^*W\to \pi_0^*V$ is
holomorphic, and thus obtain a representation as above.  This sheaf is not
unique (we can clearly replace $W$ by any subsheaf of $W$ with torsion
quotient), but the property is preserved under taking sums, and thus there
is a maximal such subsheaf, giving us a canonical representation of the
meromorphic connection as a sheaf.  (We further see as in Chapter
\ref{chap:sheaves_from_eq_comm} that $W$ is maximal iff for any sheaf of
the form $\sO_f(-1)$ (i.e., a pure 1-dimensional sheaf with Chern class $f$
and Euler characteristic 0), $\Hom(\sO_f(-1),M)=0$.)  This is dual to the
convention we used in the commutative case (being covariant rather than
contravariant), but this can be fixed by applying $R^1\ad$ and twisting by $s$.

A similar calculation in the nonsymmetric difference cases gives a discrete
analogue of a meromorphic connection, namely a meromorphic map
\[
A:\tau_q^*V\ratto V,
\]
where $\tau_q$ is the appropriate automorphism of $C$.  (When $\tau_q$ is
translation by $1\in \G_a$, this is a ``d-connection'' as in
\cite{ArinkinD/BorodinA:2006}.)  Again, there is a canonical way to take
such a discrete connection and turn it back into a sheaf on the ruled
surface.  The symmetric case is only slightly more complicated: the
discrete connection lives on $Q$ rather than $C$, and gives a meromorphic
map
\[
A:\tau_q^*\pi_0^*V\ratto \pi_0^*V
\]
such that $s_1(A)A = 1$, with the map $\pi_1^*W\to \pi_0^*V$ coming from a
factorization of $A$ as guaranteed by Hilbert's Theorem 90.

Note that although the covariant description is somewhat simpler, the
contravariant model is more natural in some respects: the point is that
because morphisms between line bundles can be represented by difference
operators, there is a representation $\Mer$ of the category of line bundles
in which the space associated to each line bundle is the space of
analytically meromorphic functions.  (In the symmetric difference case, we
must impose suitable symmetry conditions on the meromorphic functions.
There are technical difficulties in extending this approach to the
differential case, since we need to allow branch points.)  It is then
immediate that for any other module $M$ over that category, $\Hom(M,\Mer)$
is computed by solving a suitable system of difference equations associated
to $M$.  When $M$ comes from a $1$-dimensional sheaf with a presentation as
above, we can compute $\Hom(\pi_0^*V,\Mer)$ and $\Hom(\pi_1^*W,\Mer)$ by
noting that $\pi_{0*}\Mer$ and $\pi_{1*}\Mer$ are just the sheaves of
meromorphic functions, and thus, e.g., $\Hom(\pi_0^*V,\Mer)$ is just
$\rank(V)$ copies of the space of meromorphic functions; similarly, given
our explicit descriptions of $\overline S$, it is straightforward to see
that the morphisms $\pi_0^*V\to \Mer$ that descend to $M$ are nothing other
than the solutions of the equation corresponding to the dual of the above
connection.

In the symmetric difference case, since the connection lives on $Q$, we may
feel free to twist the interpretation by any automorphism of $Q$ (which
will, of course, act on all the various parameters including $q$).  This
freedom survives in a somewhat modified version in the nonsymmetric
difference case, as there is still a hidden dependence on $Q$: we had to
choose a component of $Q$ to fix the isomorphism between $C_0$ and $C_1$
and thus to determine the discrete connection.  In particular, the discrete
connection is more naturally viewed as living on that component of $Q$.  We
may thus again twist the interpretation by automorphisms of $Q$ that
preserve that component.  (Of course, if the automorphism {\em fixes} that
component, then it will not affect the interpretation.)  There are also
automorphisms that swap the two components, the most natural of which is
the deck transformation corresponding to the degree 2 map $Q\to C_0$.  In
the $q$-difference case, this transformation takes
\[
v(qz)=A(z)v(z)\qquad\text{to}\qquad v(z/q)=A(z/q)^{-1}v(z)
\]
when $V$ is trivial, and the other nonsymmetric difference cases are
analogous.  (In fact, the deck transform also has this effect in the
symmetric difference cases, this being precisely the condition to {\em be}
a symmetric difference equation.)  Similarly, in the differential case, we
most naturally have a connection on $C=Q^{\red}$, with the deck transform
acting trivially since it fixes $C$.

When the sheaf bimodule is no longer just the structure sheaf of $Q$, the
situation is more complicated.  In the invertible sheaf case, we may view
it as a twisted connection; this is not too difficult to deal with in the
difference cases, but is already quite subtle in the differential case.
And although many of the non-invertible (but torsion-free) cases are in
principle not twisted, and thus come with an explicit representation via
untwisted operators, the identification as an untwisted sheaf is only up to
nonunique isomorphism, and thus we have {\em multiple} such
representations.

Consider the differential case again.  The issue is that we want in the end
to have a meromorphic map $\pi_1^*V\ratto \pi_0^*V$ on the double diagonal,
but instead are given a meromorphic map $\pi_1^*V\ratto \pi_0^*V \otimes
{\cal E}$.  So to fix this in a consistent way, we need to choose an
invertible meromorphic map $\sO_Q\ratto {\cal E}$.  But such a map can be
specified via a suitable sheaf on the ruled surface!  In other words, to
translate meromorphic connections to sheaves on the ruled surface, we need
only specify a $\rho_{{-}1*}$-acyclic sheaf $\Sigma$ transverse to $Q$ with
$\rank(\Sigma)=0$, $\rho_{0*}\Sigma\cong \sO_C$.  This forces
$\rank(\rho_{1*}\Sigma)=1$, and thus $c_1(\Sigma)\cdot f=1$, so that we may
think of $\Sigma$ as a noncommutative analogue of a section of the ruling.
(Note that we can still get an interpretation when $\rho_{0*}\Sigma$ is a
nontrivial line bundle, but will get a connection on $V\otimes
\rho_{0*}\Sigma^{-1}$.)  This is, of course, analogous to the normalization
we used in the commutative case, with the only difference being that it is
less natural to take the structure sheaf of the section as normalization.

This of course works more generally; given a noncommutative ruled surface
$X$ and a ``section'', we obtain a correspondence between suitable
$1$-dimensional sheaves and meromorphic (discrete) connections.  Moreover,
changing the section simply tensors the connection with the inverse of the
connection corresponding to the new section relative to the old section, so
corresponds to a scalar gauge transformation.  Note that the conditions on
$M$ are the same as in the commutative case: to correspond to a (discrete)
connection (with maximal $W$), $M$ must be transverse to $Q$ and satisfy
$\Hom(\sO_f(-1),M)=\Hom(M,\sO_f(-1))=0$ for all sheaves $\sO_f(-1)$.  (The
first vanishing is needed for maximality of $W$, and implies that $V$ and
$W$ are torsion-free, while the second vanishing corresponds to
$\rho_{{-}1*}$-acyclicity of $M$.)  We of course must insist that $\Sigma$
satisfy these conditions as well.

From this perspective, the reason why the untwisted invertible case has a
canonical interpretation as connections is that it has a canonical choice
of $\Sigma$.  Even this is not unique, however: e.g., in the differential
case, we can add any holomorphic differential to the canonical flat
connection on $\sO_C$.  (In this case, the automorphism group of the ruled
surface acts transitively on such sections.)  Similarly, twisting by line
bundles ${\cal L}_i$ with $\pi_0^*{\cal L}_0\cong \pi_1^*{\cal L}_1$ has
the effect of tensoring with a holomorphic connection on ${\cal L}_0$; if
${\cal L}_0$ is degree 0, then such a connection exists, but there is no
canonical choice unless ${\cal L}_0$ is trivial.  Beyond that, we can of
course add any {\em meromorphic} differential to the canonical flat
connection, which now corresponds to a section that is no longer disjoint
from $Q$.

Once we include the section in the data, we see that the interpretation as
connections is preserved by elementary transformations.  That is, if we
start with $X_0$ with section $\Sigma$, then blowing up a point gives a
surface $X_1$ and the minimal lift $\alpha_1^{*!}\Sigma$.  This has no map
to or from $\sO_{e_1}(-1)$, and thus when we blow down $f-e_1$ to get
$X'_0$, the image of $\alpha_1^{*!}\Sigma$ will still be a section.  It
moreover agrees with the original section on the generic fiber over $C$, so
is obtained from the original meromorphic map $\sO_{\bar{Q}}\ratto {\cal
  E}$ by using the natural identification of the generic fibers of ${\cal
  E}$ and ${\cal E}'$.  This lets us give a more geometric interpretation
of how the section comes into play: we simply perform elementary
transformations until $\Sigma$ is disjoint from $Q$, which forces the
resulting surface to be the standard untwisted surface, with $\Sigma$
determined up to a holomorphic differential (in the differential case) or a
scalar factor (in the nonsymmetric elliptic difference case), as
appropriate.

Since the interpretation as a meromorphic connection essentially depends
only on the generic fiber over $C$, we find the following: a map $\phi:M\to
N$ of pure 1-dimensional sheaves such that the kernel and cokernel have
Chern class proportional to $f$ induces a gauge transformation between the
corresponding connections.  Indeed, this map is an isomorphism on the
generic fiber, and thus induces an isomorphism between the generic fibers
of $\rho_{0*}M$ and $\rho_{0*}N$ compatible with the two connections.  But
this is the same as specifying a meromorphic map $\rho_{0*}M\ratto
\rho_{0*}N$ that gauges the connection on $\rho_{0*}M$ to the connection on
$\rho_{0*}N$.  Such a gauge transformation is essentially a discrete
isomonodromy transformation, with the caveats that it may introduce or
remove apparent singularities, and that ``isomonodromy'' (though
traditional) is not really the right word when the equation has irregular
singularities.  Rather, the key property of the transformation is that any
local meromorphic solution of the equation corresponding to $\rho_{0*}M$
induces a local meromorphic solution of the equation corresponding to
$\rho_{0*}N$, and vice versa.  In the differential case, this not only
implies that their monodromy representations agree, but that they have the
same Stokes data at irregular singularities.

An important source of such transformations comes from pseudo-twists (i.e.,
taking the minimal lift, twisting by $\pm e_i$ and taking the direct
image).  Indeed, pseudo-twists have no effect on the generic fiber, so the
effect on the connection is a gauge transformation (corresponding to a
modification of $V$ at the point lying under $e_i$).  One consequence is
that, since as we have seen Lemma \ref{lem:resolving_pseudotwist} extends
to the noncommutative case, we can always arrange by a sequence of such
canonical gauge transformations to get a sheaf and a blowup such that the
minimal lift is disjoint from $Q$.  This process seems to correspond to one
of removing apparent singularities, and at the very least has the effect
that any further twists by exceptional line bundles preserve disjointness
from $Q$.  (Note that if we specify the section $\Sigma$ as the image of a
sheaf on the blowup, then twisting $\Sigma$ by exceptional line bundles
will not preserve the condition that its image on $C$ be the structure
sheaf, so will introduce an overall twist by a line bundle as mentioned
above.)  It also follows from Proposition \ref{prop:pseudo_twists_coalesce}
that {\em any} gauge transformation can be expressed as a composition of
pseudo-twists and their inverses.

This mostly describes the action of twisting by line bundles; the only
remaining detail for exceptional bundles is the description of the specific
subspace of $V|_{\pi(e_i)}$ along which one should modify $V$, but this is
relatively straightforward in most cases: for indecomposable singularities
involving $e_i$ (see below), it is straightforward, and in general the
space is given by the sum of spaces corresponding to indecomposable
singularities.

Twisting by a line bundle of Chern class $f$ is also easy to understand,
once we remember to keep track of the action on $\Sigma$: we thus find that
twisting the surface by any line bundle in $\rho_0^*\Pic(C)$ acts trivially
on the corresponding (discrete) connection.

Before considering twisting by $s$, we first consider the action of
$R^1\ad$.  This commutes with elementary transforms, so we may reduce to
the untwisted case.  We should note, of course, that the adjoint does not
preserve the canonical section $\Sigma$, so that it will be somewhat more
convenient to include a twist by a line bundle, and thus consider $M\mapsto
R^1\ad M\otimes \omega_C^{-1}$.  It moreover suffices to consider how this
behaves on a suitable localization of $C$, since we can use the known
effect on $\Sigma$ to guide the gluing on $M$.  In particular, we may as
well consider only how $M$ restricts to the {\em generic} point of $C$, so
that in particular all of the vector bundles are trivial and $M$
corresponds to a module over a $\Z$-algebra of the form $\overline{S}$.  We
know how the adjoint acts on such modules, and thus find that this acts via
the composition of the (na\"{i}ve) dual on the corresponding equation and
the deck transform over $C$.  To be precise, if $M$ corresponds to the
$q$-difference equation
\[
v(qz)=A(z)v(z),
\]
then $R^1\ad M\otimes \omega_C^{-1}$ corresponds to
\[
v(z/q)=A(z/q)^t v(z)
\]
and similarly for the other types of difference equations.  Similarly, in
the differential case, the operation takes
\[
v'(z) = A(z)v(z) \qquad\text{to}\qquad v'(z)=-A(z)^t v(z).
\]
In either case, composing by the deck transform of $Q\to C_0$ recovers the
standard duality on equations, e.g., taking $v(qz)=A(z)v(z)$ to
$v(qz)=A(z)^{-t}v(z)$.  Here we should note that although this duality
comes from the natural duality on $\GL_n$-torsors, most other functors on
$\GL_n$-torsors do not have nice translations into the sheaf framework, for
the simple reason that they tend to greatly increase the complexity of the
singularities (not to mention acting nonlinearly on the order of the
equation).

To understand the twist by $s$, we consider another contravariant symmetry,
namely taking the transpose of the morphism $B:\pi_1^*W\to \pi_0^*V$ to
get a morphism
\[
B^t:
\pi_0^*\sHom(V,\sO_C)
\to
\pi_1^*\sHom(W,\sO_C)
\]
corresponding to a sheaf on the noncommutative ruled surface arising by
swapping the roles of $C_0$ and $C_1$ in the sheaf bimodule.  We thus
obtain a (discrete) meromorphic connection on $W^*$ from the original
(discrete) meromorphic connection on $V$, and this operation is clearly a
contravariant involution.  This is closely connected to the adjoint, since
if $M$ has presentation given by $B$, then $R^1\ad M$ has a presentation of
the form
\[
0\to \rho_2^*\Hom(V,\omega_C)\to \rho_1^*\Hom(W,\omega_C)\to R^1\ad M\to 0
\]
in which the map is essentially given by $B^t$.  We thus see that (in the
untwisted case) transposing $B$ takes $M$ to $R^1\ad M\otimes L$, where $L$
is a suitable line bundle of Chern class $s+(1-g)f$ determined by the
condition that $\Sigma\cong R^1\ad \Sigma\otimes L$.  So to twist by $s$,
we need merely take the adjoint and then apply the above operation.  Note
that the combined operation is again expressible as an isomonodromy
transformation, the one caveat being that in the symmetric case, it changes
the symmetry.

The above calculations suffice to allow one to express the twist by any
line bundle as an isomonodromy transformation.\footnote{This, of course, is
not a surprise, given our discussion in Section \ref{sec:elldiff}.}  One
case worth singling out is twisting by $\theta \sO_X$.  This twist is not
quite given by $\theta$, for the simple reason that twisting by
$\theta\sO_X$ changes the identification of $Q$ with the anticanonical
curve.  We thus find that twisting by $\theta \sO_X$ is given by the
pullback through a translation, or in other words takes
\[
v(qz)=A(z)v(z)\qquad\text{to}\qquad
v(qz)=A(qz)v(z).
\]
(Note that since $A(qz)=A(qz) A(z) A(z)^{-1}$, this is indeed an
isomonodromy transformation.)

\section{Singularities}
\label{sec:sigularities_nc}

The discussion of Section \ref{sec:singularities} in the commutative case
suggests that for a $1$-dimensional sheaf $M$ on a ruled surface (which we
may suppose untwisted with $Q=\bar{Q}$), the singularities of the
corresponding connection are determined by the sheaf $M|_Q$.  This is
somewhat tricky to deal with in general, but the previous discussion
suggests that we should instead consider sheaves on iterated blowups that
are disjoint from $Q$, in which case the singularities should be determined
by the surface $X_m$ and the Chern class of $M$.  In the commutative case,
the precise correspondence to singularities was determined in Section
\ref{sec:singularities}, but this required fairly explicit calculations in
affine patches of blowups, making it nontrivial to carry out the analogous
calculation in the noncommutative setting.

In the case of a rational surface, this can be avoided in the following
way.  We first note that the blowups needed to separate a sheaf $M$ on
$X_0$ from $Q$ depend only on local information, and thus we may feel free
to modify $M$ at points not in the orbit of the support of $M|_Q$, or in
other words to apply an isomonodromy transformation which is invertibly
holomorphic on the orbit.  (This modifies $M$ by sheaves of Chern class
$\propto f$ with direct image outside the orbit.)  This makes the
separating blowup slightly more complicated, but we can recognize which
blowups were needed on the original surface from the images of the
corresponding points on $Q$, so do not lose information.  In particular, by
a suitable such modification, we may arrange for $V$ to be isomorphic to
$\sO_C^n$.  We then find that $R\Gamma(M(-f))=0$, and thus $M$ induces a
family $M_q$ of objects as we vary $q$ leaving the commutative surface
associated to $X(-f)$ unchanged (so in particular not changing the
separating blowups).  Moreover, changing $q$ has no effect on either $W$ or
the map $\rho_1^*W\to \rho_0^*V$, and thus $M_q$ is actually a family of
{\em sheaves} (injectivity does not depend on $q$), and corresponds to a
family of equations $\hbar v'(z) = A(z) v(z)$, $v(z+\hbar)=A(z)v(z)$, or
$v(qz)=A(z)v(z)$ with $A$ independent of $\hbar$ or $q$.  In addition, the
minimal lift of $M$ to the separating blowup also satisfies
$R\Gamma(M(-f))=0$, and thus the separating blowup is essentially
independent of $\hbar$ or $q$.  (We may need to exclude countably many
values of the noncommutative parameter where $M$ picks up apparent
singularities, but if $M_q$ has such singularities, we remove them by
suitable gauge transformations.) As a result, we may reduce to the case
$\hbar=0$ or $q=1$ as appropriate, and thus may apply the calculations of
Section \ref{sec:singularities}.  In particular, we find that each $e_i$ on
the separating blowup which is not a component of $Q$ corresponds to an
indecomposable singularity (depending only on the separating blowup), which
appears with multiplicity $c_1(M)\cdot e_i$ in the equation corresponding
to $M$.

Another way of viewing the above reduction is that if we localize to a
given orbit, then the commutator ideal of the corresponding ring cuts out
$Q$, and this continues after blowing up.  So if we base change to the
product of complete local rings and factor $M$ into irreducibles, each of
which is now a principal ideal, the reductions mod $Q$ of the factors are
themselves unchanged as we deform $q$.  The only caveat is that $q$
introduces an automorphism of $Q$, and thus we start to see an effect once
that automorphism becomes nontrivial, but this only happens once we have
separated to the extent that $M$ meets $Q$ in a smooth point.  In
particular, we see that the various Puiseux series representations of
Section \ref{sec:singularities} continue to hold, with the only caveat
being that there are mild additional subtleties with the lowest-order
nontrivial terms in the expressions.

The nonsymmetric elliptic difference case is also straightforward; although
we cannot as easily reduce to the commutative case, we need merely observe
that $M|_Q$ is the cokernel of the restriction to $Q$ of the map $W\to
V\otimes \bar{S}_{01}$, which we may think of as a pair of maps on the two
components of $Q$.  The matrix $A$ of the equation is simply the ratio of
those two maps, and thus we may read off the singularities of
$v(z+q)=A(z)v(z)$ from the zeros and poles of $A$, just as with finite
singularities in the $q$- and additive difference cases.

It remains to consider the higher genus differential cases.  Here we note
again that the blowup is determined by the local structure of the equation,
but now the relevant point of $C$ lies under a singular point of $Q$ (since
$Q$ is nonreduced!), and thus we may base change to the local ring at that
point.  Moreover, any sufficiently good approximation to the connection
(relative to the associated valuation) will give rise to the same sequence
of blowups, and thus we may pass to the {\em complete} local ring.  But
then the same reasoning lets us approximate the equation by one over
$k(z)$, and thus reduce to the rational case.

\section{Computing centers via singularities}

One application of the singularity classification is that it gives a
convenient way to compute centers of maximal orders in sufficiently large
characteristic.  That is, we can recognize a surface by looking at the
singularities of the connection (or generalized Higgs bundle) corresponding
to a sheaf with Chern class a sufficiently ample divisor on the surface, or
more precisely from those properties of the singularities that do not
depend on the specific choice of sheaf.  Thus to understand the structure
of the center, it suffices to understand how taking the direct image to the
center affects the separating blowup.  This in turn reduces to
understanding indecomposable singularities.  Finite singularities are
straightforward to control, so we focus on the singularities of equations
lying over singular points of $Q$, letting us work over the appropriate
complete local ring.

As an initial example, for the nonsymmetric $q$-difference case, an
indecomposable singularity is the restriction of scalars to $k[[z]]$ of an
equation of the form $v(qz) = z^{a/b} B(z) v(z)$ for $B(z)\in
\GL_m(k[[z^{1/b}]])$.  (Here we assume we are working over a field of
sufficiently large characteristic to avoid wild ramification.  Also, we
need to choose a $b$-th root of $q$ to extend the action of $z\mapsto qz$
to the given Puiseux series.)  If $q$ has order $n$, we may compute the
direct image of the corresponding sheaf as the restriction of scalars to
$k[[z^n]]$ of the equation
\[
v(q^n z)
=
\prod_{0\le i<n} (q^i z)^{a/b} B(q^i z) v(z)
=
(z^n)^{a/b} q^{an(n-1)/2b} \prod_{0\le i<n} B(q^i z) v(z).
\]
The associated separating blowup depends only on the leading term of this
matrix, and is thus determined by $a/b$ and the matrix
\[
q^{an(n-1)/2b} B(0)^n,
\]
while the separating blowup corresponding to the original sheaf determines
the Jordan block structure of $B(0)$.  (That $a/b$ does not change follows
from the fact that it is purely combinatorial: it corresponds to a sequence
of blowups in nodes of the anticanonical curve, which must therefore map to
nodes.  That an $n$th power is involved similarly follows from the fact
that the map from $Q$ to $Q'$ is the quotient by $\theta$ since $Q'$ is
reduced.)  The same calculation applies in the symmetric $q$-difference
case (as the classification of singularities is the same).

The additive difference case is somewhat trickier to deal with.  The direct
image can in general be computed as the restriction of scalars to
$k[[(z^p-z)^{-1}]]$ of $v(z+p) = A(z+p-1)\cdots A(z) v(z)$, but it is
nontrivial to compute a sufficiently good estimate of the latter product in
the case of an indecomposable singularity.  Such singularities correspond
to equations of the form
\[
v(z+1) = \alpha z^{a/b} (1+\sum_{1\le i<b} c_i z^{-1/b} + z^{-1}B(z)) v(z),
\]
where $c_1,\dots,c_{b-1}\in k$, $B\in \Mat_n(k[[z^{-1/b}]])$, and the
difficulty is that we need to compute the corresponding product to
precision $o(z^{-p})$.  This turns out to be doable, as long as the
characteristic is larger than $b$.  (If the characteristic is too small,
there are two issues: the computation of the requisite blowups may involve
wild ramification, and the anticanonical curve of the blowup may have
components of multiplicity $>p$, which complicates the action of taking the
center on the moduli stacks.)  We first note that we can factor this into a
scalar contribution and a matrix contribution, where the matrix is
$1+O(1/z)$.

\begin{lem}
  Let $k$ be a field of characteristic $p$, and let $b$ be a positive
  integer prime to $p$.  If $B\in \Mat_n(k[[z^{-1/b}]])$ is such that $B=1+B_0
  z^{-1}+o(z^{-1})$, then $B(z+p-1)\cdots B(z) = 1 +
  (B_0^p-B_0)z^{-p}+o(z^{-p}).$
\end{lem}

\begin{proof}
This is a polynomial equation in the coefficients of $B$ to order
$(z^{-1/b})^{pb+1}$, and thus it suffices to consider the generic case.  In
particular, we may assume that $B_0$ is diagonalizable (so WLOG diagonal)
with eigenvalues $\lambda_i$, such that $\lambda_i-\lambda_j\notin \F_p$
for $i\ne j$.  But then there is a matrix $C\in 1+o(1)$ such that
$C(z+1)B(z)C(z)^{-1}$ is diagonal, with
\[
B_{ii}(z) = 1+z^{-1} f_i((z^p-z)^{-1/b})
\]
for suitable power series $f_i$.  (Gauging by a matrix $C$ such that
$C-1$ has a single nonzero coefficient $C_{ij}-\delta_{ij}=\alpha z^{-a/b}$
adds $\alpha (\lambda_j-\lambda_i-a/b) z^{-a/b-1}$ to $B_{ij}$, and thus we
may set the $o(1/z)$ part of $B$ arbitrarily, except for those terms of
degree $-1$ modulo $p$.)

We then have
\begin{align}
\prod_{0\le j<p} B_{ii}(z+j)
&=
\prod_{0\le j<p} \frac{z+j+f_i((z^p-z)^{-1/b})}{z+j}\notag\\
&=
\frac{(z+f_i((z^p-z)^{-1/b})^p-(z+f_i((z^p-z)^{-1/b}))}
     {z^p-z}\notag\\
&=
1 + \frac{f_i((z^p-z)^{-1/b})^p-f_i((z^p-z)^{-1/b})}{z^p-z}\notag\\
&=
1+\frac{f_i(0)^p-f_i(0)}{z^p-z} + o(z^{-p})
\end{align}
as required.
\end{proof}

For the scalar factor, we assume $b<p$, and rewrite it using the polynomial
$\exp_p(z) = \sum_{0\le i<p} z^i/i!$, which is invertible under composition
and approximately a homomorphism.

\begin{lem}
  Let $k$ be a field of characteristic $p$, and let $1\le b<p$.  Then for
  $g\in z^{-1/b}k[[z^{-1/b}]]$ and $f=\exp_p(g)$, we have
  \[
  \prod_{0\le i<b} \exp_p(g(z+i)) =
  \exp_p(g(z)^p+g^{(p-1)}(z))|_{z=z-z^{1/p}} + o((z^p-z)^{-1}).
  \]
\end{lem}

\begin{proof}
  Since $\exp_p$ is a homomorphism to order $O(z^{-p/b})=o(1/z)$, and we
  have already seen that this holds when $\exp_p(g(z))=1+o(1/z)$, we see
  that both sides are homomorphisms to the desired order, and thus it
  suffices to prove the result for $\exp_p(g(z))=1-\alpha z^{-1/b}$.  In
  that case, we may take
  \[
  g(z) = -\sum_{1\le i\le b} i^{-1} \alpha^i z^{-i/b}+o(1/z)
  \]
  and thus
  \[
  g(z)^p+g^{(p-1)}(z)
  = -\Bigl(\sum_{1\le i\le b} i^{-1} \alpha^{ip} z^{-ip/b}\Bigr)
    +b^{-1} \alpha^b z^{-p} + o(z^{-p})
  \]
  so that we need to show
  \begin{align}
  \prod_{0\le i<p} (1-\alpha (z+i)^{-1/b})
  &=
  \exp_p\Bigl(-\sum_{1\le i\le b} i^{-1} \alpha^{ip} (z^p-z)^{-i/b}
  +b^{-1} \alpha^b (z^p-z)^{-1}\Bigr) + o((z^p-z)^{-1})\notag\\
  &=
  (1-\alpha^p (z^p-z)^{-1/b})(1+b^{-1} \alpha^b (z^p-z)^{-1}) + o((z^p-z)^{-1})\notag\\
  &=
  1-\alpha^p (z^p-z)^{-1/b}+b^{-1} \alpha^b (z^p-z)^{-1}+o((z^p-z)^{-1}).
  \end{align}
  Replacing $\alpha$ by $1/\lambda$ and multiplying both sides by
  $\lambda^p$ turns the left-hand side into the minimal polynomial of
  $z^{-1/b}$ over $k(((z^p-z)^{-1/b}))$, and thus to show that the given
  estimate holds, it suffices to plug in $z^{-1/b}$ and
  show that the value is sufficiently small.  Thus the claim follows from
  the fact that
  \[
  z^{-p/b} - (z^p-z)^{-1/b} + b^{-1} z^{(b-p)/b} (z^p-z)^{-1}
  =
  O(z^{2-2p-p/b}),
  \]
  which we may verify by expanding both powers of $z^p-z$ using the
  binomial series.
\end{proof}

We thus see that for an equation
\[
v(z+1) = \alpha z^{a/b} \Bigl(\exp_p\Bigl(\sum_{1\le i<b} c_i z^{-i/b} + B_0
z^{-1}\Bigr)+o(1/z)\Bigr) v(z),
\]
the direct image to the center is given by the the restriction of scalars
to $k[[(z^p-z)^{-1}]]$ of the matrix
\[
\alpha^p (z^p-z)^{a/b}
\Bigl(\exp_p\Bigl(\sum_{1\le i<b} c_i^p (z^p-z)^{-i/b} + (B_0^p-B_0) (z^p-z)^{-1}\Bigr)
+o((z^p-z)^{-1})\Bigr).
\]
Note that the same calculations also let us compute the requisite product
for a symmetric equation.

The differential case is similar; if we think of the equation as $Dv = Av$,
then the direct image is the restriction of scalars of the action of $D^p$.
An indecomposable equation over $k((z))$ is the restriction of scalars of
an equation
\[
v' = (f(z) + B(z))v(z),
\]
where $B\in z^{-1}\Mat_n(k[[z^{1/b}]])$ and $f\in k((z^{1/b}))$.  Since
$D\mapsto D-f$ is an automorphism, we may compute the action of $D^p$ by
first computing the action with $f(z)=0$ and adding $f(z)^p+f^{(p-1)}(z)$.
The contribution of $B$ may be computed in the generic case, when $B$ is
again diagonalizable by a suitable gauge transformation, and we find that
the contribution is $B(z)^p+B^{(p-1)}(z) + o(z^{-p})$, giving the action of
$D^p$ to sufficiently high precision.  Here we see quite explicitly that
the action of passing to the center on the parameters is more complicated
when $b>p$; any term of $f$ with exponent $<-1$ but congruent to $-1$
modulo $p$ contributes a term with exponent $<-p$ to $f^{(p-1)}$, and thus
the corresponding coefficient does not only appear as a $p$-th power.

\section{Continuous isomonodromy}
\label{sec:cont_isomonodromy}

In the discussion above, we saw that twisting by line bundles corresponds
to discrete isomonodromy transformations.  In the differential (and
ordinary difference) cases, these do not in general exhaust the full set of
monodromy- (or Stokes-, or whatever the analogue for ordinary difference
equations might be) preserving transformations: there are also typically
continuous flows between moduli spaces coming from, e.g., the fact that the
relevant fundamental group depends only on the topology of the given
punctured Riemann surface.  Although the traditional explanations of these
flows are largely analytic in nature (they are mediated by the
Riemann-Hilbert correspondence), we can hope that the {\em infinitesimal}
flow may have a more algebraic/geometric interpretation in our context.

The typical form of such a flow (we consider the differential case, but the
difference case is analogous) is (locally) given by a Lax pair, i.e., a
pair of equations
\[
v_z = Av,\qquad v_w = Bv
\]
satisfying the consistency condition ($v_{zw}=v_{wz}$)
\[
A_w = B_z + [B,A].
\]
(Here $A$ and $B$ are analytic in $w$ but algebraic in $z$.)  To first
order, this corresponds to an infinitesimal deformation of $A$, or
equivalently to the self-extension of the equation given by
\[
\begin{pmatrix}
  A & A_w\\
  0 & A
\end{pmatrix}
=
\begin{pmatrix}
  A & B_z + [B,A]\\
  0 & A
\end{pmatrix}
\]
But this is in turn the gauge transformation of the trivial self-extension
by a block unipotent matrix corresponding to $B$.  Since this is meant to
come from a flow that preserves the complexity of the singularities, we see
that the corresponding sheaf is the direct image on $X_0$ of an extension
of $M$ to a first-order deformation of the separating blowup.  Moreover, we
want this deformation to preserve the combinatorics of $Q$ (both because we
want to preserve the complexity of the singularities and because we wish to
avoid obstructed deformations), and thus in particular any singular point
of $Q^{\red}$ that must be blown up will continue to be blown up after the
deformation.  Thus we may just as well blow up such singular points before
considering any further deformations, replacing $X_0$ by some $X_{m_1}$
that still remains invariant under the deformation, but now satisfying the
condition that $M$ meets $Q^{\red}$ only in the smooth locus.

Now, the restriction to $Q$ of the corresponding sheaf on $X_{m_1}$ must
correspond to a flat deformation of $M|_Q$.  This is important because the
translation between sheaves and equations is only modulo sheaves of Chern
class $\propto f$, but any such sheaf will contribute to the restriction to
$Q$.  We thus find that an isomonodromy deformation of the kind we want
corresponds (locally on $C$) to a class in $\Ext^1(M,M)$.  Moreover, the
block upper triangular meromorphic gauge transformation that trivializes it
corresponds to an explicit trivialization of the induced extension of $M$
by $M^+$ where $M^+$ is a corresponding extension of a sheaf of Chern class
$\propto f$ by $M$.  Such an explicit trivialization corresponds to a
representation of the self-extension of $M$ as the image of a class in
$\Hom(M,M^+/M)$.  (Note that this is now actually {\em globally} defined on
$C$, since it corresponds to a class in the bottom-left corner of the
relevant hypercohomology spectral sequence.)  Assuming that $M$ has no
quotient of Chern class $\propto f$, the image of such a morphism will be
$0$-dimensional, and thus in order for the extension to exist, we must have
$M^+\subset M(D)$ for some divisor $D$ which is a nonnegative linear
combination of components of $Q$.

In other words, isomonodromy transformations that deform the separating
blowup without deforming $X_0$ correspond to morphisms $M\to M(D)/M$.
Moreover, since we are more precisely interested in isomonodromy
transformations that are consistent across an entire symplectic leaf, this
morphism should depend only on $M|_Q$.  Here there is a subtle point to
bear in mind: although every point in the symplectic leaf has isomorphic
restriction $M|_Q$, the universal family need not admit such an isomorphism
globally.  Thus our assignment of a morphism $M\to M(D)/M$ given $M|_Q$
must respect automorphisms of $M|_Q$.  In particular, we may apply this to
the induced morphism $M\to M|_Q(D)$, or equivalently $M|_Q\to M|_Q(D)$,
which must in particular commute with the local action of $\sO_Q^*$.  Since
we are working over a field of characteristic 0, ${-}(D)$ acts nontrivially
on any component of $Q$ appearing in $D$, so that the left and right
actions of $\sO_Q^*$ only agree modulo $D\cap (Q-Q^{\red})$.  (Indeed,
since $M$ only meets $Q^{\red}$ in smooth points, the localization to that
point (in the sense of \cite{vanGastelM/VandenBerghM:1997}, which suffices
to control the category of $0$-dimensional sheaves) has the form
$k\langle\langle u,v\rangle\rangle/([u,v]-v^\mu)$ where $\mu$ is the
multiplicity of the corresponding component, see Proposition
\ref{prop:local_structure_of_order_single_component}.  This lets us compute
the action of conjugation by $v^d$ for any $d$ and verify that it is
nontrivial for $d<\mu$ which is not a multiple of the characteristic.)  We
thus find that the morphism factors through $D-(D\cap Q^{\red})$, and thus
that we may reduce to the case $D=Q^{\nr}:=Q-Q^{\red}$, and thus need to
understand the {\em natural} maps $M|_Q\to M|_{Q^{\nr}}(Q^{\nr})$.

If such a map corresponds to an unobstructed deformation of the separating
blowup, then it in particular induces deformations of the support of
$M|_{Q^{\red}}$ along $Q^{\red}$.  If the deformation of a given point is
trivial (in particular if the component it lies on has multiplicity $1$ in
$Q$), then we may again blow it up, making $M|_{Q^{\red}}$ smaller.  We
thus see that we have at most one dimension of such deformations for each
time we blow up a smooth point of $Q^{\red}$ lying on a component of
multiplicity $>1$ of $Q$.  To see that this bound is tight, it suffices to
show that the map to the tangent space is surjective at every point in the
support of $M|_{Q^{\nr}\cap Q^{\red}}$.  Fixing such a point, we may work
in the category of $0$-dimensional sheaves supported at that point, or
equivalently (as we have already discussed) in the category of
finite-length modules over $R_\mu=k\langle\langle
u,v\rangle\rangle/([u,v]-v^{\mu})$, and in particular in the subcategory of
such modules annihilated by $v^{\mu-1}$.  This is the category of modules
over a {\em commutative} ring, and in that category, the functor
${-}(Q^{\nr})$ is just tensoring with the invertible bimodule generated by
$v^{1-\mu}$, so that in particular the element $v^{1-\mu}$ induces a
natural transformation $M|_{Q^{\nr}}\to M|_{Q^{\nr}}(Q^{\nr})$.  Moreover,
this element induces a derivation on $R_\mu$ taking $v$ to 0 and $u$ to
$1-\mu$, and thus in particular does not preserve the maximal ideal; it
follows immediately that the corresponding deformation moves the base point
as required.

Note that this description of the dimension of such deformations agrees
with the prediction of Section \ref{sec:singularities}, and thus in
particular the isomonodromy deformations of the various indecomposable
singularities do, in fact, glue together to form global isomonodromy
deformations.

A major simplifying assumption made above was that the surface $X_0$ itself
was not being deformed.  (This is not an issue in the ordinary difference
case, when there are no such deformations!)  If we wish to relax this
assumption, then there are some additional issues that arise.  The first is
that, as we have seen, the relation between sheaves and equations depends
not only on the surface but on a choice of section, and there is only a
canonical such choice in the untwisted case.  Thus to have an actual
isomonodromy deformation, we either need to restrict to the untwisted
case or allow the normalizing section to deform along with the surface.
Focusing on the first case (which is in any event sufficient: we can always
separate the normalizing section along with the sheaf of interest, at which
point changing the blowdown structure reduces to the case $X_0$ untwisted,
and determining which deformations act trivially on the sheaf of interest
is then straightforward), we see that the only remaining degree of freedom
is the possibility to deform $C$ itself.  Such a deformation is given by a
class in $H^1(T_C)$, or in other words by a 1-cocycle $z_{ij}$ in the space
of first-order holomorphic differential operators on $C$.  Given a
differential equation on $C$, we can then construct an extension to the
deformation in the local form
\begin{align}
\begin{pmatrix}
  1 & z_{ij} D\\
  0 & 1
\end{pmatrix}
\begin{pmatrix}
  D-A & 0\\
  0 & D-A
\end{pmatrix}
\begin{pmatrix}
  1 & -z_{ij} D\\
  0 & 1
\end{pmatrix}
&=
\begin{pmatrix}
  D-A & [D-A,z_{ij}D]\\
  0   & D-A
\end{pmatrix}\notag\\
&=
\begin{pmatrix}
  D-A & -z'_{ij}D - z_{ij}A'\\
  0   & D-A
\end{pmatrix}.
\end{align}
This is equivalent (add $z'_{ij}$ times the second row to the first) to the
equation
\[
v' = \begin{pmatrix} A & -z_{ij}A'-z'_{ij}A\\ 0 & A\end{pmatrix}v,
\]  
which is locally the isomonodromy transformation by $-z_{ij}A$.  (Note that
if $z_{ij}$ is sufficiently close to constant at a singular point, then
this is locally just the standard isomonodromy deformation moving that
singular point, while if it vanishes at the singular point, the singular
point does not move and the combinatorics of the singularity does not
change.)  To see that this is globally monodromy preserving, we note that
we obtain a local system for each open subset in the covering, and the
$z_{ij}$ induce compatible isomorphisms between the local systems on the
overlaps, and in particular ensure that at any point of $U_j\setminus
U_{ij}$ where the equation is regular, the monodromy in $U_i$ around that
point is trivial.  Thus each $U_i$ extends to a local system on the full
complement of the singularities of the equation, and these local systems
are compatibly isomorphic.  Since Stokes data is entirely local, and the
isomorphisms also ensure isomorphisms between Stokes data, the analogous
statement continues to hold for non-Fuchsian equations.

We thus see that in addition to the isomonodromy deformations that fix
$X_0$, we also have isomonodromy deformations parametrized by $H^1(T_C)$.
(When $C\cong \P^1$, we do not gain any new isomonodromy deformations, but
should recognize that three of the existing isomonodromy deformations
correspond to global automorphisms of $\P^1$; when $C$ is elliptic, we both
gain and lose an isomonodromy deformation in this way.)  In particular, it
is straightforward to determine the number of such deformations in general:
in the differential case, we start with $\chi(T_C)=3g-3$ and then as we
pass to the separating blowup, add 1 each time we blow up a smooth point of
$Q^{\red}$ lying on a multiple component.  Given the role of $Q^{\nr}$
above, we note that an easy induction shows that the expected number of
such deformations can be encapsulated in a single quantity, to wit
$-\chi(\sO_{Q^{\nr}}(Q^{\nr}))$, this time computed on the separating
blowup, with the same number holding in the ordinary difference case as
well.  (In particular, this explains the existence of a isomonodromy
deformation of the symmetric difference equations considered in
\cite{OrmerodCM/RainsEM:2017b}.)

In particular, we find that the number of isomonodromy deformations
explained by the above considerations is actually intrinsic to the surface.
This suggests that there should be a more geometrical description of these
deformations.  We have not quite been able to do this, but note the
following very suggestive facts.

First, if we consider how the continuous and discrete isomonodromy
transformations act on the given piece of the moduli stack of surfaces, we
see that they are essentially complementary: in the parametrization of
singularities, the parameters that move under continuous isomonodromy
deformations are precisely the ones that do not move under {\em discrete}
isomonodromy transformations.  Moreover, the continuously movable
parameters have a particularly nice interpretation in finite (sufficiently
large) characteristic: they are the parameters that get taken to their
$p$-th powers when passing to the center.  In particular, we find that an
infinitesimal deformation of such a parameter induces the trivial
deformation of the center!  This also applies to the transformations that
deform $C$, since the corresponding curve on the center is the image under
Frobenius.  Conversely, the parameters that move under discrete
isomonodromy transformations are acted on nontrivially by ${-}(Q)$, and
thus taking the center induces a {\em separable} map on those parameters.
We thus see that there is at least formally a correspondence between
isomonodromy deformations in characteristic 0 and center-preserving
deformations in finite characteristic.

Second, in finite characteristic, the associated deformations not only come
with extensions of $1$-dimensional sheaves disjoint from $Q$, but of {\em
  any} quasi-coherent sheaf disjoint from $Q$, for the simple reason that a
center-preserving deformation of an Azumaya algebra is trivial!  This
suggests more generally that an isomonodromy deformation should correspond
to a deformation of $X_m$ equipped with an explicit trivialization of the
induced deformation of $X_m\setminus Q$ (i.e., of the ``quasiprojective''
category obtained by inverting the natural transformation ${-}(Q)\to
\text{id}$), possibly with some additional conditions imposed to ensure the
lack of obstructions.

Even without a full geometric understanding of isomonodromy deformations,
the above considerations at least suffice to tell us that they exist when
expected, and, with limited exceptions, respect changes in blowdown
structure.  (The exceptions come from the fact that we treated deformations
of $X_0$ separately from the deformations that move base points of blowups;
the latter are treated in a way that is intrinsic to the surface!)  There
is no particular issue with elementary transformations, as those at most
have the effect of a scalar gauge transformation, so do not affect
continuous isomonodromy.  Thus the only issue is the Fourier transform.
Since this only arises when $X$ is rational, we find that the issue in that
case is not so much deformations but infinitesimal automorphisms.  There
are four cases of noncommutative $\P^1\times \P^1$s with nonreduced $Q$,
and in each case the multiple components have multiplicity 2.  By
adjunction, we find that $\sO_{Q^{\nr}}(Q^{\nr})\cong
\omega_{Q^{\nr}}^{-1}(Q^{\nr}-Q^{\red})$, and thus $\sO_{Q^{\nr}}(Q^{\nr})$
controls the deformation theory of the pair $(Q^{\nr},Q^{\nr}\cap
(Q^{\red}-Q^{\nr}))$.  In particular, there is an induced map from
infinitesimal automorphisms of the surface to global sections of
$\sO_{Q^{\nr}}(Q^{\nr})$, and in each case that map is surjective.
Moreover, in each case any automorphism of the surface acts on equations by
a change of variables, so indeed gives a trivial isomonodromy
transformation.  (Note that the automorphisms of the surface do not in
general act {\em faithfully} on equations once one takes into account the
action on the normalizing section $\Sigma$.  Moreover, the kernel of the
action on equations is not preserved by the Fourier transform!)

In particular, we find that the 2-dimensional moduli spaces arising from
Theorem \ref{thm:Painleve_moduli_spaces} admit continuous isomonodromy
transformations precisely when the curve $Q$ is nonreduced.  (A fuller
geometric understanding would presumably tell us that the consistency
conditions for those isomonodromy transformations are precisely the usual
Painlev\'e equations.  This is true in the cases involving second-order
equations, since then we can reduce to the corresponding differential
case.)

We also find that the second order difference equations admitting
continuous isomonodromy transformations controlled by Painlev\'e
transcendents have ``matrix Painlev\'e'' analogues that again admit
continuous isomonodromy transformations, which again will produce the same
nonlinear consistency equations as the corresponding differential cases,
and similarly that the equations considered in
\cite{KawakamiH/NakamuraA/SakaiH:2013} associated to 4-dimensional moduli
spaces are either of this form or have second-order difference or
differential Lax pairs.

\section{A natural group action}

It is worth noting that not every systematic family of infinitesimal
deformations preserving the complexity of singularities corresponds to an
isomonodromy transformation.  Indeed, when $g(C)>0$, we have a
$2g$-dimensional family of such deformations that do not change the
singularities at all.  The point is that on the untwisted ruled surface,
although there is a natural choice of section, the corresponding sheaf is
no longer rigid when $g>0$, and thus we obtain a $2g$-dimensional family of
regular connections on line bundles of degree 0.  Tensoring an equation by
such a connection gives a new equation with the same singularities, and
thus the tangent space at the natural section induces an infinitesimal
deformation of the equation.  By the Riemann-Hilbert correspondence, the
connection corresponding to a nontrivial section necessarily has nontrivial
monodromy, since it is not isomorphic to the connection on $\sO_C$
corresponding to the equation $v'=0$; since the tensor product of
connections has tensor product monodromy, we see that these actions do not,
in fact, preserve monodromy.  This is likely to complicate the geometric
interpretation of continuous isomonodromy transformations, especially since
on anything other than an untwisted surface, the association between
sheaves and equations is only determined modulo precisely such tensor
products!  One way to resolve this would be to consider only the {\em
  projective} monodromy (i.e., the image in $\PGL$, say by taking the
adjoint representation, of the monodromy or Galois group), at the cost of
including additional deformations (since the above $2g$-dimensional space
of deformations indeed preserve the projective monodromy).

The above action of a (symplectic) group scheme on the moduli space
interacts nicely with the natural map $\det R\rho_*$ from the moduli space
to $\Pic(C)$.  Indeed, tensoring with a connection on a line bundle simply
tensors the determinant with the $n$-th power of that bundle, assuming the
original connection is on a rank $n$ vector bundle.  We thus see in
particular that the Lagrangian subgroup $\Gamma(\omega_C)$ fixes the map to
$\Pic(C)$, which presumably plays the role of a moment map for the action.
In particular, this suggests that the quotient by $\Gamma(\omega_C)$ of any
fiber over $\Pic(C)$ should again be symplectic.

Something similar holds in the elliptic difference case, where again the
moduli space of first order ``equations'' (isomorphisms $q^*{\cal L}\cong
{\cal L}$ for ${\cal L}$ a line bundle of degree 0) is a symplectic group
scheme with a Lagrangian subgroup (now $\G_m$) respecting $\det R\rho_*$.

Of particular interest are cases in which the quotients of fibers by the
Lagrangian subgroup are 0-dimensional (so essentially rigid), or
2-dimensional (open symplectic leaves of Poisson surfaces).  The
calculations we did above for $-1$-curves and $-2$-curves can be carried
out more generally to find all possible divisors $D$ in the fundamental
chamber such that $D\cdot Q=0$ and $D^2=2g-2$ or $D^2=2g$ (since the
overall dimension is $D^2+2$, the reduction has dimension $D^2+2-2g$).  The
only possibility that exists for $g>1$ is (with respect to an even ruling)
$D=s+(g-1)f$, which is of course just the case of first-order equations
with no singularities.  For $g=1$, we obtain three additional cases, namely
the case $D=2s+f-2e_1$ with $D^2=0$ and the cases $D=2s+f-e_1-e_2$,
$D=3s+f-2e_1$ with $D^2=2$.  The first case can be ruled out because it
would correspond via an elementary transformation to a sheaf disjoint from
$Q$ on an odd ruled surface, and such sheaves cannot exist; the third case
can also be ruled out by considering the possible decompositions of $Q$ on
$X_1$.  The remaining case $D=2s+f-e_1-e_2$ was studied in the differential
case by \cite{KawaiS:2003}, where it was found that the moduli space is
rational, and the continuous isomonodromy transformations (deforming the
elliptic curve) are controlled by a special case of the usual Painlev\'e VI
equation.

Let us consider the elliptic difference analogue more closely.  We may use
elementary transformations to ensure that both singular points lie on the
same component of $Q$.  An equation then corresponds to a holomorphic map
$V\to {\cal L}_1\otimes q^*V$ where $V$ is a rank 2 vector bundle on $E$
with fixed determinant, the determinant of the map vanishes at $x_1$,
$x_2$, and ${\cal L}_1^2(-x_1-x_2)\otimes q^{\deg(V)}\cong \sO_C$.  The
action of the Lagrangian subgroup simply multiplies the map by a scalar, so
that we are really looking at the {\em projective} space of such maps,
modulo the action of $\Aut(V)$.  The simplest version takes the determinant
of $V$ to have odd degree, say 1, as then (with limited exceptions) $V$
will be the unique indecomposable bundle of that determinant.  In
particular, we find that $\dim\Hom(V,{\cal L}_1\otimes q^*V)=4$, with the
determinant giving a quadratic map
\[
\Hom(V,{\cal L}_1\otimes q^*V)
\to
\Gamma(\det(V)^{-1}\otimes {\cal L}_1^2\otimes q^*\det(V)).
\]
Since the latter line bundle has degree 2, we obtain a pencil of quadrics
inside the $\P^3$ of maps, with an open subset of the desired moduli space
being given by the complement of the base of the pencil inside the quadric
corresponding to $x_1$.  To identify the specific quadric surface and the
embedding of the curve $C$, note that any injective morphism as described
factors through a unique bundle $M_1$ such that $M_1/V\cong \sO_{x_1}$, and
similarly for $M_2$.  These bundles have degree 2, so are generically a sum
of two line bundles of degree 1, so that the set of injective morphisms
with a given $M_i$ is a $\G_m$-torsor.  (When $M_i$ is not a sum of line
bundles, it is a nontrivial self-extension of a line bundle, and we must
replace $\G_m$ by $\G_a$.)  Since each $M_i$ is classified by a quotient of
$\Pic^1(C)$ by an involution, we see that these give rulings of the
quadric.  To see how this meets the image of $C$, note that such sheaves
correspond to morphisms which are not injective, and thus to points on the
boundary of the relevant $\G_m$-torsor (or $\G_a$-torsor).  Such sheaves
are determined by giving the degree 1 line bundle through which they
factor, and each such bundle lies on a unique $M_i$.  In particular, the
map $\Pic^1(C)\to \P^1$ corresponding to $M_i$ is the quotient by the
involution ${\cal L}\mapsto \det(V)(x_i)\otimes {\cal L}^{-1}$, since
$\det(M_i)=\det(V)(x_i)$.  As long as $x_1$ and $x_2$ are distinct, these
give distinct rulings, and thus completely determine the anticanonical
quadric surface.

This is not quite the correct moduli space, for the simple reason that $V$
need not be stable in general.  If ${\cal L}_1(-x_1)$ and ${\cal
  L}_1(-x_2)$ are nontrivial, then the only unstable bundles that admit
morphisms with the correct nonzero determinant are those of the form ${\cal
  L}\oplus \det(V)\otimes {\cal L}^{-1}$ with ${\cal L}^2\cong
\det(V)\otimes {\cal L}_1$.  Modulo the action of $\Aut(V)$, each
choice of ${\cal L}$ gives rise to an $\A^1$ worth of additional points
in the moduli space, except that one sheaf on each $\A^1$ has nontrivial
stabilizer (modulo the action of the Lagrangian subgroup, that is) of order
2, and thus corresponds to a singular point of type $A_1$.

It follows that the true moduli space is obtained as follows from the data
$\det(V)$, ${\cal L}_1$, $x_1$, $x_2$ and $q$: embed $C$ in $\P^1\times
\P^1$ by the line bundles $\det(V)(x_1)$ and $\det(V)(x_2)$, then for each
point $y_i$ such that $\det(V)\otimes {\cal L}_1(-2y_i)$ is trivial, blow
up $y_i$ twice, and finally contract the four $-2$-curves that result.  In
particular, apart from the presence of singular points, this is precisely
the space on which a special case of the elliptic Painlev\'e equation acts.
Since the isomonodromy transformations induce isomorphisms of the
quasiprojective moduli space, they necessarily correspond to isomorphisms
between these surfaces.  The correspondence is somewhat more complicated
than just a direct relation between the isomonodromy equation and elliptic
Painlev\'e, for the simple reason that the $q$ parameter of the latter is
actually twice the original $q$, and thus the elliptic Painlev\'e
translation group can only give an index two subgroup of the true group of
isomonodromy transformations.  This relates to the fact that due to the
relations between the parameters, there are additional isomorphisms between
the moduli spaces; in particular, changing $\det(V)$ has no effect, and
there is even an action of $\Pic^0(C)[2]$ corresponding to twisting $V$ by
a $2$-torsion line bundle.  Modulo this residual freedom, reflecting by the
four orthogonal roots $f-e_{2i-1}-e_{2i}$ and then by $s-f$ has the desired
effect on the parameters, and thus establishes an isomorphism of moduli
spaces agreeing with the isomonodromy transformation modulo automorphisms.
It is straightforward to verify that the only automorphisms come from
$\Pic^0(C)[2]$ (an automorphism is determined by its action on the root
lattice of $\tilde{E}_8$, and must preserve the intersection form, the map
to $\Pic^0(C)$, and the set of effective roots), and thus the isomorphisms
agree for $C$ with no rational $2$-torsion, and thus in general.

\section{Painlev\'e Lax pairs}
\label{sec:Painleve_Lax}

Since the Lax pairs for Painlev\'e coming from Theorem
\ref{thm:Painleve_moduli_spaces} are of particular interest, it seems worth
spelling them out.  In the simplest case, for Painlev\'e VI, the
typical form of such a Lax pair is as a connection on a vector bundle of
rank $2r$ and degree $d$ (with $r\ge 1$ and $\gcd(r,d)=1$), with four
Fuchsian singularities and residues having two $r$-dimensional eigenspaces
with fixed eigenvalues.  The main difficulty is that although the theory
guarantees the existence of a universal sheaf in this case, and thus tells
us that there {\em is} a rational parametrization of the corresponding
family of connections, it does not give us any particular idea on how to
write down such a parametrization.

In the special case $d=1$, we can write the equation down in an alternate
form which is somewhat easier to parametrize; the cost is the introduction
of an apparent singularity, but we can use the location of that singularity
as one of the parameters.  The idea is that although the matrix form of an
equation is more natural for our purposes, we can also in general write
equations in ``straight-line'' form (i.e., as a linear relation between the
higher derivatives of a single function).  This is non-unique in general,
but in our case there is a natural choice coming from the generically
unique global section of the sheaf.  Indeed, we have noted that $M$ is
generically acyclic, and generically irreducible, and an irreducible
acyclic sheaf with the given invariants is generated by its unique global
section.  The kernel of the global section is a torsion-free sheaf of rank
1, and twisting by $c_1(M)=rQ$ makes the torsion-free sheaf a point of the
$1$-point Hilbert scheme of $X$.  A generic point of that Hilbert scheme in
turn has a unique map from the line bundle $\sO_X(-f)$, the cokernel of
which is a sheaf of the form $\sO_f(-1)$.  We thus generically obtain a
complex
\[
0\to \sO_X(-rQ-f)\to \sO_X\to M\to 0
\]
which is exact except at $\sO_X$ where the cohomology has Chern class $f$.
Thus the morphism $\sO_X(-rQ-f)\to \sO_X$ represents the same equation as
$M$, apart from introducing a single apparent singularity.  (As discussed
in Section \ref{sec:duality_of_Hilb}, this works for any $1$-dimensional
sheaf disjoint from $Q$ having a single global section, except that one
acquires $1+c_1(M)^2/2$ apparent singularities.)

In the Painlev\'e VI case, this morphism corresponds to an equation
\[
\sum_{0\le i\le 2r} c_i(z) (z(z-1)(z-\lambda))^i f^{(i)}(z) = 0
\]
such that each $c_i(z)$ is a polynomial of degree at most $4r-2i+1$.  This
equation is regular except at $0,1,\lambda,\infty$ and the unique root $v$
of $c_{2r}(z)$ (which we assume is not one of the four usual
singularities).  Each exponent at one of the true singularities induces an
arithmetic progression $e$, $e+1$,\dots,$e+r-1$ of corresponding exponents
of the straight-line equation.  The exponents at $v$ must be nonnegative
integers, and thus the global constraint on the exponents implies that they
must be $0,1,\dots,2r-2$ and $2r$; we also must have that the equation is
integrable at this point.

We start with $\sum_{0\le i\le 2r} (4r-2i+2) = (2r+1)(2r+2)$ undetermined
coefficients, and each of the eight original exponents imposes $r(r+1)/2$
linear conditions, resulting in a $2r+2$-dimensional space of equations
(that automatically satisfies the remaining condition on the exponents at
the apparent singularity; that the linear conditions are independent
follows from the fact that $\chi(\sO_X(rQ+f))=2r+2$).  Fixing the leading
term (i.e., fixing $v$ and the overall scalar) gives a $2r$-dimensional
affine space of equations, and the integrability condition at $v$ imposes
one quadratic condition for each exponent below the gap, so $2r-1$
conditions.  It follows from the general theory that the intersection of
these quadrics is a rational curve, and in fact (since the map to $v$ is
the natural ruling on the Hilbert scheme) is arithmetically $\P^1$; that
is, it admits a parametrization over the field generated by $v$ and the
exponents.  Moreover, we also find that these straight-line equations
admit a continuous isomonodromy transformation (presumably given by the
Painlev\'e VI equation), as well as a lattice of discrete isomonodromy
transformations, which are actually $r$-fold iterations of the usual
B\"acklund transformations of the space of initial conditions of PVI.  (The
full lattice can still be obtained geometrically in this setting, but now
corresponds to the translation part of the affine Weyl group action, and
thus the atomic translations are actually given by formal integral
transformations!)

There is a somewhat canonical way to parametrize the above family
of operators.  If one chooses one of the exponents at each non-apparent
singular point and gauges those exponent to be 0 (which can be done
globally in terms of twisted operators, or away from $\infty$ by
allowing that exponent to remain nonzero), then the operator can
be left-divided by $(z(z-1)(z-\lambda))^r$, and the condition to
be integrable at $v$ can be (at least generically) rephrased as
follows: there exists $u$ such that left multiplying by
$(z-u)D+(e_{\infty0}-1)$ makes the operator left-divisible by
$(z-v)$, where $e_{\infty0}$ is the chosen exponent at $\infty$.
This parametrization has the property that applying the Fourier
transform (in this case ``middle convolution'') swaps the roles of
$u$ and $v$, so in particular comes from the natural map from the
$1$-point Hilbert scheme to $\P^1\times \P^1$.  (See the discussion
in Appendix \ref{chap:fourier} on how to make the Fourier transform
work in straight-line form.  Note that gauging four exponents to be 0 makes the surface $\P^1\times \P^1$.)

\section{Morphisms}
\label{sec:morphisms_noncomm}

The relation between differential equations on an elliptic curve with a
single singularity (say at 0) and Fuchsian differential equations with four
singular points is easy enough to explain in one sense: it simply
corresponds to rewriting the equation in terms of the function $x$ on the
elliptic curve.  As long as the original equation is invariant under
$z\mapsto -z$, this change of variables can be performed, and gives a
rational equation with singular points $0$, $1$, $\infty$ and $\lambda$,
with the exponents at those singular points determined from the original
exponents at $0$.  At that point, of course, the isomonodromy
transformation deforming the elliptic curve translates directly to an
isomonodromy transformation deforming $\lambda$.

This raises the question of how the given noncommutative surfaces are
related.  More generally, given an algebraic map $\phi:C\to C'$, we should
expect some relationship between noncommutative surfaces rationally ruled
over $C$ of differential type and noncommutative surfaces rationally ruled
over $C'$ of differential type.  The most natural form such a relation
could take would be that of a morphism (i.e., an adjoint pair of functors
$\pi_*$, $\pi^*$ between the two categories).  Indeed, given an equation on
$C'$ we can certainly pull it back, while an equation on $C$ has a direct
image equation on $C'$ (of order $\deg(\phi)$ times the original order), so
that the relation we considered in the elliptic case boils down to writing
the elliptic equation as a pullback.  The analogous functors for actual
sheaves will be somewhat complicated by ramification, and indeed it seems
likely that the codomain of the morphism of noncommutative surfaces will
end up being singular in general, as it was in the elliptic difference case
considered above.  (Although we have not discussed singular noncommutative
surfaces, they are easy enough to construct: simply take the $\Z$-algebra
corresponding to a divisor on the boundary of the nef cone.  But of course
we would like to know that the result only depends on the face of the nef
cone containing that divisor, and will need to understand the categories.)

Similar sources of morphisms should arise in the difference settings; not
only should there be a direct analogue to the above coming from morphisms
between the respective curves that are equivariant with respect to $q$ or
the infinite dihedral group, one also expects morphisms taking equations
$v(qz)=A(z)v(z)$ to $v(q^k z)=A(q^{k-1}z)\cdots A(z)v(z)$, as well as
morphisms related to forgetting the symmetry of a symmetric difference
equation.  These are closely related to the notion of $G$-equivariant sheaf
discussed above in the moduli space context, with some caveats.  For a
non-symmetric difference equation, the map from
\[
v(qz)=A(q^{k-1}z)\cdots A(z)v(z)\qquad\text{to}\qquad
v(qz)=A(q^k z)\cdots A(qz)v(z)
\]
corresponds to twisting by a line bundle (assuming that $A$ has no apparent
singularities, at least!), and thus we find that translation by $q$ has the
same effect as twisting by a line bundle and permuting the blowups, so that
the sheaf is indeed isomorphic to its image under the corresponding
automorphism of order $k$ of the abelian category.  Moreover, although we
cannot recover $A$ from the new equation alone (we could, e.g., multiply
$A$ by any automorphism of order $k$ of the original equation), we {\em
  can} recover it from a specific compatible choice of isomorphism, so that
the original equation truly does represent a $\Z/k\Z$-equivariant sheaf on
the new surface.  This also applies to the case of a nonsymmetric
difference equation obtained by forgetting the symmetry of a symmetric
difference equation, where now the automorphism has order 2 and acts on $Q$
in such a way as to swap the components while preserving $q$.  The map from
$q$-difference equations to $q^k$-difference equations is more subtle in
the symmetric case, and is only the quotient by an automorphism when $k=2$;
the difficulty more generally is that the symmetry combines with the cyclic
group to form the dihedral group of order $2k$.  Of course this means that
we can still view the morphism as a map between two different quotients of
the category of nonsymmetric equations.  Similar comments apply to the case
of a $q$-difference equation equation invariant under $z\mapsto \zeta_k z$.

Also of interest in this context are the symmetries of the moduli space
coming from duality, which include cases in which the discrete connection
takes values in $\text{GO}$, $\text{GSp}$, or $\text{U}$ (or in the
corresponding Lie algebra, in the differential case).  Here it is worth
noting that second order equations {\em always} have such a symmetry, since
$\GL_2\cong \text{GSp}_2$, or more concretely since for a $2\times 2$
matrix $A$, $\det(A)A^{-t}=J A J^{-1}$, where $J$ is any nonzero
alternating matrix.  (In the differential case, this becomes
$\Tr(A)-A^t=JAJ^{-1}$.)  More precisely, given a sheaf of Chern class
$2s+df-e_1-\cdots-e_m$, the symmetry is given by composing the canonical
adjoint, the longest element of $W(D_m)$, and the twist by a suitable line
bundle, with the latter in general depending on $\det(A)$.  (Thus in the
higher genus case, this is really a symmetry of the subspace of the moduli
space on which the determinant of the equation has been fixed.)  This does
not preserve $\Sigma$, so switching back to the original $\Sigma$
introduces the requisite scalar gauge.  This description only works as
stated when none of the roots of $D_m$ are effective, but can be fixed
easily enough, at least generically.  Indeed, we can typically produce a
{\em derived} equivalence corresponding to the longest element of $W(D_m)$,
and this will take our sheaf to a sheaf unless that sheaf has a subquotient
of the form $\sO_\alpha(d)$ for $\alpha$ a root of $D_m$ (which essentially
says that the equation has apparent singularities).  The only issue (apart
from nonuniqueness of the derived equivalence) would be if we were ever
trying to reflect in an effective simple root of $D_m$ that was not
irreducibly effective.  But this cannot happen: this would imply that $Q$
had a component of the form $f-e_{i_1}-\cdots-e_{i_l}$ or
$e_{i_1}-\cdots-e_{i_l}$, but these have negative intersection with the
Chern class of our sheaf.

Although this symmetry in principle survives the Fourier transform, it does
so in a particularly obscure form, and indeed it is not clear whether the
resulting symmetry has any simpler description than as a conjugate by the
Fourier transform.  However, it can still lead to interesting consequences
when combined with other symmetries.  An interesting case is related to the
Lax pair for Painlev\'e VI of \cite{NoumiM/YamadaY:2002}, a differential
equation of the form $v' = (A_1z+A_0)v$ where $A_0,A_1\in \so_8$.  Since
$\rank(A_1)=2$, applying the Mellin transform gives a second-order
difference equation, and this equation inherits a contravariant symmetry
from the original $\so_8$ structure.  This symmetry combines with the
contravariant symmetry coming from having rank 2 to give a nontrivial
covariant symmetry, and one then finds that a suitable scalar gauge
transformation puts that symmetry in the form $A(z)=A(-z)^{-1}$.  In other
words, this second-order equation is a symmetric difference equation, which
turns out to be precisely the linear problem arising in
\cite{OrmerodCM/RainsEM:2017b}.  (One can also verify that the same thing
happens for the other two eighth-order differential equations arising via
triality.)

It may be instructive to work the above example backwards.  If we start
with the linear problem of \cite{OrmerodCM/RainsEM:2017b} and forget the
symmetry, then the result (modulo scalar gauge) is a sheaf of Chern class
$2s+4f-e_1-\cdots-e_{12}$ on a surface on which $Q$ decomposes as
\[
 (s-e_5-e_7-e_9-e_{11})
+(s-e_6-e_8-e_{10}-e_{12})
+2(f-e_1-e_2)
+(e_1-e_3)
+(e_2-e_4),
\]
with the original symmetry being reflected via the composition of swapping
$e_{2i-1}$ with $e_{2i}$ for each $i$ and a suitable involution on $Q$.
The corresponding contravariant symmetry involves the product of
reflections in the roots $f-e_5-e_6$, $f-e_7-e_8$, $f-e_9-e_{10}$,
$f-e_{11}-e_{12}$ (once we have taken into account the fact that some roots
of $D_{12}$ are components of $Q$), and thus still does not behave well
under the inverse Mellin transform.  However, if we reflect the sheaf in
$f-e_6-e_8$ and $f-e_{10}-e_{12}$, then the contravariant symmetry only
involves a permutation of the blowups.  This still does not allow a Mellin
transform--there are now too many singularities on one of the two
horizontal components of $Q$--but we can fix this by repeatedly performing
an elementary transformation in the point where that component meets the
vertical component.  (I.e., blow up that point, permute that blowup to be
the first blowup, then take the elementary transform.)  After doing that
four times, the result is a sheaf with Chern class
$2s+8f-2e_1-2e_2-2e_3-2e_4-e_5-\cdots-e_{16}$ on a surface with
anticanonical curve decomposing as
\begin{align}
Q=&\hphantom{+}(s-e_9-e_{10}-e_{11}-e_{12}-e_{13}-e_{14}-e_{15}-e_{16})
+(s-e_1-e_2-e_3-e_4)\\
&+2(f-e_1)
+2(e_1-e_2)
+2(e_2-e_3)
+2(e_3-e_4)
+2(e_4-e_5-e_6)
+(e_5-e_7)
+(e_6-e_8),\notag
\end{align}
with a contravariant symmetry involving permuting $e_9$ through $e_{16}$
and an involution on $Q$.  It follows that the inverse Mellin transform is
an 8th order equation that still has such a contravariant symmetry, and
thus after a scalar gauge to make the trace vanish preserves a bilinear form.

The above example points out two things: First, there is a fair amount of
freedom in how we obtained a sheaf having both a transform and a
contravariant symmetry surviving that transform in a recognizable form.
Indeed, although it was natural to perform the elementary transformations
in the point ``at infinity'' on the offending component, we could instead
have performed any four such transformations that respected the involution
acting on $Q$.  Most of the time this would still produce an 8th order
equation with a contravariant symmetry (albeit rather different qualitative
behavior), but including finite singularities of the difference equation in
the set also gives rise to 6th or even 4th order equations.  Second, it is
unclear how to tell a priori that the above construction gives an equation
in $\so_8$ rather than ${\mathfrak{sp}}_8$.

Geometrically, the source of this symmetry appears to be the fact that in
the commutative case, the linear system consists entirely of hyperelliptic
curves.  Another natural instance where this happens is the matrix
Painlev\'e case (on a rational surface) $D=4s+4f-2e_1-\cdots-2e_7-e_8-e_9$,
where the linear system consists of genus 2 curves.  There is indeed a
symmetry in this case as well, though the action on equations is much more
difficult to describe, as it involves acting by the longest element of
$W(E_8A_1)$ and (in the elliptic case) reflecting in the degree 2 divisor
class $x_8+x_9$.  It turns out (since $D|_Q\sim 0$) that these operations
give a covariant equivalence to the adjoint surface, and thus induce a
contravariant autoequivalence of the original surface.

This suggests that in addition to any intrinsic interest in understanding
morphisms between noncommutative surfaces, such an understanding would also
be quite fruitful in the application to special functions, e.g., by
systematically explaining quadratic (or higher-order) transformations.  In
addition, the symmetries of the moduli space associated to dualities may
give some insight into the structure of the moduli spaces of meromorphic
$G$-connections (or discrete $G$-connections) for more general structure
groups than $\GL_n$.  Note that we cannot expect there to be any reasonable
interpretation of a $G$-structure on a {\em sheaf} per se, for the simple
reason that the order of the corresponding equation depends on a choice of
ruling on the surface, so that changing the ruling may prevent the
corresponding $V$ from being a $G$-torsor.  (Indeed, this already happens
when $G=\GL_n$ and $n$ is not the minimal order of an interpretation of $M$
as an equation.)  Thus any notion of $G$-structure on a sheaf must at the
very least be taken relative to a choice of ruling, and possibly a choice
of section $\Sigma$.

\chapter{Open problems}
\label{chap:openprobs}

Although we have tried to develop the above theory as fully as possible,
there are a number of places where open problems remain.  Although we have
mostly discussed these above, it seems appropriate to list them here, with some
additional discussion.

One set of such problems relates to the classification of surfaces.  In the
commutative case, the combinatorics of the anticanonical curve gives rise
to a decomposition of the moduli stack by ``type'' of surface.  We have
seen that this fails to be a stratification in finite characteristic, but
it is open whether it is a stratification in characteristic 0.  Of course,
even if it fails to be a stratification, it still {\em induces} a
stratification, and in that case one would like to have a natural
description of the corresponding refined notion of type.  And, of course,
we would certainly like an explicit (and hopefully efficiently computable!)
description of the corresponding poset structure.  (And any relation in the
poset leads to corresponding questions about limits of equations and their
solutions.)

For noncommutative surfaces, our results show that for any reasonable
notion of ``noncommutative surface'', anything that is {\em not} of the
above form is either commutative, a deformation of an abelian or K3 surface
(which presumably does not exist) or sporadic, in that the corresponding
family does not contain any commutative fibers.\footnote{With the caveat
that our favorite pathological characteristic 2 example might have
additional noncommutative deformations, since we can currently only explain
the Poisson structure via a mixed characteristic deformation.}  Such
sporadic families exist, of course (most families of quasi-ruled surfaces),
but all examples coming from the above construction are orders in central
simple algebras.  This gives some further evidence towards a conjecture
attributed in \cite{ChanD/NymanA:2013} to Artin, namely that any
``projective noncommutative surface'' is either a sheaf of algebras on a
commutative surface or is birational (in the Van den Bergh sense) to a
ruled surface.  (Artin's actual conjecture was more precise, but also
involved a different notion of birationality.)  Of course, this may very
well depend on how one defines ``projective noncommutative surface''; in
addition to the Chan-Nyman axioms, one might reasonably assume the
existence of an ample divisor class.

One flaw in the description of noncommutative surfaces via $t$-structures
is that we still need to use the constructions in the literature to prove
several facts, including the fact that one can recover the derived category
from the heart of the $t$-structure.  Can one prove these things more
directly?  In addition to showing that the heart is large enough to recover
the dg-category, one would also like to show that the center is large in
semicommutative cases, as well as that the heart is locally Noetherian.
(For the last question, there is also the question of whether one can make
the existing argument more uniform.)

As we mentioned above, a number of natural constructions on difference(tial)
equations appear to correspond to morphisms of ruled surfaces, raising the
question of whether one can understand morphisms of noncommutative
rationally ruled surfaces in general.  One major source of morphisms comes
from quotienting by a finite group action, where there is already a problem
in the Poisson relaxation: if the finite group has a fixed point in the
open symplectic leaf, then the quotient will be singular at the image of
that point, and resolving that singularity will make the Poisson structure
on the domain have poles!  Thus it seems likely that a full description of
morphisms will either require us to deal with {\em singular} noncommutative
surfaces or with noncommutative analogues of Deligne-Mumford stacks (which
should ideally be related by some analogue of the McKay correspondence to
more usual noncommutative surfaces).

A somewhat related problem is to understand abelian and derived
equivalences of noncommutative surfaces more generally.  For abelian
equivalences, everything largely reduces to the Poisson case (an abelian
equivalence factors into a twist by a line bundle and an equivalence fixing
the structure sheaf, and the latter corresponds from an automorphism of the
corresponding anticanonical surface fixing $q$), so this is essentially a
commutative question, with the main open problem being to understand {\em
  auto}equivalences.  For derived equivalences, the natural question is: if
$X$ and $X'$ are derived equivalent but not abelian equivalent, does this
force them to be deformed elliptic surfaces?

As we mentioned above, the derived equivalences of deformed elliptic
surfaces can be viewed as generalizations of a certain instance of the
derived equivalence of geometric Langlands theory in the special case of
structure group $\GL_2$.  This more generally suggests the conjecture that
the Poisson structures on the moduli spaces are limits of noncommutative
deformations in such a way that the autoequivalences of the commutative
relaxation (coming from its structure as an Abelian fibration) extend to
the deformations, as a sort of discrete generalization of geometric
Langlands.  The one (conjectural) construction of such a noncommutative
deformation is a noncommutative deformation of the Hilbert scheme of points
\cite{elldaha} in the case of a smooth anticanonical curve, which has an
explicit description in terms of multivariate difference operators with
suitable symmetry.  (The Hilbert scheme of points is not quite an abelian
fibration, but is naturally birational to many moduli spaces of equations,
so this deformation {\em should} induce deformations of the latter.)  This
raises a couple of further questions: do the constructions of
\cite{elldaha} extend to degenerate anticanonical curves?  And can one make
the above theory work for structure groups other than $\GL_n$?

There are also a couple of holes remaining in the theory of semistable
moduli spaces.  The most glaring hole is that we have not shown that
semistable moduli spaces are projective when the rank is greater than $1$
(though luckily this only omits cases that do not correspond to
equations!).  In addition, although we can show in many cases that the
semistable moduli space is Poisson, this works by deforming the ample
bundle to give a resolution by a {\em stable} moduli space; is the
semistable moduli space still Poisson when this construction fails?
Finally, when showing algebraic integrability, we only considered a {\em
  relative} version of degree (which is arguably the right thing to look at
in the non-autonomous case!), but in the elliptic case could show {\em
  absolute} integrability.  Since the elliptic case is the one in which the
degrees of functions on the base grow the fastest, it ``ought'' to be the
hardest case to deal with, so it is natural to conjecture that algebraic
integrability holds in the absolute sense.  (Of course, we would also like
to know at least relative algebraic integrability in the case of higher
genus surfaces!  This would reduce to showing that the line bundles we know
how to construct on the moduli space includes bundles that are ample on the
stable moduli space.)

Although we have given constructions of continuous isomonodromy
deformations such that the resulting dimension is intrinsic to the surface,
the fact that we treat deformations (or automorphism!) of the base of the
ruling separately is somewhat unsatisfying, and it would be nicer to have a
fully intrinsic description.  In addition, it is not at all clear from the
present description how discrete isomonodromy transformations become
continuous isomonodromy deformations when the surface degenerates.  Another
interesting question that one might hope would be more accessible from our
algebraic perspective is to understand how the continuous isomonodromy
deformations behave in the commutative limit.  The results should be
global flows on the corresponding moduli space, and thus global sections of
the tangent bundle; how can one characterize which global sections arise in
the limit?  (And does it depend on the direction in which one takes the
limit?)

Finally, although the original question came from analysis, we have mostly
ignored analytical questions above.  For rigid equations, in addition to the
question mentioned above about limits of solutions, there is also a
question related to integral representations.  Indeed, since rigid
equations correspond to $-2$-curves, they can always be transformed into
$s-f$ by some sequence of reflections in $W(E_{m+1})$, and those
reflections are either trivial, scalar gauge transformations, or formal
integral transforms.  In the symmetric elliptic case, it was shown in
\cite{generic} that those formal integral transformations can be turned
into actual contour integrals, and thus that any rigid equation (more
precisely, any rigid sheaf with a unique global section) has a solution
given by an iterated contour integral.  Does something similar hold for
degenerations?

The other main set of analytic questions has to do with the very name
``isomonodromy''.  This is already something of a misnomer for non-Fuchsian
equations (where one must also take into account Stokes data), but as soon
as we generalize past the differential case, we find that there are
relatively few cases in which there is even a proposed replacement for
monodromy.  The main exception is the case of symmetric elliptic difference
equations, where the role of the monodromy is played by a {\em different}
symmetric elliptic difference equation (in which the step of the equation
is swapped with one of the periods).  Unfortunately, the construction is
sufficiently delicate that it is not even clear how to degenerate it to the
case of $q$-difference equations that are regular at $\infty$ (which have a
notion of monodromy introduced by Birkhoff), let alone any more degenerate
case.  One intriguing aspect of this question is the fact that although the
monodromy at the elliptic level is classified by sheaves on a
noncommutative surface, the monodromy at the differential level is given by
a local system (or a perverse sheaf of local systems in the non-Fuchsian
case), and these descriptions do not appear compatible!

\appendix
\chapter{Spaces and stacks}

The general paradigm of modern algebraic geometry---that moduli problems
should be represented by geometric objects---often encounters obstructions
to the existence of moduli schemes.  A primary example of this is moduli
spaces of sheaves; the original approach is to modify the moduli problem
slightly (by excising unstable points and weakening the equivalence
relation for semistable points), but such constructions typically break
symmetry (e.g., twisting by a line bundle does not preserve stability) and
require fairly complicated proofs.  The more modern approach involves
generalizing the notion of scheme, yielding algebraic spaces and algebraic
stacks.  Since these are not yet generally covered by basic references in
algebraic geometry, we include a brief discussion of them here.  A good
reference to working with stacks without actually {\em defining} stacks is
\cite[\S 3D]{HarrisJ/MorrisonI:1998}, while the gory details may be found in
\cite{stacksproject}.

For our purposes, it will be convenient to view algebraic spaces as a
special case of algebraic stacks.  The key observation giving rise to the
latter is that many of the pathologies of moduli problems come from
automorphisms, even in cases that {\em are} represented by schemes (e.g.,
the fact that elliptic curves with the same $j$-invariant are only {\em
  geometrically} isomorphic).  So when possible one should consider the
moduli problem as specifying a {\em groupoid} rather than an equivalence
relation; that is, given a scheme, one has a both a notion of family over
that scheme and a notion of isomorphism between families.  And, of course,
we should be able to pull families back through morphisms of schemes, and
thus pullback should be functorial for isomorphisms.  In other words, the
correct milieu for moduli problems is (contravariant) functors on schemes
(possibly over some base scheme $S$) valued in groupoids, also known as
``prestacks''.\footnote{These are normally described as ``categories
fibered in groupoids'', which is better from some technical standpoints,
but less conceptual.}

Of course, any set-valued functor gives rise to a prestack by taking the
groupoid with the given set of objects and no nonidentity morphisms.
Conversely, a ``setoid''-valued functor (i.e., in which the groupoids have
at most one morphism between any two objects) produces a set-valued functor
by taking the set of equivalence classes.  In particular, there is a
natural Yoneda embedding of the category of schemes over $S$ to the
category of prestacks over $S$.

Given two maps $F_1,F_2\to G$ of prestacks, we can define the fiber product
termwise (as the ``2-fiber product'' of groupoids).  To be precise, the
objects of $(F_1\times_G F_2)(T)$ are triples $(f_1,f_2,\phi)$ with $f_1\in
F_1(T)$, $f_2\in F_2(T)$ and $\phi$ an equivalence between their images in
$G$, while the equivalences are equivalences $f_i\cong f'_i$ making the
obvious diagram in $G$ commute, and functoriality under pullback is
straightforward.  In particular, given morphisms $F\to G$ and $T\to G$ with
$T$ (the functor corresponding to) a scheme, we may consider the base
change $F\times_G T\to T$.  A morphism $F\to G$ of prestacks is said to be
``representable (by schemes)'' if every such base change is a scheme; this
should be thought of as staying that $F$ is a family of schemes over $G$
(or a {\em relative} scheme).  Any property of morphisms of schemes that
respects base change can be defined on prestacks by asking for every base
change to a scheme to have that property; e.g., an open sub(pre)stack is a
morphism $F\to G$ such that every base change $F\times_G T\to T$ is an open
immersion.  Similarly, structures on schemes that admit pullbacks make
sense on prestacks; most notably, a {\em sheaf} on a prestack is an
assignment of a sheaf on $T$ to any object of $F(T)$, functorial with
respect to the groupoid structure on $F(T)$ and pullback through morphisms
$T'\to T$.

For set-valued functors, one typically imposes a sheaf condition.  This has
an analogue for prestacks in the form of a ``descent'' condition.  For
set-valued functors, the idea is for ``nice'' coverings $U\to X$ (e.g.,
coproducts of open subschemes in the Zariski case), the points of $F(X)$
should be in bijection with those points of $F(U)$ such that the two
pullbacks to $U\times_X U$ are the same.  Of course, for moduli problems,
two points being ``the same'' means that the corresponding families are
equivalent, which in the stack setting means that they are related by an
isomorphism.  In order to get a groupoid, we need to include a specific
isomorphism in the data, and thus we get the following notion of ``descent
datum'': (a) an object of $F(U)$ and (b) an equivalence in $F(U\times_X U)$
of the two pullbacks of the object, such that (c) the induced diagram in
$F(U\times_X U\times_X U)$ commutes.  It is straightforward to see that
descent data naturally form a groupoid, and the descent condition is that
the groupoid of descent data corresponding to a nice covering should be
equivalent to $F(X)$.

A {\em stack} is then a prestack satisfying such a descent condition
(typically w.r.to either \'etale coverings or ``fppf'' coverings).  Note
that for the stacks we consider, this will be mostly automatic: moduli
spaces of {\em sheaves} satisfy whatever descent (up to ``fpqc'') is
compatible with the conditions on the sheaves, and if $F\to G$ is
representable by schemes, then $F$ satisfies the same descent conditions
(up to fppf at least) as $G$.  (Similarly, if $T/S$ is a map of schemes and
$F$ is a stack over $T$, then the induced prestack over $S$ (with $F(T')$
the groupoid of pairs $(\phi\in \Hom_S(T',T),x\in F(\phi))$) is a stack.)
The one exception is when we consider moduli problems built from the moduli
stack of smooth genus 1 curves, as the na{\"\i}ve version of that moduli
problem (smooth proper morphisms of relative dimension 1 with reduced,
geometrically irreducible fibers of arithmetic genus 1) is only a prestack.
This can be fixed either by taking the stackification (the colimit
over \'etale coverings of the groupoid of descent data) or by considering
analogous families of algebraic spaces.

Another condition traditionally imposed on stacks is that the diagonal
should be representable.  (This is reasonable if one considers that the
diagonal of a map of schemes is generally better behaved than the original
domain: it is always separated, the diagonal of a separated morphism is a
closed immersion, and the diagonal of a closed immersion is the identity.)
This can be rephrased in several ways.  The direct rephrasing is that for
any pair of objects $x_1,x_2\in F(T)$, the ``equalizer'' should be a
scheme; more precisely, this is the functor taking $T'/T$ to the set of
isomoprhisms between the pullbacks of $x_1,x_2$ to $T'$.  An equivalent
condition is that any map $T\to F$ from a scheme should be representable,
or equivalently that any fiber product $T_1\times_F T_2$ should be a
scheme.  (This is the functor of triples $(x_1\in \Hom(T,T_1),x_2\in
\Hom(T,T_2),\phi)$ where $\phi$ is an equivalence between the pullbacks to
$F(T)$.)  In particular, this means that properties of representable maps of
prestacks make sense for maps from schemes.

This lets one state the final condition: a stack is {\em algebraic} if its
diagonal is representable and it admits a smooth covering by a scheme.
(I.e., a map $U\to F$ such that any base change to a scheme is a smooth
covering.)  There are two important subcases: a {\em Deligne-Mumford} (DM)
stack admits an {\em \'etale} covering by a scheme, while an {\em algebraic
space} is an algebraic stack in which every groupoid $F(T)$ is a setoid.
It is a nontrivial theorem that algebraic spaces are DM stacks, so in
particular also admit \'etale coverings by schemes.

The above definition is not {\em quite} general enough for our purposes;
though most of the stacks and spaces we consider not only have
representable diagonal, but {\em affine} diagonal, there is one set of
examples (involving the moduli stack of smooth genus 1 curves) where this
fails, and the diagonal is only representable by an algebraic space.  (That
is, the equalizer of two objects in $F(T)$ is the set-valued functor
associated to an algebraic space over $T$.)  We thus (as
in \cite{stacksproject}) weaken the notion of ``algebraic stack''
accordingly.  (It is still the case, luckily, that a setoid-valued
algebraic stack has scheme-representable diagonal, so there is no need to
iterate this.)  Note that in \cite{stacksproject}, one also allows
the smooth cover to be an algebraic space, but since a smooth cover of a
smooth cover is a smooth cover, this does not add any generality.

A morphism of moduli problems typically arises by imposing some additional
structure; e.g., for the moduli problem of smooth genus 0 curves, one might
impose the structure of marking and ordered triple of distinct points.
This is, in fact, an example of a smooth morphism: given a family $C\to S$
of smooth genus 0 curves, the different ways of marking an ordered triple
of distinct points are classified by an open subscheme of the smooth 3-fold
$C\times_S C\times_S C/S$.  Moreover, since this additional structure can
be defined for {\em any} family, this smooth morphism is a smooth covering
by the scheme $\Spec(\Z)$.  (Since genus 0 curves are canonically
projective (via the anticanonical bundle), this is automatically a stack,
and representability of the diagonal follows by observing that isomorphisms
are linear in the natural projective model.)  More generally, it may be
difficult to find globally defined additional structures; typically one
will have a collection of additional structures, each defined on a
different open subproblem, such that any family at least locally admits
such a structure.

Any property of a morphism $X\to Y$ of schemes that is preserved by base
change and compositions with smooth morphisms $U\to X$ can be defined for
general morphisms of algebraic stacks by asking for the property to hold
for any composition $U\to X\times_Y T\to T$ with $U\to X\times_Y T$ smooth.
For DM stacks (in particular algebraic spaces), one can replace ``smooth''
here by ``\'etale''.  For properties that are fppf-local on the base, it
suffices to consider the case that $T\to Y$ is a smooth covering.

Similarly, structures that descend from smooth (resp. \'etale) covers make
sense on algebraic stacks; most notably for our purposes, one can define
the sheaf of differentials on a DM stack or algebraic space by descent from
the sheaf of differentials on any \'etale cover.  I.e., if $F$ is a DM
stack with \'etale cover $U\to F$, and $f:T\to F$ is a morphism, the sheaf
$f^*\Omega_{F/S}$ is the sheaf which on $T\times_F U$ is given by the
pullback of $\Omega_{U/S}$, with obvious isomorphisms on $T\times_F
U\times_F U$.

In either case, it is worth noting that when $Y$ arises as the \'etale
stackification of a more natural moduli prestack $\bar{Y}$, the covering
$U\to Y$ corresponds to descent data for $\bar{Y}$ over some \'etale cover
$U'\to U$.  Since $U'$ is itself a smooth (or \'etale) cover of $Y$, we may
replace $U$ by $U'$ and thus ensure that the map $U\to Y$ comes from an
actual family.  Similarly, when working with differentials on a DM stack or
algebraic space, we need only consider maps $U\to Y$ coming from some
underlying prestack (or presheaf).  (This is under some mild assumptions on
the prestack, namely that the map from $F(X)$ to the groupoid of descent
data relative to $U\to X$ is full and faithful, but this is straightforward
to verify in the cases of interest.)

In addition, since dimension behaves well under smooth maps, we can define
dimensions on algebraic stacks.  Indeed, if $f:U\to X$ is a smooth covering
of {\em schemes}, then we can compute the dimension of $X$ at a geometric
point $x\in X(\bar{k})$ in the following way.  The base change of $f$ to
$x$ is a smooth scheme over $\bar{k}$, and thus we can choose a preimage
$\hat{x}$ of the base change.  The dimension of the base change at
$\hat{x}$ is nothing other than the relative dimension of $f$ at the image
of $\hat{x}$, and thus the dimension of $X$ at $x$ is equal to the
difference between the dimension of $U$ at the image of $\hat{x}$ and the
dimension of the base change at $\hat{x}$.  This computation continues to
make sense if we replace $X$ by an algebraic stack, and it is
straightforward to see that the dimension we get is well-defined.

An important special case of algebraic stacks comes from the action of a
smooth group scheme $G$ on a scheme $X$.  In that case, there is a
well-defined quotient prestack in which the groupoid has the same objects
as $X(T)$ and morphisms $x\to gx$ for any $g\in G(T)$.  This is not in
general a stack, but can be sheafified to give a quotient stack $[X/G]$,
which is algebraic, with smooth cover $X\to [X/G]$.  (In particular, if $X$
and $G$ are both equidimensional, then $[X/G]$ is equidimensional, with
dimension $\dim(X)-\dim(G)$.)  This is not quite the universal form of a
stack, but all of the stacks arising abover are at least {\em locally} of
this form.  Note that a map $T\to [X/G]$ is given by a $G$-equivariant map
$T'\to X$ where $T'$ is a $G$-torsor over $T$.  (I.e., $T'/T$ is
\'etale-locally isomorphic to $T\times G$, with the obvious action of $G$.)

If $G$ acts {\em freely} on $X$, then $[X/G]$ is setoid-valued, so is an
algebraic space.  Something similar holds for free actions of algebraic
group {\em stacks} on algebraic stacks.  In particular, there is a natural
action of the group stack $B\G_m:=[\Spec(\Z)/\G_m]$ classifying line
bundles on the moduli stack of coherent sheaves, which acts freely on
nonzero coherent sheaves.  The quotient stack quotients the automorphism
groups in the groupoid $\coh(X_T)$ by the normal subgroup $\G_m(T)$, and
thus the subgroupoid corresponding to {\em simple} coherent sheaves is a
setoid, giving the algebraic space classifying simple coherent sheaves.
Note here that it is difficult to give an explicit {\em \'etale} covering
of this algebraic space, while a {\em smooth} covering is relatively
straightforward (e.g., choose a basis of global sections of $M\otimes {\cal
  L}$ for sufficiently ample ${\cal L}$, where ``sufficiently ample'' is
allowed to vary over an open covering of the family of coherent sheaves);
indeed, with care, one can exhibit the moduli stack of coherent sheaves on
(projective) $X/S$ as a union of open substacks which are quotients of
schemes by smooth group schemes.

\chapter{Generalized Fourier transforms}
\label{chap:fourier}

One of the main motivations for the derived category approach we have used
above is the sheer proliferation of cases that would otherwise need to be
considered.  For instance, in the case of the Fourier transform (i.e.,
swapping the rulings of $\P^1\times \P^1$), there are 16 different cases
that naturally arise (not even including some of the issues in
characteristic 2), and one would in each case need to show that the two
representations via operators not only satisfy the same relations, but also
give rise to the same category of sheaves.  We do not attempt the latter
directly, but for applications to special functions, it is still useful to
understand the former.  The main use is that, given a new linear problem
for an integrable system of isomonodromy type, one would generally like to
know whether it is truly new, or if it can be reduced to a known linear
problem.  In particular, if we translate it into a sheaf on an appropriate
noncommutative surface, then there is a blowdown structure relative to
which its Chern class is in the fundamental domain, and thus will give a
simplest form for the equation (in particular, of lowest order).  (This may
not be unique, of course, but there will be only finitely many such forms.)
The relation to the original linear problem is via the appropriate element
of $W(E_{m+1})$, and we have already seen that $W(D_m)$ acts by scalar
gauge transformations, so that the only truly nontrivial action involves
the Fourier transform.  Thus a suitable understanding of the Fourier
transform will let us understand all of the minimal linear problems
equivalent to a given linear problem.

One tricky issue is that on the untwisted (rational) ruled surfaces for
which we have natural interpretations of equations, the divisor class $s-f$
is always effective, and thus we do not have a Fourier transform (as an
abelian equivalence, that is).  Thus there is invariably an issue with
twisting to consider.  One approach would be to choose the normalizing
section $\Sigma$ to have Chern class $s+f$; this is ample when $s-f$ is
ineffective, and thus there is always such a sheaf.  This has the
disadvantage of introducing additional parameters, and thus additional ways
for things to degenerate; we obtain a total of 48 different possibilities
for the possible structure of the anticanonical curve on the separating
blowup for such a $\Sigma$.  (To be precise, if we take into account the
order in which the points are blown up, there are 225 cases, but they fall
into 48 orbits.)  It is also somewhat cumbersome to compute the transform
in this form, as $\Sigma$ normalizes equations, but does not quite
normalize an algebra of operators: the issue is that $\Sigma$ is not
invariant under twisting by line bundles, so we can only use this
normalization to control operators acting on the trivial vector bundle.
This approach is still workable, as we can still use it to compute the
kernel of a suitable formal integral transform (or, more precisely, the
equations satisfied by the kernel), which is enough to enable the
computation of the Fourier transform on equations in straight-line form.

An alternate approach is to simply find {\em some} representation in terms
of operators in each of the 16 cases, and then check that there are
isomorphisms as required.  The main disadvantage of this approach is that
we need to represent the full category of morphisms between line bundles on
$X_0$, and there is a great deal of nonuniqueness in that representation.
In particular, if we assign to each line bundle a first-order equation,
with specified gauge equivalences between them, then we can gauge by the
resulting system to obtain a new representation.  (Note that $\Sigma$
only specifies the equation associated to the trivial bundle!)  We can also
similarly gauge by a system of automorphisms of $Q$, with the resulting
effect on matrix equations being to pull back by the automorphism
associated to the trivial bundle.  Although this nonuniqueness makes it
relatively easy to find representations, there is a significant cost when
it comes to understanding limits: to degenerate one case to another, it
may be necessary to make a suitable gauge transformation first.

If we consider applying such a representation to computing the Fourier
transform of an equation, we find that there is a considerable
simplification available.  The point is that the sheaves
$\sO_X,\sO_X(-s),\sO_X(-f),\sO_X(-s-f)$ form a strong exceptional
collection, and thus if $M$, $M(-s)$, $M(-f)$, and $M(-s-f)$ are all
acyclic, then $M$ has a resolution of the form
\[
0\to \sO_X(-s-f)^a\to \sO_X(-f)^b\oplus \sO_X(-s)^c\to \sO_X^d\to M\to 0.
\]
In other words, $M$ can be expressed as the solution of a system of $b+c$
equations in $d$ unknowns, with each equation either being a linear
equation (with coefficients in $\Hom(\sO_X(-f),\sO_X)$) or a first-order
difference/differential equation (coming from $\Hom(\sO_X(-s),\sO_X)$).  We
can recover the corresponding (discrete) connection by using the linear
equations to express the $d$ unknowns as the global sections of a vector
bundle, and observing that the remaining equations describe a (discrete)
connection on that vector bundle.  But then to understand how the Fourier
transform acts on equations, it suffices to understand how it acts on
suitable sections of $\Hom(\sO_X(-f),\sO_X)$ and $\Hom(\sO_X(-s),\sO_X)$.
Moreover, the condition for a pair of maps from those 2-dimensional spaces
to operators (with the first mapping to multiplication operators and the
second to first-order operators) to extend to a representation of the
category is basically that there be two more such maps, from
$\Hom(\sO_X(-s-f),\sO_X(-s))$ and $\Hom(\sO_X(-s-f),\sO_X(-f))$, such that
the compositions span a $4$-dimensional space.  Note that by a judicious
use of the gauge freedom, we can choose the representation so that twisting
by $s+f$ does not affect the representation: gauge by a suitable square
root of the anticanonical natural transformation.  We can furthermore
arrange in this way for all of the spaces of multiplication operators to
agree, at which point we can deduce all of the spaces of degree $s$ given
one such space and the compatibility condition.

For instance, in the (symmetric) elliptic difference case, we choose a
ramification point of $Q\to C_0$ to make it an honest elliptic curve, and
then use the remaining gauge-by-automorphisms freedom (at the cost of
choosing an element $q/2$) to ensure that all of the operators are
invariant under $z\mapsto -z$.  Then there is a representation in which the
typical operator of degree $f$ is proportional to $\vartheta(z\pm
u):=\vartheta(z+u,z-u):=\vartheta(z+u)\vartheta(z-u)$, while the typical
operator of degree $s$ is proportional to
\[
D_q(c\pm u)
:=
\frac{\vartheta(c+u+z,c-u+z)}{\vartheta(2z)} T^{1/2}
+
\frac{\vartheta(c+u-z,c-u-z)}{\vartheta(-2z)} T^{-1/2}
\]
where $(T^{\pm 1/2}f)(z)=f(z\pm q/2)$ and $c$ is a parameter depending not
only on the surface but on the domain of the morphism of line bundles.
(Here the theta function $\vartheta(z)$ is defined as in \cite[\S
  2]{elldaha}, as the canonical section of the line bundle on $E$ of
functions with a simple pole at the identity.  This can also be represented
analytically as the ratio of the usual theta function by its derivative at
the origin, see Lemma 2.1 op.~cit.)  We in particular find (by comparing
coefficients of $T^{1/2}$, say) that the spans of
\[
D_q(c\pm v)\vartheta(z\pm u)
\qquad\text{and}\qquad
\vartheta(z\pm u)D_q(c+q/2\pm v)
\]
agree and are 4-dimensional, so that these indeed extend to give a
representation of a category of the desired form.  Moreover, we see (by
comparison to \cite{generic}, or simply by noting that each $2$-dimensional
space is a space of global sections of a line bundle, and the relations are
the same as those satisfied by the global sections) that this actually
gives the general form of a relation in the elliptic case, so every
noncommutative $\P^1\times \P^1$ with smooth anticanonical curve has a
representation in this form.  The Fourier transform can be viewed as a
formal system of operators ${\cal F}_q(c)$ such that ${\cal F}_q(-c)={\cal
  F}_q(c)^{-1}$ and
\[
\vartheta(z\pm u) {\cal F}_q(c) = {\cal F}_q(c+q/2) D_q(c+q/2\pm u).
\]
(One can in fact take ${\cal F}_q(c)$ to be a certain formal difference
operator, as the univariate case of \cite[\S8]{elldaha}, though it is
perhaps more natural to view it as a formal integral operator.)
In particular, given an equation of the form
\[
\vartheta(z\pm u_1) \lambda_1\cdot v(z) +
\vartheta(z\pm u_2) \lambda_2\cdot v(z) = 0,
\]
then $w={\cal F}_q(c)v$ formally satisfies
\[
D_q(-c+q/2\pm u_1) \lambda_1\cdot w(z) +
D_q(-c+q/2\pm u_2) \lambda_2\cdot w(z) = 0,
\]
and similarly if
\[
D_q(c+q/2\pm u_1) \lambda_1\cdot v(z) +
D_q(c+q/2\pm u_2) \lambda_2\cdot v(z) = 0,
\]
then
\[
\vartheta(z\pm u_1) \lambda_1\cdot w(z) +
\vartheta(z\pm u_2) \lambda_2\cdot w(z) = 0.
\]
The compatibility condition essentially reduces to
\[
  {\cal F}_q(c) D_q(u_0,u_1,u_2,u_3)
  =
  D_q(u_0-c,u_1-c,u_2-c,u_3-c){\cal F}_q(c)
\]
for $u_0+u_1+u_2+u_3=2c+q$, where
\[
D_q(u_0,u_1,u_2,u_3)
:=
\frac{\vartheta(u_0+z,u_1+z,u_2+z,u_3+z)}{\vartheta(2z)} T^{1/2}
+
\frac{\vartheta(u_0-z,u_1-z,u_2-z,u_3-z)}{\vartheta(-2z)} T^{-1/2}
\]
is (modulo scalars) the typical element of degree $s+f$.  (In particular,
these elements generate a graded algebra which is the quotient of the
Sklyanin algebra \cite{SklyaninEK:1982,SklyaninEK:1983} by a central
element of degree 2, see also \cite{RosengrenH:2004,sklyanin_anal}.)

The top $q$-difference case (so again with $Q$ integral) is a
straightforward limit of the elliptic case: we simply replace the function
$\vartheta$ by $\exp(-\pi\sqrt{-1}z)-\exp(\pi\sqrt{-1}z)$, and replace the
various variables and parameters by suitable logarithms.  Thus the
multiplication operators become $z+1/z-u-1/u$, the difference operators
become
\[
D_q(c u^{\pm 1})
:=
\frac{cz+1/cz-u-1/u}{1/z-z} T^{1/2}
+
\frac{c/z+z/c-u-1/u}{z-1/z} T^{-1/2}
\]
and
\begin{align}
&D_q(u_0,u_1,u_2,u_3)\notag\\
&\,:=
\frac{(1-u_0z)(1-u_1z)(1-u_2z)(1-u_3z)}{\sqrt{u_0u_1u_2u_3}z(1-z^2)}T^{1/2}
+
\frac{(1-u_0/z)(1-u_1/z)(1-u_2/z)(1-u_3/z)}{\sqrt{u_0u_1u_2u_3}z^{-1}(1-z^{-2})}T^{-1/2},
\end{align}
and the Fourier transform acts by
\begin{align}
  {\cal F}_q(q^{-1/2}c) (z+1/z-u-1/u) &= D_q(q^{1/2}u^{\pm 1}/c) {\cal F}_q(c),
\\
  {\cal F}_q(q^{1/2}c) D_q(q^{1/2}cu^{\pm 1}) &= (z+1/z-u-1/u) {\cal F}_q(c),
\\
  {\cal F}_q(c) D_q(u_0,u_1,u_2,u_3)
  &=
  D_q(u_0/c,u_1/c,u_2/c,u_3/c){\cal F}_q(c)\qquad (u_0u_1u_2u_3=qc^2).
\end{align}
The top ordinary difference case is similar: just replace $\vartheta(z)=z$
and optionally set $q$ to $1$.  Note that in these cases, it is no longer
true that every operator in the appropriate space is proportional to one of
the given form (e.g., we now have a multiplication operator $1$), but this
remains true for a dense set of operators.  Also, in these cases, we may
interpret ${\cal F}_q(c)$ as a fractional power of a suitable symmetric
lowering operator (the Askey-Wilson or Wilson operator, as appropriate).
Indeed, the limit as $u\to\infty$ of the main identity of the Fourier
transform gives
\[
{\cal F}_q(c) = {\cal F}_q(q^{1/2}c) (1/z-z)^{-1} (T^{1/2}-T^{-1/2}),
\]
and thus
\[
{\cal F}_q(q^{-n/2})
=
{\cal F}_q(1)
\bigl((1/z-z)^{-1} (T^{1/2}-T^{-1/2})\bigr)^n,
\]
with ${\cal F}_q(1)$ acting trivially.  This lowering operator takes
symmetric Laurent polynomials to symmetric Laurent polynomials, decreasing
the degree, and thus in particular the equation
\[
(1/z-z)^{-1} (T^{1/2}-T^{-1/2})v=0
\]
may be viewed as the normalizing section $\Sigma$.  The additive case
similarly corresponds to powers of $(2z)^{-1}(T^{1/2}-T^{-1/2})$.
(Something similar is true in the elliptic case, see \cite[\S8]{elldaha},
although it is no longer true that the result is actually a power of the
operator for $c=q^{-1/2}$.)

The next easiest set of degenerations to consider are those for which $Q$
has two components of class $s+f$, and thus the operators are no longer
symmetric.  Here we first encounter the issue mentioned above with limits:
in order to obtain a nonsymmetric operator from a symmetric operator, we
must conjugate by an automorphism to reintroduce the symmetry as a
parameter and then take the limit in that parameter.  Thus in the
$q$-difference case, we conjugate by $z\mapsto wz$ and take a limit
$w\to\infty$.  Doing so in a na\"{\i}ve fashion makes the operators
of interest blow up, and rescaling the operators kills the parameters, so
we find that we must also rescale $u$.  There remains a further subtle
issue, in that we need to rescale the operators of degree $f$ but not of
degree $s$, which implicitly means that we need to include an additional
scale factor in the Fourier transforms.  In the end, we find that the
operators of degree $f$ become $z-u$, the operators of degree $s$ become
$(-cT^{1/2}+c^{-1}T^{-1/2})+(u/z)(T^{1/2}-T^{-1/2})$, and the Fourier
transform acts by
\[
  (z-u) {\cal F}'_q(c) = {\cal F}'_q(q^{1/2}c)
  ((u/z-q^{1/2}c)T^{1/2} +(-u/z+1/q^{1/2}c)T^{-1/2}).
\]
(We omit the action on elements of degree $s+f$.)  This again corresponds
to a fractional power of a lowering operator, namely
$z^{-1}(T^{-1/2}-T^{1/2})$.  Similarly, the additive case involves taking a
limit $w\to\infty$ after substituting $z\mapsto z+w$, $u\mapsto u+w$, and
gives (again up to rescaling ${\cal F}'$ accordingly)
\[
(z-u){\cal F}'_q(c) = {\cal F}'_q(c+q/2)
((u-z-c-q/2)T^{1/2}+(z-u-c-q/2)T^{-1/2}),
\]
corresponding to fractional powers of $T^{-1/2}-T^{1/2}$.  If we then let
$q,c\to 0$ at comparable rates, then a similar rescaling gives a transform
involving differential operators:
\[
(z-u) {\cal F}'(\alpha) = {\cal F}'(\alpha+1/2)
((z-u)D+(2\alpha+1)),
\]
again with ${\cal F}'(-\alpha)={\cal F}'(\alpha)^{-1}$.  This is
essentially the transform of \cite{KatzNM:1996} (often called ``middle
convolution'' in the later literature).  Moreover, by taking a limit of the
corresponding fractional-power-of-lowering-operator interpretation, we see
that it can be described via fractional differentiation: ${\cal F}'(\alpha)
= D^{-2\alpha}$; indeed, one has
\[
D^{2\alpha+1} (z-u) = ((z-u)D+(2\alpha+1)) D^{2\alpha}.
\]

We next turn to the transforms that map between symmetric and nonsymmetric
operators, or in geometric terms relate surfaces with $Q=(2s+f)+(f)$ to
surfaces with $Q=(s+2f)+(s)$.  Here the main complication is that by making
the symmetric operators truly symmetric, we obtain a transform without any
parameters, but still need mild dependence of the transform on the domain
of the morphism.  Luckily, this dependence is very mild, and indeed simply
involves a power of $T^{1/2}$.  The overall idea is to take the limit as
$c\to 0$ of the symmetric $q$-Fourier transform, except that we need to
include an overall shift of the parameter on one side to break the
symmetry.  It also turns out that we need to include an overall gauge on
the nonsymmetric side to make the limits work, which we do in such a way as
to ensure that the operators take polynomials to polynomials.  This gives a
pair of inverse transforms ${\cal F}^{<}_q$ and ${\cal F}^{>}_q$ such
that
\begin{align}
{\cal F}^{<}_{q}q^{-1/2}(z-u)T^{-1/2}
&=
\bigl(
\frac{z-u}{z^2-1}T^{-1/2}
-
\frac{z^{-1}-u}{z^{-2}-1}T^{1/2}
\bigr)
{\cal F}^{<}_{q}\\
{\cal F}^{<}_{q}
\bigl((z+1/z-u-1/u)T^{-1/2}-z^{-1}T^{1/2}\bigr)
T^{1/2}
&=
(z+1/z-u-1/u){\cal F}^{<}_{q}.
\end{align}
The additive equivalent is
\begin{align}
{\cal F}^{<}_{+}T^{-1/2}
(z-u)
&=
\bigl(
\frac{z-u-1/2}{2z}T^{-1/2}+
\frac{z+u+1/2}{2z}T^{1/2}
\bigr)
{\cal F}^{<}_{+}\\
{\cal F}^{<}_{+}
\bigl((z^2-u^2)T^{-1/2}-T^{1/2}\bigr)
T^{1/2}
&=
(z^2-u^2){\cal F}^{<}_{+}.
\end{align}

The next cases to consider are those for which $Q$ has components of Chern
class $s+f$, $s$, and $f$, which again comes in multiplicative and additive
flavors.  In these cases, the corresponding commutative Poisson surface is
uniquely determined, and thus we have the additional possibility of
arranging for the Fourier transform to act as an involution on that
surface.  There is, however, a technical issue that arises, particularly in
the multiplicative case: it turns out that the involution is not Poisson,
but rather anti-Poisson.  (It has non-isolated fixed points in the open
symplectic leaf, so negates the volume form at those points.)  As a result,
the simplest form of the corresponding Fourier transform turns out to
invert $q$.  Note also that because $Q$ has components of degrees $s$ and
$f$, we can use the corresponding natural transformations to determine the
gauge in place of a square root of the anticanonical natural
transformation.  This has the advantage of expressing the category of line
bundles as a Rees algebra relative to a filtration by $\N^2$.  The
corresponding filtered algebra may be taken to be generated by $z$ (of
degree $f$) and by the lowering operator $z^{-1}(1-T)$ (of degree $s$).
These generators satisfy the equivalent relations
\[
yx = qxy+(1-q)\qquad\text{and}\qquad
xy = q^{-1}yx+(1-q^{-1}),
\]
so that swapping $x$ and $y$ gives the algebra with $q$ inverted.  The
additive case is somewhat simpler, as in that case we can arrange for the
involution to be actually Poisson.  We obtain a bifiltered algebra with
relation
\[
[y,x]=y-x
\]
and involution that swaps $x$ and $y$, with a representation $x\mapsto z$,
$y\mapsto z+T$ in difference operators.

The remaining cases have similar interpretations in terms of bifiltered
algebras.  The multiplicative case (with $Q=(s)+(s)+(f)+(f)$) is
particularly simple: the filtered algebra is simply the $q$-Weyl algebra
$yx=qxy$, with representation $x\mapsto z$, $y\mapsto T$ and involution
that again inverts $q$.  The additive case ($Q=(s)+(s)+2(f)$) transforms
to a differential case ($Q=2(s)+(f)+(f)$), with filtered algebra
\[
[y,x] = y
\]
and representations $x\mapsto z$, $y\mapsto T$ and $x\mapsto -tD$,
$y\mapsto t$.  The corresponding transforms are quite familiar: from the
differential side to the discrete side is precisely the Mellin transform,
and the inverse is the $z$ transform.  (And, of course, the
differential-to-discrete direction underlies the standard power series
method for solving linear ODEs.)  The remaining case ($Q=2(s)+2(f)$) is the
usual Weyl algebra $[y,x]=1$, with representations $x\mapsto t$, $y\mapsto
D$ and $x\mapsto -D$, $y\mapsto t$.  This of course is just the classical
Fourier transform.

\medskip

As we have seen, the above calculations suffice to compute the Fourier
transform of a typical sheaf, i.e., for which $M,M(-s),M(-f),M(-s-f)$ are
all acyclic (which, of course, we can always arrange by twisting by
$\sO_X(d(s+f))$ for $d\gg 0$), as then the equation translates to a system
of equations with two forms of relations: ones of the form $m_1v=m_2w$ for
$m_1$, $m_2$ generators of the space of multiplication operators, and ones
of the form $\delta_1v=\delta_2w$ where $\delta_1,\delta_2$ span the space
of minimal first-order operators, with the Fourier transform swapping the
two.  To turn the result back into a more standard form of equation, note
that the relations involving multiplication operators define a sheaf on
$\P^1$, which is nothing other than the sheaf $W$.  One caution here is
that the map from sheaves to equations is not faithful, and the kernel is
not invariant under the Fourier transform; in particular, a
block-triangular equation with a sufficiently simple first-order block will
at best give a transform with apparent singularities, and when unlucky can
fail to give an equation at all (e.g., when the transformed sheaf has a
subsheaf of the form $\sO_f(d)$, $d\ge 0$).

This issue is particularly important when considering equations in
straight-line form.  Although in a sense the transform is particularly
straightforward here (simply express the operator as a polynomial in the
generators and transform every generator), the result will usually not be
what one wants.  The problem is that the typical source of straight-line
equations is a 1-dimensional sheaf with a global section, but the resulting
straight-line form almost always has apparent singularities, and thus the
transformed equation will be reducible!

In the simplest case, with a single apparent singularity (i.e., when
$c_1(M)^2=0$ on the separating blowup), this is easy enough to resolve:
left multiplying by operators of ``degree'' $s$ gives a $2$-dimensional
space of operators in which the spurious subsheaf has class $s+f$, which is
codimension $1$ in the space of all operators with a spurious subsheaf of
that class.  In particular, that space meets the $2$-dimensional space of
operators which are left-divisible by a multiplication operator, and
dividing that out gives an operator in which the spurious subsheaf has
class $s$.  The transform of that new operator then has spurious subsheaf
of class $f$, and is thus a straight-line form for the transform of our
original sheaf.

More generally, when $c_1(M)^2=2g-2$, one finds that the space of operators
with spurious subsheaf of class $gs+g^2f$ is $(g+1)(g^2+1)-g$-dimensional,
while the corresponding multiples of the fixed straight-line form form a
space of dimension $(g+1)(g^2-g+1)$, which has codimension $g^2$.  The
other hand, the multiples of the desired operator with a spurious subsheaf
of class $gs$ form a $g^2+1$-dimensional space, and thus again we have an
intersection, from which we may compute the straight-line form of the
transformed sheaf.

\medskip

Our purpose above was to give enough information to enable explicit
calculations of the transforms in special functions applications.  As such,
there are a number of issues we have not attempted to address.  One is that
in the elliptic case \cite{generic} and in the cases involving differential
operators (classical), the generalized Fourier transform can be expressed
as an integral operator.  Of course, in the differential cases, there are
potential issues with convergence (e.g., the usual Fourier transform is
difficult to define unless the argument decays at infinity), and in
particular it is unclear whether there is an actual operator acting on a
space containing solutions to the equations of interest.  (This was mostly
settled in the elliptic case, but even there, the actual contour integral
description requires additional constraints on the solutions, so that not
every solution appears, just a sufficiently independent set of solutions.)
The situation is worse in the discrete cases, as the solutions are only
determined up to multiplication by meromorphic functions invariant under
the shifts.  In particular, one finds that the kernel of the formal
integral transform is itself only determined up to shift-invariant
functions.  One may also be tempted to think that since the Fourier
transform takes operators of degree $s+f$ to operators of degree $s+f$,
that this means that the Fourier transform of an explicit solution of such
an equation will be easily expressible in closed form.  Although this is
true in many cases, the indeterminacy of the solution and kernel poses some
difficulties.  Indeed, all we know is that the integral is a solution of
the equation, but this only determines the integral up to a periodic
function of $z$, which need not have any nice form.

Another is that at the elliptic level, the transforms have natural
multivariate analogues related to Macdonald theory; indeed, the limit to
the most general multiplicative case is essentially a fractional power of
the Macdonald-Koornwinder lowering operator.  This suggests that all 16 of
the above transforms should have similar multivariate analogues.

An alternate approach to constructing the list of Fourier transforms
involves degenerating the canonical gauge of \cite{generic}.  In general,
any surface of type $F_1$ such that $s$ is not a component of $Q$ has such
a canonical gauge, which (up to constants) is given by gauging so that
every operator of degree $s$ annihilates $1$, while the operators of degree
$f$ remain rational functions on $Q$.  In particular, if we start with a
surface of type $F_0$ and blow up a point which is not on a component of
$Q$ of class $s$ or $f$, then either way of blowing down to $F_1$ gives a
surface with such a canonical gauge.  This makes the limits easier to take
(there is no longer a need to consider a change of gauge), at the cost of
both splitting some of the cases (e.g., at the top level, there is a
subcase in which the point being blown up is a singular point of $Q$) and
not applying to the four cases in which every component of $Q$ has class
$s$ or $f$.  There is a further variant of this approach that applies to
these cases as well: simply blow up three points in such a way that there
is no component of $Q$ having negative intersection with $s+f-e_1-e_2-e_3$,
and gauge so that the corresponding operator annihilates $1$.  Each of
these approaches is analogous to the normalization of equations via a
choice of section, except that in order to normalize morphisms of line
bundles, we must choose a normalizing sheaf for each line bundle, which in
the above approach we do by specifying the Chern class and all but one of
the points of intersection with $Q$.

\chapter{A new Lax pair for Painlev\'e II}
\label{app:PII}

As mentioned above, the derived equivalences of deformed elliptic surfaces
give rise to minimal Lax pairs associated to any level of the hierarchy of
Figure \ref{fig:sakai_good} and any rational number $d/r\in \P^1(\Q)$,
in the form of a (possibly discrete) connection on a vector bundle of rank
$2r$ and degree $d$, with all irreducible singularities appearing with
multiplicity $r$.  The discrete isomonodromy transformations correspond to
iterated discrete Painlev\'e equations, so have expected complexity growing
quadratically in $r$, but the continuous isomonodromy transformations are
flows along the same family of quasi-projective varieties as in the usual
Painlev\'e case, so one expects them to actually be the usual continuous
Painlev\'e equations.  We gave above a prescription for finding
straight-line forms of the linear problems for $d=1$, but this makes it
more difficult to obtain the corresponding differential equation in the
continuous deformation parameter.  An alternate approach uses the standard
algebraic description of a vector bundle via transition matrices between
patches: on each patch, the given connection is given by a differential
equation in the usual sense, and the gluing condition is simply that the
trnasition matrices of the vector bundle give gauge transformations between
the various equations.  For vector bundles on $\P^1$, this is particularly
straightforward, as one can take the standard coordinate patches (i.e., the
respective complements of $0$ and $\infty$) and any vector bundle is
represented by diagonal transition matrices in which the entries are powers
of $z$.

Now, let us specialize to the case of PII, or more precisely the version of
PII in which the standard linear problem is regular on $\C$.  In the
simplest new case $d/r=1/2$, we are looking for a Lax pair of the form
$v_z=A v$, $v_t=C v$ satisfying compatibility along with the following
conditions:

\begin{itemize}
  \item At $\infty$, A has irreducible singularities of the form $f'/f\in
    \{az^2+bt+\gamma/z+o(1/z),-az^2-bt-\delta/z+o(1/z)\}$ (leaving some scaling
    freedom to let us get the correct form of PII).  This essentially
    corresponds to setting $(A-1/4z)^2+O(1)$, given that the second
    condition forces $\Tr(A)=1/z+O(1)$ near 0.
  \item At $0$, if we gauge by the diagonal matrix with entries
    $(1,1,1,1/z)$, then $A$ and $C$ are both holomorphic on $\C$.
\end{itemize}

(For general $d/r$, WLOG with $0\le d<2r$, the first condition becomes that
$(A-\Tr(A)/4)^2+O(1)$ is an appropriate multiple of the identity, while the
second condition involves gauging by the diagonal matrix with $d$ copies of
$1/z$.)

The first condition ensures that $A=O(z^2)$, and determines the leading
terms of the coefficients up to changes of basis compatible with the second
condition.  Using the remaining freedom coming from the automorphism group
of the vector bundle to eliminate degrees of freedom, we end up with two
degrees of freedom for $A$.  We may furthermore identify a good candidate
for a standard Painlev\'e transcendent by observing that $v_4(z)/z$
satisfies a straight-line equation with a single finite singularity; we
take $x(t)$ to be the negative of the location of that singularity.  The
compatibility condition of the Lax pair then gives a differential equation
allowing us to express the other coordinate on the moduli space in terms of
$y(t)=x'(t)$, and giving us a final form for the Lax pair (after setting
$a,b,c,c'$ to make things match up with the form of PII appearing in the
introduction):
\begin{align}
A
=&
\hphantom{{}+{}}
\begin{pmatrix}
  0&0&1/2&0\\
  0&0&0&1/2\\
  1/2&0&0&0\\
  0&1/2&0&0
\end{pmatrix}
z^2\notag\\
&+
\begin{pmatrix}
  0 & 0 & -x/2 & 0\\
  -1 & 0 & 0 & 0\\
  x/2 & 0 & 0 & 0\\
  0 & 0 & 1 & 0
\end{pmatrix}
z\notag\\
&+
\begin{pmatrix}
  y/2 & w & x^2/2+t/2 & 1/4\\
  -x & -y/2 & 0 & t/2\\
  0 & -1/4 & -y/2 & -w\\
  0 & 0 & 0 & y/2
\end{pmatrix}\notag\\
&+
\begin{pmatrix}
  0 & 0 & 0 & -x/4\\
  0 & 0 & 0 & 2w -\beta+1/2\\
  0 & 0 & 0 & -xw\\
  0 & 0 & 0 & 1
\end{pmatrix}
z^{-1},
\end{align}
where $w=x^3/4+tx/4+\beta/2-1/4$, and
\begin{align}
C
=&
\hphantom{{}+{}}
\begin{pmatrix}
  0&0&1/4&0\\
  0&0&0&1/4\\
  1/4&0&0&0\\
  0&1/4&0&0
\end{pmatrix}
z\notag\\
&+
\begin{pmatrix}
  0 & -5x^2/8 - t/8 & -3x/4 & - y/4\\
  -1 & 0 & 0 & 0\\
  -x/4 & -y/4 & 0 & -5x^2/8 - t/8\\
  0 & 0 & 0 & 0
\end{pmatrix}\notag\\
&+
\begin{pmatrix}
  0 & 0 & 0 & xy/4+1/4\\
  0 & 0 & 0 & 5x^2/4+t/2\\
  0 & 0 & 0 & -7x^3/8 - 3tx/8 - \beta/2+1/4\\
  0 & 0 & 0 & 0
\end{pmatrix}
z^{-1}
\end{align}
with compatibility condition
\[
x_{tt} = 2x^3 + tx + \beta-1/2.
\]
If we set $u = x_t + x^2+t/2$, then this becomes the system
\[
x_t = u - x^2-t/2, \qquad u_t = 2xu + \beta,
\]
which is precisely the same form of PII corresponding to the Lax pairs of
the introduction.  Note that the irreducible singularities of $A$ at
$\infty$ are
\[
f'/f\in
\{z^2/2+t/4+(1-\beta)/2z+o(1/z),-z^2/2-t/4+\beta/2z+o(1/z)\},
\]
and more generally the Lax pair near $\infty$ is given by two copies of
\[
\log(f)_z = (z^2/2+t/4+(1-\beta)/2z+o(1/z)),
\qquad \log(f)_t = (z/4+O(1))
\]
and two copies of
\[
\log(f)_z = (-z^2/2-t/4+\beta/2z+o(1/z)),
\qquad \log(f)_t = (-z/4+O(1)).
\]
(This is why $C=A/2z+O(1)$ above.)

By comparison, the analogous rank 2 Lax pair for PII, still on a vector
bundle of degree 1, has
\[
A
=
  \begin{pmatrix} 0 & 1\\1 & 0\end{pmatrix}z^2
  +
  \begin{pmatrix} 0 & -x \\ x & 0\end{pmatrix}z
  +
  \begin{pmatrix} y & x^2+t \\ 0 & -y\end{pmatrix}
  +
  \begin{pmatrix} 0 & -x^3-tx-2\beta-1\\ 0 & 1\end{pmatrix}z^{-1}
\]
and
\[
C
=
\begin{pmatrix} 0 & 1/2\\1/2 & 0\end{pmatrix}z
+  
\begin{pmatrix} 0 & -x \\ 0 & 0\end{pmatrix}
+
\begin{pmatrix} 0 & 3x^2/2+t/2\\ 0 & 0\end{pmatrix}z^{-1}.
\]

It is worth noting, as we mentioned in the introduction, that although the
above theory gives an existence result for such Lax pairs, any explicit
version of such a Lax pair (including the Lax pairs above) gives an {\em
  independent} existence result.  Indeed, in finding the rank 4 Lax pair,
the only thing we used from the noncommutative theory {\em was} the
existence result, and that only as reassurance that the
computation would indeed produce a result; everything else was just
straightforward manipulations in (differential) algebra.

\backmatter
\bibliographystyle{plain}
\cleardoublepage
\phantomsection
\addcontentsline{toc}{chapter}{Bibliography}
\bibliography{book}

\end{document}